%% file: thesis.tex
\pdfoutput=1
\documentclass[b5paper]{book}
\usepackage{amsmath}
\usepackage{amssymb}
\usepackage{latexsym}
\usepackage{index}
\usepackage{graphicx}
\usepackage{fancyhdr}
\usepackage{hyperref}

\makeindex
\newtheorem{thm}{Theorem}[chapter]
\newtheorem{cor}[thm]{Corollary}
\newtheorem{lem}[thm]{Lemma}
\newtheorem{prop}[thm]{Proposition}

\newtheorem{conj}[thm]{Conjecture}
\newtheorem{cond}[thm]{Condition}
\providecommand{\norm}[1]{\left\| #1 \right\|}

\newcommand{\mb}{\mathbf}
\newcommand{\mh}{\mathbb}
\newcommand{\mr}{\mathrm}
\newcommand{\mc}{\mathcal}

\newcommand{\mf}{\mathfrak}

\newcommand{\ds}{\displaystyle}
\newcommand{\scs}{\scriptstyle}

\newcommand{\isom}{\xrightarrow{\;\sim\;}}
\newcommand{\mosi}{\xleftarrow{\;\sim\;}}

\newcommand{\ep}{\epsilon}
\newcommand{\es}{\emptyset}
\newcommand{\hot}{\widehat{\otimes}}
\newcommand{\oot}{\overline{\otimes}}
\newcommand{\inp}[2]{\langle #1 \,,\, #2 \rangle}
\newcommand{\inde}[1]{\index{#1} #1}
\newcommand{\prefix}[3]{\vphantom{#3}#1#2#3}
\newcommand{\hexagon}[6]{\begin{array}{ccccc}
  #1 & \to & #2 & \to & #3 \\
  \uparrow & & & & \downarrow \\
  #4 & \leftarrow & #5 & \leftarrow & #6
  \end{array}}

\begin{document}

\include{preface}

\tableofcontents

\pagestyle{fancy}
\renewcommand{\chaptermark}[1]{\markboth{\chaptername \ \thechapter.\ #1}{}}
\renewcommand{\sectionmark}[1]{\markright{\thesection.\ #1}}
\renewcommand{\headrulewidth}{0.2pt}
\fancyfoot{}
\fancyhead[LE,RO]{\thepage}
\fancyhead[LO]{\nouppercase{\rightmark}}
\fancyhead[RE]{\nouppercase{\leftmark}}

\include{chapter1}

\include{chapter2}

\include{chapter3}

\include{chapter4}

\include{chapter5}

\include{chapter6}

\appendix
\renewcommand{\chaptername}{Appendix}

\include{appendix}

\include{bibliographyb5}
\addcontentsline{toc}{chapter}{Index}
\printindex

\include{summary}

\end{document}

%% file: preface.tex
\pagestyle{empty}

$\quad$
\vspace{5cm}
\begin{center}
{\Large \bf Periodic cyclic homology of
affine Hecke algebras }
\\[3cm]
{\Large Maarten Solleveld }
\end{center}

\newpage

$\quad$
\vspace{14cm}

\noindent Periodic cyclic homology of affine Hecke algebras 
/ Maarten Solleveld, 2007 -\\ 253 p. : fig. ; 24 cm. -
Proefschrift Universiteit van Amsterdam -\\
Met samenvatting in het Nederlands.\\
ISBN 978-90-9021543-3
\\[5mm]
Cover designed with courtesy of William Wenger

\newpage

\begin{center}
{\LARGE \bf Periodic cyclic homology\\ of \\[3mm]
affine Hecke algebras}
\\[2cm]
{\large ACADEMISCH PROEFSCHRIFT} 
\\[1cm]
ter verkrijging van de graad van doctor \\
aan de Universiteit van Amsterdam \\
op gezag van de Rector Magnificus \\
prof. dr. J.W. Zwemmer \\
ten overstaan van een door het \\
college voor promoties ingestelde commissie, \\
in het openbaar te verdedigen in de Aula der Universiteit \\
op dinsdag 6 maart 2007, te 12:00 uur
\\[1cm]
door
\\[1cm]
{\large \bf Maarten Sander Solleveld}
\\[5mm]
geboren te Amsterdam
\end{center}

\newpage

\noindent\textbf{Promotiecommissie:}
\\[3mm]
\begin{tabular}[t]{l@{\qquad \qquad}l}
Promotor: & Prof. dr. E.M. Opdam \\
 & \\
Co-promotor: & Prof. dr. N.P. Landsman\\
 & \\
Overige leden: & Prof. dr. G.J. Heckman \\
 & Prof. dr. T.H. Koornwinder \\
 & Prof. dr. R. Meyer \\
 & Prof. dr. V. Nistor \\
 & Dr. M. Crainic \\
 & Dr. H.G.J. Pijls \\
 & Dr. J.V. Stokman
\end{tabular}
\\[3mm]
Faculteit der Natuurwetenschappen, Wiskunde en Informatica
\\[5cm]
Dit proefschrift werd mede mogelijk gemaakt door de Nederlandse
Organisatie voor Wetenschappelijk Onderzoek (NWO). 
Het onderzoek vond plaats in het kader van het NWO Pionier project
"Symmetry in Mathematics and Mathematical Physics".
\\[1cm]
\hspace{1cm}
\includegraphics[width=39mm,height=20mm]{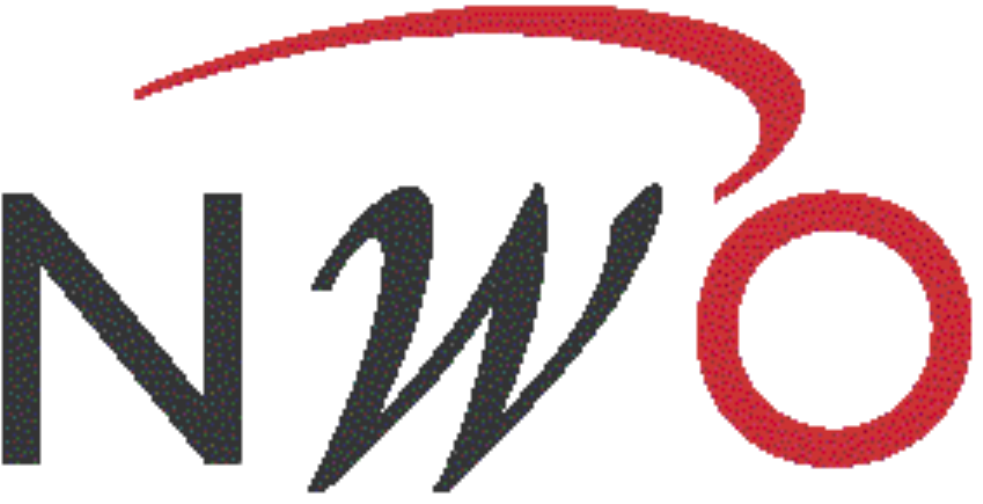}
\hspace{2cm}
\includegraphics[width=62mm,height=26mm]{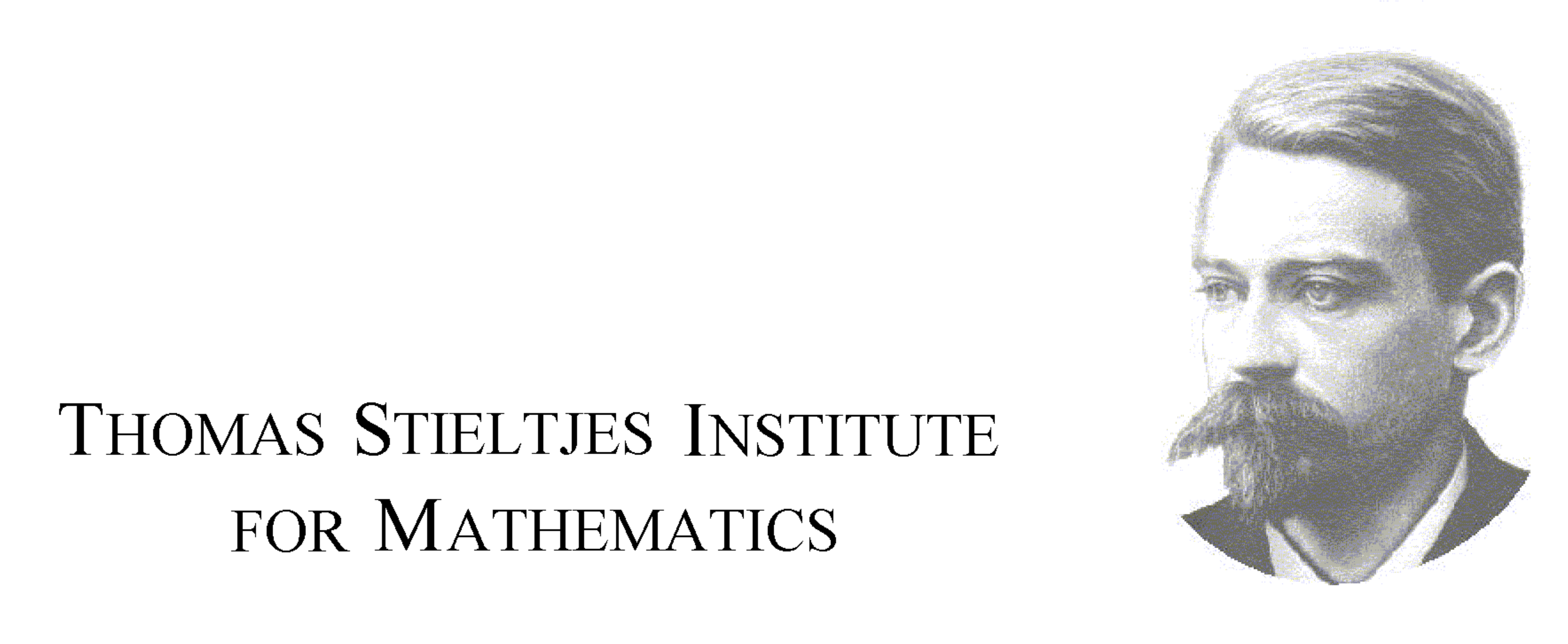}

\cleardoublepage

\chapter*{Preface}
\addcontentsline{toc}{chapter}{Preface}

The book you've just opened is the result of four years of research 
at the Korteweg-de Vries Institute for Mathematics of the Universiteit
van Amsterdam. With this thesis I hope to obtain the degree of doctor.

Due to the highly specialized nature of my research, quite some 
mathematical background is required to read this book. If you want to
get an idea of what it's about, then you may find it useful to start
with the Dutch summary. 

I would like to use this opportunity to express my gratitude to all 
those who supported me in this work, intentionally or not.
\\[1mm]

First of course my advisor, Eric Opdam, without whom almost my entire 
research would have been impossible. Over the years we had many pleasant
conversations, not only about mathematical issues, but also about chess,
movies, education, Japan, and many other things. The intensity of our 
contact fluctuated a lot. Sometimes we didn't talk for weeks, while in
other periods we spoke many hours, discussing our new findings every day.

Vividly I recall one particular evening. After a prolonged discussion I
had finally returned home and was preparing dinner. Then unexpectedly 
Eric called to tell me about some further calculations. Watching the 
boiling rice with one eye and trying to visualize an affine Hecke algebra
with the other, I quickly decided that I had to opt for dinner this time.
Nevertheless that phonecall had a profound influence on Section 6.6. I
admire Eric's deep mathematical insight, with which he managed to put me
on the right track quite a few times.
\\[1mm] 

Also I would like  to thank all the members of the promotion committee 
for the time and effort they made to read the manuscript carefully. Niels 
Kowalzig was so kind to read and comment on the second chapter.
For the summary I am indebted to Klaas Slooten. His thesis was a source 
of inspiration, even though I didn't use many theorems from it.

The KdVI is such a nice place that I come there for more than nine years 
already. Everybody who worked there is responsible for that, but in
particular Erdal, Fokko, Geertje, Harmen, Misja, Peter, Rogier and Simon.
\\[1mm]

Both my ushers, Mariska Berthol\'ee-de Mie and Ionica Smeets, are very 
dear to me. For being such good friends, for having faith in me, for 
lending me an ear, for cheering me up when things did not go as I would 
have liked. I am very happy that they will support me during the defence 
of this thesis.
\\[1mm]

Furthermore I thank all my friends at US badminton for lots of (sporting)
pleasure. Especially Paul den Hertog, who also took care of printing this
book. 

To Karel van der Weide I am grateful for sharing his laconic yet hilarious 
views on the chess world, on internetdating and on life in general.

Bill Wenger was very generous in granting me permission to use his artwork
on the cover of this book.
\\[1mm]

But above all I thank Lieske Tibbe, for being my mother, and 
everything that naturally comes with that. If there is anybody who gave
me the right scientific attitude to complete this thesis, it's her.
\\[5mm]

Amsterdam, January 2007

%% file: chapter1.tex
\chapter{Introduction}

This thesis is about a very interesting kind of algebras, Hecke 
algebras. They appear in various fields of mathematics, for 
example knot theory, harmonic analysis, special functions and 
noncommutative geometry. The motivation for the research presented
here lies mainly in the harmonic analysis of reductive $p$-adic
groups. The description and classification of smooth representations
of such groups is a long standing problem. This is motivated by
number theoretic investigations.

The category of smooth representations is divided in certain blocks
called Bernstein components. It is known that in many cases Bernstein
components can be described with the representation theory of 
certain affine Hecke algebras. This translation is a step forward,
since in a sense affine Hecke algebras are much smaller and easier 
to handle than reductive $p$-adic groups.

Hence it is desirable to obtain a good description of all irreducible
representations of an affine Hecke algebra. Such an algebra $\mc H 
(W,q)$ can be considered as a deformation of the (complex) group 
algebra of an affine Weyl group $W$, which involves a few parameters 
$q_i \in \mh C^\times$. Let us briefly mention what is known about 
the classification of its spectrum in various cases.

\begin{description}\label{p:1.1}
\item[1)] All parameters $q_i$ equal to 1.\\
In this special case $\mc H (W,q)$ is just the group algebra of $W$. 
The representations of groups like $W$ have been known explicitly 
for a long time, already from the work of Clifford \cite{Cli}.

\item[2)] All parameters $q_i$ equal to the same complex number,
not a root of unity.\\
With the use of equivariant $K$-theory, Kazhdan and Lusztig 
\cite{KaLu} gave a complete classification of the irreducible 
representations of $\mc H (W,q)$. It turns out that they are in 
bijection with the irreducible representations of $W$. This 
bijection can be made explicit with Lusztig's asymptotic Hecke 
algebra \cite{Lus2}.

\item[3)] Exceptional cases.\\
These may occur for example if there exist integers $n_i$ such 
that $\prod_i q_i^{n_i /2} \neq 1$ is a root of unity, cf. page
\pageref{eq:5.2}. The affine Hecke algebras for such parameters 
may differ significantly from those in the other cases, so we 
will not study them here.

\item[4)] General positive parameters.\\
Quite strong partial classifications are available, mainly from 
the work of Delorme and Opdam \cite{DeOp1,DeOp2}.

\item[5)] General unequal parameters.\\
These algebras have been studied in particular by Lusztig 
\cite{Lus6}. Recently Kato \cite{Kat4} parametrized the 
representations of certain affine Hecke algebras with three 
independent parameters, extending ideas that were used in case 2).
\end{description}

The affine Hecke algebras that arise from reductive $p$-adic groups
have rather special parameters: they are all powers of the cardinality
of the residue field of the $p$-adic field. Sometimes they are all
equal, and sometimes they are not, so these algebras are in case 4).

Lusztig \cite[Chapter 14]{Lus6} conjectured that 

\begin{conj}\label{conj:1.1}
For general unequal parameters there also exists an asymptotic 
Hecke algebra. It yields a natural bijection between the 
irreducible representations of $\mc H (W,q)$ and those of $W$.
\end{conj}

We are mainly concerned with a somewhat weaker version:

\begin{conj}\label{conj:1.4}
For positive parameters there is an isomorphism between the 
Grothendieck groups of finite dimensional $\mc H (W,q)$-modules and 
of finite dimensional $W$-modules.
\end{conj}

In principle the verification of this conjecture would involve two
steps 

\begin{description}
\item[a)] Assign a $W$-representation to each (irreducible) $\mc H 
(W,q)$-representation, in some natural way (or the other way round).

\item[b)] Prove that this induces an isomorphism on the above 
Grothendieck groups.
\end{description}

In our study we make use of a technique that is obviously not 
available for $p$-adic groups, we deform the parameters continuously.
We would like to do this in the context of topological algebras, 
preferably operator algebras. For this reason, and to avoid the 
exceptional cases 3), we assume throughout
most of the book that all $q_i$ are positive real numbers. It was
shown in \cite{Opd3} that for such parameters there is a nice 
natural way to complete $\mc H (W,q)$ to a Schwartz algebra $\mc S 
(W,q)$ (these notations are preliminary). We will compare this algebra 
to the Schwartz algebra $\mc S (W)$ of the group $W$. Using the 
explicit description of $\mc S (W,q)$ given by Delorme and Opdam 
\cite{DeOp1}, in Section \ref{sec:5.3} we construct a Fr\'echet 
algebra homomorphism
\begin{equation}\label{eq:1.1}
\phi_0 : \mc S (W) \to \mc S (W,q)
\end{equation}
with good properties. This provides a map from $\mc S (W,q)
$-representations to $\mc S (W)$-representations. Together with 
the Langlands classification \cite[Section 6]{DeOp2} for 
representations of $\mc H (W,q)$ and of $W$ this takes care of a).

To reformulate b) in more manageable terms we turn to noncommutative
geometry. There are (at least) three functors which are suited to 
deal with such problems: periodic cyclic homology $(HP)$, in either
the algebraic or the topological sense, and topological $K$-theory.

\begin{conj}\label{conj:1.2}
There are natural isomorphisms
\begin{enumerate}
\item $HP_* (\mc H (W,q)) \cong HP_* (\mh C [W])$
\item $HP_* (\mc S (W,q)) \cong HP_* (\mc S (W))$
\item $\;\;\: K_* (\mc S (W,q)) \cong \;\; K_* (\mc S (W))$
\end{enumerate}
\end{conj}
The relation with Conjecture \ref{conj:1.4} is as follows. Although
$HP_*$ and $K_*$ are functors on general noncommutative algebras,
in our setting they depend essentially only on the spectra of the
algebras that we are interested in. These spectra are ill-behaved
as topological spaces: the spectrum $\widehat{\mc H (W,q)}$ is a 
non-separated algebraic variety, while $\widehat{\mc S (W,q)}$ is a
compact non-Hausdorff space. Nevertheless topological $K$-theory 
and periodic cyclic homology can be considered as cohomology 
theories on such spaces. With this interpretation Conjecture 
\ref{conj:1.2} asserts that the "cohomology groups" of 
$\widehat{\mc H (W,q)}$ and of $\widehat{\mc S (W,q)}$ are invariant 
under the deformations in the parameters $q_i$. Contrarily to what 
one would expect from the results on page \pageref{p:1.1}, from 
this noncommutative geometric point of view the algebras $\mc S 
(W,q)$ actually become easier to understand when the parameters 
$q_i$ have less relations among themselves.

For equal parameters Conjecture \ref{conj:1.2} has been 
around for a while. Part 3 already appeared in the important paper
\cite{BCH}, while part 1 was proven by Baum and Nistor \cite{BaNi}.
The proof relies on the aforementioned results of Kazhdan and 
Lusztig. In the unequal parameter case Conjecture \ref{conj:1.2}.3
was formulated independently by Opdam \cite[Section 1.0.1]{Opd3}.

In this thesis we make the following progress concerning these 
conjectures. In Section \ref{sec:3.4} we prove that there are 
natural isomorphisms 
\begin{equation}\label{eq:1.2}
HP_* (\mc H (W,q)) \cong HP_* (\mc S (W,q)) \cong 
K_* (\mc S (W,q)) \otimes_{\mh Z} \mh C
\end{equation}
Hence parts 1 and 2 of Conjecture \ref{conj:1.2} are equivalent,
and both are weaker than part 3. Moreover in Section \ref{sec:5.4} 
we show that the Conjectures \ref{conj:1.4} and \ref{conj:1.2}.3 
are equivalent.

Guided by these considerations we propose the following 
refined conjecture, which has also been presented by 
Opdam \cite[Section 7.3]{Opd4}:

\begin{conj}\label{conj:1.3}
The natural map
\[
K_* (\phi_0 ) : K_* (\mc S (W)) \to K_* (\mc S (W,q))
\]
is an isomorphism for all positive parameters $q_i$.
\end{conj}

Now let us give a brief overview of this book. More detailed
outlines can be found at the beginning of each chapter.

Chapter 2 deals with noncommutative geometry. This chapter does
not depend on the rest of the book, we do not mention any Hecke
algebras. We provide a solid
foundations for the philosophy that $K$-theory and periodic cyclic
homology can be considered as cohomology theories for certain 
non-Hausdorff spaces. Among others we prove comparison theorems 
like \eqref{eq:1.2} for more general classes of algebras, all
derived from so-called finite type algebras \cite{KNS}. 

In the next chapter we introduce affine Hecke algebras. A large
part of the material presented here relies on the work of Opdam,
in collaboration with Delorme, Heckman, Reeder and Slooten.
We study the representation theory of affine Hecke algebras, which
provides a clear picture of their spectra as topological spaces.
We are especially interested in the image of $\mc S (W,q)$ under
the Fourier transform, as this turns out to be an algebra of the
type that we studied in Chapter 2. We conclude with the above
isomorphisms \eqref{eq:1.2}.

These are also interesting because they can be generalized to 
algebras associated with reductive $p$-adic groups, which we will
do in Chapter 4. Given a reductive $p$-adic group $G$ we recall 
the constructions of its Hecke algebra $\mc H (G)$, its Schwartz
algebra $\mc S (G)$ and its reduced $C^*$-algebra $C_r^* (G)$.
The main new results in this chapter are natural isomorphisms
\begin{equation}\label{eq:1.4}
HP_* (\mc H (G)) \cong HP_* (\mc S (G)) \cong
K_* (C_r^* (G)) \otimes_{\mh Z} \mh C
\end{equation}
These have some consequences in relation with the Baum-Connes
conjecture for $G$.

In Chapter 5 we really delve into the study of deformations of
affine Hecke algebras. The Fr\'echet space underlying $\mc S 
(W,q)$ is independent of $q$, and we show that all the (topological) 
algebra operations in $\mc S (W,q)$ depend continously on $q$. After 
that we focus on parameter deformations of the form $q \to q^\ep$
with $\ep \in [0,1]$. For positive $\ep$ we construct isomorphisms
\begin{equation}\label{eq:1.5}
\phi_\ep : \mc S (W,q^\ep ) \to \mc S (W,q)
\end{equation}
that depend continuously on $\ep$. The limit $\lim_{e \downarrow 
0} \phi_\ep$ is well-defined and indeed is \eqref{eq:1.1}. 
Furthermore we elaborate on the conjectures mentioned in this 
introduction.

In support of the these conjectures, and to show what the 
techniques we developed are up to, we dedicate Chapter 6 to 
calculations for affine Hecke algebras of classical type. We verify
Conjecture \ref{conj:1.3} in some low-dimensional cases and for 
types $GL_n$ and $A_n$. 

We conclude the purely scientific part of the book with a short
appendix. It contains some rather elementary results on crossed
product algebras that are used at various places.

%% file: chapter2.tex
\chapter{$K$-theory and cyclic type homology theories}

This chapter is of a more general nature than the rest of
this book. We start with the study of some important covariant
functors on the category of complex algebras. These are 
Hochschild homology, cyclic homology and periodic cyclic homology.
Contravariant versions of these functors also exist, but we will 
leave these aside. All these functors go together by the name of
cyclic theory.

It is well-known that cyclic homology is related to $K$-theory
by a natural transformation of functors called the Chern 
character. We are not satisfied with $K$-theory for Banach
algebras, but instead study its extension to the larger 
categories of Fr\'echet algebras or even m-algebras. From
these abstract considerations we will see that there are three
functors which share almost identical properties:
\begin{description}
\item[a)] periodic cyclic homology, purely algebraically
\item[b)] periodic cyclic homology, with the completed 
projective tensor product
\item[c)] $K$-theory for Fr\'echet algebras
\end{description}
These functors can be regarded as noncommutative analogues of
\begin{description}
\item[1)] De Rham cohomology in the algebraic sense, for 
complex affine varieties
\item[2)] De Rham cohomology in the differential geometric 
sense, for smooth manifolds
\item[3)] $K$-theory for topological spaces
\end{description}
By a comparison theorem of Deligne and Grothendieck 1) and 2) 
agree for a complex affine variety. For smooth manifolds 2) and 
3) (with real coefficients) give the same result, essentially
because both are generalized cohomology theories. This is also
the reason that both can be computed as
\begin{description}
\item[4)] \v Cech cohomology of a constant sheaf
\end{description}
A noncommutative analogue of 4) does not appear to exist, so we 
develop it. It will be a sheaf that depends only on the spectrum
of an algebra. Then we can also consider
\begin{description}
\item[d)] \v Cech cohomology of this sheaf
\end{description}
The main goal of this chapter is to generalize the isomorphisms
between 1) - 4) to the setting of noncommutative algebras. So
far this has been done only for b) and c).

Let us also give a more concrete overview of this chapter.
We start by recalling the definitions and properties of cyclic
theory in the purely algebraic sense. Then we specialize to
finite type algebras, mainly following \cite{KNS}. For such
algebras we define a sheaf which provides the isomorphism 
between a) and d).

After that we move on to topological algebras, especially 
Fr\'echet algebras. Most of the properties of algebraic cyclic
theory have been carried over to this topological setting, but
unfortunately these results have been scattered throughout the 
literature. We hope that bringing them together will serve the
reader. We also recall several results concerning $K$-theory
for Fr\'echet algebras, which are mostly due to Cuntz \cite{Cun1}
and Phillips \cite{Phi2}. 

In the final section we have to decide for what kind of 
topological algebras we want to compare b), c) and d). Natural
candidates are algebras that are finitely generated as modules 
over their center. For finitely generated (non-topological) 
algebras this condition leads to the aforementioned finite type
algebras. Their spectrum has the structure of a non-separated
complex affine variety.

In the topological setting we need to impose more conditions. 
Cyclic homology works best if there is a kind of smooth 
structure, so our topological analogue of a finite type algebra
is of the form 
\begin{equation}\label{eq:2.84}
C^\infty (X ;M_N (\mh C ))^G
\end{equation}
where $X$ is a smooth manifold and $G$ a finite group. The action
of $G$ is a combination of an action on $X$ and conjugation by
certain matrices.

The comparisons between 1) - 4) all rely on triangulations and 
Mayer-Vietoris sequences. We will apply these techniques to $X$ in
a suitable way. This will enable us to define d) and prove that 
it gives the same results as b) and c). Moreover we prove that if 
$X$ happens to be a complex affine variety, then there is a natural
isomorphism
\begin{equation}
HP_* \big( \mc O (X ;M_N (\mh C ))^G \big) \isom
HP_* \big( C^\infty (X ;M_N (\mh C ))^G \big)
\end{equation}
The $\mc O$ stands for algebraic functions, so the left hand side
corresponds to a) and 1) above.
\\[4mm]

\section{Algebraic cyclic theory}
\label{sec:2.1}

We give the definitions and most important properties of 
Hochschild homology, cyclic homology and periodic cyclic homology. 
We do this only in the category of algebras over $\mh C$, although
quite a big deal of the theory is also valid over arbitrary fields.
We will be rather concise, referring to Loday's monograph 
\cite{Lod} for more background and proofs.

Cyclic homology was discovered more or less independently by
several people, confer the work of Connes \cite{Con}, Loday and
Quillen \cite{LoQu} and Tsygan \cite{Tsy}. We will define it with
the so-called cyclic bicomplex. So let $n \in \mh N$, let $A$ be 
any complex algebra and let $A^{\otimes n}$ be the $n$-fold
tensor product of $A$ over $\mh C$. We abbreviate an elementary
tensor $a_1 \otimes \cdots \otimes a_n$ to
$(a_1, \ldots, a_n)$ and we define linear operators
$b, b' : A^{\otimes n+1} \to A^{\otimes n}$ and 
$\lambda, N : A^{\otimes n} \to A^{\otimes n}$ by the following 
formulas:
\begin{equation}\label{eq:2.1}
\begin{aligned}
& b' (a_0,a_1,\ldots,a_n) & =\;\; & \sum_{i=0}^{n-1} (-1)^i 
  (a_0,\ldots,a_{i-1},a_i a_{i+1},a_{i+2},\ldots,a_n) \\
& b (a_0,a_1,\ldots,a_n) & =\;\; & b' (a_0,a_1,\ldots,a_n) + 
  (-1)^n (a_n a_0,a_1,\ldots, a_{n-1}) \\
& \lambda (a_1,\ldots,a_n) & =\;\; & 
  (-1)^{n-1} (a_n,a_1,\ldots,a_{n-1}) \\
& N & =\;\; & 1 + \lambda + \cdots + \lambda^{n-1}
\end{aligned}
\end{equation}
For unital $A$ we also define 
$s,B : A^{\otimes n} \to A^{\otimes n+1}$:
\begin{equation}\label{eq:2.2}
\begin{aligned}
&s (a_1,\ldots,a_n) & =\;\; & (1,a_1,\ldots,a_n) \\
&B & =\;\; & (1-\lambda) s N
\end{aligned}
\end{equation}
These are the ingredients of a bicomplex $CC^{per}(A)$ :
\index{cyclic bicomplex}
\begin{equation}
\begin{array}{c@{\qquad\qquad}ccccccccccc}
 & & & \cdots & & \cdots & & \cdots & & \cdots & & \\
 & & & \downarrow & & \downarrow & & \downarrow & & \downarrow & & 
  \\[2mm]
{\scs 2} & \cdots & \leftarrow & A^{\otimes 3} & \xleftarrow{N} & 
  A^{\otimes 3} & \xleftarrow{1-\lambda} & A^{\otimes 3} & 
  \xleftarrow{N} & A^{\otimes 3} & \leftarrow & \cdots \\[2mm]
 & & & \; \downarrow \! -b' & & \downarrow b & & 
  \; \downarrow \! -b' & & \downarrow b & & \\[2mm]
{\scs 1} & \cdots & \leftarrow & A^{\otimes 2} & \xleftarrow{N} & 
  A^{\otimes 2} & \xleftarrow{ 1-\lambda} & A^{\otimes 2} & 
  \xleftarrow{N} & A^{\otimes 2} & \leftarrow & \cdots \\[2mm]
 & & & \; \downarrow \! -b' & & \downarrow b & & 
  \; \downarrow \! -b' & & \downarrow b & & \\[2mm]
{\scs 0} & \cdots & \leftarrow & A & \xleftarrow{N} & A & 
  \xleftarrow{1-\lambda} & A & \xleftarrow{N} & A & 
  \leftarrow & \cdots \\[7mm]
 & & & {\scs -1} & & {\scs 0} & & {\scs 1} & & {\scs 2} & &
\end{array}
\end{equation}
The indicated grading means that 
$CC^{per}_{p,q}(A) = A^{\otimes q+1}$. 

Consider also the subcomplexes $CC(A)$, consisting of all the 
columns labelled by $p \geq 0$, and $CC^{\{2\}}(A)$, which 
consists only of the columns numbered 0 and 1.
With these bicomplexes we associate differential complexes with 
a single grading. Their spaces in degree $n$ are
\begin{equation}
\begin{aligned}
&CC_n^{\{2\}}(A) &=\;\;& CC_{0,n}^{per}(A) \oplus CC_{1,n-1}^{
  per}(A) &=\;\;&  A^{\otimes n+1} \oplus A^{\otimes n} \\
&CC_n (A) &=\;\;& \bigoplus_{p=0}^n CC_{p,n-p}^{per}(A) &=\;\;&
  A^{\otimes n+1} \oplus A^{\otimes n} \oplus \cdots \oplus A \\
&CC_n^{per}(A) &=\;\;& \prod_{p+q = n} CC_{p,q}^{per}(A) &=\;\;&
  \prod_{q \geq 0} A^{\otimes q+1}
\end{aligned}
\end{equation}
This enables us to define the \inde{Hochschild homology} 
$HH_n (A)$, the \inde{cyclic homology} $HC_n (A)$ and the 
\inde{periodic cyclic homology} $HP_n (A)$ :
\begin{equation}
\begin{aligned}
HH_n (A) &= H_n ( CC_*^{\{2\}}(A) ) \\
HC_n (A) &= H_n ( CC_* (A) ) \\
HP_n (A) &= H_n ( CC_*^{per}(A) )
\end{aligned}
\end{equation}
Since all the above complexes are functorial in $A$ these 
homology theories are indeed covariant functors. The 
definitions we gave are neither the simplest possible ones, 
nor the best for explicit computations, but they do have the 
advantage that they work for every algebra, unital or not. 

By the way, we can always form the \inde{unitization} $A^+$. 
This is the vector space $\mh C \oplus A$ with multiplication
\begin{equation}\label{eq:2.6}
(z_1,a_1) (z_2,a_2) = (z_1 z_2, z_1 a_2 + z_2 a_1 + a_1 a_2)
\end{equation}
Clearly every algebra morphism $\phi : A \to B$ gives a unital
algebra morphism $\phi^+ : A^+ \to B^+$.
There are natural isomorphisms
\begin{equation}\label{eq:2.3}
HH_n (A) \cong \mr{coker} \left( HH_n (\mh C) \to HH_n (A^+) 
\right) \cong \ker \left( HH_n (A^+) \to HH_n (\mh C) \right)
\end{equation}
and similarly for $HC_n$ and $HP_n$. 

Often we shall want to consider all degrees at the same time, 
and for this purpose we write
\begin{align*}
HH_* (A) &= \bigoplus_{n \geq 0} HH_n (A) \\
HC_* (A) &= \bigoplus_{n \geq 0} HC_n (A)
\end{align*}
The map $S: CC_{p,q}^{per}(A) \to CC_{p-2,q}^{per}(A)$
simply shifting everything two columns to the left is clearly an
automorphism of $CC^{per}(A)$. Moreover it decreases the degree
by two, so it induces a natural isomorphism
\begin{equation}
HP_n (A) \isom HP_{n-2}(A)
\end{equation}
Thus we may consider periodic cyclic homology as a $\mh Z / 2 
\mh Z$-graded functor, or we may restrict $n$ to $\{0,1\}$.
In particular we shall write
\[
HP_* (A) = HP_0 (A) \oplus HP_1 (A)
\]
Regarding $CC (A)$ as a quotient of $CC^{per}(A)$, we get an
induced map $\bar S : CC(A) \to CC(A)$. This map is surjective, 
and its kernel is exactly $CC^{\{2\}}(A)$.
This leads to Connes' \inde{periodicity exact sequence} :
\begin{equation}\label{eq:2.4}
\cdots \to HH_n (A) \xrightarrow{I} HC_n (A) \xrightarrow{S} 
  HC_{n-2}(A) \xrightarrow{B} HH_{n-1}(A) \to \cdots
\end{equation}
Here $I$ comes from the inclusion of $CC^{\{2\}}(A)$ in $CC(A)$
and $B$ is induced by the map from \eqref{eq:2.2}.
In combination with the five lemma this is a very useful tool;
it enables one to prove easily that many functorial properties
of Hochschild homology also hold for cyclic homology.

Furthermore we notice that the bicomplex $CC^{per}(A)$ is the 
inverse limit of its subcomplexes $S^r (CC(A))$. In many cases
this gives an isomorphism between $HP_n (A)$ and $\varprojlim 
HC_{n+2r} (A)$. In general however it only leads to a short exact 
sequence
\begin{equation}\label{eq:2.5}
0 \to \sideset{}{^1}\lim_{\infty \leftarrow r} HC_{n+1+2r}(A) \to
HP_n (A) \to \lim_{\infty \leftarrow r} HC_{n+2r}(A) \to 0
\end{equation}
Here $\varprojlim^1$ is the first derived functor of 
$\varprojlim$, see \cite[Propostion 5.1.9]{Lod}.

Next we state some well-known features of the functors under 
consideration. 

\begin{enumerate}
\item Additivity. If $A_m \, (m \in \mh N)$ are algebras then
\index{additivity}
\begin{align*}
HH_n \left( \bigoplus_{m=1}^\infty A_m \right) &\cong
\bigoplus_{m=1}^\infty HH_n (A_m) \\
HH_n \left( \prod_{m=1}^\infty A_m \right) &\cong
\prod_{m=1}^\infty HH_n (A_m)
\end{align*}
and similarly for $HC_n$ and $\big( HP_n , \prod \big)$.

\item Stability. If $A$ is H-unital then
\index{stability}
\[
HH_n (M_m (A)) \cong HH_n (A)
\]
More generally, if $B$ and $C$ are unital and Morita-equivalent,
then
\[
HH_n (B) \cong HH_n (C)
\]
These statements hold also for $HC_n$ and $HP_n$.

\item Continuity. If $A = \lim\limits_{m \to \infty} A_m$ is an
inductive limit then
\index{continuity}
\begin{align*}
HH_n (A) &\cong \lim_{m \to \infty} HH_n (A_m) \\
HC_n (A) &\cong \lim_{m \to \infty} HC_n (A_m)
\end{align*}
However, $HP_*$ is not continuous in general. A sufficient 
condition for continuity can be found in \cite[Theorem 3]{BrPl1} :
there exists a $N \in \mh N$ such that $HH_n (A_m) = 0\;\;
\forall n > N \, \forall m$.
\end{enumerate}

To an extension of algebras
\[
0 \to A \to B \to C \to 0
\]
we would like to associate long exact sequences of homology 
groups. This is no problem for unital algebras, but in general 
it is not always possible. However we can consider unitizations 
instead, see \eqref{eq:2.6}. The sequence
\[
0 \to A \to B^+ \to C^+ \to 0
\]
is still exact, so by \eqref{eq:2.3} we may, without loss of 
generality, assume that $B$ and $C$ are unital. It was discovered 
by Wodzicki that what we need for $A$ is not so much unitality, 
but a weaker notion called homological unitality, or H-unitality 
for short. It is easily seen that for unital algebras the map $s$ 
defines a contracting homotopy for the complex 
$(A^{\otimes n},b')$, and in fact with some slight modifications 
this construction also applies to algebras that have left or 
right local units. Thus, we call a complex algebra $A$ 
\inde{H-unital} if the homology of the complex 
$(A^{\otimes n}, b')$ is 0. 
Now Wodzicki's \inde{excision} theorem \cite{Wod} says

\begin{thm}\label{thm:2.1}
Let $0 \to A \to B \to C \to 0$ be an extension of algebras, with 
$A$ H-unital. There exist long exact sequences 
\[
\begin{array}{ccccccccccc}
\cdots & \to & HH_n (A) & \to & HH_n (B) & \to & HH_n (C) & \to & 
HH_{n-1}(A) & \to & \cdots \\
\cdots & \to & HC_n (A) & \to & HC_n (B) & \to & HC_n (C) & \to & 
HC_{n-1}(A) & \to & \cdots \\
\cdots & \to & HP_n (A) & \to & HP_n (B) & \to & HP_n (C) & \to & 
HP_{n-1}(A) & \to & \cdots
\end{array}
\]
\end{thm}

It turns out \cite{CuQu} that for $HP_*$ it is not necessary to 
require H-unitality. Due to the 2-periodicity of this functor we 
get, for any extension of algebras, an exact hexagon
\begin{equation}\label{eq:2.7}
\hexagon{HP_0 (A)}{HP_0 (B)}{HP_0 (C)}{HP_1 (C)}{HP_1 (B)}{HP_1 (A)}
\end{equation}
It will be very useful to combine the excision property with the
\inde{five lemma}:

\begin{lem}\label{lem:2.5}
Suppose we have a commutative diagram of abelian groups, 
with exact rows:
\[
\begin{array}{ccccccccc}
A_1 & \to & A_2 & \to & A_3 & \to & A_4 & \to & A_5 \\[2mm]
\; \downarrow \! f_1 & & \; \downarrow \! f_2 & & 
\; \downarrow \! f_3 & & \; \downarrow \! f_4 & & 
\; \downarrow \! f_5 \\[2mm]
B_1 & \to & B_2 & \to & B_3 & \to & B_4 & \to & B_5
\end{array}
\]
If $f_1$ is surjective, $f_2$ and $f_4$ are 
isomorphisms and $f_5$ is injective, then $f_3$ is an isomorphism.
\end{lem}

Because we intend to apply the next result to several different
functors, we formulate it very abstractly.

\begin{lem}\label{lem:2.12}
Let $\mc A$ and $\mc B$ be categories of algebras, and $\mc{AG}$
the category of abelian groups. Suppose that $F_* : \mc A \to 
\mc{AG}$ and $G_* : \mc B \to \mc{AG}$ are $\mh Z$-graded, 
covariant functors satisfying excision, and that $T_* : F_* \to 
G_*$ is a natural transformation of such functors. 
Consider two sequences of ideals
\begin{equation}\label{eq:2.77}
\begin{array}{ccccccccccc}
0 & = & I_0 & \subset & I_1 & \subset \;\cdots\; \subset & I_n & 
\subset & I_{n+1} & = & A \\
0 & = & J_0 & \subset & J_1 & \subset \;\cdots\; \subset & J_n & 
\subset & J_{n+1} & = & B
\end{array}
\end{equation}
in $\mc A$ and $\mc B$ respectively. If we have an algebra 
homomorphism $\phi : A \to B$ such that $\phi (I_m ) \subset J_m$
and
\[
T(J_m / J_{m+1}) F(\phi ) = G(\phi ) T(I_m / I_{m+1}) : 
F (I_m / I_{m+1}) \to G (J_m / J_{m+1} )
\]
is an isomorphism $\forall m \leq n$, then
\[
T (\phi) := T(B) F(\phi ) = G(\phi ) T(A) : F(A) \to G(B)
\]
is an isomorphism. Similarly, consider two exact sequences
\begin{equation}\label{eq:2.78}
\begin{array}{ccccccccc}
0 & \to & A_1 & \to & A_2 & \to \;\cdots\; \to & A_n & \to & 0 \\
0 & \to & B_1 & \to & B_2 & \to \;\cdots\; \to & B_n & \to & 0
\end{array}
\end{equation}
in $\mc A$ and $\mc B$. Suppose that we have a morphism of exact
sequences $\psi = (\psi_m )_{m=1}^n$, such that 
\[
T (\psi_m) : F(A_m ) \to G (B_m )
\]
is an isomorphism for all but one $m$. Then it is an isomorphism
for all $m$.
\end{lem}
\emph{Proof.}
Consider the short exact sequences
\[
\begin{array}{ccccccccc}
0 & \to & I_{m-1} & \to & I_m & \to & I_m / I_{m-1} & \to & 0 \\
0 & \to & J_{m-1} & \to & J_m & \to & J_m / J_{m-1} & \to & 0 \\
0 & \to & \mr{im}\: (A_{m-1} \to A_m) & \to & A_m & \to & 
  \mr{im}\: (A_m \to A_{m+1}) & \to & 0  \\
0 & \to & \mr{im}\: (B_{m-1} \to B_m) & \to & B_m & \to & 
  \mr{im}\: (B_m \to B_{m+1}) & \to & 0
\end{array}
\]
They degenerate for $m=1$, so with induction we reduce the entire
lemma to the statement for exact sequences, with $m=3$. Now we 
consider only the case where $T(\psi_m )$ is an isomorphism for
$m=1$ and $m=3$, since the other cases are very similar. For any 
$k \in \mh Z$ we see from the commutative diagram
\[
\begin{array}{ccccccccc}
F_{k+1} (A_3 ) & \to & F_k (A_1 ) & \to & F_k (A_2 ) & \to & 
F_k (A_3 ) & \to & F_{k-1} (A_1 ) \\
\downarrow & & \downarrow & & \downarrow & & \downarrow & & 
\downarrow\\
G_{k+1} (B_3 ) & \to & G_k (B_1 ) & \to & G_k (B_2 ) & \to & 
G_k (B_3 ) & \to & G_{k-1} (B_1 ) 
\end{array}
\]
and Lemma \ref{lem:2.5} that $T_k (\psi_2 )$ is an isomorphism. 
$\qquad \Box$ \\[2mm]

Having elaborated a little on the functorial properties of $HH_* 
,\, HC_*$ and $HP_*$, we will show now what they look like on 
some nice algebras. First we fix the notations of some well-known
objects from algebraic geometry.

Assume for the rest of this section that $A$ is a commutative, 
unital complex algebra. The $A$-module of K\"ahler differentials 
$\Omega^1 (A)$ is generated by the symbols $da$, subject
to the following relations, for any $a,b \in A, z \in \mh C$ :
\begin{equation}\label{eq:2.19}
\begin{aligned}
d (za) &= z\, da\\
d (a+b) &= da + db\\
d (ab) &= a\, db + b\, da
\end{aligned}
\end{equation} 
\index{differential forms}
The $A$-module of differential $n$-forms is the $n$-fold exterior
product over $A$ :
\begin{equation}
\Omega^n (A) = \bigwedge\nolimits_A^n \Omega^1 (A)
\end{equation}
and, just to be sure, we decree that $\Omega^0 (A) = A$. 
The formal operator $d$ defines a differential 
$\Omega^n (A) \to \Omega^{n+1} (A)$ by
\begin{equation}
d(a_0 da_1 \wedge \cdots \wedge da_n) = 
da_0 \wedge da_1 \wedge \cdots \wedge da_n
\end{equation}
The \inde{De Rham homology} of $A$ is defined as
\begin{equation}
H^{DR}_n (A) = H_n (\Omega^* (A), d)
\end{equation}
If $A = \mc O(V)$ is the ring of regular functions on an affine 
complex algebraic variety $V$, not necessarily irreducible,  
then we also write
\[
\Omega^n (V) = \Omega^n (A) \quad \mr{and} \quad
H_{DR}^n (V) = H^{DR}_n (A)
\]
One can check that the following formulas define natural maps:
\begin{equation}\label{eq:2.13}
\begin{array}{l@{\;\to\;\;}l@{\;:\quad}l@{\;\to\;\;}l}
HH_n (A) & \Omega^n (A) & (a_0,a_1,\ldots,a_n) & 
a_0 da_1 \wedge \cdots \wedge d a_n \\
\Omega^n (A) & HH_n(A) & a_0 da_1 \wedge \cdots \wedge da_n 
 & \sum\limits_{\sigma \in S_n} \ep(\sigma) 
(a_0,a_{\sigma(1)},\ldots,a_{\sigma(n)})
\end{array}
\end{equation}
The celebrated \inde{Hochschild-Kostant-Rosenberg theorem} 
\cite{HKR} says that these maps are isomorphisms if $A$ is a 
smooth algebra. Yet the author believes that a precise definition 
of smoothness would digress too much, so we only mention that a
typical example is $\mc O(V)$ with $V$ nonsingular, and that all 
the details can be found in \cite[Appendix E]{Lod}. Anyway, 
under \eqref{eq:2.13} the differential $d$ corresponds to 
the map $B$ from \eqref{eq:2.2} and therefore the 
Hochschild-Kostant-Rosenberg theorem also gives
the (periodic) cyclic homology of smooth algebras:
\begin{align}
HH_n (A) &\cong \Omega^n (A)\\
HC_n (A) &\cong \Omega^n (A) / d \Omega^{n-1}(A) \oplus 
  H^{DR}_{n-2}(A) \oplus H^{DR}_{n-4}(A) \oplus \cdots \\
\label{eq:2.8} HP_n (A) &\cong \prod_{m \in \mh Z} 
  H^{DR}_{n+2m}(A)
\end{align}
\\[2mm]

\section{Periodic cyclic homology of finite type algebras}
\label{sec:2.2}

The theory of finite type algebras was built by Baum, Kazhdan, 
Nistor and Schneider \cite{BaNi, KNS}. This turns out to be a pleasant
playground for cyclic theory, culminating roughly speaking in the 
statement ``the periodic cyclic homology of a finite type algebra 
is an invariant of its spectrum.'' We discuss this result, and
some of its background. We also add one new ingredient to support
this point of view, namely a sheaf, depending only on the 
spectrum of $A$, whose \v Cech cohomology is isomorphic to 
$HP_* (A)$. 

All this is made possible by several extra features that $HP_*$
possesses, compared to $HH_*$ and $HC_*$.
Recall that an extension of algebras
\begin{equation}
0 \to I \to A \to A/I \to 0
\end{equation}
is called \inde{nilpotent} if the ideal $I$ is nilpotent, i.e.
if $I^n = 0$ for some $n \in \mh N$.

An algebraic \inde{homotopy} between two algebra homomorphisms 
$f,g : A \to B$ is a collection $\phi_t : A \to B$ of morphisms, 
depending polynomially on $t$, such that $\phi_0 = f$ and 
$\phi_1 = g$. This is equivalent to the existence of a morphism 
$\phi : A \to B \otimes \mh C[t]$ such that
$f = \mr{ev}_0 \circ \phi$ and $g = \mr{ev}_1 \circ \phi$.

Goodwillie \cite[Corollary II.4.4 and Theorem II.5.1]{Goo}
established two closely related features:

\begin{thm}\label{thm:2.2}
The functor $HP_*$ is homotopy invariant and turns nilpotent 
extensions into isomorphisms. Thus, with the above notation,
\begin{align*}
HP_* (f) &= HP_* (g) \\ 
HP_* (I) &= 0
\end{align*} 
and $HP_* (A) \isom HP_* (A/I)$ is an isomorphism. 
\end{thm}

Homotopy invariance can be regarded as a special case of the
K\"unneth theorem, which holds for periodic cyclic homology
under some mild conditions.

\begin{thm}\label{thm:2.28}
Suppose that $A$ is a unital algebra such that
\begin{itemize}
\item the $\varprojlim^1$-term in \eqref{eq:2.5} vanishes, i.e.
$HP_n (A) \cong \lim\limits_{\infty \leftarrow r} HC_{n+2r}(A)$
\item $HP_* (A)$ has finite dimension
\end{itemize} 
Then the K\"unneth theorem holds for $HP_* (A)$. 
This means that for any unital algebra $B$ there is a natural 
isomorphism of $\mh Z / 2 \mh Z$-graded vector spaces
\[
HP_* (A) \otimes HP_* (B) \isom HP_* (A \otimes B) 
\]
\end{thm}
\emph{Proof.} 
See \cite[Theorem 3.10]{Kas} and \cite[Theorem 4.2]{Emm}. 
$\qquad \Box$ \\[3mm]

Reconsider the Hochschild-Kostant-Rosenberg theorem \eqref{eq:2.8}
for the periodic cyclic homology of the ring of regular functions
on a nonsingular affine complex variety. It gives an isomorphism
of $\mh Z / 2 \mh Z$-graded vector spaces
\begin{equation}\label{eq:2.9}
HP_* (\mc O(V)) \cong H_{DR}^* (V)
\end{equation}
Now let $V^{an}$ be the set $V$ endowed with its natural analytic
topology. By a famous theorem of Grothendieck and Deligne (cf. 
\cite{Gro3} and \cite[Theorem IV.1.1]{Hart}) the algebraic De Rham
cohomology of $V$ is naturally isomorphic to the analytic De Rham 
cohomology of $V^{an}$:
\begin{equation}\label{eq:2.20}
H^*_{DR}(V) \cong H^*_{DR} (V^{an} ;\mh C )
\end{equation}
As is well-known, all classical cohomology theories agree on the 
category of smooth manifolds, for instance
\begin{equation}\label{eq:2.10}
H_{DR}^* (V^{an} ;\mh C ) \cong \check H^* (V^{an} ;\mh C )
\end{equation}
the latter denoting \v Cech cohomology with coefficients in 
$\mh C$. Because of the similar functorial properties, it is not
surprising that the composite isomorphism of \eqref{eq:2.9} - 
\eqref{eq:2.10} holds in greater generality. 
This was confirmed in \cite[Theorem 9]{KNS} :

\begin{thm}\label{thm:2.3}
Let $X$ be an affine complex variety, $I \subset \mc O(X)$ an
ideal and $Y \subset X$ the subvariety defined by $I$. It is 
neither assumed that $X$ is nonsingular or irreducible, nor that
$I$ is prime. There is a natural isomorphism
\[
HP_n (I) \cong \check H^{[n]} (X^{an},Y^{an}; \mh C)
:= \prod_{m \in \mh Z} \check H^{n+2m} (X^{an},Y^{an}; \mh C)
\]
\end{thm}

Recall that a \inde{primitive ideal} in a complex algebra is the 
kernel of a (nonzero) irreducible representation of $A$. The 
primitive ideal spectrum Prim$(A)$ is the set of all primitive 
ideals of $A$, and the Jacobson radical \inde{Jac$(A)$} is the 
intersection of all these primitive ideals. Note that every
nilpotent ideal is contained in Jac$(A)$.
We endow \inde{Prim$(A)$} with the \inde{Jacobson topology}, 
which means  that all closed subsets are of the form
\begin{equation}\label{eq:2.79}
\overline S := \{I \in \mr{Prim}(A) : I \supset S \}
\end{equation}
for some subset $S$ of $A$. Denote by \inde{$d_I$} the dimension 
of an irreducible representation with kernel $I \in \mr{Prim}(A)$.
If $d_I < \infty \: \forall I$ then Prim$(A)$ is a $T_1$-space, 
but in general it is only a $T_0$-space.  

For commutative $A$ the primitive ideals are precisely the
maximal ideals, and Prim$(A)$ is an algebraic variety. In this 
case there also is a natural topology on the set Prim$(A)$ that 
makes it into an analytic variety, see \cite[Section 5]{Ser}.

If $\phi : A \to B$ is an algebra homomorphism and $J \in 
\mr{Prim}(B)$ then $\phi^{-1}(J)$ is an ideal, but it is not
necessarily primitive. So Prim is not a functor, it only induces
a map $J \to \overline{\phi^{-1} (J)}$ from Prim$(B)$ to the
power set of Prim$(A)$. However, if for every $J \in \mr{Prim}(B)$
there exists exactly one $I \in \mr{Prim}(A)$ containing 
$\phi^{-1}(J)$, then $\phi$ does induce a continuous map
$\mr{Prim}(B) \to \mr{Prim}(A)$ and we call $\phi$ 
\inde{spectrum preserving}.

Now we give the definition of a \inde{finite type algebra}.
Let \inde{$\mb k$} be a finitely generated commutative unital 
complex algebra, i.e. the ring of regular functions on some 
affine complex variety. A $\mb k$-algebra is a (nonunital) 
algebra $A$ together with a unital morphism from $\mb k$ to 
\inde{$Z(\mc M(A))$}, the center of the multiplier algebra of 
$A$. An algebra $B$ is of finite type if there exists a $\mb k$ 
such that $B$ is $\mb k$-algebra which is finitely generated as 
a $\mb k$-module. 
An algebra morphism $\phi: A \to B$ is a morphism of finite type 
algebras if it is $\mb k$-linear, for some $\mb k$ over which 
both $A$ and $B$ are of finite type.

As announced, the most important theorem in this category says 
that $HP_*$ is determined by Prim, see 
\cite[Theorems 3 and 8]{BaNi}. 

\begin{thm}\label{thm:2.4}
Let $\phi : A \to B$ be a spectrum preserving morphism of finite
type algebras. Then $\mr{Prim}(B) \to \mr{Prim}(A)$ is a 
homeomorphism and 
\[
HP_* (\phi) : HP_* (A) \to HP_* (B)
\]
is an isomorphism.
\end{thm}

More generally, we might have ideals like in \eqref{eq:2.77}, 
such that the induced maps $I_m / I_{m+1} \to J_m / J_{m+1}$ are 
all spectrum preserving, but $\phi : A \to B$ is not. In that 
case $\phi$ is called weakly spectrum preserving. By Theorem 
\ref{thm:2.4} and  Lemma \ref{lem:2.12} such maps also induce
isomorphisms on periodic cyclic homology.
\index{spectrum preserving!weakly}

To understand this better we zoom in on the spectrum, relying 
heavily on \cite[Section 1]{KNS}. Until further notice we assume 
that $A$ is a unital finite type algebra. 
The \inde{central character} map \index{Thet@$\Theta$}
\begin{equation}\label{eq:2.11}
\Theta : \mr{Prim}(A) \to \mr{Prim}(Z(A)) : I \to I \cap Z(A)
\end{equation}
is a finite-to-one continuous surjection. 
For $k,p \in \mh N$ we write
\begin{equation}
\begin{array}{lll}
\mr{Prim}_k (A) & = & \{ I \in \mr{Prim}(A) : d_I = k \} \\
\mr{Prim}_{\leq p}(A) & = & \bigcup_{k=1}^p \mr{Prim}_k (A)
\end{array}
\end{equation}
The sets $\mr{Prim}_{\leq p}(A)$ are all closed and, as the 
frequent occurence of the word ``finite'' already suggests, 
there exists a $N_A \in \mh N$ such that $\mr{Prim}_{\leq N_A}(A) 
= \mr{Prim}(A)$. This leads to the so-called \inde{standard 
filtration} of $A$ :
\begin{equation}\label{eq:2.12}
\begin{split}\index{$I_p^{st}$}
&A = I_0^{st} \supset I_1^{st} \supset \cdots \supset 
  I_{N_A -1}^{st} \supset I_{N_A}^{st} = \mr{Jac}(A) \\
&I_p^{st} \:=\: \bigcap_{d_I \leq p} I = \{ a \in A : \pi(a) = 0 
  \;\mr{if}\: \pi \: \mbox{is a representation with}\: 
  \dim \pi \leq p \}
\end{split}
\end{equation}
Observe that 
\begin{equation}
\mr{Prim} \left( I_q^{st} / I_p^{st} \right) = 
\bigcup_{k = q+1}^p \mr{Prim}_k (A)
\end{equation}
From \eqref{eq:2.12} we also get a filtration of the cyclic 
bicomplex :
\begin{equation}\label{eq:2.59}
\begin{split}
&CC_*^{per}(A) = CC_*^{per}(A)_0 \supset CC_*^{per}(A)_1 \supset
  \cdots \supset CC_*^{per}(A)_{N_A -1} \supset 
  CC_*^{per}(A)_{N_A}\\
&CC_*^{per}(A)_p = \ker \left( CC_*^{per}(A) \to 
  CC_*^{per} \left( A / I_p^{st} \right) \right)
\end{split}
\end{equation}

Using standard (but involved) techniques from homological algebra
we construct a spectral sequence $E_r^{p,q}$ with
\begin{align}
E_0^{p,q} &= CC_*^{per}(A)_{p-1} / CC_*^{per}(A)_p \cong
  \ker \left( CC_{-p-q}^{per} \left( A / I_p^{st} \right) \to
  CC_{-p-q}^{per} \left( A / I_{p-1}^{st} \right) \right) \\
E_1^{p,q} &= HP_{-p-q} \left( I_{p-1}^{st} / I_p^{st} \right) \\
E_\infty^{p,q} &= HP_{-p-q} \left( I_{p-1}^{st} \right) / 
  HP_{-p-q} \left( I_p^{st} \right) 
\end{align}
Moreover $d_0^E : E_0^{p,q} \to E_0^{p,q+1}$ comes directly from
the differential in the cyclic bicomplex and
$d_1^E :  E_1^{p,q} \to E_1^{p+1,q}$ is the composition 
\[
HP_{-p-q} \left( I_{p-1}^{st} / I_p^{st} \right) \to 
HP_{-p-q} \left( A/ I_p^{st} \right) \to HP_{-p-q-1} 
\left( I_p^{st} / I_{p+1}^{st} \right)
\]
of the map induced by the inclusion 
$I_{p-1}^{st} / I_p^{st} \to A/ I_p^{st}$
and the connecting map of the extension
\[
0 \to I_p^{st} / I_{p+1}^{st} \to A / I_{p+1}^{st} \to
A / I_p^{st} \to 0
\]
A most pleasant property of the standard filtration 
\eqref{eq:2.12} is that the quotients $I_{p-1}^{st} / I_p^{st}$
behave like commutative algebras. More precisely, consider the 
analytic space \inde{$X_p$} associated to 
Prim$\left( Z \left( A / I_p^{st} \right) \right)$, 
and its subvariety \index{$Y_p$}
\begin{equation}
Y_p = \left\{ I \in X_p : Z \left( A / I_p^{st} \right)
\cap I_{p-1} / I_p \subset I \right\}
\end{equation}
The central character map for $A / I_p^{st}$ defines a bijection
\begin{equation}\label{eq:2.15}
\mr{Prim}_p (A) = \mr{Prim} \left( I_{p-1}^{st} / 
I_p^{st} \right) \to X_p \setminus Y_p
\end{equation}
and according to \cite[Theorem 1]{KNS} there is a natural 
isomorphism
\begin{equation}\label{eq:2.14}
E_1^{p,q} = HP_{-p-q} \left( I_{p-1}^{st} / I_p^{st} \right)
\cong \check H^{[p+q]} (X_p,Y_p; \mh C)
\end{equation}
Comparing this with Theorem \ref{thm:2.3} we see that 
$I_{p-1}^{st} / I_p^{st}$ is indeed ``close to commutative'' in
the sense that its periodic cyclic homology can be computed as
the \v Cech cohomology of a constant sheaf over its spectrum.

We seek to generalize this to ``less commutative'', nonunital 
finite type algebras. Let $X$ be the set Prim$(\mb k)$ with the 
analytic topology, and $V(A)$ the set Prim$(A)$ with the coarsest 
topology that makes $\Theta : \mr{Prim}(A) \to X$ continuous and 
is finer than the Jacobson topology. This topology depends only 
on the fact that $A$ is a finite type algebra, and not on the 
particular choice of $\mb k$. So if $A$ is unital we may just as 
well assume that $\mb k = Z(A)$ and $X = X_0$.

We will construct a \inde{sheaf} $\mf A$ \index{Amf@$\mf A$} 
over $X$ whose stalk at $x$ is the (finite dimensional) complex 
vector space with basis $\Theta^{-1}(x)$. By definition all 
continuous sections of this collection of stalks are constructed 
from local sections of $\Theta: V(A) \to X$. More precisely, 
given an open $Y \subset X$ we call a section $s$ of 
$\prod_{x \in Y} \mf A (x) \to Y$ continuous at $y \in Y$ if 
there exist 
\begin{itemize}
\item a neighborhood $U$ of $y$ in $Y$
\item connected components $C_1, \ldots, C_n$ of 
      $\Theta^{-1}(U)$, not necessarily different
\item for every $i$ a section $s_i$ of the quotient map from 
      $C_i$ to its Hausdorffization $C_i^H$
\item complex numbers $z_1, \ldots, z_n$
\end{itemize}
such that $\forall x \in U$
\begin{equation}\label{eq:2.91}
s(x) = \sum_{i=1}^n z_i 
\left( s_i (C_i^H) \cap \Theta^{-1}(x) \right)
\end{equation}
For example if $X'$ is a closed subvariety of $X$ and $A = \{ f 
\in \mc O (X) : f (X') = 0\}$ then $\mf A$ is the direct image 
of the constant sheaf (with stalk $\mh C$) on $X \setminus X'$.

Notice that $\mf A$ is functorial in $A$. If $\phi : A \to B$ is
a morphism of finite type $\mb k$-algebras and $V$ is a left 
$A$-module, then $B \otimes_A V$ is a $B$-module. If we consider
only the semisimple forms of these modules, then we get a 
homomorphism
\[
\mh Z [\mr{Prim}(A)] \to \mh Z [\mr{Prim}(B)]
\]
which extends naturally to a morphism $\mf A \to \mf B$ of 
sheaves over $X$.

The motivation for this sheaf comes from topological $K$-theory:
the local sections $s_i$ are supposed to model ``local'' 
idempotents in $A$. The classes of these things should
generate $HP_* (A)$, leading to 

\begin{thm}\label{thm:2.6}
There is an unnatural isomorphism of finite dimensional vector 
spaces
\[
HP_* (A) \cong \check H^* (X ; \mf A)
\]
\end{thm}
\emph{Proof.}
Assume first that $A$ is unital. 
Let $\mf A_p$ be the sheaf (over $X$) constructed from 
$A / I_p^{st}$ in the same way as we constructed $\mf A$ from 
$A$; it has stalks
\begin{equation}
\mf A_p (x) = \mh C \{ \Theta^{-1}(x) \cap 
\mr{Prim}_{\leq p}(A) \}
\end{equation}
Since $\mr{Prim}_{\leq p}(A)$ is closed in $\mr{Prim}(A)$ there 
is a natural surjection $\mf A \to \mf A_p$, which comes down
to forgetting all primitive ideals $I$ with $d_I > p$.
Thus we get filtrations of the (pre)sheaf $\mf A$:
\begin{equation}
\begin{aligned}
&\mf A = \mf I_0 \supset \mf I_1 \supset \cdots \supset
\mf I_{N_A -1} \supset \mf I_{N_A} = 0 \\
&\mf I_p = \ker \left( \mf A \to \mf A_p \right) 
\end{aligned}
\end{equation}
and of the \v Cech complex $\check C^* (X;\mf A)$ (this is a
pretty complicated object, see \cite[\S 5.8]{God}) :
\begin{equation}\label{eq:2.60}
\begin{split}
&\check C^* (X;\mf A) = \check C^* (X;\mf A)_0 \supset 
  \check C^* (X;\mf A)_1 \supset \cdots \supset 
  \check C^* (X;\mf A)_{N_A -1} \supset 
  \check C^* (X;\mf A)_{N_A} = 0\\
&\check C^* (X;\mf A)_p = \ker \left( \check C^* (X;\mf A) \to
  \check C^* (X;\mf A_p) \right) \cong \check C^* (X;\mf I_p)
\end{split}
\end{equation}
The presheaf $\mf B_p := \ker (\mf A_p \to \mf A_{p-1})$ is 
actually a sheaf, and it has stalks
\begin{equation}
\mf B_p(x) = \mh C \{ \Theta^{-1}(x) \cap \mr{Prim}_p (A) \}
\end{equation}
From these data we construct a spectral sequence $F_r^{p,q}$
with terms
\begin{align}
F_0^{p,q} &= \check C^{p+q}(X; \mf A)_{p-1} /
  \check C^{p+q}(X;\mf A)_p \cong \check C^{p+q}(X; \mf B_p) \\
F_1^{p,q} &= \check H_{p+q}(X; \mf B_p) \\
F_\infty^{p,q} &= \check H_{p+q}(X; \mf I_{p-1}) /  
  \check H_{p+q}(X; \mf I_p)
\end{align}
In this sequence $d_0^F : F_0^{p,q} \to F_0^{p,q+1}$ is the 
normal \v Cech differential, while $d_1^F : F_1^{p,q} \to 
F_1^{p+1,q}$ is induced by the inclusion $\mf B_p \to \mf A_p$
and the connecting map associated to the short exact sequence
\[
0 \to \mf B_{p+1} \to \mf A_{p+1} \to \mf A_p \to 0
\]
From \eqref{eq:2.15} and the local nature of the continuity
condition for $\mf B_p$ we see that there are natural 
isomorphisms
\begin{align*}
\check C^{p+q}(X; \mf B_p) &\cong \check 
  C^{p+q}(X_p ,Y_p ; \mh C) \\
\check H^{p+q}(X; \mf B_p) &\cong \check 
  H^{p+q}(X_p ,Y_p ; \mh C)
\end{align*}
Clearly, all this was set up to compare the spectral sequences
$E_r^{p,q}$ and $F_r^{p,q}$. On the first level we have a diagram
\begin{equation}\label{eq:2.16}
\begin{array}{ccc}
E_1^{p,q} & \xrightarrow{\quad d_1^E \quad} & E_1^{p+1,q} \\[1mm]
\cong & & \cong \\[1mm]
\check H^{[p+q]} (X_p ,Y_p ;\mh C) & & 
\check H^{[p+q+1]} (X_{p+1} ,Y_{p+1} ;\mh C) \\[1mm]
\cong & & \cong \\
\prod\limits_{n \in \mh Z} F_1^{p,q+2n} & \xrightarrow{\quad 
  d_1^F \quad} & \prod\limits_{n \in \mh Z} F_1^{p+1,q+2n}
\end{array}
\end{equation}
Since $d_1^E$ ($d_1^F$) is natural with respect to 
filtration-preserving morphisms of $\mb k$-algebras
(of presheaves over $X$), these differentials must commute with
the natural isomorphisms in the diagram \eqref{eq:2.16}. This
yields natural isomorphisms
\[
E_r^{p,q} \cong \prod_{n \in \mh Z} F_r^{p,q+2n}
\]
for all $r \geq 1$. For $r = \infty$ we see that there exist 
filtrations of finite length on $HP_* (A)$ and 
$\check H^* (X;\mf A)$, such that the associated graded objects 
are isomorphic. Hence $HP_* (A)$ and $\check H^* (X;\mf A)$,
being vector spaces, are unnaturally isomorphic.

Moreover they have finite dimension since every term
$\check H^{[p+q]} (X_p ,Y_p ;\mh C )$, being the cohomology of 
the affine algebraic variety $X_p \setminus Y_p$, has finite
dimension by \cite[Theorems 4.6 and 6.1]{Hart}.

This proves the theorem for unital finite type algebras, so let
us now assume that $J$ is an nonunital finite type 
$\mb k$-algebra. By stability $HP_* (M_2 (J)) \cong HP_* (J)$
and the sheaves corresponding to $M_2 (J)$ and $J$ are 
isomorphic, so we may assume that $J$ has no one-dimensional 
representations. Consider now the unital finite type algebra 
$A = \mb k + J$, with multiplication
\begin{equation}
(f_1 ,b_1 ) (f_2 ,b_2 ) = (f_1 f_2 ,f_1 b_2 + f_2 b_1 + b_1 b_2 )
\end{equation}
Its standard filtration is 
\begin{equation}
A = I_0^{st} \supset J = I_1^{st} \supset \cdots \supset
I_{n_A -1}^{st} \supset I_{n_A}^{st} = \mr{Jac}(A) = \mr{Jac}(J)
\end{equation}
The above considerations show that, as vector spaces, 
\begin{equation}
HP_{-m}(J) = HP_{-m}(I_1^{st}) \cong \prod_{p=2}^m 
E_\infty^{p,m-p} \cong \prod_{p=2}^m \prod_{n \in \mh Z}
F_\infty^{p,2n+m-p} \cong \check H^{[m]}(X;\mf I_1 )
\end{equation}
It only remains to see that $\mf I_1$ is isomorphic to the sheaf 
constructed from $J$, but this is clear from looking at the 
stalks. $\qquad \Box$
\\[2mm]

So we managed to describe the periodic cyclic homology of a 
finite type $\mb k$-algebra using only the following data:

\begin{itemize}
\item the spectrum Prim$(A)$ with a natural topology that makes
      it a \inde{non-Hausdorff manifold}
\item the complex analytic variety $X$
\item the continuous map $\Theta : \mr{Prim}(A) \to X$
\end{itemize}

For some time the author believed that this construction on page
\pageref{eq:2.91} could be extended to a cohomology theory on
the category of non-Hausdorff manifolds, but now it seems to him
that it only gives good results under rather restrictive 
conditions. Apparently we need the following implication of
\eqref{eq:2.15} : there exists a stratification of Prim$(A)$
such that at every level the set of non-Hausdorff points in a
component is either the whole component, or a submanifold of 
lower dimension.
\\[4mm]

\section{Topological cyclic theory}
\label{sec:2.3}

We would like to discuss the topological counterpart of the 
algebraic cyclic theory of Section \ref{sec:2.1}. 
To prepare for this, and to fix certain notations, we start 
by recalling some general results for m-algebras and 
topological tensor products.
When studying the literature, it quickly becomes clear that this
topological setting is significantly more tricky than the 
purely algebraic setting, for several reasons. Firstly, the 
category of topological vector spaces is not abelian, i.e. not 
every closed subspace has a closed complement. Secondly, the 
tensor product of two topological vector spaces is not unique, 
and the functor ``$\otimes_t A$'' (for some unambiguous choice 
of a topological tensor product) is in general not exact. 
And finally, although the appropriate results are all known to 
experts, there does not appear to be an overview available.
\\[2mm]

A \inde{topological algebra} $A$ (over $\mh C$) is an algebra 
with a topology such that addition and scalar multiplication 
are jointly continuous, while multiplication is separately 
continuous. When we talk about the spectrum of $A$, we usually
mean the set \inde{Prim$(A)$} of all closed primitive ideals of 
$A$. The closed subsets of Prim$(A)$ are as in \eqref{eq:2.79}. 

A \inde{seminorm} on $A$ is a map $p : A \to [0, \infty)$ 
with the properties
\begin{itemize}
\item $p(\lambda a) = |\lambda| p(a)$
\item $p(a+b) \leq p(a) + p(b)$
\end{itemize}
for all $a,b \in A$ and $\lambda \in \mh C$.
Moreover $p$ is called submultiplicative if
\begin{itemize}
\item $p(ab) \leq p(a) p(b)$
\end{itemize}
We say that $p'$ dominates $p$ if $p'(a) \geq p(a) \; \forall 
a \in A$. If $\{ p_i \}_{i \in I}$ is a collection of seminorms, 
then there is a coarsest topology on $A$ making all the $p_i$ 
continuous. The sets 
\[
\{ a \in A : p_i (a-b) < 1/n \} \qquad 
b \in A, n \in \mh N, i \in I
\]
form a subbasis for this topology. If it agrees with the original 
topology, then we call $A$ a \inde{locally convex algebra} and 
say that it has the topology defined by the family of seminorms 
$\{ p_i \}_{i \in I}$. Notice that the $p_i$ may have nontrivial 
nullspaces $N_i$ and that $A$ is Hausdorff if and only if 
$\cap_{i \in I} N_i = 0$. Furthermore the multiplication in
$A$ is jointly continuous if, but not only if, all the $p_i$ are 
submultiplicative.

Two families of seminorms are equivalent if every member of 
either family is dominated by a finite linear combination of
seminorms from the other family. Two families of seminorms define 
the same topologies if and only if they are equivalent.

A locally convex algebra is \inde{metrizable} if and only if its 
topology can be defined by a countable family of seminorms
$\{ p_i \}_{i=1}^\infty$ with $\cap_{i=1}^\infty N_i = 0$.
In that case a metric is given by
\begin{equation}
d(a,b) = 
\sum_{i=1}^\infty \frac{2^{-i} p_i (a-b)}{1 + p_i (a-b)}
\end{equation}
Clearly this implies a notion of completeness for such algebras, 
and it can be generalized to all locally convex algebras by
means of Cauchy filters on uniform spaces. For sequences this
comes down to calling a sequence $\{ a_n \}_{n=1}^\infty$ in $A$
Cauchy if and only if for every $i \in I$ the sequence
$\{ x_n + N_i \}_{n=1}^\infty$ is Cauchy in the normed space
$A / N_i$. \index{Fr\'echet algebra}

Combining all these notions, an \inde{m-algebra} is a complete 
Hausdorff locally convex algebra $A$ whose topology can be 
defined by a family of submultiplicative seminorms. We call $A$ 
Fr\'echet if it is metrizable on top of that. If 
$B$ is a topological algebra such that $GL_1 (B^+)$ is open in 
$B^+$, then we call $B$ a \inde{Q-algebra}. Every Banach algebra, 
but not every Fr\'echet algebra, is a Q-algebra.

Since m-algebras are not so well-known we state some important
properties. Let $A$ be a unital m-algebra and 
$A^\times = GL_1 (A)$ the set of invertible elements in $A$.
Recall that the \inde{spectrum} of an element $a \in A$ is
\[
\mr{sp}(a) = \{ \lambda \in \mh C : a - \lambda \notin A^\times \}
\]
Contrarily to the Banach algebra case, sp$(a)$ is in general
not compact.

\begin{thm}\label{thm:2.7}
\begin{enumerate}
\item M-algebras are precisely the projective limits of Banach
algebras.
\item Inverting is a continuous map from $A^\times$ to $A$.
\item Suppose that $U \subset \mh C$ is an open neighborhood of
sp$(a)$, and let $C^{an}(U)$ be the algebra of holomorphic 
functions on $U$. There exists a unique continuous algebra
homomorphism, the \inde{holomorphic functional calculus}
\[
C^{an}(U) \to A : f \to f (a)
\]
such that $1 \to 1$ and $\mr{id}_U \to a$.
\item If $\Gamma$ is a positively oriented smooth simple 
closed contour, around sp$(a)$ and in $U$, then 
\[
f(a) = {\ds \frac{1}{2 \pi i}} \int_\Gamma f(\lambda ) 
(\lambda - a)^{-1} d \lambda \qquad \forall f \in C^{an}(U)
\]
\end{enumerate}
\end{thm}
\emph{Proof.}
1 and 2 were proved by Michael \cite[Theorems 5.1 and 5.2]{Mic}.
3 and 4 are well-known for Banach algebras, see e.g. 
\cite[Proposition 2.7]{Tak}. Using 1 they can be extended to
m-algebras, as was noticed in \cite[Lemma 1.3]{Phi2}.
$\qquad \Box$ \\[3mm]

For some typical examples, consider a $C^k$-manifold $X$, 
with $k \in \{ 0,1,2, \ldots, \infty\}$. We shall always assume 
that our manifolds are $\sigma$-compact, hence in particular 
paracompact. Let $U \subset \mh R^d$ be an open set and 
$\phi : U \to X$ a chart. For a multi-index $\alpha$ with 
$|\alpha| = n$ and $g \in C^n (U)$ let \index{manifold}
\begin{equation}
\partial^\alpha g = \frac{\partial^n g}{\partial y_{\alpha_1} 
\cdots \partial y_{\alpha_n}} \in C^0 (U)
\end{equation}
be the derivative of $g$ with respect to the standard 
coordinates $y_1 ,\ldots ,y_d$ of $\mh R^d$. For $K \subset U$ 
compact and $n \in \mh N_{\leq k}$ we define a seminorm 
$\nu_{n, \phi, K}$ on $C^k (X)$ by
\begin{equation}
\nu_{n, \phi, K} (f) = \sup_{y \in K} \sum_{|\alpha| \leq n}
\frac{\left| \partial^\alpha (f \circ \phi) (y) \right|}{
|\alpha| !}
\end{equation}
Straightforward estimates show that every $\nu_{n, \phi, K}$
is submultiplicative and that $C^k (X)$ is complete with respect
to the family of such seminorms. Moreover, because $X$ is \
$\sigma$-compact, we can cover it by countably many sets 
$\phi_i (K_i)$.
\begin{equation}
\{ \nu_{n, \phi_i, K_i} : i,n \in \mh N , n \leq k \}
\end{equation}
is a countable collection of seminorms defining the topology
of $C^k (X)$, which therefore is a Fr\'echet algebra.

Finally, if $X$ is compact and $k \in \mh N$ then $C^k (X)$ is a
Banach algebra. Indeed, if we cover $X$ by finitely many sets
$\phi_i (K_i)$ then
\begin{equation}
\norm{f} = \sum_i \nu_{n, \phi_i, K_i} (f)
\end{equation}
is an appropriate norm.

We now give a quick survey of topological tensor products, 
completely due to Grothendieck \cite{Gro1}. To fix the notation, 
we agree that by $\otimes$ without any sub- or superscript we 
always mean the algebraic tensor product. By default we take it 
over $\mh C$ if both factors are complex vector spaces, and over 
$\mh Z$ if there is no field over which both factors are vector 
spaces. \index{tena@$\otimes$}

The algebraic tensor product of two vector spaces $V$ and $W$ 
solves the universal problem for bilinear maps. This means that 
every bilinear map from $V \times W$ to some vector space $Z$ 
factors as \index{tensor product!algebraic}
\[
\begin{array}{c}
V \times W \quad \longrightarrow \quad Z \qquad \\
\searrow \qquad\;\; \nearrow \; \\
V \otimes W 
\end{array}
\]
resulting in a bijection between Bil($V \times W,Z)$ and 
Lin($V \otimes W,Z)$. \index{tensor product!topological}
This procedure can be extended in several ways to the category
of locally convex spaces, corresponding to different classes of
bilinear maps and different topologies on $V \otimes W$.

For example we have the projective tensor product $V 
\otimes_\pi W$ \cite[Subsection I.1.1]{Gro1}, called so because 
it commutes with projective limits. It is $V \otimes W$ with the 
topology solving the universal problem for jointly continuous 
bilinear maps $V \times W \to Z$. If $\{ p_i \}_{i \in I}$ and 
$\{ q_j \}_{j \in J}$ are defining families of seminorms for 
$V$ and $W$, then this topology is defined by the family of 
seminorms \index{tensor product!projective}
\begin{equation}\label{eq:2.17}
\gamma_{ij}(x) = \inf \left\{ \sum_{k=1}^n p_i (v_k) q_j (w_k) 
: x = \sum_{k=1}^n v_k \otimes w_k \right\} 
\qquad i \in I, j \in J
\end{equation}
The completion $V \hot W$ of $V \otimes W$
for the associated uniform structure is called the completed 
projective tensor product. 
\index{tenp@$\hot$} \index{tensor product!inductive}

Similarly the inductive tensor product $V \otimes_i W$
\cite[Subsection I.3.1]{Gro1} solves the universal problem for
separately continuous bilinear maps, and it commutes with
inductive limits. The topology of $V \otimes_i W$ is finer than 
that of $V \otimes_\pi W$, and the associated completion is 
denoted by $V \oot W$. Typically, for a $C^k$-manifold $X$ and 
a Banach space $V$ we have \index{teni@$\oot$}
\begin{equation}
C^k (X) \oot V \cong C^k (X;V)
\end{equation}
There exists also more subtle structures on $V \otimes W$, such 
as the injective tensor product $V \otimes_\ep W$
\cite[p. I.89]{Gro1}, which in a certain sense has the weakest
reasonable topology. \index{tensor product!injective}

\index{nuclear vector space} 
If $V$ satisfies $V \otimes_\ep Z = V \otimes_\pi Z$ for every
$Z$ then it is called nuclear, and if both $V$ and $W$ are 
nuclear, then so are $V \otimes_\pi W$ and $V \hot W$.
On the other hand, if $V$ and $W$ are both Fr\'echet spaces, 
then $V \otimes_i W = V \otimes_\pi W$ \cite[p. I.74]{Gro1}
and its completion $V \oot W = V \hot W$ is again a 
Fr\'echet space \cite[Th\'eor\`eme II.2.2.9]{Gro1}.

Consequently the tensor powers of a nuclear Fr\'echet space
can be defined unambiguously. For example if $X$ and $Y$ are 
smooth manifolds, then $C^\infty (X)$ and $C^\infty (Y)$ are 
nuclear Fr\'echet spaces and 
\begin{equation}
C^\infty (X) \hot C^\infty (Y) \cong C^\infty (X \times Y)
\end{equation}
Now that we have come this far, it is logical to spend a few
words on topological tensor products over rings.
\index{tensor product!topological!over a ring}
So let $A$ be an m-algebra, $V$ a right $A$-module and
$W$ a left $A$-module. We assume that $V$ and $W$ are complete
Hausdorff locally convex spaces and that the module operations
are jointly continuous. Then the completed projective tensor
product $V \hot_A W$ is the completion of 
$V \otimes_A W$ for the topology solving the universal problem
for jointly continuous $A$-bilinear maps from $V \times W$ to
some $A$-module $Z$. Just as over $\mh C$, this topology is
defined by the family of seminorms \eqref{eq:2.17}.

Let us return to homology of algebras. In any category of
locally convex algebras with a topological tensor product
$\otimes_t$ we can form the bicomplex $CC^{per}(A, \otimes_t )$
with spaces 
\[
CC^{per}_{p,q}(A, \otimes_t) = A^{\otimes_t q+1}
\]
The maps from \eqref{eq:2.1} and \eqref{eq:2.2} are continuous
because they use only the algebra operations of $A$.
This, and the subcomplexes $CC (A,\otimes_t )$ and
$CC^{\{2\}}(A,\otimes_t )$, lead to functors 
$HH_n (A,\otimes_t ) ,\, HC_n (A,\otimes_t )$ and
$HP_n (A,\otimes_t )$. They are related by Connes' periodicity 
exact sequence, but to get more nice features it is imperative
that we use only completed tensor products and place ourselves
in one of the following categories:

\begin{itemize}\label{p:cat}
\item \inde{$\mc{CLA}$}: complete Hausdorff locally convex 
      algebras
\item \inde{$\mc{MA}$} : m-algebras
\item \inde{$\mc{FA}$} : Fr\'echet algebras
\item \inde{$\mc{BA}$} : Banach algebras
\end{itemize}

Although the objects of none these categories form a set, we 
allow ourselves to use $\in$ to indicate with what kind of
algebra we are dealing.

By Theorem \ref{thm:2.7} the completed projective tensor
product of two m-algebras is again an m-algebra, so we use
$\hot$ as our default and $\oot$
as a reserve. Just as in Section \ref{sec:2.1} we are going to
study the functorial properties of the resulting homology 
theories. Let $A, A_m \in \mc{CLA} ,\, m \in \mh N$.

Notice that the topological cyclic bicomplexes under 
consideration contain the algebraic cyclic bicomplexes. This 
yields natural transformations from the algebraic cyclic 
theories to their topological counterparts. Therefore any 
homomorphism from a complex algebra $B$ to $A$ induces maps on 
homology groups like 
\[
HH_n (B) \to HH_n (A ,\hot)
\]
These maps are compatible with all the properties below.

\begin{enumerate}
\item Additivity. \index{additivity}
\begin{align*}
HH_n \left( \bigoplus_{m=1}^\infty A_m ,\oot \right) 
&\cong \bigoplus_{m=1}^\infty HH_n (A_m ,\oot) \\
HH_n \left( \prod_{m=1}^\infty A_m ,\hot \right) 
&\cong \prod_{m=1}^\infty HH_n (A_m ,\hot)
\end{align*}
and similarly for $HC_n$. The corresponding isomorphisms for 
$\oot$ and $\prod$ hold if $A_m \in \mc{FA} \; \forall m$ and 
the isomorphisms for $\hot$ and $\bigoplus$ are valid if 
$A_m \in \mc{BA} \; \forall m$. For $HP_n$ we can only be sure 
about the case with $\prod$ and $\hot$.

The proof of all these statements can be reduced to that of 
the algebraic case, by using 
\cite[Propositions I.1.3.6 and I.3.1.14]{Gro1}.

\item Stability. \index{stability} \label{p:hp}
\[
HH_n (M_m (A), \hot ) \cong HH_n (A, \hot )
\]
and similarly with $HC_n ,\, HP_n$ and $\oot$.

This follows from the algebraic case, since all topological 
tensor products of $A$ with a finite dimensional vector space 
(such as $M_m (\mh C)$) are the same, and essentially equal to
the algebraic tensor product.

It is not known to the author whether $HH_n$ and $HC_n$ are
Morita-invariant in a more general sense, but for $HP_n$
we will soon return to this point.

\item Continuity. \index{continuity}
Here great concessions to the algebraic case must be made. Assume
that all the $A_m$ are nuclear Fr\'echet algebras and that
$A = \lim_{m \to \infty} A_m$ is a \inde{strict inductive limit}.
(Strict means that all the maps $A_m \to A_{m+1}$ are injective 
and have closed range.) In this setting Brodzki and Plymen showed 
\cite[Theorem 2]{BrPl1} that
\begin{align*}
HH_n (A, \oot) &\cong 
  \lim_{m \to \infty} HH_n (A_m ,\oot) \\
HC_n (A, \oot) &\cong 
  \lim_{m \to \infty} HC_n (A_m ,\oot)
\end{align*}
To make $HP_n$ continuous we need even more conditions. For 
example if $\exists N \in \mh N$ such that $HH_n (A_m 
,\oot) = 0 \; \forall n > N ,\, \forall m$,
then by \cite[Theorem 3]{BrPl1}
\[
HP_n (A, \oot) \cong 
\lim_{m \to \infty} HP_n (A_m ,\oot)
\]
The author knows of no continuity results for $\hot$, which is
not surprising, considering the bad compatibility of projective
tensor products with inductive limits.
\end{enumerate}

Excision is also pretty subtle for topological algebras. Let
$\mc A$ be one of the four categories from page \pageref{p:cat}.
Extending Wodzicki's terminology, we call $A \in \mc A$ strongly
H-unital if, for every $V \in \mc A$, the homology of the
differential graded complex  \index{H-unital!strongly}
$(A^{\hot n} \hot V \,,\, b' \hot \, \mr{id}_V)$ is 0.

It follows from \cite[Section 1]{Joh} that every Banach algebra
with a left or right bounded approximate identity (e.g. a 
$C^*$-algebra) is strongly H-unital, and in \cite[Section 3]{BrLy}
it is claimed that this also holds in $\mc{FA}$.

Recall that an extension of topological vector spaces
$0 \to Y \to Z \to W \to 0$ \index{admissible extension}
is admissible if it has a continuous linear splitting.
This implies in particular that $Y$ (or more precisely, its image) 
has a closed complement in $Z$. Furthermore we call an extension
\begin{equation}\label{eq:2.18}
0 \to A \to B \to C \to 0
\end{equation}
in $\mc A$ topologically pure if, for every $V \in \mc A$, 
\index{topologically pure extension}
\begin{equation}
0 \to A \hot V \to B \hot V \to C \hot V \to 0
\end{equation}
is exact. According to \cite[Section 4]{BrLy} the following types
of extensions are topologically pure in $\mc{FA}$:

\begin{enumerate}
\item admissible extensions
\item extensions \eqref{eq:2.18} such that $A$ has a bounded
left or right approximate identity
\item\label{it:2.1} extensions of nuclear Fr\'echet algebras
\end{enumerate}

With this terminology, the following is proved in
\cite[Theorems 2 and 4]{BrLy} :

\begin{thm}\label{thm:2.8}
Let $0 \to A \to B \to C$ be a topologically pure extension of
Fr\'echet algebras, with $A$ strongly H-unital. Then there are
long exact sequences
\[
\begin{array}{ccccccccc}
%\cdots & 
\to & HH_n (A,\hot) & \to & HH_n (B,\hot) & \to & 
 HH_n (C,\hot) & \to & HH_{n-1}(A,\hot) & \to \\% & \cdots \\
%\cdots & 
\to & HC_n (A,\hot) & \to & HC_n (B,\hot) & \to & 
 HC_n (C,\hot) & \to & HC_{n-1}(A,\hot) & \to \\% & \cdots \\
%\cdots & 
\to & HP_n (A,\hot) & \to & HP_n (B,\hot) & \to & 
 HP_n (C,\hot) & \to & HP_{n-1}(A,\hot) & \to % & \cdots 
\end{array}
\]
\end{thm}
With the help of Theorem \ref{thm:2.7}, all these results on 
\inde{excision} (except \ref{it:2.1}.) can be extended to the 
category of m-algebras.

Actually $HP_*$ has much more features than those listed 
above. Let $f,g : A \to B$ be morphisms in $\mc{CLA}$.
We say that they are homotopic if there exists a morphism 
$\phi : A \to C ([0,1],B)$ such that $f = \mr{ev}_0 \circ \phi$
and $g = \mr{ev}_1 \circ \phi$. They are called diffeotopic 
if there exists a morphism
\index{homotopy} \index{diffeotopy}
\[
\phi : A \to C^\infty ([0,1]) \hot B \cong C^\infty ([0,1];B)
\]
with these properties.

\begin{thm}\label{thm:2.24}
In the category $\mc{MA}$ the functor 
$HP_* ( \,\cdot\,, \hot )$ has the following properties:
\begin{enumerate}
\item If $f,g$ are diffeotopic, then $HP_* (f) = HP_* (g)$.
\item Let $E$ and $F$ be linear subspaces of an m-algebra $A$,
and let $A (EF)$, respectively $A (FE)$, be the subalgebra 
generated by all the products $e f$, respectively $f e$, with 
$e \in E \,, f \in F$. Then $HP_* (A(EF)) \cong HP_* (A(FE))$.
\item Every admissible extension \eqref{eq:2.18} gives rise to 
an exact hexagon 
\[
\hexagon{HP_0 (A,\hot)}{HP_0 (B,\hot)}{HP_0 (C,\hot)}{HP_1 
(C,\hot)}{HP_1 (B,\hot)}{HP_1 (A,\hot)}
\]
\end{enumerate}
\end{thm}
\emph{Proof.} 
1 comes from \cite[p. 125]{Con} 2 from \cite{Cun3} and 3 from 
\cite{Cun2}. $\qquad \Box$
\\[2mm]

A clear omission at this point is a K\"unneth theorem for 
topological periodic cyclic homology. It certainly exists, but
the author does not know in what generality. Fortunately, for
all the algebras that we use there is an ad hoc argument 
available to prove the K\"unneth isomorphism.

What happens to \inde{differential forms} in the presence of a 
topology?
If $A$ is a commutative unital m-algebra, then the definition of
$\Omega^1 (A)$ must be modified to retain completeness. So, 
identifying the K\"ahler differential $a\, db$ with the elementary
tensor $a \otimes b$, we define $\Omega^1 (A,\hot)$ to be the 
quotient of $A \hot A$ by the closed $A$-submodule generated by
the relations \eqref{eq:2.19}. Furthermore let $V_n$ be closed
subspace of $\left( \Omega^1 (A,\hot) \right)^{\hot_A n}$ 
generated by all the $n$-forms 
$\omega_1 \wedge \cdots \wedge \omega_n$ for which
there exist $i \neq j$ with $\omega_i = \omega_j$. Then
\begin{equation}
\Omega^n (A,\hot) = \bigwedge\nolimits ^n_A \Omega^1 (A,\hot)
:= \left( \Omega^1 (A,\hot) \right)^{\hot_A n} / V_n
\end{equation}
Thus, finally, we have the topological \inde{De Rham homology}
\begin{equation}
H^{DR}_n (A,\hot) = H_n \left( \Omega^* (A,\hot), d \right)
\end{equation}
Let us consider the topological counterpart of a smooth algebra.
It is not exactly clear what that should be, but obviously it
should be related to algebras of smooth functions. Nuclearity is
also an advantage.
So let $C^\infty (X)$ be the (nuclear Fr\'echet) algebra of
infinitely often differentiable complex valued functions on a
smooth real manifold $X$. It is well known that 
in this case we have natural isomorphisms
\begin{align}
\label{eq:2.28} \Omega^n (C^\infty (X),\hot) &\cong \Omega^n (X)\\
\label{eq:2.29} H^{DR}_n (C^\infty (X),\hot) &\cong H_{DR}^n (X)
\end{align}
These hold both with real and with complex coefficients, but we 
are mostly interested in the latter. Furthermore the maps
\eqref{eq:2.13} and the Hochschild-Kostant-Rosenberg theorem
can be extended to this topological situation 
\cite{Con,Tel,Was2}, so there are natural isomorphisms
\begin{align}
\label{eq:2.33} HH_n (C^\infty (X),\hot) &\cong \Omega^n (X)\\
\label{eq:2.34} HC_n (C^\infty (X),\hot) &\cong \Omega^n (X) / 
 d \Omega^{n-1}(X) \oplus H_{DR}^{n-2}(X) \oplus H_{DR}^{n-4}(X) 
\oplus \cdots \\
\label{eq:2.35} HP_n (C^\infty (X),\hot) 
 &\cong \prod_{m \in \mh Z} H_{DR}^{n+2m}(X)
\end{align}

We conclude the section with a warning. An algebra may be ``too 
big'' for cyclic theory to work properly. In fact the results 
are pretty noninformative for most Banach algebras. Let $A$ be 
an \inde{amenable Banach algebra} \cite{Joh}, for example $C(Y)$ 
with $Y$ a compact Hausdorff space or $L^1 (G)$ with $G$ a 
locally compact amenable group. Then we have
\begin{equation}\label{eq:2.26}
HH_n (A,\hot) =  
\left\{ \begin{array}{c@{\quad \mr{if} \quad}c}
A / [A,A] & n = 0\\
0 & n > 0 
\end{array} \right.
\end{equation}
where $[A,A]$ is the range of the commutator map
$A \hot A \to A$. Thus $\forall n \geq 0$
\begin{equation}\label{eq:2.27}
\begin{aligned}
&HC_{2n}(A,\hot)   &=\;\;& HP_{2n}(A,\hot)   &=\;\;& A / [A,A]\\
&HC_{2n+1}(A,\hot) &=\;\;& HP_{2n+1}(A,\hot) &=\;\;& 0 
\end{aligned}
\end{equation}
\\[2mm]

\section{Topological $K$-theory and the Chern character}
\label{sec:2.4}

Topological $K$-theory is at the very heart of noncommutative 
geometry. For a compact topological space it is defined roughly 
speaking as the Grothendieck group of equivalence classes of
vector bundles over $X$. By the Gelfand-Na\u{\i}mark and
Serre-Swan theorems it can be transferred to (commutative)
$C^*$-algebras, and there it becomes the Grothendieck group of
isomorphism classes of finitely generated projective modules. 
This in turn can be extended to Banach algebras, and on that 
category $K_* (A)$ is something like the Grothendieck group of
homotopy classes of idempotents or invertibles in the 
stabilization of $A$. 

In the present section we study the $K$-functor on even larger 
categories of topological algebras. We collect some important
theorems, focussing especially on those results that are fit
to compare $K$-theory with topological periodic cyclic homology.

Of course there also exists a purely algebraic $K$-theory, which
is a natural companion of the algebraic cyclic theory of 
Section \ref{sec:2.1}. However, since these algebraic $K$-groups
are notoriously difficult to compute, and since they contain 
more number-theoretic than geometric information, we will not
study them here.

The most general construction of a topological $K$-functor is due
to Cuntz \cite{Cun1}, and it realizes $K_*$ as the covariant half
of a bivariant functor on the category of m-algebras. As Cuntz's
construction is rather complicated, we will not elaborate on it.
Instead we recall the definition of Phillips \cite{Phi2}, which
works for Fr\'echet algebras and is similar to that for Banach
algebras.

Let $\mf K$ be the nuclear Fr\'echet algebra of infinite matrices 
with rapidly decreasing coefficients. It is also referred to as the 
algebra of \inde{smooth compact operators}, because it is a 
holomorphically closed dense *-subalgebra of the usual algebra of 
compact operators, and it is isomorphic, as a nuclear Fr\'echet 
space, to the algebra of smooth functions $C^\infty (\mh T^2)$ on 
the two-dimensional torus. \index{Kmf@$\mf K$}

For any Fr\'echet algebra $A$, let $(\mf K \hot A)^+$ be the
unitization of $\mf K \hot A$, and consider the Fr\'echet
algebra $M_2 \left( (\mf K \hot A)^+ \right)$. Define
$\bar P(A)$ to be the set of all idempotents $e$ in this 
algebra satisfying 
\[
e - \begin{pmatrix} 1 & 0  \\ 0 & 0 \end{pmatrix} 
\in M_2 (\mf K \hot A)
\]
Similarly $\bar U(A)$ is the set of all invertible elements 
$u \in M_2 \left( (\mf K \hot A)^+ \right)$ for which
\[
u - \begin{pmatrix} 1 & 0  \\ 0 & 1 \end{pmatrix} 
\in M_2 (\mf K \hot A)
\]
Following \cite[Definition 3.2]{Phi2} we put \index{$K$-theory}
\begin{align}
K_0 (A) &= \pi_0 \left( \bar P (A) \right)\\
K_1 (A) &= \pi_0 \left( \bar U (A) \right)
\end{align}
With the multiplication defined by the direct sum of matrices,
these turn out to be abelian groups with unit elements 
\[
\left[ \begin{pmatrix} 1 & 0  \\ 0 & 0 \end{pmatrix} \right]
\quad \mr{and} \quad
\left[ \begin{pmatrix} 1 & 0  \\ 0 & 1 \end{pmatrix} \right]
\]
Later we shall want to pick ``nice'' representants of $K$-theory
classes, so now we try to discover how much is possible in this
respect. Let $A$ be unital, $e \in M_n (A)$ idempotent and 
$u \in GL_n (A)$. Pick a rank one projector $p \in \mf K$ and
an isomorphism $M_n (\mf K) \to \mf K$ and extend it to 
\[
\lambda_n : M_n (\mf K \hot A) \isom \mf K \hot A
\]
Now consider the elements
\begin{equation}\label{eq:2.82}
\begin{pmatrix} 1 & 0 \\ 0 & \lambda_n (pep) \end{pmatrix} 
\in \bar P (A) \quad \mr{and} \quad \begin{pmatrix} 1 & 0 \\ 
0 & \lambda_n (1 - p + pup) \end{pmatrix} \in \bar U (A)
\end{equation}
The resulting classes in $K_* (A)$ do not depend on $p$ and 
$\lambda_n$, and are simply denoted $[e]$ and $[u]$. The natural 
inclusion $u \to u \oplus 1$ of $GL_n (A)$ in $GL_{n+1} (A)$ 
enables us to construct the inductive limit group
\[
\lim_{n \to \infty} \pi_0 \big( GL_n (A) \big)
\]
Similarly, the inclusion $e \to e \oplus 0$ of $M_n (A)$ in 
$M_{n+1} (A)$ leads to the inductive limit space \index{$K_0^+$}
\[
K_0^+ (A) := \lim_{n \to \infty} 
\pi_0 \big( \mr{Idem}\; M_n (A) \big)
\]
Actually this is an abelian semigroup with unit element 0.
By \cite[Lemma 7.4]{Phi2} it is naturally isomorphic to the 
monoid of equivalence classes of finitely generated projective 
$A$-modules. In this notation \cite[Theorem 7.7]{Phi2} becomes

\begin{thm}\label{thm:2.21}
Let $A$ be a unital Fr\'echet Q-algebra. The assignments 
$e \to [e]$ and $u \to [u]$ extend to natural isomorphisms 
\begin{align*}
G \left( K_0^+ (A) \right) &\isom K_0 (A) \\
\lim_{n \to \infty} \pi_0 \big( GL_n (A) \big) &\isom K_1 (A)
\end{align*}
where the G stands for \inde{Grothendieck group}.
\end{thm}

In particular $K_0 (A)$ has a natural ordering, for which 
$K_0^+ (A)$ is precisely the semigroup of positive elements.
These construction are especially important in connection
with \inde{density theorem} for $K$-theory
\cite[Th\'eor\`eme A.2.1]{Bos} :

\begin{thm}\label{thm:2.19}
Let $A$ and $B$ be Fr\'echet Q-algebras, and $\phi : A \to B$ 
a morphism with dense range. Suppose that $a \in A^+$ is 
invertible whenever $\phi^+ (a) \in B^+$ is invertible. 
Then for any $n \in \mh N$ the induced maps
\[
\begin{array}{lll}
\mr{Idem}\; M_n (A^+) & \to & \mr{Idem}\; M_n (B^+) \\
GL_n (A^+) & \to & GL_n (B^+)
\end{array}
\] 
are homotopy equivalences, and 
$K_* (\phi ) : K_* (A) \to K_* (B)$ is an isomorphism.
\end{thm}

The conditions are typically satisfied if $B$ is a unital 
Banach algebra, $A$ is a dense unital subalgebra which is 
Fr\'echet in its own finer topology, and 
$A \cap B^\times = A^\times$.

If we are working in m*-algebras then everywhere in the
above discussion we may replace invertibles by unitaries,
and idempotents by projections. This is a consequence of the
following elementary result.

\begin{lem}\label{lem:2.27} 
Let $A$ be a unital m*-algebra such that $sp (z^* z) \subset 
\mh R^+ \; \forall z \in A$. The set of unitaries in $A$ 
is a deformation retract of the set of invertibles in $A$. 
Likewise, the set of projections in $A$ is a deformation 
retract of the set of idempotents in $A$.
\end{lem}
\emph{Proof.}
Using Theorem \ref{thm:2.7}.3 write $|z| = (z^* z)^{1/2}$.
Then $z |z|^{-1}$ is unitary for every $z \in A^\times$ and
\[
[0,1] \times A^\times \to A^\times : 
(t,z) \to z |z|^{-t}
\]
is the desired deformation retraction. Similarly, there is a 
natural path from an idempotent to its associated Kaplansky 
projector, see e.g. \cite[Proposition 4.6.2]{Bla}. 
$\qquad \Box$ \\[3mm]

Quite often it is possible to find a bound on the size $n$
of matrices that we need to construct all $K_1$-classes.
To measure this we recall the notion of topological stable 
rank. Given a unital topological algebra $A$ define
\begin{align}
Lg_n (A) &:= \{ (a_1 ,\ldots,a_n ) \in A^n : 
A a_1 + \cdots + A a_n = A \}\\
tsr (A) &:= \inf\: \{ n : Lg_n (A) \mr{\;is\;dense\;in\;} A^n \}
\end{align}
Rieffel \cite{Rie1,Rie2} showed that this is useful for 
$K$-theory of $C^*$-algebras. The most general result in this
direction is \cite[Theorem 2.10]{Rie2} :

\begin{thm}\label{thm:2.22}
Let $A$ be a unital $C^*$-algebra. For any $n \geq tsr (A)$ 
we have
\[
\pi_0 \big( GL_n (A) \big) \cong K_1 (A)
\]
\end{thm}

To bound the topological stable rank of algebras that are not
too far from commutative we use the following tools, 
cf. \cite[Propostion 1.7]{Rie1} and \cite[Theorem 2.4]{OsTe} :

\begin{prop}\label{prop:2.23}
Let $X$ be a compact Hausdorff space and dim $X$ its covering
dimension. Also let $A \subset B$ be an inclusion of unital 
$C^*$-algebras, such that $B$ is a left $A$-module of rank $n$. 
Then
\begin{align*}
& tsr (C(X)) = 1 + \lfloor \dim X / 2 \rfloor\\
& tsr (B) \leq n \; tsr (A)
\end{align*}
\end{prop}

Together with Theorems \ref{thm:2.21} - \ref{thm:2.22} this
will allow us to realize the $K_1$-group of certain 
$C^*$-algebras entirely by invertible matrices, of a certain 
bounded size, with coefficients in a dense subalgebra.

Now we return to the study of the more abstract features of
the $K$-functor.

\begin{enumerate}
\item \label{it:2.2} Additivity. \index{additivity}
For any m-algebras $A_m \, (m \in \mh N)$
\[
K_n \left( \prod_{m=1}^\infty A_m \right) \cong  
\prod_{m=1}^\infty K_n (A_m)
\]

\item \label{it:2.3} Stability. \index{stability}
\[
K_n (\mf K \hot A) \cong K_n (M_m (A)) \cong  K_n (A)
\]

\item \label{it:2.4} Continuity. \index{continuity}
If $A_m \, (m \in \mh N)$ are Banach algebras and $A = 
\lim_{m \to \infty} A_m$ is their Banach inductive limit, then 
\[
K_n (A) \cong \lim_{m \to \infty} K_n (A_m)
\]

\item \label{it:2.5} Excision. \index{excision}
Let $0 \to A \to B \to C \to 0$ be an extension of m-algebras, 
admissible if not all the algebras are Fr\'echet.
There exists an exact hexagon
\[
\hexagon{K_0 (A)}{K_0 (B)}{K_0 (C)}{K_1 (C)}{K_1 (B)}{K_1 (A)}
\]

\item \label{it:2.6} Diffeotopy invariance. 
Let $f,g : A \to B$ be diffeotopic morphisms of m-algebras, or 
homotopic morphisms of Fr\'echet algebras. Then 
\[K_* (f) = K_* (g)
\]
\end{enumerate}

In this list \ref{it:2.4} is classical, but the author does not
know of any extension to Fr\'echet algebras. Proofs of
\ref{it:2.2}, \ref{it:2.3}, \ref{it:2.5} and \ref{it:2.6} can
be found in \cite{Cun1} and \cite{Phi2}. 

Obviously we will compare the features of $K_*$ with those of 
$HP_*$ given in Section \ref{sec:2.3}. Since topological 
$K$-theory is built with the completed projective tensor product,
it only makes sense to compare it with the cyclic theory with
the same topological tensor product. Hence, from now on
$HP_* (A)$ will mean $HP_* (A,\hot)$ for any m-algebra $A$, 
unless explicitly specified otherwise. 

First we deduce from the excision and diffeotopy properties 
that $K_*$ and $HP_*$ react in the same way on suspending an 
algebra. This is essentially a manifestation of Bott 
periodicity. By definition the \inde{smooth suspension} 
of $A$ is
\[
S_\infty (A) = \{ f \in C^\infty (S^1 ; A) : f(1) = 0 \}
\]

\begin{lem}\label{lem:2.25}
For $i = 0,1$ there are natural isomorphisms 
\[
\begin{array}{lll}
K_{1-i}(A) & \isom & K_i (S_\infty (A)) \\
K_0 (A) \oplus K_1 (A) & \isom & K_i (C^\infty (S^1 ; A))
\end{array}
\]
The same holds for periodic cyclic homology.
\end{lem}
\emph{Proof.}
Let (not an algebra! but easy to repair) 
\[
C_\infty (A) := \{ f \in C^\infty ([0,1];A) : f(0) = 0,
f^{(n)}(0) = f^{(n)}(1) \; \forall n > 0 \}
\]
be the smooth cone of $A$. Consider the admissible extension
\begin{equation}\label{eq:2.80}
0 \to S_\infty (A) \to C_\infty (A) \to A \to 0
\end{equation}
where the third map is evaluation at 0 and the second arrow
comes from composing a function with the surjection
\[
[0,1] \to S^1 : x \to e^{2 \pi i x}
\]
The boundary maps from the exact hexagon associated with
\eqref{eq:2.80} are the desired maps $K_{1-i}(A) \to 
K_i (S_\infty (A))$. To prove that these are isomorphisms
we will show that $C_\infty (A)$ is diffeotopy equivalent
to the algebra 0. Let $r : [0,1] \to [0,1]$ be a bijective
diffeomorphism with the properties
\begin{itemize}
\item $r'(x) > 0 \; \forall x \in (0,1)$
\item $r^{(n)}(0) = r^{(n)}(1) = 0 \; \forall n > 0$
\end{itemize} 
We have m-algebra homomorphisms
\begin{align*}
& C_\infty (A) \to C_0^\infty ([0,1],\{0\} ; A) : f \to f\\
& C_0^\infty ([0,1],\{0\} ; A) \to C_\infty (A) :
f \to f \circ r
\end{align*}
Since $r$ is diffeotopic to $\mr{id}_{[0,1]}$, both 
compositions of these algebra homomorphisms are 
diffeotopic to the respective identity homomorphisms. 
Hence we get natural isomorphisms
\[
K_* (C_\infty (A)) \isom 
K_* (C_0^\infty ([0,1],\{0\} ; A))
\]
However, $C_0^\infty ([0,1],\{0\} ; A)$ is diffeotopy 
equivalent to 0 by means of the homomorphisms
\begin{align*}
&\phi_t : C_0^\infty ([0,1],\{0\} ; A) \to  
C_0^\infty ([0,1],\{0\} ; A)\\
&\phi_t (f)(s) = f (ts)
\end{align*}
and therefore its $K$-theory vanishes.

Similarly there is an admissible extension
\begin{equation}\label{eq:2.81}
0 \to S_\infty (A) \to C^\infty (S^1 ; A) \to A \to 0
\end{equation}
But this extension splits, just send $a \in A$ to the
element $\tilde a \in C^\infty (S^1 ; A)$ with
$\tilde a (t) = a \; \forall t \in S^1$. Applying $K_i$
we get a split exact sequence of abelian groups
\[
0 \to K_i (S_\infty (A)) \to K_i (C^\infty (S^1 ; A)) 
\to K_i (A) \to 0
\]
Combined with the above this yields natural isomorphisms
\[
K_{1-i}(A) \oplus K_i (A) \isom K_i (S_\infty (A)) 
\oplus K_i (A) \isom K_i (C^\infty (S^1 ; A)) 
\]
The same proof applies with periodic cyclic homology.
$\qquad \Box$ \\[3mm]

Continuing our comparison, we see from Theorem 
\ref{thm:2.24}.2 that $HP_*$ is also $\mf K$-stable. Namely, 
we may take for $E$ all the matrices whose only nonzero entries 
are in the first column, and for $F$ all the matrices which 
have only zeros outside the first row. The isomorphism 
\begin{equation}
HP_* (A) \isom HP_* (\mf K \hot A)
\end{equation}
is induced by the algebra morphism $a \to p a p$, where 
$p \in \mf K$ is an arbitrary rank one projector. Its inverse
\begin{equation}
HP_* (\mf K \hot A) \isom HP_* (A)
\end{equation}
is a little more tricky, since it is not given by an algebra 
morphism, but a morphism of bicomplexes, the so-called 
\inde{generalized trace map}. This is the linear map
\begin{equation}
tr : \left( \mf K \hot A \right)^{\hot n} \to A^{\hot n}
\end{equation}
defined on elementary tensors by
\begin{equation}\label{eq:2.30}
tr (k_1 a_1 \otimes \cdots \otimes k_n a_n) =
tr (k_1 \cdots k_n) \, a_1 \otimes \cdots \otimes a_n
\end{equation}
Note that this works equally well if $\mf K$ is replaced 
by a finite dimensional matrix algebra $M_m (\mh C)$, cf. 
\cite[Section 1.2]{Lod}.

So now we know that $HP_*$ is halfexact, diffeotopy invariant 
and $\mf K$-stable. On the other hand, Cuntz' $kk$ 
\cite[Section 6]{Cun1} is the universal halfexact, 
diffeotopy invariant, $\mf K$-stable bivariant functor from 
m-algebras to abelian groups. This implies the existence of 
a unique natural transformation of functors
\begin{equation}\label{eq:2.31}
ch : K_* \to HP_*
\end{equation}
respecting these features. It is called the \inde{Chern 
character}, because it is a far-reaching generalization of the 
classical Chern character
\begin{equation}\label{eq:2.32}
Ch : K^* (X) \to \check H^* (X;\mh Q)
\end{equation}
that assigns to a complex vector bundle over a paracompact 
Hausdorff space $X$ a class in the even \v Cech cohomology.
Indeed, we can get \eqref{eq:2.32} for smooth manifolds by 
applying \eqref{eq:2.31} to $C^\infty (X)$ and using the 
isomorphism \eqref{eq:2.35}.

The Chern character is compatible with the countable additivity
of $K_*$ and $HP_*$, and also with excision, as was shown by
Nistor \cite[Theorem 1.6]{Nis2} :

\begin{thm}\label{thm:2.10}
Let $0 \to A \to B \to C \to 0$ be an extension of m-algebras.
The various Chern characters make a commutative diagram
\begin{center}
$\begin{array}{ccccccccccc}
K_1(A) & \!\to\! & K_1(B) & \!\to\! & K_1(C) & \!\to\! & K_0(A) 
& \!\to\! & K_0(B) & \!\to\! & K_0(C) \\
\downarrow & & \downarrow & & \downarrow & & \downarrow & & 
\downarrow & & \downarrow \\
HP_1(A) & \!\to\! & HP_1(B) & \!\to\! & HP_1(C) & \!\to\! & 
HP_0(A) & \!\to\! & HP_0(B) & \!\to\! & HP_0(C) 
\end{array}$
\end{center}
Moreover, if the extension is admissible and 
$\eta : K_0 (C) \to K_1 (A)$ and \\
$\partial : HP_0 (C) \to HP_1 (A)$ denote the connecting maps, 
then $ch \circ \eta = 2 \pi i \, \partial \circ ch$.
\end{thm}

Explicit formulas for the Chern character, from Phillips'
picture to the cyclic bicomplex, were first given by Karoubi
\cite[Chapitre II]{Kar2}. In the setting of \eqref{eq:2.82}
we may replace $A$ by $M_n (A)$ to achieve that $e, u \in A$.
Define
\begin{equation}
\begin{array}{cccl}
c_{2m+1}(e) & = & (-1)^m \ds{\frac{(2m)!}{m!}} \: 
  (e,\ldots,e) & \in A^{\otimes 2m+1} \\
c_{2m}(e) & = & (-1)^{m-1} \ds{\frac{(2m)!}{2(m!)}} \: 
  (e,\ldots,e) & \in A^{\otimes 2m} \\
c_{2m}(u) & = & (m-1)! \: (u^{-1},u,\ldots,u^{-1},u) &
  \in A^{\otimes 2m} \\
c_{2m+1}(u) & = & m! \: (2,u^{-1},u,\ldots,u^{-1},u) &
  \in A^{\otimes 2m+1}
\end{array}
\end{equation}
and place $c_n (e)$ in $CC^{per}_{1-n,n-1}(A,\hot) 
\cong A^{\hot n}$ and $c_n (u)$ in $CC^{per}_{2-n,n-1}(A,\hot) 
\cong A^{\hot n}$. By Lemma 2.1.6, Theorem 8.3.4 and Proposition
8.4.9 of \cite{Lod} we have
\begin{align}
ch \, [e] &= [ (c_n (e))_{n=1}^\infty ] \in HP_0 (A) \\
ch \, [u] &= [ (c_n (u))_{n=1}^\infty ] \in HP_1 (A) 
\end{align}

With the density theorem and the homotopy invariance of 
$K$-theory we can compute it for many Fr\'echet algebras, 
in particular commutative ones. 
The \inde{maximal ideal space} of a commutative m-algebra $A$ is 
defined like in algebraic geometry : it is the collection 
\inde{Max$(A)$} of all closed maximal ideals of $A$, endowed with 
the coarsest topology that makes all elements of $A$ into 
continuous functions on Max$(A)$. This is called the 
\inde{Gelfand topology}, and we denote it by $\mc T_G$. 

Contrarily to the $C^*$-algebra case, there may be several 
commutative Fr\'echet algebras with the same maximal ideal space.
The spectrum of a commutative Fr\'echet algebra is Hausdorff, 
$\sigma$-compact and paracompact, but it need not be locally 
compact. 
Therefore we also consider the \inde{compactly generated topology}
$\mc T_c$ on Max$(A)$. This means that we call $U \subset 
\mr{Max}(A)$ open in $\mc T_c$ if and only if $U \cap C$ is open
in $C$, for any compact $C$ with the relative topology from
$(\mr{Max}(A),\mc T_G )$. If $A^+$ is the unitization of $A$, then 
Max$(A^+) = \mr{Max}(A) \cup \{ A \}$ and we put
\begin{equation}
C_A := \{ f \in C( \mr{Max}(A^+), \mc T_c ): f(A) = 0 \}
\end{equation}

\begin{thm}\label{thm:2.11}
For any commutative Fr\'echet algebra $A$ there are natural 
isomorphisms 
\[
K_* (A) \cong K_* (C_A) \cong 
K^* \left( \mr{Max}(A^+), \{ A \} \right)\\
\]
If Max$ (A)$ is locally compact then 
\[
K_* (A) \otimes \mh Q \cong 
\check H^* (\mr{Max}(A^+),\{ A \} ) \otimes \mh Q 
\]
\end{thm}
\emph{Proof.}
The isomorphisms with integral coefficients are due to Phillips 
\cite[Theorem 7.15]{Phi2}. Here $K^*$ means \inde{representable 
$K$-theory} of topological spaces, in the sense of Karoubi 
\cite{Kar1}. He represents 
\begin{equation}\label{eq:2.37}
K^n (X) := [X,\mf F U^n]
\end{equation}
as the set of homotopy classes of continuous maps from a 
paracompact Hausdorff space $X$ to some classifying space 
$\mf F U^n$. This agrees with the usual definition if $X$ is 
compact, but in general it yields a generalized cohomology 
theory without compact supports. 
It has been known since the beginning of topological $K$-theory 
that \eqref{eq:2.32} gives an isomorphism
\begin{equation}\label{eq:2.36} 
Ch \otimes \mr{id}_{\mh Q} : K^* (X) \otimes \mh Q 
\isom \check H^* (X) \otimes \mh Q
\end{equation}
if $X$ is a finite CW-complex \cite[Section 2.4]{AtHi}. With 
spectral sequences, as in \cite{Seg1}, one can extend this to 
all compact Hausdorff spaces, since these are homotopy 
equivalent to CW-complexes. $\quad \Box$
\\[2mm]

This theorem can be considered as the counterpart in topological 
$K$-theory of Theorem \ref{thm:2.3}. If we apply it to a smooth 
manifold we get
\begin{equation}\label{eq:2.39}
K_* (C^\infty (X)) \otimes \mh C \cong 
\check H^* (X;\mh C)
\end{equation}
Since \v Cech cohomology agrees with De Rham cohomology
(both with complex coefficients) on the category of smooth 
manifolds, we deduce from \eqref{eq:2.35}, \eqref{eq:2.39} and 
the naturality of the Chern character that
\begin{equation}\label{eq:2.40}
ch \otimes \mr{id}_{\mh C} : 
K_* (C^\infty (X)) \otimes \mh C \isom HP_* (C^\infty (X))
\end{equation}
is an isomorphism. If we think a little more about this, it 
becomes clear that such an isomorphism should hold for many more 
algebras, even noncommutative ones. To make this precise, we 
introduce yet another category of topological algebras, denoted 
\inde{$\mc{CIA}$}. It is a full subcategory of the category of 
m-algebras $\mc{MA}$, and its objects are those $A \in \mc{MA}$ 
for which the Chern character induces an isomorphism
\begin{equation}\label{eq:2.41}
ch \otimes \mr{id}_{\mh C} : 
K_* (A) \otimes \mh C \isom HP_* (A)
\end{equation}
\begin{prop}\label{prop:2.8}
The category $\mc{CIA}$ is closed under the following operations:
\begin{enumerate}
\item countable direct products
\item tensoring with $M_m (\mh C)$ or $\mf K$
\item diffeotopy equivalences
\item admissible extensions, quotients and ideals
\end{enumerate}
\end{prop}
\emph{Proof.}
1,2 and 3 follow directly from the features of $K_*$ and $HP_*$
on pages \pageref{p:hp} and \pageref{it:2.4}. As concerns 4, by 
Theorem \ref{thm:2.10} we can apply Lemma \ref{lem:2.12} to the 
functors $K_* (\,\cdot\,) \otimes \mh C$ and $HP_* (\,\cdot
\,, \hot)$ on the category $\mc{MA}$ with admissible morphisms.
The factor $2 \pi i$ in Theorem \ref{thm:2.10} is inessential.
$\qquad \Box$ \\[2mm]

In view of these similarities, it is logical to try to extend 
the material from Section \ref{sec:2.2} to topological 
$K$-theory. However, this is somewhat problematic, as general
compact Hausdorff spaces are much less easy to handle then 
algebraic varieties. In the next section we will avoid
these difficulties by considering only smooth manifolds with 
a finite group action. Right now we will prove a coarse 
analogue of Theorem \ref{thm:2.4}, which applies to semisimple 
algebras which "live" on finite simplicial complexes. We 
formulate it in terms of $C^*$-algebras, but with Theorem 
\ref{thm:2.19} it can easily be generalized to certain 
Fr\'echet algebras.

\begin{prop}\label{prop:2.26}
Let $\Sigma$ be a finite simplicial complex and $\phi : A \to B$ 
a homomorphism of $C^*$-algebras. Suppose that
\begin{itemize}
\item there are unital homomorphisms from $C (\Sigma )$ to the
centers of the multiplier algebras $\mc M (A)$ and $\mc M (B)$
\item for every simplex $\sigma$ of $\Sigma$ there are finite
dimensional $C^*$-algebras $A_\sigma$ and $B_\sigma$ such that
\[
A C_0 (\sigma ,\delta \sigma ) \cong C_0 (\sigma ,\delta \sigma 
;A_\sigma ) \quad \mr{and} \quad B C_0 (\sigma ,\delta \sigma ) 
\cong C_0 (\sigma ,\delta \sigma ; B_\sigma )
\]
\item $\phi$ is $C(\Sigma )$-linear
\item for every $x_\sigma \in \sigma \setminus \delta \sigma$ 
the localization $\phi (x_\sigma ) : A_\sigma \to 
B_\sigma$ induces an isomorphism on $K$-theory
\end{itemize}
Then 
\[
K_* (\phi ) : K_* (A) \isom K_* (B)
\]
is an isomorphism.
\end{prop}
\emph{Proof.}
Let $\Sigma^n$ be the $n$-skeleton of $\Sigma$ and consider 
the ideals
\begin{equation}
\begin{aligned}
& C(\Sigma ) = I_0 \supset I_1 \supset \cdots \supset I_n 
\supset \cdots \\
& I_n = C_0 (\Sigma ,\Sigma^n )
\end{aligned}
\end{equation}
They give rise to ideals $A_n = A I_n$ and $B_n = B I_n$. 
Because $\Sigma$ is finite all these ideals are 0 for large 
$n$. We can identify
\[
A_{n-1} / A_n \cong A C_0 (\Sigma^n ,\Sigma^{n-1}) \cong
\bigoplus_{\sigma \in \Sigma \,:\, \dim \sigma = n} 
\!\!\!\! A C_0 (\sigma ,\delta \sigma) := 
\bigoplus_{\sigma \in \Sigma \,:\, \dim \sigma = n} 
\!\!\!\! C_0 (\sigma ,\delta \sigma ;A_\sigma )
\]
and similarly for $B$.
Because $\phi$ is $C(\Sigma )$-linear, it induces homomorphisms
\[
\phi (\sigma ) : C_0 (\sigma ,\delta \sigma ; A_\sigma ) \to
C_0 (\sigma ,\delta \sigma ; B_\sigma )
\]
By Lemma \ref{lem:2.12} and the additivity of $K$-theory it 
suffices to show that every $\phi (\sigma )$ induces an 
isomorphism on $K$-theory. Let $x_\sigma$ be any interior point of 
$\sigma$. Because $\sigma \setminus \delta \sigma$ is contractible, 
$\phi_\sigma$ is homotopic to $\mr{id}_{C_0 (\sigma ,\delta \sigma )}
\otimes \phi (x_\sigma )$. By assumption the latter map induces an 
isomorphism on $K$-theory. With the homotopy invariance of $K$-theory
it follows that $K_* (\phi (\sigma ))$ is an isomorphism. 
$\qquad \Box$ \\[3mm]
Note that this proof applies equally well to the functor
$K_* ( \cdot ) \otimes_{\mh Z} \mh Q$.

\section{Equivariant cohomology and algebras of invariants}
\label{sec:2.5}

This section is all about algebras that carry an action of a 
finite group, and their subalgebras of invariant elements. 
To place things in a classical context we first recall some
beautiful theorems on equivariant topological $K$-theory and 
cyclic theory for crossed product algebras. 

All the results after that are due to the author, but some of
them already appeared in \cite{Sol}. We broaden our view, to 
algebras of the form
\[
A^G = C^\infty (X ; M_N (\mh C))^G
\]
with $X$ a smooth manifold. We relate $K$-theory of such 
algebras to a $G$-equivariant cohomology theory due to Bredon 
\cite{Bre}, which can also be described as the \v Cech 
cohomology of a certain sheaf over the orbifold $X/G$.
This depends on the existence of a $G$-equivariant 
triangulation of $X$. Using the same Mayer-Vietoris type of
arguments we also prove that the Chern character for $A^G$ 
becomes an isomorphism after tensoring with $\mh C$. Finally,
if $X$ happens to be a complex affine variety, then we show
that the "polynomial" subalgebra $A^G_{alg}$ has the same 
periodic cyclic homology as $A^G$.
\\[2mm]

Let $G$ be a topological group acting continuously on a 
topological space $X$. Then $G$ also acts continuously on the 
closed subspace
\begin{equation}\label{eq:2.42}
\widetilde X := \{ (g,x) \in G \times X : g x = x \}
\end{equation}
of $G \times X$ by
\begin{equation}
g (g',x) = (g g' g^{-1}, gx)
\end{equation}
and $\widetilde X / G$ is called the \inde{extended quotient} of 
$X$ by $G$. In the literature one often encounters the notation 
$\hat X$ for \inde{$\widetilde X$}, but we avoid this because it 
might be confused with the spectrum of a topological group.

Let \inde{$\langle G \rangle$} be the set of conjugacy classes 
in $G$, and denote the class containing $g$ by 
\inde{$\langle g \rangle$}. We have a decomposition
\begin{equation}\label{eq:2.43}
\widetilde X / G \cong \bigsqcup_{\langle g \rangle \in \langle G 
\rangle} (g, X^g / Z_G (g)) \cong  \bigsqcup_{\langle g \rangle 
\in \langle G \rangle} X^g / Z_G (g)
\end{equation}
where \inde{$Z_G (g)$} is the centralizer of $g$ in $G$ and
\index{$X^g$}
\begin{equation}
X^g = \{ x \in X : g x = x \} 
\end{equation}
Notice that the components of this partition are always closed, 
and they are open if $G$ is finite.

A $G$-vector bundle over $X$ is a vector bundle $p: V \to X$ 
together with an action of $G$ on $V$, such that 
$\forall v \in V, x \in X, g \in G$ \index{$G$-vector bundle}
\begin{itemize}
\item $p (gv) = g p(v)$
\item $g : p^{-1}(x) \to p^{-1}(gx)$ is linear
\end{itemize}

If $X$ and $G$ are both compact Hausdorff, then it makes sense
to consider the Grothendieck group of equivalence classes of
complex $G$-vector bundles. This group was first studied by Atiyah 
\cite{Ati}, and it is denoted by $K_G^0 (X)$.
By the same suspension procedure as in the nonequivariant case
this leads to a sequence of functors $K^n_G$, together called
\inde{equivariant $K$-theory}. This is a equivariant cohomology
theory which shares most of the properties of ordinary 
topological $K$-theory.

For any $g \in G$, the restriction of a complex $G$-bundle
$p: V \to X$ to $X^g$ is a vector bundle on which $g$ acts
linearly in every fiber. So we can decompose it canonically 
into to its $g$-eigenspaces :
\begin{equation}\label{eq:2.44}
V \big|_{X^g} = \bigoplus_i V_i := \bigoplus_i
\left\{ v \in p^{-1}(X^g) : g v = \lambda_i v \right\}
\end{equation}
By the continuity of the action and the compactness of $X^g$ 
there are only finitely many $\lambda_i \in \mh C$ for which 
$V_i$ is nonzero. From this decompostion we cook a canonical map 
\begin{equation}\label{eq:2.45}
\rho_g : K_G^* (X) \to K^* (X^g) \otimes \mh C
\end{equation}
sending $[V]$ to $\sum_i \lambda_i [V_i]$.
All these $\rho_g$'s together combine to a map that classifies
$G$-bundles over $X$ in terms of ordinary vector bundles over
the extended quotient $\widetilde X / G$. Indeed, for finite $G$ 
the identification 
\begin{equation}
K^* (\widetilde X) \cong \bigoplus_{\langle g \rangle \in 
\langle G \rangle} K^* (X^g)
\end{equation}
gives a map
\begin{equation}\label{eq:2.46}
\rho := \sum_{g \in G} \rho_g : K_G^* (X) \to 
K^* (\widetilde X ) \otimes \mh C 
\end{equation}
It is easy to see that the image of $\rho$ is contained in the 
subspace of $G$-invariants, so if we compose it with the Chern 
character for $\widetilde X$ we land in $\check H^* (\widetilde X 
; \mh C )^G$, which by \cite[Corollaire 5.2.3]{Gro2} is naturally 
isomorphic to $\check H^*(\widetilde X / G; \mh C )$. 
By the way, this composition 
\begin{equation}\label{eq:2.47}
Ch_G := (Ch \otimes \mr{id}_{\mh C}) \circ \rho
\end{equation}
is called the \inde{equivariant Chern character}. The punchline
is of course \cite[Theorem 1.19]{BaCo} :

\begin{thm}\label{thm:2.13}
For any finite group $G$ acting on a compact Hausdorff space $X$
there are natural isomorphisms
\begin{align*}
\rho \otimes \mr{id}_{\mh C} &: K_G^* (X) \otimes \mh C
\isom \left( K^* (\widetilde X ) \otimes \mh C \right)^G\\
Ch_G \otimes \mr{id}_{\mh C} &: K_G^* (X) \otimes \mh C
\isom \left( \check H^* (\widetilde X ; \mh C ) \right)^G \cong
\check H^* (\widetilde X / G; \mh C )
\end{align*}
\end{thm}

We switch back to a more algebraic point of view. Suppose that 
the compact group $G$ acts by *-automorphisms on a 
$C^*$-algebra $A$. The above leads us to consider finitely 
generated projective $A$-modules $M$ with a $G$-action 
satisfying $g (a m) = (ga) (gm)$. The Grothendieck group of 
equivalence classes of such modules is the equivariant 
$K$-theory $K_0^G (A)$. This gives rise to sequence of functors 
$K_n^G$ which by the equivariant Serre-Swan theorem 
\cite[Theorem 2.3.1]{Phi1} are related to the above homonymous 
functors as 
\begin{equation}\label{eq:2.48}
K_*^G (C(X)) \cong K^*_G (X)
\end{equation}
We already managed to describe $G$-bundles in terms of
ordinary vector bundles over a related space, and the same is 
possible for $C^*$-algebras. Namely, Julg \cite{Jul} showed 
that there is a natural identification
\begin{equation}\label{eq:2.49}
K_*^G (A) \cong K_* (A \rtimes G)
\end{equation}
Combining \eqref{eq:2.48} and \eqref{eq:2.49} with Theorem 
\ref{thm:2.13} we see that 
\begin{equation}\label{eq:2.50}
K_* (C(X) \rtimes G) \otimes \mh C 
\cong \check H^* (\widetilde X ;\mh C )^G
\cong \check H^* (\widetilde X / G; \mh C )
\end{equation}
This is important for our purposes, since it shows that the 
homology of a crossed product algebra can be described in 
geometric terms. This can even be refined in cyclic 
theory. Let $X$ be either a nonsingular complex affine
variety or a smooth real manifold, not necessarily compact,
and $A$ the algebra of either regular or smooth functions on $X$.
Suppose that the action of $G$ preserves this structure, so that
the partially invariant subspaces $X^g$ are of the same type as 
$X$. Brylinski \cite{Bry1,Bry2} proved that
\begin{align}
\label{eq:2.23} HH_n (A \rtimes G) 
&\cong \bigoplus_{\langle g \rangle \in \langle G \rangle} 
\Omega^n \left( X^g \right)^{Z_G (g)} \;\cong\; \Omega^n \bigl( 
\widetilde X \bigr)^G \\
\label{eq:2.24} HP_n (A \rtimes G) 
&\cong \bigoplus_{\langle g \rangle \in \langle G \rangle}
\prod_{m \in \mh Z} H_{DR}^{n+2m} \left( X^g \right)^{Z_G (g)} 
\;\cong\; \prod_{m \in \mh Z} H_{DR}^{n+2m} 
\bigl( \widetilde X \bigr)^G \\
\label{eq:2.25} HC_n (A \rtimes G) 
&\cong \bigoplus_{\langle g \rangle \in \langle G \rangle} 
\!\!\!\Bigl( \Omega^n (X^g) / d \Omega^{n-1}(X^g) \oplus 
H_{DR}^{n-2}(X^g) \oplus H_{DR}^{n-4}(X^g) \oplus \cdots 
\!\Bigr)^{Z_G (g)} \\
\nonumber &\cong\; \left( \Omega^n \bigl( \widetilde X \bigr) 
  / d \Omega^{n-1} \bigl( \widetilde X \bigl) \oplus H_{DR}^{n-2} 
  \bigl( \widetilde X \bigr) \oplus H_{DR}^{n-4} \bigl( 
  \widetilde X \bigr) \oplus \cdots \right)^G
\end{align} 
Moreover Nistor \cite[Theorem 2.11]{Nis4} 
constructed an explicit map 
\begin{equation}
HH_n (A \rtimes G) \to \bigoplus_{\langle g \rangle \in 
\langle G \rangle} \Omega^n (X^g)
\end{equation}
(not the naive restriction!) and showed that it induces these 
isomorphisms.

When we combine \eqref{eq:2.50},\eqref{eq:2.24} and Theorem 
\ref{thm:2.19} with the naturality of the Chern character, 
we arrive at

\begin{thm}\label{thm:2.14}
Suppose that a finite group $G$ acts by diffeomorphisms on a
compact smooth manifold $X$. Then the Chern character gives
an isomorphism
\[
ch \otimes \mr{id}_{\mh C} : K_* (C^\infty (X) \rtimes G)
\otimes \mh C \isom HP_* (C^\infty (X) \rtimes G)
\]
\end{thm}

Later we will see that the compactness assumption in this
theorem is not necessary, so we drop it now, at least for the 
rest of this section.

Interestingly, the crossed product $C^\infty (X) \rtimes G$ 
can also be realized as an algebra of invariants :
\begin{equation}
C^\infty (X) \rtimes G \cong 
C^\infty \big( X; \mr{End}(\mh C [G]) \big)^G
\end{equation}
where the group algebra $\mh C [G]$ carries the right regular 
representation of $G$, see Lemma \ref{lem:A.3}. We propose 
to study such algebras also with other $G$-representations 
instead of $\mh C [G]$. Then $\mh C [G]$ is universal in the 
sense that it contains every irreducible $G$-representation. 
At the other extreme we have the trivial representation, 
which leads us to the Fr\'echet algebra
\begin{equation}\label{eq:2.21}
C^\infty (X/G) := C^\infty (X)^G 
\end{equation}
of smooth functions on the \inde{orbifold} $X / G$, cf. 
\cite{Sat}. Wassermann \cite[Section IV]{Was2} showed that the 
periodic cyclic homology of this algebra equals, as one would 
expect, the \v Cech cohomology of $X/G$:
\begin{equation}\label{eq:2.51}
HP_* \left( C^\infty (X)^G \right) \cong H_{DR}^* (X)^G \cong 
\check H^* (X/G;\mh C )
\end{equation}
But he also noticed that Hochschild homology does not behave so 
well in this case, as
\[
HH_* \left( C^\infty (X)^G \right) \quad \mr{and} \quad
\Omega^* (X)^G
\]
are not isomorphic in general.

Let $Z \subset Y$ be arbitrary subsets of $\mh R^n$, and $V$ a
Fr\'echet space. To include manifolds with boundary in our 
studies we adhere to the following conventions :
\index{$C_0^\infty (Y,Z)$}
\begin{equation}
\begin{aligned}
C^\infty (Y) &:= \left\{ f \big|_Y : f \in C^\infty (U) 
\;\mr{for \; some \; open} \;U\; \mr{with} \; Y \subset U \subset 
\mh R^n \right\} \\
C_0^\infty (Y,Z) &:= \left\{ f \in C^\infty (Y) : 
f \big|_Z = 0 \right\} \\
C_0^\infty (Y,Z;V) &:= C_0^\infty (Y,Z) \hot V
\end{aligned}
\end{equation}
Unfortunately this slightly ambiguous for orbifolds embedded in 
$\mh R^n$. For example if $Y = \mh R / \{ \pm 1 \}$, identified
as a topological space with $[0, \infty )$, then 
\[
C^\infty \big( [0,\infty ) \big) \supsetneq 
C^\infty (\mh R )^{\{ \pm 1 \}}
\]
since the right hand side contains only functions whose odd
derivatives vanish at 0. However, the difference is not too big,
since both algebras are diffeotopy equivalent to 
$\mh C \oplus C_0^\infty \big( \mh R ,(-\infty,0] \big)$
via $f \to f \circ \phi$, where 
$\phi \in C^\infty_0 \big( \mh R ,(-\infty,0] \big)$ 
is an automorphism of $[0,\infty)$ which is diffeotopic to the
identity of $[0, \infty)$. 
In such situations we shall usually give priority to the 
orbifold structure and use \eqref{eq:2.21} as a definition, 
at least locally. 

Now let $X$ be a smooth manifold with boundary, still 
$\sigma$-compact, and consider the Fr\'echet algebra 
\begin{equation}
A := C^\infty (X;M_N (\mh C )) \cong M_N (C^\infty (X))
\end{equation}
We assume we have elements \inde{$u_g$} $\in A^\times$ and 
diffeomorphisms $\alpha_g$ of $X$ such that 
\begin{equation}\label{eq:2.66}
g a (x) = u_g (x) a (\alpha_g^{-1} x) u_g^{-1}(x)
\end{equation}
defines an action of $G$ on $A$. Although this implies that 
$g \to \alpha_g$ is a group homomorphism, $g \to u_g$ need not
be one. The algebra 
\index{$A^G$}\index{finite type algebra!Fr\'echet}
\begin{equation}\label{eq:2.52}
A^G = C^\infty (X; M_N (\mh C ))^G
\end{equation}
will be our Fr\'echet version of a finite type algebra. Clearly 
$A^G$ is finitely generated as a module over $C^\infty (X)^G$.
It follows from a classical theorem of Newman that the set 
of points of $X$ whose $G$-stabilizer equals ker $\alpha$ is 
open and everywhere dense in $X$, see \cite[Theorem 1]{Dre}. 
Hence, if $\ker \alpha = \{e\}$ then $Z (A^G ) = C^\infty (X)^G$.

To compute its $K$-theory we will use an ``\inde{equivariant 
cohomology} theory with a \inde{local coefficient system}'', 
as defined by Bredon \cite{Bre}. This theory can be combined 
with the ideas of Segal \cite{Seg2} and S\l omi\'nska \cite{Slo1}
to describe $K_* \left( A^G \right)$ in sheaf-cohomological terms.

First we recall some of Bredon's constructions, referring to 
\cite{Bre} for more precise information. Let $\Sigma$ be a 
countable, locally finite and finite dimensional $G$-CW complex. 
Assume that all cells are oriented and that the action of $G$ 
preserves these orientations.
\index{Bredon cohomology} \index{Sigma@$\Sigma$}

We define a category $\mc K$ whose objects are the finite 
subcomplexes of $\Sigma$. The morphisms from $K$ to $K'$ are the 
maps $K \to K' : x \to gx$ for $g \in G$ such that 
$gK \subset K'$. Now a local coefficient system \inde{$\mf L$} 
on $\Sigma$ is a covariant functor \index{Kmc@$\mc K$}
from $\mc K$ to the category of abelian groups, and the group
$C^q (\Sigma;\mf L)$ of $q$-cochains is the set of all functions
$f$ on the $q$-cells of $\Sigma$ with the property that 
$f(\tau) \in \mf L(\tau) \, \forall \tau$.
Furthermore we define a coboundary map $\partial : C^q (\Sigma 
;\mf L ) \to C^{q+1} (\Sigma ;\mf L )$ by
\begin{equation}\label{eq:2.53}
(\partial f) (\sigma) = \sum_\tau [\tau : \sigma ] 
\mf L(\tau \to \sigma) f(\tau)
\end{equation}
where the sum runs over all $q$-cells $\tau$ and the incidence
number $[\tau : \sigma ]$ is the degree of the attaching map 
from $\partial \sigma$ (the boundary of $\sigma$ in the 
standard topological sense) to $\tau / \partial \tau$. 
The group $G$ acts naturally on this complex by cochain maps, 
so, for any $K \subset \Sigma \,,\, \left( C^* (K;\mf L )^G  
,\partial \right)$ is a differential complex and we can define 
the equivariant cohomology of $K$ with coefficients in $\mf L$ as
\begin{equation}
H^q_G (K;\mf L) := H^q \left( C^* (K;\mf L)^G ,\partial \right)
\end{equation}
More generally for $K' \subset K ,\: C^* (K, K'; \mf L )$ 
is the kernel of the restriction map 
$C^* (K; \mf L) \to C^*(K'; \mf L)$ and 
\begin{equation}
H^q_G (K,K';\mf L ) = 
H^q \left( C^* (K,K';\mf L )^G ,\partial\right)
\end{equation}
By construction there exists a local coefficient 
system $\mf L^G$ (more or less consisting of the $G$-invariant 
elements of $\mf L$) on the CW-complex $\Sigma / G$ such that 
the differential complexes $(C^* (K,K';\mf L )^G ,\partial)$ and
$(C^* (K/G, K'/G; \mf L^G ), \partial)$ are isomorphic. 
Notice that $\mf L^G$ defines a sheaf over $\Sigma / G$ (with
the cells as cover), so that
\begin{equation}\label{eq:2.54}
H^q_G (K,K';\mf L ) \cong 
\check H^q \left( K/G, K'/G; \mf L^G \right)
\end{equation}

Let $\Sigma^p$ be the $p$-skeleton of $\Sigma$. We capture all 
the above things in a spectral sequence $(E_r^{p,q})_{r \geq 1}$, 
degenerating already for $r \geq 2$, as follows : 
\begin{align}
E_1^{p,q} = H^{p+q}_G (\Sigma^p ,\Sigma^{p-1} ;\mf L ) &= 
 \left\{ \begin{array}{l@{\;\;\quad \mr{if} \quad}l}
   C^p (\Sigma ; \mf L )^G & q = 0 \\ 
   0                       & q > 0 
 \end{array} \right. \label{eq:2.55} \\
E_2^{p,q} &= \left\{ \begin{array}{l@{\qquad \mr{if} \quad}l}
   H^p_G (\Sigma; \mf L ) & q = 0 \\ 
   0                      & q > 0 
 \end{array} \right. \label{eq:2.56}
\end{align}
The differential $d_1^E$ is the composition
\begin{equation}\label{eq:2.57}
E_1^{p,q} \to H^{p+q}_G(\Sigma^p ; \mf L ) \to E_1^{p+1,q}
\end{equation}
of the maps induced by the inclusion $(\Sigma^p ,\emptyset ) \to 
(\Sigma^p ,\Sigma^{p-1})$ and the coboundary $\partial$.

Now let $B^G$ be an algebra like \eqref{eq:2.52}, but without
the differentiable structure. We will see later that this algebra
has the same $K$-theory as $A^G$, for a suitable triangulation
of $X$. So we put 
\begin{equation}
B = C(\Sigma ; M_N (\mh C )) = M_N (C (\Sigma )) 
\end{equation}
and we assume that we have $u_g \in B^\times$ such that 
\begin{equation}
g b (x) = u_g (x) b (g^{-1} x) u_g^{-1}(x)
\end{equation}
defines an action of $G$ on $B$. To associate a local coefficient 
system \inde{$\mf L_u$} to this algebra we first assume 
that $K$ is connected. In that case we let
\begin{equation}\label{eq:2.62}
G_K := \{ g \in G : g x = x \quad \forall x \in K \}
\index{$G_K$}
\end{equation}
be the isotropy group of $K$ and we define $\mf L_u (K)$ 
to be the free abelian group on the (equivalence classes of) 
irreducible projective $G_K$-representations contained in 
$(\pi_x, \mh C^N)$, where $\pi_x (g) = u_g (x)$ for $g \in G_K , 
x \in K$. By the continuity of the $u_g$ we get the same group 
for any $x \in K$. If $K$ is not connected, then we let 
$\{K_i \}_i$ be its connected components, and we define
\begin{equation}\label{eq:2.63}
\mf L_u (K) = \prod_i \mf L_u (K_i)
\end{equation}
Suppose that $g K \subset K'$ and that $\rho$ is a 
projective $G_K$-representation. Then we define a 
projective $G_{K'}$-representation by
\begin{equation}\label{eq:2.64}
\mf L_u (g: K \to K') \rho (g') = \rho (g^{-1} g' g) 
\quad g' \in G_{K'}
\end{equation}
If $h \in G$ gives the same map from $K$ to $K'$ as $g$ then 
$h^{-1} g \in G_K$ and
\begin{equation}
\mf L_u (h: K \to K') \rho (g') = \rho (h^{-1} g' h) =
\rho (h^{-1} g) \rho  (g^{-1} g' g) \rho (g^{-1} h)
\end{equation}
so $\mf L_u (h: K \to K') \rho$ is isomorphic to 
$\mf L_u (g: K \to K') \rho$ as a projective representation. 
This makes $\mf L_u$ into a functor.

Suppose for example that $u_g (x) = 1 \, \forall x \in \Sigma, 
g \in G$. Then $\mf L_u$ and $\mf L_u^G$ are the constant 
sheaves $\mh Z$ over $\Sigma$ and $\Sigma/G$ respectively, and 
\begin{equation}
H^*_G (\Sigma;\mf L_u) \cong \check H^*(\Sigma/G; \mh Z)
\end{equation}
is the ordinary cellular cohomology of $\Sigma/G$. 
Moreover $B^G \cong C(\Sigma /G; M_N (\mh C ))$, so 
$K_* \left( B^G \right) \cong K^* (\Sigma /G)$, which is 
isomorphic to $\check H^* (\Sigma /G ;\mh Z )$ modulo torsion.

Of most interest is the case where $B^G \cong C(\Sigma ) 
\rtimes G$ is the crossed product, as in Lemma \ref{lem:A.3}. 
Then we compare $\mf L_u^G$ to the direct image of the 
constant sheaf $\mh Z$ on $\widetilde \Sigma$, under the 
canonical map $p: \widetilde \Sigma / G \to \Sigma / G$. 
Although they may not always be isomorphic, their \v Cech 
complexes are identical, for any cover that refines the cell 
structure. Since $p$ is finite to one we can deduce, using 
even more spectral sequences \cite[Chapitre 5]{God}, that 
\begin{equation}\label{eq:2.65}
H^*_G (\Sigma;\mf L_u \otimes \mh Q) \cong 
\check H^* \big( \widetilde \Sigma /G; \mh Q \big)
\end{equation}
and by Theorem \ref{thm:2.14} this is the same as
$K_* \left( B^G \right) \otimes \mh Q$ if $\Sigma$ is 
compact, i.e. if it is a finite CW-complex.

It turns out that this close relation between 
$K_* \left( B^G \right)$ and the \v Cech cohomology 
$H^* \left( \Sigma /G;\mf L_u^G \right)$ is valid in general. 
Consider the following analogue of \eqref{eq:2.59}
\begin{equation}\label{eq:2.58} 
\begin{split}
&K_*^0 \left( B^G \right) = K_*^0 \left( B^G \right) \supset 
K_*^1 \left( B^G \right) \supset \cdots \supset 
K_*^{\dim \Sigma} \left( B^G \right) \supset 
K_*^{1 + \dim \Sigma} \left( B^G \right) = 0 \\
&K_*^p \left( B^G \right) := 
\ker \left( K_* \left( C(\Sigma ;M_N (\mh C ) )^G \right) \to 
K_* \left( C(\Sigma^{p-1}; M_N(\mh C ) )^G \right) \right)
\end{split}
\end{equation}

\begin{thm}\label{thm:2.15}
The graded group associated with the filtration \eqref{eq:2.58}
is isomorphic to $\check H^* \left( \Sigma/G ; \mf L_u^G \right)$.
In particular there is an (unnatural) isomorphism
\begin{equation}\label{eq:2.85}
K_* \left( B^G \right) \otimes \mh Q \cong \check H^* 
\left( \Sigma /G; \mf L_u^G \otimes \mh Q \right)
\end{equation}
and
\[
K_* \left( B^G \right) \cong 
\check H^* \left( \Sigma /G; \mf L_u^G \right)
\]
if one of both sides is torsion free.
\end{thm}
\emph{Proof.}
Using \cite[Section XV.7]{CaEi} we construct a spectral sequence
$(F_r^{p,q})_{r \geq 1}^{q \in \mh Z / 2 \mh Z}$ with the 
following terms:
\begin{equation}
\begin{aligned}
F_1^{p,q} &= K_{p+q} \left( C_0 (\Sigma^p / \Sigma^{p-1}; M_N
(\mh C ) )^G \right) \\
F_2^{p,q} &= H^p (\Sigma /G ;\mc K^q_u ) \\
F_\infty^{p,q} &= K_{p+q}^p \left( B^G \right) / 
K_{p+q}^{p+1} \left(B^G \right)
\end{aligned}
\end{equation}
where $\mc K^q_u (\sigma ) = K_* \left( C(G \sigma ; 
M_N (\mh C) )^G \right)$.
Now replace $\mf L$ in \eqref{eq:2.55} by $\mf L_u$ and sum over
all $q$. If we compare the result with $F_1^p = F_1^{1,0} \oplus
F_1^{1,1}$ we see that $E_1^p \cong F_1^p$. So we get a diagram
like \eqref{eq:2.16} :
\begin{equation}\label{eq:2.61}
\begin{array}{ccc}
F_1^{p,q} & \xrightarrow{\quad d_1^F \quad} & F_1^{p+1,q} \\[1mm]
\cong & & \cong \\
\prod\limits_{n \in \mh Z} E_1^{p,q+2n} & \xrightarrow{\quad 
  d_1^E \quad} & \prod\limits_{n \in \mh Z} E_1^{p+1,q+2n}
\end{array}
\end{equation}
The differential $d_1^F$ for $F_1$ is induced from the 
construction of a mapping cone of a Puppe sequence in the 
category of $C^*$-algebras. 
This is the noncommutative counterpart of the construction of 
the differential in cellular homology, so by naturality $d_1^F$ 
corresponds to $d_1^E$ under the above isomorphism. Therefore 
the spectral sequences $E_r^p$ and $F_r^p$ are isomorphic, and
in particular $F_r$ degenerates for $r \geq 2$. Now the 
isomorphism \eqref{eq:2.85} follows from \eqref{eq:2.54}.

If either $K_* \left( B^G \right)$ or $\check H^* \left( 
\Sigma /G; \mf L_u^G \right)$ is torsion free, then every term 
$E_\infty^{p,q} \cong F_\infty^{p,q}$ must be torsion free. 
Hence in this case both $K_* \left( B^G \right)$ and 
$\check H^* \left( \Sigma /G; \mf L_u^G \right)$ are free 
abelian groups, of the same rank. $\qquad \Box$ 
\\[2mm]

The main use of this theorem is really to compute 
$K_* \left( B^G \right)$, for now the extensive machinery of 
\v Cech cohomology becomes available.

The requirements \eqref{eq:2.62} - \eqref{eq:2.64} allow us to 
construct the sheaves $\mf L_u$ and $\mf L_u^G$ 
without reference to the cellular structure of $\Sigma$. If 
we do this for the algebra $A^G$ of \eqref{eq:2.52} then the 
stalk of \inde{$\mf L_u$} over $x \in X$ is the free abelian 
group on the (equivalence classes of) irreducible projective
$G_x$-representations contained in $(\pi_x, \mh C^N)$ and the 
$G$-action on $\mf L_u$ is determined by \eqref{eq:2.64}.  
A section $s$ is continuous at $x$ if there exists a 
neighborhood $U$ of $x$ such that $\forall y \in U$ :
\begin{itemize}\label{p:lug}
\item $G_y \subset G_x$
\item $s(x) = s(y)$ as virtual projective $G_y$-representations
\end{itemize}
This $\mf L_u$ is a generalization of the sheaf constructed in
\cite[\S 2]{BaCo}. Clearly, the subsheaf $\mf L_u^G$ of 
$G$-invariant continuous sections descends to a sheaf on $X/G$.  

To relate this sheaf to the $K$-theory and periodic cyclic 
homology of $A^G$ we need two preparatory results. 
The first is a weak version of Theorem \ref{thm:2.11}, 
which however does include $HP_*$.

\begin{lem}\label{lem:2.16}
Let $U \subset \mh R^n$ be an open bounded star-shaped set. 
The Fr\'echet algebra $C_0^\infty (\mh R^n, \mh R^n \setminus U)$ 
belongs to $\mc{CIA}$ and 
\[
K_* (C_0^\infty (\mh R^n, \mh R^n \setminus U)) \cong 
\check H^* (\mh R^n, \mh R^n \setminus U; \mh Z) \cong \mh Z
\] 
is concentrated in degree $n$.
\end{lem}
\emph{Proof.}
Clearly we may assume that 0 is the center of $U$. Let $P$ be 
the point of the $n$-sphere corresponding to infinity under 
the stereographic projection $S^n \to \mh R^n$. By assumption 
$C_0^\infty (\mh R^n ,\mh R^n \setminus U) \cong 
C_0^\infty (S^n ,Y)$ for some closed neighborhood $Y$ of $P$, 
and we will show that the latter algebra is diffeotopy 
equivalent to $C_0^\infty (S^n ,P)$. Let $(r_t )_{t \in [0,1]}$ 
be a diffeotopy of smooth maps $S^n \to S^n$ such that
\begin{itemize}
\item $\forall t \;: r_t (P) = P$ and $r_t (Y) \subset Y$ 
\item a neighborhood of $-P$ is fixed pointwise by all $r_t$
\item $r_1 = \mr{id}_{S^n}$ and $r_0 (Y) = P$
\end{itemize}
To construct such maps, we can require that $r_t$ stabilizes 
every geodesic from $-P$ to $P$ and declare that furthermore 
$r_t (Q)$ depends only on $t$ and on the distance from $Q$ to 
$P$. Then we only have to pick a suitable smooth function of 
$t$ and this distance. Given this, consider the Fr\'echet
algebra homomorphisms 
\begin{equation}
\begin{aligned}
& \phi : C_0^\infty (S^n ,P) \to C_0^\infty (S^n ,Y) \\
& \phi (f) = f \circ r_0 \\
& i : C_0^\infty (S^n ,Y) \to C_0^\infty (S^n ,P) \\
& i (f) = f
\end{aligned}
\end{equation}
By construction $\phi \circ i$ and $i \circ \phi$ are 
diffeotopic to the respective identity maps on 
$C_0^\infty (S^n ,Y)$ and $C_0^\infty (S^n ,P)$, so these 
algebras are indeed diffeotopy equivalent. Thus we reduced our 
task to calculating the $K$-groups and periodic cyclic homology 
of $C_0^\infty (S^n ,P)$. 
Fortunately there is an obvious split extension
\begin{equation}
0 \to C_0^\infty (S^n, P) \to C^\infty (S^n) \to \mh C \to 0
\end{equation}
which by Proposition \ref{prop:2.8} consists entirely of 
Fr\'echet algebras in the category $\mc{CIA}$. 
It is well known that
\begin{equation}
K_* (C^\infty (S^n )) \cong K^* (S^n ) \cong 
\check H^* (S^n ; \mh Z) \cong \mh Z^2
\end{equation}
with one copy of $\mh Z$ in degree 0 and the other in degree 
$n$. Since $K_* (\mh C ) = K_0 (\mh C ) \cong \mh Z$ the lemma 
follows from the excision property of the $K$-functor. 
$\quad \Box$
\\[2mm]

Next we prove an equivariant version of the Poincar\'e lemma.
\index{equivariant Poincar\'e lemma}

\begin{lem}\label{lem:2.17}
Let $X, A, G$ and $A^G$ be as in \eqref{eq:2.52}, and 
suppose that $X$ is $G$-equivariantly contractible to a point 
$x_0 \in X$. Then $A^G$ is diffeotopy equivalent to its fiber 
$\mr{End}_G \left( \mh C^N \right)$ over $x_0$. In particular 
$K_* \left( A^G \right) = K_0 \left( A^G \right)$ is a free 
abelian group of finite rank, and $A^G \in \mc{CIA}$.
\end{lem}
\emph{Proof.}
Our main task is to adjust the $u_g$ suitably. Since $X$ is 
contractible we can find for every $g \in G$ a function 
$f_g \in C^\infty (X)$ such that $f_g^{-N} = \det (u_g )$. 
The $G$-action does not change when we replace $u_g$ with 
$f_g u_g$, so we may assume that $\det (u_g ) \equiv 1 ,\, 
\forall g \in G$. The premise that \eqref{eq:2.66} is a group 
action guarantees that there is a smooth 
function $\lambda : G \times G \times X \to \mh C$ such that
\begin{equation}\label{eq:2.67}
u_g(x) u_h( \alpha_g^{-1} x) = \lambda (g,h,x) u_{gh}(x)
\end{equation}
Taking determinants we see that in fact 
$\lambda (g,h,x )^N \equiv 1$, so $\lambda$ does not depend on 
$x \in X$. All the elements of $\alpha (G)$ fix $x_0$, so the 
fiber $V_0 = \mh C^N$ over that point carries a projective 
$G$-representation $\pi_0$. Thus we are in a 
position to apply Schur's theorem \cite{Schu}, which says that 
there exists a finite central extension $G^*$ of $G$ such that 
$\pi_0$ lifts to a representation of $G^*$. This lift only 
involves scalar multiples of the $u_g (x_0) $, so it immediately 
extends to $X$. Then \eqref{eq:2.67} becomes the cocycle relation
\begin{equation}\label{eq:2.68}
u_{gh}(x) = u_g(x) u_h( \alpha_g^{-1} x)
\end{equation}
Notice that still $A^{G^*} = A^G$, so without loss of generality 
we can replace $G$ by $G^*$.

Now we want to make the $u_g (x)$ independent of $x \in X$. 
Wassermann \cite{Was1} indicated how this can be done in the 
continuous case, and his argument can easily be adapted to our 
smooth setting. The crucial observation, first made by 
Rosenberg \cite{Ros}, is that $A^G$ can be rewritten as the 
image of an idempotent in a larger algebra. 
This idempotent can then be deformed to one independent of $x$.

Indeed, let $A \rtimes_\alpha G$ be the crossed product of $A$ 
and $G$ with respect to the action $\alpha$ of $G$ on $X$, and 
$(r_t)_{t \in [0,1]}$ a smooth $G$-equivariant contraction from 
$X$ to $x_0$. (For smooth manifolds the existence of a continuous 
contraction implies the existence of a smooth one.) Define
\begin{equation}
p_t(x) := |G|^{-1} \sum_{g \in G} u_g(r_t x) g
\end{equation}
Then $p_t \in A \rtimes_\alpha G$ is an idempotent by 
\eqref{eq:2.68}, and by Lemma \ref{lem:A.2}
\begin{equation}
\begin{aligned}
& \phi_1 : A^G \to p_1 (A \rtimes_\alpha G) p_1 \\ 
& \phi_1 (\sigma ) = p_1 \sigma p_1
\end{aligned}
\end{equation}
is an isomorphism of Fr\'echet algebras. Clearly the idempotents 
$p_t$ are all homotopic, so they are conjugate in the completion 
$C(X; M_N (\mh C)) \rtimes_\alpha G$ of $A \rtimes_\alpha G$,
which is a Banach algebra if $X$ is compact. Moreover the 
standard argument for this, as 
for example in \cite[Proposition 4.3.2]{Bla}, shows that $p_0$ 
and $p_1$ are conjugate by an element of $A \rtimes_\alpha G$. 
Alternatively we can use the stronger result that homotopic 
idempotents in Fr\'echet algebras are conjugate, but this 
statement is vastly more difficult to prove than its Banach 
algebra version, cf. \cite[Lemmas 1.12 and 1.15]{Phi2}. In any 
case, we have 
\begin{equation}\label{eq:2.69}
A^G \cong p_1 (A \rtimes_\alpha G) p_1 \cong p_0 
(A \rtimes_\alpha G) p_0 
\cong C^\infty (X; \mr{End}_{\mh C } (V_0 ))^G
\end{equation}
To this last algebra we can apply the obvious diffeotopy 
$\sigma \to \sigma \circ r_t$. This shows that $A^G$ is 
diffeotopy equivalent to its fiber $\mr{End}_G (V_0 )$ over 
$x_0$, and the remaining statements on $K_* \left( A^G \right)$ 
and $HP_* \left( A^G \right)$ follow from the semisimplicity of 
the finite dimensional algebra $\mr{End}_G (V_0 ). \quad \Box$
\\[2mm]

Now we can prove the main result of this section, which extends
Lemma \ref{lem:2.17} to general $X$.

\begin{thm}\label{thm:2.18}
Let $X, A, G$ and $A^G = C^\infty (X ;M_N (\mh C ))^G$ be as 
in \eqref{eq:2.52}, and let $\mf L_u^G$ be the sheaf over 
$X/G$ constructed on page \pageref{p:lug}. Then there exists 
a filtration on $K_* \left( A^G \right)$ whose associated 
graded group is isomorphic to 
$\check H^* \left( X/G; \mf L_u^G \right)$, and 
the Chern character induces an isomorphism
\[
K_* \left( A^G \right) \otimes \mh C \isom 
HP_* \left( A^G \right)
\]
Moreover $K_* \left( A^G \right)$ is a finitely generated 
abelian group whenever $X$ is compact.
\end{thm}
\emph{Proof.}
All our arguments will depend on the existence of a specific 
cover of $X$. To construct it we use a theorem of Illman 
\cite{Ill}, which states that $X$ has a smooth \inde{equivariant 
triangulation}. In slightly more down-to-earth language this means 
(among others) that there exists a countable, locally finite 
simplicial complex $\Sigma$ in a finite dimensional orthogonal 
representation space $V$ of $G$, and a $G$-equivariant 
homeomorphism $\psi: \Sigma \to X$. Moreover $\psi$ is smooth as 
a map from a subset of $V$ to $X$, and its restriction to any 
simplex $\sigma$ of $\Sigma$ is an embedding. In particular
$\Sigma$ is a $G$-CW complex, so the assertion on the \v Cech
cohomology of $\mf L_u^G$ follows from Theorems \ref{thm:2.19}
and \ref{thm:2.15}.

For a simplex $\sigma$ we put
\begin{equation}\label{eq:2.70}
U'_\sigma := \{ v \in \Sigma : d(v,\sigma ) \leq r_\sigma \}
\end{equation}
where $d$ is the Euclidean distance in $V$. We require that the
radius $r_\sigma$ depends only on the $G$-orbit of $\sigma$ and 
that $r_\tau > r_\sigma > 0$ if $\tau$ is a face of $\sigma$. 
The orthogonality of the action of $G$ on $V$ guarantees that 
\[
g U'_\sigma = U'_{g \sigma} \quad \mr{and} \quad 
U'_\sigma \cap U'_\tau \subset U'_{\sigma \cap \tau }
\]
if we take our radii small enough. Let $D'_\sigma$ be the 
union, over all faces $\tau$ of $\sigma$, of the $U'_\tau$, 
and $G_\sigma$ the stabilizer of $\sigma$ in $G$. From the 
above we deduce that $U'_\sigma \setminus D'_\sigma$ is 
$G_\sigma$-equivariantly retractible to 
$\sigma \setminus D'_\sigma$.

Now we abbreviate $U_\sigma := \psi (U'_\sigma )$ and 
$D_\sigma := \psi (D'_\sigma )$, so that 
$\{ U_\sigma : \sigma$ simplex of $\Sigma \}$
is a closed $G$-equivariant cover of $X$.
Let $X_m$ be the union of all those $U_\sigma$ for which 
$m + \dim \sigma \leq \dim X$. It is a closed subvariety (with 
boundary and corners) of $X$ and it is stable under the action 
of $G$. Define the following $G$-stable ideals of $A$:
\begin{equation}
I_m := \{ a \in A : a \big|_{X_m} = 0 \}
\end{equation}
By \cite[Th\'eor\`eme IX.4.3]{Tou}
\begin{equation}\label{eq:2.75}
0 \to I_m \to A = C^\infty (X ; M_N (\mh C )) \to 
C^\infty(X_m ; M_N (\mh C )) \to 0
\end{equation}
is an admissible extension of Fr\'echet algebras. 
Using the finiteness of $G$ we see that $I_m^G$ is 
an admissible ideal in $I_{m+1}^G$ and that
\begin{equation}\label{eq:2.71}
I_{m+1}^G / I_m^G \cong (I_{m+1} / I_m)^G \cong 
C_0^\infty (X_m , X_{m+1} ; M_N (\mh C ))^G 
\end{equation}
In order to apply Lemma \ref{lem:2.12} to the sequence
\begin{equation}\label{eq:2.72}
0 = I_0^G \subset I_1^G \subset \cdots \subset I^G_{\dim X}
\subset I^G_{1 + \dim X} =  A^G
\end{equation}
we only have to show that the algebras in \eqref{eq:2.71} are in
the category $\mc{CIA}$. In fact, since 
$\overline{U_\sigma \setminus D_\sigma} \cap \overline{U_\tau 
\setminus D_\tau} = \emptyset$ if $\dim \sigma = \dim \tau$ and 
$\sigma \neq \tau$, we have an isomorphism
\begin{equation}
I_{m+1} / I_m \cong \prod_{m + \dim \sigma = \dim X}
  C_0^\infty ( U_\sigma, D_\sigma ; M_N (\mh C ))
\end{equation}
Now $G$ permutes the simplices in this product, so
\begin{equation}\label{eq:2.73}
I^G_{m+1} / I^G_m \cong \prod_{\sigma \in L_m} C_0^\infty 
(U_\sigma ,D_\sigma ; M_N (\mh C ))^{G_\sigma }
\end{equation}
where $L_m$ is a set of representatives of the simplices of
dimension $\dim X - m$ modulo the action of $G$. Invoking the 
additivity of $K_*$ and $HP_*$ we reduce our task to verifying 
that every factor of \eqref{eq:2.73} belongs to $\mc{CIA}$.

If $m = \dim X$ then $D_\sigma$ is empty and we see from 
Lemma \ref{lem:2.17} that\\ $C^\infty (U_\sigma ; 
M_N (\mh C ))^{G_\sigma}$ has the required property.

For smaller $m$ there also exists (for every $\sigma$) a 
$G_\sigma$-equivariant contraction of $U_\sigma$ to a point 
$x_\sigma \in \psi(\sigma)$. Thus we can follow the proof of 
Lemma \ref{lem:2.17} up to equation \eqref{eq:2.69}, where we 
find that the factor of \eqref{eq:2.73} corresponding to $\sigma$ 
is isomorphic to $C_0^\infty (U_\sigma ,D_\sigma ; 
\mr{End}_{\mh C }(V_\sigma ))^{G_\sigma }$. Here 
$(\pi_\sigma ,V_\sigma )$ denotes the projective 
$G_\sigma$-representation over the point $x_\sigma$. Using the 
$G_\sigma$-equivariant retraction of $U_\sigma \setminus 
D_\sigma$ to $\psi (\sigma \setminus D'_\sigma )$ we see that this 
algebra is diffeotopy equivalent to $C_0^\infty (\sigma ,\sigma 
\cap D'_\sigma) \otimes \mr{End}_{G_\sigma } (V_\sigma )$. The 
right hand side of this tensor product has finite dimension and 
is semisimple, so by the stability of $\mc{CIA}$ it presents no 
problems. Seen from its barycenter $\sigma \setminus D'_\sigma$ 
is star-shaped, hence by Lemma \ref{lem:2.16} the left hand side 
is also in the category $\mc{CIA}$.

We conclude that all the algebras in \eqref{eq:2.71} and 
\eqref{eq:2.73} are indeed objects of $\mc{CIA}$, so Lemma 
\ref{lem:2.12} can be applied to \eqref{eq:2.72} to prove that 
$A^G \in \mc{CIA}$.

Note that the simplicial complex $\Sigma$ has only finitely many
vertices if $X$ is compact, so then all the above direct 
products are in fact finite and $K_* \left( A^G \right)$ is 
finitely generated. $\quad \Box$
\\[2mm]

It is clear from the proofs of Lemma \ref{lem:2.17} and Theorem 
\ref{thm:2.18} that many similar Fr\'echet algebras are also in 
$\mc{CIA}$. For example if $Y$ is a closed submanifold of $X$ 
then the algebra
\begin{equation}
B = \{ f \in C^\infty (X; M_2 (\mh C )) : 
f(y) \;\mr{diagonal}\; \forall y \in Y \}
\end{equation}
is in $\mc{CIA}$, as can be seen from the admissible extension
\begin{equation}
0 \to C_0^\infty (X,Y; M_2 (\mh C )) \to 
B \to C^\infty (Y )^2 \to 0
\end{equation}
One might even study arbitrary Fr\'echet algebras
that are finitely generated over $C^\infty (Y)$, with $Y$
an orbifold. Although it is not unlikely that these are all in
the category $\mc{CIA}$, it seems that a substantial 
generalization of Lemma \ref{lem:2.17} is needed to show this. 

Since the periodic cyclic homology of $A^G$ is finite dimensional
and can be computed in terms of \v Cech cohomology, it is not
surprising that the K\"unneth formula is an isomorphism for such
algebras

\begin{cor}\label{cor:2.29}
Let $(X,A,G,u)$ and $(X',A',G',u')$ both be sets of data like we 
used in \eqref{eq:2.52}. There is a natural isomorphism of
graded vector spaces
\[
HP_* \big( A^G \big) \otimes HP_* \big( A'^{G'} \big) \isom 
HP_* \big( A^G \hot A'^{G'} \big)
\]
\end{cor}
\emph{Proof.}
By Theorem \ref{thm:2.18} we have
\begin{equation}\label{eq:2.83}
\begin{array}{lll}
HP_* \big( A^G \big) & \cong & \check H^* \left( X/G; 
\mf L_u^G \otimes_{\mh Z} \mh C \right) \\
HP_* \big( A'^{G'} \big) & \cong & \check H^* \big( X' / G'; 
\mf L_{u'}^{G'} \otimes_{\mh Z} \mh C \big) \\
HP_* \big( A^G \hot A'^{G'} \big) & \cong & \check H^* \left(
(X \times X') / (G \times G') ; (\mf L_u^G \otimes_{\mh Z}
\mf L_{u'}^{G'})^{G \times G'} \otimes_{\mh Z} \mh C \right) \\
& \cong & \check H^* \left( X/G \times X' / G' ; 
(\mf L_u^G \otimes_{\mh Z} \mh C ) \otimes_{\mh C}
(\mf L_{u'}^{G'} \otimes_{\mh Z} \mh C ) \right)
\end{array}
\end{equation}
According to \cite[\S 6.3]{God} there is a natural map of 
\v Cech complexes
\begin{multline*}
C^* \left( X/G; \mf L_u^G \otimes_{\mh Z} \mh C \right) 
\otimes_{\mh C} C^* \big( X' / G'; \mf L_{u'}^{G'} 
\otimes_{\mh Z} \mh C \big) \quad \longrightarrow \\
\check C^* \left( X/G \times X' / G' ; (\mf L_u^G 
\otimes_{\mh Z} \mh C ) \otimes_{\mh C}
(\mf L_{u'}^{G'} \otimes_{\mh Z} \mh C ) \right)
\end{multline*}
Because all three cohomology groups are finite dimensional
vector spaces the abstract K\"unneth theorem 
\cite[Theorem VI.3.1]{CaEi} tells us that
\[
HP_* \big( A^G \big) \otimes HP_* \big( A'^{G'} \big) 
\cong HP_* \big( A^G \hot A'^{G'} \big)
\]
The construction of the map $\Theta$ in \cite[p. 211]{Emm} is
also possible in the topological setting, and yields a natural
map
\[
\Theta : HP_* \big( A^G \big) \otimes HP_* \big( A'^{G'} \big) 
\to HP_* \big( A^G \hot A'^{G'} \big)
\]
Although the isomorphisms \eqref{eq:2.83} are not natural,
they come from certain filtrations of the underlying topological
spaces, which is enough to ensure that $\Theta$ is also an 
isomorphism. $\qquad \Box$ 
\\[3mm]

Notice also the similarity between $\mf L_u^G$ and the sheaf 
$\mf A$ contructed on page \pageref{eq:2.91}. Suppose that our 
manifold $X$ has the underlying structure of a complex 
nonsingular affine variety \inde{$X_{alg}$}, that the 
$\alpha_g$ are automorphisms of $X_{alg}$ and that the $u_g$ 
are invertibles in \index{$A_{alg}$}
\begin{equation}
A_{alg} := \mc O (X_{alg}) \otimes M_N (\mh C) = 
M_N (\mc O (X_{alg} ))
\end{equation}
Then 
\begin{equation}
g a(x) = u_g (x) a (\alpha_g^{-1} x) u_g^{-1}(x)
\end{equation}
defines an action of $G$ on $A_{alg}$, and $A_{alg}^G$ has the 
same irreducible representations as $A^G$. Applying the recipe 
on page \pageref{eq:2.91} to the finite type algebra $A_{alg}^G$, 
we see that the bases of the stalks $\mf L_u^G (Gx)$ and 
$\mf A (Gx)$ can be identified, and that the requirements for 
continuity of sections on page \pageref{eq:2.91} reduce to those 
on page \pageref{p:lug}. Therefore 
\begin{equation}
\mf A = \mf L_u^G \otimes \mh C
\end{equation}
This insight, together with Theorems \ref{thm:2.6} and 
\ref{thm:2.18}, was the inspiration for an explicit comparison 
theorem between algebraic and topological periodic cyclic 
homology, analogous to \eqref{eq:2.20}.

\begin{thm}\label{thm:2.20}
The inclusion $A_{alg}^G \to A^G$ induces an isomorphism of
finite dimensional vector spaces
\[
HP_* \left( A_{alg}^G \right) \isom HP_* \left( A^G \right)
\]
\end{thm}
\emph{Proof.} 
After noticing that by Theorem \ref{thm:2.6} the left hand side 
has finite dimension, we introduce some notations. 
Let $Y$ be any complex algebraic variety, $Z$ a subvariety and 
$V$ a complex vector space. Like for smooth functions we write 
\index{$\mc O_0 (Y,Z)$} \index{$\mc O_0 (Y,Z;V)$}
\begin{align}
\mc O_0 (Y,Z) &= \{ f \in \mc O (Y) : f \big|_Z = 0 \}\\
\mc O_0 (Y,Z;V) &= \mc O_0 (Y,Z) \otimes V
\end{align}
Start with the finite collection $\mc L$ of all irreducible 
components of the $X_{alg}^g$, as $g$ runs over $G$. Extend this 
to a collection $\{ V_j \}_j$ of subvarieties of $X_{alg}$ by 
including all irreducible components of the intersection of any
subset of $\mc L$. Notice that 
\begin{equation}\label{eq:2.74}
\dim \left( X^g \cap X^h \right) < 
\max \left\{ \dim X^g ,\, \dim X^h \right\}
\end{equation}
if $\alpha_g \neq \alpha_h$. Define $G$-stable subvarieties
\begin{equation}
X_p = \bigcup_{j : \dim V_j \leq p} V_j
\end{equation}
and construct the ideals 
\begin{equation}
\begin{array}{ccccr}
I_p & := & \{ a \in A_{alg}^G \!: a (X_p) = 0 \} & \cong &
  \mc O_0 (X, X_p ;M_N (\mh C ))^G \\
J_p & := & \{ a \in A^G \;: a (X_p) = 0 \} & \cong &
  C_0^\infty (X, X_p ;M_N (\mh C ))^G
\end{array}
\end{equation}
From \eqref{eq:2.75} and \eqref{eq:2.71} we see that all the 
ideals $J_p$ are admissible in $A^G$, and by Theorem 
\ref{thm:2.18} they are in $\mc{CIA}$. By Lemma \ref{lem:2.12} 
it suffices to show that for every $p$ the inclusion 
\begin{equation}
\mc O_0 (X_p ,X_{p-1} ;M_N (\mh C ))^G \cong I_{p-1} / I_p \to
C_0^\infty (X_p ,X_{p-1} ;M_N (\mh C ))^G \cong J_{p-1} / J_p
\end{equation}
induces an isomorphism on periodic cyclic homology. These 
algebras have the same primitive ideal spectrum, namely
$Y_p \setminus Z_p$, where
\begin{equation}
\begin{array}{lllll}
Y_p & = & \mr{Prim} \left( A^G / J_p \right) & = &
  \mr{Prim} \left( A_{alg}^G / I_p \right) \\
Z_p & = & \mr{Prim} \left( A^G / J_{p-1} \right) & = &
  \mr{Prim} \left( A_{alg}^G / I_{p-1} \right) 
\end{array}
\end{equation}
Because all the representations $(\pi_x, \mh C^N)$ are completely
reducible
\begin{equation}
A_{alg}^G = I_0 \supset I_1 \supset \cdots \supset I_N = 0
\end{equation}
is an abelian filtration of $A_{alg}^G$, in the sense of
\cite[Definition 3]{KNS}, and  
\begin{equation}
Z \left( A_{alg}^G / I_p \right) \cap I_{p-1} / I_p =
Z \left( I_{p-1} / I_p \right) \cong \mc O_0 (Y_p ,Z_p)
\end{equation}
This gives alternative descriptions
\begin{align}
Y_p &= \mr{Max} \left( Z \left( A_{alg}^G / I_p \right)\right)\\
Z_p &= \left\{ I \in Y_p : Z \left( I_{p-1} / I_p \right) 
\subset I \right\}
\end{align}
and the proof of \cite[Theorem 10]{KNS} shows that there are 
natural isomorphisms
\begin{equation}\label{eq:2.76}
HP_* \left( I_{p-1} / I_p \right) \cong 
HP_* \left( Z \left( I_{p-1} / I_p \right) \right) \cong
\check H^* (Y_p ,Z_p ;\mh C )
\end{equation}
To get something similar on the topological side we turn to
$K$-theory, knowing already that $J_{p-1} / J_p \in \mc{CIA}$.
Moreover this algebra is dense in $C_0 (X_p ,X_{p-1} 
;\mh C )^G$, so from theorem \ref{thm:2.19} we get 
\begin{equation}
HP_* (J_{p-1} / J_p) \cong K_* \left( C_0 (X_p ,X_{p-1} 
;\mh C )^G \right) \otimes \mh C
\end{equation}
Since $u_g \in A_{alg}$ the type of $(\pi_x ,\mh C^N )$ as a
projective $G_x$-representation cannot change along the 
(connected or irreducible) components of 
$X^{G_x} \cap X_p  \setminus X_{p-1}$. It follows from this, 
\eqref{eq:2.74} and Theorem \ref{thm:2.15} that the inclusions
\begin{equation}
C_0^\infty (Y_p ,Z_p) \to C_0 (Y_p ,Z_p) \cong
Z \left( C_0 (X_p ,X_{p-1} ;\mh C )^G \right) \to
C_0 (X_p ,X_{p-1} ;\mh C )^G
\end{equation}
induce isomorphisms on $K$-theory with rational coefficients.
Moreover 
\begin{equation}
K_* \left( C_0^\infty (Y_p ,Z_p) \right) \otimes \mh C 
\cong HP_* \left( C_0^\infty (Y_p ,Z_p) \right) 
\cong \check H^* (Y_p ,Z_p ;\mh C )
\end{equation}
We put all the above in a diagram
\begin{equation}
\begin{array}{ccccc}
HP_* (I_{p-1} / I_p ) & \leftarrow & HP_* ( \mc O_0 (Y_p ,Z_p) ) & 
\to & \check H^* (Y_p ,Z_p ;\mh C ) \\
\downarrow & & \downarrow & & \parallel \\
HP_* (J_{p-1} / J_p ) & \cong & HP_* \left( C_0^\infty (Y_p ,Z_p) 
\right) & \to & \check H^* (Y_p ,Z_p ;\mh C )
\end{array}
\end{equation}
The horizontal arrows are all natural isomorphisms, so the 
diagram commutes and the vertical arrows are isomorphisms as 
well. $\qquad \Box$

%% file: chapter3.tex
\chapter{Affine Hecke algebras}

Here we commence our study of the main subjects of this 
thesis, affine Hecke algebras. They first appeared in the 
representation theory of certain topological groups, but
that will be discussed only in the next chapter. Instead we
consider Hecke algebras as deformations of the group algebra
of a Weyl group. More precisely, Iwahori-Hecke algebras are
deformations of Coxeter groups, while affine Hecke algebras
are deformations of affine Weyl groups. This deformation
is achieved as follows. Let $s$ be a typical generator of
a Coxeter group $W$. There is a relation 
\[
(T_s - 1)(T_s + 1) = 0
\]
in $\mh Z [W]$. We replace this relation by
\[
(T_s - q(s))(T_s + 1) = 0
\]
where the label $q(s)$ can be any element of a commutative 
ring. In general it is possible to have different labels for 
different generators. Section \ref{sec:3.1} is mainly 
dedicated to making this precise, by providing the definitions 
of root data, label functions and related objects.

The affine Hecke algebra associated with these data will be 
denoted by $\mc H (\mc R ,q)$. If the labels are all positive 
then one can complete this to a $C^*$-algebra $C_r^* (\mc R ,q)$ 
or, more subtly, to a Schwartz algebra $\mc S (\mc R ,q)$.

We have two main goals in this chapter. On one hand we want to
prepare everything for a careful study of deformations in the
parameters $q$, which we will undertake in Chapter 5. This
dictates that we should provide explicit formulas whenever
possible.

On the other hand we would like to apply the ideas developed
in Chapter 2 to affine Hecke algebras. Therefore it is
imperative that we get a clear picture of representations and
the spectrum of $\mc H (\mc R ,q)$. This is provided by the
work of Opdam \cite{Opd3} on the Plancherel measure and the
Fourier transform for affine Hecke algebras. It turns out that
Prim$(\mc H (\mc R, q))$ is a non-separated variety lying over
a complex torus modulo a finite group, see Theorems 
\ref{thm:3.18} and \ref{thm:3.19}. Similarly 
Prim$(\mc S (\mc R ,q))$ is a non-Hausdorff orbifold. On these
spaces $\mc H (\mc R ,q)$ is related to polynomial functions,
$\mc S (\mc R ,q)$ to smooth functions and $C_r^* (\mc R ,q)$
to continuous functions. In fact $\mc S (\mc R ,q)$ is isomorphic 
to an algebra of the type that we studied in Section \ref{sec:2.5} 

However, this is not enough, we also need the Langlands 
classification for $\mc H (\mc R ,q)$. In Section \ref{sec:3.2} 
we explain that it says essentially that Prim$(\mc S (\mc R ,q))$ 
is a deformation retract of Prim$(\mc H (\mc R ,q))$. The final 
form in which we will actually apply this is the rather technical 
parametrization of irreducible $\mc H (\mc R ,q)$-representations 
Theorem \ref{thm:3.24}. With all these preparations, and the
interpretation of periodic cyclic homology as a cohomology theory 
on primitive ideal spectra, we can prove the main theorem of this 
chapter. It says that there are natural isomorphisms
\begin{equation}\label{eq:3.50}
HP_* (\mc H (\mc R ,q)) \cong HP_* (\mc S (\mc R ,q)) \cong
K_* (C_r^* (\mc R ,q)) \otimes \mh C
\end{equation}

\section{Definitions of Hecke algebras}
\label{sec:3.1}

We give precise definitions of (most of) the objects needed to
construct Hecke algebras. We do this both for Iwahori-Hecke
algebras associated to Coxeter groups and for affine Hecke
algebras associated to root data.

A \inde{Coxeter system} $(W,S)$ consists of a finite set $S$
such that $W$ is the group generated by $S$, subject only to
the relations
\[
(s_i s_j)^{m_{ij}} = e
\]
for certain $m_{ij} \in \{ 1,2, \ldots, \infty \}$ such that
\begin{itemize}
\item $m_{ij} = 1$ if and only if $s_i = s_j$
\item $m_{ji} = m_{ij}$
\end{itemize}
Because the most important examples are Weyl groups, we denote 
the \inde{Coxeter group} by $W$ and call the elements of $S$
\inde{simple reflections}. This simple definition still imposes 
a lot of structure, and indeed Coxeter groups have been studied 
deeply. The most relevant results for us can be found for 
example in \cite{Hum}.

A Coxeter system is completely determined by its \inde{Coxeter
graph}. This is a graph whose vertices correspond to elements
of $S$. There is an edge between $s_i$ and $s_j$ if and only
if $m_{ij} \geq 3$, and it is labelled by this number $m_{ij}$.
A Coxeter system is called irreducible if its Coxeter graph is
connected. \index{ll@$\ell$}

Just as for any finitely generated group, there is a natural
\inde{length function} $\ell$ on $W$, which assigns length 1 
to any $s \in S$. To define it, pick $w \in W$ and write it as
\[
w = s_1 \cdots s_r
\]
If $r \geq 0$ is as small as possible, then this is called a
\inde{reduced expression} for $w$ and $\ell (w) = r$.
Notice that $\ell (w^{-1}) = \ell (w)$ since all the simple
reflections have order 2. Depending on the numbers
$m_{ij} ,\, (W, \ell )$ can be finite, of polynomial growth or of
exponential growth.

If $P \subset S$ then $W_P := \langle P \rangle$ is a very special
kind of subgroup of $W$, a standard parabolic subgroup. In general 
a \inde{parabolic subgroup} of $W$ is conjugate to some $W_P$.
The pair $(W_P ,P)$ is a Coxeter system in its own right, and its
length function agrees with the restriction of $\ell$ to $W_P$.
Every right coset $w W_P$ contains a unique element of minimal 
length, so there is a canonical set of representatives $W^P$ for
$W / W_P$. Moreover, if $\{ P_j \}_j$ are the connected components 
of the Coxeter graph of $(W,S)$ then $W = \bigoplus_j W_{P_j}$. 
Hence one can learn a lot about Coxeter systems by studying only 
irreducible ones.

Let $q : W \to \mb k$ be a function from $W$ to a commutative 
unital ring $\mb k$ which is \inde{length-multiplicative}, i.e.
\begin{equation}\label{eq:3.1}
q (wv) = q(w) q(v) \quad \mr{if} \quad
\ell (wv) = \ell (w) + \ell (v)
\end{equation}
This is equivalent to giving a map $q: S \to \mb k$ such that
$q(s_i ) = q(s_j )$ whenever $m_{ij}$ is odd. 

We say that $q$ is an equal label function in the special case
that $q(s_i) = q(s_j )$ for all $s_i , s_j \in S$. In this 
situation it is customary to denote $q(s_i )$ simply by $q$, 
so that 
\[
q(w) = q^{\ell (w)}
\]
The \inde{Iwahori-Hecke algebra} \inde{$\mc H (W,q)$}
$= \mc H_{\mb k} (W,q)$ associated to the Coxeter group $W$ 
($S$ is usually suppressed from the notation) and the 
\inde{label function} \inde{$q$} is an associative 
$\mb k$-algebra which is a free $\mb k$-module with bases 
$\{T_w : w \in W \}$ and fulfills the multiplication rules
\index{$T_w$} \label{p:iha}
\begin{align}
\label{eq:3.2} & T_w T_v = T_{wv} &&
\mr{if} \quad \ell (wv) = \ell (w) + \ell (v) \\
\label{eq:3.3} & T_s T_s = (q(s) - 1) T_s + q(s) T_e &&
\mr{if} \quad s \in S
\end{align}
It is proved in \cite[Section 7.1]{Hum} that such an object
exists and is uniquely determined by these conditions.

Notice that if $q(s) = 1 \: \forall s \in S$, then $\mc H 
(W,q) = \mb k [W]$ is the group algebra of $W$ over $\mb k$.
Any standard parabolic subgroup $V$ gives rise to a 
\inde{parabolic subalgebra} $\mc H \left( V, q \big|_V \right)$, 
which as an $A$-module has bases $\{ T_v : v \in V\}$.

Two choices of $\mb k$ are especially important. The first is 
simply $\mb k = \mh C$. For the second, write \inde{$q_i$} 
$= q(s_i )$ if $s_i \in S$. Let $q_i^{1/2}$ be indeterminates 
satisfying $ \left( q_i^{1/2} \right)^2 = q(s_i )$, and put
\begin{equation}\label{eq:3.85}
\mb k = \mh Z \left[ \big\{ q_i^{1/2}, q_i^{-1/2}
\big\}_{s_i \in S} \right]
\end{equation}
In this case we have
\begin{equation}\label{eq:3.4}
T_s^{-1} = q(s)^{-1} T_s + (q(s)^{-1} - 1) T_e \quad s \in S
\end{equation}
So $T_w$ is invertible in $\mc H (W,q)$ for any $w \in W$.
Indeed, if $w = s_1 \cdots s_r$ is a reduced expression, then
\[
T_w^{-1} = T_{s_r}^{-1} \cdots T_{s_1}^{-1}
\]
To explain the introduction of the square roots we write
\index{$N_w$}
\begin{equation}\label{eq:3.5}
N_w = q(w)^{-1/2} T_w
\end{equation}
These elements form again a bases of $\mc H (W,q)$, while
\eqref{eq:3.2} and \eqref{eq:3.3} become somewhat more
manageable:
\begin{align}
\label{eq:3.6} & N_w N_v = N_{wv} &&
\mr{if} \quad \ell (wv) = \ell (w) + \ell (v)\\
\label{eq:3.7} & \big( N_{s_i} - q_i^{1/2} \big)
\big( N_{s_i} + q_i^{-1/2} \big) = 0 &&
\mr{if} \quad s_i \in S
\end{align}
We generalize the notion of an Iwahori-Hecke algebra as follows.
Let $\Omega$ be a group acting on $W$, and consider the
semidirect product $W \rtimes \Omega$. Assume that
$\forall w \in W, \omega \in \Omega$
\[
\ell (\omega w \omega^{-1}) = \ell (w) \quad \mr{and} \quad
q (\omega w \omega^{-1}) = q (w)
\]
so that we can extend $\ell$ and $q$ to $W \rtimes \Omega$ by
\[
\ell (w \omega) := \ell (w) \quad \mr{and}
\quad q (w \omega) := q(w)
\]
In this situation the rules \eqref{eq:3.2} and \eqref{eq:3.3}
again define a unique associative $\mb k$-algebra with basis
$\{ T_g : g \in W \rtimes \Omega \}$.
Such algebras are called extended Iwahori-Hecke algebras. 
\index{Iwahori-Hecke algebra!extended}
Note that this is extremely general, since $\Omega$ can be 
nearly any group. If we do not want to get too far away from
proper Iwahori-Hecke algebras we have to impose some 
restrictions on this group.

We do this in the setting of an important class of such algebras, 
namely affine Hecke algebras. We will mainly follow the notation 
of \cite{Opd3}, which implies that sometimes we attach a subscript 
0 to a finite object, to distinguish it from its affine counterpart.

First we introduce \inde{root} data, for which we need the 
following objects.
\begin{itemize}
\item $X$ and $Y$ are free abelian groups of the same finite 
rank, and $\inp{\cdot}{\cdot} : X \times Y \to \mh Z$
is a perfect pairing between them
\item $R_0 \subset X$ and $R_0^\vee \subset Y$ are finite
subsets with a given bijection $\alpha \to \alpha^\vee$
\end{itemize}
The elements of $R_0$ are called roots and the elements
of $R_0^\vee$ are called coroots. \index{coroot}
Define endomorphisms \inde{$s_\alpha$} of $X$ and
\inde{$s_\alpha^\vee$} of $Y$ by
\begin{align}
s_\alpha (x) &= x - \inp{x}{\alpha^\vee} \alpha \\
s_\alpha^\vee (y) &= y - \inp{\alpha}{y} \alpha^\vee
\end{align}
For every $\alpha \in R_0$ we impose the conditions
\begin{itemize}
\item $\inp{\alpha}{\alpha^\vee} = 2$
\item $s_\alpha \big( R_0 \big) = R_0$
\item $s_\alpha^\vee \big( R_0^\vee \big) = R_0^\vee$
\end{itemize}
A quadruple \inde{$\mc R$} $= (X,Y,R_0 ,R_0^\vee )$ with 
these properties is a \inde{root datum}. Furthermore it is
\begin{itemize}
\item reduced if $\mh Z \alpha \cap R_0 = 
\{ \alpha, -\alpha \} \; \forall \alpha \in R_0$
\item semisimple if $R_0^\perp = \{ 0 \} \subset Y$
\end{itemize}
Let \inde{$Q$} $:= \mh Z R_0 \subset X$ and
\inde{$Q^\vee$} $:= \mh Z R_0^\vee \subset Y$ be the
\inde{root lattice} and the \inde{coroot lattice}. The
\inde{weight lattice} is $\mr{Hom}_{\mh Z}(Q^\vee, \mh Z)
\supset Q$. If $\mc R$ is semisimple then it contains $X$.

These lattices do not necessarily span \inde{$\mf t^*$} $:= 
X \otimes \mh R$ and its linear dual \inde{$\mf t$} $:= 
Y \otimes \mh R$, but except for that $R_0$ and $R_0^\vee$ are 
\inde{root system}s in the classical sense, and they are dual 
to each other.

Recall that a basis of $R_0$ is a linearly independent
subset \inde{$F_0$} such that every $\alpha \in R_0$ can be
written as
\[
\alpha = \sum_{\beta \in F_0} n_\beta \beta
\]
where either
\[
n_\beta \in \mh Z_{\geq 0} \:\forall \beta \in F_0 \quad \mr{or}
\quad n_\beta \in \mh Z_{\leq 0} \:\forall \beta \in F_0
\]
\index{root!simple} \index{root!negative} \index{root!positive}
This gives a partition $R_0 = R_0^+ \cup R_0^-$. We call the
roots in $F_0$ simple, those in \inde{$R_0^+$} positive and those
in \inde{$R_0^-$} negative.
Bases always exist, and we will assume that one is given with the
root datum, which we will henceforth write as
\begin{equation}
\mc R = (X,Y, R_0 ,R_0^\vee ,F_0 )
\end{equation}
There are several ways to construct new root data from a given
one. Firstly, we can take the direct product of two root data:
\begin{equation}
\mc R \times \mc R' := (X \times X', Y \times Y', R_0 \cup R_0',
R_0^\vee \cup {R'_0}^\vee ,F_0 \cup F'_0 )
\end{equation}
In particular we always have the "trivial one dimensional extension"
\begin{equation}
\mc R \times \mh Z := (X \times \mh Z ,Y \times \mh Z, R_0 ,
R_0^\vee ,F_0 )
\end{equation}
Alternatively one can simply exchange the roles of $X$ and $Y$.
Then \index{$\mc R^\vee$} \index{root datum!dual} \index{$R_P$} 
\begin{equation}
\mc R^\vee = (Y,X, R_0^\vee ,R_0 ,F_0^\vee )
\end{equation}
is the dual root datum of $\mc R$. Furthermore, for $P \subset F_0$
write \index{$R_P^\vee$} \index{$\mc R \times \mc R'$}
\[
R_P = \mh Q P \cap R_0 \quad \mr{and} \quad
R_P^\vee = \mh Q P^\vee \cap R_0^\vee
\]
This $R_P$ is a \inde{parabolic root subsystem} of $R_0$, and 
with it we associate the root datum \index{$\mc R^P$}
\begin{equation}
\mc R^P := \big( X,Y, R_P ,R_P^\vee ,P \big)
\end{equation}
Finally, define \index{$X_P$} \index{$Y_P$} \index{$X^P$}
\index{$Y^P$} \index{$\mc R_P$}
\begin{equation}
\begin{aligned}
&X_P = X \big/ \big( X \cap (P^\vee )^\perp \big) &&
Y_P = Y \cap \mh Q P^\vee \\
&X^P = X / \left( X \cap \mh Q P \right) &&
Y^P = Y \cap P^\perp \\
&\mc R_P = \big( X_P ,Y_P ,R_P ,R_P^\vee ,P \big)
\end{aligned}
\end{equation}
Let us have a look at the various Weyl groups associated to
$\mc R$. Clearly every $s_\alpha$ induces a reflection of
$\mf t$, and they generate a finite Coxeter group \inde{$W_0$}, 
the \inde{Weyl group} of $R_0$. 
In accordance with the terminology for roots, we can take
\begin{equation}
S_0 = \left\{ s_\alpha : \alpha \in F_0 \right\}
\end{equation}
as the set of simple reflections. The action of $W_0$ on $X$ is 
by group homomorphisms, so we can construct the semidirect product 
\index{$W$} \index{$S_0$}
\begin{equation}
W = W (\mc R) := X \rtimes W_0
\end{equation}
By identifying $x \in X$ with the translation \inde{$t_x$}, we
can regard $W$ as a group of affine linear transformations of $X$.
This $W$ is the Weyl group of $\mc R$, and in the same way
we construct its normal subgroup 
\index{$W_{\mr{aff}}$} \index{Weyl group!affine}
\begin{equation}
W_{\mr{aff}} := Q \rtimes W_0
\end{equation}
which we call the affine Weyl group of either $\mc R$ or $R_0$.
Although $W$ may be isomorphic to $\mh Z^n$ (if $R_0 = \es ) 
,\: W_{\mr{aff}}$ is always a Coxeter group, and it 
is possible to extend $S_0$ to a set of Coxeter generators for 
$W_{\mr{aff}}$. To do so, observe that pairing with $F_0$ defines 
a partial ordering on $Y$, and let \inde{$F_m^\vee$} be the set 
of maximal elements of $R_0^\vee$ for this ordering. For every 
$\alpha^\vee \in F_m^\vee$ we consider the affine reflection
\begin{equation}
t_\alpha s_\alpha :
x \to x - \inp{x}{\alpha^\vee} \alpha + \alpha
\end{equation}
in the hyperplane $\{ x \in X : \inp{x}{\alpha^\vee} = 1 \}$.
The desired set of generators is \index{$S_{\mr{aff}}$}
\begin{equation}
S_{\mr{aff}} := S_0 \cup \
\left\{ t_\alpha s_\alpha : \alpha^\vee \in F_m^\vee \right\}
\end{equation}
The \inde{Coxeter graph} of $\mc R$ is that of
$\left( W_{\mr{aff}}, S_{\mr{aff}} \right)$. It is obtained from
the Coxeter graph of $\big( W_0 ,S_0 \big)$ by suitably adding
exactly one vertex to every connected component. Hence a standard
parabolic subgroup $W_P \subset W_{\mr{aff}}$ is infinite whenever 
$|P|>|F_0 |$. If $\left( W_{\mr{aff}}, S_{\mr{aff}} \right)$ is 
irreducible, then any proper parabolic subgroup of $W_{\mr{aff}}$ 
is finite. 

We have a length function $\ell$ on $W_{\mr{aff}}$, which
however does not immediately extend to $W$. To achieve that we
need a different characterization of $\ell$. It is well known that
in the finite Weyl group $W_0$ the length of $w$ is the number of
positive roots made negative by $w$ :
\begin{equation}\label{eq:3.8}
\ell (w) = \# \left\{ \alpha \in R_0^+ : w(\alpha) \in R_0^-
\right\} =  \left| R_0^+ \cap w^{-1} \big( R_0^- \big) \right|
\end{equation}
The same holds in $W_{\mr{aff}}$, one only has to replace $R_0$
by an affine root system. \cite[Propostion 1.23]{IwMa} says

\begin{prop}\label{prop:3.1}
The following formula defines a natural extension of $\ell$ to $W$:
\[
\ell (w t_x) = \sum_{\alpha \in R_0^+ \cap w^{-1} (R_0^- )}
| \inp{x}{\alpha^\vee} + 1 | \;+
\sum_{\alpha \in R_0^+ \cap w^{-1} (R_0^+ )}
| \inp{x}{\alpha^\vee} | \quad\quad w \in W_0, x \in X
\]
Moreover the set \index{Omega@$\Omega$}
\[
\Omega := \{ \omega \in W : \ell (\omega ) = 0 \}
\]
is a subgroup of $W$, complementary to $W_{\mr{aff}}$:
\[
W = W_{\mr{aff}} \rtimes \Omega
\]
\end{prop}

Hence we can say that a reduced expression for $w \in W$ is
a decomposition $w = v \omega$ or $w = \omega v$, where 
$\omega \in \Omega$ and $v \in W_{\mr{aff}}$ is written in 
a reduced way. Define \index{$X^+$} \index{$X^-$}
\begin{equation}
\begin{aligned}
&X^+ := \{ x \in X : \inp{x}{\alpha^\vee} \geq 0
\; \forall \alpha \in F_0 \}\\
&X^- := \{ x \in X : \inp{x}{\alpha^\vee} \leq 0
\; \forall \alpha \in F_0 \} = -X^+
\end{aligned}
\end{equation}
It is well known that 
\begin{equation}\label{eq:3.16}
\begin{array}{ccl}
Z(W) & = & X^+ \cap X^- \\
W_0 X^+ & = & X
\end{array}
\end{equation}
It follows immediately from Proposition \ref{prop:3.1} that
for $x \in X^+$
\begin{equation}\label{eq:3.9}
\ell (t_x) = \sum_{\alpha \in R_0^+} \inp{x}{\alpha^\vee}
\end{equation}
so $\ell$ is additive on the abelian semigroup $X^+ \subset W$.

Our definition of a \inde{label function} is more restrictive
than that for an extended Iwahori-Hecke algebra. It is a 
function $q : W \to \mh C^\times$ such that
$\forall w,v \in W ,\, \omega \in \Omega$
\begin{itemize}
\item $q ( \omega ) = 1$
\item $q (wv) = q(w) q(v) \quad \mr{if} \quad
\ell (wv) = \ell (w) + \ell (v)$
\end{itemize}
The set of label functions is in bijection with the set of 
functions $q : S_{\mr{aff}} \to \mh C^\times$ such that $q (s_i ) 
= q (s_j )$ whenever $s_i$ and $s_j$ are conjugate in $W$.

One can also describe $q$ in terms of a function on a root system
associated to $\mc R$. Assume for simplicity that $R_0$ is reduced,
and define a non-reduced root sytem in $X$ : \index{$R_{nr}$}
\begin{equation}
R_{nr} := R_0 \cup \{ 2 \alpha : \alpha^\vee \in 2 Y \}
\end{equation}
Obviously we write $(2 \alpha)^\vee = \alpha^\vee / 2$, and we
let \inde{$R_1$} be the root system of long roots in $R_{nr}$ :
\begin{equation}
R_1 := \{ \alpha \in R_{nr} : \alpha^\vee \notin 2 Y \}
\end{equation}
The set of simple long roots is
\begin{equation}
F_1 := \{ \alpha \in R_1 : \alpha \in F_0 \;\mr{or}\; 
\alpha /2 \in F_0 \}
\end{equation}
For $\alpha \in R_0 \setminus R_1$ and $\beta \in R_0 \cap R_1$
we put
\begin{align}
\label{eq:3.12} & q_{\alpha^\vee } = q (t_\alpha s_\alpha ) &&
q_{\alpha^\vee / 2} = q(s_\alpha ) q (t_\alpha s_\alpha )^{-1} \\
& q_{\beta^\vee } = q (s_\beta ) = q (t_\beta s_\beta) &&
q_{\beta^\vee / 2} = q_{2 \beta^\vee} = 1 \label{eq:3.13}
\end{align}
Since $\beta^\vee / 2 \notin Y$ and $2 \beta^\vee \notin 
R_{nr}^\vee$, we can ignore them here. However we include 
these conventions because they will simplify some future 
notations. Clearly $q : R_{nr}^\vee \to \mh C^\times$ is 
$W_0$-invariant, and conversely every $W_0$-invariant function 
$R_{nr}^\vee \to \mh C^\times$ determines a unique label 
function $W \to \mh C^\times$ as above.
This correspondence implies the following formulas
\cite[Corollaries 1.3 and 1.5]{Opd2}

\begin{cor}\label{cor:3.2}
For $w \in W_0$ and $x \in X^+$
\[
q (w) = \prod_{\alpha \in R_{nr}^+ \cap w^{-1} ( R_{nr}^- )}
q_{\alpha^\vee} \quad \mr{and} \quad 
q (t_x ) = \prod_{\alpha \in R_{nr}^+}
q_{\alpha^\vee}^{\inp{x}{\alpha^\vee}}
\]
\end{cor}

Given a reduced root datum $\mc R = (X,Y,R_0 ,R_0^\vee ,F_0 )$
and a label function $q$, we define the affine Hecke algebra
$\mc H (\mc R ,q)$ as the unique associative $\mh C$-algebra
which has basis $\{ T_w : w \in W \}$ and satisfies the
multiplication rules \index{Hecke algebra!affine}
\begin{equation}
\begin{aligned}
&T_w T_v = T_{wv} &&
\mr{if} \quad \ell (wv) = \ell (w) + \ell (v) \\
&T_s T_s = (q(s) - 1) T_s + q(s) T_e &&
\mr{if} \quad s \in S_{\mr{aff}}
\end{aligned}
\end{equation}
This algebra is canonically isomorphic to the crossed product
of an Iwahori-Hecke algebra and the group of elements of
length 0 in $W$: 
\begin{equation}
\mc H (\mc R ,q) \cong \mc H (W_{\mr{aff}},q) \rtimes \Omega
\end{equation} 
Our affine Hecke algebra has a large commutative subalgebra
$\mc A$, isomorphic to the group algebra of $X$. We will regard 
this also as $\mc O (T)$, where \index{$T$}
\begin{equation}
T = \mr{Hom}_{\mh Z} \left( X, \mh C^\times \right)
\cong \mr{Prim} \big( \mh C [X] \big) 
\end{equation}
The action of $W_0$ on $X$ induces an action on $T$ by
\begin{equation}
(w \cdot t)(x) = t (w^{-1} x)
\end{equation}
To identify this algebra $\mc A$, let $q^{1/2} : W \to \mh C$ 
be a label function such that $q^{1/2}(w)^2 = q(w) \;
\forall w \in W$. Abbreviate \index{etai@$\eta_i$}
\begin{equation}
q^{1/2}(w)^{-1} = q^{-1/2}(w) \quad \mr{and} \quad 
q^{1/2}(s_i) - q^{-1/2}(s_i) = \eta_i
\end{equation}
In terms of the new bases \index{$N_w$}
\[
\left\{ N_w = q^{-1/2}(w) T_w : w \in W \right\}
\]
the multiplication rules for $\mc H (\mc R ,q)$ become
\begin{equation}\label{eq:3.20}
\begin{aligned}
&N_w N_v = N_{wv} &&
\mr{if} \quad \ell (wv) = \ell (w) + \ell (v) \\
&\big( N_s - q^{1/2}(s) \big) \big( N_s + q^{-1/2}(s) \big) = 0
&& \mr{if} \quad s \in S_{\mr{aff}}
\end{aligned}
\end{equation}
Furthermore every $N_w$ is invertible, and
\begin{equation}
N_{s_i}^{-1} = N_{s_i} - \eta_i N_e 
\quad s_i \in S_{\mr{aff}}
\end{equation}
Because $\ell$ is additive on $X^+$, the map 
\begin{equation}
X^+ \to \mc H (\mc R ,q )^\times : x \to N_{t_x}
\end{equation}
is a monomorphism of semigroups. Extend it to a group homomorphism
$\theta : X \to \mc H (\mc R ,q )^\times$ by defining
\index{thetax@$\theta_x$}
\begin{equation}
\theta_x = N_{t_y} N^{-1}_{t_z} =  N^{-1}_{t_z} N_{t_y}
\end{equation}
if $x = y - z$ with $y,z \in X^+$. This independent of the choice
of $y$ and $z$.

The following theorem is due to Bernstein, Lusztig and Zelevinski,
see \cite[Section 3]{Lus4}. It describes what is also known as the
Bernstein presentation of $\mc H (\mc R ,q)$. 
\index{Amc@$\mc A$} \index{$a^w$}
\begin{thm}\label{thm:3.4}
\begin{enumerate}
\item The sets $\{ N_w \theta_x : w \in W_0 ,x \in X \}$ and\\
$\{ \theta_x N_w : w \in W_0 ,x \in X \}$ are both bases of 
$\mc H (\mc R ,q)$.
\item The subalgebra $\mc A := \mr{span} \{ \theta_x : x \in X \}$ 
is naturally isomorphic to $\mh C [X]$.
\end{enumerate}
\noindent
The Weyl group $W_0$ acts on $\mc A$ by 
$w \cdot \theta_x = \theta_{w x}$, or equivalently
\[
(w \cdot a)(t) = a^w (t) := a (w^{-1} t)  
\qquad t \in T, a \in \mc A \cong \mc O (T)
\]
\begin{enumerate}\setcounter{enumi}{2}
\item The center of $\mc H (\mc R ,q)$ is
\[ 
Z (\mc H (\mc R ,q)) = \mc A^{W_0} \cong 
\mc O (T)^{W_0} \cong \mc O (T / W_0)
\]
\item Take $a \in \mc A ,\, \alpha_i \in F_0$ and let $s_0$ be the
unique vertex of the Coxeter graph of $\mc R$ which is connected
to $s_i = s_{\alpha_i}$ but does not lie in $S_0$. Then
\end{enumerate}
\begin{equation}\label{eq:3.10}
a N_{s_i} - N_{s_i} a^{s_i} =
\left\{ \begin{array}{lll}
\eta_i (a - a^{s_i}) (\theta_0 - \theta_{-\alpha_i})^{-1} &
  \mr{if} & \alpha_i^\vee \notin 2 Y \\
(\eta_i + \eta_0 \, \theta_{-\alpha_i}) (a - a^{s_i})
  (\theta_0 - \theta_{-2\alpha_i})^{-1} & \mr{if} & 
  \alpha_i^\vee \in 2Y
\end{array} \right.
\end{equation}
\end{thm}

Equations \eqref{eq:3.10} are also known as the 
\inde{Bernstein-Lusztig-Zelevinski relations}. 
Since $\mc A$ is of finite rank over $\mc A^{W_0}$, it follows that
$\mc H (\mc R ,q)$ is of finite rank as a module over its center.
In particular it is a finite type algebra in the sense of Section
\ref{sec:2.2}.

For $P \subset F_0$ we can use the above to define label 
functions \inde{$q^P$} and \inde{$q_P$} on the root data 
$\mc R^P$ and $\mc R_P$. 
For $q^P$, use \eqref{eq:3.12} and \eqref{eq:3.13} as a 
definition for all $\alpha ,\beta \in R_P$, and extend this to 
$q^P : W \left( \mc R^P \right) \to \mh C^\times$. Notice that 
$q^P (t_x ) = 1$ whenever $x \perp P^\vee$. Now $q_P : 
W (\mc R_P) \to \mh C^\times$ is simply the map induced by $q^P$.
The affine Hecke algebra $\mc H \left( \mc R^P ,q^P \right)$ 
can be identified with the \inde{parabolic subalgebra} of 
$\mc H (\mc R ,q)$ generated by $\mc A$ and $\mc H \big( W_P 
,q \big|_{W_P} \big)$. Furthermore $\mc H (\mc R_P ,q_P)$ is 
naturally a quotient of $\mc H \left( \mc R^P ,q^P \right)$. 

Notice also that if $\mc R = \mc R_1 \times \mc R_2$, then $q$ 
restricts to label functions on the Weyl groups associated with 
$\mc R_1$ and $\mc R_2$, and there is a canonical identification
\begin{equation}\label{eq:3.46}
\mc H (\mc R ,q) \cong \mc H (\mc R_1 ,q) \otimes 
\mc H (\mc R_2 ,q)
\end{equation}
To study the above algebras we introduce the complex tori
\begin{equation}
\begin{array}{lllll}
T_P & = & \mr{Hom}_{\mh Z} \left( X_P , \mh C^\times \right) & =
& \big\{ t \in T : t(x) = 1 \;\mr{if}\; x \perp P^\vee \big\} \\
T^P & = & \mr{Hom}_{\mh Z} \left( X^P , \mh C^\times \right) & =
& \big\{ t \in T : t(x) = 1 \;\mr{if}\; x \in \mh Q P \big\} \\
\end{array}
\end{equation}
They decompose into a unitary and a real split part:
\index{$T_P$} \index{$T^P$} \index{$T_{P,u}$} 
\index{$T_{P,rs}$} \index{$T_u^P$} \index{$T_{rs}^P$}
\begin{equation}
\begin{array}{llcll}
T_P & = & T_{P,u} \times T_{P,rs} & = & 
\mr{Hom}_{\mh Z} \left( X_P , S^1 \right) \times
\mr{Hom}_{\mh Z} \left( X_P , \mh R^+ \right)\\
T^P & = & T_u^P \times T_{rs}^P & = &
\mr{Hom}_{\mh Z} \left( X^P , S^1 \right) \times
\mr{Hom}_{\mh Z} \left( X^P , \mh R^+ \right)
\end{array}
\end{equation}
Notice that \index{$K_P$}
\[
K_P = T^P \cap T_P = T_u^P \cap T_{P,u}
\]
is a finite group, not necessarily equal to $\{1\}$.
We make the identifications \index{$\mf t^P$}
\begin{equation}
\begin{array}{cclcc}
\mr{Lie} \left( T^P_{rs} \right) & = & \mf t^P & = & 
  Y^P \otimes_{\mh Z} \mh R\\
\mr{Lie} \left( T^P \right) & = & \mf t^P \otimes_{\mh R} 
  \mh C & = & Y^P \otimes_{\mh Z} \mh C\\
\mr{Lie} \left( T_u^P \right) & = & i \mf t^P & = &
  Y^P \otimes_{\mh Z} i \mh R
\end{array}
\end{equation}
Also define the positive parts
\index{$\mf t^{P,+}$} \index{$T_{rs}^{P,+}$}
\begin{equation}
\begin{aligned}
& \mf t^{P,+} \;=\; \left\{ \lambda \in \mf t^P : 
\inp{\alpha}{\lambda} > 0 \; \forall \alpha \in 
F_0 \setminus P \right\}\\
& T^{P,+}_{rs} \;=\; \left\{ t \in T^P_{rs} : t(\alpha ) > 1 \;
  \forall \alpha \in F_0 \setminus P \right\}
  \;=\; \exp \left( \mf t^{P,+} \right)
\end{aligned}
\end{equation}
\\[2mm]

\section{Representation theory}
\label{sec:3.2}

This section is meant as an introduction to the representation
theory of affine Hecke algebras. None of the results presented
here are original, varying from classical (Theorem \ref{thm:3.28})
to very recent (Proposition \ref{prop:3.26}).

Since an affine Hecke algebra is of finite type over its center, 
all its irreducible representations have finite dimension. 
Therefore we mainly study finite 
dimensional representations. We give two partial classifications 
of $\mc H (\mc R ,q)$-representations. One is in terms of their
restrictions to subalgebras associated with finite Coxeter groups. 
The other is more important, and in the spirit of Langlands. It 
shows that the study of $\mc H (\mc R ,q)$-representations can be 
reduced to the study of so-called tempered representations.
After that we define the $C^*$ and Schwartz completions 
$C_r^* (\mc R ,q)$ and $\mc S (\mc R ,q)$ of an affine Hecke 
algebra. It turns out that $\mc S (\mc R ,q)$-representations are
characterized among $\mc H (\mc R ,q)$-representations by the 
requirement that they are tempered.
\\[2mm]

Let \inde{$\mc H$} = $\mc H (\mc R ,q)$ be an affine Hecke 
algebra and \inde{$\mr{Rep}(\mc H )$} its category of finite
dimensional representations. Since $\mc A \cong \mh C [X]$ is 
commutative, every irreducible $\mc A$-module is onedimensional, 
of the form $\mh C_t$ for a character $t \in T$. For $(\pi ,V) 
\in \mr{Rep}(\mc H )$ and $t \in T$ we put \index{$V_t$}
\begin{equation}\label{eq:3.43}
V_t = \left\{ v \in V : \exists n \in \mh N : 
(\pi (a) - a(t))^n v = 0 \; \forall a \in \mc A \right\}
\end{equation}
If $V_t \neq 0$ then there exists an eigenvector $v \in V$ 
with $\pi (a) v = a(t) v \; \forall a \in \mc A$.
In this case $t$ is called an \inde{$\mc A$-weight} of $V$, 
and $V_t$ a generalized weight space. As an $\mc A$-module $V$ 
is the direct sum of the nonzero $V_t$.

Let us consider principal series representations. 
By definition they are the representations
\index{$I_t$} \index{representation!principal series}
\begin{equation}
I_t = \mr{Ind}_{\mc A}^{\mc H} (\mh C_t)
\qquad \mr{for}\; t \in T
\end{equation}
They have dimension $|W_0 |$, and we realize them all on the 
vector space 
\begin{equation}
\mc H (W_0 ,q) \cong \mc H / \langle \{ a - a(t) : 
a \in \mc A \} \rangle
\end{equation}
The importance of principal series representations is already 
clear from the following well-known result.

\begin{lem} \label{lem:3.12}
\begin{enumerate}
\item Every irreducible $\mc H$-representation $(\pi ,V)$ is a 
quotient of some $I_t$
\item If $I_t (h) = 0$ for all $t$ in a Zariski-dense subset
of $T$, then $h = 0$
\item Jac$ (\mc H) = 0$
\end{enumerate}
\end{lem}
\emph{Proof.}
1 was first proven in \cite{Mat2}. Take an $\mc A$-weight $t$ 
of $V$ and a corresponding eigenvector $v \in V$. Define an 
$\mc H$-module homomorphism
\begin{equation}
I_t \to V : h \to \pi (h) v
\end{equation}
This is surjective because $V$ is irreducible. \\
2 is based upon Theorem \ref{thm:3.4}. The function 
\begin{equation}\label{eq:3.44}
T \to \mr{End}_{\mh C}(\mc H (W_0 ,q)) : t \to I_t (h)
\end{equation}
is polynomial, and since it is 0 on a Zariski-dense subset it
vanishes identically. Write $h = \sum_{w \in W_0} a_w T_w$
with $a_w \in \mc A$ and suppose that $h \neq 0$. Then we can 
find $w' \in W_0$ such that $a_{w'} \neq 0$ and $\ell (w')$ is 
maximal with respect to this property. From \eqref{eq:3.10} we
see that 
\[
I_t (h) (T_e) = \sum_{w \in W_0} b_w (t) T_w
\quad \mr{with} \quad b_{w'} = a_{w'}
\]
There \eqref{eq:3.44} is not identically 0, so our assumption 
$h \neq 0$ must be false.\\
3. By \cite[(3.4.5)]{Mat1} or \cite[Theorem 2.2]{Kat1} there is 
a nonempty Zariski-open subset $T'$ of $T$ such that $I_t$ is
irreducible $\forall t \in T'$. So if $h \in \mr{Jac}(\mc H )$
then $h = 0$ by part 2. $\qquad \Box$
\\[2mm]

Similarly we can define the $\mc Z$-weights of $V$, where
\[
\mc Z = Z \big( \mc H (\mc R ,q) \big)
\]
is the center of $\mc H$. Using \ref{thm:3.4}.2 we identify 
Prim$(\mc Z)$ with $T / W_0$. Assume that \inde{$\mc Z$} acts by 
scalars on $V$, which by Schur's lemma is the case if $V$ is 
irreducible. Then the \inde{central character} of $V$ is the 
unique orbit $W_0 t \in T / W_0$ such that
\[
\pi (a) v = a(t) v \quad \forall v \in V, a \in \mc Z
\]
For any $W_0$-stable $U \subset T$ let 
\inde{$\mr{Rep}_U (\mc H (\mc R ,q))$} be the category of all 
finite dimensional $\mc H$-modules whose $\mc Z$-weights are 
contained in $U / W_0$. There are a few ways to "localize" the
algebra $\mc H$ at $U$, i.e. to construct an algebra whose 
representations are precisely $\mr{Rep}_U (\mc H (\mc R ,q))$.
One way, suitable for open $U$, is by tensoring (over $\mc A$) 
with analytic functions on $U$, as we shall see in 
\eqref{eq:3.24}. Another way, which works best if $U$ is a 
closed subvariety of $T$, is by completing $\mc H$ with
respect to the ideal 
\[
J_U = \{ h \in \mc H : I_t (h) = 0 \; \forall t \in U \}
\]
Notice that \inde{$J_U$} is generated by $\{ a \in \mc Z : a 
\big|_U = 0 \}$. Explicitly, we get the modules 
\index{$V_U$} \index{localization} \index{$\mc H_U$}
\begin{equation}\label{eq:3.45}
\begin{aligned}
& V_U = V \otimes_{\mc H} \mc H_U \\
& \mc H_U = \lim_{n \leftarrow \infty} \mc H / J_U^n
\end{aligned}
\end{equation}
For $U = W_0 t$ with $t \in T$ this was done in \cite[Section 
7]{Lus4}, and it is consistent with \eqref{eq:3.43} in the 
sense that
\begin{equation}
V_{W_0 t} = \sum_{w \in W_0} V_{w t}
\end{equation}
Finally, we can vary on this construction in a less subtle fashion,
by replacing $\mc H_U$ in \eqref{eq:3.45} with $\mc H / J_U$. This
has the advantage that we reduce things to modules over finite
dimensional algebras.

Now we start the preparations for the Langlands classification,
which can be found in \cite{DeOp2}.
We say that $V$ is a tempered $\mc H$-module if $|x(t)| \leq 1$
for every $x \in X^+$ and every $\mc A$-weight of $V$. The 
explanation will follow in Lemma \ref{lem:3.9}. Contrarily, we
call $V$ anti-tempered if $|x(t)| \geq 1$ for every such $x$ and 
$t$. \index{representation!anti-tempered}
And less restrictively we say that $V$ is essentially 
tempered if $|x(t)| \leq 1$ for every  $\mc A$-weight $t$ of $V$ 
and every $x \in \mh Z R_0 \cap X^+$. The only difference with 
tempered is that $t \big|_{Z(W)}$ need not be a unitary 
character. \index{representation!tempered}
\index{representation!essentially tempered}
\begin{lem}\label{lem:3.21}
Let $(\pi ,V)$ be an essentially tempered 
$\mc H$-representation which admits a central character
$W_0 r t$, with $t \in T_{rs}^{F_0}$ and $|r| \in T_{F_0,rs}$. 
There exists an automorphism $\psi_t$ of $\mc H$ such that 
$(\pi \circ \psi_t^{-1} ,V)$ is a tempered $\mc H$-representation
with central character $W_0 r$. \index{psit@$\psi_t$}
\end{lem}
\emph{Proof.} Define
\begin{equation}
\psi_t (N_w \theta_x) = t(x) N_w \theta_x
\end{equation}
This is an automorphism because $t$ is 1 on $\mh Z R_0$. Let
$t_1 ,\ldots ,t_d$ be the $\mc A$-weights of $(\pi ,V)$. Clearly,
the $\mc A$-weights of $(\pi \circ \psi_t^{-1} ,V)$ are
$t_1 t^{-1},\ldots ,t_d t^{-1}$. For\\ $x \in \mh Z R_0 \cap X^+$
we have 
\[
|t_i t^{-1} (x)| = |t_i (x)| \leq 1
\]
because $\pi$ is essentially tempered, and for $x \in X^+ 
\cap X^-$ we have $|t_i t^{-1} (x)| = 1$ because 
$|t_i| \in W_0 |r| \subset T_{F_0,rs}$. Hence 
$(\pi \circ \psi_t^{-1} ,V)$ is tempered with central character
$W_0 t_i t^{-1} = W_0 r . \qquad \Box$
\\[3mm]

In general, let $(\sigma ,V_\sigma )$ be any finite 
dimensional representation of\\
\inde{$\mc H^P$} = $\mc H (\mc R^P ,q^P )$. 
Recall \cite[Proposition 1.10]{Hum} that \index{$W^P$}
\begin{equation}
W^P = \{ w \in W_0 : w (P) \subset R_0^+ \}
\end{equation}
is the set of minimal length representatives of $W_0 / W_P$.
Construct the vector space \index{$\mc H (W^P )$}
\begin{equation}
\mc H \big( W^P \big) = \mr{span} \big\{ N_w : w \in W^P \big\} 
\subset \mc H (W_0 ,q)
\end{equation}
The $\mc H$-representation \index{pips@$\pi (P, \sigma)$}
\begin{equation}
\pi (P, \sigma) := \mr{Ind}^{\mc H}_{\mc H^P} (\sigma)
\end{equation}
can be realized on $\mc H \big( W^P \big) \otimes_{\mh C} 
V_\sigma$. From the proof of \cite[Proposition 4.20]{Opd3} 
we can see what the weights of this representation are:

\begin{lem}\label{lem:3.10}
The $\mc A$-weights of $\pi (P,\sigma)$ are precisely the
elements $w (t)$, where $t$ is an $\mc A$-weight 
of $(\sigma ,V_\sigma )$ and $w \in W^P$.
\end{lem}

If $\sigma$ is irreducible then it has a central character
$W_P r \in T / W_P$. Since
\begin{equation}
T_{rs} = T_{rs}^P \times T_{P,rs}
\end{equation}
and $W_P$ acts trivially on $T^P_{rs}$, there is a unique 
\inde{$r_\sigma$} $\in T^P_{rs}$ such that 
\[
W_P |r| = W_P r' r_\sigma \text{ for some } r' \in T_{P,rs}
\]
Let $\Lambda$ be the set of all pairs $(P,\sigma )$, where 
$\sigma$ is an irreducible essentially tempered representation 
of $\mc H^P$. The set of \inde{Langlands data} is 
\index{Lambda@$\Lambda$} \index{Lambda+@$\Lambda^+$}
\begin{equation}
\Lambda^+ = \left\{ (P,\sigma) \in \Lambda : 
r_\sigma \in T_{rs}^{P,+} \right\}
\end{equation}
Now we can state the \inde{Langlands classification} 
for affine Hecke algebras:

\begin{thm}\label{thm:3.11}
Let $(P,\sigma ) \in \Lambda^+$. The $\mc H$-module 
$\pi (P,\sigma )$ is indecomposable and has a unique 
irreducible quotient $L (P,\sigma )$. For every irreducible 
$\mc H$-representation $\pi$ there is precisely one Langlands 
datum $(P,\sigma) \in \Lambda^+$ such that $\pi$ is equivalent 
to $L (P,\sigma )$.
\end{thm}
\emph{Proof.}
This is entirely analogous to the corresponding statements for
graded Hecke algebras, which were proved by Evens \cite[Theorem 
2.1]{Eve}. See also \cite[Section 6]{DeOp2}. $\qquad \Box$
\\[2mm]

Another way to study $\mc H$-representations is by their 
restrictions to simpler subalgebras. We do this by constructing 
a nice (projective) resolution of an $\mc H$-module, which stems 
from joint work of Opdam and Reeder, cf. \cite[Section 8]{Opd4}.

Number the $s_i \in S_{\mr{aff}}$ such that the elements 
corresponding to one connected component of the Coxeter graph
of $(W_{\mr{aff}},S_{\mr{aff}})$ are numbered successively. For 
$I \subset \{ 1,2,\ldots ,| S_{\mr{aff}}| \} $ let $\mh C [I]$
be the vector space with basis $\{ e_i : i \in I \}$ and 
\inde{$W_I$} the standard parabolic subgroup of $W_{\mr{aff}}$ 
generated by $\{ s_i : i \in I \}$. Recall that the "length 0" 
subgroup $\Omega$ of $W$ acts on $S_{\mr{aff}}$ by conjugation. 
Transfer this to an action of $\Omega$ on the indices $i$ and put
\[
\Omega_I := \{ \omega \in \Omega : \omega (I) = I \}
\]
By definition $\Omega_I$ acts on $W_I$, so the extended 
Iwahori-Hecke algebra 
\index{Omegai@$\Omega_I$} \index{$\mc H (\mc R ,I,q)$}
\[
\mc H (\mc R ,I,q) = \mc H (W_I ,q) \rtimes \Omega_I
\]
is well-defined. Note that $W_I$ can be either finite or infinite, 
but that we always have $X^+ \cap X^- = Z(W) \subset \Omega_I$.
If $\mc R$ is semisimple and $W_I$ is finite, then 
$\mc H (\mc R ,I,q)$ has finite dimension.

Because we want to define a conjugate-linear, anti-multiplicative
involution on $\mc H (\mc R ,q)$, from now we will we assume the 
following.
\begin{cond}\label{cond:3.3}
The label function of an affine Hecke algebra only takes
values in $(0,\infty)$
\end{cond}

\begin{lem}\label{lem:3.27}
Suppose that $W_I$ is finite and let $\chi$ be a character of 
$Z(W)$. Then 
\[
\mc H (\mc R ,I,q)_\chi := \mc H (\mc R ,I,q) / \ker \left( 
\mr{Ind}_{\mh C [Z(W)]}^{\mc H (\mc R ,I,q)} \mh C_\chi \right)
\]
is a finite dimensional semisimple algebra.
\end{lem}
\emph{Proof.}
As vector spaces we may identify
\[
\mc H (\mc R ,I,q)_\chi = \mr{Ind}_{\mh C [Z(W)]}^{\mc H (\mc R ,I,q)} 
\mh C_\chi = \mc H (W_I ,q) \otimes \mh C [\Omega_I / Z(W)]
\]
We can extend $|\chi |$ canonically to $X \otimes \mh R$, making
it 1 on $R_0$. Using this extension we define an involution $*_\chi$ 
on $\mc H (\mc R, I,q)$ by
\[
(h_w T_w )^{*_\chi} = \overline{h_w} \: |\chi |(2 w (0)) T_{w^{-1}}
\]
This map is antimultiplicative by Condition \ref{cond:3.3}.
The associated bilinear form is
\[
\inp{h}{h'}_\chi = x_e \quad \mr{if} \quad 
h^{*_\chi} \cdot h' = \sum_{w \in W} x_w N_w
\]
By construction $\mr{Ind}_{\mh C [Z(W)]}^{\mc H (\mc R ,I,q)} 
\mh C_\chi$ is now a unitary representation. This makes\\ 
$\mc H (\mc R ,I,q)_\chi$ into a finite dimensional Hilbert 
algebra, so in particular it is\\ semisimple. $\qquad \Box$
\\[2mm]

For $(\pi ,V) \in \mr{Rep}(\mc H )$ and $n \in \mh N$ consider
the $\mc H$-module \index{$P_n (V)$}
\begin{equation}
P_n (V) = \bigoplus_{|I| = n ,|W_I | < \infty} \mc H 
\otimes_{\mc H (W_I ,q) \rtimes Z(W)} V \big|_{\mc H (W_I ,q) 
\rtimes Z(W)} \otimes_{\mh C} \bigwedge^n \mh C [I]
\end{equation}
Here $\rtimes Z(W)$ is just an abbreviaton of $\otimes \mh C 
[Z(W)]$. We put $P_{|F_0 | +1}(V) = V$ and $P_n (V) = 0$ if 
$n < 0$ or $n \geq |F_0 | + 2$. Define $\mc H$-module homomorphisms
\begin{equation}
\begin{aligned}
& d_n : P_n (V) \to P_{n+1}(V) \\
& d_n (h \otimes_{\mc H (W_I ,q) \rtimes Z(W)} v \otimes \lambda) =
\bigoplus_{j : |W_{I \cup \{j\} }| < \infty} h \otimes_{\mc H 
(W_{I \cup \{j\} },q) \rtimes Z(W) } v \otimes \lambda \wedge e_j
\end{aligned}
\end{equation}
Notice that actually the sum runs only over $j \notin I$, for
otherwise $\lambda \wedge e_j = 0$. To construct a fitting map 
$d_{|F_0 |}$ we need to exert ourselves a little more. 
We introduce a sign function by
\[
\mr{sign}(e_{n_1} \wedge \cdots \wedge e_{n_k}) =
\left\{ \begin{array}{clccr} 
1 & \mr{if} & e_{n_1} \wedge \cdots \wedge e_{n_k} \wedge \eta &
= & e_1 \wedge e_2 \wedge \cdots \wedge e_{|S_{\mr{aff}}|} \\
-1 & \mr{if} & e_{n_1} \wedge \cdots \wedge e_{n_k} \wedge \eta &
= & -e_1 \wedge e_2 \wedge \cdots \wedge e_{|S_{\mr{aff}}|}
\end{array} \right.
\]
where $\eta$ is the wedge product of the $e_i$ with $1 \leq i \leq
|S_{\mr{aff}}|$ and $i \neq n_j$, in standard order. Using this
convention we define
\[
\begin{aligned}
& d_{|F_0 |} : P_{|F_0 |} (V) \to V \\
& d_{|F_0 |} (h \otimes_{\mc H (W_I ,q) \rtimes Z(W)} v \otimes 
\lambda ) = \mr{sign}(\lambda ) \pi (h) v
\end{aligned}
\]
Looking at the $\wedge$-terms we see that $d_{n+1} \circ d_n = 0$,
so $(P_* (V) ,d_* )$ is an augmented differential complex. 
The group $\Omega$ acts naturally on this complex by
\begin{equation}
\omega (h \otimes_{\mc H (W_I ,q) \rtimes Z(W)} v \otimes \lambda) 
= h T_\omega^{-1} \otimes_{\mc H (W_{\omega (I)},q) \rtimes Z(W)} 
\pi (T_\omega ) v \otimes \omega (\lambda )
\end{equation}
This action commutes with the action of $\mc H$ and with the 
differentials $d_n$, so $\left( P_* (V)^\Omega ,d_* \right)$ 
is again a differential complex.

\begin{prop}\label{prop:3.26}
\[
0 \longrightarrow P_0 (V)^\Omega \xrightarrow{\; d_0 \;} 
P_1 (V)^\Omega \xrightarrow{\; d_1 \;} \cdots \longrightarrow 
P_{|F_0 |}(V)^\Omega \xrightarrow{d_{|F_0 |}} V \to 0
\]
is a natural resolution of $V$ by finitely generated modules. 
If $V$ admits the $Z(W)$-character $\chi$ then every 
module $P_n (V)^\Omega$ is projective in the category of all 
$\mc H$-modules with $Z(W)$-character $\chi$.
\end{prop}
\emph{Proof.}
This result a generalization of \cite[Proposition 8.1]{Opd4} 
to root data that are not semisimple. The proof is based upon 
\cite[Section 1]{Kat3}, as indicated by Opdam and Reeder. 

First we consider the case $\Omega = Z(W) = \{e\} ,\, W = 
W_{\mr{aff}}$. There is a linear bijection
\begin{equation}\label{eq:3.42}
\begin{aligned}
& \phi : \mh C [W] \otimes_{\mh C} V \to 
\mc H \otimes_{\mh C} V\\
& \phi (w \otimes v) = T_w \otimes \pi (T_w)^{-1} v
\end{aligned}
\end{equation}
For $s_i \in S_{\mr{aff}}$ we write
\begin{equation}\label{eq:3.40}
\begin{array}{cclcl}
L_i & := & \mr{span}\{ h T_{s_i} \otimes \pi (T_{s_i})^{-1} v 
- h \otimes v : h \in \mc H ,\, v \in V \} 
& \subset & \mc H \otimes_{\mh C} V \\
\mh C [W]_i & := & \big\{ \sum_{w \in W}
x_w w : x_{w s_i} = -x_w \; \forall w \in W \big\}
& \subset & \mh C [W]
\end{array}
\end{equation}
This $L_i$ is interesting because
\begin{equation}
\mc H \otimes_{\mc H (W_I ,q)} V = \big( \mc H \otimes_{\mh C}
V \big) \Big/ \sum_{i \in I} L_i
\end{equation}
Let $w \in W$ be such that $\ell (w s_i) > \ell (w)$. 
\begin{equation}
\begin{aligned}
\phi ((w s_i - w) \otimes v) & = T_{w s_i} \otimes \pi 
(T_{w s_i})^{-1} v - T_w \otimes \pi (T_w )^{-1} v \\
& = T_w T_{s_i} \otimes \pi (T_{s_i})^{-1} \pi (T_w )^{-1} v 
- T_w \otimes \pi (T_w )^{-1} v \in L_i
\end{aligned}
\end{equation}
so $\phi (\mh C [W]_i \otimes V) \subset L_i$. On the other hand,
$L_i$ is spanned by elements as in \eqref{eq:3.40} with 
$h = T_w$ or $h = T_{w s_i}$. 
\begin{equation}
\begin{array}{ll}
\phi^{-1}(T_{w s_i} T_{s_i} \otimes \pi (T_{s_i})^{-1} v -
T_{w s_i} \otimes v) & = \\
\phi^{-1}(q_i T_w + (q_i -1) T_{w s_i} \otimes 
\pi (T_{s_i})^{-1} v) - w s_i \otimes \pi (T_{w s_i}) v & = \\
q_i w \otimes \pi (T_w ) \pi (T_{s_i})^{-1} v + (q_i -1) w s_i
\otimes \pi (T_{w s_i}) \pi (T_{s_i})^{-1} v - 
w s_i \otimes \pi (T_{w s_i}) v & = \\
q_i (w - w s_i ) \otimes \pi (T_w T_{s_i}^{-1})v +
w s_i \otimes \pi \big( q_i T_w T_{s_i}^{-1} + (q_i -1) T_{w s_i}
T_{s_i}^{-1} - T_{w s_i} \big) v & = \\
(w - w s_i ) \otimes \pi (T_w q_i T_{s_i}^{-1})v + 
w s_i \otimes \pi \big( T_w (T_{s_i} + 1 - q_i ) + (q_i -1) T_w 
- T_w T_{s_i} \big) v \!\!\!& = \\
(w - w s_i ) \otimes \pi \big( T_w (T_{s_i}+1 - q_i ) \big) v 
\qquad \in \qquad \mh C [W]_i \otimes V
\end{array}
\end{equation}
We conclude that $\phi^{-1}(L_i ) = \mh C [W]_i \otimes V$.
Now we bring the linear bijections
\begin{equation}
\mh C [W] \Big/ \sum_{i \in I} \mh C [W]_i \to \mh C [W / W_I] :
w \to w W_I
\end{equation}
into play. Under these identifications our differential complex
becomes
\[
0 \to \mh C [W] \otimes V \otimes \mh C \to \cdots \to 
\bigoplus_{|I| = n , |W_I | < \infty} \!\!\!\! \mh C [W / W_I ] 
\otimes V \otimes \bigwedge^n \mh C [I] \to \cdots \to V \to 0
\]
We have to show that the cohomology of this complex vanishes in
all degrees. Because $d_{|F_0 |}$ is surjective it suffices to
show that the complex 
\begin{equation}
\begin{array}{lll}
C_n' & = & \bigoplus\limits_{|I| = n ,|W_I | < \infty} \mh C 
[W/W_I ] \otimes \bigwedge^n \mh C [I] \\
d_n' (w W_I \otimes \lambda ) & = & \sum\limits_{j : |W_{I \cup 
\{j\} }| < \infty} w W_{I \cup \{j\} } \otimes \lambda \wedge e_j
\end{array}
\end{equation}
has cohomology
\begin{equation}\label{eq:3.41}
H'_n = \left\{ \begin{array}{lll}
\mh C & \mr{if} & n = |F_0 | \\
0 & \mr{if} & n \neq |F_0 |
\end{array} \right.
\end{equation}
This is best seen by a geometrical interpretation. It is 
well-known that the alcove 
\[
C_\es := \big\{ x \in Q \otimes_{\mh Z} \mh R : 
\inp{x}{\alpha^\vee} \geq 0 \; \forall \alpha^\vee \in F_0^\vee ,\, 
\inp{x}{\beta^\vee} \leq 1 \; \forall \beta^\vee \in F_m^\vee \big\}
\]
is a fundamental domain for the action of $W$ on 
$Q \otimes_{\mh Z} \mh R$. The finite groups $W_I$ are naturally 
identified with the stabilizers of the faces $C_I$ of 
$C_\es$. Thus every coset $w W_I$ corresponds to a 
polysimplex $w C_I$ of dimension $|F_0 | - |I|$. It follows that 
$(C'_* ,d'_* )$ is the complex that computes $H_{|F_0 | - *}(Q 
\otimes_{\mh Z} \mh R ; \mh C)$ by means of a (polysimplicial) 
triangulation. Together with the Poincar\'e Lemma this leads to 
\eqref{eq:3.41}, proving the proposition in the special case 
$\Omega = \{e\}$.

Now the general case. Clearly $P_n (V)$ is finitely generated because 
$V$ has finite dimension and because there are only finitely many 
$I$'s involved. 

Since the action of $\Omega$ on $S_{\mr{aff}}$ 
factors through a finite group we can construct a Reynolds operator 
\[
R_\Omega := [\Omega : Z(W) ]^{-1} \sum_{\omega \in \Omega / Z(W)} 
\omega \quad \in \quad \mr{End}_{\mc H} \big( P_n (V) \big)
\]
It follows that $P_n (V )^\Omega = R_\Omega \cdot P_n (V)$ is a
direct summand of $P_n (V)$, and hence also finitely generated.

We generalize \eqref{eq:3.42} to a bijection
\begin{equation}
\begin{aligned}
& \phi : \mh C [W / Z(W)] \otimes_{\mh C} V \to 
\mc H \otimes_{\mh C [Z(W)]} V \\
& \phi (w \otimes v) = T_w \otimes \pi (T_w)^{-1} v
\end{aligned}
\end{equation}
Just as above this leads to bijections
\begin{equation}
\bigoplus\limits_{|I| = n ,|W_I | < \infty} \mh C [W/(W_I 
\times Z(W))] \otimes V \otimes \bigwedge^n \mh C [I] \to P_n (V)
\end{equation}
Since both sides are free $\Omega / Z(W)$-modules we see that
\begin{equation}
\bigoplus\limits_{|I| = n ,|W_I | < \infty} \mh C 
[W_{\mr{aff}}/W_I ] \otimes V \otimes \bigwedge^n \mh C [I] 
\cong P_n (V)^\Omega
\end{equation}
Now the above geometrical argument shows that the $P_n (V)^\Omega$
do indeed form a resolution of $V$.

Assume now that $V$ admits a $Z(W)$-character $\chi$. 
Then $\pi \big|_{\mc H (\mc R ,I,q)}$ factors through 
$\mc H (\mc R ,I,q )_\chi$ and by Lemma \ref{lem:3.27} it is 
projective as a representation of this algebra. 
This implies that $P_n (V)$ and its direct summand 
$P_n (V)^\Omega$ are projective in the category of 
$\mc H$-representations with $Z(W)$-character $\chi. \qquad \Box$

\begin{cor}\label{cor:3.29}
Suppose that $V,V' \in \mr{Rep}(\mc H )$ are such that
$V \big|_{\mc H (\mc R ,I,q)}$ and\\ 
$V' \big|_{\mc H (\mc R ,I,q)}$ are equivalent whenever $W_I$ 
is finite. Then $V$ and $V'$ define the same class in the 
Grothendieck group of finitely generated $\mc H$-modules.
\end{cor}
\emph{Proof.}
By assumption we can find a collection of $\mc H (\mc R 
,I,q)$-module isomorphisms $\alpha_I : V \to V'$ such that
\[
\pi' (T_\omega ) \alpha_I = \alpha_{\omega (I)} \pi (T_\omega )
\qquad \forall \omega \in \Omega
\]
These combine to $\mc H$-module isomorphisms
\begin{align*}
& \alpha_n : P_n (V )^\Omega \to P_n (V' )^\Omega \\
& \alpha_n (h \otimes_{\mc H (W_I ,q) \rtimes Z(W)} v 
\otimes \lambda) = h \otimes_{\mc H (W_I ,q) \rtimes Z(W)} 
\alpha_I (v) \otimes \lambda
\end{align*}
Hence by Proposition \ref{prop:3.26} $V$ and $V'$ have 
equivalent finitely generated resolutions. $\qquad \Box$
\\[2mm]

Unfortunately this result is not very strong, as the class of
$V$ is nearly always 0. This certainly is the case if $V$ is
of the form $\pi (P,\sigma )$ with $P \neq F_0$.

Let \inde{$G(\mc H )$} be the Grothendieck group of finite 
dimensional $\mc H$-modules, and \inde{$K_0 (\mc H )$} the
Grothendieck group of finitely generated projective 
$\mc H$-modules. The \inde{Euler-Poincar\'e pairing} 
\cite[Section III.4]{ScSt} on $G (\mc H )$ is defined as
\begin{equation}\label{eq:3.55}
EP \,(V,V') = \sum_{n=0}^\infty (-1)^n \dim \mr{Ext}^n_{\mc H}(V,V')
\end{equation}
where $\mr{Ext}^n_{\mc H}$ is the higher derived functor of 
$\mr{Hom}_{\mc H}$. \index{$EP$}

Let us recall some relevant observations from \cite[Section 
8]{Opd4}. Assume that $\mc R$ is semisimple. By Proposition 
\ref{prop:3.26} the \inde{Euler characteristic} can be defined
for finite dimensional $\mc H$-modules by \index{Eul}
\begin{equation}\label{eq:3.56}
\begin{aligned}
& \mr{Eul} : G (\mc H ) \to K_0 (\mc H ) \\
& \mr{Eul}\, [V] = \sum_{n=0}^{|F_0 |} (-1)^{|F_0 |-n} 
\left[ P_n (V)^\Omega \right]
\end{aligned}
\end{equation}
Moreover there is a natural pairing between $G(\mc H)$ and
$K_0 (\mc H )$. Given a representation $\pi$ of $\mc H$ and an
idempotent $p \in M_n (\mh C ) \otimes \mc H$ we put
\begin{equation}
\big[\, [p] , [\pi ] \,\big] = 
\mr{rank \,(id}\otimes \pi)(p) \in \mh Z
\end{equation}
Because $(P_* (V)^\Omega ,d_*)$ is a projective resolution of
$V$ we have the equalities
\begin{equation}\label{eq:3.58}
EP \,(V,V') = \sum_{n=0}^{|F_0 |} (-1)^{|F_0 |-n} 
\dim \mr{Hom}_{\mc H} \left( P_n (V)^\Omega , V' \right) = 
\big[\, \mr{Eul}\, [V], [V'] \,\big] 
\end{equation}
But the modules $P_n (V)$ are induced from semisimple subalgebras
of finite dimension, so this can be expressed more explicitly. 
By Frobenius reciprocity we have
\begin{equation}\label{eq:3.59}
\dim \mr{Hom}_{\mc H} (P_n (V),V') = \sum_{|I| = n ,|W_I | < \infty} 
\dim \mr{Hom}_{\mc H (W_I ,q) \rtimes Z(W)} (V,V')
\end{equation}
Let $\ep_I$ be the character of the $\Omega_I$-representation
$\bigwedge^{|I|} \mh C [I]$. Taking the $\Omega$-invariants of
$P_n (V)$ in \eqref{eq:3.59} we find
\begin{equation}\label{eq:3.57}
EP \,(V,V') \;= \sum_{I \subset S_{\mr{aff}} ,|W_I | < \infty} 
\frac{(-1)^{|F_0| - |I|}}{[\Omega : \Omega_I ]} 
\dim \mr{Hom}_{\mc H (\mc R ,I ,q)} (V \otimes \ep_I , V' )
\end{equation}
Thus the Euler-Poincar\'e pairing of two finite dimensional 
$\mc H$-modules is completely determined by their restrictions to 
the finite dimensional semisimple subalgebras $\mc H (\mc R ,I 
,q)$. Because such algebras have only finitely many inequivalent 
irreducible modules, this pairing is symmetric and invariant under 
continuous deformations of its input. However, from \eqref{eq:3.55} 
we quickly deduce that modules with different central characters 
are orthogonal for $EP$. We will see in \eqref{eq:3.54} that 
every module of the form $\pi (P,\sigma )$ with $P \subsetneq F_0$
admits continuous deformations of its central character. 
Therefore all such modules are in the radical of $EP$.
\\[2mm]

So far for the purely algebraic properties of affine Hecke 
algebras, we are also interested in their analytic structure.
Condition \ref{cond:3.3} gives us a canonical way to define the
label function $q^{1/2}$ and the basis elements $N_w$. Let 
$x = \sum_{w \in W} x_w N_w$ be an element of $\mc H (\mc R ,q)$ 
and define its adjoint and its trace by
\index{$x^*$} \index{tau@$\tau$}
\begin{equation}\label{eq:3.15}
x^* := \sum_{w \in W} \overline{x_w} N_{w^{-1}} \quad 
\mr{and} \quad \tau (x) = x_e 
\end{equation}
Condition \eqref{cond:3.3} assures that indeed $\tau$ is positive
and that $(x y)^* = y^* x^*$. This leads to a \inde{bitrace} or 
Hermitian inner product
\begin{equation}\label{eq:3.inp}
\inp{x}{y} := \tau (x^* y) \quad x,y \in \mc H (\mc R ,q)
\end{equation}
and a norm \index{$\norm{\:\cdot\:}_\tau$}
\begin{equation}
\norm{x}_\tau := \sqrt{\inp{x}{x}} = \sqrt{\tau (x^* x)}
\end{equation}
By a simple calculation one can show that $\{ N_w : w \in W \}$ is 
an orthonormal bases of $\mc H (\mc R ,q)$ for this inner product. 
The bitrace $\inp{\cdot}{\cdot}$ gives $\mc H (\mc R ,q)$ the 
structure of a Hilbert algebra, in the sense of \cite[Appendice 
A 54]{Dix}. Let \inde{$\mf H (\mc R, q)$} be its Hilbert space 
completion, and $B (\mf H (\mc R ,q))$ the associated $C^*$-algebra 
of bounded operators. Consider the multiplication maps 
\index{lambda@$\lambda$} \index{rho@$\rho$}
\begin{align*}
&\lambda (x) : y \to x y \\
&\rho    (x) : y \to y x 
\end{align*}
According to \cite[Lemma 2.3]{Opd3} these maps extend to bounded
operators on $\mf H (\mc R ,q)$ of the same norm, which we denote 
by \index{$\norm{\:\cdot\:}_o$}
\begin{equation}
\norm{x}_o := \norm{\lambda (x)}_{B (\mf H (\mc R ,q))} =
\norm{\rho (x)}_{B (\mf H (\mc R ,q))}
\end{equation}
Thus, $\mc H (\mc R ,q)$ being a *-subalgebra of 
$B (\mf H (\mc R ,q))$, we can take its closure 
\inde{$C_r^* (\mc R ,q)$} in the norm topology. By definition 
this is a separable unital $C^*$-algebra, and we call it the 
\inde{reduced $C^*$-algebra} of $\mc H (\mc R ,q)$. 
This is analogous to the reduced $C^*$-algebra of a locally 
compact group $G$, which is the completion of $C_0 (G)$ for the 
topology coming from the regular representation of $G$.

There is also a "smoother" way to complete an affine Hecke 
algebra. As a topological vector space, it will consist of 
rapidly decreasing functions on $W$ with respect to some length
function. For this purpose it is a little unsatisfactory that 
$\ell$ is 0 on the subgroup $Z(W) = X^+ \cap X^-$, since this 
can be the whole of $W$. To overcome this inconvenience, let 
$f : X \otimes \mh R \to [0, \infty )$ be a function such that
\begin{itemize}
\item $f(X) \subset \mh Z$
\item $f(x+q) = f(x) \quad \forall x \in X \otimes \mh R ,\,
q \in Q \otimes \mh R$
\item $f$ induces a norm on
$X \otimes \mh R \big/ Q \otimes \mh R \cong
Z(W) \otimes \mh R$
\end{itemize}
Now we define for $w \in W$ \index{$\mc N$}
\begin{equation}
\mc N (w) := \ell (w) + f (w(0))
\end{equation}
so that for any $w,v \in W ,\, u \in W_{\mr{aff}}, 
\omega \in \Omega$
\begin{align}
&\mc N (u\omega ) = \mc N (\omega u) = \ell (u)+ \mc N (\omega )\\
&\mc N (v w) \leq \mc N (v) + \mc N (w)
\end{align}
Because $Z(W) + Q$ is of finite index in $X$, the set $\Omega' = 
\{ w \in W : \mc N (w) = 0 \}$ is finite. Moreover, since $W$ is 
the semidirect product of a finite group and an abelian group, 
it is of polynomial growth, and different choices of $f$ lead to 
equivalent length functions $\mc N$. Now we can define for any 
$n \in \mh N$ the norm 
\index{$p_n$} \index{Omega'@$\Omega'$}
\begin{equation}\label{eq:3.48}
p_n \Big( \sum_{w \in W} x_w N_w \Big) := 
\sup_{w \in W} |x_w| ( \mc N (w) + 1 )^n
\end{equation}
The completion \inde{$\mc S (\mc R ,q)$} of $\mc H (\mc R ,q)$ 
with respect to the family of norms $\{ p_n \}_{n \in \mh N}$ is a 
nuclear Fr\'echet space. It consists of all possible infinite sums 
$x = \sum_{w \in W} x_w N_w$ such that $p_n (x) < \infty \; 
\forall n \in \mh N$.
Opdam \cite[Section 6.2]{Opd3} proved that these norms behave
reasonably with respect to multiplication:

\begin{thm}\label{thm:3.5}
There exist $C_q > 0 ,\, d \in \mh N$ such that 
$\forall x,y \in \mc S (\mc R ,q), n \in \mh N$
\begin{align}
\norm{x}_o &\leq C_q p_d (x) \\
p_n (xy) &\leq C_q p_{n+d}(x) p_{n+d}(y)
\end{align}
In particular $\mc S (\mc R ,q)$ is a unital locally convex 
*-algebra, and it is contained in $C_r^* (\mc R ,q)$.
\end{thm}
The proof of this theorem uses heavy machinery, namely the spectral
decomposition of the trace $\tau$, which we will state in Theorem
\ref{thm:3.18}.4. Closer examination of that proof shows that we
can take $d = \mr{rk}(X) + |W_0|^2 + 1$. In Section \ref{sec:5.2}
we will use more elementary tools to prove a generalization of 
this theorem, resulting in a smaller $d$.

We call $\mc S (\mc R ,q)$ the \inde{Schwartz algebra} of
$\mc H (\mc R ,q)$. In Section \ref{sec:4.2} we will see that this 
construction is analogous to the Schwartz algebra of a reductive 
$p$-adic group.
Considering Schwartz completions of parabolic subalgebras, we 
see that $\mc S (\mc R_P ,q_P )$ is still a quotient of 
$\mc S (\mc R^P ,q^P)$, but that the latter algebra is in general
no longer contained in $\mc S (\mc R ,q)$. The same holds for the 
respective reduced $C^*$-algebras.

From the work of Casselman \cite[\S 4.4]{Cas} one can 
find necessary and sufficient conditions under which 
$\mc H$-representations extend continuously to certain 
completions. Indeed it is shown in \cite[Lemma 2.20]{Opd3} that 
\index{representation!discrete series}
\begin{lem}\label{lem:3.8}
For an irreducible $\mc H$-representation $(\pi ,V)$ the following 
conditions are equivalent, and summarized by saying that $\pi$
belongs to the discrete series:
\begin{itemize}
\item $(\pi ,V)$ is a subrepresentation of the left regular 
representation $(\lambda ,\mf H (\mc R ,q))$
\item all matrix coefficients of $(\pi ,V)$ are in 
$\mf H (\mc R ,q)$
\item $|x(t)| < 1$ for every $\mc A$-weight $t$ of $V$ and 
every $x \in X^+ \setminus 0$
\end{itemize}
\end{lem}

By definition a discrete series representation is unitary, and it
extends to the reduced $C^*$-algebra $C_r^* (\mc R ,q)$. Because 
this is a Hilbert algebra, a suitable version of 
\cite[Proposition 18.4.2]{Dix} shows that $\pi$ is an isolated
point in its spectrum Prim$\big( C_r^* (\mc R ,q) \big)$. Moreover,
since $C_r^* (\mc R ,q)$ is unital, its spectrum is compact
\cite[Proposition 3.18]{Dix}, so there can be only finitely many
inequivalent discrete series representations. On the other hand,
they can only exist if $X^+ \cap X^- = 0$, i.e. if the root 
datum $\mc R$ is semisimple.

Let us mention two important examples. If $W = W_{\mr{aff}}$ and  
$q(s) > 1 \; \forall s \in S_{\mr{aff}}$ then 
\[
\pi_{St} : N_w \to (-1)^{\ell (w)} q(w)^{-1/2}
\]
defines a discrete series representation. This is called the 
Steinberg representation of $\mc H (\mc R ,q)$. Contrarily, if 
$q(s) < 1  \; \forall s \in S_{\mr{aff}}$ then the representation 
\index{representation!Steinberg} \index{representation!trivial}
\[
\pi_{triv} : N_w \to q(w)^{1/2}
\]
is discrete series. This is known as the trivial representation
of $\mc H (\mc R ,q)$, because it is a deformation of the trivial
representation of $W$.

A linear functional $f : \mc H \to \mh C$ is tempered if there
exist $C,N \in (0,\infty)$ such that for all $w \in W$
\[
|f (N_w)| \leq C (1 + \mc N (w) )^N
\]
The collection of all \inde{tempered functional}s is the linear 
dual of $\mc S (\mc R ,q)$. From \cite[Lemma 2.20]{Opd3} we get
the following characterization:
\index{representation!tempered}
\begin{lem}\label{lem:3.9}
For $(\pi, V) \in \mr{Rep}(\mc H )$ the following are equivalent:
\begin{itemize}
\item $\pi$ extends continuously to $\mc S (\mc R ,q)$
\item every matrix coefficient of $(\pi, V)$ 
is a tempered functional
\item $V$ is a tempered $\mc H$-representation
\end{itemize}
\end{lem}

Suppose that our root datum is a product $\mc R = \mc R_1 \times 
\mc R_2$. It is clear from \eqref{eq:3.48} that the Schwartz 
completion respects the decomposition \eqref{eq:3.46}, so
\begin{equation}\label{eq:3.47}
\mc S (\mc R ,q) \cong \mc S (\mc R_1 ,q) \,\hot\,
\mc S (\mc R_2 ,q)
\end{equation}
An isomorphism like
\[
C_r^* (\mc R ,q) \cong C_r^* (\mc R_1 ,q) \,\otimes_t \,
C_r^* (\mc R_2 ,q)
\]
should also hold, but here the particular topological tensor
product $\otimes_t$ is probably dictated by the very isomorphism.

Let \inde{$q^0$} be the label function that is identically 1 and
\inde{$T_u$} the unitary or real compact part of $T$. For any 
$\mc R$ the above constructions reduce to the well-known algebras
\begin{equation}\label{eq:3.51}
\begin{array}{lllllll}
\mc H (\mc R ,q^0) & = & \mh C [W] & = & \mh C [X] \rtimes W_0 & 
\cong & \mc O (T) \rtimes W_0\\
\mc S (\mc R ,q^0) & = & \mc S (W) & \cong & 
\mc S (X) \rtimes W_0  & \cong & C^\infty (T_u ) \rtimes W_0\\
C_r^* (\mc R ,q^0) & = & C_r^* (W) & \cong & 
C_r^* (X) \rtimes W_0  & \cong & C (T_u ) \rtimes W_0
\end{array}
\end{equation}
where $\mc S (X)$ denotes the space of Schwartz functions on $X$.
The representation theory of these algebras is not very 
complicated, and can be obtained from classical results which
go back to Frobenius and Clifford \cite{Cli}.

\begin{thm}\label{thm:3.28}
\begin{enumerate}
\item The $W$-module $I_t = \mr{Ind}_X^W (\mh C_t )$ is completely
reducible for any $t \in T$.
\item $I_t$ is unitary if and only if $t \in T_u$, in which case
it is also tempered and anti-tempered.
\item The number of inequivalent irreducible summands of $I_t$
is exactly the number of conjugacy classes in the isotropy group
$W_{0,t}$.
\end{enumerate}
\end{thm}
\emph{Proof.}
1. As an $X$-representation 
\[
I_t \cong \bigoplus_{w \in W_0} \mh C_{wt}
\]
and the $T_w$ act on $I_t$ by permuting these onedimensional 
subspaces. Hence there is a natural bijection between 
subrepresentations of the left regular representation of $W_{0,t}$
and subrepresentations of $I_t$. It is given explicitly by $V \to 
\mr{Ind}_{W_{0,t}}^{W_0} V$, with $X$ acting on the subspace
$T_w V$ by the character $w t$. Since $\mh C [W_{0,t}]$ is
completely reducible, so is $I_t$. \\
2 will be a special case of Proposition \ref{prop:3.13}.\\
3 follows from counting the irreducible representations of the
finite group $W_{0,t}. \qquad \Box$
\\[4mm]

\section{The Fourier transform}
\label{sec:3.3}

Now we really delve into the representation theory of affine
Hecke algebras. This is a very complicated subject and we will
barely prove anything. In fact this section is more like an 
overview of some of the work of Delorme and Opdam \cite{DeOp1, 
DeOp2, Opd3}. In principle we mention only those things that 
we use later in some way or another, but this turns out to 
be a lot. The highlights are undoubtedly the concise Theorems
\ref{thm:3.18} and \ref{thm:3.19}, which describe the Fourier 
transform of affine Hecke algebras and their completions.
Most of the other things are technicalities that can be skipped 
on first reading.
\\[2mm]

For any $t \in T^P$ there is a surjective algebra
homomorphism \index{phit@$\phi_t$}
\[
\phi_t : \mc H^P \to \mc H_P
\]
which is the identity on $\mc H (W_0 ,q)$ and sends $\theta_x$ 
to $t(x) \theta_{x_P}$, where $x_P$ is the image of $x$ in 
$X_P = X / \big( X \cap (P^\vee)^\perp \big)$. 
If $\sigma \in \mr{Rep}(\mc H_P )$ then we can construct the 
$\mc H$-representation \index{pist@$\pi (P, \sigma, t)$}
\begin{equation}\label{eq:3.54}
\pi (P, \sigma, t) = 
\mr{Ind}_{\mc H^P}^{\mc H} (\sigma \circ \phi_t )
\end{equation}
Because $\mc A$ is in general not a *-subalgebra of $\mc H$, 
it is not immediately clear whether this procedure preserves
unitarity or temperedness of representations. Recall from
\cite[Propostion 1.12]{Opd2} that

\begin{lem}\label{lem:3.6}
Let \inde{$w_0$} be the longest element of $W_0$. 
For $x \in X$ we have
\[
\theta_x^* = T_{w_0} \theta_{- w_0 x} T_{w_0}^{-1} =
N_{w_0} \theta_{- w_0 x} N_{w_0}^{-1}
\]
\end{lem}

In particular for $x \in X \cap \big( P^\vee \big)^\perp$ we have 
$\theta_x^* = \theta_{-x}$ in $\mc H^P$, so $\phi_t$ preserves
the * if and only if $t \in T_u^P$. Similarly the inclusion 
$\mc H^P \to \mc H$ is in general not *-preserving. Nevertheless

\begin{prop}\label{prop:3.13}
Let $\sigma$ be an irreducible $\mc H_P$-representation and 
$t \in T^P$.
\begin{enumerate}
\item $\pi (P,\sigma ,t)$ is unitary if and only if $\sigma$ is
unitary and $t \in T_u^P$
\item $\pi (P,\sigma ,t)$ is (anti-)tempered if and only if 
$\sigma$ is (anti-)tempered and $t \in T_u^P$
\end{enumerate}
\end{prop}
\emph{Proof.}
The "if" statements are \cite[Propositions 4.19 and 4.20]{Opd3},
so we prove the "only if" parts.\\ 
1. Clearly $\sigma \circ \phi_t$ can only be unitary if $\sigma$ 
is. If now $t \in T^P \setminus T_u^P$ then there is an $x \in X 
\cap \big( P^\vee \big)^\perp$ with $|t(x)| \neq 1$. Hence 
\[
\sigma \circ \phi_t (\theta_x^* ) = \sigma \circ \phi_t 
(\theta_{-x}) = \sigma ( t(-x) ) = t(x)^{-1} \neq \overline{t(x)} 
= \sigma ( t(x) )^* = \big( \sigma \circ \phi_t (\theta_x) \big)^*
\]
2. Again it is obvious that $\sigma$ needs to be tempered for 
$\sigma \circ \phi_t$ to be so. If $t$ and $x$ are as above, then
either $|t(x)| > 1$ or $|t(-x)| > 1$, so by Lemma \ref{lem:3.9}
$\sigma \circ \phi_t$ is not tempered. This argument also applies
in the anti-tempered case. $\qquad \Box$
\\[2mm]

Let $\Delta_P$ be the finite set of equivalence classes of 
discrete series representations of \inde{$\mc H_P$} = 
$\mc H(\mc R_P , q_P )$. We pick a representative in every class
and confuse $\Delta_P$ with this set of representations. 
Write \index{Delta@$\Delta$} \index{Deltap@$\Delta_P$}
\[
\Delta = \{ (P,\delta) : P \subset F_0 ,\, \delta \in \Delta_P \}
\]
and denote by $\Xi$ the complex algebraic variety consisting
of all triples $(P,\delta ,t)$ with $P \subset F_0 ,\, \delta \in
\Delta_P$ and $t \in T^P$. Similarly let $\Xi_u$ be the 
compact smooth manifold of all triples $(P,\delta ,t) \in \Xi$ with
$t \in T_u^P$. \index{Xi@$\Xi$} \index{Xiu@$\Xi_u$}

Let $V_\delta$ be the representation space of $\delta$, 
endowed with the inner product $\inp{\cdot}{\cdot}_\delta$, 
and define an inner product on 
$\mc H \big( W^P \big) \otimes_{\mh C} V_\delta$ by
\index{<>@$\inp{\cdot}{\cdot}_\delta$} \index{Vd@$V_\delta$}
\begin{equation}\label{eq:3.25}
\inp{x \otimes u}{y \otimes v} = \tau (x^* y) \inp{u}{v}_\delta 
\qquad x,y \in \mc H \big( W^P \big) \quad u,v \in V_\delta
\end{equation}
With $\xi = (P,\delta ,t) \in \Xi$ we associate the 
parabolically induced representation \index{pix@$\pi (\xi)$}
\index{representation!parabolically induced}
\begin{equation}
\pi (\xi) = \pi (P,\delta ,t) = 
\mr{Ind}_{\mc H^P}^{\mc H} (\delta \circ \phi_t)
\end{equation}
We realize it on $\mc H \big( W^P \big) \otimes_{\mh C} V_\delta$
with the inner product \eqref{eq:3.25}. 

For $P = \es$ we have $\mc H_P \cong \mh C$. We denote its
unique discrete series representation by $\delta_\es$. Note that 
$\pi (\es, \delta_\es ,t)$ is just the principal series 
representation $I_t$. \index{deltaes@$\delta_\es$}.

We gather all parabolically induced representations in the 
following vector bundle over $\Xi$ : \index{VX@$\mc V_{\Xi}$}
\begin{equation}
\mc V_{\Xi} = \bigcup_{(P, \delta ) \in \Delta} T^P \times
\mc H \big( W^P \big) \otimes V_\delta 
\end{equation}

Sometimes we will write $\xi \in \Xi$ as $(P,W_Pr, \delta ,t)$
to indicate that $W_P r \in T_P / W_P$ is the central character 
of $\delta$. Then the central character of $\pi (\xi)$ is 
$W_0 r t \in T / W_0$. \index{pwdt@$(P,W_P r,\delta ,t)$}

If we let $t$ run over $T^P$, we obtain a coset $r T^P$ of the
subtorus $T^P$ in $T$. The collection of all such cosets, for
different $(P,\delta) \in \Delta$ and different representatives of
their central characters, has some special properties. It can be 
described combinatorially with the help of Macdonald's $c$-function,
whose construction we recall now.
For a long root $\alpha \in R_1$ we put \index{$c_\alpha$}
\begin{equation}
c_\alpha = \left(1 + q_{\alpha^\vee}^{-1/2} \theta_{-\alpha / 2} 
\right) \left( 1 - q_{\alpha^\vee}^{-1/2} q_{2 \alpha^\vee}^{-1}
\theta_{-\alpha / 2} \right) \big( 1 - \theta_{-\alpha} 
\big)^{-1} \qquad \alpha \in R_1
\end{equation}
This is a rational function on $T$, i.e. an element of the 
quotient field \inde{$Q (\mc A)$} of $\mc A$. 
Strictly speaking, for $\alpha \in R_0 \cap R_1$ this formula 
is ill-defined, because there is no such thing as 
$\theta_{-\alpha / 2}$. However in that case 
$q_{2 \alpha^\vee} = 1$ by our convention \eqref{eq:3.13}, 
so we can interpret the above formula as
\begin{equation}
c_\alpha = \left(1 - q_{\alpha^\vee}^{-1} \theta_{-\alpha} \right) 
\big( 1 - \theta_{-\alpha} \big)^{-1} 
\hspace{4cm} \alpha \in R_0 \cap R_1
\end{equation}
\inde{Macdonald's $c$-function} is the product \index{$c(t)$}
\begin{equation}
c(t) = \prod_{\alpha \in R_1^+} c_\alpha (t)
\end{equation}
Suppose that $L \subset T$ is a coset of a subtorus \inde{$T^L$}, 
and let \index{$R_L$}
\begin{equation}
R_L = \{ \alpha \in R_1 : \alpha (T^L ) = 1 \}
\end{equation}
be the set of long roots that are constant on $T^L$. Then
\begin{equation}
T^L = \{ t \in T : x(t) = 1 \;\mr{if}\; x \in \mh Q R_L \cap X \} =
\mr{Hom}_{\mh Z} \big( X / (\mh Q R_L \cap X) , \mh C^\times \big)
\end{equation}
Hence $\theta_{-\alpha /2} (L)$ (or $\theta_{-\alpha}(L)$ if 
$\alpha \notin R_0$) is well-defined. Put
\index{$R^p_L$} \index{$R^z_L$}
\begin{equation}\label{eq:3.11}
\begin{aligned}
R_L^p &= \big\{ \alpha \in R_L : \big( 1 + q_{\alpha^\vee}^{-1/2} 
  \theta_{-\alpha / 2} (L) \big) \big( 1 - 
  q_{\alpha^\vee}^{-1/2} q_{2 \alpha^\vee}^{-1} \theta_{-\alpha /2} 
  (L) \big) = 0 \big\} \\
R_L^z &= \left\{ \alpha \in R_L : 1 - \theta_{-\alpha}(L) = 0 
  \right\}
\end{aligned}
\end{equation}
By construction the pole order of the rational function
\begin{equation}\label{eq:3.26}
t \to c^{-1}(t) c^{-1}(t^{-1})
\end{equation}
along $L$ is $|R_L^p| - |R_L^z|$. We say that $L$ is a 
\inde{residual coset} if
\begin{equation}\label{eq:3.27}
|R_L^p| - |R_L^z| = \dim (T) - \dim (L)
\end{equation}
A residual coset of dimension 0 is also called a \inde{residual 
point}. More or less everything there is to know about residual 
cosets can be found in \cite[Section 4]{HeOp} and 
\cite[Appendix A]{Opd3}. The collection of residual cosets is 
finite and $W_0$-invariant. They have been classified completely, 
and it turns out that a coset of a subtorus is residual if and 
only if its pole order for \eqref{eq:3.26} is maximal, 
given its dimension. 

From \eqref{eq:3.11} we see that the residual cosets depend
algebraically on $q$. In particular the number of residual cosets
is maximal for all $q$ in a certain Zariski open subset of the
parameter space of all possible $q$'s. Hence it makes sense to
call $q$ generic if and only if the number of residual cosets
attains its maximum at $q$. \index{label function!generic}

Define \index{$T_L$}
\begin{equation}\label{eq:3.49}
T_L = \left\{ t \in T : t(x) = 1 \;\mr{if}\; x \in 
( R_L^\vee )^\perp \cap X \right\} =
\mr{Hom}_{\mh Z} \left( X \left/ \big( ( R_L^\vee )^\perp 
\cap X \big) \right. ,\mh C^\times \right)
\end{equation}
This is a subtorus of $T$ and $K_L = T_L \cap T^L$ is a finite 
group. Notice that if $L$ is a coset of $T^P$ for some 
$P \subset F_0$, then $R_L = R_P ,\, T_L = T_P$ and $K_L = K_P$.
Pick any $r_L \in T_L \cap L$ and consider 
\[
r_L |r_L|^{-1} \in T_{L,u} \cap L
\] 
Since $r_L$ is unique up to $K_L \cap T_u^L$, it follows that the
\inde{tempered form} \index{$L^{\mr{temp}}$}
\begin{equation}
L^{\mr{temp}} = r_L |r_L|^{-1} T_u^L
\end{equation}
is independent of the choice of $r_L$. The tempered forms of 
residual cosets tend to be disjoint, cf. 
\cite[Theorem A.18]{Opd3} :

\begin{thm}\label{thm:3.14}
Let $L_1$ and $L_2$ be residual cosets. If 
$L_1^{\mr{temp}} \cap L_2^{\mr{temp}} \neq \es$ then 
$L_1 = w (L_2)$ for some $w \in W_0$.
\end{thm}

Our reason for studying residual cosets is a particular 
consequence of \cite[Theorem 4.23]{Opd3} :

\begin{thm}\label{thm:3.15}
The collection of residual cosets is precisely
\[
\left\{ r T^P : (P, W_P r, \delta ,t) \in \Xi \right\}
\]
The tempered form of $L = r T^P$ is $L^{\mr{temp}} = r T_u^P$.
\end{thm}

This theorem implies that for every $W_0$-orbit of residual 
points there is at least one discrete series representation with
that central character. In many cases there is in fact exactly
one such representation:

\begin{prop}\label{prop:3.30}
\begin{enumerate}
\item Let $t \in T$ be a residual point such that the orbit 
$W_0 t$ contains exactly $|W_0|$ points. Then there is, up to 
equivalence, a unique discrete series representation with this 
central character. 
\item If the root system $R_0$ is of type $A_n$ then 1. applies, 
and the representation in question is onedimensional.
\end{enumerate}
\end{prop}
\emph{Proof.}
For 1 see \cite[Corollary 1.2.11]{Slo2}. 
For 2 we use the classification of the residual points for type 
$A_n$ in \cite[Proposition 4.1]{HeOp}. This shows that we can
apply 1. To show that the resulting representation has dimension 
one, we just construct it. Since all simple roots are conjugate, 
we necessarily have an equal label function. If $W = W_{\mr{aff}}$
then either $\pi_{St}$ or $\pi_{\mr{triv}}$ is a onedimensional 
discrete series representation. Which one depends on whether 
$q(s) > 1$ or $q(s) < 1$ for any $s \in S_{\mr{aff}}$, see page 
\pageref{lem:3.8}. If $W \neq W_{\mr{aff}}$ then, given the 
central character, there is unique way to extend this 
representation to $C_r^* (\mc R ,q). \qquad \Box$ 
\\[3mm]

However, in general there may be more than one discrete series
representation with a given central character. This is known to 
happen already for $R_0$ of type $B_2$, for certain label functions, 
see Section \ref{sec:6.4}.

Following \cite[Definition 3.24]{Opd3} we define a radical and
a residual algebra at $t \in T$ :
\index{$\mr{Rad}_t$} \index{$\overline{\mc H^t}$}
\begin{equation}
\begin{aligned}
& \overline{\mc H^t} = \mc H (\mc R ,q) / \mr{Rad}_t \\
& \mr{Rad}_t = \cap_\pi \ker \pi \\
\end{aligned}
\end{equation}
where $\pi$ runs over the irreducible representations of 
$\mc H (\mc R ,q)$ with the following two properties :
\begin{itemize}
\item the central character of $\pi$ is $W_0 t$ 
\item $\pi$ extends to $C_r^* (\mc R ,q)$ 
\end{itemize}
Clearly $\mr{Rad}_t$ contains all elements of $\mc Z$ that vanish 
at $t$, so $\overline{\mc H^t}$ is a finite dimensional 
$C^*$-algebra whose irreducible representations correspond to the 
irreducible representations of $C_r^* (\mc R ,q)$ with central 
character $W_0 t$. By Theorem \ref{thm:3.15} the collection of 
discrete series representations of $\mc H (\mc R,q)$ is precisely 
the union, over all residual points $t$, of the irreducible 
representations of $\overline{\mc H^t}$.
\\[2mm]

As one would expect, two parabolically induced representations 
can be related if they have the same central character. Let us try 
to clarify when and how this occurs. For $P,Q \subset F_0$ write 
\begin{equation}
W (P,Q) = \{ n \in W_0 : n (P) = Q \}
\end{equation}
For such $n$ there is a *-algebra isomorphism $\psi_n$ :
$\mc H^P \to \mc H^Q$, defined by 
\index{$W (P,Q)$} \index{psin@$\psi_n$} \index{psik@$\psi_k$}
\begin{equation}
\psi_n (N_w \theta_x ) = N_{n w n^{-1}} \theta_{n x}
\qquad w \in W_P ,\, x \in X
\end{equation}
Similarly for $k \in K_P$ there is a *-algebra automorphism
$\psi_k : \mc H^P \to \mc H^P$, defined by
\begin{equation}
\psi_k (N_w \theta_x ) = k(x) N_w \theta_x 
\qquad w \in W_P ,\, x \in X
\end{equation}
These maps descend to isomorphisms $\psi_n : \mc H_P \to \mc H_Q$
and $\psi_k : \mc H_P \to \mc H_P$.

Let $\mc W$ be the finite groupoid, over the power set of $F_0$, 
with \index{Wmc@$\mc W$}
\begin{align}
\mc W_{PQ} &= K_Q \times W (P,Q) \\
(k_1 ,n_1 ) (k_2 ,n_2 ) &= (k_1 n_1 (k_2 ), n_1 n_2 )
\qquad (k_1 ,n_1 ) \in \mc W_{PQ} ,\, (k_2 n_2) \in \mc W_{P' P}
\end{align}
For $kn \in$ \inde{$\mc W_{PQ}$} and $\delta \in \Delta_P$, let 
$\Psi_{kn}(\delta )$ be the equivalence class of the discrete series 
representation $\delta \circ \psi_n^{-1} \circ \psi_k^{-1}$ of 
$\mc H_Q$. Then there exists a unitary isomorphism $\tilde 
\delta_{kn}$ : $V_\delta \to V_{\Psi_{kn}(\delta)}$ such that 
\index{deltkn@$\tilde \delta_{kn}$} \index{Psikn@$\Psi_{kn}(\delta )$}
\begin{equation}\label{eq:3.28}
\tilde \delta_{kn} \circ \delta (\psi_n^{-1} \psi_k^{-1} h) =
\Psi_{kn}(\delta ) (h) \circ \tilde \delta_{kn} \qquad h \in \mc H_Q
\end{equation}
By Schur's lemma $\tilde \delta_{kn}$ is unique up to a complex 
number of absolute value 1, but whether or not there is a canonical
way to choose these scalars is an open problem. This is a subtle 
point to which we will return later.

The above maps induce an action of $\mc W$ on $\Xi$ :
\begin{equation}\label{eq:3.29}
kn (P, W_P r, \delta ,t) = 
\big( Q, W_Q n(r), \Psi_{kn}(\delta), k n(t) \big)
\end{equation}

We would like to lift $\tilde \delta_g$ to an intertwiner 
$\pi(g, \xi )$ between $\pi (\xi )$ and $\pi (g \xi )$. 
For $k \in K_p$ this is easy, we can simply define 
\index{pikx@$\pi (k,\xi )$} \index{pigx@$\pi(g, \xi )$}
\begin{equation}\label{eq:3.21}
\pi (k,\xi ) = \mr{id}_{\mc H ( W^P )} \otimes \tilde \delta_k :
\mc H \otimes_{\mc H^P} V_\delta \to 
\mc H \otimes_{\mc H^P} V_\delta
\end{equation}
This is a well-defined intertwiner since, for all 
$h \in \mc H ,\, h' \in \mc H^P ,\, v \in V$ :
\begin{equation}\label{eq:3.30}
\begin{split}
\pi (k,\xi ) (h h' \otimes v) & = h h' \otimes \tilde \delta_k (v)\\ 
& = h \otimes \Psi_k (\delta ) (\phi_{kt} h' ) 
 \big( \tilde \delta_k v \big) \\
& = h \otimes \Psi_k (\delta ) \big( \psi_k (\phi_t h') \big) 
 \big( \tilde \delta_k v \big)\\ 
& = h \otimes \tilde \delta_k \big( \delta (\phi_t h') v \big) \\ 
& = \pi (k,\xi ) \big( h \otimes \delta (\phi_t h') v \big)
\end{split}
\end{equation}
On the other hand, for $w \in W_0$ it is a tricky business, 
which is taken care of in \cite[\S 4]{Opd3}. We shall recall the 
parts of the construction that we use later on.

Let \inde{$U$} be a $W_0$-invariant subset of $T$, open in the 
analytic topology, and let \inde{$C^{an}(U)$} and \inde{$\mc Z^{an} 
(U)$} be the algebras of, respectively, holomorphic functions on 
$U$ and $W_0$-invariant holomorphic functions on  $U$. There are 
natural injections $\mc A \to C^{an}(U)$ and 
$\mc Z \to \mc Z^{an} (U)$, so we can construct 
\begin{equation}\label{eq:3.24}
\mc H^{an} (U) = \mc Z^{an} (U) \otimes_{\mc Z} \mc H =
C^{an}(U) \otimes_{\mc A} \mc H
\end{equation}
Similarly we define the algebras \inde{$C^{me}(U)$} and 
\inde{$\mc Z^{me} (U)$} of ($W_0$-invariant) meromorphic 
functions on $U$, and
\begin{equation}\label{eq:3.14}
\mc H^{me} (U) = \mc Z^{me} (U) \otimes_{\mc Z} \mc H =
C^{me}(U) \otimes_{\mc A} \mc H
\end{equation}
It follows from Theorem \ref{thm:3.4}.1 that 
\inde{$\mc H^{an} (U)$} and \inde{$\mc H^{me} (U)$} are free 
modules over respectively $C^{an}(U)$ and $C^{me}(U)$, with 
bases $\{N_w : w \in W_0\}$. We define a * on these algebras by
\begin{equation}\label{eq:3.17}
f^* = T_{w_0} f^{-w_0} T_{w_0}^{-1} = 
N_{w_0} f^{-w_0} N_{w_0}^{-1}  \qquad f \in C^{me}(U)
\end{equation}
By Lemma \ref{lem:3.6} this extends the * on $\mc H (\mc R ,q)$.

A finite dimensional $\mc H$-representation extends to 
$\mc H^{an} (U)$ if and only if all its $\mc Z$-weights are in
$U / W_0$. Applying this to $\mc H^P$, we see that 
$\delta \circ \phi_t$ extends to $\mc H^{P,an}(U)$ if and only if
$W_P rt \in U / W_P$. The proof of these claims can be
found on \cite[p. 582]{Opd3}. 

Pick $t \in T$ and consider a set 
$B \subset \mr{Lie}(T) = \mf t \otimes_{\mh R} \mh C$
which satisfies the following conditions:
\index{$B$} \index{$U_t$}
\begin{cond}\label{cond:3.7}
\begin{enumerate}
\item $B$ is an open ball centred around 
$0 \in \mf t \otimes_{\mh R} \mh C$
\item $\forall \alpha \in R_{nr}, b \in \overline B : \; 
| \Im (\alpha (b)) | < \pi $ 
(so $\exp : B \to \exp (B)$ is a diffeomorphism)
\item Write $U_t = t \exp (B)$. If $w \in W_0$ and $w 
\overline{U_t} \cap \overline{U_t} \neq \es$, then $wt = t$
\item $\forall u \in \overline{U_t}$ we have
\[
R^p_{\{u\}} \cup R^z_{\{u\}} \subset R^p_{\{t\}} \cup R^z_{\{t\}}
\]
\end{enumerate}
\end{cond}

Recall that $R^p_{\{t\}} \cup R^z_{\{t\}}$ is the set of roots 
$\alpha \in R_1$ for which $t$ is a critical point of the 
rational function $c_\alpha$. Notice that these conditions are 
almost the same as \cite[Conditions 4.9]{Opd3}. In fact the 
only difference is that we replaced some statements with $B$ 
by slightly stronger versions with $\overline B$. This is not 
an essential difference, it only makes it easier in Section 
\ref{sec:5.3} to see that certain functions are bounded. 

By \cite[Proposition 4.7]{Opd3} Condition \ref{cond:3.7} can 
always be fulfilled. This yields an alternative description of 
the representation space of $\pi (P, W_P r, \delta, t)$ :
\begin{equation}\label{eq:3.31}
\mc H \big( W^P \big) \otimes_{\mh C} V_\delta = 
\mc H \otimes_{\mc H^P} V_\delta = \mc H^{an}(W_0 U_{rt})
\otimes_{\mc H^{P,an}(W_P U_{rt})} V_\delta
\end{equation}
Let $s \in S_0$ be the simple reflection corresponding to
$\alpha \in F_1$. Consider the following element of $Q (\mc A) 
\otimes_{\mc A} \mc H = Q (\mc Z) \otimes_{\mc Z} \mc H$ :
\index{i@$\imath^0_w$}
\begin{equation}\label{eq:3.23}
\begin{aligned}
\imath^0_s \;=\; & \big( q_{\alpha^\vee}^{1/2} + 
 \theta_{-\alpha /2} \big)^{-1} \big( q_{\alpha^\vee}^{1/2} 
 q_{2 \alpha^\vee} - \theta_{-\alpha /2} \big)^{-1} \\
   & \big( (1 - \theta_{-\alpha}) T_s + (q_{\alpha^\vee} 
  q_{2 \alpha^\vee} - 1) + q_{\alpha^\vee}^{1/2} 
  (q_{2 \alpha^\vee} - 1) \theta_{-\alpha /2} \big) \\
 \;=\; & \big( T_s (1 - \theta_\alpha) + ( q_{\alpha^\vee} 
  q_{2 \alpha^\vee} - 1) \theta_\alpha + q_{\alpha^\vee}^{1/2} 
  (q_{2 \alpha^\vee} - 1) \theta_{\alpha /2} \big) \\
  & \big( q_{\alpha^\vee}^{1/2} + \theta_{\alpha /2} 
  \big)^{-1} \big( q_{\alpha^\vee}^{1/2} q_{2 \alpha^\vee} - 
  \theta_{\alpha /2} \big)^{-1}
\end{aligned}
\end{equation}
For  $\alpha \in F_1 \cap R_0$ this simplifies to
\[
\imath^0_s \;=\; (q(s) - \theta_{-\alpha} )^{-1} ( (1- 
\theta_{-\alpha} ) T_s + q(s) - 1) \;=\; (T_s (1 - \theta_\alpha ) 
+ (q(s) - 1) \theta_\alpha ) (q(s) - \theta_\alpha )^{-1}
\]
Important properties of such elements come from
\cite[Proposition 5.5]{Lus4} and \cite[Lemma 4.1]{Opd3} :

\begin{thm}\label{thm:3.16}
The map 
\[
S_0 \to Q (\mc A ) \otimes_{\mc A} \mc H : s \to \imath^0_s
\]
extends to a group homomorphism
\[
W_0 \to \big( Q (\mc A ) \otimes_{\mc A} \mc H \big)^\times
\]
For all $w \in W_0, f \in Q (\mc A)$ we have
\begin{equation}\label{eq:3.19}
\imath^0_w f \imath^0_{w^{-1}} = f^w = 
f \circ w^{-1} \in Q (\mc A)
\end{equation}
As $Q (\mc A )$-modules, we can identify
\[
Q (\mc A ) \otimes_{\mc A} \mc H =
\bigoplus_{w \in W_0} \imath^0_w Q (\mc A ) =
\bigoplus_{w \in W_0} Q (\mc A ) \imath^0_w
\]
\end{thm}

If $t \in T$ and $c_\alpha^{-1}(t) = 0$ then by definition 
$\imath^o_{s_\alpha}(t) = 1$. If we combine this with \eqref{eq:3.19}
we see that it is not specific for simple reflections:
\begin{equation}\label{eq:3.52}
\beta \in R_0 ,\, c_\beta^{-1}(t) = 0 \quad \Rightarrow \quad 
\imath^o_{s_\beta}(t) = 1
\end{equation}
With a straightforward but tedious computation one may show that
\begin{equation}\label{eq:3.18}
\big( \imath_w^0 \big)^* = T_{w_0} \prod_{\alpha \in R_1^+ 
\cap w' R_1^- } \left( {\ds \frac{c_\alpha}{c_{-\alpha}}} 
\right) \imath^0_{w'} T_{w_0}^{-1}
\end{equation}
where $w' = w_0 w^{-1} w_0$. It follows more readily from the 
Bernstein-Lusztig-Zelevinski relations \eqref{eq:3.10} that, 
for $n \in W (P,Q)$ and $w \in W_P$,
\begin{equation}\label{eq:3.32}
\imath^0_{n^{-1}} N_w \imath^0_n = N_{n w n^{-1}}
\end{equation}
Hence we have the following equality in 
$Q (\mc A) \otimes_{\mc A} \mc H$ :
\begin{equation}\label{eq:3.33}
\imath^0_{n^{-1}} h \imath^0_n = \psi_n (h) \qquad h \in \mc H^P
\end{equation}
For $\xi = (P, W_P r, \delta, t) \in \Xi$ we define
\index{pinx@$\pi (n,\xi )$}
\begin{equation}\label{eq:3.22}
\begin{aligned}
&\pi (n,\xi ) : \mc H^{an}(W_0 U_{rt}) \otimes_{\mc H^{P,an}
(W_P U_{rt})} V_\delta \to \mc H^{an}(W_0 U_{rt}) 
\otimes_{\mc H^{Q,an} (W_Q U_{n(rt)})} V_{\Psi_n (\delta )} \\
&\pi (n,\xi ) (h \otimes v) = h \, \imath^0_{n^{-1}} \otimes 
\tilde \delta_n (v) \qquad h \in \mc H , v \in V_\delta
\end{aligned}
\end{equation}
With \eqref{eq:3.31} in mind, this map intertwines the 
$\mc H$-representations $\pi (\xi)$ and $\pi (n \xi)$ since,
for $h \in \mc H ,\, h' \in \mc H^P ,\, v \in V_\delta$ :
\begin{equation}\label{eq:3.34}
\begin{split}
\pi (n,\xi ) (h h' \otimes v) & = h h' \imath^0_{n^{-1}} 
 \otimes \tilde \delta_n (v) \\
& = h \imath^0_{n^{-1}} ( \imath^0_n h' \imath^0_{n^{-1}} ) 
 \otimes \tilde \delta_n (v) \\
& = h \, \imath^0_{n^{-1}} \otimes \Psi_n (\delta ) 
 \big( \phi_{n(t)} ( \imath^0_n h' \imath^0_{n^{-1}} ) \big) 
 \big( \tilde \delta_n v \big) \\
& = h \, \imath^0_{n^{-1}} \otimes \Psi_n (\delta ) \big( 
 \phi_{n(t)} (\psi_n h' ) \big) \big( \tilde \delta_n v \big) \\
& = h \, \imath^0_{n^{-1}} \otimes \Psi_n (\delta ) \big( \psi_n
 (\phi_t h' ) \big) \big( \tilde \delta_n v \big) \\
& = h \, \imath^0_{n^{-1}} \otimes \tilde \delta_n 
 \big( \delta (\phi_t h') v \big) \\
& = \pi (n,\xi ) \big( h \otimes \delta (\phi_t h') v \big)
\end{split}
\end{equation}
For $g = k n \in \mc W$ we define the \inde{intertwiner} 
\index{pigx@$\pi (g,\xi)$}
\begin{equation}\label{eq:3.35}
\pi (g, \xi) = \pi (k, n \xi) \pi (n, \xi)
\end{equation}
Because every $\tilde \delta_g$ is unique up to a scalar, there
is a \inde{$c (g_1 ,g_2 ,\delta )$} $\in \mh T$ such that
\begin{equation}\label{eq:3.36}
\pi (g_1 ,g_2 \xi ) \pi (g_2 ,\xi ) = 
c (g_1 ,g_2 ,\delta ) \pi (g_1 g_2 ,\xi ) 
\end{equation}
Things would simplify considerably if we could arrange that all 
$c (g_1 ,g_2 ,\delta )$ would be 1. In many cases this is indeed
possible, but there are some examples for which it cannot be done.

A priori $\pi (kn,\xi)$ is well-defined only on the Zariski-open 
subset of $(P,\delta ,T^P )$ where $ \imath^0_n$ and 
$\imath^0_{n^{-1}}$ are regular. In a worst case scenario this set 
could even be empty, but fortunately it is assured by 
\cite[Theorem 4.33 and Corollary 4.34]{Opd3} that

\begin{thm}\label{thm:3.17}
For any $g \in \mc W_{PQ}$ intertwining map 
\[
\pi (g,\xi ) : \mc H \big( W^P \big) \otimes V_\delta \to
\mc H \big( W^Q \big) \otimes V_{\Psi_g (\delta)}
\]
is rational as a function of $t \in T^P$. It is regular on 
an analytically open neighborhood of $T^P_u$. Moreover, if 
$t \in T_u^P$ then $\pi (g,\xi )$ is unitary with respect to 
the inner products defined by \eqref{eq:3.25}.
\end{thm}

With this in mind, let $\Gamma_{rr} \left( \Xi ;\mr{End} 
(\mc V_\Xi ) \right)$ be the algebra of rational sections of 
End$(\mc V_\Xi)$ that are regular on $\Xi_u$. It can be described 
more explicitly as \index{Gammarxe@$\Gamma_{rr} \left( \Xi 
;\mr{End} (\mc V_\Xi ) \right)$} 
\begin{multline}
\Gamma_{rr} \big( \Xi ;\mr{End} (\mc V_\Xi ) \big) =
\bigoplus_{(P,\delta) \in \Delta} \Gamma_{rr} \left( T^P ; \mr{End}
\big( \mc H \big( W^P \big) \otimes V_\delta \big) \right) \\
= \bigoplus_{(P,\delta) \in \Delta} \left\{ f \in Q \big( 
\mc O \big( T^P \big) \big) \otimes \mr{End} \big( 
\mc H \big( W^P \big) \otimes V_\delta \big) : 
f \mr{\;is\;regular\;on\;} T^P_u \right\}
\end{multline}
Obviously this algebra contains the algebra of polynomial sections
\index{$\mc O (\Xi ;\mr{End} (\mc V_\Xi ))$}
\begin{equation}
\mc O \big( \Xi ;\mr{End} (\mc V_\Xi ) \big) =
\bigoplus_{(P,\delta) \in \Delta} \mc O \big( T^P \big) \otimes 
\mr{End} \big( \mc H \big( W^P \big) \otimes V_\delta \big)
\end{equation}
Using the standard analytic structure on $\Xi$, we define the
algebras \inde{$C (\Xi ;\mr{End} (\mc V_\Xi ))$} and
\inde{$C^\infty (\Xi ;\mr{End} (\mc V_\Xi ))$} of continuous
(respectively smooth) sections of End$(\mc V_\Xi )$ in the same 
way. Furthermore, if $\mu$ is a sufficiently "nice" measure 
on $\Xi$ then the following formula defines a (degenerate) 
Hermitian inner product on $C (\Xi ;\mr{End} (\mc V_{\Xi}))$:
\begin{equation}
\inp{f_1}{f_2}_\mu = \int_\Xi \mr{tr} 
\big( f_1 (\xi )^* f_2 (\xi ) \big) d \mu (\xi)
\end{equation}
We denote the corresponding Hilbert space completion by 
\inde{$L^2 (\Xi ,\mu ;\mr{End} (\mc V_\Xi ))$}. \index{mu@$\mu$}

From the action of the groupoid $\mc W$ on $\Xi$ and the above 
intertwiners we get an action of $\mc W$ on 
$\Gamma_{rr} \left( \Xi ;\mr{End} (\mc V_\Xi ) \right)$
by algebra homomorphisms :
\begin{equation}\label{eq:3.37}
(g \cdot f) (\xi ) = \pi (g, g^{-1} \xi ) f (g^{-1} \xi ) 
\pi (g, g^{-1} \xi )^{-1}
\end{equation}
whenever $g^{-1} \xi$ is defined. 
The average of $f$ over $\mc W$ is
\begin{equation}
p_{\mc W} (f) (P,\delta ,t) = 
\big| \{ (Q,g) : Q \subset F_0 ,\, g \in \mc W_{QP} \} \big|^{-1} 
\sum_{Q \subset F_0} \sum_{g \in \mc W_{QP}} 
(g \cdot f)(P,\delta ,t)
\end{equation}
Notice that for $f \in \mc O (\Xi ;\mr{End} (\mc V_\Xi ))$ and 
$g \in \mc W ,\; g \cdot f$ and $p_{\mc W} (f)$ need not lie in
$C (\Xi ;\mr{End} (\mc V_\Xi ))$. However, if $M$ is a 
$\mc W$-stable subset of $\Xi$ on which all the intertwiners 
are regular, then there is a well-defined action of $\mc W$ 
on $C (M ;\mr{End} (\mc V_\Xi ))$. The same holds for smooth 
sections if $M$ is a smooth submanifold of $\Xi$ on top of this.

The \inde{Fourier transform} is the algebra homomorphism 
\begin{equation}
\begin{aligned}
& \mc F : \mc H \to \mc O \big( \Xi ;\mr{End} (\mc V_\Xi ) \big) \\
& \mc F (h)(\xi ) = \pi (\xi )(h)
\end{aligned}
\end{equation}
We can extend it continously to various completions of 
$\mc H (\mc R ,q)$. After doing so, its image can be described 
completely with our intertwiners. For the Hilbert space completion 
$\mf H (\mc R,q)$ the following was proved in \cite[Theorem 4.43 
and Corollary 4.45]{Opd3} : \index{Fmc@$\mc F$}

\begin{thm}\label{thm:3.18}
There exists a unique "Plancherel" measure $\mu_{Pl}$ on $\Xi$ 
with the properties \index{mupl@$\mu_{Pl}$} 
\begin{enumerate}
\item the support of $\mu_{Pl}$ is $\Xi_u$
\item $\mu_{Pl}$ is $W_0$-invariant
\item on every component $(P,\delta, T^P) \; \mu_{Pl}$ is 
absolutely continuous with respect to the Haar measure of $T_u^P$
\item the Fourier transform extends to a bijective isometry
\[
\mc F : \mf H(\mc R,q) \to 
L^2 (\Xi_u ,\mu_{Pl} ;\mr{End} (\mc V_\Xi ))^{\mc W}
\]
i.e. $\mu_{Pl}$ is the \inde{Plancherel measure} for $\tau$
\item the adjoint map
\[
\mc J : L^2 (\Xi_u ,\mu_{Pl} ;\mr{End} (\mc V_\Xi )) \to
\mf H (\mc R, q)
\]
satisfies $\mc J \mc F = \mr{id}_{\mf H (\mc R ,q)}$ and 
$\mc F \mc J = p_{\mc W}$.
\end{enumerate}
\end{thm}

The corresponding statements for the Schwartz and $C^*$-completions
are \cite[Theorem 4.3 and Corollary 4.7]{DeOp1} :

\begin{thm}\label{thm:3.19}
The Fourier transform induces algebra homomorphisms
\[
\begin{array}{rrr}
\mc H (\mc R ,q) & \to & 
  \mc O \big( \Xi ;\mr{End} (\mc V_\Xi ) \big)^{\mc W}\\
\mc S (\mc R ,q) & \to & 
  C^\infty \big( \Xi_u ;\mr{End} (\mc V_\Xi ) \big)^{\mc W}\\
C_r^* (\mc R ,q) & \to & 
  C \big( \Xi_u ;\mr{End} (\mc V_\Xi ) \big)^{\mc W}
\end{array}
\]
The first one is injective, the second is an isomorphism of
Fr\'echet *-algebras and the third is an isomorphism of
$C^*$-algebras.
\end{thm}

Some remarkable consequences of this theorem are

\begin{cor}\label{cor:3.20}
\begin{enumerate}
\item The centers of $\mc S (\mc R ,q)$ and $C^*_r (\mc R ,q)$ are
\begin{align*}
Z \big( \mc S (\mc R ,q) \big) &\cong C^\infty (\Xi_u )^{\mc W}\\
Z \big( C_r^* (\mc R ,q) \big) &\cong C (\Xi_u )^{\mc W}
\end{align*}
\item Every irreducible tempered $\mc H (\mc R ,q)$-representation
$(\pi, V)$  is a direct summand of $\pi (\xi)$, for some 
$\xi \in \Xi_u$. In particular we can endow $V$ with an inner 
product such that $\pi$ is unitary and extends to $C_r^* (\mc R ,q)$.
\item For any $\xi \in \Xi$ and $g \in \mc W$ such that $g \xi$
is defined, the $\mc H (\mc R ,q)$-representations $\pi (\xi)$ and
$\pi (g \xi)$ have the same irreducible subquotients, counted
with multiplicity.
\end{enumerate}
\end{cor}
\emph{Proof.} 
1 comes from \cite[Corollary 4.5]{DeOp1}.\\
2. As an $\mc S (\mc R ,q)$-representation $(\pi ,V)$ has a 
single $Z (\mc S (\mc R ,q))$-weight\\ $\mc W \xi \in \Xi_u / \mc W$.
Then it is also an irreducible representation of the finite 
dimensional $C^*$-algebra 
\[
\mc S (\mc R ,q) / \{ h \in \mc S (\mc R ,q) : \pi (g \xi)(h) = 0
\; \forall g \in \mc W \}
\]
By Theorem \ref{thm:3.19} every irreducible representation of 
this algebra is a direct summand of a parabolically induced 
representation $\pi (g \xi)$. Again by Theorem \ref{thm:3.19} this 
means that $(\pi, V)$ extends to a (unitary) representation of 
$C^*_r (\mc R ,q).$\\
3. See \cite[Proposition 6.1]{DeOp2}. We have to show that the 
characters of $\pi (\xi)$ and $\pi (g \xi)$ are equal, i.e. that
the function
\begin{equation}
\mc H \times T^P \to \mh C : (h,t) \to \mr{tr}\: 
\pi (P,\delta, t)(h) - \mr{tr}\: \pi(Q, \Psi_g (\delta), g(t)) (h)
\end{equation}
is identically 0. Because this is a polynomial function of $t$, it
suffices to show that it is 0 for all $t \in T^P_u$. This follows
directly from Theorem \ref{thm:3.19}. $\qquad \Box$
\\

\begin{cor}\label{cor:3.39}
\begin{enumerate}
\item An element $h \in \mc S (\mc R ,q)$ is invertible in
$\mc S (\mc R ,q)$ if and only if it is invertible in 
$C_r^* (\mc R ,q)$, which happens if and only if it is invertible
in $B (\mf H (\mc R ,q))$.
\item The set of invertible elements $\mc S (\mc R ,q)^\times$ 
is open in $\mc S (\mc R ,q)$, and inverting is a continuous map 
from this set to itself.
\end{enumerate}
\end{cor}
\emph{Proof.} 
By \cite[Proposition 4.8]{Tak} and Theorem \ref{thm:3.19} 
the inclusions 
\begin{equation}
\mc S (\mc R ,q) \to C_r^* (\mc R ,q) \to B (\mf H (\mc R ,q))
\end{equation}
are isospectral, which is another way to state 1. Because the set 
of invertibles is always open in a Banach algebra, so is
\begin{equation}
\mc S (\mc R ,q)^\times \cong
C^\infty \big( \Xi_u ;\mr{End} (\mc V_\Xi ) \big)^{\mc W} \cap
C \big( \Xi_u ;\mr{End} (\mc V_\Xi ) \big)^\times
\end{equation}
Finally, by Theorem \ref{thm:2.7}.2 inverting is continuous on 
$\mc S (\mc R ,q)^\times . \qquad \Box$
\\[3mm]

Now we will combine the Langlands classification with the Fourier
transform to obtain a finer classification of irreducible 
$\mc H$-modules. Although the proofs of the following results are
mostly in \cite{DeOp2}, we give them anyway, because we want to 
generalize them to Hecke algebras of reductive $p$-adic groups.

For $Q \subset F_0$, let \inde{$\mc W^Q$}, $\Xi^Q ,\, \pi^Q (\xi)$ 
etcetera denote the same as the corresponding objects without the
superscript $Q$, but now for the affine Hecke algebra $\mc H^Q$. 
For $\xi = (P, W_P r, \delta, t) \in \Xi$ we put \index{$P (\xi)$}
\index{XiQ@$\Xi^Q$} \index{piQX@$\pi^Q (\xi)$}
\begin{equation}
P (\xi) = \{ \alpha \in F_0 : |\alpha (t)| = 1 \} \supset P
\end{equation}
We will study $\pi (\xi )$ by induction in stages, first we 
construct $\pi^{P (\xi )} (\xi )$ and then we induce that 
representation to $\mc H$.
The following result from \cite[Proposition 6.17.1]{DeOp2}
provides the link with the Langlands classification.

\begin{lem}\label{lem:3.25}
The $\mc H^{P (\xi )}$-representation $\pi^{P (\xi)} (\xi )$ 
is essentially tempered and completely reducible.
\end{lem}
\emph{Proof.} 
By Proposition \ref{prop:3.13} 
$\pi^{P (\xi )} (P,\delta, t |t|^{-1} )$ is tempered and 
unitary. From the proof of Lemma \ref{lem:3.21} we see that
\begin{equation}\label{eq:3.60}
\pi^{P (\xi )} (\xi ) \circ \psi_{|t|}^{-1} = 
\pi^{P (\xi )} (P,\delta, t |t|^{-1} )
\end{equation}
so $\pi^{P (\xi )} (\xi )$ is completely reducible. 
If $r_1 ,\ldots, r_d$ are the $\mc A$-weights of \eqref{eq:3.60}, 
then $|t| r_1 ,\ldots, |t| r_d$ are the $\mc A$-weights of
$\pi^{P (\xi )}(\xi)$. But for $x \in \mh Z P (\xi )$ we have 
$|r_i t (x)| = |r_i (x)|$, so $\pi^{P (\xi )} (\xi )$ is 
essentially tempered. $\qquad \Box$
\\[2mm]

Similar to the definition of Langlands data, we need a kind of
positivity condition on $\Xi$. Thus we say that $\xi \in \Xi^+$ 
if $|t| \in \overline{T_{rs}^{P,+}} = 
\exp \overline{\mf t^{P,+}}$. \index{Xi+@$\Xi^+$}

\begin{prop}\label{prop:3.22}
Take $\xi = (P,W_P r, \delta, t) \in \Xi^+$.
\begin{enumerate}
\item Let $\sigma$ be an irreducible direct summand of 
$\pi^{P (\xi)} (\xi )$. Then $(P (\xi ), \sigma)
\in \Lambda^+$. 
\item The functor $\mr{Ind}_{\mc H^{P (\xi )}}^{\mc H}$ 
induces an isomorphism 
\[
\mr{End}_{\mc H} (\pi (\xi )) \cong 
\mr{End}_{\mc H^{P (\xi )}} \big( \pi^{P (\xi )} (\xi )\big)
\]
\item The irreducible quotients of $\pi (\xi )$ are precisely 
the modules $L (P (\xi ), \sigma)$ with $\sigma$ as above.
\item Every $L (P (\xi ), \sigma)$ has an $\mc A$-weight 
$t_\sigma$ such that there exists a root subsystem 
$R_\sigma \subset R_0$, of rank $|P|$, with the properties
\begin{equation}\label{eq:3.39} 
\begin{split}
& \forall \alpha \in R_0^+ \cap R_\sigma : 
|t_\sigma (\alpha )| < 1\\
& \forall \alpha \in R_0^+ \cap \big( R_\sigma^\vee 
\big)^\perp : |t_\sigma (\alpha )| \geq 1
\end{split}
\end{equation}
\item Suppose that $(\rho ,V)$ is an irreducible constituent 
of $\pi (\xi)$ which is not a quotient, and that $t_\rho$ 
is an $\mc A$-weight of $V$. Then every root subsystem $R_\rho$
with the properties \eqref{eq:3.39} has rank $> |P|$.
\end{enumerate}
\end{prop}
\emph{Proof.} 
1. By definition $r_\sigma = |t| \in T_{rs}^{P,+}$, so 
$(P,\sigma) \in \Lambda^+$.\\ 
2 and 3 follow from 1 and Theorem \ref{thm:3.11}. \\
4. By Lemma \ref{lem:3.10} every $\mc A$-weight of $\sigma$
is of the form $t_\sigma = w (r t)$ with $r$ an 
$\mc A_P$-weight of $\delta$ and
\[
w \in W_{P(\xi )} \cap W^P = 
\big\{ w \in W_{P (\xi )} : w (P) \subset R_0^+ \big\}
\]
By \cite[(2.7)]{Eve} $L (P (\xi ),\sigma)$ also has an
$\mc A$-weight of this form. Hence we may take
$R_\sigma = w (R_P)$.\\
5. By \cite[(2.7)]{Eve} every $\mc A$-weight of $(\rho ,V)$ 
is of the form $t_\rho = n w (r t)$, with 
$n \in W^Q \setminus \{e\}$ and $w, r, t$ as in the proof of 4. 
Clearly, for 
\[
\alpha \in n w (R_P^+) \subset n \big( R_Q^+ \big) 
\subset R_0^+
\]
we have 
\[
\big| t_\rho (\alpha )\big| = 
\big| rt \big( w^{-1} n^{-1} \alpha \big) \big| < 1
\]
Furthermore, since $n \neq e$, there exists a $\beta \in 
R_0^+$ with $n^{-1}(\beta ) \in R_0^-$ and $\beta \perp n 
(P(\xi ))^\vee$. Now $P' = n w (P) \cup \{\beta \}$ 
is a linearly independent set of positive roots such that 
$|t_\rho (\alpha )| < 1 \; \forall \alpha \in P'$. Therefore
any suitable root system $R_\rho$ must have rank at least 
$|P'| = |P| + 1. \qquad \Box$
\\[3mm]
Although it is written down differently, the proof of 
Proposition \ref{prop:3.22} is rather similar to that of the 
Langlands classification. In particular we use the same ideas 
as Langlands' geometric lemmas \cite[p. 61-63]{Lan}.

\begin{lem}\label{lem:3.23}
Every $\xi \in \Xi$ is $W_0$-associate to an element of $\Xi^+$.
If $\xi_1 ,\xi_2 \in \Xi^+$ are $\mc W$-associate, then 
$P (\xi_1 ) = P (\xi_2 ) := Q$, and $\pi^Q (\xi_1 )$ and 
$\pi^Q (\xi_2 )$ are equivalent as $\mc H^Q$-representations.
\end{lem}
\emph{Proof.} 
By \cite[Section 1.15]{Hum} every $W_0$-orbit in $\mf t$ 
contains a unique point in a positive chamber $\mf t^{Q,+}$. 
Hence $|t_1 | = |t_2 |$ and $P (\xi_1 ) = P (\xi_2 ) = Q$. 
From Lemmas \ref{lem:3.21} and \ref{lem:3.25} we see that there 
is a single automorphism $\psi_{|t_1 |} = \psi_{|t_2 |} := \psi$ 
of $\mc H^Q$ such that, for $i = 1 , 2$, 
\begin{equation}\label{eq:3.38}
\pi^Q (\xi_i ) \circ \psi^{-1} \cong \pi^Q (\xi'_i )
\quad \mr{where} \quad \xi'_i = (P_i, \delta_i, t_i |t_i|^{-1})
\in \Xi_u^Q
\end{equation}
If $g \xi_1 = \xi_2$ for some $g \in \mc W$, then also
$g \xi'_1 = \xi'_2$. Applying Theorem \ref{thm:3.19} to 
$\mc H^Q$, we see that $\pi^Q (\xi'_1 )$ and $\pi^Q (\xi'_2 )$
are unitarily equivalent. It follows from this and
\eqref{eq:3.38} that $\pi^Q (\xi_1 )$ and $\pi^Q (\xi_2 )$
are equivalent. $\qquad \Box$
\\

\begin{thm}\label{thm:3.24}
For every irreducible $\mc H$-representation $\pi$ there 
exists a unique association class $\mc W (P,\delta ,t) \in \Xi 
/ \mc W$ such that the following (equivalent) statements hold :
\begin{enumerate}
\item $\pi$ is equivalent to an irreducible quotient of
$\pi (\xi^+ )$, for some 
$\xi^+ \in \Xi^+ \cap \mc W (P,\delta ,t)$
\item $\pi$ is equivalent to an irreducible subquotient of
$\pi (P,\delta, t)$, and $|P|$ is maximal for this property
\end{enumerate}
\end{thm}
\emph{Proof.} 
For 1 we copy \cite[Corollary 6.19]{DeOp2}.
By Theorem \ref{thm:3.11} there is a unique Langlands datum
$(Q,\sigma) \in \Lambda^+$ such that $\pi \cong L (P,\sigma )$,
and by Lemma \ref{lem:3.21} $\sigma \circ \psi_{r_\sigma}^{-1}$
is tempered. Now Theorem \ref{thm:3.19} tells us that there 
exists a unique association class $\mc W^Q \xi_1 = 
\mc W^Q (P_1 , \delta_1 ,t_1 ) \in \Xi_u^Q / \mc W^Q$ such that 
$\pi (\xi_1 )$ contains $\sigma \circ \psi_{r_\sigma}^{-1}$
as an irreducible direct summand. Put 
$\xi = (P_1 ,\delta_1 ,t_1 r_\sigma ) \in \Xi^Q$ and, using
Lemma \ref{lem:3.23} for $\mc H^Q$, pick 
$\xi^+ = (P_2 , \delta_2 , t_2) \in \mc W^Q \xi \cap \Xi^+$.
Then $\sigma$ is a direct summand of $\pi^Q (\xi^+ )$, and we
see from Proposition \ref{prop:3.22}.3 that $\pi$ is an
irreducible subquotient of $\pi (\xi^+ )$.
By Lemma \ref{lem:3.23} and Theorem \ref{thm:3.11} the class 
$\mc W \xi \in \Xi / \mc W$ is unique for this property.

Suppose that $\pi$ is also an irreducible subquotient of
$\pi (\xi_3) = \pi (P_3 , \delta_3 , t_3 )$, where 
$|P_3 | \geq |P|$. By Corollary \ref{cor:3.20}.2 we may assume 
that $\xi_3 \in \Xi^+$. Comparing parts 4 and 5 of Proposition
\ref{prop:3.22} we see that in fact $\pi$ must be equivalent to
an irreducible quotient of $\pi (\xi_3 )$. But then $\xi_3$
is $\mc W$-associate to $\xi^+$ by the above.
Because the class $\mc W \xi$ is unique for both properties
1 and 2, this also shows that 1 and 2 are equivalent.
$\qquad \Box$.
\\[4mm]

\section{Periodic cyclic homology}
\label{sec:3.4}

We prove comparison theorems between the periodic cyclic 
homology of an affine Hecke algebra, that of its Schwartz 
completion, and the $K$-theory of its $C^*$-completion. 
We reap some fruits from our previous labor in the sense
that the technicalities are limited, although still 
substantial.

Let $\mc H (\mc R ,q)$ be an affine Hecke algebra. Observe 
that by Theorems \ref{thm:2.18} and \ref{thm:3.19} the Chern 
character gives an isomorphism
\begin{equation}\label{eq:3.68}
K_* (\mc S (\mc R ,q)) \otimes \mh C \isom HP_* (\mc S (\mc R ,q))
\end{equation}
However, we cannot apply Theorem \ref{thm:2.20} to 
$\mc H (\mc R ,q)$, since the action of the groupoid $\mc W$ 
on $\mc O (\Xi ; \mr{End} (\mc V_{\Xi}))$ is by rational 
intertwiners, which may have poles outside $\Xi_u$. Therefore 
the proof of the next theorem will involve several steps.

\begin{thm}\label{thm:3.38}
The inclusion $\mc H (\mc R ,q) \to \mc S (\mc R ,q)$ induces an
isomorphism 
\[
HP_* (\mc H (\mc R ,q)) \isom HP_* (\mc S (\mc R ,q))
\]
\end{thm}
\emph{Proof.} 
We start by constructing stratifications of the primitive ideal
spectra of $\mc H (\mc R ,q)$ and $\mc S (\mc R ,q)$. Choose an 
increasing chain
\[
\es = \Delta_0 \subset \Delta_1 \subset \cdots \subset 
\Delta_n = \Delta
\]
of $\mc W$-invariant subsets of $\Delta$ with the properties
\begin{itemize}
\item if $(P,\delta) \in \Delta_i$ and $|Q| > |P|$ then 
      $\Delta_Q \subset \Delta_i$
\item the elements of $\Delta_i \setminus \Delta_{i-1}$ form
      exactly one association class for the action of $\mc W$
\end{itemize}
To this correspond two decreasing chains of ideals
\begin{equation}\label{eq:3.53}
\begin{aligned}
\mc H &= I_0 \supset I_1 \supset \cdots \supset I_n = 0 \\ 
\mc S &= J_0 \supset J_1 \supset \cdots \supset J_n = 0 \\
I_i &= \{ h \in \mc H : \pi(P,t,\delta)(h) = 0 \;\mr{if}\: 
  (P,\delta) \in \Delta_i, t \in T^P \} \\
J_i &= \{ h \in \mc S : \pi(P,t,\delta)(h) = 0 \;\mr{if}\: 
  (P,\delta) \in \Delta_i, t \in T_u^P \}
\end{aligned}
\end{equation}
For every $i$ pick an element $(P_i,\delta_i) \in \Delta_i 
\setminus \Delta_{i-1}$, let $\mc W_i$ be the stabilizer of
$(P_i,\delta_i)$ in $\mc W$ and write $V_i = 
V_{\pi(P_i,t,\delta_i)}$. By Theorem \ref{thm:3.19} the
Fourier transform gives isomorphisms
\begin{equation}\label{eq:3.69}
J_{i-1} / J_i \cong C^\infty (T_u^{P_i}; \mr{End}\: V_i)^{\mc W_i}
\end{equation}
On the other hand, by Theorem \ref{thm:3.24} the primitive ideal 
spectrum of $(I_{i-1} / I_i)$ corresponds to the 
inverse image of $\Delta_i \setminus \Delta_j$ under the 
projection $\Xi \to \Delta$. Moreover the induced map 
Prim$(I_{i-1} / I_i) \to \mc W_i \setminus T^{P_i}$ is continuous.
(In fact it is the central character map for this algebra.)
By Lemma \ref{lem:2.12} it suffices to show that each inclusion
\begin{equation}\label{eq:3.73}
I_{i-1} / I_i \to J_{i-1} / J_i
\end{equation}
induces an isomorphism on periodic cyclic homology. Therefore we 
zoom in on $J_{i-1} / J_i$. By Theorem \ref{thm:3.17} we can 
extend the action of $\mc W_i$ on $C^\infty (T_u^{P_i}; \mr{End}\; 
V_i)$ to a neighborhood $T'$ of $T_u^{P_i}$. We may take $T' \;
\mc W_i$-equivariantly diffeomorphic to $T_u^{P_i} \times 
[-1,1]^{\dim T_u^{P_i}}$. Because $[-1,1]$ is compact and 
contractible, we can make the inner product on $V_i$ depend on
$t \in T'$ in a smooth way, such that the intertwiners 
$\pi (g,P_i ,t,\delta_i )$ are unitary on all of $T'$. To avoid
some technical difficulties we want to replace $J_{i-1} / J_i$ by
$C^\infty (T' ; \mr{End}\; V_i )^{\mc W}$, but this needs some
justification.

\begin{lem}\label{lem:3.40}
The inclusion $T_u^{P_i} \to T'$ and the Chern character induce 
isomorphisms
\[
HP_* (J_{i-1} / J_i ) \mosi HP_* \left( C^\infty (T'; \mr{End}\: 
V_i)^{\mc W_i} \right) \isom K_* \left( C(T'; \mr{End}\: 
V_i)^{\mc W_i} \right) \otimes \mh C
\]
\end{lem}
\emph{Proof.} 
The second isomorphism follows directly from Theorems
\ref{thm:2.18} and \ref{thm:2.19}. For the first one, we pick a
$\mc W_i$-equivariant triangulation $\Sigma \to T_u^{P_i}$ and we
construct $U_\sigma$ and $U_\sigma$ as on page \pageref{eq:2.70}.
Using the projection $p_u : T' \to T_u^{P_i}$ we get a closed cover
of $T'$ :
\begin{align*}
& \{ T'_\sigma : \sigma \text{ simplex of } \Sigma \} \\
&T'_\sigma = p_u^{-1} (U_\sigma ) \cong U_\sigma \times 
 [-1,1]^{\dim T_u^{P_i}}
\end{align*}
From the proof of Theorem \ref{thm:2.18} we see that it suffices to
show that for any simplex $\sigma$ we have
\begin{equation}\label{eq:3.70}
HP_* \left( C_0^\infty (U_\sigma, D_\sigma; \mr{End}\: 
V_i )^{\mc W_\sigma} \right) \cong HP_* \left( C_0^\infty (T'_\sigma, 
p_u^{-1}(D_\sigma ); \mr{End}\: V_i)^{\mc W_\sigma} \right)
\end{equation}
where $\mc W_\sigma$ is the stabilizer of $\sigma$ in $\mc W_i$.
Well, $U_\sigma \setminus D_\sigma$ is $\mc W_\sigma$-equivariantly\\
contractible by construction, and it is an equivariant deformation 
retract of\\ 
$T'_\sigma \setminus p_u^{-1}(D_\sigma) = p_u^{-1} (U_\sigma 
\setminus D_\sigma)$. So we are in the setting of Lemma
\ref{lem:2.17} and we may use its proof. It says that there 
exist a finite central extension $G$ of $\mc W_\sigma$ and a 
linear representation
\[
G \to GL (V_i ) : g \to u_g 
\]
such that the Fr\'echet algebras in \eqref{eq:3.70} are isomorphic to
\begin{align}
\label{eq:3.71} C_0^\infty (U_\sigma, D_\sigma; \mr{End}\: V_i)^G \\
\label{eq:3.72} C_0^\infty (T'_\sigma, p_u^{-1}(D_\sigma); 
  \mr{End}\: V_i)^G
\end{align}
The $G$-action on these algebras is given by
\[
g(f)(t) = u_g f(g^{-1} t) u_g^{-1}
\]
where we simply lifted the action of $\mc W_\sigma$ on $T'_\sigma$ 
to $G$. 

It is clear that the retraction $T'_\sigma \to U_\sigma$ induces a 
diffeotopy equivalence between \eqref{eq:3.71} and \eqref{eq:3.72}, 
so it also induces the required isomorphism \eqref{eq:3.70}. 
$\qquad \Box$
\\[2mm]

Consider the finite collection $\mc L$ of all irreducible components
of $(T^{P_i})^w$, as $w$ runs over $\mc W_i$. These are all cosets
of complex subtori of $T^{P_i}$ and they have nonempty intersections 
with $T_u^{P_i}$. Extend this to a collection $\{ L_j \}_j$ of 
cosets of subtori of $T^{P_i}$ by including all irreducible
components of intersections of any number of elements of $\mc L$.
Because the action $\alpha_i$ of $\mc W_i$ on $T^{P_i}$ is algebraic
\[
\dim \big( \left( T^{P_i} \right)^g \cap \left( T^{P_i} \right)^w 
\big) < \max \big\{ \dim \left( T^{P_i} \right)^g ,\:  
\dim \left( T^{P_i} \right)^w \big\}
\]
if $\alpha_i (w) \neq \alpha_i (g)$. Define $\mc W_i$-stable 
submanifolds
\begin{align*}
T_m &= \bigcup_{j :\; \dim L_j \leq m} L_j \\
T'_m &= T_m \cap T'
\end{align*}
and construct the following ideals
\begin{equation}
\begin{aligned}
A_m &= \{ h \in I_{i-1} / I_i : \pi(P_i ,t,\delta_i )(h) = 0 
  \;\mr{if}\; t \in T_m \} \\
B_m &= C^\infty (T', T'_m ; \mr{End}\: V_i )^{\mc W_i} \\
C_m &= C (T', T'_m ; \mr{End}\: V_i)^{\mc W_i} 
\end{aligned}
\end{equation}
Now we have $A_n = B_n = C_n = 0$ for $n \geq \dim T^{P_i}$ and
\[
A_n = I_{i-1} / I_i \quad 
B_n = C^\infty (U; \mr{End}\: V_i)^{\mc W_i} \quad 
C_n = C (U; \mr{End}\: V_i)^{\mc W_i} \quad \mr{for}\; n<0
\]
Just as in \eqref{eq:3.73} it suffices to show that the inclusions
\[
A_{m-1} / A_m \to B_{m-1} / B_m
\] 
induce isomorphisms on $HP_*$, so let us compute the periodic 
cyclic homologies of these quotient algebras.

Because $T_m$ is an algebraic subvariety of $T^{P_i}$ the spectrum 
of $A_{m-1} / A_m$ consists precisely of the 
irreducible representations of $I_{i-1} / I_i$ with tempered 
central character in $(P_i, T_m \setminus T_{m-1}, \delta_i)$.
We let $r_i(t)$ be the number of $\pi \in \mr{Prim}(I_{i-1} / I_i )$
corresponding to $(P_i,t,\delta_i)$. 
From Theorem \ref{thm:3.11} and we see that $r_i(t |t|^s) = r_i(t) 
\; \forall s > -1$, and from Theorem \ref{thm:3.19} that
$r_i(t |t|^{-1}) = r_i(t)$ if the stabilizers in $\mc W_i$ of $t$
and $t |t|^{-1}$ are equal. Choose a minimal subset  
$\{ L_{m,k} \}_k$ of $\mc L$ such that every $m$-dimensional element
of $\mc L$ is conjugate under $\mc W_i$ to a $L_{m,k}$. Let 
$\mc W_{m,k}$ be the stabilizer of $L_{m,k}$ in $\mc W_i$ and write
$r_k = r_i(t)$ for some $t \in L_{m,k} \setminus T_u^{P_i}$. By 
construction $\mc W_{m,k}$ acts freely on $L_{m,k} \setminus T_{m-1}$,
and the spectrum of $A_{m-1} / A_m$ is homeomorphic to
\begin{align*}
X_m \setminus Y_m  \;:=\; & \bigsqcup_k \bigsqcup_{l=1}^{r_k} 
\left( L_{m,k} \setminus T_{m-1} \right) / \mc W_{m,k} \\
\;=\; & \bigsqcup_k \bigsqcup_{l=1}^{r_k} \left( L_{m,k} / \mc W_{m,k} 
 \right) \setminus \big( (L_{m,k} \cap T_{m-1}) / \mc W_{m,k} \big)
\end{align*}
These are separable algebraic varieties, so the morphisms of finite
type algebras
\begin{equation}\label{eq:3.74}
A_{m-1} / A_m \leftarrow Z (A_{m-1} / A_m ) \to \mc O_0 (X_m ,Y_m )
\end{equation}
are spectrum preserving. Thus from Theorems \ref{thm:2.3} and 
\ref{thm:2.4} we get natural isomorphisms
\begin{equation}\label{eq:3.75}
HP_* (A_{m-1} / A_m ) \cong HP_* \big( \mc O_0 (X_m ,Y_m ) \big) \to
\check H^* (X_m,Y_m; \mh C)
\end{equation}
On the other hand, by \cite[Th\'eor\`eme IX.4.3]{Tou} the extension
\begin{equation*}
0 \to C^\infty (T', T'_m ; \mr{End}\: V_i) \to C^\infty 
(T'; \mr{End}\: V_i) \to C^\infty (T'_m ; \mr{End}\: V_i) \to 0
\end{equation*}
is admissible, and since $\mc W_i$ is finite the same holds for
\begin{equation*}
0 \to B_m \to C^\infty (T'; \mr{End}\: V_i)^{\mc W_i} \to 
C^\infty (T'_m ; \mr{End}\: V_i)^{\mc W_i} \to 0
\end{equation*}
So from Theorems \ref{thm:2.19} and \ref{thm:2.18} and 
Proposition \ref{prop:2.8} we get isomorphisms
\begin{equation}
HP_* (B_m ) \mosi K_* (B_m ) \otimes \mh C \isom 
K_*(C_m ) \otimes \mh C
\end{equation}
The spectrum of $C_{m-1}/C_m$ is 
\begin{align*}
X'_m \setminus Y'_m \;:=\; & \left( X_m \cap T' / \mc W_i \right) 
\setminus \left( Y_m \cap T' / \mc W_i \right) \\
\;=\; & \bigsqcup_k \bigsqcup_{l=1}^{r_k} (L_{m,k} \cap T'_m) / 
\mc W_{m,k} \setminus (L_{m,k} \cap T'_{m-1}) / \mc W_{m,k}
\end{align*}
These are locally compact Hausdorff spaces, so the $C^*$-algebra
homomorphisms
\begin{equation}\label{eq:3.83}
C_{m-1}/ C_m \leftarrow Z (C_{m-1}/C_m ) \isom C_0 (X'_m ,Y'_m )
\end{equation}
are spectrum preserving. By construction the stabilizer in 
$\mc W_i$ of $t \in T'$ is constant on the connected components of 
$T'_m \setminus T'_{m-1}$, so by the continuity of the intertwiners 
$\pi (g,P_i ,t,\delta_i )$ the vector space $C_{m-1} / C_m$ is a 
projective module over $C_0 (X'_m ,Y'_m )$. Thus by Proposition
\ref{prop:2.26} (for $K_* ( \cdot ) \otimes \mh Q$) \eqref{eq:3.83} 
induces isomorphisms on $K$-theory with rational coefficients. From 
this and Theorems \ref{thm:2.18} and \ref{thm:2.19} we obtain natural 
isomorphisms
\begin{equation}\label{eq:3.76}
\begin{split}
HP_* (B_{m-1}/B_m ) &\cong K_* (C_{m-1}/C_m ) \otimes \mh C \\
& \cong K_* (Z (C_{m-1}/C_m )) \otimes \mh C \\
& \cong K_* (C_0 (X'_m ,Y'_m )) \otimes \mh C \\
& \cong K_* (C_0^\infty (X'_m ,Y'_m )) \otimes \mh C \\
& \cong HP_* (C_0^\infty (X'_m ,Y'_m )) \\
& \cong \check H^* (X'_m ,Y'_m ; \mh C)
\end{split}
\end{equation}
From \eqref{eq:3.75} - \eqref{eq:3.76} we construct the commutative 
diagram
\begin{equation}
\begin{array}{ccccc}
HP_* (A_{m-1} / A_m ) & \cong & HP_* \big( \mc O_0 (X_m ,Y_m ) \big) 
 & \to & \check H^* (X_m,Y_m;\mh C) \\
\downarrow & & \downarrow & & \downarrow \\
HP_* (B_{m-1} / B_m ) & \cong & HP_* \big( C_0^\infty 
 (X'_m ,Y'_m ) \big) & \to & \check H^* (X'_m , Y'_m ; \mh C)
\end{array}
\end{equation}
The pair $(X'_m, Y'_m)$ is a deformation retract of $(X_m ,Y_m )$, 
so all maps in this diagram are isomorphisms. 
Working our way back up, using excision, we find that also
\[
HP_* (I_{i-1}/I_i ) \to HP_* \left( C^\infty (T'; \mr{End}\: 
V_i )^{\mc W_i} \right) \to HP_* (J_{i-1}/J_i)
\]
and finally 
\[
HP_* (\mc H (\mc R ,q)) \to HP_*(\mc S (\mc R ,q))
\] 
are isomorphisms. $\qquad \Box$
\\[1mm]

Note that Theorem \ref{thm:3.38} is in accordance with our earlier
results for direct products of root data. If $\mc R = \mc R_1 
\times \mc R_2$ then by \eqref{eq:3.46} and Theorem 
\ref{thm:2.28} we have
\begin{equation}
HP_* (\mc H (\mc R ,q)) \cong HP_* (\mc H (\mc R_1 ,q)) \otimes 
HP_* (\mc H (\mc R_2 ,q))
\end{equation}
while by \eqref{eq:3.47}, Theorem \ref{thm:3.19} and Corollary
\ref{cor:2.29}
\begin{equation}
HP_* (\mc S (\mc R ,q)) \cong HP_* (\mc S( \mc R_1 ,q)) \otimes 
HP_* (\mc S (\mc R_2 ,q))
\end{equation}
We can pursue the path of \eqref{eq:3.68} and Theorem \ref{thm:3.38}
a little further. Let $\mb k$ be any (unital) subring of $\mh C$
containing $\{ q(w) : w \in W\}$. As on page \pageref{p:iha} we
consider the extended Iwahori-Hecke algebra $\mc H_{\mb k}(\mc R ,q)$.
It makes sense to take its periodic cyclic homology in the category 
of $\mb k$-algebras. This is a $\mb k$-module which we denote by 
$HP_* (\mc H_{\mb k} (\mc R ,q) | \mb k)$.

\begin{thm}\label{thm:3.41}
There are natural isomorphisms
\begin{multline*}
HP_* (\mc H_{\mb k} (\mc R ,q) | \mb k ) \otimes_{\mb k} \mh C 
\;\isom\; HP_* (\mc H (\mc R ,q)) \;\isom\; 
HP_* (\mc S (\mc R ,q)) \;\mosi 
\\ K_* (\mc S (\mc R ,q)) \otimes_{\mh Z} \mh C 
\;\isom\; K_* (C_r^* (\mc R ,q)) \otimes_{\mh Z} \mh C
\end{multline*}
\end{thm}
\emph{Proof.} 
By \cite[Proposition IX.5.1]{CaEi} we have
\begin{equation}
HH_n (\mc H_{\mb k} (\mc R ,q) | \mb k ) \otimes_{\mb k} \mh C
\cong HH_n (\mc H (\mc R ,q))
\end{equation}
where $HH_n ( \,\cdot\, | \mb k )$ means Hochschild homology in the
category of $\mb k$-algebras. Now the first isomorphism follows
from \cite[Proposition 5.1.6]{Lod}. The second isomorphism is
Theorem \ref{thm:3.38} and the third was already noticed in 
\eqref{eq:3.68}. Finally, the fourth isomorphism is a consequence 
of Theorem \ref{thm:2.19}. $\qquad \Box$
\\[1mm]

Apparently this is an important invariant of the labelled root datum
$(\mc R ,q)$. By the way, we really need complex coefficients. It 
does not follow from Theorem \ref{thm:3.41} that $HP_* (\mc H_{\mb k} 
(\mc R ,q) | \mb k )$ and $K_* (C_r^* (\mc R ,q)) \otimes_{\mh Z} 
\mb k$ are naturally isomorphic, we merely know that they have the 
same (finite) rank as $\mb k$-modules. In general there is no reason 
why the image of a class in $K_* (C_r^* (\mc R ,q))$ should land in 
$HP_* (\mc H_{\mb k} (\mc R ,q) | \mb k )$ under the composition of
the above isomorphisms.

%% file: chapter4.tex
\chapter{Reductive $p$-adic groups}

Iwahori and Matsumoto were the first to recognize that any finite 
group with a $BN$-pair gives rise to a finite dimensional
Hecke algebra with a very nice description in terms of generators 
and relations. In particular this applies to a connected reductive
algebraic group defined over a finite field. 

On a higher level, if $G$ is a reductive algebraic group over a
non-Archimedean local field, then the Hecke algebra $\mc H (G)$ has 
infinite dimension. However, the valution of the field provides
enough extra structure to show that $\mc H (G)$ is a direct sum of
factors which tend to be Morita equivalent to affine Hecke algebras.

In the the first section we browse through the literature on reductive
groups, and we report when and how we see something that looks like
an affine Hecke algebra. Meanwhile we also recall some important
notions from the representation theory of totally disconnected 
groups (like groups over a $p$-adic field).

In Section \ref{sec:4.2} we start working towards the main new
result of this chapter, namely the construction of natural
isomorphisms
\begin{equation}\label{eq:4.1}
HP_* (\mc H (G)) \cong HP_* (\mc S (G), \oot) \cong K_* (C_r^* (G))
\end{equation}
Obviously we have to recall the definitions of the
involved algebras. The reduced $C^*$-algebra of $G$ is defined in a
standard way, but the construction of the Schwartz algebra 
$\mc S (G)$, originally due to Harish-Chandra \cite{HC1},
is much more difficult. The greater part of Section \ref{sec:4.2}
is used to give a proper definition of this algebra, and to 
characterize its representations among all $G$-representations.

Acknowledging the difference between $\mc H (G)$ and affine Hecke
algebras, we still proceed like we did in Chapter 3. Thus for 
information about the primitive ideal spectra of $\mc H (G)$ and
$\mc S (G)$ we turn to the Fourier transform and the Plancherel 
theorem for reductive $p$-adic groups, both of which are due to
Harish-Chandra \cite{HC2}. These will show that Prim$(\mc H (G))$
is a countable union of non-separated complex affine varieties. 
The algebras $\mc S (G)$ and $C_r^* (G)$ have the same spectrum,
which turns out to be a countable union of compact non-Hausdorff
spaces. The Langlands classification tells us that 
Prim$(\mc S (G))$ is in a sense a deformation retract of 
Prim$(\mc H (G))$.

Nearly all the new material of this chapter is contained
in the final section. There we use all the above to lift the 
result \eqref{eq:3.50} to Hecke algebras of reductive $p$-adic 
groups, which yields \eqref{eq:4.1}. We remark that we are 
careful with topological periodic cyclic homology, here we 
have to take it with respect to the completed inductive tensor 
product $\oot$.

Moreover \eqref{eq:4.1} is related to the Baum-Connes conjecture
for reductive $p$-adic groups by means of the diagram
\[
\begin{array}{ccc}
K_*^G (\beta G ) & \to & K_* (C_r^* (G)) \\
\downarrow & & \downarrow \\
HP_* (\mc H (G)) & \to & HP_* (\mc S (G))
\end{array}
\]
We conclude the chapter with a discussion of some subtleties
of this diagram.
\\[4mm]

\section{Hecke algebras of reductive groups}
\label{sec:4.1}

In Section \ref{sec:3.1} we defined Iwahori-Hecke algebras in terms 
of generators and relations, but this is hardly the way in which 
they emerged. Iwahori and Matsumoto \cite{Iwa1,Iwa2,IwMa,Mat2} 
discovered that convolution algebras associated to a reductive group 
and a suitable subgroup are of the type we described. We have a look 
at these and then we extend our view to more general convolution
algebras, of smooth functions on reductive $p$-adic groups. We 
recall everything that is needed to state all the known cases in
which such convolution algebras yield affine Hecke algebras or
closely related structures.

The most direct way to arrive at Iwahori-Hecke algebras is through
groups with a $BN$-pair. Recall that a group $G$ has a 
$BN$-pair if it satisfies \index{$BN$-pair}
\begin{enumerate}
\item $G$ is generated by two subgroups $B$ and $N$
\item $B \cap N$ is normal in $N$
\item $W := N / B \cap N$ is generated by a set 
      $S = \{ s_i : i \in I \}$ of elements of order 2
\item if $n_i \in N$ and $n_i (B \cap N) = s_i$ then 
      $n_i B n_i \neq B$
\item \label{enum:4.v} for all such $n_i$ and $n \in N$ we have
      $n_i B n \subset B n_i n B \cap B n B$
\end{enumerate}

These axioms were first formalized by Tits, cf. \cite{BrTi1}. Some
important consequences are proven in \cite[Chapitre IV.2]{Bou}.
For example, it turns out that $(W,S)$ is a Coxeter system and that 
the group $G$ has a \inde{Bruhat decomposition}, i.e. there is a 
bijection between $W$ and the double cosets of $B$ in $G$, given by
\begin{equation}\label{eq:4.16}
N \ni n \to B n B \in B \setminus G / B
\end{equation}

Any connected reductive group $\mc G$ over an algebraically closed 
field $\mh K$ has a $BN$-pair. To be precise, in this case $\mc B$ 
is a Borel subgroup of $\mc G ,\, \mc B \cap \mc N = \mc T$ is a 
maximal torus and $\mc N$ is the normalizer of $\mc T$ in $\mc G$, 
see \cite[Chapter 8]{Spr2}. With $\mc G$ and $\mc T$ one can 
associate a root datum 
\[
\mc R (\mc G ,\mc T ) = (X,Y,R_0 ,R_0^\vee ) \quad \mr{where} 
\]
\begin{itemize}
\item $X$ is the character lattice of $\mc T$
\item $Y$ is the cocharacter lattice of $\mc T$
\item $R_0$ is the set of roots of $(\mc G ,\mc T )$
\item $R_0^\vee$ is the set of coroots of $(\mc G ,\mc T )$
\item $W = \mc N / \mc T$ is isomorphic to the Weyl group 
$W_0$ of $R_0$
\end{itemize}
This results in a bijection between isomorphism classes of 
connected reductive algebraic groups and root data, see 
\cite[Expos\'e 24]{Che} and \cite[Expos\'e XXV]{DeGr}. Under this
bijection semisimple groups correspond to semisimple root data.

The most important example is of course $GL (n, \mh K)$. In 
this group we may take for $\mc B$ the subgroup of upper triangular 
matrices and for $\mc T$ the subgroup of diagonal matrices. Then 
$\mc N$ consists of the matrices that have exactly one nonzero 
entry in every row and every column. In the root datum 
$\mc R (\mc G ,\mc T )$ we have $X \cong Y \cong \mh Z^n$ and
\[
R_0 = R_0^\vee = \{ e_i - e_j : 1 \leq i,j \leq n ,\, i \neq j \}
\]
where $\{ e_i \}_{i=1}^n$ is the standard basis of $\mh Z^n$.
So $R_0$ is the root system $A_{n-1}$, and the Weyl group $W_0$ 
is isomorphic to the symmetric group $S_n$.
\\[2mm]

Assume now that we have a group with a $BN$-pair such that $B$ is
\inde{almost normal} in $G$. This means that every double coset of 
$B$ is a finite union of left cosets. (Clearly $B$ would be almost 
normal if it were normal in $G$, but this can only happen in the
degenerate situation $B = G ,\, N = \{ e \}$.)
Let $\mb k$ be a unital commutative ring and consider the 
$\mb k$-module $\mc H (G,B)$ of $\mb k$-valued $B$-biinvariant 
functions on $G$ that are nonzero on only finitely many left 
$B$-cosets. Clearly this is a free module with basis 
$\{ T_w : w \in W \}$, where $T_w$ is the characteristic function 
of $B w B$. (This sloppy notation is justified by \eqref{eq:4.16}.)
Define a measure $\mu$ on $B$-left invariant subsets by 
$\mu (H) = | B \setminus H |$. The product on $\mc H (G,B)$ is 
\inde{convolution} with respect to $\mu$:
\begin{equation}
(f_1 * f_2 ) (w) = \int_{B \setminus G} f_1 (w x^{-1}) f_2 (x) 
d\mu (x)
\end{equation}
This notion of a Hecke algebra stems from Shimura \cite[\S 7]{Shi}.
Fortunately it agrees with the definition of an Iwahori-Hecke 
algebra as given on page \pageref{p:iha}:

\begin{thm}\label{thm:4.7}
For $w \in W$ write 
\[
q(w) = \mu (B w B) = | B \setminus B w B |
\]
Then $q$ is a label function and the relations \eqref{eq:3.2} and 
\eqref{eq:3.3} hold in $\mc H (G,B)$.
\end{thm}
\emph{Proof.} 
This result is due to Iwahori \cite[Theorem 3.2]{Iwa1} and 
Matsumoto \cite[Th\'eor\`eme 4]{Mat2}. A full proof can be found in
\cite[Theorem 8.4.6]{GePf}. $\qquad \Box$
\\[2mm]

In the above situation of a reductive group $\mc G$ over an 
algebraically closed field $\mh K$, $\mc B$ cannot be almost 
normal, because it has lower dimension than $\mc G$ and $\mh K$ 
is infinite. However, suppose that the characteristic $p$ of 
$\mh K$ is nonzero and that $\mc G$ is defined over a finite 
field $\mh F_q$. The group $\mc G (\mh F_q)$ of 
$\mh F_q$-rational points still has a $BN$-pair, where 
$\mc B (\mh F_q)$ is a Borel subgroup. The associated Hecke 
algebras were studied by Iwahori \cite{Iwa1} (for Chevalley
groups) and by Howlett and Lehrer, see \cite{Car1}.
In these cases $W$ is a subgroup of the finite Weyl group 
$W_0$ of $\mc R (\mc G ,\mc T)$, and the numbers $q(s)$ are 
certain powers of $p$.
\\[2mm] 

Now we turn to $p$-adic groups. Recall that a \inde{non-Archimedean
local field} is a field $\mh F$ with a discrete valuation
\begin{equation}
v : \mh F \to \mh Z \cup \{ \infty \}
\end{equation}
such that $\mh F$ is complete with respect to the induced norm
$\norm{x}_{\mh F} = q^{-v (x)}$. \index{$v$}
Here $q$ is the cardinality of the residue field $\mc O / \mf P$,
\begin{equation}
\mc O = \{ x \in \mh F : v(x) \geq 0 \}
\end{equation}
being the ring of integers of $\mh F$ and 
\index{Omc@$\mc O$} \index{Pmf@$\mf P$} \index{Fmh@$\mh F$}
\begin{equation}
\mf P = \{ x \in \mh F : v(x) > 0 \}
\end{equation}
its unique maximal ideal. This implies that $\mh F$ is a 
totally disconnected, nondiscrete, locally compact Hausdorff 
space. If its characteristic is zero then $\mh F$ is isomorphic 
to a finite algebraic extension of the field of $p$-adic numbers 
$\mh Q_p$. On the other hand, if char$(\mh F ) > 0$ then 
$\mc O \cong \mh F_q [[t]]$, the ring of formal power series over 
the finite field $\mh F_q$. With a slight abuse of terminology, 
non-Archimedean local fields are also known as $p$-adic fields.
\index{$p$-adic field}

So let $\mh F$ be a non-Archimedean local field and $\mc G$ a 
connected reductive algebraic group that is defined over $\mh F$.
Consider the group $G = \mc G (\mh F)$ of $\mh F$-rational points
of $\mc G$. Affine Hecke algebras play an important role in the
representation theory of such groups, as we will try to explain.
We are mainly interested in smooth representations, i.e.
representations of $G$ on a complex vector space $V$ such that for
every $v \in V$ the group $\{ g \in G : g v = v \}$ is open.
Recall that, because $\mh F$ is non-Archimedean, the identity 
element $e$ of $G$ has a countable neighborhood basis consisting of 
compact open subgroups. Hence smooth representations can also 
be characterized by the condition \index{representation!smooth}
\[
V = \bigcup_K V^K
\]
where $K$ runs over all compact open subgroups of $G$.
We denote the category of smooth $G$-representations by 
\inde{Rep$(G)$}, and the set of equivalence classes of irreducible 
smooth $G$-representations by \inde{Irr$(G)$}. We call a 
map from $G$ to a Hausdorff space smooth if it is uniformly locally 
constant, i.e. if it is bi-invariant for some compact open subgroup
of $G$. \index{smooth map}

Furthermore, because $G$ is reductive, it has an
\inde{affine building} $\beta G$, also known as the 
\inde{Bruhat-Tits building} of $G$. We quickly recall the 
construction and terminology of this polysimplicial complex,
referring to \cite{BrTi1,Tit} for more detailed information.
Let \inde{$A_0$} $= \mc A_0 (\mh F)$ be a maximal $\mh F$-split 
torus of $\mc G$, $X^* (A_0 ) = X^* (\mc A_0 )$ its character 
lattice, and put \inde{$\mf a_0^*$} $= X^* (A_0 ) \otimes_{\mh Z} 
\mh R$. Then $\beta G$ is $G \times \mf a_0^*$ modulo a certain
equivalence relation. \index{betaG@$\beta G$}

The affine building is a universal space for proper $G$-actions. 
Such a universal space always exists, based on general 
categorical considerations, but it is unique only up to homotopy. 
On the other hand, the proof that $\beta G$ really has the 
required properties is very tricky, and ultimately relies on the 
existence of ``valuated root data'' \cite{BrTi2}.

The images of $\mf a_0^*$ under $G$ are the appartments of 
$\beta G$, and a polysimplex of maximal dimension in $\beta G$ 
is called a chamber.\index{appartment}\index{chamber} 
The stabilizer $I$ of such a chamber (or equivalently
of an interior point of a chamber) is an \inde{Iwahori subgroup} 
of $G$. More generally the stabilizer of an arbitrary point 
of $\beta G$ is known as a \inde{parahoric subgroup}.
If $x_0 \in \beta G$ is a ``special'' point (in particular it must 
lie in a polysimplex of minimal dimension) then its stabilizer
\inde{$K_0$} is a "good" maximal compact subgroup of $G$ in the 
sense of \cite[\S 0.6]{Sil1}. This implies that 
\begin{equation}\label{eq:4.24}
G = P K_0 = K_0 P
\end{equation} 
for any parabolic subgroup $P$ of $G$ containing $A_0$.
\\[2mm]

Normalize the Haar measure $\mu$ on $G$ by $\mu (K_0 ) = 1$.
Any compact open $K < G$ is almost normal, so we can consider the 
convolution algebra \inde{$\mc H (G,K)$} of compactly supported 
$\mh C$-valued $K$-biinvariant functions on $G$. 
For example, if $G = GL_2 (\mh Q_p )$ and $K = GL_2 (\mh Z_p )$, 
then $\mc H (G,K)$ is the classical algebra of Hecke operators, 
hence the "$\mc H$" for "Hecke" algebra.

If $K' \subset K$ is another compact subgroup then there is a 
natural inclusion $\mc H (G,K) \to \mc H (G,K')$. The inductive limit 
of this system of inclusions (i.e. the union), over all compact open 
subgroups, is called the \inde{Hecke algebra} \inde{$\mc H (G)$} of 
$G$. It consists of all compactly supported smooth functions on $G$.
For every $K$ we define the idempotent \inde{$e_K$}$\in \mc H (G)$ by
\begin{equation}
e_K (g) = \left\{ \begin{array}{lcc}
\mu (K)^{-1} & \mr{if} & g \in K \\
0 & \mr{if} & g \notin K
\end{array} \right.
\end{equation}
This gives the useful identification
\begin{equation}\label{eq:4.42}
e_K \mc H (G) e_K = \mc H (G)^{K \times K} = \mc H (G,K)
\end{equation}
By construction a smooth $G$-representation is the same thing as
a nondegenerate representation of $\mc H (G)$. There are natural maps
\begin{equation}
\begin{array}{lllllll}
\mr{Rep}(G) & \to & \mr{Rep}\big( \mc H (G,K) \big) & : & 
V & \to & V^K \; = \; \pi (e_K ) V \\
\mr{Rep}\big( \mc H (G,K) \big) & \to & \mr{Rep}(G) & : &
W & \to & \mr{Ind}^{\mc H (G)}_{\mc H (G,K)} W 
%& = & \mc H (G) \otimes_{\mc H (G,K)} W
\end{array}
\end{equation}
Bernstein \cite[Corollaire 3.9]{BeDe} showed that there exist
arbitrarily small $K$ for which these maps define equivalences 
between the category of nondegenerate $\mc H (G,K)$-representations 
and the category of those smooth $G$-representations that are 
generated by their $K$-fixed vectors. Thus $\mc H (G,K)$ covers a 
clear part of the representation theory of $G$.

It is quite possible that $\mc H (G,K)$ is an extended 
Iwahori-Hecke algebra. For example, suppose that $\mc G$ is split 
over $\mh F$, let $A_0 = \mc A_0 (\mh F )$ be a maximal split torus 
and $B = \mc B (\mh F )$ a Borel subgroup containing $\mc T$. 
Assume that $\mc B$ is defined over $\mc O$, so that 
$\mc B (\mc O / \mf P )$ is a Borel subgroup of 
$\mc G (\mc O / \mf P )$. The inverse image $I$ of
$\mc B (\mc O / \mf P )$ under the quotient map 
$\mc G (\mc O ) \to \mc G (\mc O / \mf P )$ is an Iwahori subgroup
of $G$. It is known \cite{IwMa} that $\mc H (G,I)$ is an 
affine Hecke algebra. The Weyl group of the associated root datum 
is $W = N_G (A_0 ) / \mc A_0 (\mc O )$ and it decomposes as 
$W = W_0 \ltimes X^* (A_0 )$ where 
$W_0 = N_{\mc G(\mc O )} (A_0 ) / \mc A_0 (\mc O )$. Finally, 
the value of the label function on any simple reflection is 
$q = |\mc O / \mf P |$.

Or suppose that $\mc G$ is simply connected, but not necessarily 
split over $\mh F$. Let $N$ be the stabilizer of an appartment in 
the affine building of $G$, and $I$ the stabilizer of a chamber 
of this appartment. According to \cite[Proposition 5.2.10]{BrTi2} 
$(I,N)$ is a $BN$-pair in $G$, so $\mc H (G,I)$ is an Iwahori-Hecke 
algebra. The Coxeter group $W = N / I \cap N$ is an affine
Weyl group coming from a root datum that is contained in
$\mc R (\mc G ,\mc T )$, for a suitable torus $\mc T$.
\\[2mm]

An important decomposition of the category Rep$(G)$ was 
discovered by Bernstein \cite{BeDe}. To describe it we introduce 
several classes of smooth representations.

We call $(\pi ,V) \in \mr{Rep}(G)$
\begin{itemize}
\item admissible if $V^K$ has finite dimension for every compact
      open subgroup $K$ \index{representation!admissible}
\item supercuspidal if it is admissible and all matrix coefficients
      of $\pi$ have compact support modulo the center of $G$
\index{representation!supercuspidal}
\end{itemize}
By \cite[Corollaire 3.4]{BeDe} for every $K$ the algebra $\mc H 
(G,K)$ is of finite type, so all its irreducible representations 
have finite dimension. In combination with \cite[Corollaire 3.9
]{BeDe} this shows that every irreducible smooth
representation is automatically admissible. Supercuspidal
representations are also known as absolutely cuspidal (or just
cuspidal) representations, but there seems to be no agreement in 
the terminology here.

There is a natural notion of the contragredient of a smooth 
representation. Let \index{representation!contragredient}
\begin{equation}\label{eq:4.40}
\breve{V}^K = \{ f \in V^* : f \circ \pi (e_K ) = f \}
\end{equation}
be the dual space of $V^K$ and define 
\begin{equation}
\breve V = \bigcup_K \breve{V}^K
\end{equation}
Then $(\breve \pi ,\breve V )$ is the contragredient representation 
of $(\pi ,V)$. By construction it is smooth, and it is admissible 
whenever $V$ is. \index{piv@$(\breve \pi ,\breve V )$}

Suppose that $P$ is a parabolic subgroup of $G$ and that $M$ is a
Levi subgroup of $P$. Although $G$ and $M$ are unimodular, the 
modular function $\delta_P$ of $P$ is in general not constant. 
Let $\sigma$ be an irreducible supercuspidal representation of $M$.
Under these conditions we call $(M,\sigma )$ a \inde{cuspidal 
pair}. From this we construct a parabolically induced 
$G$-representation \index{$I_P^G (\sigma)$}
\index{representation!parabolically induced}
\begin{equation}
I_P^G (\sigma) = 
\mr{Ind}_P^G \big( \delta_P^{1/2} \otimes \sigma \big)
\end{equation}
This means that we first inflate $\sigma$ to $P$, and then we 
apply the \inde{normalized induction} functor, i.e. we twist it by
$\delta_P^{1/2}$ and take the smooth induction to $G$. This twist 
is useful to preserve unitarity, cf. \cite[Theorem 3.2]{Car2}.

For every $(\pi ,V) \in \mr{Irr}(G)$ there is a cuspidal pair
$(M, \sigma)$, uniquely determined up to $G$-conjugacy, such that
$V$ is a subquotient of $I_P^G (\sigma)$. If $P'$ is another 
parabolic subgroup of $G$ containing $M$ then $I_{P'}^G (\sigma)$ 
and $I_P^G (\sigma)$ have the same irreducible subquotients, but 
they need not be isomorphic.

To define a suitable equivalence relation on the set of cuspidal 
pairs, we now recall a particular kind of characters. Let $H$ be 
any algebraic group over $\mh F$, and consider the subgroup
\begin{equation}
\prefix^{0}{H} = \left\{ h \in H : 
v (\gamma (h)) = 0 \; \forall \gamma \in X^* (H) \right\}
\end{equation}
This is an open normal subgroup of $G$ which contains every compact
subgroup, and $H / \prefix^{0}{H}$ is a free abelian group.
An \inde{unramified character} of $H$ is a homomorphism 
$\chi : H \to \mh C^\times$ whose kernel contains $\prefix^{0}{H}$. 
The group of these forms a complex torus \inde{$X_{nr}(H)$} and the 
map $X^* (H) \otimes_{\mh Z} \mh C^\times \to X_{nr}(H)$ defined by
\[
\gamma \otimes z \;\longrightarrow\; 
\left( h \to z^{\ds v(\gamma (h))} \right)
\]
is an isomorphism. We will denote the compact torus of unitary
unramified characters by \index{$X_{unr}(H)$}
\begin{equation}
X_{unr}(H) = \mr{Hom}\big( H / \prefix^{0}{H}, S^1 \big)
\end{equation}
We say that two cuspidal pairs $(M ,\sigma )$ and $(M' ,\sigma' )$
are inertially equivalent \index{inertial equivalence} 
if there exist $\chi \in X_{nr}(M')$ and $g \in G$ such that
$M' = g^{-1} M g$ and $\sigma' \otimes \chi \cong \sigma^g$.

With an inertial equivalence class $\mf s = [M,\sigma ]_G$ we
associate a subcategory of \inde{Rep$(G)^{\mf s}$} of Rep$(G)$. By 
definition its objects are those smooth representations $\pi$ 
with the following property: for every irreducible subquotient 
$\rho$ of $\pi$ there is a $(M,\sigma ) \in \mf s$ such that 
$\rho$ is a subrepresentation of $I_P^G (\sigma )$. 

These Rep$(G)^{\mf s}$ give rise to 
the \inde{Bernstein decomposition} \cite[Proposition 2.10]{BeDe}:
\begin{equation}\label{eq:4.18}
\mr{Rep}(G) = \prod_{\mf s \in \mf B (G)} \mr{Rep}(G)^{\mf s}
\end{equation}
The set \inde{$\mf B (G)$} of Bernstein components is countably 
infinite. We have a corresponding decomposition of the Hecke 
algebra of $G$ into two-sided ideals: \index{$\mc H (G)^{\mf s}$}
\begin{equation}\label{eq:4.19}
\mc H (G) = \bigoplus_{\mf s \in \mf B (G)} \mc H (G)^{\mf s}
\end{equation}
with Rep$\big( \mc H (G)^{\mf s} \big) = \mr{Rep}(G)^{\mf s}$.
By \cite[Proposition 3.3]{BuKu2} there exists an idempotent 
$e_{\mf s} \in \mc H (G)$ such that \index{$e_{\mf s}$}
\begin{itemize}
\item $\mc H (G)^{\mf s} = \mc H (G) e_{\mf s} \mc H (G)$
\item Rep$(G)^{\mf s}$ is equivalent to 
      Rep$\big( e_{\mf s} \mc H (G) e_{\mf s} \big)$
\end{itemize}
Under these conditions $e_{\mf s} \mc H (G) e_{\mf s}$ is a
finite type algebra, whose center was already described by 
Bernstein. The set \inde{$D_\sigma$} of all cuspidal pairs
of the form $(M, \sigma \otimes \chi)$ is in bijection with
$X_{nr}(M)$, so it has the structure of a complex torus. Put
\begin{equation}
N (M, \sigma) = \left\{ g \in G : g M g^{-1} = M \;\mr{and}\;
[M, \sigma^g ]_M = [M, \sigma ]_M \right\}
\end{equation}
Then \inde{$W_\sigma$}$= N (M, \sigma) / M$ is a finite group 
acting on $D_\sigma$, so $D_\sigma / W_\sigma$ is an irreducible 
algebraic variety. \cite[Th\'eor\`eme 2.13]{BeDe} tells us that
\begin{equation}\label{eq:4.20}
Z \big( e_{\mf s} \mc H (G) e_{\mf s} \big) \cong 
\mc O ( D_\sigma / W_\sigma ) = \mc O (D_\sigma )^{W_\sigma}
\end{equation}
For any compact open $K < G$ write 
\index{$\mf B (G,K)$} \index{$\mc H (G,K)^{\mf s}$}
\begin{align*}
& \mf B (G,K) = \big\{ \mf s \in \mf B (G) : 
\mc H (G,K)^{\mf s} \neq 0 \big\} \\
& \mc H (G,K)^{\mf s} = \mc H (G)^{\mf s} \cap \mc H (G,K)
\end{align*}

\begin{prop}\label{prop:4.42}
\begin{enumerate}
\item $\mf B (G,K)$ is finite for any compact open $K < G$.
\item If $K$ is a normal subgroup of a good maximal compact subgroup 
$K_0$ then $\mc H (G,K)$ is Morita equivalent
to $\oplus_{\mf s \in \mf B (G,K)} \mc H (G)^{\mf s}$ and
\[
Z \big( \mc H (G,K) \big) \cong \bigoplus_{\mf s \in \mf B (G,K)} 
Z \big( e_{\mf s} \mc H (G) e_{\mf s}\big)
\]
\item For any $\mf s \in \mf B (G)$ there exists a compact open
$K_{\mf s} < G$ such that for every compact open 
$K \subset K_{\mf s}$ the algebras
\[
\mc H (G,K_{\mf s})^{\mf s} \;,\; \mc H (G,K)^{\mf s} \text{ and }
\mc H (G)^{\mf s}
\]
are Morita equivalent.
\end{enumerate}
\end{prop}
\emph{Proof.}
All these results are due to Bernstein. 2 is a direct consequence of 
\cite[Corollaire 3.9]{BeDe} and \eqref{eq:4.20}. 1 and 3 follow from 
this in combination with \cite[(3.7)]{BeDe}, as was remarked in
\cite[p. 143]{BHP2}. $\qquad \Box$
\\[2mm]

We may also consider more general algebras associated to $(G,K)$.
Let $(\rho, V)$ be an irreducible smooth representation of $K$, 
and $(\breve \rho ,\breve V )$ its contragredient. Notice that 
$\rho$ is smooth and has finite dimension. Define 
\index{$\mc H (G,K,\rho )$}
\begin{equation}
\mc H (G,K,\rho ) = \left\{ f : G \to \mr{End}_{\mh C} (\breve V ) 
: f (k_1 g k_2 ) = \breve \rho (k_1) f(g) \breve \rho (k_2)
\;\forall k_1 ,k_2 \in K, g \in G \right\}
\end{equation}
This is a unital algebra under the convolution product, its
elements being smooth functions on $G$. Consider the idempotent
$e_\rho \in \mc H (G)$ defined by \index{$e_\rho$}
\begin{equation}
e_\rho (g) = \left\{ \begin{array}{lcc}
\mu (K)^{-1} \dim (V) \,\mr{tr} \big( \rho (g^{-1}) \big) & 
 \mr{if} & g \in K \\
0 & \mr{if} & g \notin K
\end{array} \right.
\end{equation}
By \cite[Proposition 4.2.4]{BuKu1} there is a natural isomorphism
\begin{equation}\label{eq:4.21}
\mc H (G,K,\rho ) \otimes_{\mh C} \mr{End}_{\mh C} (V) \cong 
e_\rho \mc H (G) e_\rho
\end{equation}
Let \inde{$\mr{Rep}_\rho (G)$} be the subcategory of Rep$(G)$ 
consisting of all representations $(\pi ,U)$ for which 
$\mc H (G) e_\rho U = U$. According to 
\cite[Proposition 4.2.3]{BuKu1} there are natural 
bijections between the sets of irreducible objects of
\begin{itemize}
\item $\mr{Rep}_\rho (G)$
\item Rep$\big( e_\rho \mc H (G) e_\rho \big)$
\item Rep$\big( \mc H (G,K,\rho ) \big)$
\end{itemize}
If moreover $\mr{Rep}_\rho (G)$ is closed under taking subquotients
(of $G$-representations) then there exists a finite subset
$\mf S \subset \mf B (G)$ such that 
\begin{equation}
\mr{Rep}_\rho (G) = \bigoplus_{\mf s \in \mf S} \mr{Rep}(G)^{\mf s}
\end{equation}
In the terminology of Bushnell and Kutzko $(K,\rho)$ is an 
$\mf S$-type \cite[(3.12)]{BuKu2}. Of special interest is the 
case when $\mf S$ consists of a single element $\mf s$, for then
$\mc H (G,K, \rho )$ is Morita equivalent to $\mc H (G)^{\mf s}$.
It is known that under this and certain extra conditions
$\mc H (G,K, \rho )$ is isomorphic to an affine Hecke algebra
\cite[Theorem 6.3]{Roc}, sometimes with unequal labels
\cite[Section 1]{Lus5}, or to the twisted crossed product of
such a thing with a (twisted) group algebra 
\cite[Theorem 7.12]{Mor}. \index{type}

Using this approach it has been shown that $\mc H (G)^{\mf s}$ is
Morita equivalent to an affine Hecke algebra for every
$\mf s \in \mf B (GL (n,\mh F ))$ \cite{BuKu1,BuKu3}, and to a
``twisted'' affine Hecke algebra for every 
$\mf s \in \mf B (SL (n,\mh F ))$ \cite[Theorem 11.1]{GoRo}.
\\[4mm]

\section{Harish-Chandra's Schwartz algebra}
\label{sec:4.2}

In this section $G$ will be a connected reductive algebraic group
defined over a non-Archimedean local field $\mh F$. The reduced
$C^*$-algebra of $G$ is defined in a standard way, but we need to
go to some lengths to construct Harish-Chandra's Schwartz algebra.
Once this is done we characterize its representations among 
admissible $G$-representations and formulate the Langlands 
classification for reductive $p$-adic groups.

Define the adjoint and the trace of $f \in \mc H (G)$ by 
\index{$f^*$} \index{tau@$\tau$}
\[
f^* (g) = \overline{f(g^{-1})} \qquad \tau (f) = f(e)
\]
This gives rise to a bitrace \index{$(f,f')$}
\[
(f,f') = \tau (f^* f')
\]
making $\mc H (G)$ into a Hilbert algebra. The Hilbert space 
completion of $\mc H (G)$ is the space \inde{$L^2 (G)$} of all
square-integrable functions on $G$. It carries two natural 
$G$-actions, left and right translation:
\index{lambda@$\lambda$} \index{rho@$\rho$}
\begin{align*}
\big( \lambda (g) f \big) (h) &= f (g^{-1} h)\\
\big( \rho (g) f \big ) (h) &= f (h g)
\end{align*}
Now $\lambda (g)$ and $\rho (g)$ are bounded operators of the same 
norm. They extend naturally to representations of $\mc H (G)$, so 
we get an injection
\[
\lambda : \mc H (G) \to B (L^2 (G))
\]
The \inde{reduced $C^*$-algebra} \inde{$C_r^* (G)$} is the closure 
of $\lambda (\mc H (G))$ in $B (L^2 (G))$. It is a separable 
nonunital $C^*$-algebra whose representations correspond to the 
unitary $G$-representations that are weakly contained in the left
regular representation $(\lambda ,L^2 (G))$ of $G$. 

Usually the reduced $C^*$-algebra of a locally compact group $H$ 
is defined as the completion of $C_c (H)$ or $L^1 (H)$ in
$B (L^2 (H))$. However if $H$ is totally disconnected we may just 
as well start with smooth functions only.

For a compact open $K < G$ we let \inde{$C_r^* (G,K)$} 
be the completion of $\mc H (G,K)$ in $B (L^2 (G))$. It is a 
unital type I $C^*$-algebra and it equals
\begin{equation}\label{eq:4.78}
e_K C_r^* (G) e_K = C_r^* (G)^{K \times K} = C_r^* (G,K)
\end{equation}
Let us mention some general facts about the structure of 
$C_r^* (G)$. They can be read off from Theorem \ref{thm:4.42},
but it seems appropriate to formulate them here already.
This algebra can be recovered as the inductive limit of the above 
subalgebras, over all compact open subgroups, partially ordered by inclusion:
\begin{equation}\label{eq:4.22}
C_r^* (G) = \varinjlim C_r^* (G,K)
\end{equation}
Moreover it has a \inde{Bernstein decomposition}, analogous 
to \eqref{eq:4.19}, with a direct sum in the $C^*$-algebra sense:
\begin{equation}\label{eq:4.2}
C_r^* (G) = \varinjlim_{\mf S} \bigoplus_{\mf s \in \mf S} 
C_r^* (G)^{\mf s}
\end{equation}
where $\mf S$ runs over all finite subsets of $\mf B (G)$.
Here \inde{$C_r^* (G)^{\mf s}$} is the two-sided ideal of 
$C_r^* (G)$ generated by $\mc H (G)^{\mf s}$. Every subalgebra 
$C_r^* (G,K)$ lives in only finitely many Bernstein components.
\\[1mm]

The construction of the Schwartz algebra of $G$ is more 
complicated, we need to introduce a lot of things to achieve it.

A \inde{p-pair} is a pair $(P,A)$ consisting of a parabolic subgroup 
$P$ of $G$, and the identity component $A$ of the maximal split 
torus in the center of some Levi subgroup $M$ of $P$. Then
\index{P,A@$(P,A)$}
\begin{align}
& M = Z_G (A) = A \times \prefix^{0}{M} \\
& P = M N = A \prefix^{0}{M} N = Z_G (A) N
\end{align}
where $N$ is the unipotent radical of $P$.
For example $(G, A_G)$ is a p-pair, where \inde{$A_G$} is the 
maximal split torus of $Z(G)$. 

There is a unique p-pair $(\bar P,A)$ such that $\bar P \cap P = M$. 
The parabolic subgroup $\bar P$ is called the opposite of $P$. 
Clearly $\bar P = M \bar N$ where $\bar N \cap N = \{1\}$.
\index{Pbar@$\bar P$} \index{PAQB@$(P,A) \geq (Q,B)$}

Let $(Q,B)$ be another parabolic pair. Write $W (A|G|B)$
for the set of all homomorphisms $B \to A$ induced by inner
automorphisms of $G$. If $B=A$ then this is a group :
\index{WG|A@$W (G|A)$} \index{WAGB@$W (A|G|B)$}
\begin{equation}
W (G|A) := W (A|G|A) = N_G (A) / Z_G (A) = N_G (A) / M
\end{equation}
We say that $(P,A)$ dominates $(Q,B)$, written $(P,A) \geq (Q,B)$, 
if $P \supset Q$ and $A \subset B$.

Recall that we have chosen a maximal split torus $A_0$ of $G$, and 
let \inde{$P_0$} be a minimal parabolic subgroup containing it. 
We call a p-pair $(P,A)$ and its Levi factor $M$ semi-standard
if $A \subset A_0$, or equivalently $A_0 \subset M$. If moreover\\
$(P,A) \geq (P_0 ,A_0 )$, then we say that $(P,A)$ is standard.
Every p-pair is conjugate to a standard p-pair.

Let $X^* (A)$ be the character lattice of $A$ and put
\index{$\mf a$} \index{$\mf a^*$}
\begin{equation}
\begin{aligned}
& \mf a \; = \mr{Hom}_{\mh Z} (X^* (A), \mh R)\\
& \mf a^* = X^* (A) \otimes_{\mh Z} \mh R
\end{aligned}
\end{equation}
There is a natural homomorphism $H_M : M \to \mf a$, defined
by the equivalent conditions \index{$H_M$}
\begin{equation}
\begin{aligned}
& \inp{\chi}{H_M (m)} = -v (\chi (m))\\
& q^{\ds \inp{\chi}{H_M (m)}} = \norm{\chi (m)}_{\mh F}
\end{aligned}
\end{equation}
where $\chi \in X^* (A)$ and $q$ is the module of $\mh F$.
Conversely, if $\nu \in \mf a^*$ then we define an unramified
character $\chi_\nu$ of $M$ by \index{xchinu@$\chi_\nu$}
\begin{equation}
\chi_\nu (m) = q^{\ds \inp{\nu}{H_M (m)}} 
\end{equation}
For a parabolic subgroup $Q$ with $P \subset Q \subset G$, let
$\Sigma (Q,A) \subset \mf a^*$ be the set of roots of 
$Q$ with respect to $A$. By this we mean the set of 
$\alpha \in X^* (A) \setminus \{ 1 \}$ such that 
\inde{$\mf q_\alpha$} is nonzero, where $\mf q$ is the Lie 
algebra of $Q$ and \index{SigmaQA@$\Sigma (Q,A)$}
\[
\mf q_\alpha := \{ x \in \mf q : \mr{Ad}(a) x = \alpha (a) x
\quad \forall a \in A \}
\]
Then $\Sigma (G,A)$ is a root system and $\Sigma (P,A)$ is a
positive system of roots. Let $\Delta (P,A)$ be the corresponding 
set of simple roots. The Weyl group of $\Sigma (G,A)$ is 
naturally isomorphic to $W (G|A)$. Notice that $\Sigma (\bar P 
,A) = -\Sigma (P,A)$. \index{DeltaPA@$\Delta (P,A)$}

The minimal p-pair $(P_0 ,A_0)$ gives us a root system
\index{Sigma0@$\Sigma_0$} \index{Delta0@$\Delta_0$}
\[
\Sigma_0 = \Sigma (G,A_0 ) \subset \mf a_0^*
\]
with simple roots $\Delta_0 = \Delta (P_0 ,A_0 )$ and
Weyl group \inde{$W_0$}$\:= W (G|A_0 )$. Fix a $W_0$-invariant
inner product on $\mf a_0^*$, so that we may identify this 
vector space with its dual $\mf a_0$.

If $(P,A)$ is standard then $\Delta (P,A)$ is the set of
nonzero projections of $\Delta_0$ on $\mf a^*$, and $W (M|A_0 )$ 
is the parabolic subgroup of $W_0$ generated by 
\[
\{ s_\alpha : \alpha \in \Delta_0 ,\, \alpha (A) = 1 \}
\]
Let us also introduce the associated sets of positive
elements in $\mf a^*$ and $A$ :
\index{$\mf a^{*,+}$} \index{abar+@$\bar{\mf a}^{*,+}$} 
\index{$A^+$} \index{Abar+@$\bar{A}^+$}
\begin{align*}
\mf a^{*,+} &= \{ \nu \in \mf a^* : \inp{\nu}{\alpha} > 0 \; 
\forall \alpha \in \Delta (P,A) \} \\
\bar{\mf a}^{*,+} &= \{ \nu \in \mf a^* : \inp{\nu}{\alpha}
\geq 0 \; \forall \alpha \in \Delta (P,A) \} \\
A^+ &= \{ a \in A : \| \alpha (a) \|_{\mh F} > 1 
\; \forall \alpha \in \Delta (P,A) \} \\
\bar{A}^+ &= \{ a \in A : \| \alpha (a) \|_{\mh F} \geq 1 
\; \forall \alpha \in \Delta (P,A) \}
\end{align*}
In order to say when a function on $G$ is rapidly decreasing,
we need a length function on this group.
For $x \in GL (m,\mh F )$ let $x_{ij}$ and $x^{ij}$
be the entries of $x$ and $x^{-1}$, and define \index{$\mc N$}
\begin{equation}
\mc N (x) = \max \{ -v(x_{ij}), -v(x^{ij}) : 1 \leq i,j \leq m \}
\end{equation}
Notice that for all $x,y \in GL(m,\mh F)$
\begin{equation} 
0 \leq \mc N(xy) \leq \mc N(x) + \mc N (y)
\end{equation}
Pick an injective homomorphism $\tau : G \to GL (m,\mh F)$ and put
\begin{equation}
\sigma = \mc N \circ \gamma : G \to \mh Z_{\geq 0}
\end{equation}
Then $\sigma$ is a continuous length function on $G$. Let 
$\delta_{P_0}$ be the modular function of $P_0$. Using the 
decomposition \eqref{eq:4.24} we extend this to a function 
$\delta_0$ on $G$ satisfying 
\index{sigma@$\sigma$} \index{delta0@$\delta_0$}
\begin{equation}
\delta_0 (pk) = \delta_{P_0}(p) \qquad p \in P_0 , k \in K_0
\end{equation}
Harish-Chandra's spherical $\Xi$-function is \index{Xi@$\Xi$}
\begin{equation}
\Xi (g) = \int_{K_0} \delta_0 (kg) d \mu (k)
\end{equation}
Important properties of this function can be found in 
\cite[Paragraphe II]{Wal} and \cite[\S 4.2]{Sil1}.
For $n \in \mh N$ consider the following norm on $\mc H (G)$ :
\begin{equation}
\nu_n (f) = \sup_{g \in G} |f(g)| \Xi (g)^{-1} (\sigma (g) + 1)^n
\end{equation}
We say that $f \in C (G)$ decreases rapidly if 
$\nu_n (f) < \infty \quad \forall \in \mh N$.
Clearly the $\nu_n$ depend on the choice of \index{nun@$\nu_n$}
$\tau : G \to GL (m,\mh F)$, but the topology defined by the family 
$\{ \nu_n \}_{n=1}^\infty$ does not. For any compact open $K < G$
let \inde{$S (G,K)$} be the completion of $\mc H (G,K)$ for this
family of norms. According to Vign\'eras \cite[Theorem 29]{Vig}
this is a unital, nuclear Fr\'echet *-algebra, 
and a dense subalgebra of $C_r^* (G,K)$. Moreover an element of 
$\mc S (G,K)$ is invertible if and only if it is invertible in 
$C_r^* (G,K)$, so $\mc S (G,K)$ is closed under the holomorphic 
functional calculus of $C_r^* (G,K)$. 

For $K' \subset K$ there is still an inclusion $S (G,K) \to 
S (G,K')$, so we can take the inductive limit over all compact
open subgroups of $G$. This yields Harish-Chandra's 
\inde{Schwartz algebra}: 
\begin{equation}\label{eq:4.23}
\mc S (G) = \varinjlim \mc S (G,K)
\end{equation}
By definition it consists of all rapidly decreasing smooth 
functions on $G$. The obvious analogue of \eqref{eq:4.42} is
\begin{equation}\label{eq:4.62}
e_K \mc S (G) e_K = \mc S (G)^{K \times K} = \mc S (G,K)
\end{equation}
Compared to the above $C^*$-algebras, \inde{$\mc S (G)$} 
inherits fewer topological properties from its subalgebras. 
Namely, it is a complete Hausdorff locally convex algebra, 
but it is not metrizable, and its multiplication is only 
separately continuous \cite[\S III.6]{Wal}.

It does have a \inde{Bernstein decomposition}
\begin{equation}
\mc S (G) = \bigoplus_{\mf s \in \mf B (G)} S (G)^{\mf s}
\end{equation}
where \inde{$S (G)^{\mf s}$} is the completion of $\mc H 
(G)^{\mf s}$, a two-sided ideal in $\mc S (G)$. This follows 
from Theorem \ref{thm:4.31}, but of course it can be proved 
without using the full strength of that result. 
\\[2mm]
To characterize those $G$-representations that extend to 
$\mc S (G)$ we need to know more about smooth representations.

Let $(\pi ,V)$ be a smooth $G$-representation, and $P$ a parabolic
subgroup with unipotent radical $N$ and a Levi factor $M$. The
\inde{Jacquet module} associated to these data is
\index{$V_P$} \index{$V(N)$}
\begin{equation}
\begin{aligned}
& V_P = V / V(N) \\
& V(N) = \mr{span} \{ \pi (n) v - v : n \in N, v \in V \}
\end{aligned}
\end{equation}
We make it into an $M$-representation $(\pi_P ,V_P )$ by
\begin{equation}\label{eq:4.43}
\pi_P (m) j_P (v) = \delta_P^{-1/2}(m) j_P (\pi (m) v)
\end{equation}
where $j_P : V \to V / V(N)$ is the natural projection. By 
Frobenius reciprocity we get, for any smooth $M$-representation
$\sigma$ :
\begin{equation}\label{eq:4.44}
\mr{Hom}_G (\pi ,I_P^G \sigma ) \cong \mr{Hom}_M (\pi_P ,\sigma )
\end{equation}
For $\chi \in \mr{Hom}(A_G ,\mh C^\times )$ define the
generalized weight space \index{Vchi@$V_\chi$}
\begin{equation}
V_\chi = \big\{ v \in V : \exists n \in \mh N : 
(\pi (a) - \chi (a))^n v = 0 \; \forall a \in A_G \big\}
\end{equation}
If $V_\chi \neq 0$ then we call $\chi$ an \inde{exponent} 
of $(\pi ,V)$. If 
\[
\pi (a) v = \chi (a) v \; \forall v \in V , a \in A_G
\]
then we say that $V$ admits the \inde{central character} $\chi$. 
Similarly for a p-pair $(P,A)$ and a smooth $M$-representation
we have exponents in Hom$(A,\mh C^\times )$. The set of exponents
of the Jacquet module $V_{\bar P}$ is \index{$\mc X_\pi (P,A)$}
\begin{equation}
\mc X_\pi (P,A) = \{ \chi \in \mr{Hom} (A, \mh C^\times ) : 
V_{\bar P ,\chi} \neq 0 \}
\end{equation}
Notice the curious shift in the notation, from $\bar P$ to $P$.
This is designed to make a nicer formulation of the Langlands 
classification possible.

Let us characterize square-integrable, discrete series and 
tempered\\ $G$-representations. Our first description is due to 
Casselman \cite[Theorem 4.4.6]{Cas}.

\begin{prop}\label{prop:4.26}
Let $\pi$ be an admissible $G$-representation which admits a
unitary central character. The following are equivalent :
\begin{itemize}
\item Every matrix coefficient of $\pi$ is square-integrable
on $G / A_G$
\item If $(P,A)$ is a semi-standard p-pair, $\chi \in \mc X_\pi 
(P,A)$ and $a \in \bar{A}^+$ is such that there is an $\alpha 
\in \Delta_0$ with $|\alpha (a)| \neq 1$, then $|\chi (a)| < 1$
\item For every semi-standard p-pair $(P,A)$ and every 
$\chi \in \mc X_\pi (P,A)$ we can write 
\[
\log |\chi | = \sum_{\alpha \in \Delta (P,A)} \chi_\alpha \alpha
\quad \mr{with}\; \chi_\alpha < 0
\]
\end{itemize}
We say that $\pi$ is square-integrable if it satisfies these 
conditions.\index{representation!square-integrable}
\end{prop}

Every square-integrable representation is unitary and completely
reducible, see \cite[Corollary 1.11.8]{Sil1} or 
\cite[Lemme III.1.3]{Wal}. A more restrictive notion is that of
a discrete series representation.

\begin{prop}\label{prop:4.27}
Let $(\pi ,V)$ be an irreducible admissible $G$-representation.
The following are equivalent :
\begin{itemize}
\item $(\pi ,V)$ is a subrepresentation of $(\lambda ,L^2 (G))$
\item $G$ is semi-simple and $\pi$ is square-integrable
\end{itemize}
If $(\pi ,V)$ satisfies these conditions then it is called a
discrete series representation. 
\index{representation!discrete series}
\end{prop}

By \cite[Proposition 18.4.2]{Dix} such a representation does 
indeed give an isolated point in Prim$(C_r^* (G))$.

A (smooth) function $f$ on $G$ is tempered if there exist 
$C,N \in (0,\infty )$ such that \index{tempered function}
\begin{equation}
|f (g)| \leq C \Xi (g) (1 + \sigma (g))^N \quad \forall g \in G
\end{equation}
\begin{prop}\label{prop:4.28}
Let $\pi$ be an admissible $G$-representation. The following 
are equivalent : \index{representation!tempered}
\begin{itemize}
\item $\pi$ extends continuously to $\mc S (G)$
\item Every matrix coefficient of $\pi$ is a tempered function
\item If $(P,A)$ is a semi-standard p-pair, $\chi \in 
\mc X_\pi (P,A)$ and $a \in A^+$, then \\ $|\chi (a)| \leq 1$
\item For every semi-standard p-pair $(P,A)$ and every 
$\chi \in \mc X_\pi (P,A)$ we can write 
\[
\log |\chi | = \sum_{\alpha \in \Delta (P,A)} \chi_\alpha \alpha
\quad \mr{with}\; \chi_\alpha \leq 0
\]
\end{itemize}
The representation $\pi$ is said to be tempered if and only if
these conditions hold.
\end{prop}
\emph{Proof.} 
Almost everything follows from 
\cite[Proposition III.2.2 and \S III.7]{Wal}. The only thing 
left is to show that all matrix coefficients of an admissible 
$\mc S (G)$-representation are tempered. This follows from the
admissibility, combined with the observation that the collection 
of tempered $K$-biinvariant functions on $G$ is the linear dual 
of $\mc S (G,K). \qquad \Box$
\\[2mm]

The properties temperedness and pre-unitarity are preserved under
normalized induction:

\begin{prop}\label{prop:4.34}
Let $(P,A)$ be a semi-standard p-pair and $(\pi ,V)$ an admissible
$M$-representation. Then
\begin{enumerate}
\item $I_P^G (\pi )$ is tempered if and only if $\pi$ is tempered
\item $I_P^G (\pi )$ is pre-unitary if and only if $\pi$ 
is pre-unitary
\end{enumerate}
\end{prop}
\emph{Proof.} 
1. comes from \cite[Lemme III.2.3]{Wal}.\\ 
2. It is clear that $I_P^G (\pi )$ cannot be pre-unitary if $\pi$ 
is not. It remains to produce a $G$-invariant inner product on 
$I_P^G (V)$, given an $M$-invariant inner product on $V$. This is 
achieved by setting
\begin{equation}\label{eq:4.65}
\inp{f}{f'} = \int_{K_0} \inp{f(k)}{f'(k)} d \mu (k)
\end{equation}
Notice that $I_P^G (V)$ is usually not complete with respect to
this inner product, even is $V$ is. $\qquad \Box$
\\[2mm]

Like in Section \ref{sec:3.2}, let $\Lambda$ be the set of triples
$(P,\sigma ,\nu )$, where $(P,A)$ is a standard p-pair, $\sigma$ is 
an irreducible tempered representation of $M = Z_G (A)$ and 
$\nu \in \mf a^*$. With such a triple we associate the admissible
$G$-representation \index{$I (P,\sigma ,\nu)$}
\index{Lambda@$\Lambda$} \index{Psnu@$(P,\sigma ,\nu )$}
\begin{equation}
I (P,\sigma ,\nu ) = I_P^G (\sigma \otimes \chi_\nu ) =
\mr{Ind}_P^G \big( \sigma \otimes \chi_\nu \otimes 
\delta_P^{1/2} \big)
\end{equation}
The set of Langlands data is \index{Lambda+@$\Lambda^+$}
\begin{equation}
\Lambda^+ = \big\{ (P,\sigma ,\nu ) \in \Lambda : 
\nu \in \mf a^{*,+} \big\}
\end{equation}
A somewhat extended version of the \inde{Langlands classification} 
for reductive $p$-adic groups reads :\index{$J (P,\sigma ,\nu )$}

\begin{thm}\label{thm:4.29}
Let $(P,\sigma ,\nu ) ,\, (P',\sigma' ,\nu' ) \in \Lambda^+$. 
\begin{enumerate}
\item The $G$-representation $I (P,\sigma ,\nu )$ is indecomposable 
and has a unique irreducible quotient, which we call 
$J (P,\sigma ,\nu )$. 
\item For every $\pi \in \mr{Irr}(G)$ there is precisely one 
Langlands datum $(P,\sigma ,\nu )$ such that $\pi$ is equivalent 
to $J (P,\sigma ,\nu )$.
\item If $J (P,\sigma ,\nu )$ is equivalent to a subquotient of
$I (P',\sigma' ,\nu' )$, then $\nu' - \nu \in \bar{\mf a}_0^{*,+}$ and
$P' \subset P$. If $P' = P$ then also $\sigma' = \sigma$ and
$\nu' = \nu$.
\end{enumerate}
\end{thm}
\emph{Proof.} 
For 1 and 2 see \cite{Sil2} or \cite[\S XI.2]{BoWa}.
As concerns 3, by \cite[Lemma XI.2.13]{BoWa} we have 
$\nu' - \nu \in \bar{\mf a}_0^{*,+}$. Now it follows from the
definition of $\Lambda^+$ that $P' \subset P$. Suppose that 
$(P,\sigma ,\nu) \neq (P',\sigma',\nu' )$ while $P = P'$.
Then, again by \cite[Lemma XI.2.13]{BoWa}, $\nu \neq \nu'$.
But by Frobenius reciprocity $\sigma$ is equivalent to a subquotient 
of $I_P^G (\sigma' \otimes \chi_{\nu'} ) (N)$. Hence $\nu = \nu' 
\circ w$ for some $w \in W(G|A) \setminus \{1\}$. Since 
$\nu' \in {\mf a'}^{*,+}$ there is an $\alpha \in \Delta_0$
with $\inp{\nu'}{\alpha} > 0$ but $\inp{\nu}{\alpha} < 0$. This 
contradicts the positivity of $\nu$ with respect to $P' = P. 
\qquad \Box$ 
\\[4mm]

\section{The Plancherel theorem}
\label{sec:4.3}

The Plancherel formula for $G$ is an explicit decomposition 
of the trace $\tau$ in terms of the traces of irreducible 
$G$-representations. Closely related is the Plancherel theorem,
which describes the image of $\mc S (G)$ under the Fourier 
transform. The crucial theorems in this section are due to 
Harish-Chandra \cite{HC2}, but unfortunately he never 
published the proofs. Based upon Harish-Chandra's notes, 
Waldspurger \cite{Wal} provided full proofs of these results, 
which we will describe in as much detail as we need. We try 
to set up a complete analogy with Section \ref{sec:3.3}. 
In particular we refine the Langlands classification using 
parabolic induction in stages.

We start with a semi-standard p-pair $(P,A)$ and an irreducible
square-integrable $M$-representation $(\omega ,E)$. Let
$(\breve \omega ,\breve E )$ be its contragredient, and construct
the admissible $G \times G$-representation 
\index{$L (\omega ,P)$} \index{omegaE@$(\omega ,E)$}
\begin{equation}
L (\omega ,P) = I_{P \times P}^{G \times G} (E \otimes \breve E ) =
I_P^G (E) \otimes I_P^G (\breve E)
\end{equation}
Using \eqref{eq:4.65} we make $I_P^G (E)^K$ into a finite 
dimensional Hilbert space, for every compact open $K < G$. This 
allows us to identify $\breve{I_P^G (E)}$ with $I_P^G (\breve E )$ 
as representations, and with $I_P^G (E)$ as vector spaces. Thus we 
can turn $L (\omega ,P)$ into a nonunital *-algebra with 
\begin{align}
& (f_1 \otimes f_2 ) (f_3 \otimes f_4 ) = \label{eq:4.63}
\inp{f_2}{f_3} f_1 \otimes f_4\\
& (f_1 \otimes f_2 )^* = f_2 \otimes f_1 \label{eq:4.64}
\end{align}
Notice that for every $\chi \in X_{unr}(M)$ the representation
$\chi \otimes \omega$ is still square-integrable, and that $L 
(\chi \otimes \omega, P)$ can be identified with $L (\omega ,P)$.
Let \inde{$K_\omega$} be the set of $k \in X_{nr}(M)$ such that
$k \otimes \omega$ is equivalent to $\omega$. This is a finite
subgroup of $X_{unr}(M)$. For every $k \in K_\omega$ we pick a 
unitary intertwiner \index{$I (k,\omega)$}
\begin{equation}
\tilde \omega_k : (k \otimes \omega ,E) \to (\omega ,E)
\end{equation}
This induces an automorphism of $L (\omega ,P)$ by 
\begin{equation}
I (k,\omega) = I_{P \times P}^{G \times G} (\tilde \omega_k 
\otimes \tilde \omega_k^{-t}) = I_P^G (\tilde \omega_k ) \otimes
I_P^G (\tilde \omega_k^{-t})
\end{equation}
where $\tilde \omega_k^{-t}$ is the inverse transpose of 
$\tilde \omega_k$. Then $I (k, \omega) \in \mr{Aut}_{G 
\times G} \big( L (\omega ,P) \big)$ is independent of the 
choice of $\tilde \omega_k$, and in general nontrivial, cf. 
\cite[\S VI.1]{Wal}.

It is more difficult to define intertwiners corresponding to
elements of the various Weyl groups. First we notice that for any
p-pair $(P',A)$ with the same Levi factor $M ,\, \omega$ can also 
be lifted to a representation of $P'$ that is trivial on $N'$.
Let $(Q, A^g)$, with $g \in G$, be yet another semi-standard 
p-pair, and put 
\[
n = [g] \in W (A^g |G|A)
\]
The equivalence class of the $M^g$ representation 
$(\omega^{g^{-1}},E)$ depends only on $n$, and is therefore 
denoted by $n \omega$.

In \cite[Paragraphe V]{Wal} certain normalized intertwiners
$\prefix^{o}{c_{Q|P}} (n, \omega)$ are constructed. 
Preferring the simpler notation \inde{$I (n,\omega)$}, we recall 
their properties.

\begin{thm}\label{thm:4.30}
Let $(P,A) ,\, (P',A')$ and $(Q,B)$ be semi-standard p-pairs,
and $n \in W (B|G|A)$. There exists an intertwiner
\[
I (n,\chi \otimes \omega) \in \mr{Hom}_{G \times G} 
(L (\omega ,P) , L (n \omega ,Q) )
\]
with the following properties :
\begin{itemize}
\item $\chi \to I (n, \chi \otimes \omega )$ is a rational 
function on $X_{nr}(M)$
\item $I (n ,\chi \otimes \omega )$ is unitary and regular
for $\chi \in X_{unr}(M)$.
\item If $n' \in W (A' |G|B)$ then
\[
I (n',n (\chi \otimes \omega) ) I (n,\chi \otimes \omega ) = 
I (n'n, \chi \otimes \omega )
\]
\end{itemize}
\end{thm}

Let $\Gamma_{rr} \big( X_{nr}(M) ; L (\omega ,P) \big)$ be the 
algebra of rational sections that are regular on $X_{unr}(M)$. 
We define an action of $K_\omega$ on this algebra by 
\begin{equation}\label{eq:4.59}
k f (\chi ) = I(k,\omega ) f(k^{-1} \chi )
\end{equation}
Similarly, for $n$ as in Theorem \ref{thm:4.30} we define an
algebra homomorphism 
\begin{equation}\label{eq:4.60}
\begin{aligned}
& n : \Gamma_{rr} \big( X_{nr}(M)) ; L (\omega ,P) \big) \to
\Gamma_{rr} \big( X_{nr}(Z_G (B)) ; L (n \omega, Q) \big) \\
& n f (\chi) = I (n, \omega) f ( \chi \circ n )
\end{aligned}
\end{equation}
To define the Fourier transform we construct a scheme containing
all tempered $G$-representations. For every Levi subgroup $M$ of 
$G$ choose a set $\Delta_M$ of irreducible square-integrable
$M$-representations, with the property that for every 
square-integrable $\pi \in \mr{Irr}(M)$ there exists exactly one
$\omega \in \Delta_M$ such that $\pi$ is equivalent to
$\chi \otimes \omega$, for some $\chi \in X_{nr}(M)$.

Let $\Xi$ be the scheme consisting of all quadruples 
$(P,A,\omega ,\chi )$, with $(P,A)$ a semi-standard p-pair, 
$\omega \in \Delta_M$ and $\chi \in X_{nr}(M)$. This is
a countable disjoint union of complex algebraic tori. Let $\Xi_u$ 
be the smooth submanifold obtained by the restriction
$\chi \in X_{unr}(M)$. Notice that $\Xi$ is naturally a finite 
cover of the set $\Theta$ defined in \cite[p. 305]{Wal}. For
$\xi = (P,A, \omega, \chi) \in \Xi$ we put $\pi (\xi) = I_P^G 
(\chi \otimes \omega)$. Let $\mc L_{\Xi}$ be the vector
bundle over $\Xi$ which is trivial on every component and whose 
fiber at $\xi$ is $L (\omega ,P)$. We say that a section of this 
bundle is algebraic or rational if it is supported on only finitely 
many components, and has the required property on every component. 
Now we define the \inde{Fourier transform} \index{pix@$\pi (\xi)$}
\index{Xi@$\Xi$} \index{PAoc@$(P,A,\omega ,\chi )$}
\begin{equation}\label{eq:4.86}
\begin{aligned}
&\mc F : \mc H (G) \to \mc O (\Xi ; \mc L_{\Xi} ) \\
&\mc F (f) (P,A, \omega, \chi) = 
I (P,A, \omega, \chi)(f) \in L (\omega ,P)
\end{aligned}
\end{equation}
This is not the same as $\hat f (\chi \otimes \omega, P)$, as
in \cite[\S VII.1]{Wal}! We adjusted the latter to make $\mc F$ 
multiplicative. Fortunately the difference is not too big, so 
most results remain valid. \index{Fmc@$\mc F$}

We construct a locally finite groupoid $\mc W$ as follows. 
The objects of $\mc W$ are triples $(P,A,\omega)$ with $(P,A)$ a
semi-standard p-pair and $\omega \in \Delta_M$. The morphisms
from $(Q,B,\eta )$ to $(P,A,\omega )$ are pairs $(k,n)$ with
the following properties \index{Wmc@$\mc W$}
\begin{itemize}
\item $k \in K_\omega$
\item $n \in W (A|G|B)$ and $n B = A$
\item $n \eta$ is equivalent to $\chi \otimes \omega$, for some
$\chi \in X_{nr}(M)$
\end{itemize}
The multiplication in $\mc W$, if possible, is given by
\begin{equation}
(k,n) (k',n') = (k (k' \circ n), n n')
\end{equation}
Let $\Gamma_{rr} (\Xi ; \mc L_\Xi)$ be the algebra of rational 
sections of $L_\Xi$ that are regular on $\Xi_u$ :
\index{Gammarxl@$\Gamma_{rr} (\Xi ; \mc L_\Xi)$}
\begin{multline}
\Gamma_{rr} (\Xi ; \mc L_\Xi) = \bigoplus_{(P,A,\omega )} 
\Gamma_{rr} \big( X_{nr}(M) ; L (\omega ,P) \big) \\
= \bigoplus_{(P,A,\omega )} \big\{ f \in Q (\mc O (X_{nr}(M)))
\otimes L (\omega, P) : f \;\mr{is \; regular \; on}\; X_{unr}(M) 
\big\}
\end{multline}
From \eqref{eq:4.59} and \eqref{eq:4.60} we get an action of
the groupoid $\mc W$ on this algebra. By construction the image 
of $\mc H (G)$ under the Fourier transform consists of 
$\mc W$-invariant sections.

Because $(\omega,E)$ is admissible
\begin{equation}
C^\infty (X_{unr}(M)) \otimes L (\omega ,P)^{K \times K}
\cong C^\infty (X_{unr}(M)) \otimes I_P^G (E)^K \otimes 
I_P^G (\breve E)^K
\end{equation}
is in a natural way a Fr\'echet space, for every compact open 
$K < G$. Endow $C^\infty (X_{unr}(M)) \otimes L (\omega ,P)$ 
with the inductive limit topology. This also gives a topology
on $C^\infty_c (\Xi_u ; \mc L_\Xi )$, as the inductive limit of
finite direct sums of such algebras. Notice that these are all 
*-algebras by \eqref{eq:4.61}. Clearly the action of $\mc W$
extends continuously to $C^\infty_c (\Xi_u ; \mc L_\Xi )$. 

Now the Plancherel theorem for reductive $p$-adic groups
\cite[p. 320]{Wal} tells us that

\begin{thm}\label{thm:4.31}
The Fourier transform
\[
\mc F : \mc S (G) \to C^\infty_c (\Xi_u ; \mc L_\Xi )^{\mc W}
\]
is an isomorphism of topological *-algebras.
\end{thm}

This guides us to the Fourier transform of $C_r^* (G)$. 
For $(\omega ,E)$ as on page \pageref{eq:4.63}, let 
\inde{$\mc K (\omega ,P)$} be the algebra of compact operators 
on the Hilbert space completion of $I_P^G (E)$. Notice that
\begin{equation}
\mc K (\omega ,P) = \varinjlim L (\omega ,P)^{K \times K}
\end{equation}
in the $C^*$-algebra sense, and that the intertwiner $I (n, \omega)$
extends to $\mc K (\omega ,P)$ because it is unitary. Let 
$\mc K_\Xi$ be the vector bundle over $\Xi$ whose fiber at 
$(P,A,\omega ,\chi )$ is $\mc K (\omega ,P)$, and 
\inde{$C_0 (\Xi_u ; \mc K_\Xi )$} the $C^*$-completion of
\[
\bigoplus_{(P,A,\omega)} C (X_{unr}(M) ; \mc K (\omega ,P))
\]
Plymen \cite[Theorem 2.5]{Ply2} proved that \index{KXi@$\mc K_\Xi$}

\begin{thm}\label{thm:4.42}
The Fourier transform extends to an isomorphism of $C^*$-algebras
\[
C_r^* (G) \isom C_0 (\Xi_u ; \mc K_\Xi )^\mc W
\]
\end{thm}

The subalgebras $\mc S (G,K)$ are more manageable than $\mc S (G)$,
so it pays off to describe their images under $\mc F$.

\begin{thm}\label{thm:4.32}
Let $K$ be a compact open subgroup of $G$. There exists a finite
set of triples $(P_i ,A_i ,\omega_i ) ,\; i = 1,\ldots ,n_K$, such
that the Fourier transform induces algebra homomorphisms
\[
\begin{array}{rrr}
\mc H (G,K) & \to & \bigoplus_{i=1}^{n_K} 
  \big( \mc O (X_{nr}(M_i )) \otimes L(\omega_i ,P_i 
  )^{K \times K} \big)^{\mc W_i} \\
\mc S (G,K) & \to & \bigoplus_{i=1}^{n_K} 
  \big( C^\infty (X_{unr}(M_i )) \otimes L(\omega_i ,P_i 
  )^{K \times K} \big)^{\mc W_i} \\
C_r^* (G,K) & \to & \bigoplus_{i=1}^{n_K} 
  \big( C (X_{unr}(M_i )) \otimes L(\omega_i ,P_i )^{K \times K}
  \big)^{\mc W_i}
\end{array}
\]
where $\mc W_i$ is the isotropy group of $(P_i ,A_i ,\omega_i )$ in 
$\mc W$. The first map is injective, the second is an isomorphism
of Fr\'echet *-algebras and the third is an isomorphism of
$C^*$-algebras.
For every $w \in \mc W_i$ there is a rational, unitary element
\[
u_w \in C^\infty (X_{unr}(M_i )) \otimes L(\omega_i ,P_i 
)^{K \times K}
\]
such that for every $f \in C (X_{unr}(M_i )) \otimes 
L(\omega_i ,P_i )^{K \times K}$
\begin{equation}\label{eq:4.61}
w f (\chi ) = u_w (\chi ) f (w^{-1} \chi ) u_w^{-1}(\chi )
\end{equation}
\end{thm}
\emph{Proof.} 
By \cite[Th\'eor\`eme VIII.1.2]{Wal} there is only a finite number 
of association classes among the objects of $\mc W$ on which the
idempotent $e_K$ does not act as zero. Pick one representant 
$(P_i ,A_i, \omega_i )$ in every such association class. From 
\eqref{eq:4.86} and Theorems \ref{thm:4.31} and \ref{thm:4.42} we 
immediately get the required description of the Fourier transforms 
of $\mc H (G,K) \,,\mc S (G,K)$ and $C_r^* (G,K)$. 

Every automorphism of $L (\omega_i ,P_i )^{K \times K} \cong
\mr{End}\big(I_P^G (E)^K \big)$ is inner, so the
formula \eqref{eq:4.61} holds for some $u_w$. Using Theorem 
\ref{thm:4.30} we can arrange that $u_w$ is rational on 
$X_{nr}(M_i )$ and unitary on $X_{unr}(M_i ). \qquad \Box$
\\[2mm]

Purely representation-theoretic consequences of the above 
isomorphisms are:

\begin{cor}\label{cor:4.33}
\begin{enumerate}
\item Every irreducible tempered $G$-representation is a direct
summand of $I(\xi )$, for some $\xi \in \Xi_u$.
\item For any $w \in \mc W$ and $\xi \in \Xi$ such that $w \xi$ is
defined, the $G$-representations $I(\xi )$ and $I(w \xi )$ have the
same irreducible subquotients, counted with multiplicity.
\end{enumerate}
\end{cor}
\emph{Proof.}
1. Let $K$ be a compact open subgroup of $G$ such that $V^K \neq 0$,
and let $V'$ be an irreducible submodule of $V^K$, considered
as a $\mc S (G,K)$-representation. With Theorem \ref{thm:4.32} and 
the same argument as in the proof of Corollary \ref{cor:3.20}.2 we 
deduce that $V'$ is a direct summand of $I (P_i ,A_i ,\omega_i 
,\chi)^K$ for some $\chi \in X_{unr}(M)$. Because $V'$ generates $V$
as a $G$-module, $V$ is a constituent of $I (P_i ,A_i ,\omega_i 
,\chi)$. By Proposition \ref{prop:4.34}.2 the latter is completely
reducible, so $V$ is in fact equivalent to a direct summand.of
$I (P_i ,A_i ,\omega_i ,\chi)$.\\
2. By \cite[Corollary 2.3.3]{Cas} we have to show that the 
characters of $I(\xi )$ and $I(w \xi )$ are the same, i.e. that
the function
\begin{equation}
\mc H (G) \times X_{nr}(M) \to \mh C : (f ,\chi) \to
\mr{tr}\, I(P,A,\omega ,\chi )(f) - 
\mr{tr}\, I (wP,wA,w \omega ,\chi \circ w^{-1})(f)
\end{equation}
is identically 0. Because this is a polynomial function of $\chi$,
it suffices to show that it is 0 on $\mc H (G) \times X_{unr}(M)$.
This is an immediate consequence of Theorem \ref{thm:4.31}.
$\qquad \Box$
\\[2mm]

For $\xi = (P,A,\omega ,\chi ) \in \Xi$ we define 
\index{$A(\xi )$} \index{omegax@$\omega (\xi)$} 
\index{nux@$\nu (\xi )$} \index{$P(\xi )$} \index{$M(\xi )$} 
\begin{equation}\label{eq:4.67}
\begin{aligned}
A(\xi ) &= \{ a \in A : \alpha (a) = 1 \;\mr{if}\; \inp{\log 
|\chi |}{\alpha} = 0 \; \forall \alpha \in \Sigma (P,A) \}\\
M(\xi ) &= Z_G (A(\xi ))\\
P(\xi ) &= P M(\xi )\\
\omega (\xi ) &= I_{M(\xi ) \cap P}^{M (\xi )} 
\big(\omega \otimes \chi |\chi |^{-1} \big)\\
\nu (\xi ) &= \log |\chi |
\end{aligned}
\end{equation}
By Proposition \ref{prop:4.34} $\omega (\xi )$ is a pre-unitary 
tempered $M(\xi )$-representation. Like in \cite[\S XI.9]{KnVo} 
these objects are designed to divide parabolic induction into 
stages:
\begin{equation}\label{eq:4.66}
\begin{aligned}
I_{P(\xi )}^G \big( |\chi | \otimes \omega (\xi ) \big) &\cong 
\mr{Ind}_{P(\xi )}^G \big( \delta_{P(\xi )}^{1/2} \otimes |\chi | 
\otimes \omega (\xi ) \big) \\
&\cong \mr{Ind}_{P(\xi )}^G \big( \delta_{P(\xi )}^{1/2} \otimes 
|\chi | \otimes \mr{Ind}_{M(\xi ) \cap P}^{M(\xi )} \big( 
\delta_{P \cap M(\xi)}^{1/2} \otimes \chi |\chi|^{-1} \otimes 
\omega \big) \big) \\
&\cong \mr{Ind}_{P M(\xi )}^G \big( \delta_{P M(\xi )}^{1/2} \otimes 
\mr{Ind}_{M(\xi ) \cap P}^{M(\xi )} \big( \delta_{P \cap M(\xi)}^{1/2} 
\otimes \chi \otimes \omega \big) \big) \\
&\cong \mr{Ind}_{P M(\xi )}^G \big(  \mr{Ind}_{M(\xi ) 
\cap P}^{M(\xi )} \big( \delta_{P M(\xi )}^{1/2} \otimes \delta_{P 
\cap M(\xi)}^{1/2} \otimes \chi \otimes \omega \big) \big) \\
&\cong \mr{Ind}_P^G \big( \delta_P^{1/2} \otimes \chi \otimes 
\omega \big) \quad = \quad I(\xi )
\end{aligned}
\end{equation}
We say that $(P,A,\omega ,\chi ) \in \Xi^+$ if $(P,A)$ is standard and
$\log |\chi | \in \bar{a}^{*,+}$. This choice of a "positive cone" in
$\Xi$ is justified by the next result.

\begin{lem}\label{lem:4.35}
Every $\xi \in \Xi$ is $\mc W$-associate to an element of $\Xi^+$. If
$\xi_1 , \xi_2 \in \Xi^+$ are $\mc W$-associate, then the objects
$A(\xi_i ) ,\, M(\xi_i ) ,\, P(\xi_i )$ and $\nu (\xi_i )$ are the 
same for $i=1$ and $i=2$, while $\omega (\xi_1 )$ and 
$\omega (\xi_2 )$ are equivalent $M(\xi_i )$-representations. 
\end{lem}
\emph{Proof.} 
Every p-pair is conjugate to a standard p-pair, and by 
\cite[Section 1.15]{Hum} every $W_0$-orbit in $\mf a_0^*$ contains 
a unique point in positive chamber $\mf a^{*,+}$. This proves the 
first claim, and it also shows that 
\begin{equation}
\log |\chi_1 | = \log |\chi_2 | \in \mf a_0^*
\end{equation}
Hence the $\nu$'s, $A$'s and $M$'s are the same for $i=1$ and $i=2$.
Because 
\begin{equation}
\Delta (P_i ,A_i ) = \{ \alpha \big|_{\ds \mf a^*} : 
\alpha \in \Delta_0 ,\, \inp{\log |\chi_i |}{\alpha} > 0 \}
\end{equation}
we must also have $P (\xi_1 ) = P(\xi_2 )$. If now $w \in \mc W$ is 
such that $w \omega_1 \cong \omega_2$, then by Theorem \ref{thm:4.30},
applied to $M(\xi_i )$, there is a unitary intertwiner between
$\omega (\xi_1 )$ and $\omega (\xi_2 ). \qquad \Box$
\\[2mm]

Some immediate consequences of the above definitions and Theorem 
\ref{thm:4.29} are:

\begin{prop}\label{prop:4.36}
Take $\xi = (P,A,\omega ,\chi ) \in \Xi^+$.
\begin{enumerate}
\item Let $\tau$ be an irreducible direct summand of $\omega (\xi )$.
Then $(P(\xi ),\tau ,\nu (\xi )) \in \Lambda^+$.
\item The functor $I_{P(\xi )}^G$ induces an isomorphism 
\[
\mr{End}_G (I(\xi )) \cong \mr{End}_{M(\xi )}(\omega (\xi ))
\]
\item The irreducible quotients of $I(\xi )$ are precisely the modules 
$J (P(\xi ),\tau ,\nu (\xi ))$ with $\tau$ as above.
\end{enumerate}
\end{prop}

\begin{thm}\label{thm:4.37}
For every $\pi \in \mr{Irr}(G)$ there exists a unique association 
class $\mc W (P,A,\omega ,\chi ) \in \Xi / \mc W$ such that the 
following equivalent statements hold :
\begin{enumerate}
\item $\pi$ is equivalent to an irreducible quotient of $I (\xi^+ )$,
for some\\ $\xi^+ \in \mc W (P,A,\omega ,\chi ) \cap \Xi^+$.
\item $\pi$ is equivalent to an irreducible subquotient of 
$I (P,A,\omega ,\chi )$, and $P$ is maximal for this property.
\end{enumerate}
\end{thm}
\emph{Proof.} 
1. Let $(Q,\sigma ,\nu )$ be the Langlands datum associated to $\pi$. 
Write $\mc W^M, \Xi^M$ etcetera for $\mc W, \Xi$, but now 
corresponding to $M$ instead of $G$. By Theorem \ref{thm:4.31} 
there exists a unique association class 
\[
\mc W^M \xi = \mc W^M (P,A,\omega ,\chi ) \in \Xi^M_u / \mc W^M
\]
such that $\sigma$ is a direct summand of $I^M (\xi ) = I_P^M 
(\omega \otimes \chi)$. Pick $\xi^+ \in \mc W^M \xi \cap \Xi^+$. By 
Proposition \ref{prop:4.36}.3 $\pi$ is equivalent to an irreducible 
quotient of $I (\xi^+ )$, and by Lemma \ref{lem:4.35} and Theorem
\ref{thm:4.31} the class $\mc W \xi^+ = \mc W \xi \in \Xi / \mc W$
is unique for this property.

2. Suppose that $\xi' = (P',A',\omega' ,\chi' ) \in \Xi^+$ and that
$\pi$ is equivalent to a\\ subquotient of $I (\xi' )$ which is not a 
quotient. By Theorem \ref{thm:4.29}.3 we have\\ $\nu (\xi' ) - 
\nu(\xi ) \in \mf a (\xi' )^{*,+}$ and $A(\xi ) \subsetneq A(\xi' )$. 
For $\alpha \in \Delta_0$ we have
\[
\inp{\nu (\xi' )}{\alpha} = 0 \; \Rightarrow \;
\inp{\nu (\xi )}{\alpha} = 0 
\]
so by \eqref{eq:4.67} $A \subsetneq A'$ and $P \supsetneq P'$.
Therefore the conditions 1 and 2 are equivalent. $\qquad \Box$
\\[4mm]

\section{Noncommutative geometry}
\label{sec:4.4}

Now that we know quite something about the representation theory
of reductive $p$-adic groups, we can turn to the study of their 
noncommutative geometry with more confidence. More specifically,
we compare the periodic cyclic homologies of $\mc H (G)$ and
$\mc S (G)$, and the $K$-theory of $C_r^* (G)$. Although these 
results are new and technical, this section remains very short, 
because we already did most of the hard work. We will discuss these
comparison theorems in relation with the Baum-Connes conjecture.

First we will prove an analogue of \eqref{eq:3.68} for the Hecke 
algebra of a reductive $p$-adic group. Since $\mc S (G)$ is defined 
as an inductive limit of Fr\'echet algebras, we take its periodic 
cyclic homology with respect to the completed inductive tensor 
product $\oot$.

\begin{thm}\label{thm:4.43}
Let $\mf s \in \mf B (G)$ be a Bernstein component and $K_{\mf s}$ 
a compact open subgroup of $G$ as in Proposition \ref{prop:4.42}.3.
The Chern character for $\mc S (G,K_{\mf s})^{\mf s}$ induces an 
isomorphism 
\[
K_* \big( C_r^* (G)^{\mf s} \big) \otimes \mh C \isom 
HP_* \big( \mc S (G)^{\mf s}, \oot \big)
\]
The direct sum of these maps, over all $\mf s \in \mf B (G)$, is
a natural isomorphism
\[
K_* (C_r^* (G)) \otimes \mh C \isom 
HP_* (\mc S (G), \oot)
\]
\end{thm}
\emph{Proof.} 
Since $\mc S (G,K_{\mf s})^{\mf s}$ is a direct summand of
$\mc S (G,K_{\mf s})$, by Theorem \ref{thm:4.32} there are 
finitely many components $(P_i,A_i,\omega_i )$ of $\Xi$ such that
\begin{equation}\label{eq:4.79}
\begin{array}{lll}
\mc S (G,K_{\mf s})^{\mf s} & \cong & \bigoplus\limits_i 
\big( C^\infty (X_{unr}(M_i )) \otimes L(\omega_i ,P_i 
)^{K_{\mf s} \times K_{\mf s}} \big)^{\mc W_i} \\
C_r^* (G,K_{\mf s})^{\mf s} & \cong & \bigoplus\limits_i 
\big( C (X_{unr}(M_i )) \otimes L(\omega_i ,P_i )^{K_{\mf s} 
\times K_{\mf s}} \big)^{\mc W_i}
\end{array}
\end{equation}
Note that the single Bernstein component $\mf s$ generally contains 
more than component of $\Xi$. According to Theorem \ref{thm:2.19} 
the inclusion induces an isomorphism
\begin{equation}\label{eq:4.80}
K_* \big( \mc S (G,K_{\mf s})^{\mf s} \big) \isom 
K_* \big( C_r^* (G,K_{\mf s})^{\mf s} \big)
\end{equation}
and by Theorem \ref{thm:2.18} the Chern character induces an
isomorphism
\begin{equation}\label{eq:4.81}
K_* \big( \mc S (G,K)^{\mf s} \big) \otimes \mh C \isom
HP_* \big( \mc S (G,K)^{\mf s} \big)
\end{equation}
For any compact open $K \subset K_{\mf s}$ the algebra
$L(\omega_i ,P_i )^{K \times K}$ is finite dimensional and 
simple, so the inclusion 
\[
L(\omega_i ,P_i )^{K^{\mf s} \times K^{\mf s}} \to
L(\omega_i ,P_i )^{K \times K}
\]
is of the type $M_n (\mh C ) \subset M_m (\mh C )$. Therefore
\begin{equation}\label{eq:4.84}
\mc S (G,K_{\mf s} )^{\mf s} \to
\mc S (G,K)^{\mf s}
\end{equation}
induces an isomorphism on $HH_* ,\, HC_* ,\, HP_*$ and $K_*$.
From this, \eqref{eq:4.2} and the properties of the topological
$K$-functor (cf. page \pageref{it:2.4}) we see that
\begin{equation}\label{eq:4.82}
\begin{aligned}
K_* (C_r^* (G)) \, &\cong K_* \left( \varinjlim_{\mf S} 
\bigoplus_{\mf s \in \mf S} C_r^* (G)^{\mf s} \right) \\
& \cong \bigoplus_{\mf s \in \mf B (G)} K_* \left( 
C_r^* (G)^{\mf s} \right) \\
& \cong \bigoplus_{\mf s \in \mf B (G)} \varinjlim K_* 
\left( C_r^* (G,K)^{\mf s} \right) \\
& \cong \bigoplus_{\mf s \in \mf B (G)} 
K_* \big( C_r^* (G,K_{\mf s} )^{\mf s} \big)\\
&\cong \bigoplus_{\mf s \in \mf B (G)} 
K_* \big( \mc S (G,K_{\mf s} )^{\mf s} \big)
\end{aligned}
\end{equation}
Since the algebras $\mc S (G,K)^{\mf s}$ are all nuclear
Fr\'echet, and have the same Hochschild homology for fixed
$\mf s \in \mf B (G)$ and $K \subset K_{\mf s}$, we may use 
the continuity and additivity of $HP_* ( \, \cdot \,, 
\oot)$, as described on page \pageref{p:hp}.
(The cohomological dimension of $\mc H (G)^{\mf s}$ and
$\mc S (G)^{\mf s}$ is bounded independently of $\mf s$.)
\begin{equation}\label{eq:4.83}
\begin{aligned}
HP_* (\mc S (G), \oot) \: &\cong 
\bigoplus_{\mf s \in \mf B (G)} 
HP_* (\mc S (G)^{\mf s}, \oot)\\
&\cong \bigoplus_{\mf s \in \mf B (G)} \varinjlim
HP_* (\mc S (G,K)^{\mf s}, \oot)\\
&\cong \bigoplus_{\mf s \in \mf B (G)} 
HP_* (\mc S (G, K_{\mf s})^{\mf s}, \oot)\\
&= \bigoplus_{\mf s \in \mf B (G)} 
HP_* (\mc S (G, K_{\mf s})^{\mf s}, \hot)\\
\end{aligned}
\end{equation}
Now the theorem follows from the combination of \eqref{eq:4.81}, 
\eqref{eq:4.82} and \eqref{eq:4.83}. $\qquad \Box$
\\[3mm]

The analogue of Theorem \ref{thm:3.38} was suggested in
\cite[Conjecture 8.9]{BHP2} and in \cite[Conjecture 1]{ABP} :

\begin{thm}\label{thm:4.44}
The inclusions $\mc H (G)^{\mf s} \to \mc S (G)^{\mf s}$ induce
isomorphisms
\begin{align*}
HP_* ( \mc H (G)^{\mf s} ) &\isom HP_* ( \mc S (G)^{\mf s},
\oot) \\
HP_* ( \mc H (G) ) &\isom HP_* ( \mc S (G),\oot) \\
\end{align*}
\end{thm}
\emph{Proof.}
Just as in \eqref{eq:4.83} we have 
\begin{equation}\label{eq:4.85}
\begin{aligned}
HP_* (\mc H (G) ) \: &\cong 
\bigoplus_{\mf s \in \mf B (G)} 
HP_* (\mc H (G)^{\mf s} )\\
&\cong \bigoplus_{\mf s \in \mf B (G)} \varinjlim
HP_* (\mc H (G,K)^{\mf s} )\\
&\cong \bigoplus_{\mf s \in \mf B (G)} 
HP_* (\mc H (G, K_{\mf s})^{\mf s} )\\
\end{aligned}
\end{equation}
Therefore we only have to show that every inclusion
\[
\mc H (G, K_{\mf s})^{\mf s} \to  
\mc S (G, K_{\mf s})^{\mf s} 
\]
induces an isomorphism on periodic cyclic homology. Number 
the direct summands in \eqref{eq:4.79}, such that
\begin{equation}\label{eq:4.87}
\mc S (G,K_{\mf s} )^{\mf s} \cong \bigoplus_{j=1}^{n_{\mf s}} 
\big( C^\infty (X_{unr}(M_j )) \otimes L(\omega_j ,P_j 
)^{K_{\mf s} \times K_{\mf s}} \big)^{\mc W_j} 
\end{equation}
and $\dim M_i \leq \dim M_j$ if $i \leq j$. Now we construct 
two chains of ideals
\begin{align*}
&\mc H (G, K_{\mf s})^{\mf s} = I_0 \supset I_1 \supset \cdots 
\supset I_{n_{\mf s}} = 0 \\ 
&\mc S (G, K_{\mf s})^{\mf s} = J_0 \supset J_1 \supset \cdots 
\supset J_{n_{\mf s}} = 0 \\
&I_j = \{ h \in \mc H (G, K_{\mf s})^{\mf s} : I(P_j ,A_j ,\omega_j 
,\chi )(h) = 0 \;\mr{if}\: \chi \in X_{nr}(M_j) 
\;\mr{and}\; j \leq i \} \\
&J_j = \{ h \in \mc S (G, K_{\mf s})^{\mf s} : I(P_j ,A_j ,\omega_j 
,\chi )(h) = 0 \;\mr{if}\: \chi \in X_{unr}(M_j) 
\;\mr{and}\; j \leq i \}
\end{align*}
Writing $V_i = I^G_{P_i} (E_i )^{K_{\mf s}}$, we clearly have
\[
J_{i-1} / J_i \cong \big( C^\infty (X_{unr}(M_i )) \otimes 
\mr{End} \: V_i \big)^{\mc W_i} 
\]
From now on we can follow the proof of Theorem \ref{thm:3.38}. 
We have to substitute Theorems \ref{thm:4.29}, \ref{thm:4.30},
\ref{thm:4.32} and \ref{thm:4.37} for, respectively, Theorems
\ref{thm:3.11}, \ref{thm:3.17}, \ref{thm:3.19} and \ref{thm:3.24}.
$\qquad \Box$. 
\\[2mm]

This theorem should be compared with the work of Meyer \cite{Mey}.

It is also interesting to compare Theorems \ref{thm:4.43} and 
\ref{thm:4.44} with other homological invariants of reductive
$p$-adic groups. One such is the chamber homology of the 
Bruhat-Tits building $\beta G$, equivariant with respect to 
$G$. This is a sequence of complex vector spaces \inde{$H^G_n 
(\beta G)$} $ ,n = 0,1,2, \ldots$, first defined in \cite[\S 
6]{BCH}. It was proved simultaneously in \cite{HiNi} and 
\cite{Schn} that there are natural isomorphisms 
\begin{equation}\label{eq:4.88} 
HP_i (\mc H (G)) \cong \bigoplus_{n \in \mh Z} H_{i+2n}^G (\beta G)
\end{equation}
Closely related is the $G$-equivariant $K$-homology of 
$\beta G$, as defined in \cite[\S 3]{BCH}. According to 
Voigt \cite[Theorem 6.8]{Voi} there exists a natural equivariant 
Chern character
\begin{equation}
ch_*^G : K_*^G (\beta G) \to H_*^G (\beta G)
\end{equation} 
which becomes an isomorphism after tensoring with $\mh C$. 
Notice that Voigt uses an alternative but equivalent definition
of $H_*^G (\beta G)$. The \inde{Baum-Connes conjecture} for 
reductive $p$-adic groups asserts that the so-called 
\inde{assembly map} \index{$K_j^G (\beta G)$}
\begin{equation}\label{eq:4.89}
\mu : K_j^G (\beta G) \to K_j (C_r^* (G))
\end{equation}
is an isomorphism. This was proved by Lafforgue \cite{Laf}, as a 
part of a much more general result. Putting all these things 
together we more or less arrive at \cite[Proposition 9.4]{BHP2} :

\begin{thm}\label{thm:4.45}
In following diagram the horizontal maps are natural
isomorphisms, and the vertical maps become natural isomorphisms 
after tensoring with $\mh C$.
\[
\begin{array}{ccccc}
K_*^G (\beta G) & & \longrightarrow & & K_* (C_r^* (G)) \\
\downarrow & & & & \downarrow \\
H_*^G (\beta G) & \cong & HP_* (\mc H (G)) & \to & 
  HP_* (\mc S (G), \oot)
\end{array}
\]
\end{thm}

Because of the naturality, the diagram is probably commutative, but 
the author does not know how to prove this.
Unfortunately the definitions are so complicated that it already is 
difficult to find any element of $K^G_* (\beta G)$ for which the
diagram can be seen to commute by direct computation. 

For the groups $GL_n (\mh F)$ and $SL_n (\mh F)$ partial results
in the direction of Theorem \ref{thm:4.45} were proved by Baum, 
Higson and Plymen in \cite{BHP1}. In fact, in \cite{BHP1} the 
Baum-Connes conjecture for these groups is proved precisely with 
the above diagram. However, the argument uses the commutativity of 
the diagram in an essential way, and unfortunately the authors do 
not provide any support of the (implicit) claim that it does commute.

%% file: chapter5.tex
\chapter{Parameter deformations in affine Hecke algebras}

So far we have always written affine Hecke algebras as 
deformations of a group algebra, but we have not really
done anything with this. Ideally speaking, several 
properties of an affine Hecke algebra $\mc H (\mc R ,q)$
should be independent of the parameters $q(s)$. This
intuitive idea comes from finite dimensional algebras, where
it is very clear from Tits' deformation theorem. Roughly
speaking, it tells us that if two semisimple algebras can be
continuously deformed into each other, then they are
isomorphic. This is so because there are only countably
many isomorphism classes of such algebras, and they lie
discrete in some sense. For infinite dimensional algebras
nothing similar holds, so there we have to find more
subtle invariants and arguments.

This has been done for affine Hecke algebras with equal
labels. Kazhdan and Lusztig \cite{KaLu} gave a complete
geometric parametrization of the irreducible representations 
of such algebras. This parametrization is independent of 
$q \in \mh C^\times$, except in a few tricky cases where $q$ 
is a proper root of unity. Baum and Nistor \cite{BaNi} 
showed that this leads to an isomorphism
\begin{equation}\label{eq:5.41}
HP_* \big( \mc H (\mc R ,q) \big) \cong 
HP_* \big( \mh C [W] \big)
\end{equation}
Acknowledging that these results cannot readily be carried 
over to the unequal label case, we follow another path, more 
analytic in nature. By careful estimates in the Schwartz algebra 
$\mc S (\mc R ,q)$ we show that the following things all depend 
continuously on $q$: the operator norm, multiplication, 
inverting and the holomorphic functional calculus.

Equipped with these tools and the knowledge from Chapter 3 we
attack a special kind of parameter deformation, scaling the
label function. Such deformations were studied first by Opdam
\cite{Opd3}. On the level of central characters this amounts to 
scaling the absolute value by a real factor $\ep$.
For $\ep > 0$ we construct isomorphisms of pre-$C^*$-algebras
\begin{equation}\label{eq:5.42}
\phi_\ep : \mc S (\mc R ,q^\ep ) \isom \mc S (\mc R ,q)
\end{equation}
They depend continuously on $\ep$, but the limit
\begin{equation}\label{eq:5.43}
\phi_0 : \mc S (W) \to \mc S (\mc R ,q)
\end{equation}
is no longer surjective. Nevertheless this map seems to behave
well. In view of \eqref{eq:5.41} one is naturally led 
to conjecture that
\begin{equation}\label{eq:5.44}
HP_* (\phi_0 ) : HP_* \big( \mc S (W) \big) \to
HP_* \big( \mc S (\mc R ,q) \big) 
\end{equation}
is an isomorphism for any positive label function $q$. We provide
various equivalent reformulations of this statement.

Conjecture \eqref{eq:5.44} can be derived from the stronger 
conjecture, originally due to Baum, Connes and Higson \cite{BCH}, 
that
\begin{equation}\label{eq:5.45}
K_* \big( C_r^* (\mc R ,q) \big) \cong K_* \big( C_r^* (W) \big)
\end{equation}
We show that \eqref{eq:5.45} is equivalent to the existence of a 
natural bijection between the Grothendieck groups of irreducible
representations of $C_r^* (W)$ and of $C_r^* (\mc R ,q)$. At the end 
of the chapter we give some clues in support of these conjectures.
\\[4mm]

\section{The finite dimensional and equal label cases}
\label{sec:5.1}

We recall what is already known about deformations of
Iwahori-Hecke algebras obtained by varying the label
function $q$. For Hecke algebras of finite type this is 
very clear: as long as they are semisimple they are rigid 
under deformations. But this is specific for the finite 
case, as it relies on the classification of finite 
dimensional semisimple algebras.

For any extended Iwahori-Hecke algebra $\mc H (\mc R ,q)$ 
with equal labels a complete parametrization of irreducible
representations is available. This is a refinement of the 
Langlands classification, and it is essentially independent 
of $q$. The link between different $q$'s is made via Lusztig's
asymptotic Hecke algebra $J$, which allows a weakly 
spectrum preserving morphism $\mc H (\mc R ,q) \to J$.
From this we will see that the periodic cyclic homology of
$\mc H (\mc R ,q)$ is independent of $q$, as long as it is
not a proper root of unity.

Recall that an algebra $A$ is semisimple if its Jacobson 
radical is 0, which means that for every nonzero $a \in A$
there is an irreducible $A$-representation $\pi$ such that
$\pi (a) \neq 0$. For example, by \cite[Th\'eor\`eme 
2.7.3]{Dix} every $C^*$-algebra is semisimple. 
The structure of finite dimensional semisimple algebras is
described in a famous theorem of Wedderburn \cite{Wed} :

\begin{thm}\label{thm:5.1}
Let $A$ be a finite dimensional \inde{semisimple algebra} 
over a field $\mh F$. There exist natural numbers $n_i$ 
and division algebras $D_i$ over $\mh F$ such that
\[
A \cong \bigoplus_{i=1}^r M_{n_i} (D_i)
\]
If $\mh F$ is algebraically closed then 
$D_i = \mh F \; \forall i$.
\end{thm}

Let $G$ be any finite group. By Maschke's theorem the 
group algebra $\mh C [G]$ is semisimple. Let 
$\{T_g : g \in G\}$ be its canonical basis, and 
$\mb k = \mh C [x_1 ,\ldots, x_r ]$ a polynomial ring
over $\mh C$. Let $A$ be a $\mb k$-algebra whose 
underlying $\mb k$-module is $\mb k [G]$ and whose 
multiplication is defined by 
\begin{equation}\label{eq:5.1}
T_g \cdot T_h = \sum_{w \in G} a_{g,h,w} T_w
\end{equation}
for certain $a_{g,h,w} \in \mb k$. For any point 
$q \in \mh C^r$ we can endow the vector space $\mh C [G]$
with the structure of an associative algebra by
\begin{equation}\label{eq:5.5}
T_g \cdot_q T_h = \sum_{w \in G} a_{g,h,w}(q) T_w
\end{equation}
We denote the resulting algebra by $\mc H (G,q)$. It is 
isomorphic to the tensor product $A \otimes_{\mb k} \mh C$
where $\mh C$ has the $\mb k$-module structure obtained 
from evaluating at $q$. Assume moreover that there exists a
$q^0 \in \mh C^r$ such that 
\[
\mc H (G ,q^0 ) = \mh C [G]
\]
We express the rigidity of finite dimensional semisimple
algebras by the following special case of \inde{Tits's 
deformation theorem} \cite[p. 357 - 359]{Car1}:

\begin{thm}\label{thm:5.2}
There exists a polynomial $P \in \mb k$ such that the 
following are equivalent :
\begin{itemize}
\item $P(q) \neq 0$
\item $\mc H (G,q)$ is semisimple
\item $\mc H (G,q) \cong \mh C [G]$
\end{itemize}
\end{thm}

Now let $(W,S)$ be a finite Coxeter system, $q$ a label 
function om $W$ and $\mc H (W,q)$ the associated Iwahori-Hecke
algebra, as in Section \ref{sec:3.1}. This is consistent with
the above notation. We want to know under which conditions 
this algebra is semisimple. Clearly this is the case if 
$q(w) > 0 \; \forall w \in W$, for then $\mc H (W,q)$ is a 
$C^*$-algebra by \eqref{eq:3.inp}. 

But the polynomials $P(q)$ of Theorem \ref{thm:5.2} have also
been determined explicitly. If we are in the equal label case 
$q(s) = q \; \forall s \in S$ then we may take
\begin{equation}
P(q) = q \sum_{w \in W} q^{\ell (w)} 
\end{equation}
except that we must omit the factor $q$ if $W$ is of type 
$(A_1 )^n$, see \cite{GyUn}. More generally, Gyoja \cite[p. 
569]{Gyo} showed that if $(W,S)$ is irreducible and $S$ consists 
of two conjugacy classes, then in most cases we may take
\begin{equation}\label{eq:5.2}
\begin{array}{ccl}
P(q_1 ,q_2 ) & = & q_1^{|W|} q_2 W(q_1 ,q_2 ) W(q_1^{-1},q_2 )\\
W(q_1 ,q_2 ) & = & \sum\limits_{w \in W} q(w)
\end{array}
\end{equation}
So generically there is an isomorphism 
\begin{equation}\label{eq:5.3}
\mc H (W,q) \cong \mh C [W]
\end{equation}
We will see later how it can be constructed explicitly. From 
our somewhat simple point of view this is all there is to 
say about parameter deformations of finite dimensional Hecke 
algebras. If they are semisimple then they are isomorphic, 
and if not, then they have nilpotent ideals and look very 
different from $\mh C [W]$.
\\[2mm]

Let $\mc R = (X,Y,R_0 ,R_0^\vee ,F_0 )$ be a root datum, let
$q \in \mh C^\times$, and consider the affine Hecke algebra
with equal labels $\mc H (\mc R ,q)$. The irreducible 
representations of this algebra have been classified 
completely by Kazhdan and Lusztig \cite{KaLu}. For this very
deep result they showed among others that $\mc H (\mc R ,q)$
is isomorphic to the equivariant algebraic $K$-theory of a
certain variety.

Let $G$ be the unique complex reductive algebraic group with
root datum $\mc R^\vee = (Y,X,R_0^\vee ,R_0 )$, and $\mf g$
its Lie algebra. For reasons of a much more general nature
$G$ is called the \inde{Langlands dual group}.
Then $T = \mr{Hom}_{\mh Z} (X, \mh C^\times )$ can be 
identified with a maximal torus of $G$. Since every semisimple
element of $G$ is conjugate to an element of $T$, and since
$N_G (T) / Z_G (T) \cong W_0$, we can parametrize the central
character of an (irreducible) $\mc H (\mc R ,q)$-module by a
unique conjugacy class of semisimple elements in $G$. So, let
$s \in G$ be semisimple and write
\begin{equation}
\mf n (s,q) = \{ N \in \mf g : N \:\mr{nilpotent} \,, 
\mr{Ad}(s) N = q N \}
\end{equation}
The $G$-conjugacy classes of pairs $(s,N)$ with 
$N \in \mf n (s,q)$ are called Deligne-Langlands parameters. 
They give an almost complete description of
Prim$(\mc H (\mc R ,q))$. For instance it works perfectly if 
$\mc R$ is of type $GL_n$, see \cite{Zel}. In a sense this is
equivalent  to the Langlands classification in Theorem 
\ref{thm:4.29}. To find really all irreducible representations 
we must add one extra ingredient. Let
\begin{equation}\label{eq:5.38}
Z (s,N) = \{ g \in G : g s = s g ,\, \mr{Ad}(g)N = N \}
\end{equation}
be the simultaneous centralizer of $s$ and $N$, and $Z^0 (s,N)$
its identity component. Assume that $q \in \mh C^\times $ is 
not a \inde{proper root of unity}, i.e. either $q=1$ or $q$ is 
not a root of unity. In these cases there is a bijection between 
Prim$(\mc H (\mc R ,q))$ and $G$-conjugacy classes of triples
$(s,N,\rho )$, where $s \in G$ is semisimple, $N \in \mf n 
(s,q)$ and $\rho$ is a "geometric" irreducible representation 
of the finite group $Z (s,N) / Z^0 (s,N)$. This was proved 
for $q=1$ in \cite[Theorem 4.1]{Kat2} and for $q$ not a root of
unity and $X$ equal to the weight lattice of $R_0^\vee$ in 
\cite[Theorem 7.12]{KaLu}. Later it was shown in \cite[Theorem 
2]{Ree} that this condition on $X$ is not necessary.

Another construction which is particular for the equal label 
case is Lusztig's asymptotic Hecke algebra \cite{Lus2,Lus3}. 
This is a finite type algebra $J$ with a basis $\{ t_w : w \in W 
\}$ over $\mh C$. It decomposes as a finite direct sum of 
two-sided ideals \index{Hecke algebra!asymptotic}
\begin{equation}\label{eq:5.4}
\begin{array}{lll}
J   & = & \bigoplus_{i=0}^{|R_0^+|} J^i\\
J^i & = & \mr{span} \{ t_w : a(w) = i \}
\end{array}
\end{equation}
where $a$ is Lusztig's $a$-function. 
For every $q \in \mh C^\times$ there is an 
injective morphism of finite type algebras
\begin{equation}
\phi_q : \mc H (\mc R ,q) \to J
\end{equation}
If $q$ is not a proper root of unity, then $\phi_q$ induces 
a bijection on irreducible representations. Namely, for any 
irreducible $J$-representation $\pi$ the $\mc H (\mc R 
,q)$-representation $\pi \circ \phi_q$ has a unique irreducible
constituent of minimal "$a$-weight". This implies that the 
morphisms of finite type algebras 
\begin{equation}
\phi_q^{-1} \Big( \bigoplus_{i \geq k} J^i \Big) \Big/ 
\phi_q^{-1} \Big( \bigoplus_{i > k} J^i \Big) 
\; \longrightarrow \; J^k
\end{equation}
are spectrum preserving. Lusztig \cite[Corollary 3.6]{Lus3} 
proved this in the case $W = W_{\mr{aff}}$, but using the 
aforementioned result of Reeder \cite{Ree} his proof can be 
extended to general root data. Combining this with Theorem 
\ref{thm:2.4} and Lemma \ref{lem:2.12} we arrive at an 
extended version of \cite[Theorem 11]{BaNi} :

\begin{thm}\label{thm:5.4}
Assume that $q$ is not a proper root of unity. Then
\[
HP_* (\phi_q ) : HP_* (\mc H (\mc R ,q)) \to HP_* (J)
\]
is an isomorphism. Consequently
\[
HP_* (\mc H (\mc R ,q)) \cong HP_* ( \mh C [W] )
\]
\end{thm}

It is expected that an asymptotic Hecke algebra can also be 
constructed for finite or affine Coxeter systems with unequal 
labels \cite[Chapter 18]{Lus6}. Assuming certain conjectures 
\cite[Chapter 15]{Lus6} one can construct algebra 
homomorphisms $\phi_q : \mc H (W ,q) \to J$ for any label 
function with the following property: there exist 
$v \in \mh C^\times$ and $n_s \in \mh N$ such that
$q(s) = v^{n_s} \: \forall s \in S$.

For finite $W$ the map $\phi_q$ is an isomorphism if and only 
if $\mc H (W,q)$ is semisimple \cite[(20.1.e)]{Lus6}. In this 
way one can find explicit formulas for the isomorphisms from 
Theorem \ref{thm:5.2}.

For affine $W \;\phi_q$ has a nilpotent kernel \cite[Proposition
18.12]{Lus6} and in general it is not surjective. It is unknown
whether $\phi_q$ is spectrum preserving in any sense. The problem 
is that in general there is no definite classification of all
irreducible representations of an affine Hecke algebra. 
Apparently the link with the Langlands dual group is much weaker
for unequal parameters.

Nevertheless in some cases the Deligne-Langlands philosophy
outlined above can be generalized. Namely, along these lines a 
classification of the irreducible representations of 
$\mc H (\mc R ,q)$ has been obtained in \cite{Kat4} for $\mc R$ 
of type $B_n / C_n$, for almost all label functions $q$.

\section{Estimating norms}
\label{sec:5.2}

In this section we lay the analytic foundations for all our
coming results on parameter deformations.
In Section \ref{sec:3.2} we defined various norms on an affine
Hecke algebra $\mc H (\mc R ,q)$: the norm $\norm{\:\cdot
\:}_\tau$ associated with the trace $\tau$, the operator norm 
$\norm{\:\cdot\:}_o$ and the Schwartz norms $p_n$. We will 
show that the operator norm, the multiplication and the 
inverse of an element depend continuously on $q$. From this
we deduce that the holomorphic functional calculus on affine Hecke
algebras is continuous in a very general sense.

We also reconstruct the Schwartz algebra $\mc S (\mc R ,q)$ 
in a different way. With this construction we can show in a
straightforward fashion that $\mc S (\mc R ,q)$ is 
holomorphically closed in $C_r^* (\mc R ,q)$. 

With respect to the bases $\{ N_w : w \in W \}$ the norms 
$\norm{\cdot}_\tau$ and $p_n$ are independent of $q$. Therefore 
we can identify all the Hilbert spaces $\mf H (\mc R ,q)$ and all 
the Schwartz spaces $\mc S (\mc R ,q)$ by means of this basis. 
When we want to consider them in this way, only as topological 
vector spaces and without a specified label function, we write
\inde{$\mf H (\mc R )$} and \inde{$\mc S (\mc R )$}. To 
indicate that $x \in \mc S (\mc R )$ should be considered as 
an element of $\mc S (\mc R ,q)$ we sometimes denote it by 
\inde{$(x,q)$}. Furthermore, to distinguish the various products 
we add a subscript, so \inde{$\cdot_q$} is the multiplication in 
$C_r^* (\mc R ,q)$.

Let \inde{$L_{\mc R}$} be the space of label functions on 
$\mc R$ satisfying the positivity Condition \ref{cond:3.3}. 
Recall that for a simple reflection $s_i \in S_{\mr{aff}}$ 
we put \index{etai@$\eta_i$} \index{$q_i$}
\begin{equation}
q_i = q (s_i ) \quad \mr{and} \quad 
\eta_i = q_i^{1/2} - q_i^{-1/2}
\end{equation}
These numbers determine $q$ uniquely, and their domain is only 
limited by the conditions $q_i > 0$ and $q_i = q_j$ whenever
$s_i$ and $s_j$ are conjugate in $W$. Hence the parameter space
$L_{\mc R}$ is homeomorphic to $\mh R^n$ for a certain $n$. We 
will use the standard topology on $L_{\mc R}$, induced by the
metric \index{rho@$\rho$}
\begin{equation}
\rho (q,q') = \max_{s_i \in S_{\mr{aff}}} | \eta_i - \eta'_i |
\end{equation}
We already know that the group $W$ with the length function 
$\mc N$ is of polynomial growth, but we need a more explicit
estimate on the number of elements of a fixed length.

\begin{lem}\label{lem:5.3}
There exists a real number $C_{\mc N}$ such that 
$\forall n \in \mh N$
\[
\# \big\{ w \in W : \mc N (w) = n \big\} < 
C_{\mc N} \, (n+1)^{\mr{rk}(X) -1}
\]
\end{lem}
\emph{Proof.} 
Let $r$ denote the rank of $X$, and pick a linear bijection
$L : X \otimes \mh R \to \mh R^r$ such that
\begin{itemize}
\item $L (X) \subset \mh Z^r$
\item $\forall x \in X^+ : \mc N (x) = \norm{L(x)}_1$
\end{itemize}
This is possible because, on $X^+, \mc N$ is additive and
takes values in $\mh N$. By \eqref{eq:3.16} we can write any
$w \in W$ as 
\[
w = u x v \quad \mr{with} \quad u,v \in W_0 \,, x \in X^+ 
\]
If $\mc N (w) = n$, then clearly 
\begin{equation}
n - |R_0 | \leq n - \mc N (u) - \mc N (v) \leq \mc N (x) \leq
n + \mc N (u) + \mc N (v) \leq n + |R_0 |
\end{equation}
Therefore we can estimate
\begin{align*}
|W_0 |^{-2} \# \big\{ w \in W : \mc N (w) = n \big\} 
& \leq \# \big\{ x \in X^+ : n - |R_0| \leq \mc N(x) 
\leq n + |R_0 | \big\} \\
& = \# \big\{ y \in L(X^+ ) : n - |R_0| \leq \norm{y}_1 
\leq n + |R_0 | \big\} \\
& \leq \# \big\{ y \in \mh Z^r : n - |R_0| \leq \norm{y}_1 
\leq n + |R_0 | \big\}
\end{align*}
For the sake of calculation we assume now that $n > |R_0 |$. 
This is allowed because there are only finitely many 
$w \in W$ of smaller length. We continue our estimate:
\begin{align*}
& \leq (n + |R_0 | + 1)^r - (n - |R_0 | - 1)^r\\
& = \sum_{i = 0}^r \left( r \atop i \right) n^{r-i} 
\big( |R_0 | + 1 \big)^i \big( 1 - (-1)^i \big)\\
& \leq n^{r-1} \sum_{i \leq r ,\, i \; \mr{odd}} 2 
\left( r \atop i \right) \big( |R_0 | + 1 \big)^i \\
& < n^{r-1} 2 \big( |R_0 | + 2 \big)^r
\end{align*}
So our candidate for $C_{\mc N}$ is 
\[
|W_0 |^2 2 \big( |R_0 | + 2 \big)^r
\]
We only have to check whether it works also for 
$n \leq |R_0|$ and, if not, increase it accordingly.
$\qquad \Box$ \\[2mm]

Put \inde{$b$} $= \mr{rk}(X) + 1$. By Lemma \ref{lem:5.3} the 
following sum converges to a limit \inde{$C_b$}:
\begin{equation}
C_b := \sum_{w \in W} (\mc N (w) + 1)^{-b} < 
\sum_{n=0}^\infty C_{\mc N} (n+1)^{\mr{rk}(X) - 1} 
(n+1)^{-\mr{rk}(X) - 1} < \infty
\end{equation}
This implies that for any $x = \sum_u x_u N_u \in \mc S 
(\mc R)$ and $n \in \mh N$
\begin{align}\label{eq:5.6}
& \sum_u |x_u| (\mc N (u) + 1)^n \leq \sum_u \sup_v \big\{ 
|x_v| (\mc N (v) + 1)^{n+b} \big\} (\mc N (u) + 1 )^{-b} =
C_b p_{n+b} (x) \\
& \norm{x}_\tau \leq \sum_u |x_u | \leq C_b p_b (x)
\end{align}
Since $\ell (w) \leq \mc N (w) \; \forall w \in W$, these 
inequalities a fortiori remain valid if we replace $\mc N$
by $\ell$.

Let $u,v,w \in W$ and let $u = \omega s_1 \cdots s_{\ell (u)}$
be a reduced expression, where $\ell (\omega ) = 0$ and
$s_i \in S_{\mr{aff}}$. (The $s_i$ need not all be different.)
For $I \subset \{ 1,2, \ldots ,\ell (u) \}$ we put 
$\eta_I = \prod_{i \in I} \eta_i$ and 
\index{etaI@$\eta_I$} \index{$u_I$}
\[
u_I = \omega \tilde s_1 \cdots \tilde s_{\ell (u)}
\qquad \mr{where} \qquad
\tilde s_i = \left\{ \begin{array}{ccc}
s_i & \mr{if} & i \notin I\\
e   & \mr{if} & i \in I
\end{array} \right.
\]
\begin{thm}\label{thm:5.5}
\[
N_u \cdot_q N_v = \sum_{I \subset \{1,2, \ldots ,\ell (u)\}}
\eta_I D_v^u (I) N_{u_I v} \qquad \mr{where}
\]
\begin{itemize}
\item \inde{$D_v^u (I)$} is either 0 or 1
\item $D_u^v (\es) = 1$ and $D_v^u (I) = 0$ if 
$|I| > |R_0^+ |$
\item $\sum\limits_{I \subset \{ 1,2, \ldots ,\ell (u) \} } 
D_v^u (I) < 3 (\ell (u) + 1)^{|R_0^+ |}$
\end{itemize}
\end{thm}
\emph{Proof.}
It follows from the multiplicaton rules \eqref{eq:3.20} that
\begin{equation}
N_{s_i} \cdot_q N_v = N_{s_i v} + D_v^{s_i} (i) \eta_i N_v 
\;\; \mr{where} \;\; D_v^{s_i} (i) = 
\left\{ \begin{array}{ccc}
0 & \mr{if} & \ell (s_i v) > \ell (v) \\
1 & \mr{if} & \ell (s_i v) < \ell (v)
\end{array} \right.
\end{equation}
The expression for $N_u \cdot_q N_v$, with $D_u^v (I)$ being 
0 or 1 and $D_u^v (\es) = 1$, follows from this, with 
induction to $\ell (u)$. By \cite[Theorem 7.2]{Lus1} for 
fixed $w \in W$ the sum
\[
\sum_{I : u_I = w} \eta_I D_v^u (I)
\] 
is a polynomial of degree at most $|R_0^+ |$ in the $\eta_i$. 
Therefore $D_v^u (I) = 0$ whenever $|I| > |R_0^+ |$. Consequently
\begin{multline}\label{eq:5.7}
\sum\limits_{I \subset \{ 1,2, \ldots ,\ell (u) \}} D_v^u (I) 
\;\; \leq \;\; \# \big\{ I \subset \{ 1,2, \ldots ,\ell (u) \} 
: |I| \leq |R_0^+ | \big\} \;\; \leq \;\; \sum\limits_{j=0
}^{|R_0^+ |} \left( \ds{\ell (u) \atop j} \right) \\
\leq \;\; {\ds \frac{\ell (u)!}{\big( \ell (u) - |R_0^+ | \big) !}} 
\sum\limits_{j=0}^{|R_0^+ |} {\ds \frac{1}{j!}} \;\; < \;\; 
3 (\ell (u) + 1)^{|R_0^+ |}
\end{multline}
where we should interpret $\big( \ell (u) - |R_0^+ | \big) !$ as 1 if
$|R_0^+ | \geq \ell (u). \qquad \Box$
\\[3mm]

Let $\eta > 0$ and put \inde{$C_\eta$} $= 3 C_b \max 
\big\{ 1,\eta^{|R_0^+ |} \big\}$.

\begin{prop}\label{prop:5.6}
For all $q,q' \in B_\rho (q^0 ,\eta ) ,\, x \in \mc S (\mc R)$
the following estimates hold.
\[
\begin{array}{lrr}
\norm{\lambda (x,q)}_{B (\mf H (\mc R))} \;=\; \norm{(x,q)}_o &
\leq & C_\eta p_{b + |R_0^+ |} (x) \\
\norm{\lambda (x,q) - \lambda (x,q')}_{B (\mf H (\mc R))} &
\leq & \rho (q,q') C_\eta p_{b + |R_0^+ |} (x) 
\end{array}
\]
In particular $\mc S (\mc R, q) \subset C_r^* (\mc R ,q)$ and 
every finite dimensional\\ $C_r^* (\mc R ,q)$-representation is
tempered.
\end{prop}
\emph{Proof.}
Let $y = \sum_v y_v N_v \in \mc S (\mc R)$. By \eqref{eq:5.6}
and Theorem \ref{thm:5.5} we have
\begin{equation}
\begin{aligned}
\norm{x \cdot_q y}_\tau & = \norm{\sum_{u,v} x_u y_v N_u 
\cdot_q N_v }_\tau \\
& = \norm{\sum_{u,v} x_u y_v \sum_I \eta_I D_v^u (I) 
N_{u_I v} }_\tau \\
& \leq \sum_u |x_u | \sum_{I : |I| \leq |R_0^+ |} |\eta_I |
\norm{\sum_v |y_v | N_{u_I v} }_\tau \\
& \leq \sum_u |x_u | (\ell (u) + 1)^{|R_0^+ |} 
3 \max \big\{ 1,\eta^{|R_0^+ |} \big\} \norm{y}_\tau \\
& \leq C_\eta p_{b + |R_0^+ |} (x) \norm{y}_\tau
\end{aligned}
\end{equation}
Since $\mc S (\mc R)$ is dense in $\mf H (\mc R)$ this gives 
the estimate, by the very definition of the operator norm on
$B(\mf H (\mc R))$. In particular we get a continuous 
embedding $\mc S (\mc R ,q) \to C_r^* (\mc R ,q)$, so every 
finite dimensional representation of the latter algebra is 
tempered by Lemma \ref{lem:3.9}.
\begin{equation}
\begin{aligned}
\norm{x \cdot_q y - x \cdot_{q'} y}_\tau & = \norm{\sum_{u,v} 
x_u y_v (N_u \cdot_q N_v - N_u \cdot_{q'} N_v) }\tau \\
& = \norm{\sum_{u,v} x_u y_v \sum_I (\eta_I -\eta'_I ) D_v^u (I) 
N_{u_I v} }_\tau \\
& \leq \norm{\sum_{u,v} x_u y_v \sum_I \rho (q,q') |I| 
\eta^{|I| - 1} D_v^u (I) N_{u_I v} }_\tau \\
& \leq \rho (q,q') \sum_u |x_u | \sum_{I : |I| \leq |R_0^+ |}
|I| \eta^{|I| - 1} \norm{\sum_v |y_v | N_{u_I v} }_\tau \\
& \leq \rho (q,q') \sum_u |x_u | (\ell (u) + 1)^{|R_0^+ |} 
\,3 \max \big\{ 1,\eta^{|R_0^+ |} \big\} \norm{y}_\tau \\
& \leq \rho (q,q') C_\eta p_{b + |R_0^+ |} (x) \norm{y}_\tau
\end{aligned}
\end{equation}
Between lines 4 and 5 we used a small calculation like 
\eqref{eq:5.7} :
\begin{equation}
\begin{aligned}
\sum_{I : |I| \leq |R_0^+ |} |I| \eta^{|I| - 1} & \leq
\sum_{j=0}^{|R_0^+ |} \left( \ds{\ell (u) \atop j} \right) 
j \eta^{j-1} \\
& \leq {\ds \frac{\ell (u)!}{\big( \ell (u) - |R_0^+ | \big) !}} 
\sum\limits_{j=0}^{|R_0^+ |} {\ds \frac{j}{j!}} 
\max \big\{ 1,\eta^{|R_0^+ | - 1} \big\} \\
& \leq (\ell (u) + 1)^{|R_0^+ |} 
\,3 \max \big\{ 1,\eta^{|R_0^+ |} \big\}
\end{aligned}
\end{equation}
Note that we did not really use that $y$ lies in the subspace
$\mc S (\mc R)$ of $\mf H (\mc R)$, it only helps to ensure 
that all intermediate expressions are well-defined. 
$\qquad \Box$ \\[3mm]

Now we can also estimate the behaviour of the Schwartz norms
$p_n$ under multiplication. Put 
\inde{$b'$} $= 2b + |R_0^+ | = 2 \:\mr{rk}(X) + |R_0^+ | + 2$.

\begin{prop}\label{prop:5.7}
Let $n \in \mh N ,\, q,q' \in B_\rho (q^0 ,\eta )$ and
$x_i = \sum_{u \in W} x_{iu} N_u \in \mc S (\mc R ,q)$. Then
\[
\begin{array}{lcr}
p_n (x_1 \cdot_q \cdots \cdot_q x_m ) & \leq & 
 \prod\limits_{i=1}^m C_\eta C_b p_{n + b'}(x_i ) \\
p_n (x_1 \cdot_q \cdots \cdot_q x_m - x_1 \cdot_{q'} \cdots 
  \cdot_{q'} x_m) & \leq &
 \rho(q,q') \prod\limits_{i=1}^m C_\eta C_b p_{n + b'}(x_i)
\end{array}
\]
\end{prop}
\emph{Proof.}
This can be deduced with a piece of careful bookkeeping:
\[
\begin{array}{ll}
p_n (x_1 \cdot_q \cdots \cdot_q x_m ) & \leq \\
  p_n( \sum_{u_i \in W} x_{1 u_1} \cdots x_{m u_m}
  N_{u_1} \cdot_q \cdots \cdot_q N_{u_m}) & \leq \\
\sum_{u_i \in W} | x_{1 u_1} \cdots x_{m u_m} | (\mc N (u_1) + 
  \cdots + \mc N (u_m) + 1)^n \prod_{i=1}^m 
  \norm{(N_{u_i},q)}_o & \leq \\
\sum_{u_i \in W} | x_{1 u_1} \cdots x_{m u_m} | \prod_{i=1}^m 
  C_\eta (\mc N (u_i) + 1)^{n + b + |R_0^+ |} & = \\
\prod_{i=1}^m C_\eta \sum_{u \in W} |x_{i u}| 
  (\mc N(u) +1)^{n + b + |R_0^+ |} & \leq \\
\prod_{i=1}^m C_\eta C_b p_{n + b'}(x_i) \\
 & \\
p_n (N_{u_1} \cdot_q \cdots \cdot_q N_{u_m} - N_{u_1} 
  \cdot_{q'} \cdots \cdot_{q'} N_{u_m}) & \leq \\
\sum_{j=1}^{m-1} p_n(N_{u_1} \cdot_q \cdots \cdot_q N_{u_j} 
  \cdot_q N_{u_{j+1}} \cdot_{q'} \cdots \cdot_{q'} N_{u_m} -& \\ 
\qquad \qquad N_{u_1} \cdot_q \cdots \cdot_q N_{u_j} \cdot_{q'} 
  N_{u_{j+1}} \cdot_{q'} \cdots \cdot_{q'} N_{u_m}) & \leq \\
\sum_{j=1}^{m-1} \rho(q,q') \prod_{i=1}^m C_\eta 
  (\mc N (u_i) + 1)^{n + b + |R_0^+ |} & \leq \\
\rho(q,q') \prod_{i=1}^m C_\eta 
  (\mc N (u_i) + 1)^{n + b + |R_0^+ |}  & \\
 & \\
p_n (x_1 \cdot_q \cdots \cdot_q x_m - x_1 \cdot_{q'} \cdots 
  \cdot_{q'} x_m) & \leq \\
\sum_{u_i \in W} | x_{1 u_1} \cdots x_{m u_m} | p_n(N_{u_1} 
  \cdot_q \cdots \cdot_q N_{u_m} - N_{u_1} \cdot_{q'} 
  \cdots \cdot_{q'} N_{u_m}) & \leq \\
\sum_{u_i \in W} | x_{1 u_1} \cdots x_{m u_m} | \rho(q,q') 
  \prod_{i=1}^m C_\eta (\mc N(u_i) + 1)^{n + b + |R_0^+ |} 
 & = \\
\rho(q,q') \prod_{i=1}^m C_\eta \sum_{u \in W}
  |x_{i u}| (\mc N(u) +1)^ {n + b + |R_0^+ |} & \leq \\
\rho(q,q') \prod_{i=1}^m C_\eta p_{n + b'} (x_i) & 
\end{array}
\]
In these calculations we used \eqref{eq:5.6} and Proposition
\ref{prop:5.6} several times. $\qquad \Box$
\\[3mm]

Knowing how to handle multiple products in $\mc S (\mc R ,q)$, 
we can even make some rough estimates for the holomorphic 
functional calculus.

\begin{cor}\label{cor:5.8}
Let $f : z \to \sum_{m=0}^\infty a_m z^m$ be a holomorphic 
function on a neighborhood of $0 \in \mh C$ and define another 
holomorphic function $\tilde f$ (with the same radius of 
convergence) by $\tilde f (z) = \sum_{m=0}^\infty |a_m| z^m$.
For any $n \in \mh N ,\;x \in \mc S(\mc R,q)$ and $q,q' \in 
B_\rho (q^0,\eta)$ such that $f(x,q)$ and $f(x,q')$ 
are defined we have
\[
\begin{array}{lcr}
p_n (f(x,q)) & \leq & 
  \tilde f \big( C_\eta C_b p_{n + b'} (x) \big) \\
p_n (f(x,q) - f(x,q')) & \leq & \rho(q,q') 
  \tilde f \big( C_\eta C_b p_{n + b'} (x) \big)
\end{array}
\]
\end{cor}
\emph{Proof.}
By Proposition \ref{prop:5.7} we have
\begin{align*}
p_n (f(x,q)) & =
  p_n \left( \sum_{m=0}^\infty a_m (x,q)^m \right) \\
 & \leq \sum_{m=0}^\infty |a_m| p_n ((x,q)^m ) \\
 & \leq \sum_{m=0}^\infty |a_m| \big( C_\eta C_b 
p_{n + b'}(x) \big)^m \\
 & = \tilde f \big( C_\eta C_b p_{n + b'}(x) \big)
\end{align*}
\begin{align*}
p_n (f(x,q) - f(x,q')) & = p_n \left( \sum_{m=0}^\infty a_m 
  ((x,q)^m - (x,q')^m) \right) \\
 & \leq \sum_{m=0}^\infty |a_m| p_n ((x,q)^m - (x,q')^m) \\
 & \leq \sum_{m=0}^\infty |a_m| \rho(q,q') 
  \big( C_\eta C_b p_{n + b'}(x) \big)^m \\
 & = \rho(q,q') \tilde f \big( C_\eta C_b p_{n + b'}(x) \big)
\end{align*}
The right hand sides of these inequalities might be infinite,
but that is no problem. $\qquad \Box$
\\[3mm]

With this result we can show that inverting is continuous as a
function of $x$ and $q$.

\begin{prop}\label{prop:5.9}
The set of invertible elements $\bigcup_{q \in L_{\mc R}} \mc S 
(\mc R ,q)^\times$ is open in $\mc S (\mc R ) \times L_{\mc R}$, 
and inverting is a continuous map from this set to itself.
\end{prop}
\emph{Proof.}
First we recall that if $\norm{(z,q)}_o < 1$, then $z$ is
invertible in $C_r^* (\mc R ,q)$, with inverse
$\sum_{n=0}^\infty (1-z )^n$. Take $q,q' \in B_\rho (q^0, \eta)
,\, y \in \mc S (\mc R) ,\, x \in \mc S (\mc R ,q)^\times$ and
write $a = (x,q)^{-1}$. If the sum converges, then
\begin{equation}\label{eq:5.8}
a \cdot_{q'} \sum_{m=1}^\infty (1 - (x+y) \cdot_{q'} a,q')^m 
= a \cdot_{q'} ((x+y) \cdot_{q'} a,q')^{-1} - a \cdot_{q'} 1
= (x+y,q')^{-1} - a
\end{equation}
By Proposition \ref{prop:5.7} 
\begin{multline}\label{eq:5.9}
p_n( (x+y) \cdot_{q'} a - 1 ) \leq
p_n( x \cdot_{q'} a - x \cdot_q a ) + p_n( y \cdot_{q'} a ) \\
\leq \rho(q,q') C_\eta^2 C_b^2 p_{n+b'}(x) p_{n+b'}(a) + 
C_\eta^2 C_b^2 p_{n+b'}(y) p_{n+b'}(a)
\end{multline}
Let $U$ be the open neighborhood of $(x,q)$ consisting of those\\
$(x+y,q') \in \mc S (\mc R) \times B_\rho (q^0 ,\eta)$ for which 
\begin{align*}
\rho(q,q') C_\eta^3 C_b^2 p_{3b + |R_0 |}(x) p_{3b + |R_0 |}(a) 
& < 1/2 \\
C_\eta^3 C_b^2 p_{3b + |R_0 |}(y) p_{3b + |R_0 |}(a) & < 1/2
\end{align*}
By \eqref{eq:5.9} and Proposition \ref{prop:5.6} we have
\[
\norm{((x+y) \cdot_{q'} a - 1, q')}_o < 1 \qquad
\forall (x+y,q') \in U
\]
so every element of $U$ is invertible. To show the continuity 
of inverting we consider the function
\[
f(z) = \sum_{m=1}^\infty z^m = z / (1-z)
\]
By \eqref{eq:5.8} and Corollary \ref{cor:5.8} we have
\begin{align*}
p_n( (x+y,q')^{-1} - a) & \leq C_b^2 C_\eta^2 p_{n+b'}(a)
p_{n+b'}\big( f(1 - (x+y) \cdot_{q'} a,q') \big) \\
 & \leq C_b^2 C_\eta^2 p_{n+b'}(a) f \big( C_b C_\eta
  p_{n+2b'} (1 - (x+y) \cdot_{q'} a) \big)
\end{align*}
Since $f(0) = 0$ we deduce from \eqref{eq:5.9} that this 
expression is small whenever $\rho (q,q')$ and $y$ are small.
$\qquad \Box$ \\[3mm]

With this result we can see that the \inde{holomorphic 
functional calculus} is continuous in the most general sense.

\begin{cor}\label{cor:5.22}
Let $V \subset \mc S (\mc R ) ,\, Q \subset L_{\mc R}$ and
$U \subset \mh C$ be open subsets such that the spectrum
of every $(x,q) \in V \times Q$ is contained in $U$. Then 
the map \index{$f(x,q)$}
\[
C^{an}(U) \times V \times Q \to \mc S (\mc R ): 
(f,x,q) \to f(x,q)
\]
is continuous.
\end{cor}
\emph{Proof.}
Recall from Theorem \ref{thm:2.7}.4 that
\[
f(x,q) = {\ds \frac{1}{2 \pi i}} \int_\Gamma f (\lambda ) 
(\lambda - x,q)^{-1} d \lambda
\]
for a suitable contour $\Gamma \subset U$ around sp$(x,q)$. 
By Corollary \ref{cor:3.39}.2 sp$(x,q)$ is compact, and by 
Proposition \ref{prop:5.9} it depends continuously on $x$ and
$q$. Therefore we can find a contour which is suitable for
every $(x',q')$ in a small neighborhood of $(x,q)$. Now apply
Proposition \ref{prop:5.9} and Theorem \ref{thm:2.7}.3.
$\qquad \Box$ \\[3mm]

Let \inde{$\mc H (\mc R)^*$} be the algebraic dual of 
$\mc H (\mc R)$, which we identify, using the bitrace $\tau$, 
with the space of all formal infinite sums $\sum_{w \in W} 
x_w N_w$. The length function $\mc N$ may also be considered 
as an endomorphism of $\mc H (\mc R )^*$: 
\begin{equation}
\lambda (\mc N) : \sum_{w \in W} x_w N_w \to 
\sum_{w \in W} \mc N (w) x_w N_w
\end{equation}
This is an unbounded operator on $\mf H (\mc R)$, but it does
restrict to a continuous endomorphism of $\mc S (\mc R)$. For
$T \in B (\mf H (\mc R))$ put \inde{$D(T)$}$= [\lambda (\mc N),
T]$. Inspired by the work of Vign\'eras \cite[Section 7]{Vig} 
we study the space \index{lambdan@$\lambda (\mc N)$}
\index{$V_{\mc N}^\infty (\mc R ,q)$}
\begin{equation}
V_{\mc N}^\infty (\mc R ,q) = \{ x \in \mc H (\mc R)^* : 
D^n (\lambda (x)) \in B (\mf H (\mc R)) \; \forall n \in 
\mh Z_{\geq 0} \}
\end{equation}
We use the topology defined by the collection of seminorms
\begin{equation}\label{eq:5.10}
\big\{ \norm{D^n (\lambda ( \:\cdot\: ) )}_{B (\mf H (\mc R))}
,\; n \in \mh Z_{\geq 0} \big\}
\end{equation}
In fact we already know this space:

\begin{lem}\label{lem:5.10}
\[
V_{\mc N}^\infty (\mc R ,q) = \mc S (\mc R ,q)
\]
\end{lem}
\emph{Proof.}
From the proof of Proposition \ref{prop:5.6} we see that for 
any $y = \sum_{v \in W} y_v N_v \in \mf H (\mc R) ,\,n \in 
\mh Z_{\geq 0} ,\, u \in W$
\begin{align*}
\norm{D^n (\lambda (N_u )) y}_\tau & = \norm{ \sum_v y_v 
\sum_{i=0}^n (-1)^i \left( \ds{n \atop i} \right) \lambda 
(\mc N )^{n-i} \lambda (N_u ) \lambda (\mc N )^i N_v }_\tau \\
& = \norm{ \sum_v y_v \sum_{i=0}^n (-1)^i \left( \ds{n \atop i} 
\right) \sum_I \eta_I D_v^u (I) \mc N (u_I v)^{n-i} \mc N (v)^i
N_{u_I v} }_\tau \\
& = \norm{ \sum_v y_v \sum_I \eta_I D_v^u (I) (\mc N (u_I v) - 
\mc N (v))^n N_{u_I v} }_\tau \\
& \leq \mc N (u)^n \norm{ \sum_v |y_v| \sum_I |\eta_I | 
D_v^u (I) N_{u_I v} }_\tau \\
& \leq \mc N (u)^n 3 (\mc N (u) + 1)^{|R_0^+ |} \max \big\{1, 
\eta^{|R_0^+ |} \big\} \norm{ \sum_v |y_v| N_{u_I v} }_\tau \\
& = \mc N (u)^n 3 (\mc N (u) + 1)^{|R_0^+ |} \max \big\{1, 
\eta^{|R_0^+ |} \big\} \norm{y}_\tau
\end{align*}
where $\eta = \rho (q,q^0 )$. Hence, for 
$x = \sum_u x_u N_u \in \mc H (\mc R)^*$
\begin{align*}
\norm{D^n (\lambda (x)) }_o & = \norm{\sum_u x_u D^n (\lambda 
(N_u)) }_o \\
& \leq \sum_u |x_u | 3 (\mc N (u) + 1)^{n+ |R_o^+ |} \max 
\big\{1, \eta^{|R_0^+ |} \big\} \\
& \leq C_\eta p_{n + b + |R_0^+ |} (x)
\end{align*}
On the other hand, since $\Omega' = \{ \omega \in W : 
\mc N (\omega ) = 0 \}$ is finite,
\begin{align*}
p_n (x)^2 & \leq \sum_{u \in W} (\mc N (u) + 1)^{2n} |x_u |^2 \\
& \leq \sum_{\omega \in \Omega'} |x_\omega |^2 + 
  4^n \sum_{u \in W} \mc N (u)^{2n}  |x_u |^2 \\
& \leq | \Omega' | \norm{x}^2_\tau + 
  4^n \norm{\lambda (N)^n x}_\tau^2 \\
& = |\Omega' | \norm{ \lambda (x) N_e }_\tau^2 + 
  4^n \norm{ D^n (\lambda (x)) N_e }_\tau^2 \\
& \leq |\Omega' | \norm{ \lambda (x) }_o^2 + 
  4^n \norm{ D^n (\lambda (x)) }_o^2 \\
& \leq \big( |\Omega' |^{1/2} \norm{ \lambda (x) }_o + 
  2^n \norm{ D^n (\lambda (x)) }_o \big)^2
\end{align*}
Therefore the collections of seminorms $\{p_n : n \in 
\mh Z_{\geq 0} \}$ and \eqref{eq:5.10} are equivalent.
$\qquad \Box$ \\[3mm]

So we found a different way to construct the Schwartz algebra of 
an affine Hecke algebra. An advantage of this construction is 
that it allows us to prove Corollary \ref{cor:3.39} in a more
elementary way, relying only on the density of $V_{\mc N}^\infty 
(\mc R ,q)$ in $C_r^* (\mc R ,q)$ and not on any representation 
theory.

\begin{thm}\label{thm:5.11}
\begin{enumerate}
\item $V_{\mc N}^\infty (\mc R ,q)$ is a complete locally 
convex algebra with jointly continuous multiplication.
\item $V_{\mc N}^\infty (\mc R ,q)^\times$ is open in 
$V_{\mc N}^\infty (\mc R ,q)$, and inverting is a continuous 
map from this set to itself.
\item An element of $V_{\mc N}^\infty (\mc R ,q)$ is invertible
if and only if it is invertible in $C_r^* (\mc R ,q)$.
\end{enumerate}
\end{thm}
\emph{Proof.}
1. By Lemma \ref{lem:5.10} $V_{\mc N}^\infty (\mc R ,q)$ is a
Fr\'echet space. Since $D$ is a derivation, it is also a 
topological algebra with jointly continuous multiplication.\\
2. See \cite[Lemma 16]{Vig}. Suppose that $x \in 
V_{\mc N}^\infty (\mc R ,q)$ and $\norm{x}_o < 1$. Then
$1 - x \in C_r^* (\mc R ,q)^\times$ and $(1-x)^{-1} = 
\sum_{n=0}^\infty x^n \in C_r^* (\mc R ,q)$. We have to show 
that this sum converges in $V_{\mc N}^\infty (\mc R ,q)$. 
For $n,r \in \mh N$
\[
D^r (\lambda (x)) = \sum_{r_1 + \cdots + r_n = r} \ds{ \frac{
r! D^{r_1} (\lambda (x)) \cdots D^{r_n} (\lambda (x)) }{r_1 !
\cdots r_n !}}
\]
Every product $D^{r_1} (\lambda (x)) \cdots D^{r_n} 
(\lambda (x))$ contains at least $n-r$ factors $\lambda (x)$
and at most $r$ factors of the form $D^i (\lambda (x))$ with
$i > 0$. Therefore
\begin{align*}
& \norm{D^r (\lambda (x)) }_{B (\mf H (\mc R))} \leq n^r M^r
\norm{x}_o^n \\
& M = \max \left\{ \norm{D^i (\lambda (x)) }_{B (\mf H (\mc R))}
\norm{x}_o^{-1} : 0 < i \leq r \right\}
\end{align*}
This gives
\[
\norm{D^r (\lambda (1-x)^{-1}) }_{B (\mf H (\mc R))} \leq
\sum_{n=0}^\infty n^r M^r \norm{x}_o^n \leq 
r! M^r (1 - \norm{x}_o )^{-r-1}
\]
from which we conclude that indeed $(1-x)^{-1} \in 
V_{\mc N}^\infty (\mc R ,q)$ and that inverting in 
$V_{\mc N}^\infty (\mc R ,q)$ is continuous around 1. In 
general, if $y \in V_{\mc N}^\infty (\mc R ,q)^\times$ then
we can use the "translation" $\lambda (y^{-1})$ to show that 
$V_{\mc N}^\infty (\mc R ,q)^\times$ contains an open 
neighborhood of $y$ and that inverting is continuous on this 
set.\\
3. Suppose that $z \in V_{\mc N}^\infty (\mc R ,q) \cap
C_r^* (\mc R ,q)^\times$. By Lemma \ref{lem:5.10}
$V_{\mc N}^\infty (\mc R ,q)$ is dense in $C_r^* (\mc R ,q)$, 
so we can find $y \in z^{-1} B \cap B z^{-1}$, where 
\[
B = \{ x \in C_r^* (\mc R ,q) : \norm{1-x}_o < 1 \}
\]
By the above 
\[
yz, zy \in B \cap V_{\mc N}^\infty (\mc R ,q) \subset 
V_{\mc N}^\infty (\mc R ,q)^\times
\]
so $z$ is also invertible in $V_{\mc N}^\infty (\mc R ,q). 
\qquad \Box$ \\[3mm]

Note that it does not follow from these considerations that
$V_{\mc N}^\infty (\mc R ,q)$ is a m-algebra. To prove that we
still have to use Theorem \ref{thm:3.19}.

Consider the bundle of Banach spaces $\bigsqcup_{q \in 
L_{\mc R}} C_r^* (\mc R ,q)$ over $L_{\mc R}$. For any fixed
$x \in \mc S (\mc R)$ the constant function $q \to x$ is 
a section of this bundle, and by Proposition \ref{prop:5.6}
the function $q \to \norm{(x,q)}_o$ is continuous on 
$L_{\mc R}$. So by \cite[Proposition 10.2.3]{Dix} there is a 
unique collection $\Gamma$ of sections of $\bigsqcup_{q 
\in L_{\mc R}} C_r^* (\mc R ,q)$ containing all these constant
sections, which makes this into a field of $C^*$-algebras, in 
the sense of Dixmier \cite[Section 10.3]{Dix}. By construction 
$\Gamma$ contains all continuous maps $L_{\mc R} \to \mc S (\mc R)$. 
In particular for any compact set $Q \subset \mc L_{\mc R}$ we 
can construct the unital $C^*$-algebra \inde{$C_r^* (\mc R ,Q)$}$
= \Gamma \big|_Q$, which contains $C (Q ; \mc S (\mc R))$ as a 
dense subalgebra. \index{Gamma@$\Gamma$}

For a related construction, assume that $Q$ is a smooth 
submanifold of $L_{\mc R}$, not necessarily compact. Although
the author is not aware of any precise definition, he believes
it makes sense to call $\bigsqcup_{q \in L_{\mc R}} \mc S 
(\mc R ,q)$ a field of Fr\'echet algebras. By Proposition 
\ref{prop:5.7} the set of smooth sections 
\index{$\mc S (\mc R ,Q)$}
\begin{equation}
\mc S (\mc R ,Q) = C^\infty (Q ; \mc S (\mc R))
\end{equation}
is a Fr\'echet space and a topological algebra with jointly
continuous multiplication. However, the author does not know
whether $\mc S (\mc R ,Q)$ is a m-algebra in general. For every
$q \in Q$ we can obtain submultiplicative seminorms on 
$\mc S (\mc R ,q)$ from the Fourier transform, but these are not
easily expressible in terms of the elements $N_w$. Therefore it
is not clear whether we can choose them in some sense continuous
in $q$ and "glue" such seminorms to a family of seminorms that 
defines the topology of $\mc S (\mc R ,Q)$.

Seminorms that we can construct are related to the principal 
series representations, which exist for every $q \in L_{\mc R}$. 
From Theorem \ref{thm:3.4}.1 we get a linear bijection 
\[
\mc H (\mc R ,q) \to \mc H (W_0 ,q) \otimes \mh C [X]
\]
This can be extended to a continuous map 
\index{fitq@$\phi_{\theta ,q}$}
\begin{equation}
\phi_{\theta ,q} : \mc S (\mc R ,q) \to 
\mc H (W_0 ,q) \otimes \mc S (X)
\end{equation}
which is essentially the direct integral of all unitary principal
series representations. In general $\phi_{\theta ,q}$ is
neither injective nor surjective. For $n \in \mh N$ we define the
following norm on $\mc H (W_0 ,q) \otimes \mc S (X)$ :
\index{sigman@$\sigma_n$} \index{sigmanq@$\sigma_{n,q}$}
\begin{equation}
\sigma_n \Big( \sum_{w \in W_0 , x \in X} y_{w,x} N_w \theta_x 
\Big) = \sum_{w \in W_0 , x \in X} (\mc N (x) + 1)^n
\end{equation}
We compose it with $\phi_{\theta ,q}$ to get the seminorm
$\sigma_{n,q} = \sigma_n \circ \phi_{\theta ,q}$ on 
$\mc S (\mc R ,q)$. These seminorms are continuous in $q$:

\begin{lem}\label{lem:5.12}
Let $n \in \mh N ,\, \eta > 0$ and $q,q'\in B_\rho (q^0, \eta)$.
There exists a real number $C_{n,\eta}$ such that for all
$z \in \mc S (\mc R)$
\[
\begin{array}{lcr}
\sigma_{n,q}(z) & \leq & C_{n,\eta} p_{n + b + |R_0^+ |}(z) \\
| \sigma_{n,q}(z) - \sigma_{n,q'}(z) | & \leq &
\rho (q,q') C_{n,\eta} p_{n + b + |R_0^+ |}(z)
\end{array}
\]
\end{lem}
\emph{Proof.}
Let $x^+ \in X^+$ and put 
\[
P = \{ \alpha \in F_0 : \alpha^\vee (x^+ ) = 0 \}
\]
By \cite[Proposition 1.15]{Hum} we have
\[
W_P = \{ w \in W_0 : w x^+ = x^+ \}
\]
and hence
\[
W_0 x^+ W_0 = W^P x^+ W_0 = W_0 x^+ \big( W^P \big)^{-1}
\]
From \eqref{eq:3.10} we see that the $N_u$ with $u \in W_P$
commute with $\theta_{x^+} \in \mc H (\mc R ,q)$. Pick any
$w \in W_0 x^+ W_0$. From Theorem \ref{thm:5.5} we see that
there is a unique way to write
\begin{equation}\label{eq:5.11}
N_w = \sum_{u \in W_0 , v \in W^P} c_{u,v}^w N_u 
\theta_{x^+} N_{v^{-1}}
\end{equation}
where every $c_{u,v}^w$ is a polynomial in the $\eta_i$ of
degree at most $|R_0 |$. From the length formula in Proposition
\ref{prop:3.1} we see that the $c_{u,v}^w$ depend only on $P$,
in the following sense. If $w = w_1 x^+ w_2 \; (w_1 ,w_2 \in 
W_0 )$ and $w' = w_1 x' w_2$ where $x' \in X^+$ and $\{ \alpha 
\in F_0 : \alpha^\vee (x') = 0 \} = P$, then $c_{u,v}^{w'} = 
c_{u,v}^w$. In particular there are only finitely many different 
$c_{u,v}^w$, less than $|W_0 |^2 2^{|F_0 |}$.

Assume for simplicity that $\eta \geq 1$, so 
\[
\norm{(N_v ,q)}_o \leq (2 \eta)^{|R_0^+ |} 
\quad \forall v \in W_0 
\]
Let $K$ be an upper bound for the absolute values of all 
$c_{u,v}^w$, also under the condition $q \in B_\rho (q^0 ,\eta )$. 
From a repeated application of \eqref{eq:3.10} we see that 
$\theta_{x^+} N_{v^{-1}}$ equals a sum of at most $(2 \mc N(x^+ ) 
+2 )^{|R_0^+ |}$ terms of the form $\eta_I N_{v'} \theta_{x}$, 
where $\mc N (x) \leq \mc N (x^+ ) ,\, v'\in W_0$ and $I$ is a 
multi-index with $|I| \leq |R_0^+ |$. This leads to the estimate
\begin{align*}
\sigma_{n,q}(N_w ) & = \sigma_n \Big( \sum_{u,v} c_{u,v}^w N_u 
\sum_{I,v',x} \eta_I N_{v'} \theta_x \Big) \\
& \leq \sum_v |W_0 |^{1/2} \norm{\sum_u c_{u,v}^w N_u }_\tau
\sum_{I,v',x} \norm{(N_{v'},q)}_o |\eta_I | (\mc N (x) + 1)^n \\
& \leq |W_0 | K (2 \eta )^{|R_0^+ |} \eta^{|R_0^+ |} 
(2 \mc N (x^+ ) + 2 )^{|R_0^+ |} (\mc N (x^+ ) + 1)^n \\
& = |W_0 | K (2 \eta )^{|R_0 |} 
  (\mc N (x^+ ) + 1)^{n + |R_0^+ |} \\
& \leq |W_0 | K (2 \eta )^{|R_0 |} (|R_0^+ | + 1)^{n + |R_0^+ |}
(\mc N (w) + 1)^{n + |R_0^+ |}
\end{align*}
For $z = \sum_{w \in W} z_w N_w \in \mc S (\mc R)$ we obtain
\begin{align*}
\sigma_{n,q}(z) & \leq \sum_{w \in W} |z_w | \sigma_{n,q}(N_w) \\
& \leq |W_0 | K (2 \eta )^{|R_0 |} (|R_0^+ | + 1)^{n + |R_0^+ |}
\sum_{w \in W} |z_w | (\mc N (w) + 1)^{n + |R_0^+ |} \\
& \leq |W_0 | K (2 \eta )^{|R_0 |} (|R_0^+ | + 1)^{n + |R_0^+ |}
C_b p_{n + b + |R_0^+ |} (z)
\end{align*}
Plugging the description of $\theta_{x^+} N_{v^{-1}}$ into 
\eqref{eq:5.11} we see that
\begin{equation}
N_w = \sum_{u \in W_0 , x \in X} y_{u,x}^w N_u \theta_x
\end{equation}
where the $y_{u,x}^w$ are polynomials in the $\eta_i$ of 
degree at most $2 |R_0 |$. Therefore we can write
\begin{equation}
\phi_{\theta ,q}(z) = \sum_{I : |I| \leq 2 |R_0 |} 
\sum_{u \in W_0 , x \in X} \eta_I z_{I,u,x} N_u \theta_x :=
\sum_I \eta_I z_I 
\end{equation}
Since the finite collection $\{ \eta_I : |I| \leq 2 |R_0 | \}$
is linearly independent, considered as functions of the 
$\eta_i$, we can find constants $K_n$ such that 
\[
\sigma_n (z_I) \leq K_n p_{n + b + |R_0^+ |}(z) \quad \forall I
\]
Hence
\begin{align*}
| \sigma_{n,q}(z) - \sigma_{n,q'}(z) | & \leq \sum_{I : |I| 
\leq 2 |R_0 |} |\eta_I - \eta'_I| \sigma_n (z_I) \\
& \leq \rho (q,q') 2 |R_0 | \eta^{2 |R_0 |} K_n 
p_{n + b + |R_0^+ |}(z)
\end{align*}
To finish the proof we take for $C_{n,\eta}$ the maximum of\\
$|W_0 | K (2 \eta )^{|R_0 |} (|R_0^+ | + 1)^{n + |R_0^+ |} C_b$ 
and $2 |R_0 | \eta^{2 |R_0 |} K_n . \qquad \Box$
\\[3mm] 

We conclude this section with an important remark. All the 
estimates obtained here can be generalized to $M_k (\mc S 
(\mc R ,q))$ for any $k \in \mh N$. Of course we first have to
(re)define \index{$p_n^{(k)}$}
\begin{equation}\label{eq:5.35}
p_n^{(k)} \big (z_{i,j} )_{i,j=1}^k \big) = 
\max_{1 \leq i,j \leq k} p_n (z_{i,j} )
\end{equation}
but then all our results can be extended by standard techniques.
In the coming sections we assume that this has been done, and we
attach a superscript $(k)$ to the modified constants.
\\[4mm]

\section{Scaling the labels}
\label{sec:5.3}

Fix a root datum $\mc R$ and a label function $q$. Instead of 
considering general deformations of $q$, we concentrate on the
scaled label functions $q^\ep$ for $\ep \in [-1,1]$.
Opdam \cite[Section 5]{Opd3} was the first to realize that a 
great deal of the representation theory of $\mc H (\mc R ,q^\ep )$ 
can be "scaled" accordingly. We will 
construct isomorphisms $\mc S (\mc R ,q^\ep ) \to \mc S (R ,q)$
for $\ep > 0$ and an injection $\mc S (W) = S (\mc R ,q^0 ) 
\to \mc S (\mc R ,q)$, all depending continuously on $\ep$.
This requires a lot of long calculations, which rely on the 
technical parts of Chapter 3 and Section \ref{sec:5.2}.

Recall from \cite[Section 6.3]{Opd3} that we always have a 
canonical isomorphism
\begin{equation}\label{eq:5.40}
\mc H (\mc R ,q) \isom \mc H (\mc R ,q^{-1}) :
N_w \to (-1)^{\ell (w)} N_w
\end{equation}
This map preserves * and $\tau$, and it extends to the Schwartz 
and $C^*$-completions. However, our scaling maps will not lead 
to this map for $\ep = -1$.

To facilitate the study of $\mc S (\mc R ,q)$ we want to regard
it not only as an algebra of invariant sections (via the Fourier
transform), but also as the image of a projector in a larger
algebra. Let $\mc P' (F_0 )$ be a complete set of representatives 
for the action of $\mc W$ on the power set of $F_0$. By Theorem 
\ref{thm:3.19} there are canonical direct sum decompositions
\begin{equation}\label{eq:5.46}
\begin{array}{lllll}
\mc S (\mc R ,q) & = & 
  \bigoplus\limits_{P \in \mc P' (F_0 )} \mc S (\mc R ,q)_P & = &
  \bigoplus\limits_{P \in \mc P' (F_0 )} \mc S (\mc R ,q) \: e_P \\
C_r^* (\mc R ,q) & = & 
  \bigoplus\limits_{P \in \mc P' (F_0 )} C_r^* (\mc R ,q)_P & = &
  \bigoplus\limits_{P \in \mc P' (F_0 )} C_r^* (\mc R ,q) \: e_P 
\end{array}
\end{equation}
where the \inde{$e_P$} are central idempotents in $\mc S (\mc R 
,q)$. But this can be refined. Let $\Delta'$ be a collection 
of representatives for the action of $\mc W$ on $\Delta$. For 
$(P, \delta ) \in \Delta$, we let $\mc W_\delta$ be the isotropy 
group of $(P, \delta )$ in $\mc W$. The Fourier transform gives 
an isomorphism \index{$A_\delta$} \index{$\mc F'$}
\begin{equation}\label{eq:5.12}
\mc F' : \mc S (\mc R ,q) \isom \bigoplus_{(P,\delta ) \in 
\Delta'} C^\infty \left( T_u^P ; \mr{End} \big(\mc H (W^P) 
\otimes V_\delta \big) \right)^{\mc W_\delta} := 
\bigoplus_{(P,\delta ) \in \Delta'} A_\delta^{\mc W_\delta}
\end{equation}
We must be careful when taking invariants, since
\begin{equation}\label{eq:5.13}
\mc W_\delta \to A_\delta^\times : 
g \to \pi (g,P,\delta, \cdot )
\end{equation}
is not necessarily a group homomorphism. In fact, by 
\eqref{eq:3.36} it is a projective representation. By Schur's
theorem \cite{Schu} there exists a finite central extension
\begin{equation}\label{eq:5.14}
\{e\} \to N_\delta \to \Gamma_\delta \to \mc W_\delta \to \{e\} 
\end{equation}
such that every projective representation of $\mc W_\delta$ lifts 
to a unique linear representation of $\Gamma_\delta$.
This lift does not depend on the $\tilde \delta_g$ that we 
chose in \eqref{eq:3.28} to construct $\pi (g,\xi)$ for
$\xi = (P,\delta ,t)$. In fact, the problems with \eqref{eq:5.13} 
arise only from the ambiguity in the definition of $\tilde \delta_g$.
Lifting things to a linear representation of $\Gamma_\delta$ is 
therefore equivalent to picking, for every lift $\gamma \in 
\Gamma_\delta$ of $g \in \mc W_\delta$, a multiple $\tilde 
\delta_\gamma$ of $\tilde \delta_g$ such that
 \index{deltatg@$\tilde \delta_\gamma$} \index{$u_\gamma$}
\[
\Gamma_\delta \to A_\delta^\times :
\gamma \to \pi (\gamma,P,\delta, \cdot)
\]
becomes multiplicative. Writing \index{Gammad@$\Gamma_\delta$}
\[
u_\gamma (\xi) = \pi (\gamma , \gamma^{-1} \xi)
\]
\eqref{eq:3.36} becomes the cocycle relation
\begin{equation}\label{eq:5.15}
u_{\gamma \gamma'} = u_\gamma u_{\gamma'}^\gamma
\end{equation}
Notice the similarity with the proof of Lemma \ref{lem:2.17}.
Consider the crossed product 
\begin{equation}
A_\delta \rtimes \Gamma_\delta \cong 
\mr{End}\big(\mc H (W^P) \otimes V_\delta \big) \otimes 
C^\infty \big( T_u^P \big) \rtimes \Gamma_\delta
\end{equation}
with respect to the action of $\Gamma_\delta$ on $T_u^P$.
Lemma \ref{lem:A.1} gives some information about this 
algebra. We still need to determine $A_\delta^{\mc W_\delta} 
= A_\delta^{\Gamma_\delta}$, but using the multiplication in
$A_\delta \rtimes \Gamma_\delta$ we can write the group 
action as
\begin{equation}\label{eq:5.21}
\gamma (a) = u_\gamma \gamma a \gamma^{-1} u_\gamma^{-1}
\end{equation}
Moreover by \eqref{eq:5.15} 
\[
\Gamma_\delta \to \big( A_\delta \rtimes \Gamma_\delta 
\big)^\times : \gamma \to u_\gamma \gamma
\]
is a unitary representation, so \index{$p_\delta (u)$}
\begin{equation}
p_\delta (u) := | \Gamma_\delta |^{-1} 
\sum_{\gamma \in \Gamma_\delta} u_\gamma \gamma 
\in A_\delta \rtimes \Gamma_\delta 
\end{equation}
is a projection. By Lemma \ref{lem:A.2} the map 
\begin{equation}\label{eq:5.24}
A_\delta^{\Gamma_\delta} \to p_\delta (u) \big( 
A_\delta \rtimes \Gamma_\delta \big) p_\delta (u) :
a \to p_\delta (u) a p_\delta (u)
\end{equation}
is an isomorphism of pre-$C^*$-algebras. 

Consider for $\ep \in [-1,1]$ the affine Hecke algebras
\inde{$\mc H_\ep$} $= \mc H (\mc R ,q^\ep )$ with label 
functions \inde{$q^\ep$}$(w) = q(w)^\ep$. Let 
\inde{$c_{\alpha,\ep}$} and $\imath^o_{w,\ep}$ be the 
$c$-functions and normalized intertwiners for these algebras. 
From the formula \eqref{eq:3.11} we see that the residual 
cosets are also related by scaling in suitable directions. 
Their tempered forms all approach $T_u$ when $\ep \to 0$. 
\index{iowe@$\imath^o_{w,\ep}$}

Take $r \in T$ and write $r = r_u \exp (r_s )$ with $r_u \in 
T_u$ and $r_s \in \mf t_{rs}$. Let $B \subset \mr{Lie}(T)$ 
be a ball satisfying Conditions \ref{cond:3.7}. Then $\ep B$ 
satisfies these conditions with respect to $r_u \exp (\ep r_s 
)$ and $q^\ep$, except that it is not open for $\ep = 0$. 
(Here we use $|\ep | \leq 1$.) We put \index{$U_\ep$}
\[
U_\ep = W_0 \big( r_u \exp (\ep (r_s + B)) \big)
\]
and we define a $W_0$-equivariant scaling map 
\begin{equation}\label{eq:5.16}
\begin{aligned}
& \sigma_\ep : U \to U_\ep\\
& \sigma_\ep \big( w (r \exp (b)) \big) = 
w \big( r_u \exp (\ep (r_s + b)) \big)
\end{aligned}
\end{equation}
Assume now that $0 \neq \ep \in [-1,1]$. As was noted in 
\cite[Lemma 5.1]{Opd3}, $\sigma_\ep$ is an analytic 
diffeomorphism. We can combine it with \eqref{eq:3.14} and 
Theorem \ref{thm:3.16} to construct algebra isomorphisms 
\index{sigmae@$\sigma_\ep$}
\begin{equation}\label{eq:5.17}
\begin{aligned} 
& \rho_\ep : \mc H_\ep^{me}(U_\ep ) \isom \mc H^{me}(U) \\
& \rho_\ep \left( \sum_{w \in W_0} a_w \imath^o_{w,\ep} \right)
= \sum_{w \in W_0 } (a_w \circ \sigma_\ep ) \imath^o_w
\qquad a_w \in C^{me}(U_\ep )
\end{aligned}
\end{equation}
We intend to show that these maps depend analytically on $\ep$
and have a well-defined limit as $\ep \to 0$. Notice
that $\sigma_0 = \lim_{\ep \to 0} \sigma_\ep$ is a locally
constant map with range $W_0 r_u$. \index{rhoe@$\rho_\ep$}

\begin{lem}\label{lem:5.13}
For $\ep \neq 0$ and $\alpha \in R_1$ write \inde{$d_{\alpha 
,\ep}$} $= (c_{\alpha, \ep} \circ \sigma_\ep) c_\alpha^{-1}$. 
This defines a bounded invertible analytic function of $u$ 
and $\ep$ which extends to a function on $\overline U \times 
[-1,1]$ with the same properties.
\end{lem}
\emph{Proof.}
This is an extended version of \cite[Lemma 5.2]{Opd3}.
Let us write
\begin{multline*}
d_{\alpha, \ep}(u) \quad = \quad {\ds \frac{f_1 f_2 f_3 f_4
}{g_1 g_2 g_3 g_4}(u) \quad = \quad
\frac{1 + \theta_{-\alpha /2}(u)}{1 + \theta_{-\alpha /2}
(\sigma_\ep (u))} \: \times } \\
{\ds \frac{1 + q_{\alpha^\vee}^{-\ep /2} \theta_{-\alpha /2}
(\sigma_\ep (u))}{1 + q_{\alpha^\vee}^{-1/2} \theta_{-\alpha 
/2}(u)} \frac{1 - \theta_{-\alpha /2}(u)}{1 - \theta_{-\alpha 
/2}(\sigma_\ep (u))} 
\frac{1 - q_{\alpha^\vee}^{-\ep /2} q_{2 \alpha^\vee}^{-\ep} 
\theta_{-\alpha /2}(\sigma_\ep (u))}{1 - q_{\alpha^\vee}^{-1/2} 
q_{2 \alpha^\vee}^{-1} \theta_{-\alpha /2}(u)} } 
\end{multline*}
We see that $d_{\alpha ,\ep}(u)$ extends to an invertible
analytic function on $\overline U \times [-1,1]$ if none of the 
quotients $f_n / g_n$ has a zero or a pole on this domain. By 
Condition \ref{cond:3.7}.2 there is a unique $b \in w (r_s + 
\overline B ) / 2$ such that $u = w(r_u ) \exp (2b)$. This forms 
a coordinate system on $w ( r \exp (\overline B))$, and 
$\sigma_\ep (u) = w(r_u ) \exp (2b \ep)$. By Condition 
\ref{cond:3.7}.4 if either $f_n (u) = 0$ or $g_n (u) = 0$ for 
some $u \in w ( r \exp (\overline B)) \subset \overline U$, then 
$f_n (w(r)) = g_n (w(r)) = 0$. One can easily check that in
this situation
\[
{\ds \frac{f_n (u)}{g_n (u)} = \left( \frac{1 - e^{-\alpha (b) 
\ep}}{1 - e^{-\alpha (b)}} \right)^{(-1)^n} }
\]
Again by Condition \ref{cond:3.7}.2 the only critical points 
of this function are those for which $\alpha (b) = 0$. If 
$\ep \neq 0$ then both the numerator and the denominator have 
a zero of order 1 at such points, so the singularity is 
removable. For the case $\ep = 0$ we need to have a closer 
look. In our new coordinate system we can write
\begin{multline*}
c_{\alpha ,\ep}(\sigma_\ep (u)) \quad = \quad {\ds \frac{f_2 
(u) f_4 (u)}{g_1 (u) g_3 (u)}} \\
= \quad {\ds \frac{r_u (w^{-1} \alpha /2) + \big( 
q_{\alpha^\vee}^{-1/2} e^{-\alpha (b)} \big)^\ep}{r_u (w^{-1} 
\alpha /2) + e^{-\alpha (b) \ep}} \frac{r_u (w^{-1} \alpha /2) - 
\big( q_{\alpha^\vee}^{-1/2} q_{2 \alpha^\vee}^{-1} e^{-\alpha 
(b)} \big)^\ep}{r_u (w^{-1} \alpha /2) - e^{-\alpha (b) \ep}} }
\end{multline*}
Standard calculations using L'Hospital's rule show that
\[
\lim_{\ep \to 0} c_{\alpha ,\ep}(\sigma_\ep (u)) = 
\left\{ \begin{array}{c@{\qquad\mr{if}\quad}ccc}
1 & r_u (w^{-1} \alpha ) & \neq & 1 \\
{\ds \frac{\alpha (b) + \log (q_{\alpha^\vee}) /2}{\alpha (b)} } 
& r_u (w^{-1} \alpha / 2) & = & -1 \\
{\ds \frac{\alpha (b) + \log (q_{2\alpha^\vee}) + 
\log (q_{\alpha^\vee}) / 2}{\alpha (b)} } & 
r_u (w^{-1} \alpha / 2) & = & 1 
\end{array} \right.
\]
Thus at least $d_{\alpha ,0} = \lim_{\ep \to 0} d_{\alpha ,\ep}$
exists as a meromorphic function on $\overline U$. 
For $r_u (w^{-1} \alpha ) \neq 1 ,\,
d_{\alpha ,0} = c_\alpha^{-1}$ is invertible by Condition 
\ref{cond:3.7}.4. For $r_u (w^{-1} \alpha / 2) = -1$ we have
\[
{\ds d_{\alpha ,0}(u) = \frac{1 - e^{-\alpha (b)}}{\alpha (b)}
\frac{\alpha (b) + \log (q_{\alpha^\vee}) / 2}{1 - 
q_{\alpha^\vee}^{-1/2} e^{-\alpha (b)}} \frac{1 + e^{-\alpha 
(b)}}{1 + q_{\alpha^\vee}^{-1/2} q_{2\alpha^\vee}^{-1} 
e^{-\alpha (b)}} }
\]
while for $r_u (w^{-1} \alpha / 2) = 1$
\[
{\ds d_{\alpha ,0}(u) = \frac{1 - e^{-\alpha (b)}}{\alpha (b)}
\frac{1 + e^{-\alpha (b)}}{1 + q_{\alpha^\vee}^{-1/2} e^{-\alpha 
(b)}} \frac{\alpha (b) + \log (q_{2\alpha^\vee}) + \log (q_{
\alpha^\vee}) / 2}{1 - q_{\alpha^\vee}^{-1/2} q_{2\alpha^\vee}^{
-1} e^{-\alpha (b)}} }
\]
These expressions define invertible functions by Condition 
\ref{cond:3.7}.2. We conclude that indeed $d_{\alpha ,\ep}(u)$ 
and $d_{\alpha ,\ep}^{-1}(u)$ are analytic functions on 
$\overline U \times [-1,1]$. Since this domain is compact, they
are bounded. $\Box$
\\[2mm]

Note that $d_{\alpha ,1} = 1$ and that 
\begin{multline}\label{eq:5.20}
d_{\alpha ,-1}(u) = {\ds \frac{r_u (w^{-1} \alpha) - e^{-2 
\alpha (b)}}{r_u (w^{-1} \alpha) - e^{2 \alpha (b)}} 
\frac{r_u (w^{-1} \alpha /2) + q_{\alpha^\vee}^{1/2} e^{\alpha 
(b)}}{r_u (w^{-1} \alpha /2) + q_{\alpha^\vee}^{-1/2} e^{-\alpha 
(b)}} \: \times } \\
{\ds \frac{r_u (w^{-1} \alpha /2) - q_{\alpha^\vee}^{1/2} 
q_{2 \alpha^\vee} e^{\alpha (b)}}{r_u (w^{-1} \alpha /2) - 
q_{\alpha^\vee}^{-1/2} q_{2 \alpha^\vee}^{-1} e^{-\alpha (b)}} }
\end{multline}
If either $r_u (w^{-1} \alpha) = \theta_{-\alpha}(w(r_u )) = 1$ 
or $q_{\alpha^\vee} = q_{2 \alpha^\vee} = 1$ this simplifies 
considerably to
\begin{equation}\label{eq:5.18}
d_{\alpha ,-1}(u) = q_{\alpha^\vee} q_{2 \alpha^\vee}
\end{equation}
We can use Lemma \ref{lem:5.13} to show that the maps $\rho_\ep$ 
preserve analyticity: \index{rho0@$\rho_0$}

\begin{lem}\label{lem:5.14}
The isomorphisms \eqref{eq:5.17} restrict to isomorphisms
\[
\rho_\ep : \mc H_\ep^{an}(U_\ep ) \isom \mc H^{an}(U)
\] 
These maps have a well-defined limit homomorphism
\[
\rho_0 = \lim_{\ep \to 0} \rho_\ep : \mh C [W] \to \mc H^{an}(U)
\] 
such that for every $w \in W$ the function 
\[ 
[-1,1] \to \mc H^{an}(U) : \ep \to \rho_\ep (N_w )
\] 
is analytic.
\end{lem}
\emph{Proof.} The first statement is \cite[Theorem 5.3]{Opd3},
but for the remainder we need to prove this anyway. It is clear
that $\rho_\ep$ restricts to an isomorphism between 
$C^{an}(U_\ep )$ and $C^{an}(U)$. For a simple reflection 
$s \in S_0$ corresponding to $\alpha \in F_1$ we have
\begin{equation}\label{eq:5.19}
\begin{array}{lll}
N_s + q(s)^{-\ep /2} & = & q^{-\ep /2} c_{\alpha ,\ep} 
(\imath^o_{s,\ep} + 1) \\
\rho_\ep (N_s ) & = & q^{\ep /2} (c_{\alpha ,\ep} \circ 
\sigma_\ep )(\imath^o_s + 1) - q(s)^{-\ep /2} \\
& = & q(s)^{(\ep -1)/2}(c_{\alpha ,\ep} \circ \sigma_\ep ) 
c_\alpha^{-1} \big( N_s + q(s)^{-1/2} \big) - q(s)^{-\ep /2} \\
& = & q(s)^{(\ep -1)/2} d_{\alpha ,\ep} \big( N_s + q(s)^{-1/2}
\big) - q(s)^{-\ep /2} \\
\end{array}
\end{equation}
By Lemma \ref{lem:5.13} such elements are analytic in $\ep \in 
[-1,1]$ and $u \in U$, so in particular they belong to 
$\mc H^{an}(U)$. Moreover, since every $d_{\alpha ,\ep}$ is 
invertible, the set $\{ \rho_\ep (N_w ) : w \in W_0 \}$ is a 
basis for $\mc H^{an}(U)$ as a $C^{an}(U)$-module. Therefore 
$\rho_\ep$ restricts to an isomorphism between $\mc H_\ep^{an}
(U_\ep)$ and $\mc H^{an}(U)$ for $\ep \neq 0$. 

For any $x \in X ,\, \rho_\ep (\theta_x ) = \theta_x \circ 
\sigma_\ep$ depends analytically on $\ep$, as a function on 
$U$. Combined with \eqref{eq:5.11} this shows that $\rho_\ep 
(N_w )$ is analytic in $\ep \in [-1,1]$ for any $w \in W$. 
Thus $\rho_0$ exists as a linear map. But, being a limit of 
algebra homomorphisms, it must also be multiplicative. 
$\qquad \Box$ \\[2mm]

Note that $\rho_0$ is never injective or surjective because
$\sigma_\ep$ is not. Moreover $\rho_{-1}$ is not the 
isomorphism \eqref{eq:5.40}. In general \eqref{eq:5.20} is not 
even rational, so $\rho_{-1}(N_s )$ cannot lie in $\mc H$. 
In the simple case $r_u (w^{-1} \alpha) = 1$ we have 
$\rho_{-1}(\theta_\alpha) = \theta_{-\alpha}$ and, by
\eqref{eq:5.18} and \eqref{eq:5.19}
\[
\rho_{-1}(N_s ) = N_s + q(s)^{-1/2} - q(s)^{1/2} = N_s^{-1}
\]
Usually the maps $\rho_\ep$ do not preserve the *, but this
can be fixed. For $\ep \in [-1,1]$ consider the element
\[
M_\ep = \rho_\ep (N_{w_0 ,\ep}^{-1})^* N_{w_0} 
\prod_{\alpha \in R_1^+} d_{\alpha ,\ep} \in \mc H^{an}(U)
\]
We will use $M_\ep$ to extend \cite[Corollary 5.7]{Opd3}. 
However, this result contained a small mistake: the 
construction of the element $A_\ep$ in \cite{Opd3} was 
unfortunately influenced by an inessential oversight. 
To correct this we replace it by $M_\ep$.

\begin{thm}\label{thm:5.15}
For all $\ep \in [-1,1]$ \inde{$M_\ep$} is invertible, positive 
and bounded. It has a positive square root in $\mc H^{an}(U)$ 
and the map $\ep \to M_\ep^{1/2}$ is analytic. 
\begin{align*}
& \tilde \rho_\ep : \mc H_\ep^{an}(U_\ep ) \to \mc H^{an}(U) \\
& \tilde \rho_\ep (h) = M_\ep^{1/2} \rho_\ep (h) M_\ep^{-1/2}
\end{align*}
is a homomorphism of topological *-algebras, and an isomorphism
if $\ep \neq 0$. For any $w \in W$ the function 
\[ 
[-1,1] \to \mc H^{an}(U) : \ep \to \tilde \rho_\ep (N_w )
\] 
is analytic.
\end{thm}
\emph{Proof.} 
By Lemmas \ref{lem:5.13} and \ref{lem:5.14} the $M_\ep$ are 
invertible, bounded and analytic in $\ep$. Consider, for 
$\ep \neq 0$, the automorphism $\mu_\ep$ of $\mc H^{me}(U)$ 
given by
\[
\mu_\ep (h) = \rho_\ep (\rho_\ep^{-1} (h)^* )^*
\]
On one hand, for $f \in C^{me}(U)$ we have by definition 
\eqref{eq:3.17} and the $W_0$-equivariance of $\sigma_\ep$
\begin{equation}
\begin{array}{lll}
\mu_\ep (f) & = & \rho_\ep ( (f \circ \sigma_\ep )^* )^* \\ 
& = & \rho_\ep \big( N_{w_0} (f^{-w_0} \circ \sigma_\ep )
N_{w_0 ,\ep}^{-1} \big)^* \\
& = & \rho_\ep (N_{w_0 ,\ep}^{-1})^* \big( f^{-w_0} \big)^*
\rho_\ep (N_{w_0 ,\ep})^* \\
& = & \rho_\ep (N_{w_0 ,\ep}^{-1})^* N_{w_0} f N_{w_0}^{-1}
\rho_\ep (N_{w_0 ,\ep})^* \\
& = & \rho_\ep (N_{w_0 ,\ep}^{-1})^* N_{w_0} \prod_{\alpha \in 
R_1^+} d_{\alpha ,\ep} f \prod_{\alpha \in R_1^+} 
d_{\alpha ,\ep}^{-1} N_{w_0}^{-1} \rho_\ep (N_{w_0 ,\ep})^* \\
& = & M_\ep f M_\ep^{-1}
\end{array}
\end{equation}
On the other hand, suppose that the simple reflections 
$s$ and $s' = w_0 s w_0 \in S_0$ correspond to $\alpha$ and 
$\alpha' = -w_0 \alpha \in F_1$. Using \eqref{eq:3.18} and 
\eqref{eq:3.19} we find
\begin{equation}
\begin{array}{lll}
M_\ep \imath_s^o M_\ep^{-1} & = & \rho_\ep (N_{w_0 ,\ep}^{-1}
)^* N_{w_0} \prod\limits_{\alpha \in R_1^+} d_{\alpha ,\ep} 
\imath_s^o \prod\limits_{\alpha \in R_1^+} d_{\alpha 
,\ep}^{-1} N_{w_0}^{-1} \rho_\ep (N_{w_0 ,\ep})^* \\
& = & \rho_\ep (N_{w_0 ,\ep}^{-1})^* N_{w_0} \imath_s^o
d_{\alpha, \ep}^{-1} d_{-\alpha ,\ep} N_{w_0}^{-1} 
\rho_\ep (N_{w_0 ,\ep})^* \vspace{2mm}\\
& = & \rho_\ep (N_{w_0 ,\ep}^{-1})^* N_{w_0} \imath_s^o {\ds 
\frac{c_{-\alpha}}{c_\alpha}  N_{w_0}^{-1} N_{w_0} 
\frac{c_{\alpha ,\ep} \circ \sigma_\ep}{c_{-\alpha ,\ep} \circ 
\sigma_\ep} } N_{w_0}^{-1} \rho_\ep (N_{w_0 ,\ep})^* \\
& = & \rho_\ep (N_{w_0 ,\ep}^{-1})^* (\imath^0_{s'})^* 
\left({\ds \frac{c_{\alpha' ,\ep} \circ \sigma_\ep}{c_{-\alpha' 
,\ep} \circ \sigma_\ep} }\right)^* \rho_\ep (N_{w_0 ,\ep})^* \\
& = & \left( \rho_\ep (N_{w_0 ,\ep}) {\ds \frac{c_{\alpha' ,\ep} 
\circ \sigma_\ep}{c_{-\alpha' ,\ep} \circ \sigma_\ep} } 
\imath^0_{s'} \rho_\ep (N_{w_0 ,\ep}^{-1}) \right)^* \\
& = & \rho_\ep \left( N_{w_0} {\ds \frac{c_{\alpha'}}{c_{
-\alpha'}} } \imath^o_{s',\ep} N_{w_0,\ep}^{-1} \right)^* \\
& = & \rho_\ep \left( (\imath^o_{s,\ep})^* \right)^*
\;=\; \mu_\ep (\imath^o_s )
\end{array}
\end{equation}
Since $C^{me}(U)$ and the $\imath^0_s$ generate 
$\mc H^{me}(U)$, we conclude that 
\[
\mu_\ep (h) = M_\ep h M_\ep^{-1} \qquad 
\forall h \in \mc H^{me}(U)
\]
Now we can see that
\[ 
\begin{array}{lll}
\rho_\ep \big( N_{w_0 ,\ep}^{-1} \big)^* & 
= & \rho_\ep \big( (N_{w_0 ,\ep}^*)^{-1} \big)^*
\;=\; \rho_\ep \big( (N_{w_0 ,\ep}^{-1})^* \big)^* \\
& = & \mu_\ep \big( \rho_\ep (N_{w_0 ,\ep}^{-1} ) \big) 
\vspace{2mm}
\;=\; M_\ep \rho_\ep (N_{w_0 ,\ep}^{-1} ) M_\ep^{-1} \\
N_e & = & M_\ep^{-1} \rho_\ep (N_{w_0 ,\ep}^{-1})^* N_{w_0} 
\prod\limits_{\alpha \in R_1^+} d_{\alpha ,\ep} 
\;=\; \rho_\ep (N_{w_0 ,\ep}^{-1}) M_\ep^{-1} N_{w_0} 
\prod\limits_{\alpha \in R_1^+} d_{\alpha ,\ep} \\
M_\ep & = & N_{w_0} \prod\limits_{\alpha \in R_1^+} d_{\alpha 
,\ep} \rho_\ep (N_{w_0 ,\ep}^{-1}) \;=\; \big( \rho_\ep 
(N_{w_0 ,\ep}^{-1})^* \big( N_{w_0} \prod\limits_{\alpha 
\in R_1^+} d_{\alpha ,\ep} \big)^* \big)^* \\
& = & \big( \rho_\ep (N_{w_0 ,\ep}^{-1})^*  N_{w_0} 
\prod\limits_{\alpha \in R_1^+} d_{\alpha ,\ep} \big)^* 
\;=\; M_\ep^*
\end{array}
\]
Thus the elements $M_\ep$ are Hermitian $\forall \ep \neq 0$. By
continuity in $\ep \; M_0$ is also Hermitian. Moreover they are 
all invertible, and $M_1 = N_e$, so they are in fact strictly 
positive. We already knew that the element $\ep \to M_\ep$ of
\[
C^{an} ([-1,1] ; \mc H^{an}(U)) \cong C^{an} ([-1,1] \times U)
\otimes_{\mc A} \mc H
\]
is bounded, so we can construct its square root using 
holomorphic functional calculus in the Fr\'echet Q-algebra 
$C_b^{an}([-1,1] \times U) \otimes_{\mc A} \mc H$. 
This ensures that $\ep \to M_\ep^{1/2}$ is still analytic.
Finally, for $\ep \neq 0$
\begin{equation}
\begin{array}{lll}
\tilde \rho_\ep (h)^* & = & 
\left( M_\ep^{1/2} \rho_\ep (h) M_\ep^{-1/2} \right)^* \\
& = & M_\ep^{-1/2} \rho_\ep (h)^* M_\ep^{1/2} \\
& = & M_\ep^{-1/2} \mu_\ep (\rho_\ep (h^*)) M_\ep^{1/2} \\
& = & M_\ep^{1/2} \rho_\ep (h^* ) M_\ep^{-1/2} 
\;=\; \tilde \rho_\ep (h^* )
\end{array}
\end{equation}
Again this extends to $\ep = 0$ by continuity. $\qquad \Box$
\\[2mm]

Let \inde{Rep$\,(C_r^* (\mc R ,q))$} and \inde{Rep$\,(\mc S 
(\mc R ,q))$} be the categories of finite dimensional 
representations of $C_r^* (\mc R ,q)$ and of $\mc S (\mc R ,q))$. 
We define \index{$\mr{Rep}_U (C_r^* (\mc R ,q))$}
\index{$\mr{Rep}_U (\mc S (\mc R ,q))$}
\begin{equation}
\begin{array}{lllll}
\mr{Rep}_U (C_r^* (\mc R ,q)) & = & \mr{Rep}\, (C_r^* (\mc R ,q)) 
& \cap & \mr{Rep}_U (\mc H (\mc R ,q)) \\
\mr{Rep}_U (\mc S (\mc R ,q)) & = & \mr{Rep}\, (\mc S (\mc R ,q)) 
& \cap & \mr{Rep}_U (\mc H (\mc R ,q))
\end{array}
\end{equation}
Recall that $\overline{\mc H_\ep^t}$ is the residual algebra 
of $\mc H_\ep$ at $t \in T$, whose (finite dimensional) 
representations are precisely 
$\mr{Rep}_{W_0 t}(C_r^* (\mc R ,q^\ep))$.

\begin{lem}\label{lem:5.16}
For $\ep \in [-1,1]$ the map $\tilde \rho_\ep$ induces a
"scaling map" \index{rhooe@$\overline{\rho_\ep}$} 
\index{sigmate@$\tilde \sigma_\ep$}
\[
\tilde \sigma_\ep : \mr{Rep}_{W_0 u}(\mc H ) \to
\mr{Rep}_{W_0 \sigma_\ep (u)} (\mc H_\ep )
\]
which preserves unitarity and is a bijection if $\ep \neq 0$.

For $\ep < 0 \; \tilde \sigma_\ep$ exchanges tempered and
anti-tempered modules. For $\ep \geq 0 \; \tilde \sigma_\ep$ 
preserves (anti-)temperedness and $\tilde \rho_\ep$ descends 
to a *-homomorphism
\[
\overline{\rho_\ep} : \overline{\mc H_\ep^{\sigma_\ep (u)}} 
\to \overline{\mc H^u}
\]
which is an isomorphism if $\ep > 0$.
\end{lem}
\emph{Proof.}
If $\pi \in \mr{Rep}(\mc H )$ and $\ep \neq 0$ then by 
construction $t \in T$ is an $\mc A$-weight of $\pi$ if and 
only if $\sigma_\ep (t)$ is an $\mc A_\ep$-weight of $\pi \circ 
\tilde \rho_\ep$. The $\mc A_0$-weights of $\pi \circ \tilde 
\rho_0$ are all contained in $W_0 r \subset T_u$. Therefore
\begin{equation}
\tilde \sigma_\ep : \pi \to \pi \circ \tilde \rho_\ep
\end{equation}
maps $\mr{Rep}_{W_0 u}(\mc H )$ to $\mr{Rep}_{W_0 \sigma_\ep 
(u)}(\mc H_\ep )$, for any $u \in U$ and $\ep \in [-1,1]$.
Since $\tilde \rho_\ep$ is a *-homomorphism and a bijection
(for $\ep \neq 0) \: \tilde \sigma_\ep$ has the desired 
properties for such $\ep$.

Moreover for $x \in X$ we have $| \sigma_\ep (t)(x)| = 
|t(x) |^\ep$, which proves the assertions about temperedness.
If $\ep \geq 0$ and $\pi$ extends to $C_r^* (\mc R ,q)$ then 
it is unitary and tempered by Proposition \ref{prop:5.6}. 
Therefore $\pi \circ \tilde \rho_\ep$ is tempered and 
completely reducible, and Corollary \ref{cor:3.20} assures 
that it extends to $C_r^* (\mc R ,q^\ep )$. This implies
\[
\tilde \rho_\ep \big( \mr{Rad}_{\sigma_\ep (u),\ep} \big) 
\subset \mr{Rad}_u
\]
so we get a *-homomorphism
\[
\overline{\rho_\ep} : \mc H_\ep / \mr{Rad}_{\sigma_\ep (u),\ep} 
= \overline{\mc H_\ep^{\sigma_\ep (u)}} \to \mc H / \mr{Rad}_u
= \overline{\mc H^u}
\]
If $\ep > 0$ then we can follow the same reasoning for 
${\tilde \rho_\ep}^{-1}$, so $\overline {\rho_\ep}$ is an 
isomorphism. $\qquad \Box$.
\\[2mm]

Assume once more that $0 \neq \ep \in [-1,1]$. Then $r \to 
\sigma_\ep (r)$ is a bijection between the residual points 
for $(\mc R ,q)$ and those for $(\mc R ,q^\ep)$. The groupoid 
of intertwiners $\mc W$ is independent of $q$, and it acts 
on the set of $\mc H_\ep$-representations of the form 
$\pi (P,\delta' ,t)$ with $\delta'$ irreducible but not 
necessarily discrete series. The definitions \eqref{eq:3.21} 
and \eqref{eq:3.22} also makes sense in this case. If we 
realize the representation
\[
\pi \big( P,\tilde \sigma_\ep (\delta ),t \big) \quad\mr{on}
\quad \mc H \big( W^P \big) \otimes V_\delta \quad\mr{as}\quad 
\mr{Ind}_{\mc H_\ep^P}^{\mc H_\ep}\big( \delta \circ \tilde 
\rho_\ep \circ \phi_{t,\ep} \big)
\]
then we get homomorphisms \index{$\mc F'_\ep$}
\begin{equation}\label{eq:5.26}
\begin{aligned}
& \mc F'_\ep : \mc H (\mc R ,q^\ep) \to \bigoplus_{(P,\delta) 
\in \Delta'} \mc O \big( T^P \big) \otimes \mr{End} 
\big(\mc H (W^P ) \otimes V_\delta \big) \\
& \mc F'_\ep (h) (P,\delta ,t) = \pi (P, \tilde \sigma_\ep 
(\delta ),t)(h)
\end{aligned}
\end{equation}
(Notice that this is also defined for $\ep = 0$.)
The isotropy groups $\mc W_{\tilde \sigma_\ep (\delta)}$ for
various $\ep$'s may be identified, so the image of 
\eqref{eq:5.26} consists of sections that are invariant under 
a certain action of $\mc W_\delta$. We add a subscript $\ep$ 
to indicate which action we consider.

Recall from \eqref{eq:5.14} that $\Gamma_\delta$ is a Schur 
extension of $\mc W_\delta$. There are uniquely determined 
\inde{$u_{\gamma ,\ep}$}$\in A_\delta^\times$ such that, just 
as in \eqref{eq:5.21}, we can write the associated 
$\Gamma_\delta$-action as
\begin{equation}
\gamma_\ep (a) = u_{\gamma ,\ep} \gamma a \gamma^{-1} 
u_{\gamma ,\ep}^{-1} \qquad 
\mr{in}\; A_\delta \rtimes \Gamma_\delta
\end{equation}

\begin{lem}\label{lem:5.17}
The elements $u_{\gamma ,\ep}$ depend analytically on $\ep$ and
$u_{\gamma ,0} = \lim_{\ep \to 0} u_{\gamma ,\ep}$ exists.
Moreover $u_{\gamma ,\ep}$ is unitary $\forall \ep \in [-1,1] 
,\, \gamma \in \Gamma_\delta$.
\end{lem}
\emph{Proof.}
Let $\gamma$ be a lift of $k n \in \mc W_\delta$ and write
$\xi = \gamma (P,\tilde \sigma_\ep (\delta),t) \in \Xi_{u,\ep}$.
By definition \eqref{eq:3.21} for $h_0 \in \mc H (W^P ) ,\, 
v \in V_\delta$
\begin{equation}\label{eq:5.22}
u_{\gamma ,\ep}(\xi)(h_0 \otimes v) = h_0 \cdot_{q^\ep} 
\imath^o_{n^{-1},\ep}(t) \otimes \tilde \delta_{\gamma ,\ep}(v)
\end{equation}
where $\tilde \delta_{\gamma ,\ep} \in U(V_\delta )$. More 
precisely, $\tilde \delta_{\gamma ,\ep}$ is a multiple of a 
map $L_\ep := \tilde \delta_{kn,\ep} \in U(V_\delta )$ that
satisfies
\begin{equation}\label{eq:5.23}
\delta (\tilde \rho_\ep h_1 ) = L_\ep^{-1} \delta (\rho_\ep 
\psi_k \psi_n h_1 ) L_\ep \qquad h_1 \in \mc H_{P,\ep}
\end{equation}
This $L_\ep$ is only defined up to scalars, but we will show 
that these can be chosen such that $L_0 = \lim_{\ep \to 0} 
L_\ep$ exists. Since $V_\delta$ is a finite dimensional vector 
space, every automorphism of End$(V_\delta )$ is inner, and 
$\mr{Aut}(\mr{End}(V_\delta)) \cong PGL(V_\delta)$. Because 
$GL (V_\delta ) \to PGL(V_\delta)$ is a fiber bundle it 
suffices to show that $L_0$ exists as a \inde{projective linear 
map}. Writing $h_1 = \rho_\ep^{-1} h_2$ and rearranging some 
terms transforms \eqref{eq:5.23} into
\[
\delta (M_\ep^{-1/2}) L_\ep \delta (M_\ep^{1/2}) \delta (h_2 ) 
\delta (M_\ep^{-1/2}) L_\ep^{-1} \delta (M_\ep^{1/2}) 
= \delta (\rho_\ep \psi_k \psi_n \rho_\ep^{-1} h_2 )
\]
If we can show that everything else in this equation is 
analytic in $\ep$ and has a well-defined limit as $\ep \to 0$, 
then the same must hold for $L_\ep \in PGL (V_\delta )$. By 
Theorem \ref{thm:5.15} we know this already for 
$M_\ep^{\pm 1/2}$. For any $f \in C^{an}(U)$ we have 
$\rho_\ep \psi_k \psi_n \rho_\ep^{-1} f = \psi_k \psi_n f$ 
by the $W_0$-equivariance of $\sigma_\ep$. For the simple 
reflection $s$ associated to $\alpha \in F_1$ we have 
\[
\begin{aligned}
\rho_\ep \psi_k \psi_n \rho_\ep^{-1} & \big( N_s + q(s)^{-1/2} 
\big) \; = \; \rho_\ep \psi_k \psi_n \big( q(s)^{(1-\ep )/2} 
(c_\alpha \circ \sigma_{1/\ep} c_{\alpha ,\ep}^{-1} 
(N_s + q(s)^{-\ep /2}) \big) \\
& = \; \rho_\ep \psi_k \big( q(s)^{(1-\ep )/2} (c_{n \alpha} 
\circ \sigma_{1/\ep}) c_{n \alpha ,\ep}^{-1} (N_{n s n^{-1}} + 
q(s)^{-\ep /2}) \big) \\
& = \; \psi_k \big( q(s)^{(1-\ep )/2} c_{n \alpha} (c_{n \alpha 
,\ep}^{-1} \circ \sigma_\ep ) \big) \rho_\ep \big( N_{n s n^{-1}} 
+ q(s)^{-\ep /2} \big)
\end{aligned}
\]
By Lemmas \ref{lem:5.13} and \ref{lem:5.14} the last expression 
has the required properties.

Now we turn our attention to the other parts of \eqref{eq:5.22}.
For $\ep > 0$ Lemma \ref{lem:5.16} guarantees that $\tilde 
\sigma_\ep (\delta )$ is discrete series. Therefore $u_{\gamma 
,\ep}$ is unitary and $\imath^o_{n^{-1},\ep}$ cannot have a pole 
at $t$. From the explicit definition \eqref{eq:3.23} we see that 
$\imath^o_{n^{-1},\ep}$ is regular at $t$ for any $\ep \in 
[-1,1]$, and that $\lim_{\ep \to 0}\imath^o_{n^{-1},\ep}(t) = 
N_{n^{-1}}$. By Proposition \ref{prop:5.7} this implies
\[ 
\lim_{\ep \to 0} h_0 \cdot_{q^\ep} \imath^o_{n^{-1},\ep}(t) =
h_0 \cdot_{q^0} N_{n^{-1}}
\]
Putting things together we conclude that $L_\ep ,\, \tilde 
\delta_{\gamma ,\ep}$ and $u_{\gamma ,\ep}$ are analytic in 
$\ep$ and have well-defined limits as $\ep \to 0$, all as 
projective linear maps. However, we already agreed that we may 
assume that this even holds for $L_\ep$ as a linear map. But by 
\cite[\S 53]{CuRe} the way to lift this representation from
$\mc W_\delta$ to $\Gamma_\delta$ is completely determined by 
the cocycle
\[
\mc W_\delta \times \mc W_\delta \to \mh C^\times :
(g_1 ,g_2 ) \to \tilde \delta_{g_1 ,\ep} \tilde \delta_{g_2 
,\ep} \tilde \delta_{g_2^{-1} g_1^{-1},\ep}^{-1}
\]
This cocycle is continuous in $\ep$ and $\Gamma_\delta$ is 
finite, so the way to lift is independent of $\ep$. Therefore 
$\lim_{\ep \to 0} \tilde \delta_{\gamma ,\ep}$ and $\lim_{\ep 
\to 0} u_{\gamma ,\ep}$ exist even as linear maps. It also 
follows that $u_{\gamma ,\ep}$ is analytic in $\ep$ which, in 
combination with its unitarity $\forall \ep > 0$, shows that 
it is unitary $\forall \ep \in [-1,1]. \qquad \Box$
\\

\begin{lem}\label{lem:5.18}
The pre-$C^*$-algebras $A_\delta^{\mc W_{\delta ,\ep}} = 
A_\delta^{\Gamma_{\delta ,\ep}}$ are all isomorphic, by
isomorphisms that are piecewise analytic in 
$\ep ,\ep' \in [-1,1]$.
\end{lem}
\emph{Proof.}
From Lemma \ref{lem:5.17} we get an analytic path of projections
\begin{equation}
[-1,1] \to A_\delta \rtimes \Gamma_\delta : \ep \to
p_\delta (u_\ep ) := |\Gamma_\delta |^{-1} \sum_{\gamma \in 
\Gamma_\delta} u_{\gamma ,\ep} \gamma
\end{equation}
Like in \eqref{eq:5.24} the map 
\[
A_\delta^{\Gamma_{\delta ,\ep}} \to p_\delta (u_\ep )
(A_\delta \rtimes \Gamma_\delta ) p_\delta (u_\ep ) :
a \to p_\delta (u_\ep ) a p_\delta (u_\ep )
\]
is an isomorphism of pre-$C^*$-algebras. If we apply 
\cite[Lemma 1.15]{Phi2} we see that the $p_\delta (u_\ep )$ 
are all conjugate, by elements depending continuously on $\ep$.
To show analyticity we construct these elements explicitly,
using the recipe of \cite[Proposition 4.32]{Bla}.
For $\ep ,\ep' \in [-1,1]$ consider the element 
\[
z(\delta ,\ep ,\ep' ) := (2 p_\delta (u_\ep' ) - 1)(2 p_\delta 
(u_\ep ) - 1) + 1 \in A_\delta \rtimes \Gamma_\delta
\]
Clearly this is analytic in $\ep$ and $\ep'$ and
\[
p_\delta (u_\ep )z(\delta ,\ep ,\ep' ) = 2 p_\delta (u_\ep' )
p_\delta (u_\ep ) = z(\delta ,\ep ,\ep' )p_\delta (u_\ep )
\]
Moreover if $\norm{\cdot}$ is the norm of the enveloping 
$C^*$-algebra 
\[
C_\delta := C \big( T_u^P ; \mr{End} \big(\mc H (W^P ) 
\otimes V_\delta \big) \big) \rtimes \Gamma_\delta
\]
of $A_\delta \rtimes \Gamma_\delta$ and 
$\norm{p_\delta (u_\ep ) - p_\delta (u_\ep' )} < 2$ then
\[
\begin{array}{lll}
\norm{z(\delta ,\ep ,\ep') - 2} & = & \norm{4 p_\delta 
(u_\ep' ) p_\delta (u_\ep ) - 2 p_\delta (u_\ep ) - 
2 p_\delta (u_\ep )} \\ 
& = & \norm{-2 \big( p_\delta (u_\ep ) - p_\delta (u_\ep' ) 
\big)^2} \quad < \quad 2
\end{array}
\]
so $z(\delta ,\ep ,\ep')$ is invertible in $C_\delta$. But 
$A_\delta \rtimes \Gamma_\delta$ is holomorphically closed in 
$C_\delta$, so $z(\delta ,\ep ,\ep')$ is also invertible in 
this Fr\'echet algebra. Moreover, because the Fr\'echet 
topology on $A_\delta \rtimes \Gamma_\delta$ is finer than 
the topology coming from $\norm{\cdot}$, there exists an open 
interval $I_\ep$ containing $\ep$ such that $z(\delta ,\ep 
,\ep')$ is invertible $\forall \ep' \in I_\ep$. For such $\ep'$ 
we construct the unitary element
\[
u(\delta ,\ep ,\ep') := z(\delta ,\ep ,\ep') 
|z(\delta ,\ep ,\ep')|^{-1}
\]
By construction the map
\[
p_\delta (u_\ep ) (A_\delta \rtimes \Gamma_\delta ) 
p_\delta (u_\ep ) \to p_\delta (u_\ep' )
(A_\delta \rtimes \Gamma_\delta ) p_\delta (u_\ep' ) :
x \to u(\delta ,\ep ,\ep') x u(\delta ,\ep ,\ep')^{-1}
\]
is an isomorphism of pre-$C^*$-algebras. The composite map
$A_\delta^{\Gamma_{\delta ,\ep}} \to 
A_\delta^{\Gamma_{\delta ,\ep'}}$ is given by
\begin{equation}\label{eq:5.25}
x \to |\Gamma_\delta | \left[ u(\delta ,\ep ,\ep') p_\delta 
(u_\ep' ) x p_\delta (u_\ep' ) u(\delta ,\ep ,\ep')^{-1}
\right]_e
\end{equation}
which is analytic in $\ep$ and $\ep'$ because $p_\delta (u_\ep 
)$ is. Now we pick a finite cover $\{I_{\ep_i}\}_{i=1}^m$ of 
$[-1,1]$. Then for any $\ep ,\, \ep' \in [-1,1]$ an isomorphism
between $A_\delta^{\Gamma_{\delta ,\ep}}$ and $A_\delta^{
\Gamma_{\delta ,\ep'}}$ can be obtained by composing at most 
$m$ isomorphisms of the type \eqref{eq:5.25}. $\qquad \Box$
\\[2mm]

The constructions in this section lead to the following

\begin{cor}\label{cor:5.19}
There exists a collection of injective *-homomorphisms
\[
\phi_\ep : \mc H (\mc R ,q^\ep ) \to \mc S (\mc R ,q)
\qquad \ep \in [-1,1]
\]
such that 
\begin{enumerate}
\item for $\ep < 0$ the map \index{fie@$\phi_\ep$}
\[
\mr{Rep}(\mc S (\mc R ,q)) \to \mr{Rep}(\mc H (\mc R ,q^\ep )) 
: \pi \to \pi \circ \phi_\ep
\]
is a bijection from tempered $\mc H$-representations to 
anti-tempered\\ $\mc H_\ep$-representations
\item $\forall (P,\delta ,t) \in \Xi_u$ the representation
$\pi (P,\delta ,t) \circ \phi_\ep$ is equivalent with
$\pi_\ep (P,\tilde \sigma_\ep (\delta),t)$
\item $\phi_1$ is the canonical embedding
\item $\phi_\ep (N_w ) = N_w \; \forall w \in Z(W)$
\item $\forall h \in \mc H (\mc R )$ the function 
\[
[-1,1] \to \mc S (\mc R ,q) :\ep \to \phi_\ep (h)
\]
is piecewise analytic, and in particular analytic at $0$.
\end{enumerate}
\end{cor}
\emph{Proof.}
By Lemma \ref{lem:5.17} the image of $\mc F'_0$ is invariant
under the action of $\mc W_{\delta ,0}$. So we can define
\begin{equation}
\phi_\ep = \mc F'^{-1} \circ \zeta_\ep \circ \mc F'_\ep
\end{equation}
where $\zeta_\ep = \bigoplus_{(P,\delta ) \in \Delta'} 
\zeta_{\ep ,\delta}$ and 
\[
\zeta_{\ep ,\delta} : A_\delta^{\mc W_{\delta ,\ep}}
\to A_\delta^{\mc W_\delta}
\]
is the isomorphism from Lemma \ref{lem:5.18}. Now 2 and 3
are valid by construction, 4 follows from the 
observation that the $Z(W)$-character of $\pi_\ep (P,\tilde 
\sigma_\ep (\delta),t)$ is equal to $t \big|_{Z(W)}$, for 
every $\ep \in [-1,1]$, and 1 is a consequence of
Lemma \ref{lem:5.16} and Proposition \ref{prop:3.13}. 
Finally, for 5 we use \eqref{eq:5.11}, Theorem \ref{thm:5.15} 
and Lemma \ref{lem:5.18}. From the proof of this lemma we see
that we can arrange that $\ep \to \phi_\ep (h)$ is analytic
at 0.

As concerns the injectivity of $\phi_\ep$, note that $\pi_\ep 
(P, \tilde \sigma_\ep (\delta_\es ),t) = I_{t,\ep}$ 
is a principal series representation for $\mc H_\ep$. 
So if $h \in \ker (\phi_\ep )$, then $h$ acts as 0 on all 
unitary principal series. By Lemma \ref{lem:3.12} we must 
have $h = 0. \qquad \Box$ 
\\[2mm]

We do not know whether $\phi_\ep (\mc H_\ep ) \subset \mc H$,
for two reasons : $\zeta_{\ep ,\delta}$ need not preserve
polynomiality, and not every polynomial section is in the 
image of $\mc F'$.

\begin{thm}\label{thm:5.20}
For $\ep \in [0,1]$ there exist homomorphisms of 
pre-$C^*$-algebras 
\[
\begin{array}{ccccc}
\phi_\ep & : & \mc S (\mc R ,q^\ep ) & \to & \mc S (\mc R ,q) \\
\phi_\ep & : & C_r^* (\mc R ,q^\ep ) & \to & C_r^* (\mc R ,q)
\end{array}
\]
with the properties
\begin{enumerate}
\item $\phi_\ep$ is an isomorphism if $\ep > 0$
and $\phi_0$ is injective
\item $\forall (P,\delta ,t) \in \Xi_u$ the representation
$\pi (P,\delta ,t) \circ \phi_\ep$ is equivalent with
$\pi_\ep (P,\tilde \sigma_\ep (\delta),t)$
\item $\phi_1$ is the identity
\item $\phi_\ep (h) = h \; \forall h \in \mc S (Z(W))$
\item $\ep \to \phi_\ep (h)$ is continuous $\forall h \in 
\mc S (\mc R)$
\end{enumerate}
\end{thm}
\emph{Proof.}
For any $(P,\delta ) \in \Delta$ the representation $\tilde 
\sigma_\ep (\delta )$, although not necessarily irreducible if
$\ep = 0$, is certainly completely reducible, being unitary.
Hence by Lemma \ref{lem:5.16} every irreducible constituent 
$\pi_1$ of $\tilde \sigma_\ep (\delta )$ is a direct summand of 
\[
\mr{Ind}_{\mc H_{\ep ,P}^{P_1}}^{\mc H_{\ep ,P}} (P_1 ,\delta_1
,\phi_{t_1 ,\ep})
\]
for certain $P_1 \subset P ,\, \delta_1 \in \Delta_\ep$ and
\[
t_1 \in \mr{Hom}_{\mh Z} \left( (X_P)^{P_1} , S^1 \right) =
\mr{Hom}_{\mh Z} \left( X / (X \cap (P^\vee )^\perp + \mh Q P_1 
), S^1 \right) \subset T_u
\]
Consequently $\pi (P, \pi_1 ,t)$ is a direct summand of
\[
\mr{Ind}_{\mc H_\ep^P}^{\mc H_\ep} \Big( 
\mr{Ind}_{\mc H_{\ep ,P}^{P_1}}^{\mc H_{\ep ,P}} (\delta_1 
\circ \phi_{t_1 ,\ep}) 
\circ \phi_{t,\ep} \Big) = \pi_\ep (P_1 ,\delta_1 ,t t_1 )
\]
In particular every matrix coefficient of $\pi_\ep (P, \tilde
\sigma_\ep (\delta ),t)$ appears in the Fourier transform of
$\mc H_\ep$, and \eqref{eq:5.26} extends to the respective
Schwartz and $C^*$-completions, as required. By Lemma 
\ref{lem:5.16} and \eqref{eq:5.12} these maps are isomorphisms 
if $\ep > 0$. 

In the same way as in the proof of Corollary \ref{cor:5.19} 
we can see that $\phi_0$ remains injective: every irreducible 
tempered $\mc H (\mc R ,q^0 )$-representation is a quotient of 
a unitary principal series, so any element of 
$C_r^* (\mc R ,q^0 )$ that vanishes on all unitary principal 
series is 0. Furthermore properties 2. and 4. are direct
consequences of Corollary \ref{cor:5.19}.

Finally, if $x = \sum_{w \in W} x_w N_w \in \mc S (\mc R)$ 
then this sum converges uniformly to $x$. For every partial sum 
$x'$ the map $\ep \to \phi_\ep (x')$ is continuous by Corollary
\ref{cor:5.19}, so this also holds for $x$ itself. 
$\qquad \Box$ \\[3mm]

Although there were quite a few arbitrary choices involved in
constructing $\phi_\ep$, the homotopy class of these maps is 
well-defined:

\begin{lem}\label{lem:5.21}
The construction of $(\phi_\ep )_{\ep \in [0,1]}$ is unique up
to a homotopy of objects with the properties of Theorem 
\ref{thm:5.20}.
\end{lem}
\emph{Proof.}
Let us inventorize all the choices we made in the above 
construction. Already in \eqref{eq:5.12} we chose a set of 
representatives $\Delta'$ of $\Delta / \mc W$. This implicitly 
fixed realizations $V_\delta$ of the discrete series 
representations $\delta \in \Delta'$, but since we never used 
a basis of $V_\delta$, the construction does not really depend 
on this vector space. Then we used intertwiners $\pi(g,\xi) 
,\:g \in \mc W_\delta$ that were defined only up to scalars, 
but this ambiguity dropped out when we lifted things to 
$\Gamma_\delta$. Finally we chose a (monotonuous) sequence 
$(\ep_i)_{i=0}^m$ such that $\ep_0 = \ep ,\: \ep_m = 1$ and 
every isomorphism $\mc S (\mc R,q^{\ep_i}) \to 
\mc S(\mc R,q^{\ep_{i+1}})$ involved only one map of the type 
\eqref{eq:5.25} for every $\delta$.

Suppose we take another sequence $(\ep'_i)_{i=0}^{m'}$. 
Allowing some of the $\ep_i$ and $\ep'_j$ to coincide, we may 
assume that $m' = m$. Then we can continuously deform the first 
sequence into the second. Since the elements $u(\delta ,\ep 
,\ep')$ depend analytically on $\ep$ and $\ep'$, this will give 
us a path of isomorphisms.

To see what happens when we use some $\Delta''$ instead of 
$\Delta'$ is more difficult. It is clearly enough to investigate 
the effect on only one component of $\Xi_u$, corresponding to 
$(P,\delta) \in \Delta'$ and to $(P',\delta') \in \Delta''$. 
By the previous paragraph we may also restrict ourselves to 
$\ep$ and $\ep'$ with $|\ep - \ep'|$ "sufficiently" small. 
Then the original isomorphism is given by
\begin{equation}\label{eq:5.27}
\phi_{\ep ,\ep'} : a \to |\Gamma_\delta | 
\left[ p_\delta (u_{\ep'}) u( \delta ,\ep ,\ep' ) p_\delta 
(u_\ep ) a p_\delta (u_\ep ) u( \delta ,\ep ,\ep' )^{-1} 
p_\delta (u_{\ep'}) \right]_e
\end{equation}
and the alternative by
\begin{equation}
\phi'_{\ep,\ep'} : a' \to |\Gamma_{\delta'} | 
\left[ u(\delta',\ep ,\ep') p_{\delta'}(u_\ep ) a' 
p_{\delta'}(u_\ep ) u(\delta' ,\ep ,\ep')^{-1} \right]_e
\end{equation}
In both formulas $[\cdot ]_e$ means taking the coefficient 
at the identity element in some group algebra. To compare 
$\phi_{\ep ,\ep'}$ and $\phi'_{\ep ,\ep'}$ we take 
$g \in \mc W_{\delta \delta'} ,\: \xi = (P,\delta,t) \in \Xi_u$ 
and evaluate $\phi'_{\ep ,\ep'} a(\xi )$, assuming that $a$ 
and $\phi'_{\ep ,\ep'}(a)$ are $g$-invariant:
\begin{equation}\label{eq:5.28}
\begin{split}
& \phi'_{\ep ,\ep'} a(\xi ) \quad = \quad g^{-1}(\phi'_{\ep ,\ep'} 
(g a))(\xi) \quad = \\
& \pi_\ep' (g^{-1},g \xi) \: \big( \phi'_{\ep ,\ep'} (g a) \big) 
(g \xi) \: \pi_\ep' (g^{-1},g \xi)^{-1} \quad =\\
& |\Gamma_{\delta'}| \: \pi_\ep' (g^{-1} , g \xi) \;
\big[ u(\delta' ,\ep ,\ep' ) p_{\delta'}(u_\eta) g(a) 
p_{\delta'}(u_\ep ) u(\delta' ,\ep ,\ep')^{-1} \big]_e (g \xi) \;
\pi_\ep' (g^{-1} , g \xi)^{-1} \quad = \\
& |\Gamma_{\delta'}| \: \pi_\ep' (g^{-1},g\xi) \; 
\big[ u(\delta' ,\ep ,\ep' ) (g \xi) p_{\delta'}(u_\ep )(g \xi) 
\pi_\ep (g,\xi ) a(\xi ) \; \circ \\
& \quad \pi_\ep (g,\xi )^{-1} p_{\delta'}(u_\ep )(g \xi ) u(\delta' 
,\ep ,\ep' )^{-1}(g \xi ) \big]_e \; \pi_\ep' (g^{-1},g \xi )^{-1}
\end{split}
\end{equation}
To compare $\Gamma_\delta$ and $\Gamma_{\delta'}$ we have to make 
the Schur extension \eqref{eq:5.14} functorial. In general it 
is not even unique, but the recipe in \cite[\S 53]{CuRe} always 
works. For the sake of the argument we might temporarily redefine 
the Schur extension to be the result of this construction. As a 
bonus it is easily seen that for $\gamma' \in \Gamma_{\delta'}$ 
the element $\gamma = g^{-1} \gamma' g \in \Gamma_\delta$ is 
well-defined. But then necessarily
\[
\pi (g,\xi )^{-1} u_{\gamma'}(g \xi ) \pi (g,\xi ) 
= u_\gamma (\xi )
\]
since both sides are proportional and satisfy the cocycle 
relations \eqref{eq:5.15}. In particular
\[
\begin{array}{lll}
p_{\delta'}(u)(g \xi ) \pi (g,\xi ) & \!\!\!=\!\! & 
\pi (g,\xi ) p_\delta (u) (\xi )\\
\pi_\ep' (g^{-1},g \xi ) u(\delta' ,\ep ,\ep') (g \xi ) 
\pi_\eta (g,\xi ) & \!\!\!=\!\! & (2 p_{\delta'}(u_{\ep'})(\xi ) 
- 1) b (\xi ) (2 p_{\delta'}(u_\ep )(\xi ) - 1) + b (\xi )
\end{array}
\]
We denote the last expression by $u(b)(\xi )$, where 
$b(\xi ) = \pi_\ep' (g^{-1}, g \xi ) \pi_\ep (g,\xi )$.
In general elements of the type $u(b)$ in a $C^*$-algebra 
are invertible if 
\[
\norm{p_{\delta'}(u_\ep ) - p_{\delta'}(u_{\ep'})} + 
4 \norm{b-1} < 1
\] 
and they satisfy 
\[
p_{\delta'}(u_{\ep'}) u(b) = p_{\delta'}(u_{\ep'}) u(b) 
p_{\delta'}(u_\ep ) = u(b) p_{\delta'}(u_\ep )
\]
With this knowledge we can rewrite \eqref{eq:5.28} as
\[
\phi'_{\ep ,\ep'} a(\xi ) = |\Gamma_\delta | 
[p_\delta (u_{\ep'}) u(b) p_\delta (u_\ep ) ap_\delta (u_\ep ) 
u(b)^{-1} p_\delta (u_{\ep'}) ]_e (\xi )
\]
but now $e$ is an element of $\Gamma_\delta$ instead of 
$\Gamma_{\delta'}$. Comparing this to \eqref{eq:5.27} we find 
that the only difference is that $u(\delta ,\ep ,\ep')$ has 
been replaced by this $u(b)$. The intertwiners $\pi_{\ep'}(g,\xi)$ 
depend analytically on $\ep'$, so we can always find a path from 
$u(b)$ to $u(\delta ,\ep ,\ep' )$ consisting only of elements of 
the same type. (This step might require a subdivision of the 
interval between $\ep$ and $\ep'$, but that is no problem.) 
The corresponding isomorphisms form a path from 
$\phi'_{\ep ,\ep'}$ to $\phi_{\ep ,\ep'}. \qquad \Box$
\\[2mm]

\section{$K$-theoretic conjectures}
\label{sec:5.4}

We saw in Section \ref{sec:5.2} that the multiplication in
$\mc S (\mc R ,q)$ varies continuously with $q$. Since a class 
in $K$-theory is rigid under small perturbations, it is natural 
to expect that the $K$-groups of $\mc S (\mc R ,q)$ are
independent of $q$. We reformulate this conjecture in terms of
the map $\phi_0$ and show that it implies some other important
conjectures.

Our main tools are the scaling maps $\phi_\ep$ from Theorem
\ref{thm:5.20}. From Theorems \ref{thm:5.20} and \ref{thm:2.19}, 
Lemma \ref{lem:5.21} and the homotopy invariance of $K$-theory 
we see that for all $\ep \in [0,1]$ the map
\begin{equation}\label{eq:5.29}
K_* (\phi_\ep ) : K_* (C_r^* (\mc R ,q^\ep )) \cong 
K_* (\mc S (\mc R ,q^\ep )) \to K_* (C_r^* (\mc R ,q )) \cong 
K_* (\mc S (\mc R ,q ))
\end{equation}
is natural. By \eqref{eq:3.68} the same goes for
\begin{equation}\label{eq:5.30}
HP_* (\phi_\ep ) : HP_* (\mc S (\mc R ,q^\ep ))  
\to HP_* (\mc S (\mc R , q))
\end{equation}
Obviously, by Theorem \ref{thm:5.20} these maps are 
isomorphisms for $\ep > 0$, but whether this holds in 
general for $\ep = 0$ is not known. 

Let \inde{$G (C_r^* (\mc R ,q))$} be the Grothendieck group
of the additive category of finite dimensional $C_r^* (\mc R 
,q)$-modules. There is a natural map \index{$G (\phi_0 )$}
\begin{equation}
\begin{aligned}
& G (\phi_0 ) : G \big( C_r^* (\mc R ,q) \big) \to 
G \big( C_r^* (W) \big) \\
& G (\phi_0 ) (\pi ,V) = (\pi \circ \phi_0 ,V)
\end{aligned}
\end{equation}
For $U \subset T / W_0$ we introduce the two-sided ideals
\index{$J_U^s$} \index{$J_U^c$} \index{fiwt@$\phi_{W_0 t_0}$} 
\begin{equation}
\begin{aligned}
& J_U^s = \{ x \in \mc S (\mc R ,q) : \pi(P,W_P r,\delta ,t)(x) 
= 0 \text{  if } t \in T_u^P ,\, W_0 r t \in U \} \\
& J_U^c = \{ x \in C_r^* (\mc R ,q) : \pi(P,W_P r,\delta ,t)(x) 
= 0 \text{  if } t \in T_u^P ,\, W_0 r t \in U \} 
\end{aligned}
\end{equation}
By Theorem \ref{thm:5.20}.2 $\phi_0$ factors through
\begin{equation}
\phi_{W_0 t_0} : \mc S (W) / J^s_{W_0 t_0} \to 
\mc S (\mc R ,q) / J^s_{W_0 t_0 T_{rs}}
\end{equation}
The induced map on $K_0^+$ can be regarded as a morphism
of semigroups
\begin{equation}\label{eq:5.34}
K_0^+ \big( \phi_{W_0 t_0} \big) : 
\mr{Rep}_{W_0 t_0} \big( C_r^* (W) \big) \to
\mr{Rep}_{W_0 t_0 T_{rs}} \big( C_r^* (\mc R ,q) \big)
\end{equation}
The direct sum of these maps, over all $W_0 t_0 \in T_u / W_0$, 
is a homomorphism \index{$K_0^{\mr{Rep}}(\phi_0 )$}
\begin{equation}
K_{\mr{Rep}}(\phi_0 ) : G \big( C_r^* (W) \big) 
\to G \big( C_r^* (\mc R ,q) \big)
\end{equation}
This map is a bit weird, it does not always preserve 
dimensions, and it certainly is not an inverse of $G (\phi_0 )$. 
For example, consider the case $R_0 = A_1 ,\, X$ 
the root lattice and $q(s_0 ) = q(s_1 ) > 1$. Then
\[
\big( K_{\mr{Rep}}(\phi_0 ) \pi \big) \circ \phi_0 = \pi
\]
if $\pi$ admits a central character $t \neq 1$. The only 
irreducible $\mh C W$-representations with central character
$t = 1$ are the trivial and sign representations of $W$.
However, there we see something strange: $K_{\mr{Rep}}(\phi_0 )$ 
sends the trivial representation to the principal series 
$I_1 = \pi (\es, \delta_\es ,1)$, and it sends the 
sign representation to the direct sum of $I_1$ and the 
Steinberg representation of $\mc H (\mc R ,q)$.

\begin{thm}\label{thm:5.25}
The following are equivalent:
\begin{enumerate}
\item $G (\phi_0 ) \otimes \mr{id}_{\mh Q} : 
G \big( C_r^* (\mc R ,q) \big) \otimes \mh Q \to 
G \big( C_r^* (W) \big) \otimes \mh Q$ is a bijection
\item $K_{\mr{Rep}}(\phi_0 ) \otimes \mr{id}_{\mh Q} : G \big( 
C_r^* (W) \big) \otimes \mh Q \to G \big( C_r^* (\mc R ,q) 
\big) \otimes \mh Q$ is a bijection
\item $K_* (\phi_0 ) \otimes \mr{id}_{\mh Q} : K_* \big( 
C_r^* (W) \big) \otimes \mh Q \to K_* \big( C_r^* (\mc R ,q) 
\big) \otimes \mh Q$ is a bijection
\item $HP_* (\phi_0 ) : HP_* (\mc S (W)) \to HP_* (\mc S 
(\mc R ,q))$ is a bijection
\item $HH_0 (\phi_0 ) : \mc S (W) / [\mc S (W),\mc S (W)] \to
\mc S (\mc R ,q) / [\mc S (\mc R ,q),\mc S (\mc R ,q)]$
is an isomorphism of Fr\'echet spaces
\end{enumerate}
\end{thm}
\emph{Proof.}
$1 \Leftrightarrow 2$. By construction it suffices to show this
for $\phi_{W_0 t_0}$, for arbitrary $W_0 t_0 \in T_u /W_0$.
But $\phi_{W_0 t_0}$ is just a homomorphism between finite
dimensional semisimple algebras, so $K_0 (\phi_0 ) \otimes 
\mr{id}_{\mh Q}$ and $G(\phi_0 ) \otimes \mr{id}_{\mh Q}$ are
linear maps between finite dimensional vector spaces. With 
respect to the bases formed by irreducible representations the
matrices of these two maps are each others transpose. In 
particular one of them is bijective if and only if the other is.

$2 \Leftrightarrow 3$. Consider the projection
\begin{equation}\label{eq:5.37}
\begin{aligned}
& \mr{pr} : \Xi_u / \mc W \to T_u / W_0 \\
& \mr{pr}\,(\mc W (P,W_P r ,\delta ,t)) = W_0 r_u t
\end{aligned}
\end{equation}
With this map we make $C_r^* (\mc R ,q))$ into a $C (T_u /W_0) 
$-algebra. By \eqref{eq:5.17} $\phi_0$ is $C (T_u /W_0)$-linear.
Triangulate $T_u /W_0$ such that every subset $T_u^G / W_0$ with
$G \subset W_0$ becomes a subcomplex. In view of Proposition
\ref{prop:2.26}, 2 implies 3.

Contrarily, suppose that $K_{\mr{Rep}}(\phi_0 ) \otimes 
\mr{id}_{\mh Q}$ is not surjective. By definition there is a 
$t_0 \in T_u$ such that $K_0 (\phi_{W_0 t_0} ) \otimes 
\mr{id}_{\mh Q}$ is not surjective. However the canonical map 
\[
K_0 (C_r^* (\mc R ,q)) \to K_0 \big( C_r^* (\mc R ,q) / 
J^c_{W_0 t_0 T_{rs}} \big)
\]
is always surjective. This can be seen as follows. Every 
component of $\Xi_u$ intersects $\mr{pr}^{-1} ( W_0 t_0 )$ in at 
most one $W_0$-orbit. If $[p_1 ] \in K_0 \big( C_r^* (\mc R ,q) 
/ J^c_{W_0 t_0 T_{rs}} \big)$ then $[p_2 ] \in K_0 (C_r^* 
(\mc R ,q))$ maps to $[p_1 ]$ if and only if the value of $p_2$
on $\Xi_u \cap \mr{pr}^{-1} (W_0 t_0 )$ is as prescribed by $p_1$.
It follows from Theorem \ref{thm:3.19} that such a $p_2$ can 
always can be found. 
This shows that $K_* (\phi_0 ) \otimes \mr{id}_{\mh Q}$ and
$K_* (\phi_0 )$ are not surjective.

Now suppose that $K_{\mr{Rep}}(\phi_0 ) \otimes \mr{id}_{\mh Q}$ 
is not injective. Pick $W_0 t_0$ such that $K_0 (\phi_{W_0 t_0})$ 
is not injective, with $| W_{0,t_0} |$ minimal for this property. 
Next pick $V,V' \in \mr{Rep}_{W_0 t_0}(C_r^* (W))$ such that 
$[V] - [V'] \in \ker K_0 \big( \phi_{W_0 t_0} \big)$. Put 
\[
T' := \{t \in T_u : W_{0,t} \not\subset W_{0,t_0} \}
\]
and introduce the ideals 
\begin{align*}
& I_0 := J^c_{T'} \subset C_r^* (W) \\
& I_1 := J^c_{T' T_{rs}} \subset C_r^* (\mc R ,q)
\end{align*}
Note that $\phi_0 (I_0 ) \subset I_1$. Recall the description
\begin{equation}
C_r^* (W) \cong C \big( T_u ;\mr{End} \:\mh C [W_0 ] \big)^{W_0}
\end{equation}
from Lemma \ref{lem:A.3}, in combination with \eqref{eq:2.50}
and Theorem \ref{thm:3.28}. These show that it is possible to 
find $m,n \in \mh N$ and projections $p,p' \in M_n (I_0^+ )$ 
such that $p(t)$ and $p'(t)$ yield the $W_{0,t_0}$-modules
$m V$ and $m V'$, for all $t \in T_u^{W_{0,t_0}} \setminus T'$.
Now we insert $[p] - [p']$ in the commutative diagram
\[
\begin{array}{ccc}
K_0 (I_0 ) & \to & K_0 (C_r^* (W)) \\
\downarrow & & \downarrow \\
K_0 (I_1 ) & \to & K_0 (C_r^* (\mc R ,q)
\end{array}
\]
By assumption 
\[
[\phi_0 (p)] - [\phi_0 (p')] = 0 \in K_0 (I_1 )
\]
On the other hand, $[p]$ and $[p']$ are different on 
$T_u^{W_{0,t_0}}$, so by Theorem \ref{thm:2.13} 
\[
Ch_{W_0} \big( [p] - [p'] \big) \neq 0 \in 
\check H^* \big( \widetilde{T_u} / W_0 ; \mh C \big)
\]
Therefore $K_0 (\phi_0 ) \otimes \mh Q$ and $K_0 (\phi_0 )$ 
are not injective.

$3 \Leftrightarrow 4$ by Theorem \ref{thm:3.41}.

$1 \Leftrightarrow 5$. The localization of 
\[
HH_0 (\mc S (\mc R ,q)) = \mc S (\mc R ,q) / 
[\mc S (\mc R ,q),\mc S (\mc R ,q)] 
\]
at $\mc W \xi \in \Xi_u / \mc W$ is a complex vector space
whose dimension is the number of inequivalent irreducible
constituents of $\pi (\xi )$. This gives a fibration of
$\Xi_u / \mc W$, and by Theorem \ref{thm:3.19} $HH_0 (\mc S
(\mc R ,q))$ can be regarded as the set of global sections
of the sheaf $\mf F_q$ of smooth sections of this fibration. 
The direct image of $\mf F_q$ under \eqref{eq:5.37} is the 
sheaf of smooth sections of a fibration of $T_u / W_0$. The 
fiber at $W_0 t_0$ is a vector space whose dimension is the 
number of irreducibles in $\mr{Rep}_{W_0 t_0 T_{rs}}(C_r^* 
(\mc R ,q))$.

The Fr\'echet space $HH_0 (\mc S (W))$ admits a similar
description it terms of a sheaf $\mf F_1$, but with fibers of
dimension the number of irreducibles in $\mr{Rep}_{W_0 t_0}
(C_r^* (W))$. Now $\phi_0$ induces a morphism 
\[
\mf F (\phi_0 ) : \mf F_1 \to \mr{pr}_* (\mf F_q )
\]
of sheaves over $T_u / W_0$. If $G (\phi_0 ) \otimes 
\mr{id}_{\mh Q}$ is not bijective, then $\mf F_1$ and $\mr{pr}_* 
(\mf F_q )$ have different stalks, so $\mf F (\phi_0 )$ and 
$HH_0 (\phi_0 )$ cannot be isomorphisms. On the other hand, if 
$G (\phi_0 ) \otimes \mr{id}_{\mh Q}$ is bijective, then
\[
\mf F (\phi_0 )(W_0 t_0 ) : \mf F_1 (W_0 t_0 ) \to 
\mr{pr}_* (\mf F_q )(W_0 t_0 )
\]
is a bijection for every $W_0 t_0 \in T_u / W_0$. This implies
that $\mf F (\phi_0 )$ is an isomorphism, see e.g. 
\cite[Section II.1.6]{God}. In particular $HH_0 (\phi_0 ) =
\mf F (\phi_0 ) (T_u / W_0 )$ is an isomorphism of Fr\'echet
spaces. $\qquad \Box$

\begin{conj}\label{conj:5.26}
The equivalent statements of Theorem \ref{thm:5.25} hold 
for every root datum and every positive label function.
\end{conj}

If $q$ is an equal label function then by the Kazhdan-Lusztig
classification (see page \pageref{eq:5.38}) and by Theorem 
\ref{thm:5.4} none of the vector spaces in Theorem 
\ref{thm:5.25} depends on $q$, up to natural isomorphisms. 
It is probable, though not a priori certain, that these 
isomorphisms can be realized with $\phi_0$.

Actually, even for unequal labels strange things happen if
Conjecture \ref{conj:5.26} does not hold. 
Suppose for example that $HH_0 (\phi_0 )$ is not injective. 
In that case there would exist an 
\[
x \in \mc S (W) \setminus [\mc S (W) ,\mc S (W)]
\quad \text{such that} \quad 
\phi_0 (x) \in [\mc S (\mc R ,q),\mc S (\mc R ,q)]
\]
However, since $\phi_\ep$ is an isomorphism $\forall \ep > 0$,
we would have
\begin{equation}
\phi_\ep^{-1} \phi_0 (x) \in [\mc S (\mc R ,q^\ep),\mc S 
(\mc R ,q^\ep)] \quad \forall \ep > 0
\end{equation}
But this is remarkable since 
\[
\ep \to \phi_\ep^{-1} \phi_0 (x) \quad \mr{and} \quad
\ep \to z \cdot_{q^\ep} y - y \cdot_{q^\ep} z 
\]
are both continuous on $[0,1] \;, \forall x,y,z \in \mc S 
(\mc R)$. Maybe one can show that it is outright impossible, 
by a thorough study of conjugacy classes in affine Weyl groups. 
Related results for finite Coxeter groups can be found in
\cite[Sections 3.2 and 8.2]{GePf}.

Or suppose that $G (\phi_0 ) \otimes \mr{id}_{\mh Q}$ is not 
injective. Then there would exist $t_0 \in T_u$ and $(\pi ,V) 
,(\pi' ,V') \in \mr{Rep}_{W_0 t_0 T_{rs}} (C_r^* (\mc R ,q))$ 
such that $\pi \circ \phi_0$ and $\pi' \circ \phi_0$ are equivalent
$C_r^* (W)$-representations. 
The Euler-Poincar\'e pairing from \eqref{eq:3.55} can help in 
this situation. Assume for the moment that $\mc R$ is semisimple. 
By Theorem \ref{thm:5.2} the finite dimensional semisimple 
subalgebras $\mc H (\mc R ,I,q)$ are rigid under $q \to q^\ep$, so 
\begin{equation}
\pi \big|_{\mc H (\mc R ,I,q)} \cong (\pi \circ \phi_0 
)\big|_{\mc H (\mc R ,I,q^0)} \cong (\pi' \circ \phi_0 
)\big|_{\mc H (\mc R ,I,q^0)} \cong \pi' \big|_{\mc H (\mc R ,I,q)}
\end{equation}
Therefore $P_n (V)^\Omega \cong P_n (V')^\Omega$ for all $n$, and
by Corollary \ref{cor:3.29}
\[
\mr{Eul}\, [V] = \mr{Eul}\, [V'] \in K_0 (\mc H )
\]
Together with \eqref{eq:3.58} this shows that $[V] - [V']$ is in 
the radical of the Euler-Poincar\'e pairing. However, we noticed on 
page \pageref{eq:3.57} that the radical of $EP$ is very large, so 
this does certainly not imply that $\pi$ and $\pi'$ are equivalent.
The next theorem is very useful to overcome this problem. It is
analogous to results of Meyer \cite[Theorems 21 and 38]{Mey} 
for Schwartz algebras of reductive $p$-adic groups. 

\begin{thm}\label{thm:5.29}
Let $\mc R$ be any root datum and $V ,V' \in \mr{Rep}(\mc S 
(\mc R ,q))$. Then
\begin{equation}\label{eq:5.47}
\left( \mc S (\mc R ,q) \otimes_{\mc H (\mc R ,q)} P_* (V)^\Omega
\,,\, \mr{id}_{\mc S (\mc R ,q)} \otimes d_* \right)
\end{equation}
is a finitely generated resolution of $V$. This resolution 
consists of projective modules if $\mc R$ is semisimple. Moreover
for such root data there are natural isomorphisms
\[
\mr{Ext}^n_{\mc H (\mc R ,q)} (V,V') \cong
\mr{Ext}^n_{\mc S (\mc R ,q)} (V,V') \qquad n \in \mh N
\]
\end{thm}
\emph{Proof.}
This has been proved recently by Opdam and the author. We can 
construct a suitable contracting homotopy operator of the 
differential complex
\[
\left( P_* (V)^\Omega \,,\,  d_* \right)
\]
Using the temperedness of $V$ we can show that this operator 
extends continuously to the complex \eqref{eq:5.47}.
The details will be published elsewhere. $\qquad \Box$.
\\[3mm]

Suppose that both $V$ and $V'$ are direct sums of discrete series
representations. By Theorem \ref{thm:3.19} a discrete series 
module is projective in Rep$\,(\mc S (\mc R ,q))$, so with 
Theorem \ref{thm:5.29}
\begin{equation}
\begin{array}{lll}
EP \, \big( [V] - [V'] , [V] - [V'] \big) & = &
\sum\limits_{n=0}^\infty (-1)^n \dim \mr{Ext}^n_{\mc S (\mc R ,q)}
\big( [V] - [V'] , [V] - [V'] \big) \\
& = & \dim \mr{Hom}_{\mc S (\mc R ,q)} \big( [V] - [V'] , 
[V] - [V'] \big)
\end{array}
\end{equation}
Clearly this is positive whenever $V$ and $V'$ are inequivalent.
This does not only imply 
$\pi \circ \phi_0 \not\cong \pi' \circ \phi_0$, but also that
\begin{equation}\label{eq:5.48}
G (\phi_0 ) \big( [V] - [V'] \big)
\text{ is not in the radical of } EP.
\end{equation}
Hence this element cannot be written as a sum of virtual 
representations that are induced from proper parabolic subalgebras. 
Because every irreducible tempered representation is a direct 
summand of a representation that is parabolically induced from a
discrete series representation, this is an essential part of 
Conjecture \ref{conj:5.26}. 

In some important cases this actually suffices to prove the 
conjecture. Namely, using the theory of R-groups \cite{DeOp2,Opd4} 
it can be shown that for certain labelled root data all 
representations of the form $\pi (P,\delta ,t)$ with $(P,\delta 
,t) \in \Xi_u$ are irreducible.
Using Theorems \ref{thm:3.15} and \ref{thm:5.29} we can
apply an inductive argument to verify Conjecture \ref{conj:5.26}
in such cases. See Section \ref{sec:6.7} for more details.

Let us have another look at Theorem \ref{thm:5.25}. Clearly the 
first three statements can also be formulated with integral 
coefficients:

\begin{prop}\label{prop:5.24}
The following are equivalent:
\begin{enumerate}
\item $G (\phi_0 ) : G \big( C_r^* (\mc R ,q) \big) \to 
G \big( C_r^* (W) \big)$ is an isomorphism
\item $K_{\mr{Rep}}(\phi_0 ) : G \big( C_r^* (W) \big) \to 
G \big( C_r^* (\mc R ,q) \big)$ is an isomorphism
\item $K_* (\phi_0 ) : K_* \big( C_r^* (W) \big) \to 
K_* \big( C_r^* (\mc R ,q) \big)$ is an isomorphism
\end{enumerate}
\end{prop}
\emph{Proof.}
This is completely analogous to the corresponding part of the
proof of Theorem \ref{thm:5.25}. $\qquad \Box$
\\[3mm]

One of the motivations for considering these maps is that 
$K_* (\phi_0 )$ is natural, in the sense that it can be 
constructed without really using 
$\phi_0$. The idea is that small \inde{perturbations} of 
invertibles or idempotents have no effect on classes in 
$K$-theory. By Theorem \ref{thm:2.18} $K_* (\mc S (W))$ is 
finitely generated. Using Theorem \ref{thm:2.21} we can find 
$k \in \mh N$ and a finite set of idempotents and invertibles 
in $M_k (\mc S (W))$ which generates $K_* (\mc S (W))$. Let 
$(u,q^0 ) \in M_k (\mc S (\mc R ,q^0))$ be such an invertible. 
In view of Proposition \ref{prop:5.9} and the remark on page
\pageref{eq:5.35} there exists an $\ep_u > 0$ such that 
\[
(u,q^\ep ) \in GL_k (\mc S( \mc R ,q^\ep )) 
\quad \forall \ep < \ep_u
\]
To handle idempotents in a similar way we need holomorphic
functional calculus. Define the holomorphic function 
\inde{$f_p$} on $\{ z \in \mh C : \Re (z) \neq 1/2 \}$ by
\[
f_p (z) = \left\{ \begin{array}{ccc}
1 & \mr{if} & \Re (z) > 1/2 \\
0 & \mr{if} & \Re (z) < 1/2
\end{array} \right.
\]
Note that $f_p (x)$ is idempotent whenever it is defined. Let
$(e,q^0 ) \in M_k (\mc S (\mc R ,q^0 ))$ be an idempotent. It
is clear from Proposition \ref{prop:5.6} that $\exists \ep_e > 0$
such that 
\[
(e,q^\ep ) - 1/2 - a i \in GL_k (\mc S (\mc R ,q^\ep ))
\quad \forall \ep < \ep_e ,\, \forall a \in \mh R
\]
In fact, by direct calculation one can show that this holds for
all $\ep$ such that 
\[
\norm{\lambda (e,q^\ep ) - \lambda (e,q^0 )}_{B (\mf H (\mc R))} 
< \frac{1}{\ds 2 + 4 \norm{\lambda (e,q^0 )}_{B (\mf H (\mc R))} } 
\]
For such $\ep$ the idempotent $f_p (e,q^\ep )$ is well-defined.

\begin{lem}\label{lem:5.27}
The following equalities of $K$-theory classes hold for $u$ and 
$e$ as above.
\[
\begin{array}{rcccc@{\qquad}c}
[(u,q^\ep )] & = & K_1 (\phi_\ep^{-1} \phi_0 ) [(u,q^0 )] & \in
& K_1 (\mc S (\mc R ,q^\ep )) & \forall \ep \in (0,\ep_u ) \\
{[f_p (e,q^\ep )]} & = & K_0 (\phi_\ep^{-1} \phi_0 ) [(e,q^0 )] & 
\in & K_0 (\mc S (\mc R ,q^\ep )) & \forall \ep \in (0,\ep_e )
\end{array}
\]
\end{lem}
\emph{Proof.}
For any $x \in M_k (\mh C ) \otimes \mc S (\mc R )$ the map
\[
[0,1] \to M_k (\mh C ) \otimes \mc S (\mc R ) : 
\ep \to \phi_\ep^{-1} \phi_0 (x)
\]
is continuous. Hence $(u,q^{\ep_1})$ and $\phi_{\ep_1}^{-1} \phi_0 
(u,q^0 )$ are homotopic in $GL_k (\mc S (\mc R ,q^{\ep_1} ))$, for
some small $\ep_1 > 0$. Clearly this implies that $\phi_\ep^{-1}
\phi_{\ep_1} (u,q^{\ep_1} )$ and $\phi_\ep^{-1} \phi_0 (u,q^0 )$
are homotopic $\forall \ep \in (0,\ep_u )$. But there also is a
path from $\phi_\ep^{-1} \phi_{\ep_1} (u,q^{\ep_1} )$ to 
$(u,q^\ep )$ along elements of the form $\phi_\ep^{-1} 
\phi_{\ep_2} (u,q^{\ep_2})$.

Similarly, by Corollary \ref{cor:5.22} 
\[
[0,\ep_e ) \to  M_k (\mh C ) \otimes \mc S (\mc R ) : 
\ep \to f_p (e,q^\ep )
\]
is continuous. According to \cite[Proposition 4.3.2]{Bla} there
is a small $\ep_3 > 0$ such that $f_p (e,q^{\ep_3})$ and 
$\phi_{\ep_3}^{-1} \phi_0 (e,q^0 )$ are homotopic in Idem$(M_k
(\mc S (\mc R ,q^{\ep_3}))$. But then, as above for $u ,\, 
\phi_\ep^{-1} \phi_0 (e,q^0 )$ and $f_p (e,q^\ep )$ are 
homotopic via $\phi_\ep^{-1} \phi_{\ep_3} (f_p (e,q^{\ep_3})).
\qquad \Box$ \\[3mm]

So we have a family of pre-$C^*$-algebras $\mc S (\mc R ,q))$ 
which are independent of $q$ as Fr\'echet spaces, and whose 
multiplication depends continuously on $q$. Moreover, replacing 
$q$ by $q^\ep$ with $\ep > 0$ sufficiently small, we may assume 
that the natural group homomorphism $K_* (\phi_0 )$ can be
constructed without using $\phi_0$. Therefore it is not 
unreasonable to suspect the following. 

\begin{conj}\label{conj:5.28}
For any root datum $\mc R$ and positive label function $q$
the map
\[
K_* (\phi_0 ) : K_* (\mc S (W)) \to K_* (\mc S (\mc R ,q))
\]
is an isomorphism.
\end{conj}

As mentioned, this conjecture stems from Baum, Connes and Higson
\cite{BCH}, at least in the equal label case. Independently, Opdam
\cite[p. 533]{Opd3} stated it for unequal labels. In Chapter 6 we
will verify Conjecture \ref{conj:5.28} for some classical root data.

Consider the root datum $\mc R \times \mh Z$ with the unique 
label function that extends $q$. By \eqref{eq:3.47}
\begin{equation}
\mc S (\mc R \times \mh Z ,q) \cong 
\mc S (\mc R ,q) \hot \mc S (\mh Z )
\cong C^\infty \big( S^1 ; \mc S (\mc R ,q) \big)
\end{equation}
In Lemma \ref{lem:2.25} we constructed natural isomorphisms
\begin{equation}\label{eq:5.36}
K_0 (\mc S (\mc R \times \mh Z ,q)) \mosi K_* (\mc S (\mc R ,q))
\isom K_1 (\mc S (\mc R \times \mh Z ,q))
\end{equation} 
Therefore it suffices to prove Conjecture
\ref{conj:5.26} either for $K_0$ and every $(\mc R ,q)$, or for 
$K_1$ and every $(\mc R ,q)$. Probably the $K_1$-case is easier, for 
two reasons. Firstly, invertibles are more flexible than idempotents. 
If we perturb them a little they remain invertible, so we can do 
without holomorphic functional calculus. Secondly, we can find a 
bound, uniform in $q$, on the size of the matrices that we need to 
represent all $K_1$-classes. In fact, by Proposition \ref{prop:2.23} 
we can bound the topological stable rank by
\begin{equation}
tsr (C_r^* (\mc R ,q)) \leq |W_0 |^2 \big( 1 + \lfloor
\dim T_u / 2 \rfloor \big)
\end{equation}
Now Theorems \ref{thm:2.19} and \ref{thm:2.22} show that 
\begin{equation}
K_1 (\mc S (\mc R ,q)) \cong \pi_0 \big( GL_n (\mc S (\mc R ,q)) 
\big) \quad \forall n \geq |W_0 |^2 (1 + \mr{rk}(X) / 2)
\end{equation}
One way to attack Conjecture \ref{conj:5.28} goes approximately as
follows. Pick a finite set of generators of $K_1 (\mc S (\mc R ,q))$.
Find suitable representants $u_i \in GL_n (\mc S (\mc R ,q))$, i.e.
$u_i$ should lie in 
\[
M_n (\mh C) \otimes \mr{span}\{ N_w : \mc N (w) \leq M \}
\]
with $M$ "small" and sp$(u_i)$ should lie in a "small" neighborhood
of the unit circle in $\mh C$. If $q$ is close to $q^0$ one may 
hope that to every $u_i$ one can associate in an unambiguous way a 
unique homotopy class in $GL_n (\mc S (W))$. This class should be 
constructed by applying the isomorphisms $\phi_\ep^{-1}$ and by 
perturbing $u_i$ a little. In particular it should contain elements
of $M_n (\mh C) \otimes \mc S (\mc R )$ that are homotopic to $u_i$ 
in $GL_n (\mc S (\mc R ,q))$. 

An analogue of Conjecture \ref{conj:5.28} does hold for 
\inde{noncommutative tori}. Let \inde{$\mh T^n$} $= \mh R^n / 
\mh Z^n \cong (S^1 )^n$ be the standard compact $n$-dimensional torus. 
In this setting the underlying group $W$ becomes $\mh Z^n$, and $q$ 
is replaced by a skew-symmetric bilinear form $\theta$ on $\mh Z^n$. 
By taking iterated crossed products with $\mh Z$, one constructs a 
$C^*$-algebra which is a deformation of $C (\mh T^n )$ and is commonly 
denoted by $A_\theta$. It has a holomorphically closed dense subalgebra 
$\mc S (\mh Z^n ,\theta)$ which, as a Fr\'echet space, is naturally 
isomorphic to $\mc S (\mh Z^n )$. By deforming $\theta$ Elliott 
\cite[Theorem 2.2]{Ell} proved that there exists a natural group 
isomorphism
\begin{equation}
K_* (A_\theta ) \isom \bigwedge \mh Z^n
\end{equation}
For $\theta = 0$ this can be interpreted geometrically as the classical
Chern character
\[
Ch : K^* (\mh T^n ) \isom H^* (\mh T^n ; \mh Z)
\]
However, there are also quite big differences between noncommutative 
tori and affine Hecke algebras. Namely, the structure of $A_\theta$ 
is very different for $\theta$ rational or irrational, and an 
essential ingredient in Elliott's proof is the Pimsner-Voiculescu 
exact sequence, which is not available for crossed products with
non-cyclic groups.

%% file: chapter6.tex
\chapter{Examples and calculations}

The final chapter of this book is completely different from the
others. We barely state or prove theorems here, we mostly make 
calculations.

In Chapters 3 and 5 we studied affine Hecke algebras in a very
abstract way, almost without mentioning any examples. However
the results in those chapters could hardly have been obtained
without first checking, by hand, what happens in some simple
exemplary cases. Basically we have two goals: we want have to 
examples of all the objects we introduced in Chapter 3, and we 
want to verify the conjectures made in Section \ref{sec:5.4} in 
some cases.

So we must devise a strategy to calculate the $K$-theory of the
$C^*$-completion $C_r^* (\mc R ,q)$ of an affine Hecke algebra
with root datum $\mc R$ and label function $q$, and to find
the homomorphisms (hopefully isomorphisms)
\begin{equation}\label{eq:6.26}
K_* (\phi_0 ) : K_* \big( C_r^* (\mc R ,q^0 ) \big) \to
K_* \big( C_r^* (\mc R ,q) \big)
\end{equation}
For $q^0$ we can use \eqref{eq:2.50} and \eqref{eq:3.51}, 
which say that 
\begin{equation}\label{eq:6.37}
K_* \big( C_r^* (W) \big) \otimes \mh C \cong K_* (C (T_u ) 
\rtimes W_0 ) \otimes \mh C \cong \check H^* \big( 
\widetilde{T_u} ; \mh C \big)^{W_0} \cong \check H^* \big( 
\widetilde{T_u} \big/ W_0; \mh C \big)
\end{equation}
Moreover, if $\check H^* \big( \widetilde{T_u} \big/ W_0; 
\mh Z \big)$ is torsion free, then by Theorem \ref{thm:2.15}
\begin{equation}\label{eq:6.38}
K_* \big( C_r^* (W) \big) = K_* \Big( C \big( T_u ; 
\mr{End}\:\mh C [W_0 ] \big)^{W_0} \Big) \cong
\check H^* \big( \widetilde{T_u} \big/ W_0; \mh Z \big)
\end{equation}
In general our procedure will involve the following steps.
\begin{enumerate}
\item Explicitly write down the root datum and the associated
Weyl groups.
\item Determine the residual cosets and distinguish the 
different "genericity classes" of label functions.

For every $q$ there is, as noticed in \eqref{eq:5.46}, a
canonical decomposition
\begin{equation}\label{eq:6.27}
C_r^* (\mc R ,q) = \bigoplus_P C_r^* (\mc R ,q)_P
\end{equation}
where $P$ runs over certain sets of simple roots.
\item List a good set of $P$'s.

For every chosen $P$ we do the following:
\item Determine the root datum $\mc R_P$ 
\item Determine the discrete series of $\mc H (\mc R_P ,q_P )$,
and all the relevant intertwining operators.

This is the only step where we can still encounter theoretical 
difficulties. The problem is that in general it is not known how 
many inequivalent discrete series representations there are. 
Fortunately we can decide this in the cases we consider.

\item Describe $C_r^* (\mc R ,q)_P$ and its primitive ideal 
spectrum.
\item Calculate $K_* \big( C_r^* (\mc R ,q) \big)_P$.

Our main tools for this are excision and Proposition 
\ref{prop:2.26}. In principle these will always lead to the 
answer, but it is tedious work that becomes unpractical in higher 
dimensions. On the other hand things become easier if we forget 
about torsion elements, for then we can use sheaf cohomology. 
For this reason will sometimes be satisfied with 
$K_* \big( C_r^* (\mc R ,q)_P \big) \otimes_{\mh Z} \mh C$.

\item Find generating idempotents and invertibles, as explicit
as possible.
\item Compare $K_* \big( C_r^* (\mc R ,q) \big)$ and
$K_* \big( C_r^* (\mc R ,q^0) \big)$.
\item Determine $K_* (\phi_0 )$ in terms of the given generators.
\\[1mm]
\end{enumerate}

The results of these calculations are
\begin{itemize}
\item There are no examples for which \eqref{eq:6.26} is known
to be no isomorphism.
\item \eqref{eq:6.26} is an isomorphism for some root data of
low rank
\item \eqref{eq:6.26} is an isomorphism the root data
$\mc R (GL_n )$ and $\mc R (A_{n-1})^\vee$. 
\end{itemize}

We will also see that these $K$-groups tend to be torsion 
free. This is due to the fact that the primitive ideal spaces
of affine Hecke algebras look like quotients of tori by reflection 
groups, or direct products and unions of those. The root data
we study all have free abelian $K$-groups, but whether this 
holds in general is hard to say.

\section{$A_1$}
\label{sec:6.1}

Like in the theory of semisimple Lie algebras, the rank one
root system $A_1$ plays an important role in the realm of Hecke
algebras. Every Iwahori-Hecke algebra is in a sense built from
such rank one Hecke algebras. There are two semisimple root
data with $R_0 = A_1$, corresponding to the root lattice and
the weight lattice. We will study the associated affine Hecke
algebras in detail, and show that Conjecture \ref{conj:5.28}
holds for these root data.

The notation $\mc R (A_1 )$ will be reserved for $X$ the root
lattice. We start with the other case. Thus we consider the root 
datum \inde{$\mc R (A_1 )^\vee$} $= (X,Y,R_0 ,R_0^\vee ,F_0 )$ with
\[
\begin{aligned}
& X = \mh Z \qquad Q = 2 \mh Z \qquad X^+ = \mh Z_{\geq 0} \\
& Y = Q^\vee = \mh Z \\
& T = \mh C^\times \\
& R_0 = R_1 = \{ \pm \alpha \} = \{ \pm 2 \} \\
& R_0^\vee = R_1^\vee = \big\{ \pm \alpha^\vee \big\} = 
  \{\pm 1 \} \\
& F_0 = \{ \alpha \} \qquad W_0 = \{ e , s_\alpha \} \\
& s_1 = s_\alpha : x \to -x \qquad 
  s_0 = t_2 s_1 = t_1 s_1 t_{-1} : x \to 2 - x \\
& S_{\mr{aff}} = \{ s_0 ,s_1 \} \qquad W \neq W_{\mr{aff}} = 
  \langle s_0 ,s_1 | s_0^2 = s_1^2 = e \rangle \\
& \Omega = \{e ,\omega \} = \{e ,t_1 s_1 \} 
\end{aligned}
\]
For any label function $q$ we have $q (s_0 ) = q(s_1 ) =
q_{\alpha^\vee}$, and we denote this value simply by $q$. Then
\[
c_\alpha = (1 - q^{-1} \theta_{-2})(1 - \theta_{-2})^{-1}
\]
so generically the residual points are
\begin{equation}\label{eq:6.1}
q^{1/2} ,\, q^{-1/2} ,\, -q^{1/2} ,\,-q^{-1/2}
\end{equation}
Hence there are only two essentially different cases,
depending on whether $q$ equals 1 or not.
\\[2mm]

\begin{itemize}
\item {\large \textbf{group case} $\mb{q = 1}$}
\end{itemize}
From Theorem \ref{thm:3.28} we know that every irreducible
representation is a direct summand of a unitary principal
series $I_t = \pi (\es ,\delta_\es ,t)$. The
underlying vector space of $I_t$ is 
\[
\mh C [W_0 ] = \mh C T_e + \mh C T_{s_1}
\]
and the intertwiner $\pi (s_1 ,\es ,\delta_\es
,t) : I_t \to I_{t^{-1}}$ is simply right multiplication
by $T_{s_1}$. So with respect to the orthonormal basis 
\begin{equation}\label{eq:6.2}
\big\{ 2^{-1/2} (T_e + T_{s_1}) ,2^{-1/2}(T_e - T_{s_1}) \big\}
\end{equation}
we have
\begin{equation}\label{eq:6.3}
\begin{array}{lll}
\mc S (W) & \cong & \left\{ f \in C^\infty \big( S^1 ;M_2 (\mh C )
\big) : f(t^{-1}) = \begin{pmatrix} 1 & 0 \\ 0 & -1 \end{pmatrix} 
f(t) \begin{pmatrix} 1 & 0 \\ 0 & -1 \end{pmatrix} \right\} \\
C_r^* (W) & \cong & \{ f \in C \big( [-1,1];M_2 (\mh C ) \big) :
f(\pm 1) \text{ is diagonal} \}
\end{array}
\end{equation}
The spectrum of these algebras is the non-Hausdorff space

\raisebox{3mm}{Prim$(C_r^* (W)) \qquad \cong$}
\begin{picture}(5,1)
\linethickness{1mm}
\put(1,0.4){\line(1,0){3.5}}
\put(1,0.3){\circle{0.14}}
\put(1,0.5){\circle{0.14}}
\put(4.5,0.5){\circle{0.14}}
\put(4.5,0.3){\circle{0.14}}
\end{picture}

To calculate the $K$-theory we use the extension 
\[
0 \to C_0 ([-1,1], \{0,1\} ;M_2 (\mh C)) \to C_r^* (W) \to \mh C^4 \to 0
\]
defined by evaluation at 1 and $-1$. The associated exact hexagon
is 
\[
\hexagon{0}{K_0 (C_r^* (W))}{\mh Z^4}{0}{K_1 (C_r^* (W))}{\mh Z}
\]
By comparing this with the standard extension
\[
0 \to C_0 ([-1,1],\{0,1\} ) \to C (S^1 ) \to \mh C \to 0
\]
we see that the vertical maps are surjective, so
\begin{equation}\label{eq:6.4}
\begin{aligned}
& K_0 (C_r^* (W)) \cong \mh Z^3 \\
& K_1 (C_r^* (W)) = 0
\end{aligned}
\end{equation}
We can even find explicit generating projections, namely
\begin{equation}\label{eq:6.5}
\begin{aligned}
& p_a = (T_e + T_{s_1}) / 2\\
& p_b = (T_e - T_{s_1}) / 2\\
& p_c = T_e / 2 - ( (\theta_1 + \theta_{-1}) T_{s_1} + i
  (\theta_{-1} - \theta_1 ) T_e ) / 4\\
& p_d = T_e / 2 + ( (\theta_1 + \theta_{-1}) T_{s_1} + i
  (\theta_{-1} - \theta_1 ) T_e ) / 4
\end{aligned}
\end{equation}
The only relation between these generators is
\[
[p_a ] + [p_b ] = [p_c ] + [p_d ] = [1] \in K_0 (C_r^* (W))
\]
\\[1mm]

\begin{itemize}
\item {\large \textbf{generic, equal label case} $\mb{q \neq 1}$}
\item $P = \es$
\end{itemize}
\[
\begin{aligned}
& R_P = \es \qquad R_P^\vee = \es \\
& X^P = X \qquad X_P = 0 \qquad Y^P = Y \qquad Y_P = 0 \\
& T^P = T \qquad T_P = \{ 1 \} \qquad K_P = \{ 1 \} \\
& W^P = W(P,P) = \mc W_{PP} = W_0  \qquad W_P = \{ e \} 
\end{aligned}
\]
Again there is a single intertwiner $\pi (s_1 ,\es 
,\delta_\es ,t) : I_t \to I_{t^{-1}}$ and
\[
\begin{aligned}
& \imath^o_{s_1} = (T_{s_1} (1- \theta_2 ) + (q-1) \theta_2) 
(q - \theta_2 )^{-1} \\
& \pi (s_1 ,\es ,\delta_\es ,1) = \pi (s_1 
  ,\es ,\delta_\es ,-1) = 1 \\
& C_r^* (\mc R ,q)_P \cong C \big( [-1,1] ;M_2 (\mh C ) \big) \\ 
& \mr{Prim} \big( C_r^* (\mc R ,q)_P \big) 
\cong S^1 / W_0 \cong [-1,1] 
\end{aligned}
\]

\begin{itemize}
\item $P = \{ \alpha \}$
\end{itemize}
\[
\begin{aligned}
& R_P = R_0 \qquad R_P^\vee = R_0^\vee \\
& W^P = W(P,P) = \mc W_{PP} = \{ e \}  \qquad W_P = W_0 \\
& X^P = 0 \qquad X_P = X \qquad Y^P = 0 \qquad Y_P = Y \\
& T^P = \{ 1 \} \qquad T_P = T \qquad K_P = \{ 1 \}
\end{aligned}
\]
From Proposition \ref{prop:3.30}.2 we know that there is exactly 
one discrete series representation for every orbit of residual 
points. So the spectrum of $C_r^* (\mc R ,q)$ contains two 
isolated points, which we call $\delta_1$ and $\delta_{-1}$, 
by the sign of their central character. 
\[
C_r^* (\mc R ,q)_P \cong \mh C^2
\]
By brute calculation one finds the associated projectors
\begin{equation}\label{eq:6.6}
\begin{array}{lllll}
p_1 & = & \sum\limits_{w \in W_{\mr{aff}}} (-q)^{\ell (w)} T_w 
(T_e - T_\omega ) \Big( \sum\limits_{w \in W} q(w)^{-1} 
\Big)^{-1} & \mr{if} & q > 1 \\
p_1 & = & \sum\limits_{w \in W} T_w \Big( \sum\limits_{w \in W} 
q(w)^{-1} \Big)^{-1} & \mr{if} & q < 1 \\
p_{-1} & = & \sum\limits_{w \in W} (-q)^{\ell (w)} T_w \Big( 
\sum\limits_{w \in W} q(w) \Big)^{-1} & \mr{if} & q > 1 \\
p_{-1} & = & \sum\limits_{w \in W_{\mr{aff}}} T_w (T_e - 
T_\omega ) \Big( \sum\limits_{w \in W} q(w) \Big)^{-1} & 
\mr{if} & q < 1
\end{array}
\end{equation}
Combining the results for $P = \es$ and $P = \{ \alpha \}$,
the spectrum of $C_r^* (\mc R ,q)$ becomes the Hausdorff space

\raisebox{3mm}{Prim$(C_r^* (W)) \qquad \cong$}
\begin{picture}(6,1)
\linethickness{1mm}
\put(1.5,0.4){\line(1,0){3.5}}
\put(1,0.4){\circle*{0.15}}
\put(5.5,0.4){\circle*{0.15}}
\end{picture}

This implies 
\begin{equation}\label{eq:6.7}
\begin{aligned}
& K_0 (C_r^* (\mc R ,q)) \cong \mh Z^3 \\
& K_1 (C_r^* (\mc R ,q)) = 0
\end{aligned}
\end{equation}
Let $p_0$ be any rank one projector in $C([-1,1];M_2 (\mh C))$.
Then (the classes of) $p_{-1} ,\, p_0$ and $p_1$ generate
$K_0 (C_r^* (\mc R ,q))$, and we have
\[
[p_{-1}] + [p_1 ] + 2 [p_0 ] = [1] \in K_0 (C_r^* (\mc R ,q))
\]
Now we can compare the cases $q=1$ and $q \neq 1$. Looking 
carefully at the behaviour near $t = 1$ and $t = -1$, we find that
$K_0 (\phi_0 )$ is always an isomorphism. We describe this map
completely with the following table.

\begin{equation}\label{eq:6.8}
\begin{array}{|l|l|l|}
\hline q < 1 & q = 1 & q > 1 \\ 
\hline p_0 + p_1 + p_{-1} & p_a & p_0 \\
p_0 & p_b & p_0 + p_1 + p_{-1} \\
p_0 + p_{-1} & p_c & p_0 + p_1 \\
p_0 + p_1 & p_d & p_0 + p_{-1} \\ 
\hline
\end{array}
\end{equation}
\\[2mm]
Let us move on to the type $A_1$ Hecke algebras with $X$ the
root lattice. This means that we work with the root datum
\inde{$\mc R (A_1 )$}:
\[
\begin{aligned}
& X = Q = \mh Z \qquad X^+ = \mh Z_{\geq 0} \\
& Y = \mh Z \qquad Q^\vee = 2 \mh Z \\
& T = \mh C^\times \\
& R_0 = \{ \pm \alpha \} = \{ \pm 1 \} \qquad R_1 = \{ \pm 2 \} \\
& R_0^\vee = \{ \pm \alpha^\vee \} = \{ \pm 2 \} \qquad 
R_1^\vee = \{ \pm 1 \}  \\
& F_0 = \{ \alpha \} \qquad W_0 = \{ e , s_\alpha \} \\
& s_1 = s_\alpha : x \to -x \qquad s_0 = t_1 s_1 : x \to 1-x \\
& S_{\mr{aff}} = \{ s_0 ,s_1 \} \qquad W = W_{\mr{aff}} =
  \langle s_0 ,s_1 | s_0^2 = s_1^2 = e \rangle
\end{aligned}
\]
Now $s_0$ and $s_1$ are no longer conjugate, so $q_0 = q (s_0 )$
and $q_1 = q(s_1 )$ may be different. By definition
\begin{align*}
& q_{\alpha^\vee} = q_0 \qquad q_{\alpha^\vee / 2} = q_1 q_0^{-1} \\
& c_\alpha = (1 + q_1^{-1/2} q_0^{1/2} \theta_{-1}) (1 - q_1^{-1/2}
q_0^{-1/2} \theta_{-1}) (1 - \theta_{-2})^{-1}
\end{align*}
Generically the residual points are
\begin{equation}\label{eq:6.9}
q_0^{1/2} q_1^{1/2} \qquad q_0^{-1/2} q_1^{-1/2} \qquad
-q_0^{1/2} q_1^{-1/2} \qquad -q_0^{-1/2} q_1^{1/2}
\end{equation}
From this we see that there are four cases to study: the group case
$q_0 = q_1 = 1$, the equal label case $q_0 = q_1 \neq 1$, the
"inverse label" case $q_0 = q_1^{-1} \neq 1$ and the generic case.

\begin{itemize}
\item {\large \textbf{group case} $\mb{q_0 = q_1 = 1}$}
\end{itemize}
This is the same as the group case for $X$ equal to the weight
lattice of $A_1$.

\begin{itemize}
\item {\large \textbf{equal label case} $\mb{q_0 = q_1 = q \neq 1}$}
\item $P = \es$
\end{itemize}
\[
\begin{aligned}
& R_P = \es \qquad R_P^\vee = \es \\
& X^P = X \qquad X_P = 0 \qquad Y^P = Y \qquad Y_P = 0 \\
& T^P = T \qquad T_P = \{ 1 \} \qquad K_P = \{ 1 \} \\
& W^P = W(P,P) = \mc W_{PP} = W_0  \qquad W_P = \{ e \} \\
& \imath^o_{s_1} = (T_{s_1} (1- \theta_1 ) + (q-1) \theta_1) 
(q - \theta_1 )^{-1} 
\end{aligned}
\]
With respect to the orthonormal basis
\begin{equation}\label{eq:6.10}
\big\{ (T_e + T_{s_1})(1+q)^{-1/2} ,(q T_e - T_{s_1})(q+q^2)^{-1/2} \big\}
\end{equation}
of $\mc H (W_0 ,q)$ we have
\[
\begin{aligned}
& \pi (s_1 ,\es ,\delta_\es ,1) = \begin{pmatrix}
1 & 0 \\ 0 & 1 \end{pmatrix} \qquad
\pi (s_1 ,\es ,\delta_\es ,-1) = \begin{pmatrix}
1 & 0 \\ 0 & -1 \end{pmatrix} \\
& C_r^* (\mc R ,q)_P \cong \big\{ f \in C \big( [-1,1] ;M_2 (\mh C) 
\big) : f(-1) \text{ is diagonal} \big\} \\
& \raisebox{3mm}{Prim$\big( C_r^* (\mc R ,q)_P \big) \qquad \cong$} 
\begin{picture}(5,1)
\linethickness{1mm}
\put(1,0.4){\line(1,0){3.5}}
\put(1,0.3){\circle{0.14}}
\put(1,0.5){\circle{0.14}}
\end{picture}
\end{aligned}
\]
\begin{itemize}
\item $P = \{ \alpha \}$
\end{itemize}
\[
\begin{aligned}
& R_P = R_0 \qquad R_P^\vee = R_0^\vee \\
& X^P = 0 \qquad X_P = X \qquad Y^P = 0 \qquad Y_P = Y \\
& T^P = \{ 1 \} \qquad T_P = T \qquad K_P = \{ 1 \} \\
& W^P = W(P,P) = \mc W_{PP} = \{ e \} \qquad W_P = W_0 
\end{aligned}
\]
Obviously $-q_0^{1/2} q_1^{-1/2} = -q_0^{-1/2} q_1^{1/2} = -1$,
so these points are not residual for this particular label function.
On the other hand, by Proposition \ref{prop:3.30}.2 there is a 
unique discrete series representation $\delta_1$ with central 
character $q^{\pm 1}$. It has dimension 1, and the corresponding 
projection is
\begin{equation}\label{eq:6.11}
\begin{array}{lllll}
p_1 & = & \sum\limits_{w \in W} (-q)^{\ell (w)} T_w \Big( 
\sum\limits_{w \in W} q(w)^{-1} \Big)^{-1} & \mr{if} & q > 1 \\
p_1 & = & \sum\limits_{w \in W} T_w \Big( \sum\limits_{w \in W} 
q(w) \Big)^{-1} & \mr{if} & q < 1
\end{array}
\end{equation}

We conclude that
\begin{align*}
& C_r^* (\mc R ,q) \cong \big\{ f \in C \big( [-1,1] ;M_2 (\mh C) 
\big) : f(-1) \text{ is diagonal} \big\} \oplus \mh C \\
& \raisebox{3mm}{Prim$\big( C_r^* (\mc R ,q) \big) \qquad \cong$}
\begin{picture}(6,1)
\linethickness{1mm}
\put(1,0.4){\line(1,0){3.5}}
\put(1,0.3){\circle{0.14}}
\put(1,0.5){\circle{0.14}}
\put(5.5,0.4){\circle*{0.15}}
\end{picture}
\end{align*}
Evaluating at $-1$ and at $q^{\pm 1}$ yields an extension
\[
0 \to C_0 \big( [-1,1] ,\{-1\};M_2 (\mh C) \big) \to
C_r^* (\mc R ,q) \to \mh C^3 \to 0
\]
whose exact hexagon in $K$-theory is
\[
\hexagon{0}{K_0 (C_r^* (\mc R ,q))}{\mh Z^3}{0}{K_1 
(C_r^* (\mc R ,q))}{0}
\]
This shows that
\begin{equation}\label{eq:6.12}
\begin{aligned}
& K_0 (C_r^* (\mc R ,q)) \cong \mh Z^3 \\
& K_1 (C_r^* (\mc R ,q)) = 0
\end{aligned}
\end{equation}
Generating projections are
\[
p_1 = \left( \begin{pmatrix} 0 & 0 \\ 0 & 0 \end{pmatrix} 
, 1 \right) \qquad p_a = \left( \begin{pmatrix} 1 & 0 \\ 
0 & 0 \end{pmatrix} ,0 \right) \qquad p_b = \left( 
\begin{pmatrix} 0 & 0 \\ 0 & 1 \end{pmatrix} ,0 \right)
\]
Explicitly this works out to
\begin{equation}\label{eq:6.13}
\begin{array}{lllll}
p_a & = & (T_e + T_{s_1})(1+q)^{-1} & \mr{if} & q > 1 \\
p_a & = & (T_e + T_{s_1})(1+q)^{-1} - p_1 & \mr{if} & q < 1 \\
p_b & = & (q T_e - T_{s_1})(1+q)^{-1} - p_1 & \mr{if} & q > 1 \\
p_a & = & (T_e + T_{s_1})(1+q)^{-1} & \mr{if} & q < 1 \\
\end{array}
\end{equation}

\begin{itemize}
\item {\large \textbf{inverse label case} 
  $\mb{q = q_1 = q_0^{-1} \neq 1}$}
\item $P = \es$
\end{itemize}
\[
\begin{aligned}
& R_P = \es \qquad R_P^\vee = \es \\
& X^P = X \qquad X_P = 0 \qquad Y^P = Y \qquad Y_P = 0 \\
& T^P = T \qquad T_P = \{ 1 \} \qquad K_P = \{ 1 \} \\
& W^P = W(P,P) = \mc W_{PP} = W_0  \qquad W_P = \{ e \} \\
& \imath^o_{s_1} = (T_{s_1} (1+ \theta_1 ) + (1-q) \theta_1) 
(q + \theta_1 )^{-1} 
\end{aligned}
\]
With respect to the basis \eqref{eq:6.10} we have
\[
\begin{aligned}
& \pi (s_1 ,\es ,\delta_\es ,1) = \begin{pmatrix}
1 & 0 \\ 0 & -1 \end{pmatrix} \qquad
\pi (s_1 ,\es ,\delta_\es ,-1) = \begin{pmatrix}
1 & 0 \\ 0 & 1 \end{pmatrix} \\
& C_r^* (\mc R ,q)_P \cong \big\{ f \in C \big( [-1,1] ;M_2 (\mh C) 
\big) : f(1) \text{ is diagonal} \big\} \\
& \raisebox{3mm}{Prim$\big( C_r^* (\mc R ,q)_P \big) \qquad \cong$} 
\begin{picture}(5,1)
\linethickness{1mm}
\put(1,0.4){\line(1,0){3.5}}
\put(4.5,0.3){\circle{0.14}}
\put(4.5,0.5){\circle{0.14}}
\end{picture}
\end{aligned}
\]
\begin{itemize}
\item $P = \{ \alpha \}$
\end{itemize}
\[
\begin{aligned}
& R_P = R_0 \qquad R_P^\vee = R_0^\vee \\
& X^P = 0 \qquad X_P = X \qquad Y^P = 0 \qquad Y_P = Y \\
& T^P = \{ 1 \} \qquad T_P = T \qquad K_P = \{ 1 \} \\
& W^P = W(P,P) = \mc W_{PP} = \{ e \} \qquad W_P = W_0 
\end{aligned}
\]
Now we have $q_0^{1/2} q_1^{1/2} = q_0^{-1/2} q_1^{-1/2} = 1$, so
these points are not residual. There is a unique discrete series
representation $\delta_{-1}$ with central character $-q^{\pm 1}$.
Its projector is already a little more difficult to describe. For
$w \in W$ write $\ell (w) = \ell_0 (w) + \ell_1 (w)$, where 
$\ell_1$ counts the number of factors $s_i$ in an reduced 
expression for $w$. Notice that this is well-defined only because
there are no relations between $s_0$ and $s_1$ in $W_{\mr{aff}}$.
\begin{equation}
\begin{array}{ccccc}
p_{-1} & = & \sum\limits_{w \in W} (-q_1 )^{\ell_1 (w)/2} T_w \Big( 
\sum\limits_{w \in W} q^{-\ell (w)} \Big)^{-1} & \mr{if} & q > 1 \\
p_{-1} & = & \sum\limits_{w \in W} (-q_0 )^{\ell_0 (w)/2} T_w \Big( 
\sum\limits_{w \in W} q^{\ell (w)} \Big)^{-1} & \mr{if} & q < 1
\end{array}
\end{equation}

Summarizing, we have
\begin{align*}
& C_r^* (\mc R ,q) \cong \mh C \oplus \big\{ f \in C \big( [-1,1] 
;M_2 (\mh C) \big) : f(-1) \text{ is diagonal} \big\} \\
& \raisebox{3mm}{Prim$\big( C_r^* (\mc R ,q) \big) \qquad \cong$}
\begin{picture}(6,1)
\linethickness{1mm}
\put(2,0.4){\line(1,0){3.5}}
\put(5.5,0.3){\circle{0.14}}
\put(5.5,0.5){\circle{0.14}}
\put(1,0.4){\circle*{0.15}}
\end{picture}
\end{align*}
Like in the equal case this leads to
\begin{equation}\label{eq:6.14}
\begin{aligned}
& K_0 (C_r^* (\mc R ,q)) \cong \mh Z^3 \\
& K_1 (C_r^* (\mc R ,q)) = 0
\end{aligned}
\end{equation}
and generating projections are
\[
p_{-1} = \left( 1, \begin{pmatrix} 0 & 0 \\ 0 & 0 \end{pmatrix} 
\right) \qquad p_a = \left( 0, \begin{pmatrix} 1 & 0 \\ 
0 & 0 \end{pmatrix} \right) \qquad p_b = \left( 0, 
\begin{pmatrix} 0 & 0 \\ 0 & 1 \end{pmatrix} \right)
\]
Now they are given by
\begin{equation}\label{eq:6.15}
\begin{array}{lllll}
p_a & = & (T_e + T_{s_1})(1+q)^{-1} & \mr{if} & q > 1 \\
p_a & = & (T_e + T_{s_1})(1+q)^{-1} - p_{-1} & \mr{if} & q < 1 \\
p_b & = & (q T_e - T_{s_1})(1+q)^{-1} - p_{-1} & \mr{if} & q > 1 \\
p_b & = & (q T_e - T_{s_1})(1+q)^{-1} & \mr{if} & q < 1 \\
\end{array}
\end{equation}

\begin{itemize}
\item {\large \textbf{generic case} $\mb{q_0 \neq q_1 \neq q_0^{-1}}$}
\item $P = \es$
\end{itemize}
\[
\begin{aligned}
& R_P = \es \qquad R_P^\vee = \es \\
& X^P = X \qquad X_P = 0 \qquad Y^P = Y \qquad Y_P = 0 \\
& T^P = T \qquad T_P = \{ 1 \} \qquad K_P = \{ 1 \} \\
& W^P = W(P,P) = \mc W_{PP} = W_0  \qquad W_P = \{ e \} \\
& \imath^o_{s_1} = (T_{s_1} (1+ \theta_1 ) + (1-q) \theta_1) 
(q + \theta_1 )^{-1} 
\end{aligned}
\]
In this case all unitary principal series representations are
irreducible:
\[
\begin{aligned}
& \pi (s_1 ,\es ,\delta_\es ,1) = 1 \qquad
\pi (s_1 ,\es ,\delta_\es ,-1) = 1 \\
& C_r^* (\mc R ,q)_P \cong C \big( [-1,1] ;M_2 (\mh C) \big) \\
& \mr{Prim} \big( C_r^* (\mc R ,q)_P \big) \cong S^1 / W_0 
\cong [-1,1]
\end{aligned}
\]
\begin{itemize}
\item $P = \{ \alpha \}$
\end{itemize}
\[
\begin{aligned}
& R_P = R_0 \qquad R_P^\vee = R_0^\vee \\
& X^P = 0 \qquad X_P = X \qquad Y^P = 0 \qquad Y_P = Y \\
& T^P = \{ 1 \} \qquad T_P = T \qquad K_P = \{ 1 \} \\
& W^P = W(P,P) = \mc W_{PP} = \{ e \} \qquad W_P = W_0 
\end{aligned}
\]
The residual points were already listed in \eqref{eq:6.9}. By
Proposition \ref{prop:3.30}.2 there are exactly two inequivalent
discrete series representations, $\delta_1$ and $\delta_{-1}$.
\[
C_r^* (\mc R ,q)_P \cong \mh C^2
\]
To write down the corresponding projections we have to 
distinguish four cases.
\begin{align}
& p_1 = \sum_{w \in W} (-1)^{\ell (w)} q(w)^{-1} T_w 
\Big( \sum_{w \in W} q(w)^{-1} \Big)^{-1} \qquad 
& \mr{if} \; q_0 > q_1^{-1} \label{eq:6.16} \\
& p_1 =  \sum_{w \in W} T_w \Big( \sum_{w \in W} q(w) 
\Big)^{-1} \qquad & \mr{if} \; q_0 < q_1^{-1} \label{eq:6.17} \\
& p_{-1} = \sum_{w \in W} (-q_1 )^{\ell_1 (w)/2} T_w \Big( 
\sum_{w \in W} q_1^{-\ell_1 (w)} q_0^{\ell_0 (w)} \Big)^{-1} 
\qquad & \mr{if} \; q_0 < q_1 \;\; \label{eq:6.18} \\
& p_{-1} = \sum_{w \in W} (-q_0 )^{\ell_0 (w)/2} T_w \Big( 
\sum\limits_{w \in W} q_0^{-\ell_0 (w)} q_1^{\ell_1 (w)} 
\Big)^{-1} \qquad & \mr{if} \; q_0 > q_1 \;\; \label{eq:6.19}
\end{align}
If we consider $\delta_1$ and $\delta_{-1}$ only as representations
of $\mc H (W_0 ,q)$, then in this list \eqref{eq:6.16} and 
\eqref{eq:6.18} are deformations of the sign representation, while 
\eqref{eq:6.17} and \eqref{eq:6.19} are deformations of the trivial 
representation.
\\[2mm]
We conclude that
\begin{align*}
& C_r^* (\mc R ,q) \cong \mh C \oplus C \big( [-1,1] ;M_2 (\mh C) 
\big) \oplus \mh C \\
& \raisebox{3mm}{Prim$\big( C_r^* (\mc R ,q) \big) \qquad \cong$}
\begin{picture}(6,1)
\linethickness{1mm}
\put(2,0.4){\line(1,0){3.5}}
\put(1,0.4){\circle*{0.15}}
\put(6.5,0.4){\circle*{0.15}}
\end{picture}
\end{align*}
Hence in the generic case also 
\begin{equation}\label{eq:6.20}
\begin{aligned}
& K_0 (C_r^* (\mc R ,q)) \cong \mh Z^3 \\
& K_1 (C_r^* (\mc R ,q)) = 0
\end{aligned}
\end{equation}
Canonical generators are $[p_{-1}] ,\, [p_1 ]$ and $[p_0 ]$, 
where $p_0$ is any rank one projector in 
$C \big( [-1,1] ;M_2 (\mh C) \big)$.
\\[2mm]
So for this root datum the $K$-groups are independent of the label
function. Moreover the various maps $K_0 (\phi_0 )$ all turn out
to be isomorphisms. We list the images of the projections 
$p_a ,\, p_b ,\, p_c$ and $p_d$ below, in that order.

\begin{equation}\label{eq:6.21}
\begin{array}{|c|c|c|}
\hline q_0 = q_1 < 1 & q_0^{-1} < q_1 > q_0 & q_0 = q_1 > 1 \\
\hline p_a & p_0 & p_a \\
p_b + p_{-1} & p_0 + p_1 + p_{-1} & p_b + p_1 \\
p_b & p_0 + p_1 & p_a + p_1 \\
p_a + p_{-1} & p_0 + p_{-1} & p_b \\
\hline q_1^{-1} > q_0 < q_1 & q_0 = q_1 = 1 & q_1^{-1} < q_0 > q_1 \\
\hline p_0 + p_1 & p_a & p_0 + p_{-1} \\
p_0 + p_{-1} & p_b & p_0 + p_1 \\
p_0 & p_c & p_0 + p_1 + p_{-1} \\
p_0 + p_1 + p_{-1} & p_d & p_0 \\
\hline q_0 = q_1 < 1 & q_0^{-1} > q_1 < q_0 & q_0^{-1} = q_1 < 1 \\
\hline p_a + p_1 & p_0 + p_1 + p_{-1} & p_a + p_{-1} \\
p_b & p_0 & p_b \\
p_a & p_0 + p_{-1} & p_b + p_{-1} \\
p_b + p_1 & p_0 + p_1 & p_a \\
\hline
\end{array}
\end{equation}
\\[2mm]

\section{$GL_2$}
\label{sec:6.2}

The simplest twodimensional root datum which is not a product of
two onedimensionale root data is \inde{$\mc R (GL_2 )$}.
\[
\begin{aligned}
& X = \mh Z^2 \qquad Q = \{ (n,-n) : n \in \mh Z \} \qquad 
  X^+ = \{ (m,n) \in \mh Z^2 : m \geq n \} \\
& Y = \mh Z^2 \qquad Q^\vee = \{ (n,-n) : n \in \mh Z \} \\
& T = (\mh C^\times )^2 \qquad t= (t_1 ,t_2 ) = (t(1,0),t(0,1)) \\
& R_0 = \{ \pm \alpha \} = \{ \pm (1,-1) \} = R_1 \\
& R_0^\vee = \{ \pm \alpha^\vee \} = \{ \pm (1,-1) \} = R_1^\vee \\
& F_0 = \alpha \qquad W_0 = \{ e ,s_\alpha \} \\
& s_1 = s_\alpha : (m,n) \to (n,m) \qquad s_0 = t_\alpha s_\alpha 
  = t_{(1,0)} s_1 t_{(-1,0)} : (m,n) \to (n+1 ,m-1) \\
& S_{\mr{aff}} = \{ s_0 ,s_1 \} \qquad W \neq W_{\mr{aff}} = 
  \langle s_0 , s_1 | s_0^2 = s_1^2 = e \rangle \\
& \Omega = \langle \omega \rangle = 
  \langle t_{(1,0)} s_1 \rangle \cong \mh Z
\end{aligned}
\]
For any label function $q$ we have $q(s_0 ) = q(s_1 ) = 
q_{\alpha^\vee}$, so we call this value $q$. There are no residual
point because $\mc R$ is not semisimple. We do have two residual
cosets of dimension one, namely
\[
\{ t \in T : t_1^{-1} t_2 = q \} \quad \mr{and} \quad
\{ t \in T : t_1^{-1} t_2 = q^{-1} \}
\]
So there are only two really different cases, $q=1$ and $q \neq 1$.

\begin{itemize}
\item {\large \textbf{group case} $\mb{q = 1}$}
\end{itemize}
As said before, we only need to look at unitary principal series 
representations. There is a single nonscalar intertwiner 
$\pi (s_1 ,\es ,\delta_\es ,t) : I_{(t_1 ,t_2 )} \to
I_{(t_2 ,t_1 )}$ It is given by right multiplication with $T_{s_1}$,
so with respect to the basis \eqref{eq:6.2} we have
\[
\mc S (W) \cong \left\{ f \in C^\infty \big( \mh T^2 ;M_2 (\mh C )
\big) : f(t_2 ,t_1 ) = \begin{pmatrix} 1 & 0 \\ 0 & -1 \end{pmatrix}
f (t_2 ,t_1 ) \begin{pmatrix} 1 & 0 \\ 0 & -1 \end{pmatrix} \right\}
\]
Let $M$ be the closed M\"obius strip and $\partial M$ its boundary.
We see that 
\[
C_r^* (W) \cong \{ f \in C (M;M_2 (\mh C )) : f(m) 
\text{ is diagonal if } m \in \partial M \}
\]
Consider the ideal
\[
A_1 = \{ f \in C (M;M_2 (\mh C )) : f(m) = \begin{pmatrix} * & 0 \\
0 & 0 \end{pmatrix} \; \mr{if} \; m \in \partial M \}
\]
According to Proposition \ref{prop:2.26} the inclusion
$A_1 \to C(M;M_2 (\mh C ))$ induces an isomorphism on $K$-theory.
However, $M$ is homotopy equivalent to a circle, so 
\[
K_0 (A_1) \cong K_1 (A_1 ) \cong \mh Z
\]
Moreover $\partial M \cong S^1$, so from
\[
0 \to A_1 \to C_r^* (W) \to C (\partial M) \to 0
\]
we get
\[
\hexagon{\mh Z}{K_0 (C_r^* (W))}{\mh Z}{\mh Z
}{K_1 (C_r^* (W))}{\mh Z}
\]
Now it is not difficult to see that the upper part of this hexagon
is exact, and hence
\begin{equation}
\begin{aligned}
K_0 \big( C_r^* (W) \big) \cong \mh Z^2 \\
K_1 \big( C_r^* (W) \big) \cong \mh Z^2 
\end{aligned}
\end{equation}
Generating projections and unitaries are
\begin{equation}\label{eq:6.22}
\begin{aligned}
& p_a = (T_e + T_{s_1})/2 \\
& p_b = (T_e - T_{s_1})/2 \\
& \theta_{(0,1)} \\
& u = \theta_{(0,1)} p_a + \theta_{(0,-1)} p_b
\end{aligned}
\end{equation}

\begin{itemize}
\item {\large \textbf{generic, equal label case} $\mb{q \neq 1}$}
\item $P = \es$
\end{itemize}
\[
\begin{aligned}
& R_P = \es \qquad R_P^\vee = \es \\
& X^P = X \qquad X_P = 0 \qquad Y^P = Y \qquad Y_P = 0 \\
& T^P = T \qquad T_P = \{ 1 \} \qquad K_P = \{ 1 \} \\
& W^P = W(P,P) = \mc W_{PP} = W_0  \qquad W_P = \{ e \} \\
& \imath^o_{s_1} = (T_{s_1} (1- \theta_{(1,-1)} ) + (1-q) 
\theta_{(1,-1)}) (q + \theta_{(1,-1)} )^{-1} 
\end{aligned}
\]
If $s_1 (t) = t$ then $\pi (s_1 ,\es ,\delta_\es 
,t) = 1$, so 
\[
\begin{aligned}
& C_r^* (\mc R ,q)_P \cong C(M;M_2 (\mh C)) \\
& \mr{Prim} \big( C_r^* (\mc R ,q)_P \big) \cong T_u / W_0 \cong M
\end{aligned}
\]
\begin{itemize}
\item $P = \{ \alpha \}$
\end{itemize}
\[
\begin{aligned}
& R_P = R_0 \qquad R_P^\vee = R_P \\
& X^P = X / \mh Z \alpha \cong \mh Z \qquad 
  X_P = X /(R_P^\vee )^\perp \cong \mh Z \alpha / 2 \\
& Y^P = Y \cap R_P^\perp = \mh Z (1,1) \qquad 
  Y_P = Y \cap \mh Q R_P^\vee = \mh Z \alpha^\vee \\
& T^P = \{ (t_1 ,t_1 ) : t_1 \in \mh C^\times \} \qquad
  T_P = \{ (t_1 ,t_1^{-1}) : t_1 \in \mh C^\times \} \\
& K_P = \{ (1,1) ,(-1,-1) \} = \{ 1, k_P \}
\end{aligned}
\]
The root datum $\mc R_P$ is isomorphic to $\mc R (A_1 )^\vee$, 
so we can use the description of the discrete series on page 
\pageref{eq:6.11}. The representations  $\pi (P ,\delta_1 
,(t_1 ,t_1 ))$ and\\ $\pi (P ,\delta_1 ,(-t_1 ,-t_1 ))$ are 
intertwined by $\pi (k_P )$.
\[
\begin{aligned}
& C_r^* (\mc R ,q)_P \cong C \big( S^1 \big) \\
& \mr{Prim} \big( C_r^* (\mc R ,q)_P \big) \cong 
\big( P, W_P (q^{1/2},q^{-1/2}), \delta_1 ,T_u^P \big) \cong S^1
\end{aligned}
\]
The associated projector is
\begin{equation}\label{eq:6.23}
\begin{array}{lllll}
p_\alpha & = & \sum\limits_{w \in W_{\mr{aff}}} (-q)^{\ell (w)} 
T_w \Big( \sum\limits_{w \in W_{\mr{aff}}} q(w)^{-1} \Big)^{-1} 
& \mr{if} & q > 1 \\
p_\alpha & = & \sum\limits_{w \in W_{\mr{aff}}} T_w \Big( 
\sum\limits_{w \in W_{\mr{aff}}} q(w) \Big)^{-1} & \mr{if} & 
q < 1 \\
\end{array}
\end{equation}
Adding these two summands we get 
\begin{equation}
\begin{aligned}
K_0 \big( C_r^* (\mc R ,q) \big) \cong \mh Z^2 \\
K_1 \big( C_r^* (\mc R ,q) \big) \cong \mh Z^2 
\end{aligned}
\end{equation}
These groups are generated by the classes of the projections
\[
p_\alpha = \left( \begin{pmatrix} 0 & 0 \\ 0 & 0 \end{pmatrix}
, 1 \right) \quad \mr{and} \quad p_0 = \left( \begin{pmatrix} 
1 & 0 \\ 0 & 0 \end{pmatrix}, 1 \right) 
\]
and of the invertibles
\[
u_\alpha = \left( \begin{pmatrix} 1 & 0 \\ 0 & 1 \end{pmatrix}
, \mr{id}_{S^1} \right) \quad \mr{and} \quad u_0 = \left( 
\begin{pmatrix} r & 0 \\ 0 & 1 \end{pmatrix}, 1 \right) 
\]
where $r: M \to S^1$ is a homotopy equivalence. Explicitly we
may take
\begin{equation}\label{eq:6.24}
\begin{array}{lllll}
u_\alpha & = & p_\alpha \theta_{(1,0)} + 1 - p_\alpha \\
u_0 & = & \theta_{(1,0)} (1 - p_\alpha ) + p_\alpha \\
p_0 & = & (T_e + T_{s_1})(1+q)^{-1} & \mr{if} & q > 1 \\
p_0 & = & (T_{s_1} - q T_e )(1+q)^{-1} & \mr{if} & q < 1
\end{array}
\end{equation}

Hence for this root system the $K$-theory of $\mc S (\mc R ,q))$
is independent of $q$. The group isomorphisms $K_0 (\phi_0 )$
are as follows:

\begin{equation}
\begin{array}{|c|c|c|}
\hline q < 1 & q = 1 & q > 1 \\
\hline p_0 + p_\alpha & p_a & p_0 \\
p_0 & p_b & p_0 + p_\alpha \\
u_0 & \theta_{(1,0)} & u_0 \\
u_0 u_\alpha & u & u_0 u_\alpha^{-1} \\
\hline
\end{array}
\end{equation}
\\[2mm]

\section{$A_2$}
\label{sec:6.3}

Somewhat unusually we will not consider $A_2$ as embedded in 
$\mh R^3$, but only as a twodimensional object. There are two 
semisimple root data with $R_0$ of type $A_2$, depending on 
whether $X$ is the root lattice or the weight lattice. The latter
case is easier, so let us draw this root system together with the
fundamental weights $x_3$ and $x_4$.

\begin{center}
\includegraphics[width=4cm,height=4cm]{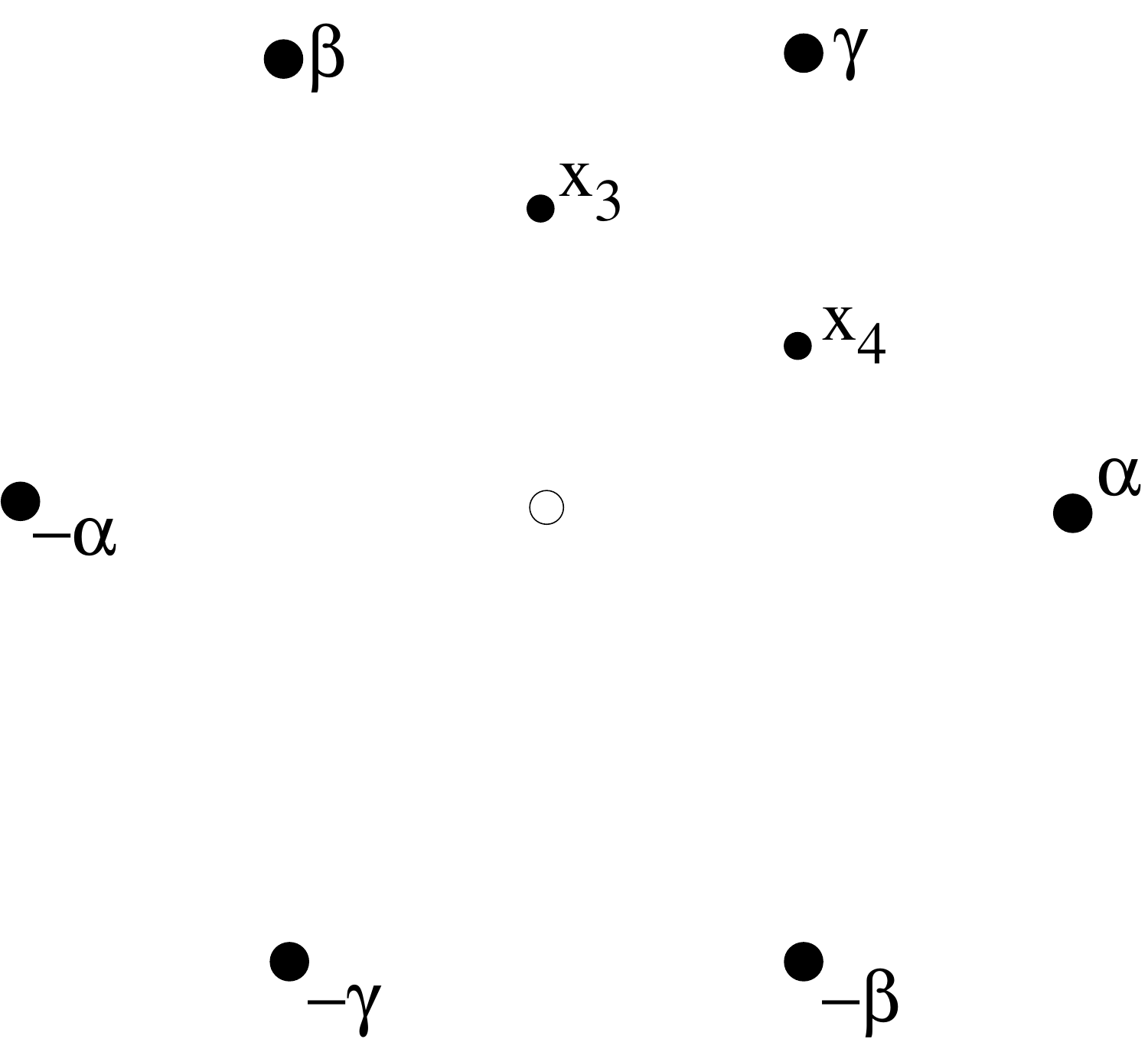}
\end{center}
The root datum \inde{$\mc R (A_2 )^\vee$} is described by
\[
\begin{aligned}
& X = \mh Z x_3 + \mh Z x_4 = \mh Z \big( 0, 3^{-1/2} \big) + 
  \mh Z \big( 1/2, 3^{-1/2} /2 \big) \qquad 
  X^+ = \mh N x_3 + \mh N x_4 \\
& Q = \mh Z (1,0) + \mh Z \big( -1/2, \sqrt 3 /2 \big) \qquad
  Y = Q^\vee = \mh Z (2,0) + \mh Z \big( 1,\sqrt 3 \big) \\
& T = ( \mh C^\times )^2 \qquad t = (t_3 ,t_4 ) = 
  (t (x_3 ), t(x_4 )) \\
& R_0 = \{ \pm \alpha , \pm \beta , \pm \gamma \} = \{ \pm (1,0) , 
  \pm \big( -1/2 , \sqrt 3 /2 \big) , \pm \big( 1/2, \sqrt 3 /2 
  \big) \} = R_1 \\
& R_0^\vee = \big\{ \pm \alpha^\vee , \pm \beta^\vee , \pm 
  \gamma^\vee \big\} = \{ \pm (1,0) , \pm \big( -1/2 , \sqrt 3 
  /2 \big) , \pm \big( 1/2, \sqrt 3 /2 \big) \} = R_1^\vee \\
& F_0 = \{ \alpha ,\beta \} \qquad 
  W_0 = \langle s_\alpha , s_\beta | s_\alpha^2 = s_\beta^2 = 
  (s_\alpha s_\beta )^3 = e \rangle \cong S_3 \\
& s_1 = s_\alpha : (n,m) \to (-n,m) \quad s_2 = s_\beta : 
  \big( n + m/2 ,m \sqrt 3 /2 \big) \to 
  \big( m + n/2, n \sqrt 3 /2 \big) \\
& s_0 = t_\gamma s_\gamma : \big( n + m/2 ,m \sqrt 3 /2 \big) \to 
  \big( (1+n-m)/2, (1-n-m) \sqrt 3 /2 \big) \\
& S_{\mr{aff}} = \{ s_0 ,s_1 ,s_2 \} \qquad 
  W \neq W_{\mr{aff}} = \langle s_0 ,W_0 | s_0^2 = (s_0 s_1 )^3 
  = (s_0 s_2 )^3 = e \rangle \\
& \Omega = \{ e ,\omega_1 ,\omega_2 \} = \{ e, t_{x_3} s_\beta 
  s_\alpha , t_{x_3} s_\alpha s_\beta \}
\end{aligned}
\]
For any label function $q$ we have
\[
q(s_0 ) = q(s_1 ) = q(s_2 ) = q(s_3 ) = q_{\alpha^\vee} =
q_{\beta^\vee} = q_{\gamma^\vee}
\]
so we denote this value simply by $q$.
\[
c_\eta = (1 - q^{-1} \theta_{-\eta} )(1 - \theta_{-\eta})^{-1}
\]
for $\eta \in \{ \alpha ,\beta ,\gamma \}$. Generically there are
6 tempered residual circles and 18 residual points. Representatives 
for the $W_0$-conjugacy classes are \index{zeta@$\zeta$}
\begin{equation}\label{eq:6.25}
(q,q) \quad (q \zeta ,q \zeta^2) \quad (q \zeta^2 ,q \zeta) 
\quad \mr{and} \quad 
\big\{ (q^{1/3} t_4^2 ,q^{2/3} t_4) : t_4 \in \mh T \big\} 
\end{equation}
where $\zeta = e^{2 \pi i /3}$ is a root of unity. 

\begin{itemize}
\item {\large \textbf{group case} $\mb{q = 1}$}
\end{itemize}

In the compact torus $T_u$ there are three $W_0$-invariant points:
\begin{equation}\label{eq:6.28}
(1,1) \qquad (\zeta ,\zeta^2 ) \qquad (\zeta^2 ,\zeta )
\end{equation}
Furthermore we have the following circles with nontrivial 
stabilizers:
\[
\begin{aligned}
& \{ (t_4^2 ,t_4 ) : t_4 \in \mh T \} & s_\alpha \\
& \{ (t_3 ,t_3^2 ) : t_3 \in \mh T \} & s_\beta \\
& \{ (t_3 ,t_3^{-1}) : t_3 \in \mh T \} & s_\gamma 
\end{aligned}
\]
These circles are conjugate under $W_0$. Therefore 
Prim$(\mc S (\mc R ,q))$ is a non-Hausdorff triangle whose interior 
is Hausdorff, whose edges are doubled and whose vertices are triple 
points. Let us call the underlying Hausdorff quotient space 
$T_u / W_0 \; D$, and its edges $D_\alpha ,\, D_\beta$ and 
$D_\gamma$, indicating their stabilizer.
\begin{center}
\includegraphics[width=6cm,height=4cm]{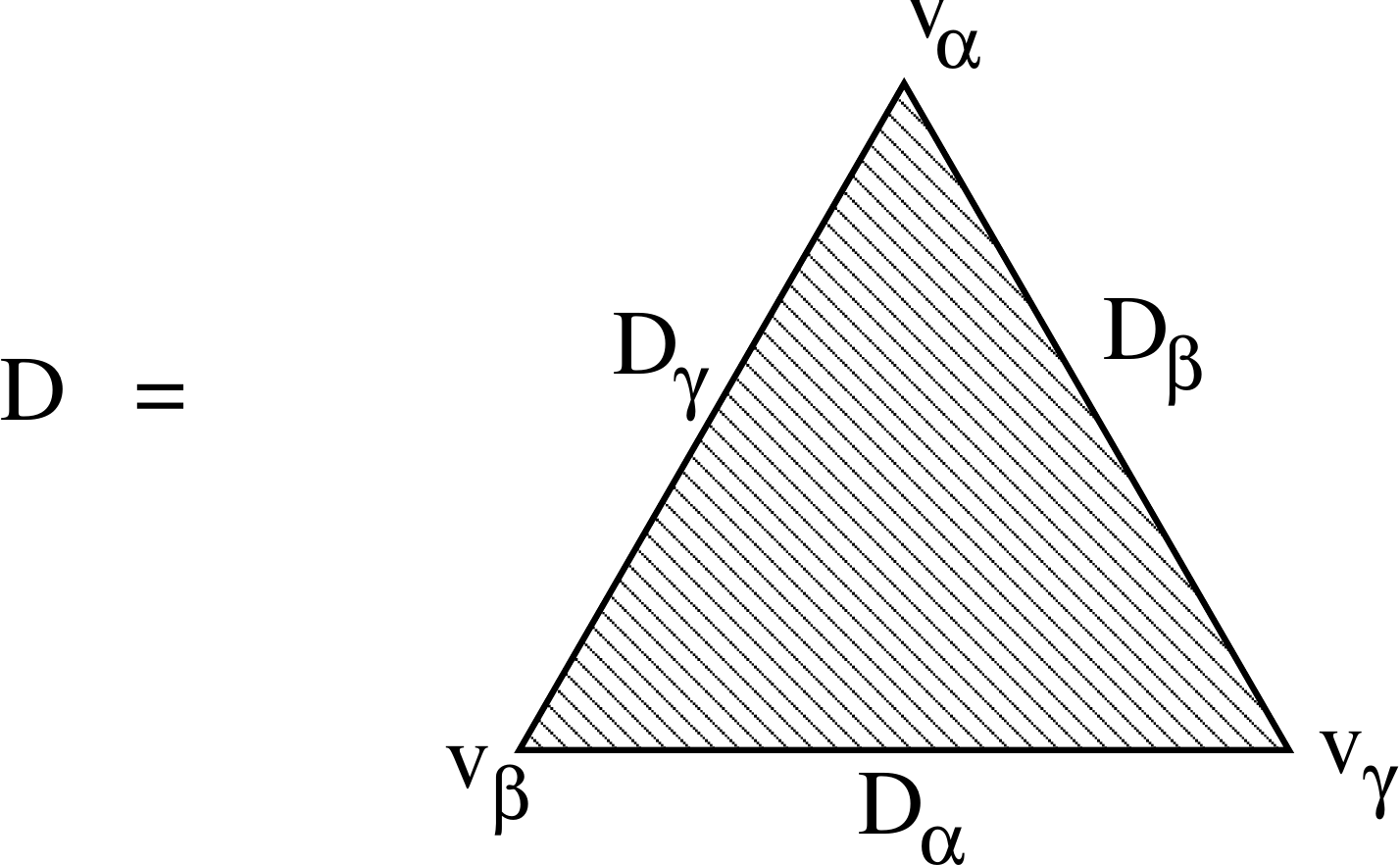}
\end{center}
There is no need to indicate which of the points \eqref{eq:6.28} is
$v_\alpha$, as all configurations occur, for a suitable choice of a
fundamental domain. 
\begin{multline*}
C_r^* (W) \cong \left\{ f \in C \big( D ; \mr{End} (\mh C [W_0 ]) 
\big): T_{s_\eta} f(t) T_{s_\eta} = f(t) \; \forall t \in D_\eta 
\: , \right. \\
\left. f(v_\eta ) \in \mr{End} (\mh C [W_0 ])^{W_0} \; 
\forall \eta \in \{ \alpha ,\beta, \gamma \} \right\}
\end{multline*}
Consider the extension
\begin{equation}\label{eq:6.29}
\begin{aligned}
& 0 \to C_0 \big( D ,\partial D ;\mr{End} (\mh C [W_0 ]) \big) \to
C_r^* (W) \to A_1 \to 0 \\
& A_1 = \left\{ f \in C \big( \partial D ;\mr{End} (\mh C [W_0 ]) 
\big): T_{s_\eta} f(t) T_{s_\eta} = f(t) \; \forall t \in D_\eta 
\: , \right. \\
& \qquad \qquad \left. f(v_\eta ) \in \mr{End} (\mh C [W_0 ])^{W_0} 
\; \forall \eta \in \{ \alpha ,\beta, \gamma \} \right\}
\end{aligned}
\end{equation}
In the vertices $v_\eta$ the $A_1$-representation $\mh C [W_0 ]$ is 
a direct sum of three irreducibles: the trivial $W_0$-representation,
the sign representation of $W_0$ and the defining (reflection) 
representation of $W_0$ (with multiplicity two). Likewise, on the 
edge $D_\eta \; \mh C [W_0 ]$ is the direct sum of a part 
corresponding to the trivial representation of $\{ e, s_\eta \}$ 
and a part corresponding to the trivial representation of 
$\{ e ,s_\eta \}$. Both summands have dimension three. Evaluating 
the reflection representations at the vertices gives an extension
\begin{equation}\label{eq:6.30}
\begin{aligned}
& 0 \to A_2^2 \to A_1 \to M_2 (\mh C)^3 \to 0 \\
& A_2 = \{ f \in C (\partial D ;M_3 (\mh C) ) : f(v_\eta ) \in 
\mh C \oplus O_2\}
\end{aligned}
\end{equation}
where $O_2$ is the $2 \times 2$ zero matrix. By Proposition 
\ref{prop:2.26} the inclusion $A_2 \to C (\partial D ;M_3 
(\mh C ))$ induces an isomorphism on $K$-theory, so we have 
an exact hexagon
\[
\hexagon{\mh Z^2}{K_0 (A_1 )}{\mh Z^3}{0}{K_1 (A_1 )}{\mh Z^2}
\]
The vertical maps are 0, so
\begin{align*}
& K_0 (A_1 ) \cong \mh Z^5 \\
& K_1 (A_1 ) \cong \mh Z^2 
\end{align*}
Plugging this into \eqref{eq:6.29} we get
\[
\hexagon{\mh Z}{K_0 \big( C_r^* (W) \big)}{\mh Z^5}{\mh Z^2}{K_1 
\big(C_r^* (W) \big)}{0}
\]
In terms of suitable generators the left vertical map becomes 
addition, and therefore
\begin{equation}\label{eq:6.31}
\begin{aligned}
& K_0 \big( C_r^* (W) \big) \cong \mh Z^5 \\
& K_1 \big( C_r^* (W) \big) \cong \mh Z
\end{aligned}
\end{equation}
Define the projections
\[
\begin{array}{lll}
p_{\mr{triv}} & = & \frac{1}{6} \sum\limits_{w \in W_0} T_w \\
p_{\mr{sign}} & = & \frac{1}{6} 
\sum\limits_{w \in W_0} (-1)^{\ell (w)} T_w 
\end{array}
\]
We can unambiguously define classes of projections 
$p_3 ,\, p_4 ,\, p_5$ by the requirements
\[
p_i (\zeta^j ,\zeta^{3-j}) \left\{ \begin{array}{lllll}
= & p_{\mr{triv}} + p_{\mr{sign}} & \mr{if} & i \neq j & 
\mr{mod}\: 3 \\ \perp & p_{\mr{triv}} + p_{\mr{sign}} & \mr{if} & 
i = j & \mr{mod}\: 3 \end{array} \right.
\]
Then $K_0 \big( C_r^* (W) \big)$ is generated by
\[
[p_{\mr{triv}}] \qquad [p_{\mr{sign}}] \qquad [p_3 ] \qquad [p_4 ]
\qquad [p_5 ]
\]
and a generator for $K_1 \big( C_r^* (W) \big)$ is 
\[
[u] = \left[ p_{\mr{triv}} N_{x_3} p_{\mr{triv}} + p_{\mr{sign}} 
N_{-x_3} p_{\mr{sign}} + T_e - p_{\mr{triv}} - p_{\mr{sign}} \right]
\]
\begin{itemize}
\item {\large \textbf{generic case} $\mb{q \neq 1}$}
\item $P = \es$
\end{itemize}
\[
\begin{aligned}
& R_P = \es \qquad  R_P^\vee = \es \\
& X^P = X \qquad X_P = 0 \qquad Y^P = Y \qquad Y_P = 0 \\
& T^P = T \qquad T_P = \{ 1 \} \qquad K_P = \{ 1 \} \\
& W^P = W(P,P) = \mc W_{PP} = W_0  \qquad W_P = \{ e \} \\
& \imath^o_{s_\eta} = (T_{s_\eta} (1- \theta_\eta ) + (q-1) 
\theta_\eta ) (q - \theta_\eta )^{-1} \qquad 
\eta \in \{ \alpha ,\beta ,\gamma \}
\end{aligned}
\]
If $s_\eta (t) = t$ then $\imath^o_{s_\eta} = 1$, and there are no
points with stabilizer $\{ e, s_1 s_2 ,s_2 s_1 \}$, so
\[
\begin{aligned}
& C_r^* (\mc R ,q)_P \cong C (D ;M_6 (\mh C)) \\
& \mr{Prim} \big( C_r^* (\mc R ,q)_P \big) \cong T_u / W_0 \cong D
\end{aligned}
\]
\begin{itemize}
\item $P = \{ \alpha \}$
\end{itemize}
\[
\begin{aligned}
& R_P = \{ \pm \alpha \} \qquad R_P^\vee = \{ \pm \alpha^\vee \} \\
& X^P = X / \mh Z \alpha \cong \mh Z x_4 \qquad 
  X_P = X / \mh Z x_3 \cong \mh Z \alpha / 2 \\
& Y^P = \mh Z (0,2 \sqrt 3 ) \qquad Y_P = \mh Z \alpha^\vee \\
& T^P = \{ (t_4^2 ,t_4 ) : t_4 \in \mh T \} \quad T_P = \{ (1,t_4 ) : 
  t_4 \in \mh T \} \quad K_P = \{ (1,1) ,(1,-1) \} = \{ 1,k_P \} \\
& W_p = \{ e ,s_\alpha \} \quad W^P = \{ e, s_\beta ,s_\alpha 
  s_\beta \} \quad W(P,P) = \{ e \} \quad \mc W_{PP} = K_P
\end{aligned}
\]
The root datum $\mc R_P$ is isomorphic to $\mc R (A_1 )^\vee$, 
so we can use the analysis on page \pageref{eq:6.8}. The 
representations $\pi (P,\delta_1 (1,t_4 ))$ and $\pi (P,\delta_{-1}
,(1,-t_4 ))$ are intertwined by $\pi (k_P )$. 
\[
\begin{aligned}
& C_r^* (\mc R ,q)_P \cong C (S^1 ;M_3 (\mh C )) \\
& \mr{Prim}( C_r^* (\mc R ,q)_P ) \cong \big( P,W_P (q^{1/3},
q^{2/3}), \delta_1 ,T_u^P \big) \cong S^1
\end{aligned}
\]
The corresponding central idempotent is 
$p_\alpha + p_\beta + p_\gamma$, where 
\begin{equation}
\begin{array}{lllll}
p_\eta & = & \sum\limits_{w \in W_\eta} (-q)^{\ell (w)} T_w \Big( 
\sum\limits_{w \in W_\eta} q(w)^{-1} \Big)^{-1} & \mr{if} & q > 1 \\
p_\eta & = & \sum\limits_{w \in W_\eta} T_w \Big( \sum\limits_{w \in 
W_\eta} q(w) \Big)^{-1} & \mr{if} & q < 1
\end{array}
\end{equation}
and $W_\eta = \langle s_\eta ,t_\eta \rangle \subset W_{\mr{aff}}$.

\begin{itemize}
\item $P = \{ \beta \}$
\end{itemize}
This subset of $F_0$ is conjugate to $\{ \alpha \}$ by 
$s_\alpha s_\beta$.

\begin{itemize}
\item $P = \{ \alpha, \beta \}$
\end{itemize}
\[
\begin{aligned}
& R_P = R_0 \qquad R_P^\vee = R_0^\vee \\
& X^P = 0 \qquad X_P = X \qquad Y^P = 0 \qquad Y_P = Y \\
& T^P = \{ 1 \} \qquad T_P = T \qquad K_P = \{ 1 \} \\
& W^P = W(P,P) = \mc W_{PP} = \{ e \} \qquad W_P = W_0 
\end{aligned}
\]
The residual points \eqref{eq:6.25} all carry exactly one 
inequivalent discrete series representation. We have 
$C_r^* (\mc R ,q) \cong \mh C^3$ with central idempotent
\begin{equation}\label{eq:6.35}
\begin{array}{lllll}
p_{\alpha ,\beta} & = & \sum\limits_{w \in W_{\mr{aff}}} 
(-q)^{\ell (w)} T_w \Big( \sum\limits_{w \in W_{\mr{aff}}} 
q(w)^{-1} \Big)^{-1} & \mr{if} & q > 1 \\
p_{\alpha ,\beta} & = & \sum\limits_{w \in W_{\mr{aff}}} T_w \Big( 
\sum\limits_{w \in W_{\mr{aff}}} q(w) \Big)^{-1} & \mr{if} & q < 1
\end{array}
\end{equation}
The projections for the specific points are
\begin{equation}
\begin{array}{lll}
p_{(1,1)} & = & (T_e + T_{\omega_1} + T_{\omega_2}) \,
  p_{\alpha ,\beta} / 3 \\
p_{(\zeta ,\zeta^2 )} & = &  (T_e + \zeta T_{\omega_1} + 
  \zeta^2 T_{\omega_2}) \, p_{\alpha ,\beta} / 3 \\
p_{(\zeta^2 ,\zeta )} & = &  (T_e + \zeta^2 T_{\omega_1} + 
  \zeta T_{\omega_2}) \, p_{\alpha ,\beta} / 3 \\
\end{array}
\end{equation}
Notice that they agree on $\mc H (W_{\mr{aff}} ,q)$.

We conclude that the spectrum of $C_r^* (\mc R ,q)$ is the 
Hausdorff space 
\[
\mr{Prim}\big( C_r^* (\mc R ,q) \big) \cong D \sqcup S^1 \sqcup
3 \;\mr{points}
\]
Since $D$ is compact and contractible this implies 
\begin{equation}
\begin{aligned}
& K_0 \big( C_r^* (\mc R ,q) \big) \cong \mh Z^5 \\
& K_1 \big( C_r^* (\mc R ,q) \big) \cong \mh Z
\end{aligned}
\end{equation}
Let $p_0$ and $p_1$ be rank one projectors in 
$C_r^* (\mc R ,q)_\es$ and $C_r^* (\mc R ,q)_{\{\alpha \}}$.
Then $p_0 ,p_1 , \\ p_{(1,1)} ,p_{(\zeta ,\zeta^2 )}$ and 
$p_{(\zeta^2 ,\zeta )}$ generate $K_0 \big( C_r^* (\mc R ,q) \big)$,
while $T_{x_3}$ generates $K_1 \big( C_r^* (\mc R ,q) \big)$.
\\[3mm]

Once again, the $K$-theory turns out to be independent of the
parameters. In terms of all the above generators, the group 
isomorphisms $K_* (\phi_0 )$ are as follows.
\begin{equation}
\begin{array}{|c|c|c|}
\hline q < 1 & q = 1 & q > 1 \\
\hline p_0 + p_1 + p_{(1,1)} + p_{(\zeta ,\zeta^2 )} + 
p_{(\zeta^2 ,\zeta )} & p_{\mr{triv}} & p_0 \\
p_0 & p_{\mr{sign}} & p_0 + p_1 + p_{(1,1)} + 
p_{(\zeta ,\zeta^2 )} + p_{(\zeta^2 ,\zeta )} \\
2 p_0 + p_1 + p_{(\zeta ,\zeta^2 )} + p_{(\zeta^2 ,\zeta )} &
p_3 & 2 p_0 + p_1 + p_{(\zeta ,\zeta^2 )} + p_{(\zeta^2 ,\zeta )} \\
2 p_0 + p_1 + p_{(1,1)} + p_{(\zeta^2 ,\zeta )} & p_4 &
2 p_0 + p_1 + p_{(1,1)} + p_{(\zeta^2 ,\zeta )} \\
2 p_0 + p_1 + p_{(1,1)} + p_{(\zeta ,\zeta^2 )} & p_5 &
2 p_0 + p_1 + p_{(1,1)} + p_{(\zeta ,\zeta^2 )} \\
T_{x_3} & u & T_{x_3}^{-1} \\
\hline
\end{array}
\end{equation}
\\[3mm]

As promised, we also discuss the root datum $\mc R (A_2 )$, where
$X$ is the root lattice.
\[
\begin{aligned}
& X = Q = \mh Z (1,0) + \mh Z \big( 1/2 ,\sqrt 3 /2 \big) \\
&  X^+ = \{ \big( n + m/2 ,m \sqrt 3 /2 \big) : 
  0 \leq m \leq 2n \leq 4m \} \\
& Y = \mh Z \big( 0,2 / \sqrt 3 \big) + \mh Z \big(1 , 1/ \sqrt 3 
  \big) \qquad 
  Q^\vee = \mh Z (1,0) + \mh Z \big( 1/2 ,\sqrt 3 /2 \big) \\
& R_0 = \{ \pm \alpha , \pm \beta , \pm \gamma \} = \{ \pm (1,0) , 
  \pm \big( -1/2 , \sqrt 3 /2 \big) , \pm \big( 1/2, \sqrt 3 /2 
  \big) \} = R_1 \\
& R_0^\vee = \big\{ \pm \alpha^\vee , \pm \beta^\vee , 
  \pm \gamma^\vee \big\} = \{ \pm (1,0) , \pm \big( -1/2 , \sqrt 3 
  /2 \big) , \pm \big( 1/2, \sqrt 3 /2 \big) \} = R_1^\vee \\
& T = (\mh C^\times )^2 \qquad 
  t = (t_1 ,t_2 ) = (t (\alpha ), t (\beta )) \\
& F_0 = \{ \alpha ,\beta \} \qquad 
  W_0 = \langle s_\alpha , s_\beta | s_\alpha^2 = s_\beta^2 = 
  (s_\alpha s_\beta )^3 = e \rangle \cong S_3 \\
& s_1 = s_\alpha : (n,m) \to (-n,m) \quad s_2 = s_\beta : 
  (n + m/2 ,m \sqrt 3 /2) \to (m + n/2, n \sqrt 3 /2) \\
& s_0 = t_\gamma s_\gamma : \big( n + m/2 ,m \sqrt 3 /2 \big) \to 
  \big( (1+n-m)/2, (1-n-m) \sqrt 3 /2 \big) \\
& S_{\mr{aff}} = \{ s_0 ,s_1 ,s_2 \} \qquad 
  W = W_{\mr{aff}} = \langle s_0 ,W_0 | s_0^2 = (s_0 s_1 )^3 
  = (s_0 s_2 )^3 = e \rangle \\
& q(s_0 ) = q(s_1 ) = q(s_2 ) = q(s_3 ) = q_{\alpha^\vee} =
  q_{\beta^\vee} = q_{\gamma^\vee} := q
\end{aligned}
\]
Generically there are 6 residual points and 6 tempered residual
circles. Both form a single $W_0$-conjugacy class, typical 
examples being
\[
(q^{-1},q^{-1}) \quad \mr{and} \quad 
\{ (q^{-1},t_2 : t_2 \in \mh T \}
\]
As usual we distinguish the cases $q = 1$ and $q \neq 1$.

\begin{itemize}
\item {\large \textbf{group case} $\mb{q = 1}$}
\end{itemize}
The following subtori of $T_u$ have nontrivial stabilizers:
\[
\begin{array}{ll}
\{ (1,t_2 ) : t_2 \in \mh T \} & \{ e ,s_\alpha \} \\
\{ (t_1 ,1) : t_1 \in \mh T \} & \{ e ,s_\beta \} \\
\{ (t_1 ,t_1^{-1}) :t_1 \in \mh T \} & \{ e ,s_\gamma \} \\
(\zeta ,\zeta) , (\zeta^2 ,\zeta^2) & 
  \{ e, s_\alpha s_\beta, s_\beta s_\alpha \} \\
(1,1) & W_0
\end{array}
\]
The following part $T'$ of $T_u$ is a fundamental domain for
the action of $W_0$. 
\begin{center}
\includegraphics[width=4cm,height=5cm]{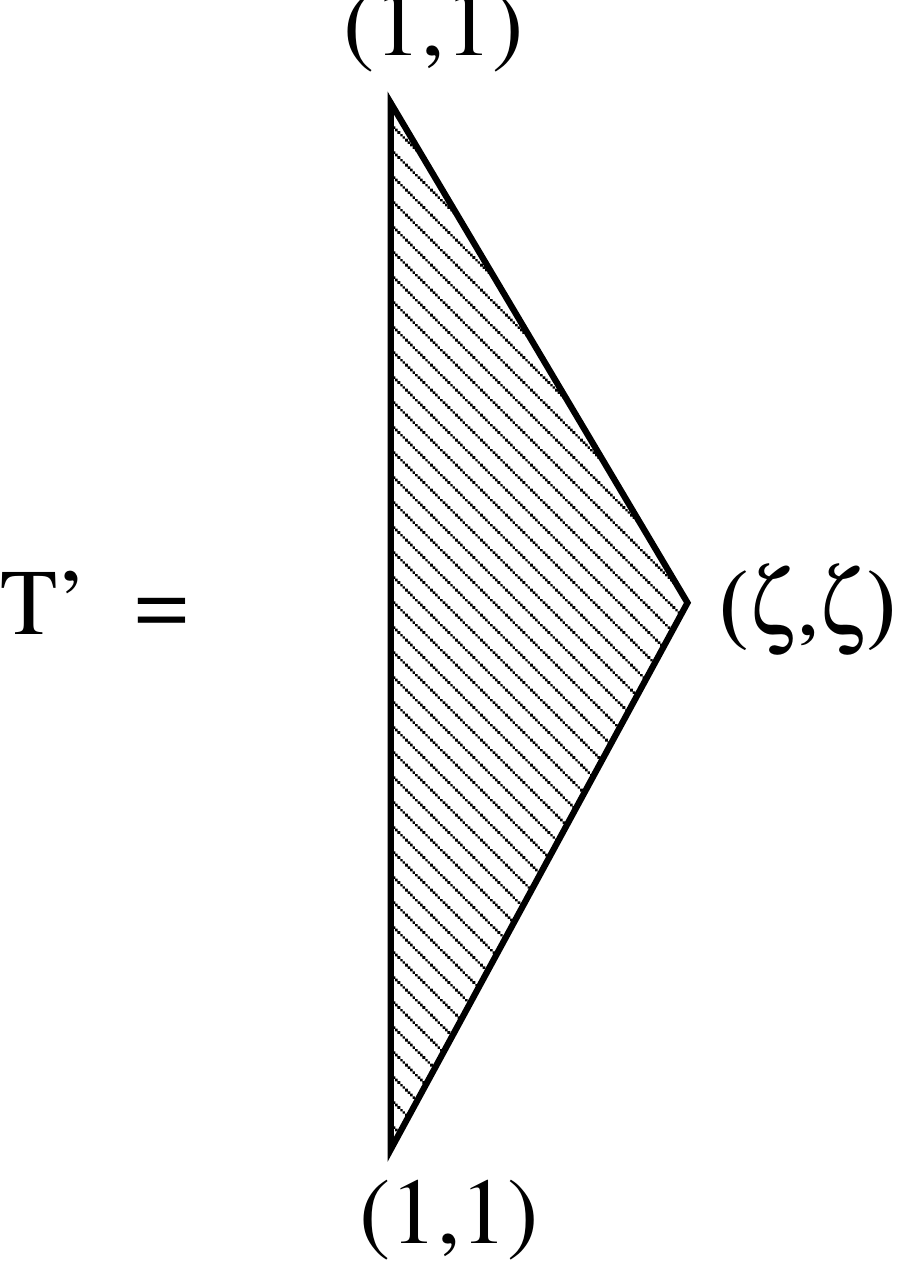}
\end{center}
The left edge is $T_u^{s_\alpha}$, and to get $T_u / W_0$ we
only have to identify the other two edges by means of a rotation
around $(\zeta ,\zeta)$. Then $C_r^* (W)$ consists of all
$f \in C \big( T' ; \mr{End}( \mh C [W_0 ]) \big)$ such that
\begin{enumerate}
\item $T_{s_\alpha} f(t) T_{s_\alpha} = f(t)$ if 
$t \in T_u^{s_\alpha}$
\item $f(1,1) \in \big( \mr{End}( \mh C [W_0 ]) \big)^{W_0}$
\item $f(s_\beta s_\alpha (t_1 ,t_1 )) = T_{s_\beta s_\alpha }
f(t_1 t_1 ) T_{s_\beta s_\alpha } \quad \forall t_1 \in
\exp (\pi i [0,2/3] )$
\end{enumerate}
If $t \in T_u^{s_\alpha}$ then $f(t)$ stabilizes 
$\mh C [W_0 ]^{s_\alpha}$, so there are extensions
\begin{align}
& \label{eq:6.32} 0 \to A_2 \to C_r^* (W) \to A_1 \to 0 \\
& \label{eq:6.33} 0 \to A_3 \to A_1 \to \mh C \to 0 \\
& \nonumber A_1 = \left\{ f \in C \Big( T_u^{s_\alpha} ; \mr{End} 
\big( \mh C [W_0 ]^{s_\alpha} \big) \Big) : f(1,1) = 
\begin{pmatrix} * & 0 & 0 \\ 0 & * & * \\ 0 & * & * \end{pmatrix}
\right\} \\
& \nonumber A_2 = \left\{ f \in C_r^* (W) : f(1,1) \in \mh C 
\oplus O_5 ,\, f(t) \in M_3 (\mh C) \oplus O_3 \: \forall t \in 
T_u^{s_\alpha} \right\} \\
& \nonumber A_3 = \{ f \in A_1 : f(1,1) \in O_1 \oplus 
M_2 (\mh C ) \}
\end{align}
Here $O_n$ denotes the $n \times n$ zero matrix. By Proposition
\ref{prop:2.26} the inclusions 
\begin{equation}\label{eq:6.34}
\begin{aligned}
& A_2 \to A_4 := \{ f \in C \big( T' ; \mr{End}( \mh C [W_0 ]) 
\big) : \text{ 3. holds } \} \\
& A_3 \to C \big( T_u^{s_\alpha} ; \mr{End} 
\big( \mh C [W_0 ]^{s_\alpha} \big) \big) 
\end{aligned}
\end{equation}
induce isomorphisms on $K$-theory. With the help of Lemma 
\ref{lem:2.17} we find that 
\[
\begin{array}{lll@{\qquad}lll}
K_0 (A_4 ) & \cong & \mh Z & K_1 (A_4 ) & = & 0 \\
K_0 (A_3 ) & \cong & \mh Z & K_1 (A_3 ) & \cong & \mh Z \\
K_0 (A_2 ) & \cong & \mh Z^3 & K_1 (A_2 ) & = & 0 \\
K_0 (A_1 ) & \cong & \mh Z^2 & K_1 (A_1 ) & \cong & \mh Z
\end{array}
\]
From \eqref{eq:6.32} we get an exact hexagon
\[
\hexagon{\mh Z^3}{K_0 \big( C_r^* (W) \big)}{\mh Z^2}{\mh Z}{
K_1 \big( C_r^* (W) \big)}{0}
\]
The upper row is exact, so
\begin{equation}
\begin{aligned}
& K_0 \big( C_r^* (W) \big) \cong \mh Z^5 \\
& K_1 \big( C_r^* (W) \big) \cong \mh Z
\end{aligned}
\end{equation}
It is rather difficult to write down explicit generators, so we
only indicate what they look like. Consider the projections
\[
\begin{array}{lll}
p_{\mr{triv}} & = & \frac{1}{6} \sum\limits_{w \in W_0} T_w \\
p_{\mr{sign}} & = & \frac{1}{6} 
\sum\limits_{w \in W_0} (-1)^{\ell (w)} T_w \\
p_{\mr{rot}} & = &  (T_e + \zeta T_{s_\alpha s_\beta} + 
  \zeta^2 T_{s_\beta s_\alpha}) / 3
\end{array}
\]
With these we can define classes of projections $[p_\zeta ]$ 
and $[p_{\zeta^2} ]$ by the conditions
\begin{align*}
& p_\zeta \big( T_u^{s_\alpha} \big) = 
p_{\zeta^2} \big( T_u^{s_\alpha} \big) = p_{\mr{triv}} \\
& p_\zeta (\zeta ,\zeta) (p_{\mr{triv}} + p_{\mr{sign}} + 
p_{\mr{rot}} ) = 0 \\
& p_{\zeta^2} (\zeta ,\zeta) p_{\mr{rot}} = 
p_{\zeta^2} (\zeta ,\zeta)
\end{align*}
These five classes of projections generate 
$K_0 \big( C_r^* (W) \big)$, and a generator for\\ 
$K_1 \big( C_r^* (W) \big)$ is
\[
u = p_{\mr{triv}} \theta_\beta \, p_{\mr{triv}} + p_{\mr{sign}}
\theta_{-\beta} \, p_{\mr{sign}} + T_e - p_{\mr{triv}} -
 p_{\mr{sign}}
\]

\begin{itemize}
\item {\large \textbf{generic case} $\mb{q \neq 1}$}
\item $P = \es$
\end{itemize}
\[
\begin{aligned}
& R_P = \es \qquad  R_P^\vee = \es \\
& X^P = X \qquad X_P = 0 \qquad Y^P = Y \qquad Y_P = 0 \\
& T^P = T \qquad T_P = \{ 1 \} \qquad K_P = \{ 1 \} \\
& W^P = W(P,P) = \mc W_{PP} = W_0  \qquad W_P = \{ e \} \\
& \imath^o_{s_\eta} = (T_{s_\eta} (1- \theta_\eta ) + (q-1) 
\theta_\eta ) (q - \theta_\eta )^{-1} \qquad 
\eta \in \{ \alpha ,\beta ,\gamma \}
\end{aligned}
\]
If $s_\eta (t) = t$ then $\imath^o_{s_\eta}(t) = 1$.
There are two points in $T_u$ whose stabilizer is not generated
by reflections: $(\zeta ,\zeta)$ and $(\zeta^2 ,\zeta^2 )$.
Let $A_4$ be as in \eqref{eq:6.34}. Then
\[
\begin{aligned}
& C_r^* (\mc R ,q)_P \cong A_4 \\
& \raisebox{9mm}{$\mr{Prim} \big( C_r^* (\mc R ,q)_P \big) 
\quad \cong \quad$} 
\includegraphics[width=2cm,height=2cm]{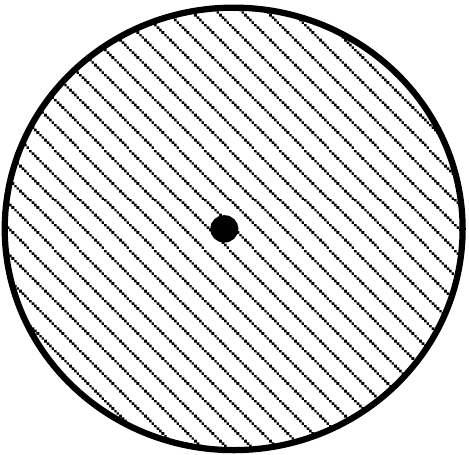}
\end{aligned}
\]
which is supposed to depict $T_u / W_0$ with in the center, 
instead of just $W_0 (\zeta ,\zeta)$, a triple non-Hausdorff 
point. By Lemma \ref{lem:2.17} we have
\begin{align*}
& K_0 \big( C_r^* (\mc R ,q)_P \big) \cong \mh Z^3 \\
& K_1 \big( C_r^* (\mc R ,q)_P \big) = 0
\end{align*}
The class of a projection in this algebra is completely 
determined by its value at $(\zeta ,\zeta )$, and over there 
we only have to say which representation of $\{ e, s_1 s_2, 
s_2 s_1 \}$ it gives. So we have generators $p_0 ,p_1 ,p_2$ with
\[
p_i (\zeta ,\zeta ) \imath^o_{s_1 s_2}(\zeta ,\zeta) = 
\zeta^i p_i (\zeta ,\zeta ) 
\]

\begin{itemize}
\item $P = \{ \alpha \}$
\end{itemize}
\[
\begin{aligned}
& R_P = \{ \pm \alpha \} \quad R_P^\vee = \{ \pm \alpha^\vee \} \\
& X^P = X / \mh Z \alpha \cong \mh Z \beta \quad
  X_P = X / \mh Z (0,2) \cong \mh Z (\alpha / 2) \\
& Y^P = \mh Z (0, 2 / \sqrt 3 ) \qquad 
  Y_P = \mh Z \alpha^\vee = \mh Z (2,0) \\
& T^P = \{ (1,t_2 ) : t_2 \in \mh C^\times \} \qquad
  T_P = \{ (t_2^{-2},t_2 ) : t_2 \in \mh C^\times \} \\
& K_P = \{ (1,1) ,(1,-1) \} = \{ 1 , k_P \} \\
& W_P = \{e ,s_\alpha \} \quad W^P = \{e ,s_\beta 
  ,s_\alpha s_\beta \} \quad W (P,P) = \{e\} \quad 
  \mc W_{PP} = K_P 
\end{aligned}
\]
The root datum $\mc R_P$ is isomorphic to $\mc R (A_1 )^\vee$,
so we already know its discrete series. The representations
$\pi (P,\delta_1 ,(1,t_2 ))$ and $\pi (P,\delta_{-1},(1,-t_2 ))$
are interwined by $\pi (k_P )$, so
\[
\begin{aligned}
& \mr{Prim} \big( C_r^* (\mc R ,q)_P \big) \cong 
\big( P, W_P (q,q^{1/2}), \delta_1 ,T_u^P \big) \cong S^1 \\
& C_r^* (\mc R ,q)_P \cong C (S^1 ; M_3 (\mh C ))
\end{aligned}
\]
The central idempotent for this component is 
$p_\alpha + p_\beta + p_\gamma$, as given by \eqref{eq:6.35}.

\begin{itemize}
\item $P = \{ \beta \}$
\end{itemize}
This subset of $F_0$ is conjugate to $\{ \alpha \}$ by
$s_\alpha s_\beta$.

\begin{itemize}
\item $P = \{ \alpha ,\beta \}$
\end{itemize}
\[
\begin{aligned}
& R_P = R_0 \quad R_P^\vee = R_0^\vee \\
& X^P = 0 \quad X_P = X \\
& Y^P = 0 \quad Y_P = Y \\
& T^P = \{1\} \quad T_P = T \quad K_P = \{1\} \\
& W_P = W_0 \quad W^P = W(P,P) = \mc W_{PP} = \{e\}
\end{aligned}
\]
Up to $W_0$-conjugacy there is a single residual point, so by
Proposition \ref{prop:3.30}.2 there is a unique discrete series
representation $\delta$. We have $C_r^* (\mc R ,q)_P \cong \mh C$,
and the corresponding projector is
\begin{equation}\label{eq:6.36}
\begin{array}{lllll}
p_\delta & = & \sum\limits_{w \in W} (-q)^{\ell (w)} T_w \Big( 
\sum\limits_{w \in W} q(w)^{-1} \Big)^{-1} & \mr{if} & q > 1 \\
p_\delta & = & \sum\limits_{w \in W} T_w \Big( 
\sum\limits_{w \in W} q(w) \Big)^{-1} & \mr{if} & q < 1
\end{array}
\end{equation}

On the whole we found that 
\begin{equation}
\begin{aligned}
& K_0 \big( C_r^* (\mc R ,q) \big) \cong \mh Z^5 \\
& K_1 \big( C_r^* (\mc R ,q) \big) \cong \mh Z
\end{aligned}
\end{equation}
which is the same as for $q=1$. Generators are the invertible
$\theta_\beta$ and the rank one projectors $p_0 ,p_1 ,p_2 
,p_\delta$ and $p_{\{ \alpha \}}$.

The isomorphisms $K_* (\phi_0 )$ take the following form:
\begin{equation}
\begin{array}{|c|c|c|}
\hline q < 1 & q = 1 & q > 1 \\
\hline p_0 + p_{\{ \alpha \}} + p_\delta & p_{\mr{triv}} & p_0 \\
p_0 & p_{\mr{sign}} & p_0 + p_{\{ \alpha \}} + p_\delta \\
2 p_2 + p_{\{ \alpha \}} & p_{\mr{rot}} & 2 p_2 + p_{\{ \alpha \}} \\
p_1 + p_{\{ \alpha \}} + p_\delta & p_\zeta & p_1 \\
p_2 + p_{\{ \alpha \}} + p_\delta & p_{\zeta^2} & p_2 \\
\theta_\beta & u & \theta_\beta^{-1} \\
\hline
\end{array}
\end{equation}

\section{$B_2$}
\label{sec:6.4}

The rank two root system $B_2$ is probably the best testcase for
Conjecture \ref{conj:5.28}. Because there are roots of different
lengths the conjecture is not yet known, and at the same time the 
calculations are still manageable. Moreover some interesting 
phenomena already occur for these affine Hecke algebras, like 
residual points that carry several inequivalent discrete series 
representations.

We will only consider the root datum $\mc R (B_2 )^\vee$ where
$X$ is the weight lattice:
\[
\begin{aligned}
& X = \mh Z^2 \quad Q = \{ (m,n) \in \mh Z^2 : n+m \text{ is even } 
\} \quad X^+ = \{ (m,n) \in \mh Z^2 : n \geq m \geq 0 \} \\
& Y = Q^\vee = \mh Z^2 \\
& T = ( \mh C^\times )^2 \qquad t = (t_1 ,t_2 ) = (t(1,0),t(0,1)) \\
& R_0 = R_1 = \{ \pm \alpha_1 , \pm \alpha_2 , \pm \alpha_3 ,
\pm \alpha_4 \} = \{ \pm (2,0) , \pm (-1,1) , \pm (1,1) , \pm (0,2) 
\} \\
& R_0^\vee = R_1^\vee = \{ \pm \alpha_1^\vee , \pm \alpha_2^\vee , 
\pm \alpha_3^\vee , \pm \alpha_4^\vee \} = \{ \pm (1,0) , \pm (-1,1) 
, \pm (1,1) , \pm (0,1) \} \\
& F_0 = \{ \alpha_1 , \alpha_2 \} \quad W_0 = \langle s_1 ,s_2 | 
s_1^2 = s_2^2 = (s_1 s_2 )^4 = e \rangle \cong D_4 \\
& s_i = s_{\alpha_i} \quad s_0 = t_{\alpha_3} s_{\alpha_3} =
(t_{(1,0)} s_2 s_1 ) s_2 (t_{(1,0)} s_2 s_1 )^{-1} :
(m,n) \to (1-n ,1-m) \\
& S_{\mr{aff}} = \{ s_0 ,s_1 ,s_2 \} \quad W \neq W_{\mr{aff}} =
\langle s_0 ,W_0 | s_0^2 = (s_0 s_2 )^2 = (s_0 s_1 )^4 = e \rangle \\
& \Omega = \{ e, t_{(1,0)} s_1 \} = \{ e , \omega \} \\
& q_1 := q(s_1 ) = q_{\alpha_1^\vee} = q_{\alpha_3^\vee} \qquad
q_2 := q(s_2 ) = q(s_0 ) = q_{\alpha_2^\vee} = q_{\alpha_4^\vee} \\
& c_{\alpha_i} = (1 - q_{\alpha_i^\vee}^{-1} \theta_{- \alpha_i} )
(1 - \theta_{- \alpha_i})^{-1} \quad i = 1,2,3,4
\end{aligned}
\]
Generically there are 24 tempered residual circles and 40
residual points. Representatives are
\begin{align}
& \label{eq:6.63} \big\{ (q_1^{1/2},t_2 ) : t_2 \in \mh T \big\} \\
& \label{eq:6.64} \big\{ (-q_1^{1/2},t_2 ) : t_2 \in \mh T \big\} \\
& \label{eq:6.65} \big\{ (q_2^{1/2} t_1 ,q_2^{-1/2} t_1 ) : 
t_1 \in \mh T \big\} \\
& (q_1^{-1/2},q_1^{1/2} q_2^{-1}) \,,\; (q_1^{1/2}, q_1^{-1/2} 
q_2^{-1}) \,,\; (-q_1^{-1/2},-q_1^{1/2} q_2^{-1}) \,,\; \nonumber \\
& \label{eq:6.59} (-q_1^{1/2},-q_1^{-1/2} q_2^{-1}) \,,\;
(-q_1^{-1/2},q_1^{1/2})
\end{align}
It turns out that there are five classes of parameters with the
same level of genericity. The first three are easily found from
\eqref{eq:6.63} - \eqref{eq:6.65}: $q_1 = 1 = q_2 ,\, q_1 \neq 
1 = q_2$ and $q_1 = 1 \neq q_2$. Furthermore we have the generic
class and the four special lines
\begin{equation}\label{eq:6.66}
q_1 = q_2 \neq 1 ,\; q_1 = q_2^{-1} \neq 1 ,\; q_1 = q_2^2 \neq 1
,\; q_1 = q_2^{-2} \neq 1
\end{equation}

\begin{itemize}
\item {\large \textbf{group case} $\mb{q_1 = q_2 = 1}$}
\end{itemize}
From \eqref{eq:6.37} and \eqref{eq:6.38} we see that it pays off 
to determine the extended quotient 
\[
\widetilde{T_u} /W_0 \cong \bigsqcup_{\langle w \rangle \in \langle 
W_0 \rangle} T_u^w \big/ Z_{W_0} (w)
\]
There are five conjugacy classes, namely
\[
\{ e \} ,\, \{ s_1 ,s_4 \} ,\, \{ s_2 ,s_3 \} ,\, \{ s_1 s_2 ,s_2 
s_1 \} ,\, \{ s_1 s_2 s_1 s_2 \}
\]
We pick the elements written first as representatives.
\[
\begin{array}{llll}
w & Z_{W_0}(w) & T_u^w & T_u^w / Z_{W_0}(w) \\
\hline 
e & W_0 & T_u & \{ (t_1 , t_2) \in [0,1]^2 : t_1 \geq t_2 \} \\
s_1 & \{ e ,s_1 ,s_4 ,s_1 s_4 \} & \{ (\pm 1,t_2) \in T_u \}& 
\{ \pm 1 \} \times [-1,1] \\
s_2  & \{ e ,s_2 ,s_3 ,s_2 s_3 \} & \{ (t_1 ,t_1 ) \in T_u \} &
[-1,1] \\
s_1 s_2 & \langle s_1 s_2 \rangle & \{(1,1) ,(-1,-1)\}
& \{(1,1) ,(-1,-1)\} \\
(s_1 s_2 )^2 & W_0 & \{ (\pm 1,\pm 1),(\pm 1,\mp 1) \} &
\{ (1,1),(-1,-1),(1,-1) \} 
\end{array}
\]
Since all components of this space are contractible,
with \eqref{eq:6.38} we find that 
\begin{equation}
\begin{aligned} 
& K_0 \big( C_r^* (W) \big) \cong \check H^* \big( \widetilde{T_u} 
/ W_0 ; \mh Z \big) \cong \mh Z^9 \\
& K_1 \big( C_r^* (W) \big) = 0
\end{aligned}
\end{equation}
We may visualize the extended quotient as
\[
\raisebox{2cm}{$\hspace{1cm} \widetilde{T_u} / W_0 \quad \cong \quad$}
\includegraphics[width=45mm,height=45mm]{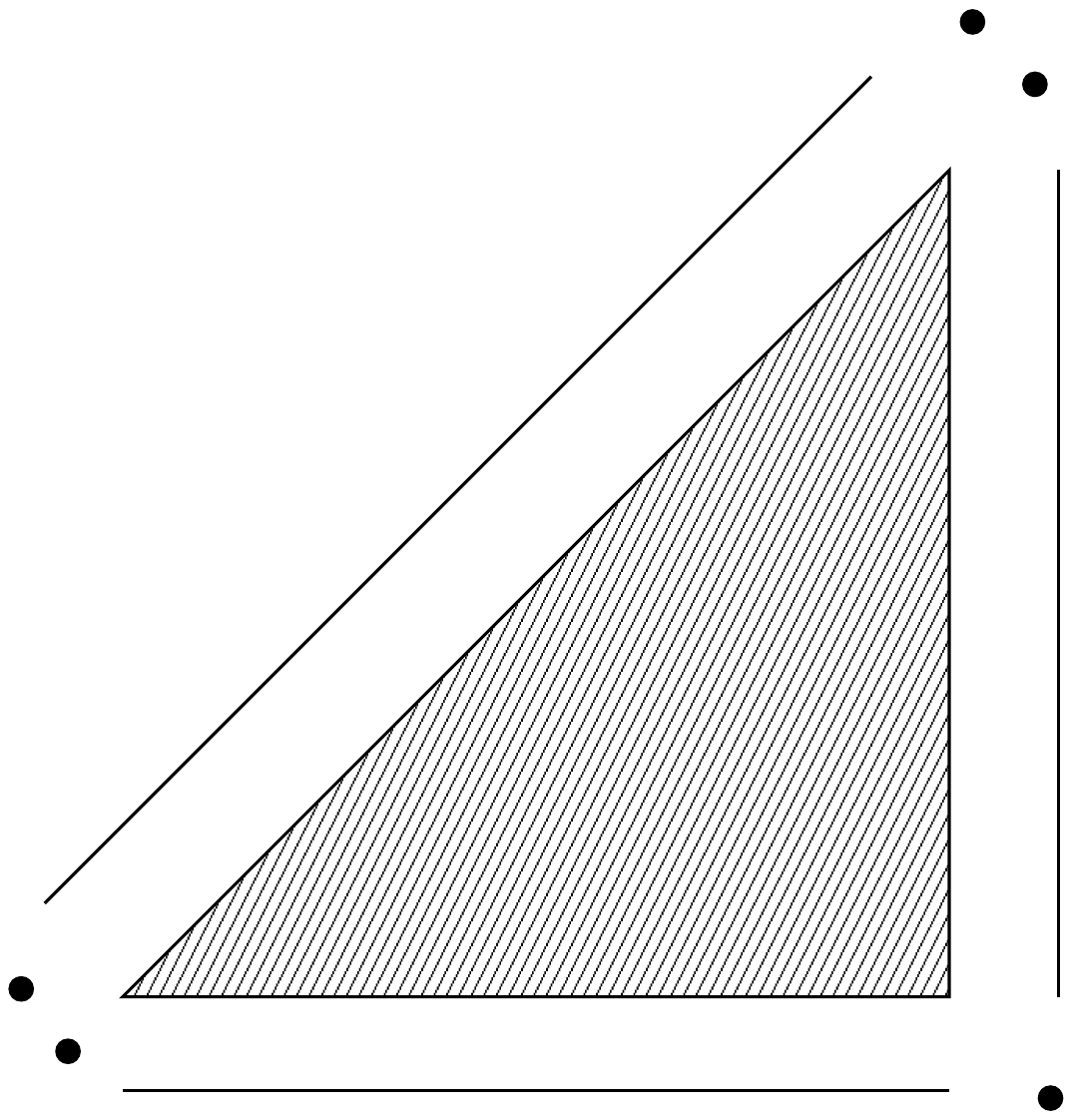}
\]
In the same way we have
\[
\raisebox{2cm}{Prim$\big( C_r^* (W) \big) \quad \cong \quad$}
\includegraphics[width=45mm,height=45mm]{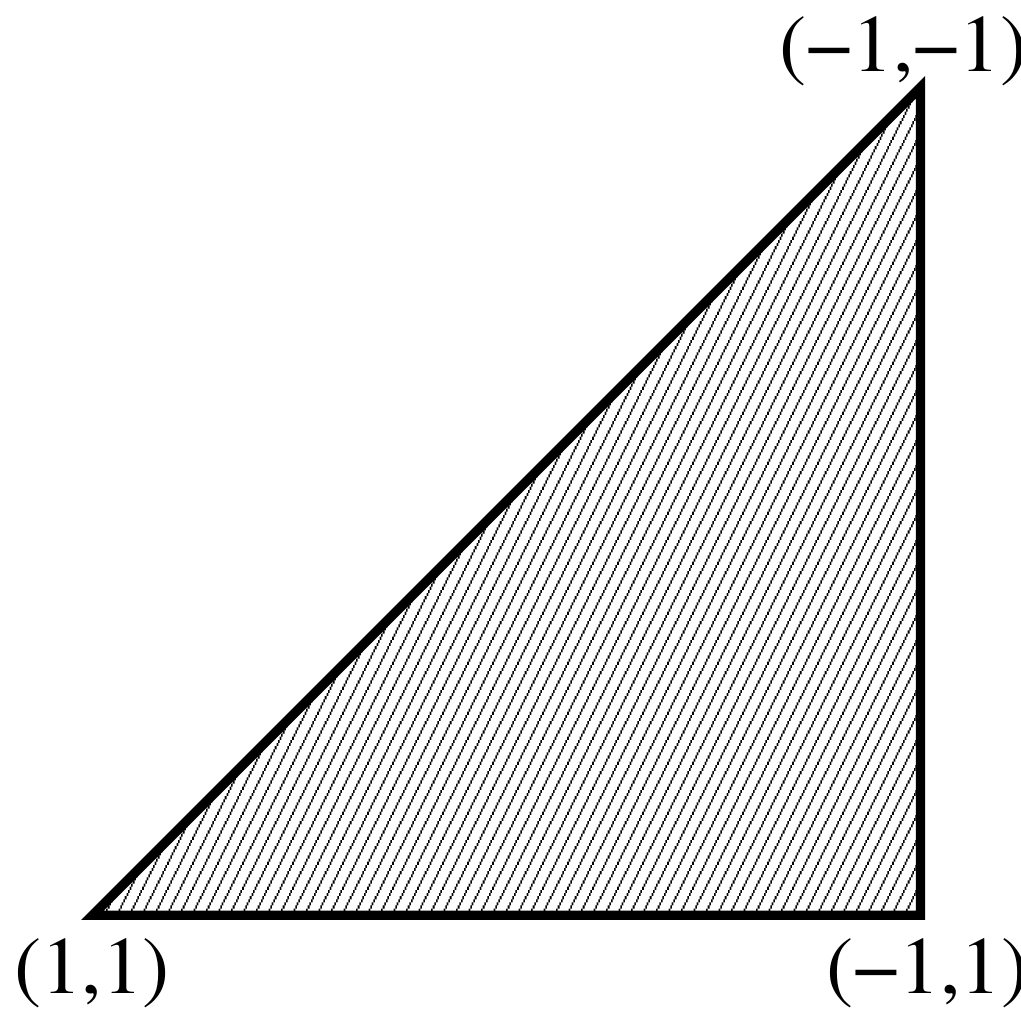}
\]
where the multiplicities of the non-Hausdorff points at the edge 
can be read off from the previous picture, by collapsing everything 
on the triangle. 

To find generating projections for the $K$-theory we first have 
a closer look at $W_0$. This group has four one-dimensional 
representations, let us call them $\ep_i$. These representations
and the corresponding projections are 
\begin{equation}\label{eq:6.62}
\begin{array}{lrrl}
\ep_i & \ep_i (s_1 ) & \ep (s_2 ) & p_i \\
\hline \ep_0 & 1 & 1 & \frac{1}{8} ( T_e + T_{s_1 s_2} + T_{s_2 s_1} 
+ T_{(s_1 s_2)^2} + T_{s_1} + T_{s_2} + T_{s_3} + T_{s_4} ) \\
\ep_1 & -1 & 1 & \frac{1}{8} ( T_e + T_{s_2} + T_{s_3} + 
T_{(s_1 s_2)^2} - T_{s_1} - T_{s_4} - T_{s_1 s_2} - T_{s_2 s_1} ) \\
\ep_2 & 1 & -1 &  \frac{1}{8} ( T_e + T_{s_1} + T_{s_4} + 
T_{(s_1 s_2)^2} - T_{s_2} - T_{s_3} - T_{s_1 s_2} - T_{s_2 s_1} ) \\
\ep_3 & -1 & -1 & \frac{1}{8} ( T_e + T_{s_1 s_2} + T_{s_2 s_1} +
T_{(s_1 s_2)^2} - T_{s_1} - T_{s_2} - T_{s_3} - T_{s_4} )
\end{array}
\end{equation}
The remaining irreducible $W_0$-representation $\rho_4$ has dimension
two, it defines $W_0$ as a reflection group. Let $p_4 \in \mh C [W_0]$
be rank one projector for that representation, i.e. 
\[
\rho_4 (p_4 ) = \begin{pmatrix} 1 & 0 \\ 0 & 0 \end{pmatrix} \quad 
\mr{and} \quad e_i (p_4 ) = 0 \quad i = 0,1,2,3 
\]
Note that $p_4$ has rank two in 
$\mr{End}_{\mh C} \big( \mh C [W_0 ] \big)$.

We also have to consider the stabilizer of $(-1,1) \in T$, the
subgroup\\ $\{ e ,s_1 ,s_4 ,(s_1 s_2 )^2 \} \cong D_2$. Let us list
its irreducible representations, which all have dimension one:
\begin{equation}\label{eq:6.60}
\begin{array}{lrrll}
\ep & \ep (s_1 ) & \ep (s_4 ) & \mr{Ind}_{\langle s_1 ,s_4 
\rangle}^{W_0} (\ep ) & p (\ep ) \\
\hline \ep_{++} & 1 & 1 & \ep_0 \oplus \ep_2 & p_0 \\
\ep_{+-} & 1 & -1 & \rho_4 & p_{+-} \\
\ep_{-+} & -1 & 1 & \rho_4 & p_{-+} \\
\ep_{--} & -1 & -1 & \ep_1 \oplus \ep_3 & p_3
\end{array}
\end{equation}
In the last column we indicate an element of $\mr{End}_{\mh C} \big( 
\mh C [W_0 ] \big)$ that acts as a rank one projector for that 
representation of $\langle s_1 ,s_4 \rangle$. For $\ep_{+-}$ and
$\ep_{-+}$ such an element cannot be found in $\mh C [W_0 ] \cong
\mr{End}_{W_0} \big( \mh C [W_0 ] \big)$, so we refrain from giving
an explicit formula.

Now we can indicate generators $p_0 ,\ldots ,p_8$ for $K_0 \big( C_r^* 
(W) \big)$. The last four classes of projections are defined by
their values in three special points.
\[
\begin{array}{llll}
p & p(1,1) & p(-1,-1) & p(-1,1) \\
\hline p_5 & p_0 & p_1 & p_{-+} \\
p_6 & p_2 & p_3 & p_{+-} \\
p_7 & p_0 + p_3 & p_4 & p_0 + p_3 \\
p_8 & p_0 + p_3 & p_4 & p_{+-} + p_{-+}
\end{array}
\]

\begin{itemize}
\item \textbf{\large generic case}
\item $P = \es$
\end{itemize}
\[
\begin{aligned}
& R_P = \es \qquad R_P^\vee = \es \\
& X^P = X \qquad X_P = 0 \qquad Y^P = Y \qquad Y_P = 0 \\
& T^P = T \qquad T_P = \{ 1 \} \qquad K_P = \{ 1 \} \\
& W^P = W(P,P) = \mc W_{PP} = W_0  \qquad W_P = \{ e \} \\
& \imath^o_{s_i} = (T_{s_i} (1 - \theta_{\alpha_i}) + 
(q(s_i ) - 1) \theta_{\alpha_i}) (q(s_i ) - \theta_{\alpha_i}
)^{-1} \qquad i = 1,2,3,4 \\
& \imath^o_{s_i}(t) = 1 \text{  if  } t \in T^{s_i}
\end{aligned}
\]
The $W_0$-stabilizer of any $t \in T$ is generated by reflections,
so all the unitary representations $\pi (\es 
,\delta_\es ,t)$ are irreducible.
\[
\begin{aligned}
& C_r^* (\mc R ,q)_P \cong C (T_u / W_0 ; M_8 (\mh C )) \\
& \raisebox{2cm}{Prim$\big( C_r^* (\mc R ,q)_P \big) \;=\;
T_u / W_0 \quad \cong \quad$}
\includegraphics[width=45mm,height=45mm]{b22}
\end{aligned}
\]
Since $T_u / W_0$ is contractible we conclude that
\[
\begin{array}{lll}
K_0 \big( C_r^* (\mc R ,q)_P \big) & \cong & \mh Z \\
K_1 \big( C_r^* (\mc R ,q)_P \big) & = & 0
\end{array}
\]
A generator is 
\[
p_\es = \frac{e_\es}{8} \sum_{w \in W_0} q(w) T_w
\]

\begin{itemize}
\item $P = \{ \alpha_1 \}$
\end{itemize}
\[
\begin{aligned}
& R_P = \{ \pm \alpha_1 \} \qquad 
  R_P^\vee = \{ \pm \alpha_1^\vee \} \\
& X^P = X / \mh Z (\alpha_1 / 2) \cong \mh Z \alpha_4 / 2 \quad
  X_P = X / \mh Z (\alpha_4 / 2) \cong \mh Z \alpha_1 / 2 \\
& Y^P = \mh Z \alpha_4^\vee \qquad Y_P = \mh Z \alpha_1^\vee \\
& T^P = \{ (1,t_2 ) : t_2 \in \mh C^\times \} \quad T_P = \{ (t_1 
  ,1) : t_1 \in \mh C^\times \} \quad K_P = \{ 1 \} \\
& W_P = \{ e ,s_1 \} \quad W^P = \{ e, s_2 ,s_4 ,s_1 s_3 \} \quad
  W(P,P) = \mc W_{PP} = \{ e ,s_4 \} \\
& \imath^o_{s_4} = (T_{s_4} (1 - \theta_{\alpha_4}) + (q_1 - 1)
  \theta_{\alpha_4}) (q_1 - \theta_{\alpha_4})^{-1}
\end{aligned}
\]
The algebra $\mc H_P$ is isomorphic to $\mc H (\mc R (A_1 )^\vee 
,q_1 )$. Hence it has two discrete series representations, with
central characters $\big( q_1^{\pm 1/2} ,1 \big)$ and 
$\big( -q_1^{\pm 1/2} ,1 \big)$. Clearly\\ 
$\imath^o_{s_4}(t) = 1$ if $t \in \big( T_u^P \big)^{s_4}$, so
\[
\begin{aligned}
& C_r^* (\mc R ,q)_P \cong C ( [0,1] ; M_4 (\mh C ) )^2 \\
& \raisebox{8mm}{Prim$\big( C_r^* (\mc R ,q)_P \big) \; \cong \;$}
\begin{picture}(2,2)
\linethickness{0.5mm}
\put(0.3,0.5){\line(1,0){1.0}}
\put(1.5,0.7){\line(0,1){1.0}}
\end{picture} \\
& K_0 \big( C_r^* (\mc R ,q)_P \big) \cong \mh Z^2 \\
& K_1 \big( C_r^* (\mc R ,q)_P \big) = 0
\end{aligned}
\]
Let us denote the canonical generators by $[p^+_{\alpha_1}]$ and
$[p^-_{\alpha_1}]$. Notice that as $\mc H (W_0 ,q)$-representations
the $\pi (P,\delta ,t)$ are deformations of
$\mr{Ind}_{W_P}^{W_0}(\mr{sign}_{W_P})$ if $q_1 > 1$ and of
$\mr{Ind}_{W_P}^{W_0}(\mr{triv}_{W_P})$ if $q_1 < 1$.
\\[2mm]
\begin{itemize}
\item $P = \{ \alpha_2 \}$
\end{itemize}
\[
\begin{aligned}
& R_P = \{ \pm \alpha_2 \} \qquad 
  R_P^\vee = \{ \pm \alpha_2^\vee \} \\
& X^P = X / \mh Z \alpha_2 \cong \mh Z \alpha_1 / 2 \qquad 
  X_P = X / \mh Z \alpha_3 \cong \mh Z \alpha_2 / 2 \\
& Y^P = \mh Z \alpha_3^\vee \qquad Y_P = \mh Z \alpha_2^\vee \\
& T^P = \{ (t_1 ,t_1 ) : t_1 \in \mh C^\times \} \qquad 
  T_P = \{ (t_2^{-1},t_2 ) : t_2 \in \mh C^\times \} \\
& K_P = \{ (1,1) ,(-1,-1) \} = \{ 1, k_P \} \\
& W_P = \{ e, s_2 \} \quad W^P = \{ e, s_1 ,s_4 ,s_3 s_1 \} \quad
  W(P,P) = \{ e ,s_3 \} \\
& \imath^o_{s_3} = (T_{s_3} (1 - \theta_{\alpha_3}) + (q_2 - 1)
  \theta_{\alpha_3}) (q_2 - \theta_{\alpha_3})^{-1} \\
& \imath^o_{s_3}(t) = 1 \text{  if  } t \in \big( T_u^P \big)^{s_3}
\end{aligned}
\]
The algebra $\mc H_P$ is isomorphic to $\mc H (\mc R (A_1 )^\vee 
,q_2 )$, so it has two discrete series representations 
$\delta_+$ and $\delta_-$. Their central characters are 
respectively $(q_2^{\pm 1/2},q_2^{\mp 1/2})$ and $(-q_2^{\pm1/2}
,-q_2^{\mp 1/2})$. However $\pi (k_P )$ intertwines $\pi (P,
\delta_+ ,(t_1 ,t_1 ))$ and\\ $\pi (P,\delta_- ,(-t_1 ,-t_1 ))$, so 
in the spectrum we get only one component $(P,\delta_+ ,T_u^P )$. 
The intertwiner $\pi (s_3 )$ acts as a reflection on $T_u^P$ and
$\pi (s_3 ,P,\delta_+ ,t) = 1$ whenever $s_3 (t) = t$. Therefore
\[
\begin{aligned}
& C_r^* (\mc R ,q)_P \cong C ([0,1] ; M_4 (\mh C )) \\
& \mr{Prim}\big( C_r^* (\mc R ,q)_P \big) \cong (P,\delta_+ 
  ,T_u^P ) / W(P,P) \cong [0,1] \\
& K_0 \big( C_r^* (\mc R ,q)_P \big) \cong \mh Z \\
& K_1 \big( C_r^* (\mc R ,q)_P \big) = 0
\end{aligned}
\]
There is a canonical generator $[p_{\alpha_2}] \in K_0 \big( C_r^*
(\mc R ,q)_P \big)$. The type of $\pi (P,\delta_+ ,t)$ as a
representation of $\mc H (W_0 ,q)$ is easily determined: for
$q_2 < 1$ it is a deformation of $\mr{Ind}_{W_P}^{W_0}(
\mr{triv}_{W_P})$ while for $q_2 > 1$ it is a deformation of
$\mr{Ind}_{W_P}^{W_0}(\mr{sign}_{W_P})$.

\begin{itemize}
\item $P = \{ \alpha_1 ,\alpha_2 \}$
\end{itemize}
\[
\begin{aligned}
& R_P = R_0 \qquad R_P^\vee = R_0^\vee \\
& X^P = 0 \qquad X_P = X \qquad Y^P = 0 \qquad Y_P = Y \\
& T^P = \{ 1 \} \qquad T_P = T \qquad K_P = \{ 1 \} \\
& W^P = W(P,P) = \mc W_{PP} = \{ e \} \qquad W_P = W_0 
\end{aligned}
\]
Generically all the residual points \eqref{eq:6.59} are in orbits
consisting of $|W_0 | = 8$ points. Proposition \ref{prop:3.30}.1 
tells us that every such $W_0$-orbit is the central character of
precisely one discrete series representation. The discrete series
representation $\delta_5$ with central character $W_0 (q_1^{-1/2},
-q_1^{-1/2})$ is the easiest to describe. It has dimension two,
and as a $W_0$-representation $\delta_5 \circ \phi$ is equivalent 
to $\ep_1 \oplus \ep_3$ (for $q_1 > 1$) or to $\ep_0 \oplus \ep_2$ 
(for $q_1 < 1$). For the other residual points we have to 
distinguish more relations between the parameters. Let $\delta_1 ,
\delta_2 ,\delta_3 ,\delta_4$ be the discrete series 
representations with respective central characters
\[
W_0 \big( q_1^{1/2},q_1^{-1/2}q_2 \big) \,, W_0
\big( q_1^{1/2},q_1^{1/2}q_2 \big) \,, W_0 
\big(-q_1^{1/2},-q_1^{-1/2}q_2 \big) \,, W_0 
\big( -q_1^{1/2},-q_1^{1/2}q_2 \big) 
\]
We list the type of $\delta_i \circ \phi_0$ as a representation of 
$W_0 \subset \mc S (W)$, for different $q$'s :
\begin{equation}\label{eq:6.61}
\begin{array}{c|cccc}
 & \delta_1 & \delta_2 & \delta_3 & \delta_4 \\
\hline 1 < q_1^{1/2} < q_2 < q_1 & 
\rho_4 & \ep_3 & \rho_4 & \ep_3 \\
q_1^{-1} < q_2 < q_1^{-1/2} < 1 & \ep_1 & \rho_4 & \ep_1 &
\rho_4 \\
1 < q_1^{-1/2} < q_2 < q_1^{-1} & \ep_2 & \rho_4 & \ep_2 &
\rho_4 \\
q_1 < q_2 < q_1^{1/2} < 1 & \rho_4 & \ep_0 & \rho_4 & 
\ep_0 \\
q_1^{-1/2} < q_2 < q_1^{1/2} > 1 & \ep_1 & \ep_3 & \ep_1 & \ep_3 \\
q_2^{-1} < q_1 < q_2 > 1 & \ep_2 & \ep_3 & \ep_2 & \ep_3 \\
1 > q_1^{1/2} < q_2 < q_1^{-1/2} & \ep_2 & \ep_0 & \ep_2 & \ep_0 \\
1 > q_2 < q_1 < q_2^{-1} & \ep_1 & \ep_0 & \ep_1 & \ep_0
\end{array}
\end{equation}
It can be shown by direct calculation that this table still gives 
all discrete series representations if either $q_1$ or $q_2$ (but
not both) equals 1. This not true for the special parameters
\eqref{eq:6.66} however. Nevertheless, for all the parameters 
under consideration here
\[
\begin{array}{lll}
K_0 \big( C_r^* (\mc R ,q)_P \big) & \cong & \mh Z^5 \\
K_1 \big( C_r^* (\mc R ,q)_P \big) & = & 0
\end{array}
\]
We denote the generating projections by $p(\delta_i ) ,\,
i = 1,2,3,4,5$. 

On the whole we found that $C_r^* (\mc R ,q)$ is Morita-equivalent
to the commutative $C^*$-algebra with spectrum
\[
\raisebox{2cm}{Prim$\big( C_r^* (\mc R ,q) \big) \quad \cong \quad$}
\includegraphics[width=45mm,height=45mm]{b21}
\]
and that 
\begin{equation}
\begin{array}{lll}
K_0 \big( C_r^* (\mc R ,q) \big) & \cong & \mh Z^9 \\
K_1 \big( C_r^* (\mc R ,q) \big) & = & 0
\end{array}
\end{equation}

\begin{itemize}
\item {\large $\mb{q_1 \neq 1 = q_2}$}
\item $P = \es$
\end{itemize}
\[
\begin{aligned}
& R_P = \es \qquad R_P^\vee = \es \\
& X^P = X \qquad X_P = 0 \qquad Y^P = Y \qquad Y_P = 0 \\
& T^P = T \qquad T_P = \{ 1 \} \qquad K_P = \{ 1 \} \\
& W^P = W(P,P) = \mc W_{PP} = W_0  \qquad W_P = \{ e \} \\
& \imath^o_{s_1} = (T_{s_1} (1 - \theta_{\alpha_1}) + (q_1 - 1)
\theta_{\alpha_1}) (q_1 - \theta_{\alpha_1})^{-1} \quad
\imath^o_{s_2} = T_{s_2}
\end{aligned}
\]
Since $\imath^o_{s_1}(t) = 1$ if $s_1 (t) = t$, all the nonscalar
selfintertwiners of unitary principal series representations come 
from $s_2$ or its conjugates. We see that $\pi (\es, \delta_\es ,t)$ 
is irreducible unless $t \in T_u^{s_2} \cup T_u^{s_3}$, in which 
case it is the direct sum of two inequivalent subrepresentations. 
Hence
\[
\begin{aligned}
& C_r^* (\mc R ,q)_P \; \cong \; \left\{ f \in C \big( T_u / W_0 
;M_8 (\mh C ) \big) : f \big( T_u^{s_2} \big) \in M_4 (\mh C ) 
\oplus M_4 (\mh C ) \right\} \\
& \raisebox{17mm}{Prim$\big( C_r^* (\mc R ,q)_P \big) 
\quad \cong \quad$}
\includegraphics[width=35mm,height=35mm]{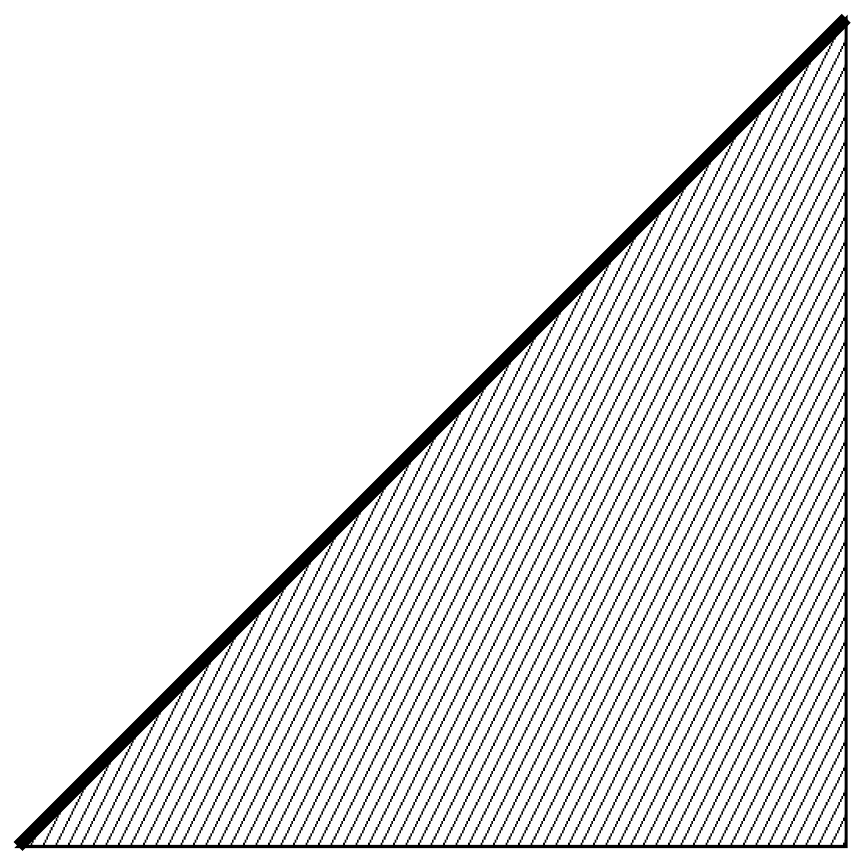}
\end{aligned}
\]
In this picture the diagonal edge should be regarded as 
consisting of double points. The algebra is diffeotopy equivalent 
to $M_4 (\mh C )^2$, so 
\[
\begin{array}{lll}
K_0 \big( C_r^* (\mc R ,q)_P \big) & \cong & \mh Z^2 \\
K_1 \big( C_r^* (\mc R ,q)_P \big) & = & 0
\end{array}
\]
Generators are for example
\[
p_0 := \frac{1}{8} e_\es \sum_{w \in W_0} q(w) T_w 
\text{  and  } p_3 :=
\frac{1}{8} e_\es \sum_{w \in W_0} (-1)^{\ell (w)} T_w
\]

\begin{itemize}
\item $P = \{ \alpha_1 \}$
\end{itemize}
This is identical to $P = \{ \alpha_1 \}$ in the generic case.

\begin{itemize}
\item $P = \{ \alpha_2 \}$
\end{itemize}
Here $\mc H_P \cong \mc H ( \mc R (A_1 )^\vee ,q_2 ) = 
\mh C [W (A_1 )]$. As we saw before, this algebra has no discrete
series representations, so there is no component in the spectrum
of $\mc S (\mc R ,q)$ corresponding to this $P$.

\begin{itemize}
\item $P = \{ \alpha_1 ,\alpha_2 \}$
\end{itemize}
\[
\begin{aligned}
& R_P = R_0 \qquad R_P^\vee = R_0^\vee \\
& X^P = 0 \qquad X_P = X \qquad Y^P = 0 \qquad Y_P = Y \\
& T^P = \{ 1 \} \qquad T_P = T \qquad K_P = \{ 1 \} \\
& W^P = W(P,P) = \mc W_{PP} = \{ e \} \qquad W_P = W_0 
\end{aligned}
\]
Some residual points confluence when $q_2 \to 1$, and only three
orbits remain. Two of those consist of four points, represented 
by $\big( q_1^{-1/2},q_1^{1/2} \big)$ and $\big( -q_1^{-1/2}
,-q_1^{-1/2} \big)$. The last orbit still contains 8 different 
points, for example $\big( -q_1^{-1/2},q_1^{1/2} \big)$. By
Proposition \ref{prop:3.30}.1 there is exactly one discrete 
series representation $\delta_5$ with central character
$W_0 \big( -q_1^{-1/2},q_1^{1/2} \big)$. It has dimension two and
restricted to $\mc H (W_0 ,q)$ it is a deformation of the 
$W_0$-representations $\ep_1 \oplus \ep_3$ or $\ep_0 \oplus \ep_2$,
depending on whether $q_1 > 1$ or $q_1 < 1$. We calculated the 
other discrete series representations already in \eqref{eq:6.61}.
For $q_1 > 1$ we have the onedimensional representations
\begin{equation}
\begin{array}{rrrr}
\delta & \delta (T_{s_1}) & \delta (T_{s_2}) & \delta (\theta_x ) \\
\hline
\delta_a & -1 & 1 & (q_1^{-1/2},q_1^{-1/2})(x) \\
\delta_b & -1 & -1 & (q_1^{-1/2},q_1^{-1/2})(x) \\
\delta_c & -1 & 1 & (-q_1^{-1/2},-q_1^{-1/2})(x) \\
\delta_d & -1 & -1 & (-q_1^{-1/2},-q_1^{-1/2})(x) \\
\end{array}
\end{equation}
On the other hand, for $q_1 < 1$ we can define onedimensional
representations by
\begin{equation}
\begin{array}{rrrr}
\delta & \delta (T_{s_1}) & \delta (T_{s_2}) & \delta (\theta_x ) \\
\hline
\delta_a & q_1 & 1 & (q_1^{1/2},q_1^{1/2})(x) \\
\delta_b & q_1 & -1 & (q_1^{1/2},q_1^{1/2})(x) \\
\delta_c & q_1 & 1 & (-q_1^{1/2},-q_1^{1/2})(x) \\
\delta_d & q_1 & -1 & (-q_1^{1/2},-q_1^{1/2})(x) 
\end{array}
\end{equation}
This leads to
\[
\begin{array}{lll}
C_r^* (\mc R ,q)_P & \cong & M_2 (\mh C ) \oplus \mh C^4 \\
K_0 \big( C_r^* (\mc R ,q)_P \big) & \cong & \mh Z^5 \\
K_1 \big( C_r^* (\mc R ,q)_P \big) & = & 0
\end{array}
\]
The generators of $K_0 \big( C_r^* (\mc R ,q)_P \big)$ 
corresponding to rank one projectors in these representations are
denoted by $p (\delta_v ) ,\, v \in \{ a,b,c,d,5 \}$.
\\[2mm]

Combining all $P$'s we find
\[
\begin{aligned}
& K_0 \big( C_r^* (\mc R ,q)_P \big) \cong \mh Z^5 \\
& K_1 \big( C_r^* (\mc R ,q)_P \big) = 0 \\
& \raisebox{2cm}{Prim$\big( C_r^* (\mc R ,q) \big) \quad \cong \quad$}
\includegraphics[width=43mm,height=43mm]{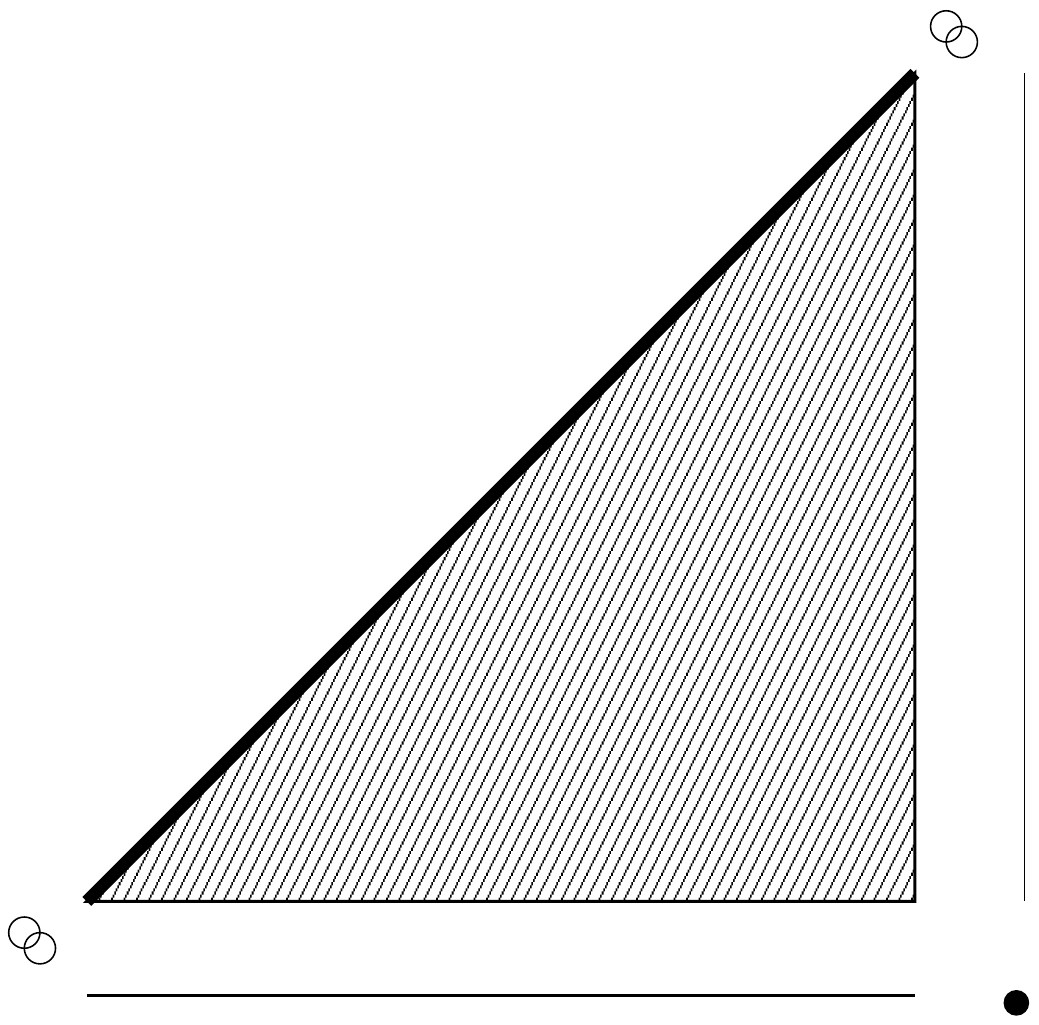}
\end{aligned}
\]

\begin{itemize}
\item {\large $\mb{q_1 = 1 \neq q_2}$}
\item $P = \es$
\end{itemize}
\[
\begin{aligned}
& R_P = \es \qquad R_P^\vee = \es \\
& X^P = X \qquad X_P = 0 \qquad Y^P = Y \qquad Y_P = 0 \\
& T^P = T \qquad T_P = \{ 1 \} \qquad K_P = \{ 1 \} \\
& W^P = W(P,P) = \mc W_{PP} = W_0  \qquad W_P = \{ e \} \\
& \imath^o_{s_1} = T_{s_1} \quad \imath^o_{s_2} = (T_{s_2} 
(1 - \theta_{\alpha_2}) + (q_2 - 1) \theta_{\alpha_2}) 
(q_2 - \theta_{\alpha_2})^{-1} 
\end{aligned}
\]
Since $\imath^o_{s_2}(t) = 1$ whenever $t \in T^{s_2}$, all the
nonscalar self-intertwiners of principal series representations 
come from $s_1 ,\, s_4 = s_2 s_1 s_2$ and $s_1 s_4$. This
implies that we should divide the points $t \in T_u$ in five 
classes.

\begin{enumerate}
\item $s_1 (t) \neq t \neq s_4 (t)$\\
Here we do not encounter nonscalar self-intertwiners, so
$\pi (\es ,\delta_\es ,t)$ is irreducible.

\item $s_1 (t) = t \neq s_4 (t)$\\
For such $t \; \pi (\es ,\delta_\es ,t)$ splits
into two summands, which as $\mc H (W_0 ,q)$-representations are
deformations of $\mr{Ind}_{\{ e,s_1 \}}^{W_0}(\mr{triv})$ and of
$\mr{Ind}_{\{ e,s_1 \}}^{W_0}(\mr{sign})$.

\item $s_1 (t) \neq t = s_4 (t)$\\
Just as in 2. $\!\pi (\es ,\delta_\es ,t)$ is the
direct sum of two irreducible subrepresentations, which can be
regarded as deformations of $\mr{Ind}_{\{ e,s_4 \}}^{W_0}
(\mr{triv})$ and of $\mr{Ind}_{\{ e,s_4 \}}^{W_0}(\mr{sign})$.

\item $(1,1)$ and $(-1,-1)$\\
These points are $W_0$-invariant, so $s_1$ and $s_4$ are 
conjugate in $W_{0,t}$. Hence $\pi (s_1, \es 
,\delta_\es ,t)$ is essentially the only independent 
intertwiner, and $\pi (\es ,\delta_\es ,t)$ can be
decomposed in only two irreducible parts. These are
deformations of the $W_0$-representations $\ep_1 \oplus \ep_3 
\oplus \rho_4$ and $\ep_0 \oplus \ep_2 \oplus \rho_4$.

\item $(1,-1)$ and $(-1,1)$\\
These points make up one $W_0$-orbit in $T_u$, their stabilizer
being $\{ e,s_1 ,s_4 ,s_1 s_4 \}$. Here the intertwiners 
$\pi (e), \pi (s_1 ), \pi (s_4 ) ,\pi (s_1 s_4 )$ are all
linearly independent, so $\pi (\es ,\delta_\es ,t)$
splits into no less than four irreducible summands, corresponding
to the irreducible representations of $\{ e,s_1 ,s_4 ,s_1 s_4 \}
\cong (\mh Z / 2 \mh Z )^2$.
\end{enumerate}

We visualize this as 
\[
\raisebox{2cm}{Prim$\big( C_r^* (\mc R ,q)_P \big) \; \cong \;$}
\includegraphics[width=35mm,height=35mm]{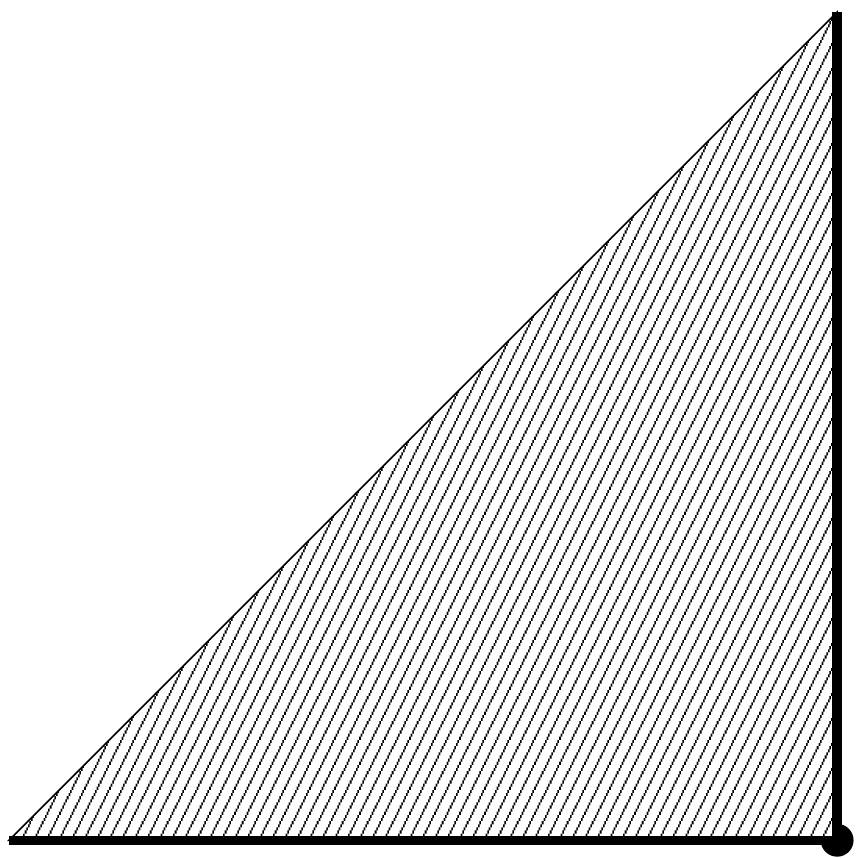}
\]
where the lower right corner depicts the fourfold non-Hausdorff
point $(-1,1)$. This algebra is diffeotopy-equivalent to its fiber 
over $(-1,1)$, so
\[
\begin{array}{lll}
K_0 \big( C_r^* (\mc R ,q)_P \big) & \cong & \mh Z^4 \\
K_1 \big( C_r^* (\mc R ,q)_P \big) & = & 0
\end{array}
\]
Generators are the rank one projectors 
\[
p_0 = \frac{1}{8} \sum_{w \in W_0} q(w) T_w \;,\;
p_3 = \frac{1}{8} \sum_{w \in W_0} (-1)^{\ell (w)} T_w
\;,\; p_{-+} \;,\; p_{+-}
\]
the last two being defined in the same way as the homonymous 
(classes of) projections in \eqref{eq:6.60}.

\begin{itemize}
\item $P = \{ \alpha_1 \}$
\end{itemize}
Here $\mc H_P \cong \mc H (\mc R (A_1 )^\vee ,q_1 ) \cong
\mh C [ W(A_1 )]$, so we do not find any discrete series
representations to induce.

\begin{itemize}
\item $P = \{ \alpha_2 \}$
\end{itemize}
This is identical to $P = \{ \alpha_2 \}$ for generic labels.

\begin{itemize}
\item $P = \{ \alpha_1 ,\alpha_2 \}$
\end{itemize}
\[
\begin{aligned}
& R_P = R_0 \qquad R_P^\vee = R_0^\vee \\
& X^P = 0 \qquad X_P = X \qquad Y^P = 0 \qquad Y_P = Y \\
& T^P = \{ 1 \} \qquad T_P = T \qquad K_P = \{ 1 \} \\
& W^P = W(P,P) = \mc W_{PP} = \{ e \} \qquad W_P = W_0 
\end{aligned}
\]
In the limit $q_1 \to 1$ many residual points confluence, and
some lose the residuality. There remain only two $W_0$-orbits of 
four points, from which we pick the representatives $(1,q_2 )$ 
and $(-1,-q_2 )$. The other discrete series representations were
already constructed in \eqref{eq:6.61}. They have dimension one, 
and for $q_2 > 1$:
\begin{equation}
\begin{array}{rrrr}
\delta & \delta (T_{s_1}) & \delta (T_{s_2}) & \delta (\theta_x ) \\
\hline
\delta_a & 1 & -1 & (1,q_2^{-1})(x) \\
\delta_b & -1 & -1 & (1,q_2^{-1})(x) \\
\delta_c & 1 & -1 & (-1,-q_2^{-1})(x) \\
\delta_d & -1 & -1 & (-1,-q_2^{-1})(x)
\end{array}
\end{equation}
On the other hand, for $q_2 < 1$ we find:
\begin{equation}
\begin{array}{rrrr}
\delta & \delta (T_{s_1}) & \delta (T_{s_2}) & \delta (\theta_x ) \\
\hline
\delta_a & 1 & q_2 & (1,q_2^{-1})(x) \\
\delta_b & -1 & q_2 & (1,q_2^{-1})(x) \\
\delta_c & 1 & q_2 & (-1,-q_2^{-1})(x) \\
\delta_d & -1 & q_2 & (-1,-q_2^{-1})(x)
\end{array}
\end{equation}
This summand of $C_r^* (\mc R ,q)$ is actually commutative:
\[
\begin{array}{lll}
C_r^* (\mc R ,q)_P & \cong & \mh C^4 \\
K_0 \big( C_r^* (\mc R ,q)_P \big) & \cong & \mh Z^4 \\
K_1 \big( C_r^* (\mc R ,q)_P \big) & = & 0
\end{array}
\]
Thus the spectrum of $C_r^* (\mc R ,q)$ is the non-Hausdorff space
\[
\raisebox{2cm}{Prim$\big( C_r^* (\mc R ,q) \big) \; \cong \;$}
\includegraphics[width=38mm,height=38mm]{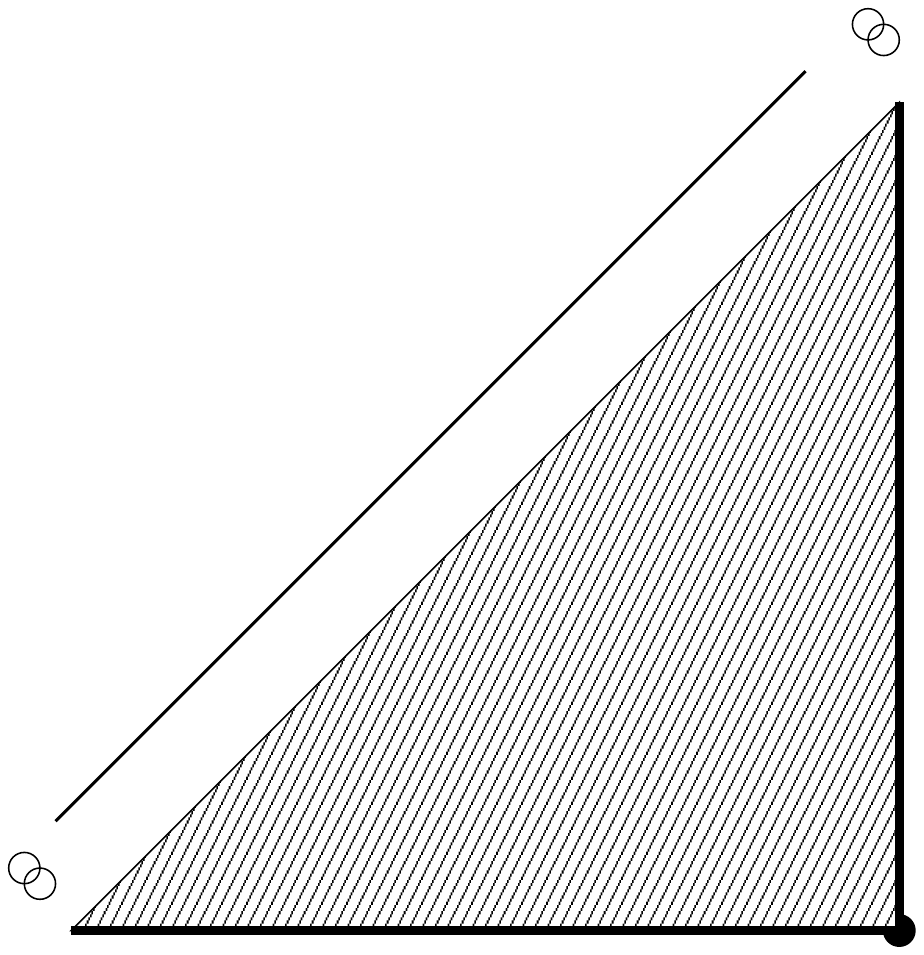}
\]
and its $K$-groups are
\begin{equation}
\begin{array}{lll}
K_0 \big( C_r^* (\mc R ,q) \big) & \cong & \mh Z^9 \\
K_1 \big( C_r^* (\mc R ,q) \big) & = & 0
\end{array}
\end{equation}

\newpage
\begin{itemize}
\item {\large $\mb{q_1 = q_2 \neq 1 ,\; q_1 = q_2^{-1} \neq 1 ,\; 
q_1 = q_2^2 \neq 1 ,\; q_1 = q_2^{-2} \neq 1}$}
\end{itemize}
These parameters are a bit too tricky to analyse with techniques
we used so far. To determine the spectra of the corresponding
$C^*$-algebras one can use the precise results on R-groups in 
\cite{Slo3}. It turns out that there are only three inequivalent 
discrete series representations. On the other hand, compared to 
the generic case two of the representations $\pi (P,\delta,t)$ 
with $|P| = 1$ become reducible. 
\\[5mm]
Let us make up the balance for the root datum $\mc R (B_2 )^\vee$.
The $K_1$-groups vanish for all label functions, and the 
$K_0$-groups are all free abelian of rank 9. We did not give all
the generating projections explicitly, but we have enough 
information to determine the maps $K_0 (\phi_0 )$. In the next table 
we list the images of the generators $p_i$ of $K_0 (\mc S (W))$. 
Assuming that all the calculations in this section are correct,
the table shows that Conjecture \ref{conj:5.28} is valid
for $\mc R (B_2 )^\vee$. 
\\[3mm]
We will not discuss the root datum $\mc R (B_2 )$ exhaustively, 
because there too are many different label functions. The group 
case is as in this section, and the equal label case is well 
understood, as described in Section \ref{sec:5.1}. We will study 
the generic label case in Section \ref{sec:6.7}. 
\newpage
\small
\[
\begin{array}{|c|c|c|}
\hline 1 < q_1^{-1/2} < q_2 < q_1^{-1} & q_1 = 1 < q_2 & 
1 < q_1^{1/2} < q_2 < q_1 \\
\hline p_\es + p_{\alpha_1}^+ + p_{\alpha_1}^- + p (\delta_5 )
& p_0 & p_\es \\
p_\es & p_3 & 
p_\es + p_{\alpha_1}^+ + p_{\alpha_1}^- + p (\delta_5 ) \\
\begin{array}{c}
p_\es + p_{\alpha_1}^+ + p_{\alpha_1}^- + p_{\alpha_2} + \\
p (\delta_1 ) + p (\delta_3 ) + p (\delta_5 ) 
\end{array} & 
p_0 + p_{\alpha_2} + p (\delta_a ) + p (\delta_c ) & 
p_\es + p_{\alpha_2} \\
\hline
p_\es + p_{\alpha_2} & p_3 + p_{\alpha_2} + p (\delta_b ) + 
p (\delta_d ) & 
\begin{array}{c}
p_\es + p_{\alpha_1}^+ + p_{\alpha_1}^- + p_{\alpha_2} \\
+ p (\delta_2 ) + p (\delta_4 ) + p (\delta_5 ) 
\end{array} \\
\begin{array}{c} 
2 p_\es + p_{\alpha_1}^+ + p_{\alpha_1}^- + \\
p (\delta_2 ) + p (\delta_4 ) 
\end{array} & 
p_{-+} + p_{+-} + p_{\alpha_2} & 
\begin{array}{c}
2 p_\es + p_{\alpha_1}^+ + p_{\alpha_1}^- + \\
p (\delta_1 ) + p (\delta_3 ) 
\end{array} \\
p_\es + p_{\alpha_1}^+ & p_{-+} & p_\es + p_{\alpha_1}^- \\
\hline
p_\es + p_{\alpha_1}^+ + p_{\alpha_2} + p(\delta_1 ) & 
p_{+-} + p_{\alpha_2} + p(\delta_a ) + p(\delta_d ) &
p_\es + p_{\alpha_1}^- + p_{\alpha_2} + p(\delta_4 ) \\
\begin{array}{c}
2 p_\es + p_{\alpha_1}^+ + p_{\alpha_1}^- + \\ 
p_{\alpha_2} + p (\delta_4 ) + p (\delta_5 ) 
\end{array} &
p_0 + p_3 + p_{\alpha_2} + p(\delta_b ) &
\begin{array}{c}
2 p_\es + p_{\alpha_1}^+ + p_{\alpha_1}^- + \\ 
p_{\alpha_2} + p_{\delta_2} + p (\delta_3 ) + p (\delta_5 ) 
\end{array} \\
\begin{array}{c}
2 p_\es + p_{\alpha_1}^+ + p_{\alpha_1}^- + \\
p_{\alpha_2} + p (\delta_4 ) 
\end{array} & 
p_{+-} + p_{-+} + p_{\alpha_2} + p(\delta_b ) &
\begin{array}{c}
2 p_\es + p_{\alpha_1}^+ + p_{\alpha_1}^- + \\ 
p_{\alpha_2} + p_{\delta_2} + p (\delta_3 ) 
\end{array} \\
% \hline \\
\hline q_1 < 1 = q_2 & q_1 = 1 = q_2 & q_1 > 1 = q_2 \\
\hline \begin{array}{c}
p_0 + p_{\alpha_1}^+ + p_{\alpha_1}^- + \\ 
p (\delta_a ) + p (\delta_c ) + p(\delta_5 ) 
\end{array} &
p_0 & p_0 \\
p_0 & p_1 & 
\begin{array}{c}
p_0 + p_{\alpha_1}^+ + p_{\alpha_1}^- + \\ 
p (\delta_a ) + p (\delta_c ) + p(\delta_5 ) 
\end{array} \\
\begin{array}{c}
p_3 + p_{\alpha_1}^+ + p_{\alpha_1}^- + \\ 
p (\delta_b ) + p (\delta_d ) + p(\delta_5 ) 
\end{array} &
p_2 & p_3 \\
\hline p_3 & p_3 & 
\begin{array}{c}
p_3 + p_{\alpha_1}^+ + p_{\alpha_1}^- + \\ 
p (\delta_b ) + p (\delta_d ) + p(\delta_5 ) 
\end{array} \\
p_0 + p_3 + p_{\alpha_1}^+ + p_{\alpha_1}^- & p_4 &
p_0 + p_3 + p_{\alpha_1}^+ + p_{\alpha_1}^- \\
p_0 + p_{\alpha_1}^+ + p(\delta_a ) & p_5 &
p_0 + p_{\alpha_1}^- + p(\delta_c ) \\
\hline p_3 + p_{\alpha_1}^+ + p(\delta_b ) & p_6 &
p_3 + p_{\alpha_1}^- + p(\delta_d ) \\
\begin{array}{c}
p_0 + p_3 + p_{\alpha_1}^+ + \\
p_{\alpha_1}^- + p(\delta_a ) + p(\delta_5 ) 
\end{array} & 
p_7 &
\begin{array}{c}
p_0 + p_3 + p_{\alpha_1}^+ + \\
p_{\alpha_1}^- + p(\delta_b ) + p(\delta_5 ) 
\end{array} \\
\begin{array}{c}
p_0 + p_3 + p_{\alpha_1}^+ + \\
p_{\alpha_1}^- + p(\delta_a ) 
\end{array} & 
p_8 & 
\begin{array}{c}
p_0 + p_3 + p_{\alpha_1}^+ + \\
p_{\alpha_1}^- + p(\delta_b ) 
\end{array} \\
% \hline \\
\hline q_1 < q_2 < q_1^{1/2} < 1 & q_1 = 1 > q_2 & 
q_1^{-1} < q_2 < q_1^{-1/2} < 1 \\
\hline \begin{array}{c}
p_\es + p_{\alpha_1}^+ + p_{\alpha_1}^- + p_{\alpha_2} \\
+ p (\delta_1 ) + p (\delta_4 ) + p (\delta_5 ) 
\end{array} &
p_0 + p_{\alpha_2} + p(\delta_a ) + p(\delta_c ) &
p_\es + p_{\alpha_2} \\
p_\es + p_{\alpha_2} & 
p_3 + p_{\alpha_2} + p(\delta_b ) + p(\delta_d ) & 
\begin{array}{c}
p_\es + p_{\alpha_1}^+ + p_{\alpha_1}^- + p_{\alpha_2} \\
+ p (\delta_1 ) + p (\delta_3 ) + p (\delta_5 ) 
\end{array} \\
p_\es + p_{\alpha_1}^+ + p_{\alpha_1}^- + p(\delta_5 ) & 
p_0 & p_\es \\
\hline p_\es & p_3 & 
p_\es + p_{\alpha_1}^+ + p_{\alpha_1}^- + p(\delta_5 ) \\
\begin{array}{c}
2 p_\es + p_{\alpha_1}^+ + p_{\alpha_1}^- + \\ 
p_{\alpha_2} + p (\delta_1 ) + p (\delta_3 ) 
\end{array} &
p_{+-} + p_{-+} & 
\begin{array}{c}
2 p_\es + p_{\alpha_1}^+ + p_{\alpha_1}^- + \\ 
p_{\alpha_2} + p (\delta_2 ) + p (\delta_4 ) 
\end{array} \\
p_\es + p_{\alpha_1}^+ + p_{\alpha_2} + p(\delta_1 ) &
p_{-+} + p_{\alpha_2} + p_(\delta_a ) + p(\delta_d ) & 
p_\es + p_{\alpha_1}^- + p_{\alpha_2} + p(\delta_3 ) \\
\hline p_\es + p_{\alpha_1}^+ & p_{+-} & p_\es + p_{\alpha_1}^- \\
\begin{array}{c}
2 p_\es + p_{\alpha_1}^+ + p_{\alpha_1}^- + p_{\alpha_2} \\
+ p_(\delta_1 ) + p (\delta_3 ) + p (\delta_5 ) 
\end{array} &
p_0 + p_3 + p_{\alpha_2} + p(\delta_a ) &
\begin{array}{c}
2 p_\es + p_{\alpha_1}^+ + p_{\alpha_1}^- + \\ 
p_{\alpha_2} + p (\delta_4 ) + p (\delta_5 ) 
\end{array} \\
\begin{array}{c}
2 p_\es + p_{\alpha_1}^+ + p_{\alpha_1}^- + \\
p_{\alpha_2} +  p_(\delta_1 ) + p (\delta_3 ) 
\end{array} &
p_{+-} + p_{-+} + p_{\alpha_2} + p(\delta_a ) &
\begin{array}{c}
2 p_\es + p_{\alpha_1}^+ + p_{\alpha_1}^- + \\
p_{\alpha_2} + p (\delta_4 ) 
\end{array} \\
\hline
\end{array}
\]
\normalsize

\section{$GL_n$}
\label{sec:6.5}

After the calculations with twodimensional root data we move on to
higher ranks. The easiest root data to study are those associated 
with the reductive group $GL_n$. The way right way to do this was 
indicated in \cite{BrPl2}. From \cite[Lemma 5.3]{Ply1} we know that 
the topological $K$-groups of these affine Hecke algebras are free 
abelian, and according to Theorem \ref{thm:5.4} they do not depend
on $q$. Because of the higher dimensionality we do not provide 
explicit generators for these $K$-groups anymore. Nevertheless, by
a different argument we show that Conjecture \ref{conj:5.28} holds
for these root data.

From now on many things will be parametrized by partitions 
and permutations, so let us agree on some notations. We write 
partitions in decreasing order and abbreviate $(x)^3 = (x,x,x)$. 
A typical partition looks like \index{mu@$\mu$}
\begin{equation}\label{eq:6.39}
\mu = (\mu_1 ,\mu_2 ,\ldots , \mu_d ) = 
(n)^{m_n} \cdots (2)^{m_2} (1)^{m_1}
\end{equation}
where some of the multiplicities $m_i$ may be 0. By 
$\mu \vdash n$ we mean that the weight of $\mu$ is 
\index{$b (\mu )$} \index{muvee@$\mu^\vee$}
\[
| \mu | = \mu_1 + \cdots + \mu_d = n
\]
The number of different $\mu_i$'s (i.e. the number of blocks in 
the diagram of $\mu$) will be denoted by $b (\mu )$ and the 
dual partition (obtained by reflecting the diagram of $\mu$) 
by $\mu^\vee$. Sometimes we abbreviate \index{gcd($\mu$)}
\begin{equation}
\begin{aligned}
& g = \mr{gcd}(\mu) = \mr{gcd}(\mu_1 ,\ldots ,\mu_d ) \\
& \mu ! = \mu_1 ! \mu_2 ! \cdots \mu_d !
\end{aligned}
\end{equation}
With a such partition $\mu$ of $n$ we associate the permutation 
\index{sigmamu@$\sigma (\mu )$}
\[
\sigma (\mu ) = (1 2 \cdots \mu_1 ) (\mu_1 + 1 \cdots \mu_1 + 
\mu_2 ) \cdots (n+1- \mu_d \cdots n) \in S_n
\]
As is well known, this gives a bijection between partitions of 
$n$ and conjugacy classes in the symmetric group $S_n$. The
centralizer $Z_{S_n}(\sigma (\mu ))$ is generated by the cycles
\[
((\mu_1 + \cdots + \mu_i + 1) (\mu_1 + \cdots + \mu_i + 2) 
  \cdots (\mu_1 + \cdots + \mu_i + \mu_{i+1}))
\]
and the ``permutations of cycles of equal length'', 
for example if $\mu_1 = \mu_2$:
\begin{equation}\label{eq:6.40}
(1 \, \mu_1 + 1) (2 \, \mu_1 + 2) \cdots (\mu_1 \, 2\mu_1)
\end{equation}
Using the second presentation of $\mu$ this means that
\[
Z_{S_n}(\sigma(\mu)) \cong \prod_{l=1}^n 
  (\mh Z / l \mh Z)^{m_l} \rtimes S_{m_l}
\]

Let us recall the definition of $\mc R (GL_n )$ :
\[
\begin{aligned}
& X = \mh Z^n \quad Q = \{ x \in X : x_1 + \cdots + x_n = 0 \} \\
& X^+ = \{ x \in \mh Z^n : x_1 \geq x_2 \geq \cdots \geq x_n \} \\
& Y = \mh Z^n \quad 
  Q^\vee = \{ y \in Y : y_1 + \cdots + y_n = 0 \} \\
& T = (\mh C^\times )^n \quad  t = (t(e_1),\ldots, t(e_n)) = 
  (t_1, \ldots, t_n) \\
& R_0 = R_1 = \{ e_i - e_j \in X : i \neq j \} \\
& R_0^\vee = R_1^\vee = \{ e_i - e_j \in Y : i \neq j \} \\
& F_{\mr{aff}} = \{ \alpha_i^\vee = e_i - e_{i+1} \} 
  \cup \{ 1 - \alpha_0^\vee = 1 - (e_1-e_n) \} \\
& s_i = s_{\alpha_i} \quad s_0 = t_{\alpha_0} s_{\alpha_0} = 
  t_{-\alpha_1} s_{\alpha_0} t_{\alpha_1} : x \to x + \alpha_0 - 
  \inp{\alpha_0^\vee}{x} \alpha_0 \\
& W_0 = \langle s_1 ,\cdots, s_{n-1} | s_i^2 = (s_i s_{i+1})^3 = 
  (s_i s_j)^2 = e : |i - j| > 1 \rangle \cong S_n \\
& S_{\mr{aff}} = \{ s_0 ,s_1 ,\ldots s_{n-1} \} \\
& W_{\mr{aff}} = \langle s_0, W_0 | s_0^2 = (s_0 s_i)^2 = 
  (s_0 s_1)^3 = (s_0 s_{n-1})^3 = e \;\mr{if}\; 2 \leq i 
  \leq n-2 \rangle \\
& W = W_{\mr{aff}} \rtimes \Omega \quad \Omega = \langle t_{e_1} 
  (1 \, 2 \cdots n) \rangle \cong \mh Z
\end{aligned}
\]
Because all roots of $R_0$ are conjugate, $s_0$ is conjugate to any
$s_i \in S_{\mr{aff}}$. Hence for any label function we have
\[
q(s_0 ) = q(s_i ) = q_{\alpha_i^\vee} := q
\]
The $W_0$-stabilizer of a point $\big( (t_1)^{\mu_1} (t_{\mu_1
+1})^{\mu_2} \cdots (t_n)^{\mu_d} \big) \in T$ is isomorphic to\\ 
$S_{\mu_1} \times S_{\mu_2} \times \cdots \times S_{\mu_d}$.

\begin{itemize}
\item {\large \textbf{group case} $\mb{q = 1}$}
\end{itemize}

By \eqref{eq:6.37} we have
\[
K_* \big( C_r^* (W) \big) \otimes \mh C \cong \check H^* \big( 
\widetilde{T_u} \big/ S_n ; \mh C \big) \cong 
\bigoplus_{\mu \vdash n} \check H^* \big( T_u^{\sigma (\mu)} 
\big/ Z_{S_n}(\sigma (\mu )) ; \mh C \big)
\]
Therefore we want to determine $T_u^{\sigma (\mu)} / 
Z_{S_n}(\sigma (\mu ))$. If $\mu$ is as in \eqref{eq:6.39} then
\begin{equation}
\begin{aligned}
& T^{\sigma (\mu )} = \{ (t_1)^{\mu_1} (t_{\mu_1+1})^{\mu_2} 
\cdots (t_n)^{\mu_d} \in T \} \\
& T^{\sigma (\mu)} \big/ Z_{S_n}(\sigma (\mu )) \cong
  (\mh C^\times)^{m_n} / S_{m_n} \times \cdots \times 
  (\mh C^\times)^{m_1} / S_{m_1}
\end{aligned}
\end{equation}
where $S_{m_l}$ acts on $(\mh C^\times)^{m_l}$ by permuting
the coordinates. To handle this space we use the following nice,
elementary result, a proof of which can be found for example in 
\cite[p. 97]{BrPl2}.

\begin{lem}\label{lem:6.1}
For any $m \in \mh N$ there is an isomorphism of algebraic
varieties
\[
(\mh C^\times )^m \big/ S_m \cong \mh C^{m-1} \times \mh C^\times
\]
\end{lem}

It follows that $T_u^{\sigma (\mu)} \big/ Z_{S_n}(\sigma (\mu ))$
has the homotopy type of $\mh T^{b (\mu )}$. The latter space
has torsion free cohomology, so by \eqref{eq:6.38}
\begin{equation}
K_* \big( C_r^* (W) \big) \cong \check H^* \big( \widetilde{T_u} 
/ S_n ; \mh Z \big) \cong \bigoplus_{\mu \vdash n}
\check H^* \big( \mh T^{b (\mu )} ; \mh Z \big) \cong
\bigoplus_{\mu \vdash n} \mh Z^{2^{b (\mu )}}
\end{equation}

\begin{itemize}
\item {\large \textbf{generic, equal label case} $\mb q \neq 1$}
\item $P_\mu = F_0 \setminus \{\alpha_{\mu_1}, \alpha_{\mu_1 + 
\mu_2}, \ldots, \alpha_{n - \mu_d} \}$
\end{itemize}

Inequivalent subsets of $F_0$ are parametrized by partitions 
$\mu$ of $n$. For the typical partition \eqref{eq:6.39} we put 
\[
\begin{aligned}
& R_{P_\mu} \cong (A_{n-1})^{m_n} \times \cdots \times (A_1)^{m_2} 
  \cong R_{P_\mu}^\vee \\
& X^{P_\mu} \cong \mh Z (e_1 + \cdots + e_{\mu_1)} / \mu_1 + 
  \cdots + \mh Z (e_{n+1 - \mu_d} + \cdots + e_n) / \mu_d \\
& X_{P_\mu} \cong \big( \mh Z^n / \mh Z (e_1 + \cdots + e_n )
  \big)^{m_n} \times  
  \cdots \times  \big( \mh Z^2 / \mh Z (e_1 + e_2 ) \big)^{m_2} \\
& Y^{P_\mu} = \mh Z (e_1 + \cdots + e_{\mu_1}) + \cdots + \mh Z 
  (e_{n+1 - \mu_d} + \cdots + e_n) \\
& Y_{P_\mu} = \{ y \in \mh Z^n : y_1 + \cdots + y_{\mu_1} = \cdots 
  = y_{n+1- \mu_d} + \cdots + y_n = 0 \} \\
& T^{P_\mu} = \{ (t_1)^{\mu_1} \cdots (t_n)^{\mu_d} \in T \} \\
& T_{P_\mu} = \{ t \in T : t_1 t_2 \cdots t_{\mu_1} = \cdots =
  t_{n+1 -\mu_d} \cdots t_n = 1 \} \\
& K_{P_\mu} = \{ t \in T^{P_\mu} : t_1^{\mu_1} = \cdots =  
  t_n^{\mu_d} = 1 \} \\
& W_{P_\mu} \cong (S_n )^{m_n} \times \cdots \times (S_2 )^{m_2}
  \qquad W(P_\mu ,P_\mu ) \cong S_{m_n} \times \cdots \times 
  S_{m_2} \times S_{m_1} \\
& \mc W_{P_\mu P_\mu} = K_{P_\mu} \rtimes W(P_\mu ,P_\mu ) \qquad
  Z_{S_n}(\sigma (\mu )) = W(P_\mu ,P_\mu ) \ltimes \prod_{l=1}^n 
  (\mh Z / l \mh Z)^{m_l}
\end{aligned}
\] 
The $W_{P_\mu}$-orbits of residual points for $\mc H_{P_\mu}$ 
are parametrized by 
\[
K_{P_\mu} \big( ( q^{(\mu_1 - 1)/2} ,q^{(\mu_1 - 3)/2}, \ldots ,
q^{(1 - \mu_1 )/2}) \cdots (q^{(\mu_d - 1)/2} ,q^{(\mu_d - 3)/2},
\ldots ,q^{(1 - \mu_d )/2}) \big)
\]
This set is obviously in bijection with $K_{P_\mu}$, and indeed
the intertwiners $\pi (k), k \in K_{P_\mu}$ act on it by 
multiplication. Together with Proposition \ref{prop:3.30}.2 this
implies
\[
\begin{aligned}
& \bigcup_{\Delta_{P_\mu}} \big( P_\mu, \delta ,T^{P_\mu} \big) 
\big/ K_{P_\mu} \cong T^{P_\mu} \\
& \bigcup_{\Delta_{P_\mu}} \big( P_\mu, \delta ,T^{P_\mu} \big) 
\big/ \mc W_{P_\mu P_\mu} \cong T^{P_\mu} \big/ W(P_\mu ,P_\mu ) =
T^{\sigma (\mu )} \big/ Z_{S_n}(\sigma (\mu ))
\end{aligned}
\]
It a point $t \in T^{P_\mu}$ has a nontrivial stabilizer in 
$W(P_\mu P_\mu) \cong \prod_{l=1}^n S_{m_l}$, then this isotropy
group is generated by transpositions. From \eqref{eq:6.40} we
see that every such transposition $w \in W_0$ can be written as
a product of mutually commuting reflections $s_\alpha$ with
$c_\alpha^{-1}(t) = 0$. By \eqref{eq:3.52} this gives 
$\imath^o_w (t) = 1$, and since dim $\delta = 1$ also
\begin{equation}\label{eq:6.48}
\pi (w,P_\mu ,\delta ,t) = 1 \quad \mr{if} \quad w (t) = t
\end{equation}
So the action of $\mc W_{P_\mu P_\mu}$ on 
\[
C \Big( \bigsqcup_{\Delta_{P_\mu}} T_u^{P_\mu} ; M_{n! / \mu !}
(\mh C ) \Big)
\]
is essentially only on $\bigsqcup_{\Delta_{P_\mu}} T_u^{P_\mu}$ 
and the conjugation part doesn't really matter. In particular we
deduce that 
\begin{equation}\label{eq:6.41}
C_r^* (\mc R ,q) \cong \bigoplus_{\mu \vdash n} M_{n! / \mu !} 
\Big( C \Big( \bigsqcup_{\Delta_{P_\mu}} T_u^{P_\mu} \Big) \Big)
\cong \bigoplus_{\mu \vdash n} M_{n! / \mu !} \big(
T_u^{\sigma (\mu )} \big/ Z_{S_n}(\sigma (\mu )) \big)
\end{equation}
Similar results were obtained by completely different methods in
\cite{Mis}. Just as in the group case it follows that 
\begin{equation}
K_* \big( C_r^* (\mc R ,q) \big) \cong \bigoplus_{\mu \vdash n}
K^* \big( T_u^{\sigma (\mu )} \big/ Z_{S_n}(\sigma (\mu )) \big)
\cong \bigoplus_{\mu \vdash n} K^* \big( \mh T^{b(\mu )} \big) 
\cong \bigoplus_{\mu \vdash n} \mh Z^{2^{b(\mu )}}
\end{equation}

As promised, we show that Conjecture \ref{conj:5.28} holds 
in this particular case.

\begin{thm}\label{thm:6.2}
\[
K_* (\phi_0 ) : K_* \big( C_r^* (\mc R (GL_n ) ,q^0) \big) \to
K_* \big( C_r^* (\mc R (GL_n ) ,q) \big)
\]
is an isomorphism.
\end{thm}
\emph{Proof.}
We assume that $q>1$. If instead $q<1$ then we only have to 
modify our argument by replacing the sign representations 
everywhere by trivial representations. 

Let $J_i^c \subset C_r^* (\mc R ,q)$ be the norm completion of
the ideal $J_i \subset \mc S (\mc R ,q)$ from \eqref{eq:3.53}.
For every $i$ there is a unique partition $\mu_i$ such that 
\[
J_{i-1}^c / J_i^c \cong C ( T_u^{P_i} ; \mr{End} V_i )^{\mc W_i}
\cong M_{n! / \mu_i !} C \big( T_u^{\sigma (\mu_i )} \big/ 
Z_{S_n}(\sigma (\mu_i )) \big)
\]
The induced homomorphism
\[
f_i : \phi_0^{-1}(J_{i-1}^c ) / \phi_0^{-1}(J_i^c ) 
\to J_{i-1} / J_i
\]
is injective. If we would know that $K_* (f_i )$ is an 
isomorphism for every $i$, then the theorem would follow with 
Lemma \ref{lem:2.12}. 

We know already that the number of components of 
Prim$(C_R^* (\mc R ,q))$ equals the number of components of 
Prim$(C_r^* (W)) \cong \widetilde{T_u} / W_0$. Therefore it 
suffices to construct, for every $i$, a projection 
$p \in \phi_0^{-1} (J_{i-1}^c ) / \phi_0^{-1}(J_i^c )$ 
such that\\ $f_i (p) \in J_{i-1}^c / J_i^c$ has rank one. 

We can do this because we know precisely what all discrete 
series look like. We may assume that the central character of
$(\delta_i ,V_i )$ lies in $T_{rs}$, so that the central 
character of $\pi (P_i, \delta_i ,t) \circ \phi_0$ is $W_0 t$.
As a $W_0$-representation $\pi (P_i, \delta_i ,t) \circ \phi_0$ 
is equivalent with $\mr{Ind}_{W_{P_i}}^{W_0} (\ep_{W_{P_i}})$,
where $\ep_{W_P}$ denotes the sign representation of $W_P$. 
The ideal $\phi_0^{-1}(J_{i-1}^c )$ does not annihilate this
representation, but it is contained in the kernel of the
representation 
\[
\pi (P,\delta ,t') \text{ for} \quad  P \supset P_i = P_{\mu_i} 
\text{ and} \quad t' \in T_u^P \subset 
T^{P_i} = T^{\sigma (\mu_i )} .
\]
As $W_0$-representation we have $\pi (P,\delta ,t') \circ \phi_0
\cong \mr{Ind}_{W_P}^{W_0} (\ep_{W_P})$. Hence the "stalk" of
$\phi_0^{-1}(J_{i-1}^c ) / \phi_0^{-1}(J_i^c )$ at $W_0 t$ 
contains
\[
\bigcap_{P \supsetneq P_i ,t \in T_u^P} 
\ker \mr{Ind}_{W_P}^{W_0} (\ep_{W_P})
\big/ \ker \mr{Ind}_{W_{P_i}}^{W_0} (\ep_{W_{P_i}})
\]
which is a subquotient algebra of $\mh C [W_0 ]$. Therefore we
can find a suitable $p$ already in $\mh C [W_0 ]$: pick a 
projection of minimal rank in $\mh C [W_0 ]$ which is in
\[
\bigcap_{P \supsetneq P_i} \ker \mr{Ind}_{W_P}^{W_0} (\ep_{W_P})
\text{ , but not in} \quad 
\ker \mr{Ind}_{W_{P_i}}^{W_0} (\ep_{W_{P_i}}) .
\] 
Then $\phi_0 (p)$ will act as a rank one projector on 
$\pi (P_i ,\delta_i, t). \qquad \Box$ 
\\[1mm]

\section{$A_{n-1}$}
\label{sec:6.6}

It is known from Theorem \ref{thm:5.4} that the periodic 
cyclic homology of type $A$ affine Hecke algebras does not depend
on $q$, but it will still be insightful to determine the spectra
of these algebras. The title of this section is $A_{n-1}$ instead
of $A_n$ because we want to consider everything as a quotient or
a subset of $\mh Z^n$. If we require that our root datum is
semisimple, then the easiest case is when $X$ is the weight 
lattice. This is completely analogous to the $GL_n$-case, we 
can even show in the same way that Conjecture \ref{conj:5.28} 
holds. The calculations are also manageable when $X$ is the root
lattice. Intermediate lattices however would require quite some 
extra bookkeeping, so we do not study those. Throughout this 
section we assume that $n > 2$, because the root system $A_1$ has
slightly different properties. 

The root datum $\mc R (A_{n-1})^\vee$ is defined as
\[
\begin{aligned}
& X = \mh Z^n / \mh Z (e_1 + \cdots e_n) \cong 
  Q + ((e_1 + \cdots + e_n) /n - e_n) \\
& Q = \{ x \in \mh Z^n : x_1 + \cdots + x_n = 0 \} \quad 
  X^+ = \{ x \in X : x_1 \geq x_2 \geq \cdots \geq x_n \} \\
& Y = Q^\vee = \{ y \in \mh Z^n : y_1 + \cdots + y_n = 0 \} \\ 
& T = \{ t \in (\mh C^\times)^n : t_1 \cdots t_n = 1 \} \quad
  t = (t(e_1),\ldots, t(e_n)) = (t_1, \ldots, t_n) \\
& R_0 = R_1 = \{ e_i - e_j \in X : i \neq j \} \\
& R_0^\vee = R_1^\vee = \{ e_i - e_j \in Y : i \neq j \} \\
& F_0 = \{ \alpha_i = e_i - e_{i+1} \} 
  \quad \alpha_0 = e_1 - e_n \\
& s_i = s_{\alpha_i} \quad s_0 = t_{\alpha_0} s_{\alpha_0} = 
  t_{-\alpha_1} s_{\alpha_0} t_{\alpha_1} : x \to x + \alpha_0 - 
  \inp{\alpha_0^\vee}{x} \alpha_0 \\
& W_0 = \langle s_1 ,\cdots , s_{n-1} | s_i^2 = (s_i s_{i+1})^3 = 
  (s_i s_j)^2 = e \; \mr{if} \; \mr |i - j| > 1 \rangle \cong S_n \\
& S_{\mr{aff}} = \{ s_0 ,s_1 ,\ldots ,s_{n-1} \} \\
& W_{\mr{aff}} = \langle s_0, W_0 | s_0^2 = (s_0 s_i)^2 = 
  (s_0 s_1)^3 = (s_0 s_{n-1})^3 = e \;\mr{if}\; 2 \leq i 
  \leq n-2 \rangle \\
& W = W_{\mr{aff}} \rtimes \Omega \quad 
  \Omega = \langle t_{e_1 - (e_1 + \cdots e_n)/n} (1 2 \cdots n)
  \rangle \cong \mh Z / n \mh Z 
\end{aligned}
\]
Because all roots are conjugate, $s_0$ is conjugate to any 
$s_i \in S_{\mr{aff}}$, and for any label function
\[
q(s_0 ) = q(s_i ) = q_{\alpha_i^\vee} = q
\]
The $W_0$-stabilizer of $\big( (t_1)^{\mu_1} (t_{\mu_1+1})^{\mu_2}
\cdots (t_n)^{\mu_d} \big)$ is isomorphic to
$S_{\mu_1} \times \cdots \times S_{\mu_d}$. Generically there are
$n! n$ residual points, and they all satisfy $t(\alpha_i ) = q$
or $t(\alpha_i ) = q^{-1}$ for $1 \leq i < n$. There residual
points form $n$ conjugacy classes unless $q = 1$, in which case
$T$ itself is the only residual coset. 

\begin{itemize}
\item {\large \textbf{group case} $\mb{q = 1}$}
\end{itemize}

According to \eqref{eq:6.37} we have 
\[
K_* \big( C_r^* (W) \big) \otimes \mh C \cong \check H^* \big( 
\widetilde{T_u} \big/ S_n ; \mh C \big) \cong \bigoplus_{\mu \vdash n}
\check H^* \big( T_u^{\sigma (\mu)} \big/ Z_{S_n}(\sigma (\mu )) ; 
\mh C \big)
\]
Pick a partition $\mu$ of $n$ and write it as in \eqref{eq:6.39}.
\[
\begin{array}{lll}
T^{\sigma (\mu )} & = & \{ (t_1)^{\mu_1} (t_{\mu_1+1})^{\mu_2} 
  \cdots (t_n)^{\mu_d} \in T \} \\
& \cong & \{ (t_1)^{\mu_1} (t_{\mu_1+1})^{\mu_2} \cdots 
  (t_n)^{\mu_d} \in (\mh C^\times )^n \} / \mh C^\times \times 
  \{ (e^{2 \pi i k / n})^n : 0 \leq k < g \} \\
\multicolumn{3}{l}{T^{\sigma (\mu)} \big/ Z_{S_n}(\sigma (\mu )) 
\quad \cong \quad
\left( (\mh C^\times)^{m_n} / S_{m_n} \times \cdots \times (
\mh C^\times)^{m_1} / S_{m_1} \right) \big/ \mh C^\times \times }\\
 & & \hspace{2cm} \{ (e^{2 \pi i k / n})^n : 0 \leq k < g \}
\end{array}
\]
where $\mh C^\times $ acts diagonally. By Lemma \ref{lem:6.1}
each factor $(\mh C^\times)^{m_i} / S_{m_i}$ is homotopy 
equivalent to a circle. The induced action of $S^1 \subset 
\mh C^\times$ on this direct product of circles identifies with 
a direct product of rotations. Hence $T^{\sigma (\mu)} /
Z_{S_n}(\sigma (\mu ))$ is homotopy equivalent with 
$\mh T^{b(\mu ) - 1} \times \{\mr{gcd}(\mu ) \; \mr{points} \}$.
The cohomology of this space has no torsion, so by \eqref{eq:6.38}
\begin{equation}\label{eq:6.46}
\begin{aligned}
& K_* \big( C_r^* (W) \big) \cong \check H^* \big( \widetilde{T_u} 
\big/ S_n ; \mh Z \big) \cong \mh Z^{d(n)} \\
& d(n) = \sum_{\mu \vdash n} \mr{gcd}(\mu) 2^{b(\mu )-1}
\end{aligned}
\end{equation}

\begin{itemize}
\item {\large \textbf{generic, equal label case} $\mb q \neq 1$}
\item $P_\mu = F_0 \setminus \{\alpha_{\mu_1}, \alpha_{\mu_1 + 
\mu_2}, \ldots, \alpha_{n - \mu_d} \}$
\end{itemize}

Inequivalent subsets of $F_0$ are parametrized by partitions 
$\mu$ of $n$. For the typical partition \eqref{eq:6.39} we put 
\[
\begin{aligned}
& R_{P_\mu} \cong (A_{n-1})^{m_n} \times \cdots \times (A_1)^{m_2}
  \cong R_{P_\mu}^\vee \\
& X^{P_\mu} \!\cong \left( \mh Z (e_1 + \!\cdots\! + e_{\mu_1)} / 
  \mu_1 + \!\cdots\! + \mh Z (e_{n+1 - \mu_d} + \!\cdots\! + e_n) 
  / \mu_d \right) \!\big/ \mh Z (e_1 + \!\cdots\! + e_n) / g \\
& X_{P_\mu} \cong \big( \mh Z^n / \mh Z (e_1 + \cdots + e_n )
  \big)^{m_n} \times
  \cdots \times \big( \mh Z^2 / \mh Z (e_1 + e_2 ) \big)^{m_2} \\
& Y^{P_\mu} = \{ y \in \mh Z (e_1 + \cdots + e_{\mu_1}) + \cdots 
  + \mh Z (e_{n+1 - \mu_d} + \cdots + e_n) : y_1 + \cdots + y_n = 
  0 \} \\
& Y_{P_\mu} = \{ y \in Y : y_1 + \cdots + y_{\mu_1} = \cdots = 
  y_{n+1- \mu_d} + \cdots + y_n = 0 \} \\
& T^{P_\mu} = \{ (t_1)^{\mu_1} \cdots (t_n)^{\mu_d} \in T :
  t_1^{\mu_1 / g} \cdots t_n^{\mu_d / g} = 1 \} \\
& T_{P_\mu} = \{ t \in T : t_1 t_2 \cdots t_{\mu_1} = \cdots =
  t_{n+1 -\mu_d} \cdots t_n = 1 \} \\
& K_{P_\mu} = \{ t \in T^{P_\mu} : t_1^{\mu_1} = \cdots =  
  t_n^{\mu_d} = 1 \} \\
& W_{P_\mu} \cong (S_n )^{m_n} \times \cdots \times (S_2 )^{m_2}
  \qquad W(P_\mu ,P_\mu ) \cong S_{m_n} \times \cdots \times 
  S_{m_2} \times S_{m_1} \\
& \mc W_{P_\mu P_\mu} = K_{P_\mu} \rtimes W(P_\mu ,P_\mu ) \qquad
  Z_{S_n}(\sigma (\mu )) = W(P_\mu ,P_\mu ) \ltimes \prod_{l=1}^n 
  (\mh Z / l \mh Z)^{m_l}
\end{aligned}
\]
The $W_{P_\mu}$-orbits of residual points for $\mc H_{P_\mu}$ are
represented by the points
\begin{multline}\label{eq:6.50}
\big( ( q^{(\mu_1 - 1)/2} ,q^{(\mu_1 - 3)/2}, \ldots ,
q^{(1 - \mu_1 )/2}) \cdots (q^{(\mu_d - 1)/2} ,q^{(\mu_d - 3)/2},
\ldots ,q^{(1 - \mu_d )/2}) \big) \, \cdot \\ 
\big( (e^{2 \pi i k_1 / \mu_1} )^{\mu_1} \cdots (e^{2 \pi i k_d /
\mu_d} )^{\mu_d} \big) \quad , \quad 0 \leq k_i < \mu_i
\end{multline}
These points are in bijection with $K_{P_\mu} \times \mh Z / 
\mr{gcd}(\mu ) \mh Z$. Also $T^{\sigma (\mu )}$ consists of
exactly gcg$(\mu )$ components, one of which is $T^{P_\mu}$.
Together with Proposition \ref{prop:3.30}.2 this leads to
\[
\begin{aligned}
& \bigcup_{\Delta_{P_\mu}} \big( P_\mu, \delta ,T^{P_\mu} \big) 
\big/ K_{P_\mu} \cong T^{P_\mu} \times \mh Z / \mr{gcd}(\mu ) \mh Z 
\cong T^{\sigma (\mu )} \\
& \bigcup_{\Delta_{P_\mu}} \big( P_\mu, \delta ,T^{P_\mu} \big) 
\big/ \mc W_{P_\mu P_\mu} \cong T^{\sigma (\mu )} \big/ 
Z_{S_n}(\sigma (\mu ))
\end{aligned}
\]
If a point $t \in T^{P_\mu}$ has a nontrivial stabilizer in
$W(P_\mu ,P_\mu ) \cong \prod_{l=1}^n S_{m_l}$ then this 
stabilizer is generated by transpositions. From \eqref{eq:6.40}
we see that every such transposition $w \in W(P_\mu ,P_\mu )$ 
can be written as a product of mutually commuting reflections
$s_\alpha$ with $\alpha \in R_0$ and $c_\alpha^{-1}(t) = 0$. 
So by \eqref{eq:3.52} $\imath^o_w (t) = 1$ and since the
discrete series representation $\delta$ is onedimensional,
also $\pi (w,P_\mu ,\delta ,t) = 1$. On the 
other hand, $K_{P_\mu}$ permutes the components of 
$\bigcup_{\Delta_{P_\mu}} \big( P_\mu, \delta ,T^{P_\mu} \big)$ 
faithfully, so the action of $\mc W_{P_\mu P_\mu}$ on 
\newpage
\[
C \bigg( \bigcup_{\Delta_{P_\mu}} T^{P_\mu} ; M_{n! / \mu !} 
(\mh C ) \bigg)
\]
is essentially only on the underlying space. Therefore
\begin{equation}
\begin{aligned}
C_r^* (\mc R ,q) & \cong \bigoplus_{\mu \vdash n} M_{n! / \mu !} 
\Big( C \Big( \bigsqcup_{\Delta_{P_\mu}} T_u^{P_\mu} \Big) \Big)
\cong \bigoplus_{\mu \vdash n} M_{n! / \mu !} \big(
T_u^{\sigma (\mu )} \big/ Z_{S_n}(\sigma (\mu )) \big) \\
K_* \big( C_r^* (\mc R ,q) \big) & \cong \bigoplus_{\mu \vdash n}
K^* \big( T_u^{\sigma (\mu )} \big/ Z_{S_n}(\sigma (\mu )) \big)
\cong \bigoplus_{\mu \vdash n} K^* \big( \mh T^{b(\mu ) -1} 
\big) \cong \mh Z^d
\end{aligned}
\end{equation}
where $d = \sum_{\mu \vdash n} \mr{gcd}(\mu) 2^{b(\mu )-1}$.

We conclude that for $\mc R (A_{n-1})^\vee$ the $K$-theory
of $C_r^* (\mc R ,q)$ does not depend on $q$, and is a free
abelian group.

\begin{thm}\label{thm:6.3}
\[
K_* (\phi_0 ) : 
K_* \big( C_r^* \big( \mc R (A_{n-1} )^\vee ,q^0 \big) \big) 
\to K_* \big( C_r^* \big( \mc R (A_{n-1} )^\vee ,q \big) \big)
\]
is an isomorphism.
\end{thm}
\emph{Proof.}
This is completely analogous to Theorem \ref{thm:6.2}. The 
essential common properties of these two root data are that the
reduced $C^*$-algebras are Morita equivalent to commutative
$C^*$-algebras, and that we have explicit descriptions of the 
discrete series representations of the algebras 
$\mc H_P. \qquad \Box$
\\[3mm]

The analogy with $\mc R (GL_n )$ is significantly weaker for
the root datum $\mc R (A_{n-1})$:
\[
\begin{aligned}
& X = Q = \{ x \in \mh Z^n : x_1 + \cdots + x_n = 0 \} \quad
  X^+ = \{ x \in X : x_1 \geq x_2 \geq \cdots \geq x_n \} \\
& Q^\vee = \{ y \in \mh Z^n : y_1 + \cdots + y_n = 0 \} \\
& Y = \mh Z^n / \mh Z (e_1 + \cdots +e_n ) \cong 
  Q^\vee + ((e_1 + \cdots + e_n)/n - e_1) \\
& T = (\mh C^\times)^n / \mh C^\times \quad t = (t_1, \ldots, t_n) 
  = (t(e_1), \ldots, t(e_n)) \\
& R_0 = R_1 = \{ e_i - e_j \in X : i \neq j\} \\
& R_0^\vee = R_1^\vee =  \{ e_i - e_j \in Y : i \neq j\} \\
& F_0 = \{ \alpha_i = e_i - e_{i+1} : 1 \leq i < n \} 
  \quad \alpha_0 = e_1 - e_n \\
& s_i = s_{\alpha_i} \quad s_0 = t_{\alpha_0} s_{\alpha_0} =  
  t_{-\alpha_1} s_{\alpha_0} t_{\alpha_1} : 
  x \to x + \alpha_0 - \inp{\alpha_0^\vee}{x} \alpha_0 \\
& W_0 = \langle s_1 ,\cdots , s_{n-1} | s_i^2 = (s_i s_{i+1})^3 = 
  (s_i s_j)^2 = e \; \mr{if} \; \mr |i - j| > 1 \rangle \cong S_n \\
& S_{\mr{aff}} = \{ s_0 ,s_1 ,\ldots ,s_{n-1} \} \quad
  \Omega = \{ e \} \\
& W = W_{\mr{aff}} = \langle s_0, W_0 | s_0^2 = (s_0 s_i)^2 = 
  (s_0 s_1)^3 = (s_0 s_{n-1})^3 = e \;\mr{if}\; 2 \leq i 
  \leq n-2 \rangle \\
& q(s_0 ) = q(s_i ) = q_{\alpha_i^\vee} := q
\end{aligned}
\]
For $q \neq 1$ there are $n!$ residual points. They form one 
$W_0$-orbit, and a typical residual point is
\newpage
\[
\big( q^{(1-n)/2}, q^{(3-n)/2}, \ldots ,q^{(n-1)/2} \big)
\]
To determine the isotropy group of points of $T$ we have to be
careful. In general the $W_0$-stabilizer of 
\[
\big( (t_1)^{\mu_1} (t_{\mu_1 + 1})^{\mu_2} \cdots (t_n 
)^{\mu_d} \big) \in T
\] 
is isomorphic to 
\[
S_{\mu_1} \times S_{\mu_2} \times \cdots \times S_{\mu_d}
\subset W_0
\]
However, in some special cases the diagonal action of $\mh 
C^\times$ on $(\mh C^\times )^n$ gives rise to extra stabilizers.
Let $r$ be a divisor of $n ,\, k \in (\mh Z / r \mh Z)^\times$
and $\lambda = (\lambda_1 ,\ldots ,\lambda_l )$ a partition of
$n/r$. The isotropy group of
\begin{multline}\label{eq:6.42}
\big( (t_1)^{\lambda_1} (e^{2 \pi i k / r} t_1 )^{\lambda_1} \cdots
(e^{-2 \pi i k / r} t_1 )^{\lambda_1} (t_{r \lambda_1 + 1}
)^{\lambda_2} \cdots \\ (e^{-2 \pi i k / r} t_{r \lambda_1 + 1} 
)^{\lambda_2} \cdots (e^{-2 \pi i k / r} t_n )^{\lambda_l} \big)
\end{multline}
is isomorphic to 
\begin{equation}\label{eq:6.43}
S^r_{\lambda_1} \times S^r_{\lambda_2} \times \cdots \times
S^r_{\lambda_l} \rtimes \mh Z / r \mh Z
\end{equation}
Explicitly the subgroup $\mh Z / r \mh Z$ is generated by
\begin{multline}\label{eq:6.44}
(1 \; \lambda_1 + 1 \; 2 \lambda_1 + 1 \cdots (r-1)\lambda_1 + 1)
(2 \; \lambda_1 + 2 \; 2 \lambda_1 + 2 \cdots (r-1)\lambda_1 + 2)
\cdots (\lambda_1 \; 2 \lambda_1 \cdots r \lambda_1 ) \\
\cdots (n+1 - r \lambda_d \; n + 1 + (1-r) \lambda_d \cdots 
n+1+ (r-1) \lambda_d ) ( n+ (1-r) \lambda_d \, n +(2-r) 
\lambda_d \cdots n)
\end{multline}
and it acts on every factor $S^r_{\lambda_j}$ in \eqref{eq:6.43}
by cyclic permutations. 

\begin{itemize}
\item {\large \textbf{group case} $\mb{q = 1}$}
\end{itemize}
As we saw before
\[
K_* \big( C_r^* (W) \big) \otimes \mh C \cong \check H^* \big( 
\widetilde{T_u} \big/ S_n ; \mh C \big) \cong \bigoplus_{\mu \vdash n}
\check H^* \big( T_u^{\sigma (\mu)} \big/ Z_{S_n}(\sigma (\mu )) ; 
\mh C \big)
\]
For the typical partition $\mu$ we have
\begin{equation}
T^{\sigma (\mu )} = \{ (t_1)^{\mu_1} (t_{\mu_1+1})^{\mu_2} 
\cdots (t_n)^{\mu_d} \} \big/ \mh C^\times \times 
\{ t : t(e_j) = e^{2 \pi i j k / g} ,\, 0 \leq k < g \}
\end{equation}
which is the disjoint union of $g = \mr{gcd}(\mu )$ complex tori
of dimension\\ $m_n + m_{n-1} + \cdots + m_1 - 1$.
\begin{multline}\label{eq:6.45}
T^{\sigma (\mu )} \big/ Z_{S_n}(\sigma (\mu )) \cong
\left( (\mh C^\times)^{m_n} / S_{m_n} \times \cdots \times 
(\mh C^\times)^{m_1} / S_{m_1} \right) \big/ \mh C^\times \times \\
\{ t : t(e_j) = e^{2 \pi i j k / g} ,\, 0 \leq k < g \}
\end{multline}
Curiously enough these sets are diffeomorphic to the corresponding
sets for $\mc R (A_{n-1} )^\vee$, a phenomenon for which the
author does not have a good explanation. Anyway, we do take
advantage of this by reusing our deduction that \eqref{eq:6.45}
is homotopy equivalent with $\mh T^{b (\mu ) - 1} \times \{
\mr{gcd}(\mu ) \; \mr{points} \}$. Just as in \eqref{eq:6.46} we
conclude that 
\begin{equation}\label{eq:6.47}
\begin{aligned}
& K_* \big( C_r^* (W) \big) \cong \check H^* \big( \widetilde{T_u}
\big/ S_n ; \mh Z \big)\cong \mh Z^{d(n)} \\
& d(n) = \sum_{\mu \vdash n} \mr{gcd}(\mu) 2^{b(\mu )-1}
\end{aligned}
\end{equation}

\begin{itemize}
\item {\large \textbf{generic case} $\mb{q \neq 1}$}
\end{itemize}
This is noticeably different from the generic cases for
$\mc R (GL_n )$ and $\mc R (A_{n-1}^\vee )$ because
$C_r^* (\mc R (A_{n-1} ,q) )$ is not Morita equivalent to a
commutative $C^*$-algebra. Of course the inequivalent subsets
of $F_0$ are still parametrized by partitions $\mu$ of $n$.

\begin{itemize}
\item $P_\mu = F_0 \setminus \{\alpha_{\mu_1}, \alpha_{\mu_1 + 
\mu_2}, \ldots, \alpha_{n - \mu_d} \}$
\end{itemize}
\[
\begin{aligned}
& R_{P_\mu} \cong (A_{n-1})^{m_n} \times \cdots \times 
  (A_1)^{m_2} \cong R_{P_\mu}^\vee \\
& X^{P_\mu} \cong \{ x \in \mh Z (e_1 + \cdots + e_{\mu_1)} / \mu_1 
  + \cdots + \mh Z (e_{n+1 - \mu_d} + \cdots + e_n) / \mu_d : \\
& \quad x_1 + \cdots + x_n = 0 \} \\ 
& X_{P_\mu} \cong \{ x \in \mh Z^{\mu_1} / \mh Z (e_1 + \cdots + 
  e_{\mu_1} ) + \cdots + \mh Z^{\mu_d} / \mh Z (e_{n+1-\mu_d} + 
  \cdots + e_n ) : \\
& \quad x_1 + \cdots + x_n \in g \mh Z / g \mh Z \} \\
& Y^{P_\mu} \cong \mh Z (e_1 + \cdots + e_{\mu_1}) + \cdots + \mh Z 
  (e_{n+1 - \mu_d} + \cdots + e_n) / \mh Z (e_1 + \cdots + e_n) \\
& Y_{P_\mu} \cong \{ y : y_1 + \cdots + y_{\mu_1} = \cdots = y_{n+
  1- \mu_d} + \cdots + y_n = 0 \}  / \mh Z (e_1 + \cdots e_n) \\
& T^{P_\mu} = \{ (t_1)^{\mu_1} \cdots (t_n)^{\mu_d} \} /
  \mh C^\times \\
& T_{P_\mu} = \{ t : t_1 t_2 \cdots t_{\mu_1} = \cdots = t_{n+1
  -\mu_d} \cdots t_n = 1 \} / \{ z \in \mh C : z^g = 1 \} \\
& K_{P_\mu} = \{ (t_1)^{\mu_1} \cdots (t_n)^{\mu_d} : t_1^{\mu_1}
  = \cdots = t_n^{\mu_d} = 1 \} / \{ z \in \mh C : z^g = 1 \} \\
& W_{P_\mu} \cong S_n^{m_n} \times S_{n-1}^{m_{n-1}} \times \cdots
  \times S_2^{m_2} \quad W(P_\mu ,P_\mu ) \cong
  S_{m_n} \times \cdots \times S_{m_2} \times S_{m_1} 
\end{aligned}
\]
Note that 
\[
T^{\sigma (\mu )} = T^{P_\mu} \times \{ t : t(e_j) = 
e^{2 \pi i j k / g} ,\, 0 \leq k < g \}
\]
The $W_{P_\mu}$-orbits of residual points for $\mc H_{P_\mu}$ are
represented by the points of
\[
K_{P_\mu} \big( q^{(\mu_1 -1)/2} ,q^{(\mu_1 -3)/2}, \ldots
q^{(1-\mu_1 )/2} ,q^{(\mu_2 -1)/2} ,\ldots ,q^{(\mu_d -1)/2} ,
\ldots ,q^{(1-\mu_d)/2} \big)
\]
Hence the intertwiners $\pi (k)$ with $k \in K_{P_\mu}$ permute
the elements of $\Delta_{P_\mu}$ faithfully, and
\[
\bigsqcup_{\Delta_{P_\mu}} \big( P_\mu ,\delta ,T^{P_\mu} \big) / 
K_{P_\mu} \cong T^{P_\mu} = \big( T^{\sigma(\mu)} \big)_1
\]
where $( \cdot )_1$ means the connected component containing 
$(1,1, \ldots, 1) = 1 \in T$. In \eqref{eq:6.48} we saw that the
intertwiners for $\mc R (GL_n) ,q \neq 1$ have the property
\[
w(t) = t \Rightarrow \pi (w, P_\mu ,\delta ,t) = 1
\]
This implies that in our present setting we can have $w(t) = t$
and $\pi (w, P_\mu ,\delta ,t) \neq 1$ only if $w(t) = t$ does
not hold without taking the action of $\mh C^\times$ into account.
Let us classify such $w \in W(P_\mu ,P_\mu )$ and $t \in
T^{P_\mu}$ up to conjugacy. For a divisor $r$ of $g^\vee := 
\mr{gcd}(\mu^\vee )$ we have the partition \index{mu1r@$\mu^{1/r}$}
\[
\mu^{1/r} := (n r)^{m_n / r} \cdots (2 r)^{m_2 /r} (r)^{m_1 /r}
\]
Notice that
\[
b (\mu^{1/r} ) = b (\mu ) = b (\mu^\vee )
\]
There exists a $\sigma \in S_n$ which is conjugate to $\sigma
(\mu^{1/r})$ and satisfies $\sigma^r = \sigma (\mu )$. We 
construct a particular such $\sigma$ as follows. If $r = g^\vee$
then (starting from the left) replace every block 
\[
(d+1 \; d+2 \cdots d + m )(d+1+m \cdots d + 2 m)
\cdots (d+(g^\vee -1)m \cdots d + g^\vee m )
\]
of $\sigma (\mu )$ by
\[
(d+1 \; d+1+m \cdots d+1+(g^\vee - 1)m \; 2 \; d+2+m \cdots
d+2+(g^\vee -1)m \; d+3 \cdots d+g^\vee m)
\]
We denote the resulting element by $\sigma (\mu )^{1 / g^\vee}$
and for general $r | g^\vee$ we define
\[
\sigma (\mu )^{1/r} := 
\big( \sigma (\mu )^{1 / g^\vee} \big)^{g^\vee /r}
\]
Consider the cosets of subtori
\[
T^{P_\mu}_{r,k} := \big( T^{\sigma (\mu )^{1/r}} \big)_1 
\big( (1)^{g^\vee \mu_1 /r} (e^{2 \pi i k /r} )^{g^\vee 
\mu_{1 + g^\vee /r} /r} \cdots (e^{-2 \pi i k /r} )^{g^\vee \mu_d
/r} \big) \quad k \in \mh Z
\]
If gcd$(k,r)$ = 1 then the generic points of $T^{P_\mu}_{r,k}$
have $W(P_\mu ,P_\mu )$-stabilizer 
\[
\inp{W_{P_\mu}}{\sigma (\mu )^{1/r}} \cap W(P_\mu ,P_\mu ) \cong
\mh Z / r \mh Z
\]
Note that for $r' | g^\vee$
\begin{equation}\label{eq:6.49}
T^{P_\mu}_{r',k} \subset T^{P_\mu}_{r,k} \qquad \mr{if}\; r | r'
\end{equation}
If a point $t \in T^{P_\mu}_{r,k}$ does not lie on any 
$T^{P_\mu}_{r',k'}$ with $r' > r$, then its $W(P_\mu ,P_\mu 
)$-stabilizer may still be larger than $\mh Z / r \mh Z$. 
However, it is always of the form
\[
S^r_{\lambda_1} \times \cdots \times S^r_{\lambda_l} \rtimes 
\mh Z / r \mh Z
\]
By an argument like on pages \pageref{eq:6.48} and 
\pageref{eq:6.50} one can show that the intertwiners 
$\pi (w, P_\mu ,\delta ,t)$ are scalar for 
$w \in S^r_{\lambda_1} \times \cdots \times S^r_{\lambda_l}$ 
and nonscalar for $w \in (\mh Z / r \mh Z) \setminus \{ e \}$. 
Because $\mh Z / r \mh Z$ is cyclic this implies that 
$\pi (P_\mu ,\delta ,t)$ is the direct sum of exactly $r$
inequivalent irreducible representations. For a more systematic
discussion of such matters we refer to \cite{DeOp2}. 

Different choices of $\sigma (\mu )^{1/r}$ or of $k \in (\mh Z / 
r \mh Z )^\times$ lead to conjugate subvarieties of $T^{P_\mu}$, 
so we have a complete discription of Prim$\big(C_r^* (\mc R ,q
)_{P_\mu}\big)$. To calculate the $K$-theory of this algebra we 
use Theorem \ref{thm:2.15}, which says that (modulo torsion) it 
is isomorphic to 
\[
H^*_{W(P_\mu ,P_\mu )} \big( T_u^{P_\mu};\mc L_u \big) \cong
\check H^* \big( T^{P_\mu} \big/ W(P_\mu ,P_\mu ) ;\mc L_u^{W(
P_\mu ,P_\mu )} \big)
\]
We know from \cite{Ill} that we can endow $T_u^{P_\mu}$ with 
the structure of a finite $W(P_\mu ,P_\mu )$-CW-complex, such
that every $T_{u,r,k}^{P_\mu}$ is a subcomplex. The local 
coefficient system $\mc L_u$ is not very complicated: 
$\mc L_u (B) \cong \mh Z^r$ if and only if $B \setminus \partial B$
consists of generic points in a conjugate of $T_{u,r,k}^{P_\mu}$.
In suitable coordinates the maps $\mc L_u (B \to B')$ are all
of the form 
\[
\mh Z^r \to \mh Z^{r/d} : (x_1 ,\ldots ,x_r ) \to 
(x_1 + x_2 + \cdots + x_d ,\ldots ,x_{1+r-d} + \cdots + x_r )
\]
Hence the associated sheaf is the direct sum of several subsheaves
$\mf F^\mu_r$, one for each divisor $r$ of gcd$(\mu^\vee )$. The
support of $\mf F^\mu_r$ is
\[
W(P_\mu ,P_\mu ) T_{u,r,1}^{P_\mu} \big/ W(P_\mu ,P_\mu ) \cong
T_u^{P_{\mu^{1/r}}} \big/ Z_{S_n} (\sigma (\mu^{1/r}))
\]
and on that space it has constant stalk $\mh Z^{\phi (r)}$.
Here $\phi$ is the Euler $\phi$-function, i.e.
\[
\phi (r) = \# \{ m \in \mh Z : 0 \leq m < r : \mr{gcd}(m,r) = 1 \}
= \# (\mh Z / r \mh Z )^\times 
\]
This the rank of $\mf F^\mu_r$ because in every point of
$T_{u,r,1}$ we have $r$ irreducible representations, but the ones
corresponding to numbers that are not coprime with $r$ are already
accounted for by the sheaves $\mf F^\mu_{r'}$ with $r' | r$.
Now we can calculate
\begin{equation}
\begin{aligned}
\check H^* \big( T^{P_\mu} / W(P_\mu ,P_\mu ) ;\mc L_u^{W(
P_\mu ,P_\mu )} \big) & 
\cong \bigoplus_{r | \mr{gcd}(\mu^\vee )} \check H^* \big( 
T^{P_\mu} / W(P_\mu ,P_\mu ) ; \mf F^\mu_r \big) \\
& \cong \bigoplus_{r | \mr{gcd}(\mu^\vee )} \check H^* \big(
T_u^{P_{\mu^{1/r}}} \big/ Z_{S_n} (\sigma (\mu^{1/r})) ; 
\mh Z^{\phi (r)} \big) \\
& \cong \bigoplus_{r | \mr{gcd}(\mu^\vee )} \check H^* \big( 
\mh T^{b(\mu^{1/r}) - 1} ; \mh Z^{\phi (r)} \big) \\
& \cong \bigoplus_{r | \mr{gcd}(\mu^\vee )} \mh Z^{\phi (r)
2^{b(\mu^{1/r}) - 1}} \\
& = \bigoplus_{r | \mr{gcd}(\mu^\vee )} \mh Z^{\phi (r) 
2^{b(\mu^\vee ) -1}} \: = \; 
\mh Z^{\mr{gcd}(\mu^\vee ) 2^{b(\mu^\vee ) -1}}
\end{aligned}
\end{equation}
By Theorem \ref{thm:2.15} $K_* \big(C_r^* (\mc R ,q)_{P_\mu}\big)$
must also be a free abelian group of rank $\mr{gcd}(\mu^\vee)
2^{b(\mu^\vee ) -1}$. 

Summing over partitions $\mu$ of $n$ we find that $K_* \big( 
C_r^* (\mc R ,q) \big)$ is a free abelian group of rank
\[
\sum_{\mu \vdash n} \mr{gcd}(\mu^\vee ) 2^{b(\mu^\vee ) -1} =
\sum_{\mu \vdash n} \mr{gcd}(\mu ) 2^{b(\mu ) -1} 
\]
We conclude that for the root datum $\mc R (A_{n-1})$ 
\begin{equation}\label{eq:6.51}
K_* \big( C_r^* (\mc R ,q) \big) \cong K_* \big( C_r^* (W) \big)
\end{equation}

\section{$B_n$}
\label{sec:6.7}

The root systems of type $B_n$ are more complicated than those 
of type $A_n$ because there are roots of different lengths. This 
implies that the associated root data allow label functions 
which have three independent parameters. Detailed information
about the representations of type $B_n$ affine Hecke algebras
is available from \cite{Slo2}.

We will compare the $C^*$-algebras for generic labelled root 
data with the reduced $C^*$-algebra of the affine Weyl group 
of type $B_n$. Consider the root datum $\mc R (B_n )$ where 
$X$ is the root lattice:
\[
\begin{aligned}
& X = Q = \mh Z^n \quad X^+ = \{ x \in X : x_1 \geq x_2 \geq
  \cdots \geq x_n \geq 0 \} \\
& Y = \mh Z^n \quad Q^\vee = \{ y \in Y : y_1 + \cdots + y_n 
  \; \mr{even} \} \\
& T = (\mh C^\times)^n \quad t = (t_1, \ldots, t_n) = 
  (t(e_1), \ldots, t(e_n)) \\
& R_0 = \{ x \in X : \norm{x} = 1 \;\mr{or}\; \norm{x} = \sqrt 2 \}
  \quad R_1 = \{ x \in X : \norm{x} = 2 \;\mr{or}\; \norm{x} = 
  \sqrt 2 \} \\
& R_0^\vee = \{ x \in X : \norm{x} = 2 \;\mr{or}\; \norm{x} = 
  \sqrt 2 \} \quad R_1^\vee = \{ x \in X : \norm{x} = 1 \;\mr{or}\;
  \norm{x} = \sqrt 2 \} \\
& F_0 = \{ \alpha_i = e_i - e_{i+1} : i=1, \ldots, n-1 \} \cup
  \{\alpha_n = e_n\} \quad \alpha_0 = e_1 \\
& s_i = s_{\alpha_i} \quad s_0 = t_{\alpha_0} s_{\alpha_0} : x \to
  x + \inp{\alpha_0^\vee}{x} \alpha_0  \\
& W_0 = \langle s_1, \ldots, s_n | s_j^2 = (s_i s_j )^2 = (s_i
  s_{i+1})^3 = (s_{n-1} s_n )^4 = e : i \leq n-2, |i-j| > 1 \rangle \\
& S_{\mr{aff}} = \{ s_0 ,s_1, \ldots ,s_{n-1} ,s_n \} \quad
  \Omega = \{ e\} \\
& W = W_{\mr{aff}} = \langle W_0, s_0 | s_0^2 = 
  (s_0 s_i)^2 = (s_0 s_1)^4 = e : i \geq 2 \rangle 
\end{aligned}
\]
For a generic label function we have different labels
$q_0 = q(s_0 ) ,\\ q_1 = q(s_i ) ,\, 1 \leq i < n$ and 
$q_2 = q(s_n )$. For completeness we mention that
\[
q_{\alpha_i^\vee} = q_1 \quad q_{\alpha_n^\vee} = q_0 \quad
q_{\alpha_n^\vee / 2} = q_2 q_0^{-1}
\]
The finite reflection group $W_0 = W_0 (B_n )$ is naturally 
isomorphic to $(\mh Z / 2 \mh Z )^n \rtimes S_n$. If $\mu \vdash n$
then the $W_0$-stabilizer of 
\[
\big( (1)^{\mu_1} (-1)^{\mu_2} (t_{\mu_1 + \mu_2 + 1})^{\mu_3}
\cdots (t_n )^{\mu_d} \big) \in T
\]
is isomorphic to 
\[
W_0 (B_{\mu_1}) \times W_0 (B_{\mu_2}) \times S_{\mu_3} \times
S_{\mu_d}
\]

\begin{itemize}
\item \textbf{\large group case $\mb{q_0 = q_1 = q_2 = 1}$}
\end{itemize}
In view of \eqref{eq:6.37} we want to determine the extended
quotient $\widetilde{T} / W_0$. Therefore we start with the
classification of conjugacy classes in $W_0$. We already know that
the quotient of $W_0$ by the normal subgroup $(\mh Z / 2 \mh Z )^n$
of sign changes is isomorphic to $S_n$, and that conjugacy classes
in $S_n$ are parametrized by partitions of $n$. So we wonder what
the different conjugacy classes in $(\mh Z / 2 \mh Z )^n \sigma 
(\mu )$ are for $\mu \vdash n$. 

To handle this we introduce the some notations. Assume that 
$|\mu | + |\lambda | = n$ and $|\mu | + |\lambda | + |\rho | = n'$.
\begin{equation}\label{eq:6.52}
\begin{array}{lll}
\nu_I & = & 
  \prod_{i \in I} s_{e_i} \quad I \subset \{1, \ldots, n\} \\
I_\lambda & = & 
  \{ 1, 1 + \lambda_1, 1 + \lambda_1 + \lambda_2, \cdots \} 
  \quad \lambda = (\lambda_1, \lambda_2, \lambda_3 ,\ldots ) \\
\sigma'(\lambda) & = & \nu_{I_\lambda} \sigma(\lambda) \in
  W_0 (B_{|\lambda |}) \\
\sigma(\mu,\lambda) & = &
  \sigma(\mu) \: (m \to m - |\lambda| \: \mr{mod} \:n) \,\sigma'
  (\lambda)\, (m \to m + |\lambda| \: \mr{mod} \:n) \\
\sigma(\mu, \lambda, \rho) & = &
  \sigma(\mu,\lambda) \: (m \to m - |\rho | \: \mr{mod} \:n') 
  \,\sigma'(\rho )\, (m \to m + |\rho | \: \mr{mod} \:n')  
\end{array}
\end{equation}
Let $I \subset \{1, \ldots, m \}$ and $J \subset \{m+1, \ldots, 2m 
\}$. It is easily verified that $\nu_I (1 \, 2 \cdots m)$ is 
conjugate to $\mu_J (m+1 \, m+2 \cdots 2m)$ if and only if $|I| + 
|J|$ is even. Therefore the conjugacy classes in $W_0$ are
parametrized by ordered pairs of partitions of total weight $n$.
Explicitly $(\mu ,\lambda )$ corresponds to $\sigma (\mu ,\lambda )$
as in \eqref{eq:6.52}. The set $T^{\sigma (\mu ,\lambda )}$ and the
group $Z_{W_0 (B_n )}(\sigma (\mu ,\lambda ))$ are both the direct
product of the corresponding objects for the blocks of $\mu$ and
$\lambda$, i.e. for the parts $(m,m,\ldots ,m)$. The centralizer of
$\sigma ((m)^k )$ in $W_0 (B_{km})$ is generated by
$(1 \; 2 \cdots m) ,\, \nu_{\{1,2,\ldots ,m \}}$ and the 
transpositions of cycles. 
\begin{equation}\label{eq:6.53}
(am+1 \; am+m+1) (am+2 \; am+m+2) \cdots (am+m \; am+2m)
\qquad 0 \leq a \leq k-2
\end{equation}
It follows that
\begin{equation}\label{eq:6.54}
\begin{array}{l}
Z_{W_0 (B_{km})} (\sigma ((m)^k )) \quad \cong \quad W_0 (B_k ) \\ 
\big( (\mh C^\times)^{km} \big)^{\sigma ((m)^k)} \quad = \quad
\big\{ \big( (t_1)^m (t_{m+1})^m \cdots (t_{km+1-m})^m \big) : 
t_i \in \mh C^\times \big\} \\
\big( \mh T^{km} \big)^{\sigma ((m)^k)} \big/ Z_{W_0 (B_{km})} 
(\sigma ((m)^k )) \quad \cong \quad [-1,1]^k \big/ S_k
\end{array}
\end{equation}
Now consider the following element of $W_0 (B_{km})$:
\[
\sigma' ((m)^k ) = \nu_{\{1,m+1,\ldots,km+1-m\}} \: (1\; 2 \cdots m) 
(m+1 \cdots 2m) \cdots (km+1-m \cdots km)
\]
It has only $2^k$ fixpoints, namely 
\[
\big( (\pm 1)^m (\pm 1)^m \cdots (\pm 1)^m \big)
\]
The centralizer of $\sigma' ((m)^k )$ is generated by $\nu_{\{1\}} 
(1 \; 2 \cdots m) ,\, \nu_{\{1,2, \ldots, m\}}$ and the
elements \eqref{eq:6.53}. Hence
\begin{equation}\label{eq:6.55}
\begin{array}{lll}
Z_{W_0 (B_{mk})} (\sigma' ((m)^k )) & \cong & W_0 (B_k ) \\
\big( \mh T^{km} \big)^{\sigma' ((m)^k)} \big/ Z_{W_0 (B_{mk})} 
(\sigma' ((m)^k )) & \cong & \{ (1)^{am} (-1)^{(k-a)m} : 
0 \leq a \leq k \}
\end{array}
\end{equation}
Now we can see what $T_u^{\sigma (\mu ,\lambda)} \big/ Z_{W_0} 
(\sigma (\mu, \lambda))$ looks like. Its number of components 
$N(\lambda )$ depends only on $\lambda$, and all these components 
are mutually homeomorphic contractible orbifolds, the shape and
dimension being determined by $\mu$. More precisely, for every 
block of $\mu$ of width $k$ we get a factor $[-1,1]^k / S_k$, and 
for every block of $\lambda$ of width $l$ we must multiply the 
number components by $l+1$. Alternatively, we can obtain the same
space (modulo the action of $W_0$) as
\begin{equation}\label{eq:6.56}
\begin{split}
T_u^{\sigma(\mu,\lambda)} \! \big/ Z_{W_0} (\sigma(\mu,\lambda))\; 
 & =\; \bigsqcup_{\lambda_1 \cup \lambda_2 = \lambda} 
  \big( T_u^{\sigma(\mu,\lambda_1,\lambda_2)} \big/ Z_{W_0} 
  (\sigma(\mu,\lambda_1,\lambda_2)) \big)_c \\
 & =\; \bigsqcup_{\lambda_1 \cup \lambda_2 = \lambda} \big( 
  \mh T^{|\mu|} \big)^{\sigma(\mu)} \!\!\big/ Z_{W_0(B_{|\mu|})}
  (\sigma (\mu)) \: (-1)^{|\lambda_1|} \: (1)^{|\lambda_2|} 
\end{split}
\end{equation}
where the subscript $c$ is supposed to indicate that we take only
the connected component containing the point
$\big( (1)^{|\mu|} (-1)^{|\lambda_1|} (1)^{|\lambda_2|} \big)$.

In effect we parametrized the components of the extended quotient
$\widetilde{T_u} / W_0$ by ordered triples of partitions 
$(\mu,\lambda_1,\lambda_2)$ of total weight $n$, and every such
components is contractible. Denote the number of ordered $k$-tuples 
of partitions of total weight $n$ by \inde{$\mc P (k,n)$}.
\[
\check H^* \big( \widetilde{T_u} / W_0 ;\mh Z \big) =
\check H^0 \big( \widetilde{T_u} / W_0 ;\mh Z \big) \cong
\mh Z^{\mc P (3,n)}
\]
By \eqref{eq:6.38} also
\begin{equation}
K_* \big( C_r^* (W) \big) = K_0 \big( C_r^* (W) \big) \cong
\mh Z^{\mc P (3,n)}
\end{equation}

\begin{itemize}
\item \textbf{\large generic case}
\item $P_\mu = F_0  \setminus \{\alpha_{\mu_1}, \alpha_{\mu_1 +
 \mu_2}, \ldots, \alpha_{|\mu|} \}$
\end{itemize}
The inequivalent subsets of $F_0$ are parametrized by partitions
$\mu$ of weight at most $n$.
\[
\begin{aligned}
& R_{P_\mu} \cong (A_{n-1})^{m_n} \times \cdots \times (A_1)^{m_2} 
\times B_{n - |\mu|} \\
& R_{P_\mu}^\vee \cong (A_{n-1})^{m_n} \times \cdots \times 
(A_1)^{m_2} \times C_{n - |\mu|} \\
& X^{P_\mu} \cong \mh Z (e_1 + \cdots + e_{\mu_1}) /\mu_1 + \cdots 
  + \mh Z (e_{|\mu|+1- \mu_d} + \cdots + e_{|\mu|}) /\mu_d \\
& X_{P_\mu} \cong ( \mh Z^n / \mh Z (e_1 + \cdots + e_n ) )^{m_n}
  \times \cdots \times (\mh Z^2 / (\mh Z (e_1 + e_2 ) )^{m_2} 
  \times \mh Z^{n - |\mu |} \\
& Y^{P_\mu} = \mh Z(e_1 + \cdots + e_{\mu_1}) + \cdots + \mh Z 
  (e_{|\mu|+1-\mu_d} + \cdots +e_{|\mu|}) \\
& Y_{P_\mu} = \{ y \in \mh Z^n : y_1 + \cdots + y_{\mu_1} = \cdots 
  = y_{|\mu|+1-\mu_d} + \cdots + y_{|\mu|} = 0 \} \\
& T^{P_\mu} = \{ (t_1 )^{\mu_1} (t_{\mu +1})^{\mu_2} \cdots
  (t_{|\mu |})^{\mu_d} (1)^{n - |\mu |} : t_i \in \mh C^\times \} \\
& T_{P_\mu} = \{ t \in (\mh C^\times )^n : t_1 \cdots t_{\mu_1} =
  t_{\mu_1} \cdots t_{\mu_1 + \mu_2} = \cdots = t_{|\mu | + 1 - 
  \mu_d} \cdots t_{|\mu |} = 1 \} \\
& K_{P_\mu} = \{ t \in T^{P_\mu} : t_1^{\mu_1} = \cdots = 
  t_{|\mu|}^{\mu_d} = 1 \} \\
& W_{P_\mu} \cong S_n^{m_n} \times \cdots \times S_2^{m_2} \times
  W_0 (B_{n - |\mu |}) \\
& W(P_\mu ,P_\mu ) \cong W_0 (B_{m_n}) \times \cdots \times 
  W_0 (B_{m_1}) 
\end{aligned}
\]
We see that $\mc R_{P_\mu}$ is the product of various root data
of type $A_m$ and one factor $\mc R (B_{n - |\mu |})$. Hence by
\eqref{eq:3.46} $\mc H_{P_\mu}$ is the tensor product of a type
$A$ part and a type $B$ part. From our study of $\mc R (A_{n-1})$ 
we recall that the discrete series representations of the type 
$A$ part of $\mc H_{P_\mu}$ are in bijection with $K_{P_\mu}$. 
From \cite[Proposition 4.3]{HeOp} and \cite[Appendix A.2]{Opd3}
we know that the
residual points for $B_{n - |\mu |}$ are parametrized by ordered
pairs $(\lambda_1 ,\lambda_2 )$ of total weight $n - |\mu |$.
The unitary part of such a residual point is in the component we
indicated in \eqref{eq:6.56}. Let $RP (\mc R ,q)$ denote the 
collection of residual points for the pair $(\mc R ,q)$.
\begin{multline}\label{eq:6.57}
\bigsqcup_{t \in RP (\mc R_{P_\mu}, q_{P_\mu})} t T_u^{P_\mu} \big/ 
\mc W_{P_\mu P_\mu} \cong \; \bigsqcup_{t \in RP (\mc R 
(B_{n - |\mu|}, q))} t T_u^{P_\mu} \big/ W (P_\mu ,P_\mu ) \\
\cong \; T_u^{P_\mu} \big/ Z_{W_0 (B_{|\mu |})} (\sigma (\mu )) 
  \times \bigsqcup_{(\lambda_1 ,\lambda_2) : |\lambda_1 | + 
  |\lambda_2 | = n} (-1)^{|\lambda_1 |} (1)^{|\lambda_2 |}
\end{multline}
This space is diffeomorphic to the extended quotient described 
on page \pageref{eq:6.56}. By Theorem \ref{thm:3.15} every point 
is the central character of at least one irreducible $C_r^* (\mc R
,q)$-representation. Therefore it is natural to compare Conjecture
\ref{conj:5.26} with the following statements.

\begin{description}
\item[1)] Every parabolically induced $C_r^* (\mc R (B_n )^\vee 
,q)$-representation $\pi (P,\delta ,t)$ is irreducible.
\item[2)] Prim$\big( C_r^* (\mc R (B_n )^\vee ,q) \big)$ 
is naturally in bijection with \eqref{eq:6.57}.
\end{description}

Opdam and Slooten \cite{Opd4} have announced a proof of 1). 
Extending the results from \cite{Slo3} they can show that all 
R-groups are trivial in this situation. In view of the above 
calculations and \eqref{eq:5.48}, 1) and Theorem \ref{thm:5.29} 
together imply 2). Probably 2) can also be derived from the 
recent work of Kato \cite{Kat4,Kat5}. 

\begin{thm}\label{thm:6.4}
Let $q$ be a generic positive label function on the root datum
$\mc R (B_n )$. Then Conjecture \ref{conj:5.26} holds. In 
particular there are precisely $\mc P (2,n)$ inequivalent
discrete series representations, one for each $W_0$-orbit of
residual points.
\end{thm}
\emph{Proof.}
This theorem was predicted in \cite[\S 8.1]{Opd4}. The proof 
relies on 1) and on Theorem \ref{thm:5.29}, whose proofs
will appear elsewhere. As said, via \eqref{eq:5.48} these 
results imply 2) and in particular the statement about the 
discrete series.

For every $W_0 t_0 \in T_u /W_0$ the image of 
$G (\phi_{W_0 t_0})$ is spanned by the modules 
\begin{equation}\label{eq:6.67}
G(\phi_0 ) \pi (P,W_P r,\delta ,t) = \pi_0 (P,W_P r_u ,
\tilde \sigma_0 (\delta ),t) \quad \mr{with} \quad 
(P,\delta )\in \Delta ,\,  W_0 r_u t = W_0 t_0
\end{equation}
By 2) the number of such modules equals the number of irreducible\\
$W$-representations with central character $W_0 t_0$.
So we only have to show that the elements \eqref{eq:6.67} are
linearly independent in $G (C_r^* (W))$. Abbreviate
\[
E_0 = \mr{Rep}_{W_0 t_0} \big( C_r^* (W) \big) 
\otimes_{\mh Z} \mh C
\]
Fix a set $\mc P'$ of representatives for the action of $W_0$ 
on the power set of $F_0$. For $P,Q \in \mc P'$ we write 
\[
P \leq Q \text{ if } \exists w \in W_0 : w P \subset Q
\]
Pick $t_1 \in W_0 t_0$ such that
\[
P_1 := \{ \alpha \in F_1 : \theta_\alpha (t_1 ) = 1 \}
\]
is maximal. We may assume that the set $P_0$ of short roots
corresponding to $P_1$ is an element of $\mc P'$. From the 
proof of Theorem \ref{thm:3.28} we see that the standard 
pairing between representations of the isotropy group 
$W_{0,t_1}$ gives an inner product on $E_0$. 

Now we can decompose $E_0$ into subspaces corresponding to the 
$P \in \mc P'$. Let $E_1^P$ be the span of the virtual 
representations in $E_0$ whose character only allows continuous
deformations along $T_u^P / W_0$, not along other directions on 
$T_u / W_0$. We put
\[
E_0^P = E_1^P \cap \Big( \sum_{Q > P} E_1^Q \Big)^\perp
\]
From Theorem \ref{thm:3.28}, \eqref{eq:6.56}, \eqref{eq:6.57} 
and \eqref{eq:2.23} (for $HH^0$) we deduce
\begin{equation}\label{eq:6.68}
\dim E_0^P = \# \{ (P,W_P r,\delta ,t) \in \Xi_u : 
W_0 r_u t = W_0 t_0 \}
\end{equation}
Notice that this is a simple version of \cite[Proposition 
6.6]{Opd4}. Upon applying \eqref{eq:5.48} to the affine Hecke 
algebra $\mc H_P$ (which has generic labels), we see that 
\begin{equation}\label{eq:6.69}
G(\phi_0 ) \otimes \mr{id}_{\mh C} \big( \mr{span} 
\{ (P,W_P r,\delta ,t) \in \Xi_u : W_0 r_u t = W_0 t_0 \} 
\big) \subset E_1^P
\end{equation}
and that $G (\phi_0 ) \otimes \mr{id}_{\mh C}$ is injective 
on this domain. Moreover by \eqref{eq:6.68} the image in 
\eqref{eq:6.69} intersects $\sum_{Q > P} E_1^Q$ only in 0.
Therefore $G (\phi_0 )$ is injective. $\qquad \Box$

%% file: appendix.tex
\chapter{Crossed products}

We collect some well-known results on crossed products of 
algebras by compact groups. We do this in the large category 
of $m$-algebras, 
but it is straightforward to see that they also hold for 
$C^*$-algebras, if we assume that the action is *-preserving. 
So let $A$ be an $m$-algebra, $G$ a compact group (with its 
normalized Haar measure) and 
\begin{equation}
\begin{aligned}
& \alpha: G \times A \to A \\
& \alpha (g,a) = \alpha_g(a)
\end{aligned}
\end{equation}
a continuous action of $G$ on $A$ by algebra homomorphisms. 
Recall that the \inde{crossed product} $A \rtimes_\alpha G$ 
is the vector space $C(G;A)$ with multiplication defined by
\begin{equation}
(f \cdot f') (g') = \int_G f(g) \alpha_g 
  (f'(g^{-1} g')) dg
\end{equation}
This is again an $m$-algebra. Notice that if $a \in A$ and 
$\delta_g$ is the $\delta$-function concentrated at $g \in G$, 
then $a \delta_g$ is an element of the multiplier algebra 
$\mc M (A \rtimes_\alpha G)$.
Explicitly, multiplying by this element is defined as
\begin{align}
(f \cdot a \delta_g )(g') &= f(g' g^{-1}) \alpha_{g' g^{-1}}(a)\\
(a \delta_g \cdot f) (g') &= a \alpha_g (f(g^{-1} g'))
\end{align}
Similarly, if $N$ is a closed subgroup of $G$ then any 
$\psi \in C^*(N)$ is also a multiplier of $A \rtimes_\alpha G$, 
defined naturally by
\begin{align}
(f \cdot \psi)(g') &= \int_N f(g' n^{-1}) \psi(n) dn \\
(\psi \cdot f)(g') &= \int_N \psi(n) \alpha_n (f(n^{-1} g')) dn
\end{align}
The next result is useful in connection with projective 
representations.

\begin{lem}\label{lem:A.1}
Let $\{ e \} \to N \to G \to H \to \{ e \}$ be a short exact 
sequence of compact groups, all equipped with their normalized 
Haar measures. Assume that $N \subset \ker \alpha$, so that 
$\alpha$ descends to $H$. Let $p_N \in C^*(N)$ be the constant 
function with value 1, considered as an idempotent in 
$\mc M(A \rtimes_\alpha G)$.
 
a) $p_N$ is central and
$p_N (A \rtimes_\alpha G) \cong (A \rtimes_\alpha H)$\\
Assume now that moreover $N$ is finite and the extension from 
$H$ to $G$ by $N$ central. Then we let $\widehat N$ be the dual 
of the abelian group $N$ and we consider every (one-dimensional) 
character $\chi$ of $N$ as an idempotent 
$p_\chi \in \mc M(A \rtimes_\alpha G)$. In particular 
$p_N = p_\chi$ for the trivial character $\chi \equiv 1$. 

b) $p_\chi$ is also central and 
$A \rtimes_\alpha G = \bigoplus_{\chi \in \widehat N} 
p_\chi \left( A \rtimes_\alpha G \right)$
\end{lem}
\emph{Proof.}
For any $b \in A \rtimes_\alpha G$ and $g' \in G$ we have
\begin{equation}
\begin{aligned}
(b \cdot p_\chi)(g') & =\; 
 \int_G b(g) \alpha_g (p_\chi(g^{-1} g')) dg \\
 & =\; \int_N b(g' g^{-1}) \chi(g^{-1}) |N|^{-1} dg \\
 & =\; \int_N |N|^{-1} \chi(g^{-1}) b(g^{-1} g') dg \\
 & =\; \int_G p_\chi(g) \alpha_g (b(g^{-1} g')) dg 
\;=\; (p_\chi \cdot b)(g')
\end{aligned}
\end{equation}
The third equality holds if $N$ is central or if $\chi \equiv 1$ 
and $N$ only normal, that is, in all the cases we need. 
So $p_\chi$ indeed commutes with $A \rtimes_\alpha G$ and 
$p_\chi (A \rtimes_\alpha G)$ is a subalgebra. Now the 
statement b) follows from two standard equalities in the 
representation theory of finite groups: 
\begin{equation}
\sum_{\chi \in \hat N} p_\chi = 1 \qquad \mr{and} \qquad
p_\chi \cdot p_{\chi'} = \delta_{\chi \chi'} p_\chi 
\end{equation}
Writing $\pi: G \to H$, it is easily verified that
\begin{equation}
\pi^* : A \rtimes_\alpha H \to A \rtimes_\alpha G :
f \to f \circ \pi
\end{equation}
is a monomorphism of $m$-algebras, so let us determine its image. 
Since
\begin{equation}
\begin{split}
(p_N\cdot \pi^* f)(g') &= 
  \int_G p_N(g) \alpha_g (\pi^* f(g^{-1}g')) dg \\
 &= \int_N p_N(g) f(\pi(g^{-1}g')) dg \\
 &= \int_N |N|^{-1} f(\pi g') dg \;=\; \pi^* f(g') 
\end{split}
\end{equation}
the image is contained in $p_N (A \rtimes_\alpha G)$. From the 
above expressions for $b p_N$ it follows immediately that 
anything of this type is $N$-biinvariant, so that it descends 
to an element of $A \rtimes_\alpha H$. Hence the image of 
$\pi^*$ is exactly $p_N (A \rtimes_\alpha G) \qquad \Box$
\\[2mm]

Suppose now that $A$ is unital and we are given 
$p \in C(G; A^\times)$ with the properties
\begin{align}
p(g g') &= p(g) \alpha_g (p(g')) \\
p(e) &= 1 
\end{align}
Then $p$ is an idempotent in $A \rtimes_\alpha G$ and we can 
define a new action $\beta$ of $G$ on $A$ by 
\begin{equation}
\beta_g (a) = p(g) \alpha_g (a) p(g)^{-1}
\end{equation}
The following description of the invariant algebra
\begin{equation}
A^{\beta (G)} := \{a \in A : \beta_g (a) = a \; \forall g \in G\}
\end{equation}
is essentially due to Rosenberg \cite{Ros}.

\begin{lem}\label{lem:A.2}
\[
A^{\beta (G)} \cong p (A \rtimes_\alpha G) p
\]
\end{lem}
\emph{Proof.}
We will show that the obvious map 
\begin{equation}
\begin{aligned}
& \phi: A^{\beta (G)} \to p (A \rtimes_\alpha G) p \\
& \phi (a) = p (a \delta_e )p
\end{aligned}
\end{equation}
is an isomorphism. Clearly $\phi$ is linear and continuous.
If $a, a' \in A^\beta$ then 
\begin{align}
\phi(a) &= p (a \delta_e) p = p (a \delta_e) = (a \delta_e) p \\
\phi(a) \phi(a') &= p (a \delta_e) p p (a' \delta_e) p = p (a 
  \delta_e ) (a' \delta_e ) p =  p (a a' \delta_e ) p = \phi(a a')
\end{align}
so $\phi$ turns out to be injective and multiplicative.
Next we observe that for any $g \in G$
\begin{equation}
p(g) \delta_g \cdot p = p = p \cdot p(g) \delta_g
\end{equation}
This gives, for $b \in p (A \rtimes_\alpha G) p$,
\begin{align}
b(g) &= (p(g) \delta_g \cdot b)(g) = p(g) \alpha_g (b(e)) \\
b(g) &= (b \cdot p(g) \delta_g)(g) = b(e) p(g)
\end{align} 
Comparing these expressions we see that $b(e) \in A^\beta$. 
Therefore $\phi$ is bijective, with continuous inverse 
$b \to b(e). \quad \Box$
\\[2mm]

Assume now that $G$ is finite. Then evaluating integrals of 
$A$-valued functions on $G$ is easy, so we do not need any 
topology on $A$ to define $A \rtimes_\alpha G$. Moreover we agree 
to use the counting measure on $G$, even though it is not 
normalized. The crossed product thus obtained also appears in 
another way:

\begin{lem}\label{lem:A.3}
Let $(\mh C [G], \rho)$ be the \inde{right regular 
representation} of $G$ and $A$ any complex algebra. Endow 
$A \otimes \mr{End} (\mh C [G] \otimes \mh C^n)$ with the 
$G$-action $g (a) := \rho(g) \alpha_g (a) \rho(g^{-1})$, where 
$\alpha$ stands also for the action of $G$ on $A$ tensored with 
the identity. Then
\[
\big( A \otimes \mr{End}(\mh C [G] \otimes \mh C^n) \big)^G 
\cong M_n (A \rtimes_\alpha G)
\]
\end{lem}
\emph{Proof.}
The left hand side is isomorphic to $M_n \left( \big( A \otimes 
\mr{End}(\mh C [G]) \big)^G \right)$, 
so it suffices to prove the case $n = 1$.

For $a \in A$ and $g,h \in G$ define 
$L(a \otimes g) (h) = \alpha_{h^{-1} g^{-1}}(a) \otimes g h$.
It is easily checked that this extends to an algebra homomorphism
$L: A \rtimes_\alpha G \to \big( A \otimes \mr{End}(\mh C [G]) 
\big)^G$. We claim that 
$L' : b \to \sum_{g \in G} b(g^{-1})_e \otimes g$
is the inverse of $L$. It is clear that $L' L = \,$Id, so we only 
have to show that $L(L'b) = b$ for any
$b \in \big( B \otimes \mr{End}(\mh C [G]) \big)^G$.
\begin{equation}
\begin{split}
L(L'b)(h) &= \sum_{g \in G} \alpha_{h^{-1} g^{-1}} (b(g^{-1})_e) 
  \otimes gh \\
 &= \alpha_{h^{-1}} \alpha_{g^{-1}}(b)(g^{-1})_e \otimes gh \\
 &= \sum_{g \in G} \alpha_{h^{-1}} \left( \rho(g) b \rho(g^{-1}) 
  (g^{-1}) \right)_e \otimes gh \\
 &= \sum_{g \in G} \alpha_{h^{-1}} \left( b(e) g^{-1} \right)_e 
  \otimes gh \\
 &= \alpha_{h^{-1}} (b)(e) h \\
 &= \rho(h^{-1}) \alpha_{h^{-1}}(b) \rho(h)(h) \;=\; b(h)
\end{split}
\end{equation}
This holds for any $h \in G$, so indeed $L' = L^{-1} \qquad \Box$

%% file: summary.tex
\chapter*{Samenvatting}
\addcontentsline{toc}{chapter}{Samenvatting}
\fancyhead[LO,RE]{Samenvatting}

De afgelopen jaren is mij vaak gevraagd wat ik nou eigenlijk 
onderzoek. Op deze vraag heb ik inmiddels een voorraadje antwoorden
uitgeprobeerd, die in zekere zin allemaal wel correct waren. 
Niettemin bleek dat sommige antwoorden aanzienlijk
meer begrip en waardering oogsten dan andere. E\'en persoon
besloot zelfs geheel af te zien van verdere communicatie nadat
ik haar de titel van dit proefschrift had verteld.

Daarom lijkt het me wel een goed idee om in ieder geval de
wellicht enigszins cryptische zinsnede \emph{Periodiek cyclische
homologie van affiene Hecke algebra's} toe te lichten. Dit
wordt dan ook niet zozeer een samenvatting van mijn onderzoek,
als wel een relatief gezellige wandeling langs de randen van de
oneindig dimensionale ruimten waarin ik mij gewoonlijk begeef.
\\[1mm]

Laten we beginnen bij groepen. Met groepen kun je de symmetrie\"en
beschrijven van uiteenlopende dingen, zoals een voetbal, een
velletje papier, een molecuul, een differentiaalvergelijking of 
ruimte-tijd, maar ook van simpele figuren als een lijn, een kubus 
of een zevenhoek.
Een eenvoudige groep, die een rol speelt in dit boek, bestaat
uit de zes symmetrie\"en van een gelijkzijdige driehoek.

\begin{center}
\includegraphics[width=36mm,height=34mm]{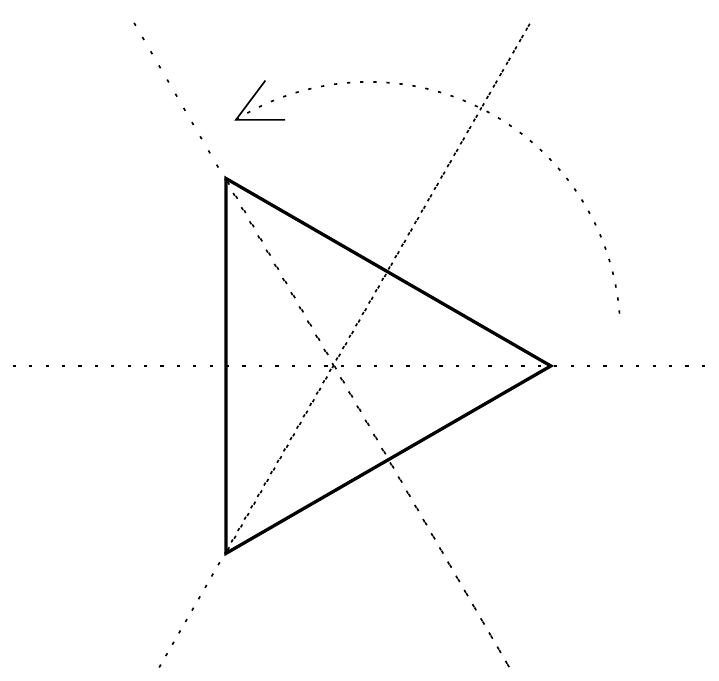}
\end{center}
We zien drie spiegelingen (in de stippellijnen) en rotaties over
$120^\circ$ en over $240^\circ$. De laatste symmetrie is de 
afbeelding die alles op zijn plek laat. We zien direct dat 
\emph{groep} binnen de wiskunde iets heel anders betekent dan in het
dagelijks leven. In abstracto is het een verzameling waarvan je
de elementen kunt samenstellen, waarbij aan bepaalde voorwaarden
voldaan moet zijn. 

De bovenstaande groep heeft allerlei bijzondere eigenschappen.
Bijvoorbeeld, als je een willekeurige lijn neemt die door het 
middelpunt van de driehoek gaat, en je past daarop een symmetrie 
van de driehoek toe, dan krijg je weer een lijn die door dat 
middelpunt gaat. Men zegt daarom wel dat deze groep bestaat uit 
lineaire afbeeldingen.

Beschouw het rooster dat is opgebouwd uit gelijkzijdige driehoeken.
(Zie de kaft voor een artistiekere impressie van dit rooster,
waarvoor mijn dank uitgaat naar Bill Wenger.)

\begin{center}
\includegraphics[width=46mm,height=40mm]{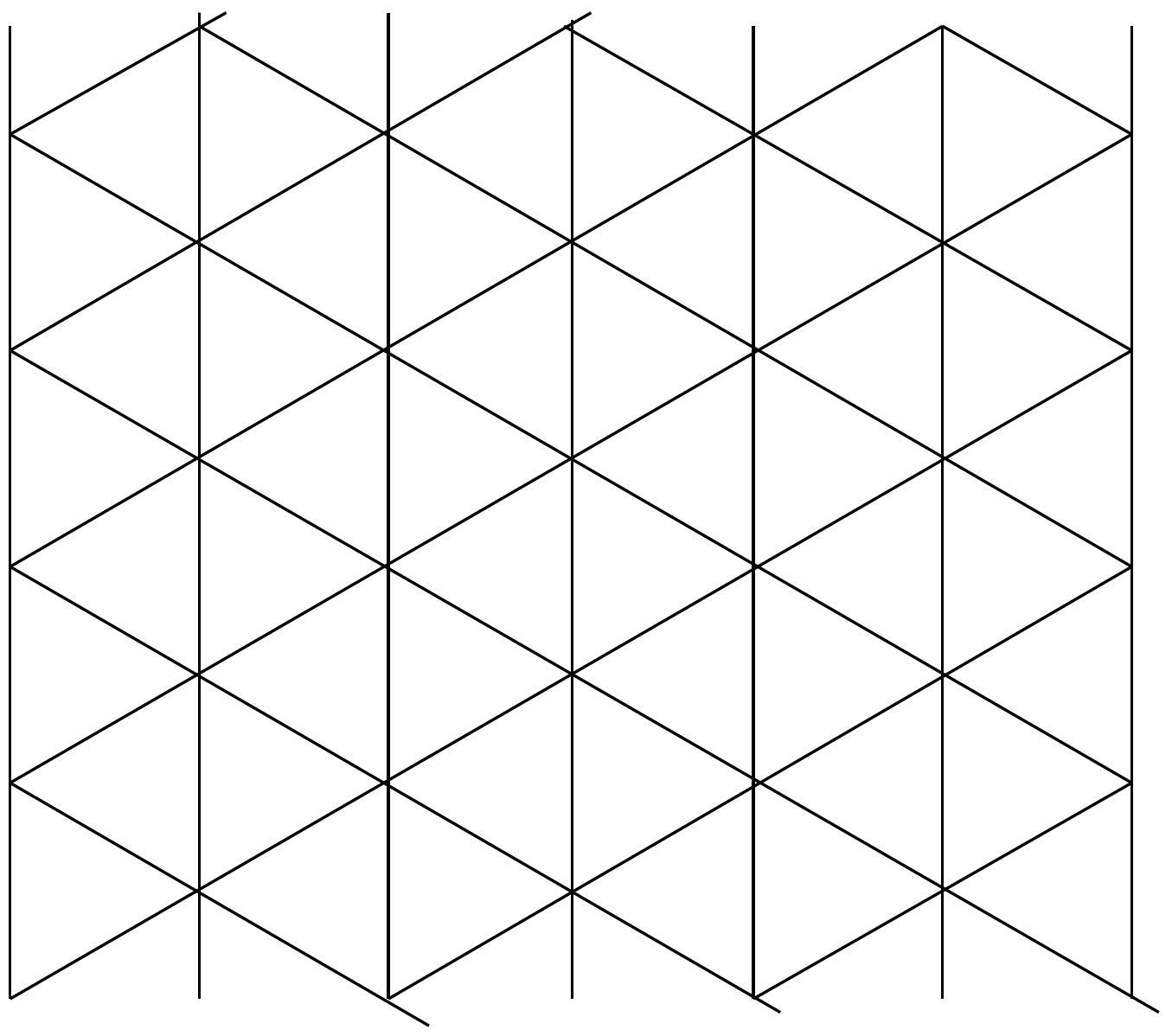}
\end{center}
De symmetriegroep van dit (oneindig grote) figuur bevat oneindig 
veel elementen, onder andere spiegelingen en rotaties, maar ook 
verschuivingen. Stel dat je een punt in dit rooster kiest, en een 
lijn door dat punt. Als je zomaar een symmetrie van het rooster 
toepast krijg je zeker weer een lijn, maar er is geen garantie dat 
die nog steeds door het gekozen punt gaat. Daarom noemt men dit soort 
symmetrie\"en geen lineaire afbeelingen, maar \emph{affiene} 
afbeeldingen.

Met deze hele opzet kunnen we iets doen waar wiskundigen dol op 
zijn, we kunnen het zaakje generaliseren. In dat geval vervangen we
de driehoek door een ingewikkelder kristal, en het driehoekige 
rooster door een rooster van hogere dimensie. De symmetriegroep
van het kristal is een (eindige) Weyl groep, en de symmetriegroep
van het rooster heet een affiene Weyl groep.
\\[1mm]

Als je een groep goed wilt begrijpen is het van groot belang om zijn
zogeheten representaties te kennen. Dit illustreren we aan de hand
van ander voorbeeld, de cirkel. Enerzijds kan deze worden beschouwd 
als de groep van rotaties om zijn eigen middelpunt. Anderszijds 
kunnen we de cirkel opvatten als een snaar, en dan kan hij trillen.
De representaties van de cirkelgroep corresponderen precies met de 
trillingen van de cirkelvormige snaar, waarbij we \'e\'en specifiek
punt P vasthouden. 

\begin{center}
\includegraphics[width=110mm,height=21mm]{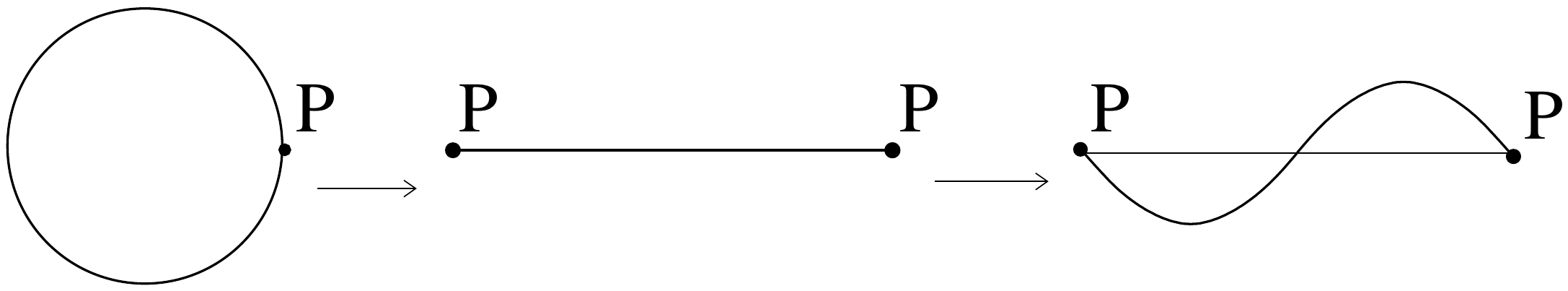}
\end{center}
Er is een grondtoon, waarvan de golflengte exact de lengte van 
de snaar is. De boventonen hebben een golflengte die een geheel 
aantal keren in de snaar past. Elke trilling is te maken als een 
geschikte combinatie van dergelijke harmonische trillingen.

In wiskundig jargon betekent dit dat de harmonische trillingen
corresponderen met de irreducibele representaties. Irreducibele
representaties zijn een soort bouwsteentjes waarmee je elke 
representatie kunt maken. Ze zijn zelf niet verder op te delen.
\\[1mm]

Grof gezegd is een algebra een verzameling waarbinnen je kunt 
optellen en vermenigvuldigen. Zo vormen de gehele getallen een 
algebra. Een wat lastiger voorbeeld zijn de $2 \times 2 $-matrices 
met re\"ele co\"efficienten:
\[
M_2 (\mh R ) := \left\{ \begin{pmatrix} x_1 & x_2 \\
x_3 & x_4 \end{pmatrix} : x_1 ,x_2 ,x_3 ,x_4 \in \mh R \right\}
\]
Dit soort matrices kan je co\"ordinaatsgewijs optellen:
\[
\begin{pmatrix} 1 & 0 \\ \pi & 4 \end{pmatrix} +
\begin{pmatrix} 3/4 & -5/9 \\ -1 & -4 \end{pmatrix} =
\begin{pmatrix} 7/4 & -5/9 \\ \pi - 1 & 0 \end{pmatrix}
\]
Een matrix is op te vatten als een lineaire afbeelding. Als we 
punten in het vlak schrijven als vectoren 
$\begin{pmatrix} x \\ y \end{pmatrix}$, dan stuurt een matrix
dus een vector naar een andere vector:
\[
\begin{pmatrix} 3/4 & -5/9 \\ -1 & -4 \end{pmatrix} 
\begin{pmatrix} a \\ 0 \end{pmatrix} = 
\begin{pmatrix} 3 a / 4 \\ -a \end{pmatrix} 
\qquad \qquad
\begin{pmatrix} 3/4 & -5/9 \\ -1 & -4 \end{pmatrix} 
\begin{pmatrix} 0 \\ b \end{pmatrix} = 
\begin{pmatrix} -5 b / 9 \\ -4b \end{pmatrix}
\]
De standaardmanier om matrices te vermenigvuldigen correspondeert
met het samenstellen van afbeeldingen, bijvoorbeeld
\[
\begin{pmatrix} 0 & 0 \\ 3 & 0 \end{pmatrix} \cdot
\begin{pmatrix} 0 & 2 \\ 0 & 0 \end{pmatrix} =
\begin{pmatrix} 0 & 0 \\ 0 & 6 \end{pmatrix} 
\qquad \qquad
\begin{pmatrix} 0 & 2 \\ 0 & 0 \end{pmatrix} \cdot
\begin{pmatrix} 0 & 0 \\ 3 & 0 \end{pmatrix} =
\begin{pmatrix} 6 & 0 \\ 0 & 0 \end{pmatrix} 
\]
Al met al maakt dit $M_2 (\mh R )$ tot een algebra over de re\"ele
getallen $\mh R$. Merk op dat bepaalde gebruikelijke rekenregels 
voor getallen niet meer gelden voor matrices. We zijn gewend dat
$7 \cdot 9 = 9 \cdot 7 = 63$. Zelfs zonder dat we de uitkomst
weten kunnen we met zekerheid zeggen dat 
\[
4676013 \cdot 2369655 = 2369655 \cdot 4676013
\]
In feite geldt $x \cdot y = y \cdot x$ voor alle re\"ele getallen
$x$ en $y$. Niettemin zien we hierboven dat er matrices $A$ en $B$ 
bestaan zodat 
\[ 
A \cdot B \neq B \cdot A
\]
Men zegt dan dat $A$ en $B$ niet commuteren en dat de algebra
$M_2 (\mh R )$ niet commutatief is.

Wat is nu het verband tussen algebra's en groepen? Vanuit een
groep kunnen we een algebra construeren die alle eigenschappen van
de groep reflecteert. Zo hebben algebra's ook representaties, en de
representaties van een groep komen overeen met de representaties 
van de bijbehorende groepsalgebra.

Een aanzet tot de algebra's in de titel van dit boek werd gegeven
door de Duitse wiskundige Erich Hecke, die leefde van 1887 tot 1947.
Hecke hield zich vooral bezig met getaltheorie, bijvoorbeeld met
$p$-adische getallen. Hier is $p$ een priemgetal, bijvoorbeeld 5. 
De verzameling 5-adische getallen geven we aan met $\mh Q_5$.
Een typisch 5-adisch getal ziet eruit als een decimale expansie in 
de verkeerde richting:
\[
x = \cdots 42130012044.113 \in \mh Q_5
\]
Omdat $p = 5$ komen alleen de symbolen 0, 1, 2, 3 en 4 in deze 
schrijfwijze van $x$ voor. We dienen $x$ te interpreteren als
\[
x = 3 \cdot 5^{-3} + 1 \cdot 5^{-2} + 1 \cdot 5^{-1} + 4 \cdot 5^0 +
4 \cdot 5^1 + 0 \cdot 5^2 + 2 \cdot 5^3 + \cdots
\]
Als $y$ een ander 5-adisch getal is, dan zijn $x+y$ en $x \cdot y$
zonder al te veel problemen te bepalen met behulp van de regeltjes
\[
\begin{array}{ccccl}
a 5^n & + & b 5^n & = & (a+b) 5^n \\
a 5^n & \cdot & c 5^m & = & (a \cdot c) 5^{n+m}
\end{array}
\]
voor gehele getallen $a,b,c,n$ en $m$. Om de co\"efficient van 
$x \cdot y$ bij $5^n$ uit te rekenen hebben we slechts kleine 
stukjes van de expansies van $x$ en $y$ nodig. Zelfs delen is 
mogelijk met $p$-adische getallen, omdat $p$ een priemgetal is.

Met $p$-adische getallen kunnen we weer groepen bouwen. Een simpel
voorbeeld van zo'n $p$-adische groep is
\[
\left\{ \begin{pmatrix} a & b \\ c & d \end{pmatrix} : 
a,b,c,d \in \mh Q_5 \:, a \cdot d - b \cdot c = 1 \right\}
\]
De samenstelling in deze groep is de gebruikelijke 
matrixvermeningvuldiging, maar dan uitgevoerd met 5-adische
getallen. Dergelijke $p$-adische groepen spelen een belangrijke 
rol in verschillende gebieden van de wiskunde. Men zou graag alle 
irreducibele representaties van zo'n groep klassificeren, maar dat
is erg lastig. Het blijkt dat een belangrijk deel van de 
representatietheorie van een $p$-adische groep valt uit te drukken
met een zekere algebra. Zo een algebra is een generalisatie van een
type algebra's dat Hecke indertijd vanuit een iets andere hoek heeft
gedefinieerd en bestudeerd, vandaar dat ze onder de naam 
\emph{Hecke algebra's} door het leven gaan.
\\[1mm]

Hecke algebra's zijn er in soorten en maten. Ik ben vooral 
ge\"interesseerd in Hecke algebra's die sterk gerelateerd zijn aan
Weyl groepen. In feite zijn deze Hecke algebra's te beschouwen als
deformaties van Weyl groepen. De groepsalgebra die bij een Weyl
groep hoort heeft zoals zoals gezegd grotendeels dezelfde 
eigenschappen als de groep zelf. Als je die groepsalgebra op een
geschikte manier vervormt krijg je een Hecke algebra.

Doe je dit met een eindige Weyl groep, dan krijg je een Hecke
algebra van eindige dimensie. Hoewel \emph{eindig dimensionaal} in
eerste instantie nog niet zo eenvoudig klinkt, is het vast
makkelijker dan \emph{oneindig dimensionaal}. In ieder geval 
begrijpt men eindig dimensionale Hecke algebra's heel goed.

Echter, als je een affiene Weyl groep deformeert krijg je een 
affiene Hecke algebra, en die heeft oneindige dimensie. Affiene 
Hecke algebra's zijn veel ingewikkelder dan Hecke algebra's 
van eindige dimensie. Toch is dat niet zo'n ramp. Men gebruikt 
affiene Hecke algebra's onder andere om tot een beter begrip te 
komen van de gecompliceerde representatietheorie van $p$-adische
groepen, en daar zou weinig van te verwachten zijn als ze te 
eenvoudig waren. Het blijkt dat affiene Hecke algebra's aan de ene
kant diepzinnig genoeg zijn om tot nieuwe inzichten te leiden, en
aan de andere kant makkelijk genoeg om er prettig mee te kunnen
werken. Dus een affiene Hecke algebra heeft precies de goede 
moeilijkheid om hem tot een interessant studieobject te maken.
\\[1mm]

Zonet zagen we de analogie tussen trillingen van een snaar en
representaties van een groep. Als we dit uitbreiden correspondeert
een algebra niet meer met \'e\'en snaar, maar met een 
snaarinstrument, bijvoorbeeld een piano. Een representatie van die
algebra wordt dan een toon die je met die piano kunt voortbrengen.
Op deze manier kunnen we een irreducibele representatie 
identificeren met een zuivere toon van de piano.

Het uiteindelijke doel van mijn promotieonderzoek was om alle
irreducibele representaties van een algemene affiene Hecke algebra
te bepalen. Het ligt voor de hand om ze eerst maar eens te tellen.
Helaas laten ze zich niet zo gemakkelijk tellen, want het zijn er
oneindig veel. Net zo heeft het weinig zin om alle zuivere tonen 
van een piano te tellen, want dat zijn er ook oneindig veel. Een 
beter idee is daarom om alle grondtonen van de snaren van de piano 
te tellen, dat vertelt je bijvoorbeeld al hoeveel snaren je piano 
heeft.

In de context van algebra's ligt dat wat subtieler, daar heet de
geschikte manier om grondtonen te tellen \emph{periodiek cyclische 
homologie}. "Homologie" komt uit het Grieks en betekent zoveel als
"studie van gelijkheid". Dat gaat ongeveer als volgt. Stel dat je
twee objecten, bijvoorbeeld twee algebra's, wilt vergelijken.
Kies een geschikte methode (een homologietheorie) om aan een 
algebra iets relatief eenvoudigs toe te kennen, bijvoorbeeld een
simpel type groep, een rijtje getallen of zelfs een rijtje
groepen. Dat heet dan de homologie van de algebra. Als je twee
algebra's in essentie hetzelfde zijn zullen ze dezelfde homologie
hebben. Daarentegen, als ze verschillende homologie hebben dan
zijn de algebra's niet hetzelfde, en ook niet ongeveer.

Op de kaft van dit boek zie je duidelijk periodieke en cyclische 
verschijnselen. Dat is geen toeval, maar de etymologische 
achtergrond van de term \emph{periodiek cyclische homologie} is 
anders. De periodiek cyclische homologie van een algebra $A$ is 
een rijtje groepen:
\[
\ldots , HP_{-2}(A) ,\, HP_{-1}(A) ,\, HP_0 (A) ,\, HP_1 (A) ,\,
HP_2 (A) ,\, HP_3 (A) , \ldots
\]
Dit is periodiek in de zin dat voor elk geheel getal $n$ geldt
\[
HP_n (A) = HP_{n+2} (A)
\]
Het cyclische zit wat dieper verstopt, dat heeft te maken met
hoe de groepen $HP_n (A)$ expliciet geconstrueerd worden. Stel dat
we zeven hokjes hebben, die allemaal gevuld zijn met een letter:

\begin{center}
\begin{tabular}{|c|c|c|c|c|c|c|}
\hline f & b & d & c & r & t & z \\ \hline
\end{tabular}
\end{center}
Nu schuiven we elke letter \'e\'en hokje naar rechts, en de meest
rechtse letter stoppen we in het vrijgekomen linker hokje:

\begin{center}
\begin{tabular}{|c|c|c|c|c|c|c|}
\hline z &f & b & d & c & r & t \\ \hline
\end{tabular}
\end{center}
Dit kunnen we nog wat suggestiever tekenen:

\begin{center}
\includegraphics[width=73mm,height=25mm]{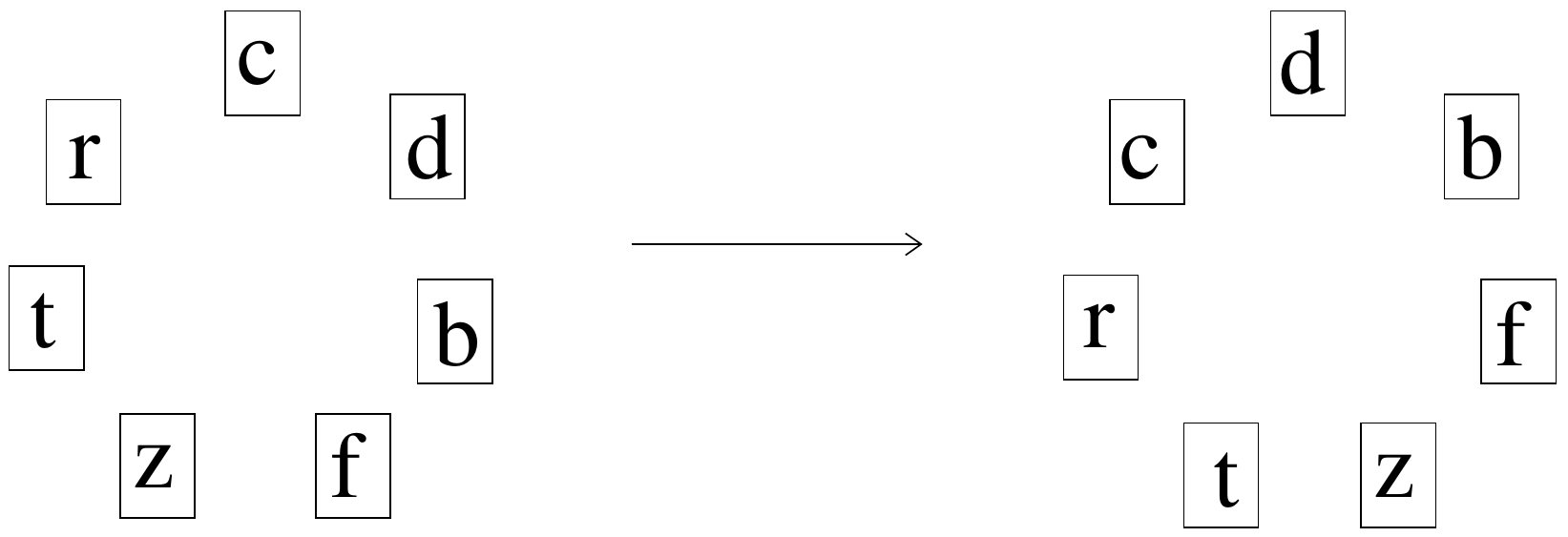}
\end{center}
Nu is het wel duidelijk waarom dit een cyclische permutatie heet. 
Zulke permutaties worden gebruikt in de definitie van periodiek 
cyclische homologie.

Het is nog niet zo eenvoudig om de periodiek cyclische homologie
van een affiene Hecke algebra ook daadwerkelijk uit te rekenen.
Daartoe grijpen we terug op de affiene Weyl groep waarvan het een
deformatie is. Als we weer denken aan piano's en snaren betekent deze
deformatie dat we wat gaan sleutelen aan de snaren: wat langer of
iets korter, een stukje dichter bij elkaar. Hoewel het in muzikaal
opzicht barbaars is zouden we zelfs sommige snaren aan elkaar vast
kunnen knopen. Zo'n deformatie kan er schematisch uitzien als

\begin{center}
\includegraphics[width=120mm,height=24mm]{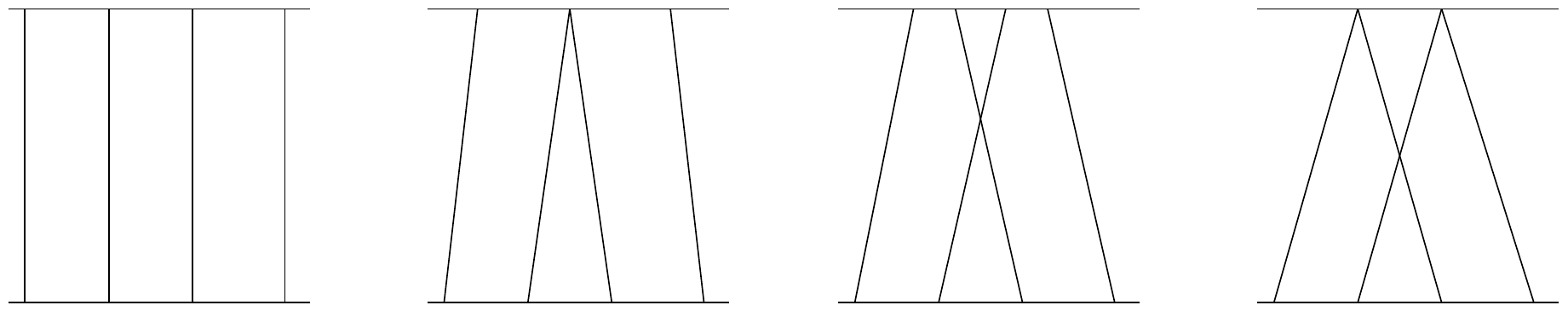}
\end{center}

Men vermoedt dat de representatietheorie van een affiene Weyl groep
niet essentieel verandert onder deze deformaties. Dit vermoeden
wordt ondersteund door diepzinnige stellingen die zeggen dat het
in bepaalde belangrijke gevallen klopt. Het is een belangrijk 
vermoeden, want hiermee kan men representaties van een affiene
Hecke algebra herleiden tot representaties van een affiene Weyl
groep, en die zijn allemaal al lang bekend.

In dit proefschrift heb ik bewezen dat dit vermoeden equivalent is 
met een ogenschijnlijk zwakkere uitspraak, namelijk dat de periodiek
cyclische homologie van de groepsalgebra van een affiene Weyl
groep niet verandert als je die algebra vervormt tot een affiene
Hecke algebra. Grof gezegd betekent de sterke versie van dit 
vermoeden dat je met elk van de boven getekende pianootjes 
evenveel verschillende tonen kan voortbrengen. De zwakke versie
zegt zo ongeveer dat al die piano's evenveel grondtonen hebben.

Verder worden in dit boek onder andere een aantal nieuwe technieken 
ge\"intro\-du\-ceerd om de periodiek cyclische homologie van een 
redelijk algemeen type algebra uit te rekenen. Mede daardoor kunnen 
de bovenstaande vermoedens nu bewezen worden in veel nieuwe gevallen.

\chapter*{Curriculum vitae}
\addcontentsline{toc}{chapter}{Curriculum vitae}
\pagestyle{empty}

De auteur werd geboren op 5 februari 1979 in Amsterdam.
Zijn middelbare schooltijd bracht hij door aan het Montessori
Lyceum Amsterdam, alwaar hij in juni 1997 voor zijn VWO-examen
slaagde. 

In september 1997 begon hij als student aan de Universiteit 
van Amsterdam. In 1998 met behaalde hij cum laude een dubbele 
propedeuse wiskunde en na\-tuur\-kunde, waarop hij besloot 
zich verder te concentreren op de wiskunde. Hij schreef zijn 
scriptie, getiteld \emph{Lie algebra cohomology and Macdonald's 
conjectures}, onder de hoede van prof. dr. Eric Opdam. Hierop
studeerde hij in september 2002 cum laude af. In november 2002
trad hij in dienst bij de Universiteit van Amsterdam als 
assistent in opleiding (AIO), eveneens onder bege\-leiding van
Eric Opdam. In deze hoedanigheid heeft hij vier jaar lang
gewerkt aan dit proefschrift en werkcolleges verzorgd voor
studenten van verschillende exacte studies. Momenteel is de
auteur werkzaam aan de Universiteit van Amsterdam als 
onderzoeker en docent wiskunde.
\\[1mm]

Min of meer parallel met zijn wiskundige activiteiten liep de
schaakcarri\`ere van de auteur. Hij reikte diverse keren tot een
gedeelde tweede plaats in Nederlandse jeugdkampioenschappen schaken:
t/m 12 jaar (1991), t/m 14 jaar (1992), t/m 16 jaar (1994) en
t/m 20 jaar (1998). In 2000 voldeed hij aan de voorwaarden voor
de titel Internationaal Meester. Hij speelde schaaktoernooien in
Nederland, Belgi\"e, Frankrijk, Andorra, Engeland, Ierland, 
Denemarken, Duitsland, Zwitserland, Liechtenstein en in de 
Dominicaanse Republiek. Een hoogtepunt kwam in 2003, toen hij
topscorer werd van de hoofdklasse van de Nederlandse 
schaakcompetitie. 

Reeds enige jaren houdt de auteur zich ook actief bezig met een
andere sport, badminton. Sinds december 2005 is hij 
penningmeester van de vereniging US badminton.

%% file: thesis.bbl
\begin{thebibliography}{999}
\addcontentsline{toc}{chapter}{Bibliography}


\bibitem{Ati} M.F. Atiyah,
\emph{$K$-theory},
Mathematics Lecture Note Series,
W.A. Benjamin, New York, 1967

\bibitem{AtHi} M.F. Atiyah, F.E.P. Hirzebruch,
``Vector bundles and homogeneous spaces'', 
pp. 7-38 in: \emph{Differential geometry},
Proc. Sympos. Pure Math. \textbf{3},
American Mathematical Society, Providence RI, 1961

\bibitem{ABP} A.-M. Aubert, P.F. Baum, R.J. Plymen, 
''The Hecke algebra of a reductive $p$-adic group:
a view from noncommutative geometry'',
pp. 1-34 in: \emph{Noncommutative geometry and number theory},
Aspects of Mathematics \textbf{E37},
Vieweg Verlag, Wiesbaden, 2006

\bibitem{BaCo} P.F. Baum, A. Connes,
``Chern character for discrete groups'', 
pp. 163-232 in: \emph{A f\^ete of topology},
Academic Press, Boston MA, 1988

\bibitem{BCH} P.F. Baum, A. Connes, N. Higson,
``Classifying space for proper actions and $K$-theory of 
group $C^*$-algebras'', pp. 240-291 in: 
\emph{$C^*$-algebras: 1943-1993, A fifty year celebration},
Contemp. Math. \textbf{167},
American Mathematical Society, Providence RI, 1994

\bibitem{BHP1} P.F. Baum, N. Higson, R.J. Plymen,
``Une d\'emonstration de la conjecture de Baum-Connes pour 
le groupe $p$-adique $GL (n)$'',
C.R. Acad. Sci. Paris \textbf{325} (1997), 171-176

\bibitem{BHP2} P.F. Baum, N. Higson, R.J. Plymen,
``Representation theory of $p$-adic groups:
a view from operator algebras'',
Proc. Sympos. Pure. Math. \textbf{68} (2000), 111-149

\bibitem{BaNi} P.F. Baum, V. Nistor,
``Periodic cyclic homology of Iwahori-Hecke algebras'',
K-Theory \textbf{27.4} (2002), 329-357

\bibitem{BeDe} J.N. Bernstein, P. Deligne,
``Le "centre" de Bernstein'', 
pp. 1-32 in: \emph{Repr\'esentations des groupes 
r\'eductifs sur un corps local},
Travaux en cours,
Hermann, Paris, 1984

\bibitem{Bla} B. Blackadar,
\emph{$K$-theory for operator algebras 2nd ed.},
Mathematical Sciences Research Institute Publications \textbf{5},
Cambridge University Press, Cambridge, 1998

\bibitem{BoWa} A. Borel, N.R. Wallach,
\emph{Continuous cohomology, discrete subgroups,
and representations of reductive groups},
Annals of Mathematics Studies \textbf{94},
Princeton University Press, Princeton NJ, 1980

\bibitem{Bos} J.-B. Bost,
``Principe d'Oka, $K$-th\'eorie et syst\`emes dynamiques 
non commutatifs'',
Invent. Math. \textbf{101} (1990), 261-333

\bibitem{Bou} N. Bourbaki,
\emph{Groupes et alg\`ebres de Lie. Chapitres IV, V et VI},
Hermann, Paris, 1968

\bibitem{Bre} G.E. Bredon,
\emph{Equivariant cohomology theories},
Lecture Notes in Mathematics \textbf{34},
Springer-Verlag, Heidelberg, 1967

\bibitem{BrLy} J. Brodzki, Z.A. Lykova,
``Excision in cyclic type homology of Fr\'echet algebras'',
Bull. London Math. Soc. \textbf{33} (2000), 283-291

\bibitem{BrPl1} J. Brodzki, R.J. Plymen,
``Periodic cyclic homology of certain nuclear algebras'',
C.R. Acad. Sci. Paris \textbf{329} (1999), 671-676

\bibitem{BrPl2} J. Brodzki, R.J. Plymen,
``Complex structure on the smooth dual of $GL (n)$'',
Doc. Math. \textbf{7} (2002), 91-112

\bibitem{Bro} L.G. Brown,
``Stable isomorphism of hereditary subalgebras of 
$C^*$-algebras'',
Pacific J. Math. \textbf{71.2} (1977), 335-348

\bibitem{BrTi1} F. Bruhat, J. Tits,		
``Groupes r\'eductifs sur un corps local I.
Donn\'ees radicielles valu\'ees'',
Publ. Math. Inst. Hautes \'Etudes Sci. \textbf{41} (1972), 5-251

\bibitem{BrTi2} F. Bruhat, J. Tits,		
``Groupes r\'eductifs sur un corps local II., Sch\'emas en 
groupes. Existence d'une donn\'ee radicielle valu\'ee'',
Publ. Math. Inst. Hautes \'Etudes Sci. \textbf{60} (1984), 5-184

\bibitem{Bry1} J.L. Brylinski,
\emph{Algebras associated with group actions and their homology},
Brown University preprint, Providence RI, 1987

\bibitem{Bry2} J.L. Brylinski,
``Cyclic homology and equivariant theories'',
Ann. Inst. Fourier \textbf{37.4} (1987), 15-28

\bibitem{BuKu1} C.J. Bushnell, P.C. Kutzko,
\emph{The admissible dual of GL(N) via compact open subgroups},
Annals of Mathematics Studies \textbf{129},
Princeton University Press, Princeton NJ, 1993

\bibitem{BuKu2} C.J. Bushnell, P.C. Kutzko,
``Smooth representations of reductive $p$-adic groups:
structure theory via types'',
Proc. London Math. Soc. \textbf{77.3} (1998), 582-634

\bibitem{BuKu3} C.J. Bushnell, P.C. Kutzko,  
``Semisimple types in $GL_n$'',
Compositio Math. \textbf{119.1} (1999), 53-97

\bibitem{CaEi} H. Cartan, S. Eilenberg,
\emph{Homological algebra},
Princeton University Press, Princeton NJ, 1956

\bibitem{Car1} R.W. Carter, 
\emph{Finite groups of Lie type.
Conjugacy classes and complex characters},
Pure and Applied Mathematics (New York),
John Wiley \& Sons, New York, 1985

\bibitem{Car2} P. Cartier,
``Representations of $p$-adic groups : a survey'', 
pp. 111-155 in: \emph{Automorphic forms, representations 
and L-functions. Part 1},
Proc. Sympos. Pure Math. \textbf{33},
American Mathematical Society, Providence RI, 1979

\bibitem{Cas} W. Casselman,
``Introduction to the theory of admissible representations
of $p$-adic reductive groups'',
draft, 1995

\bibitem{Che} C.C. Chevalley,
\emph{Classification des groupes de Lie alg\'ebriques},
S\'eminaire Ecole Normale Sup\'erieure 1956-1958,
Secr\'etariat math\'ematique, Paris, 1958

\bibitem{Cli} A.H. Clifford,
``Representations induced in an invariant subgroup'',
Ann. of Math. \textbf{38} (1937), 533-550

\bibitem{Con} A. Connes,
``Noncommutative differential geometry'',
Publ. Math. Inst. Hautes \'Etudes Sci. \textbf{62} (1985), 41-144

\bibitem{Cun1} J.R. Cuntz, 
``Bivariante $K$-Theorie f\"ur lokalkonvexe Algebren und
der Chern-Connes-Charakter'',
Doc. Math \textbf{2} (1997), 139-182

\bibitem{Cun2} J.R. Cuntz,
``Excision in periodic cyclic theory for topological algebras'',
pp. 43-53 in: \emph{Cyclic cohomology and noncommutative geometry},
Fields Inst. Commun. \textbf{17},
American Mathematical Society, Providence RI, 1997

\bibitem{Cun3} J.R. Cuntz,
``Morita invariance in cyclic homology for nonunital algebras'',
K-Theory \textbf{15} (1998), 301-305

\bibitem{CuQu} J.R. Cuntz, D. Quillen,
``Excision in bivariant periodic cyclic cohomology'',	
Inv. Math. \textbf{127} (1997), 67-98

\bibitem{CuRe} C.W. Curtis, I. Reiner,
\emph{Representation theory of finite groups and associative algebras},
Pure and Applied Mathematics \textbf{11},
John Wiley \& Sons, New York - London, 1962

\bibitem{DeGr} M. Demazure, A. Grothendieck,
\emph{Sch\'emas en groupes III. 
Structure des sch\'emas en groupes r\'eductifs},
Lecture Notes in Mathematics \textbf{153},
Springer-Verlag, Berlin - New York, 1964

\bibitem{DeOp1} P. Delorme, E.M. Opdam,
``The Schwartz algebra of an affine Hecke algebra'',
arXiv:math.RT/0312517, 2004

\bibitem{DeOp2} P. Delorme, E.M. Opdam,
``Analytic R-groups of affine Hecke algebras'',
preprint, 2005

\bibitem{Dix} J. Dixmier,
\emph{Les $C^*$-alg\`ebres et leurs representations},
Cahiers Scientifiques \textbf{29},
Gauthier-Villars \'Editeur, Paris, 1969

\bibitem{Dre} A. Dress,
``Newman's theorems on transformation groups'',
Topology \textbf{8} (1969), 203-207

\bibitem{Ell} G.A. Elliott,
``On the $K$-theory of the $C^*$-algebra generated by a projective
representation of a torsion-free discrete abelian group'',
pp. 157-184 in: \emph{Operator algebras and group representations 
Vol. I}, Monogr. Stud. Math. \textbf{17},
Pitman, Boston MA, 1984

\bibitem{Emm} I. Emmanouil,
``The K\"unneth formula in periodic cyclic homology'',
K-Theory \textbf{10.2} (1996), 197-214

\bibitem{Eve} S. Evens,
``The Langlands classification for graded Hecke algebras'',
Proc. Amer. Math. Soc. \textbf{124.4} (1996), 1285-1290

\bibitem{GePf} M. Geck, G. Pfeiffer,	
\emph{Characters of finite Coxeter groups and Iwahori-Hecke algebras},
London Mathematical Society Monographs, New Series \textbf{21},
Oxford University Press, New York, 2000

\bibitem{God} R. Godement,
\emph{Topologie alg\`ebrique et th\'eorie des faisceaux},
Publ. Math. Univ. Strasbourg \textbf{13},
Hermann, Paris, 1958

\bibitem{GoRo} D. Goldberg, A. Roche, 
``Hecke algebras and $SL_n$-types'',
Proc. London Math. Soc. \textbf{90.1} (2005), 87-131

\bibitem{Goo} T.G. Goodwillie,
``Cyclic homology, derivations, and the free loopspace'',
Topology \textbf{24.2} (1985), 187-215

\bibitem{Gro1} A. Grothendieck,
\emph{Produits tensoriels topologiques et espaces nucl\'eaires},
Mem. Amer. Math. Soc. \textbf{16},
American Mathematical Society, Providence RI, 1955

\bibitem{Gro2} A. Grothendieck,
``Sur quelques points d'alg\`ebre homologique'',
T\^ohoku Math. J. II \textbf{9} (1957), 119-221

\bibitem{Gro3} A. Grothendieck,
``On the De Rham cohomology of algebraic varieties'',
Publ. Math. Inst. Hautes \'Etudes Sci. \textbf{29} (1966), 95-103

\bibitem{Gyo} A. Gyoja,
``Modular representation theory over a ring of higher dimension, 
with applications to Hecke algebras'',
J. Algebra \textbf{174.2} (1995), 553-572

\bibitem{GyUn} A. Gyoja, K. Uno,
``On the semisimplicity of Hecke algebras'',
J. Math. Soc. Japan \textbf{41.1} (1989), 75-79

\bibitem{HC1} Harish-Chandra,
''Harmonic analysis on reductive $p$-adic groups'',
pp. 167-192 in: \emph{Harmonic analysis on homogeneous spaces},
Proc. Sympos. Pure Math. \textbf{26},
American Mathematical Society, Providence RI, 1973

\bibitem{HC2} Harish-Chandra,
``The Plancherel formula for reductive $p$-adic groups'',
pp. 353-367 in: \emph{Collected papers Vol. IV},
Springer-Verlag, New York, 1984

\bibitem{Hart} R. Hartshorne,
``On the De Rham cohomology of algebraic varieties'',
Publ. Math. Inst. Hautes \'Etudes Sci. \textbf{45} (1975), 5-99

\bibitem{HeOp} G.J. Heckman, E.M. Opdam,
``Yang's system of particles and Hecke algebras'',
Ann. of Math. \textbf{145.1} (1997), 139-173

\bibitem{HiNi} N. Higson, V. Nistor,
``Cyclic homology of totally disconnected
groups acting on buildings'',
J. Funct. Anal. \textbf{141.2} (1996), 466-495

\bibitem{HKR} G. Hochschild, B. Kostant, A. Rosenberg,	
``Differential forms on regular affine algebras'',
Trans. Amer. Math. Soc. \textbf{102.3} (1962), 383-408

\bibitem{Hum} J.E. Humphreys,
\emph{Reflection groups and Coxeter groups},
Cambridge Studies in Advanced Mathematics \textbf{29},
Cambridge University Press, Cambridge, 1990

\bibitem{Ill} S. Illman,
``Smooth equivariant triangulations of $G$-manifolds 
for $G$ a finite group'', 
Math. Ann. \textbf{233} (1978), 199-220

\bibitem{Iwa1} N. Iwahori,
``On the structure of a Hecke ring of a Chevalley
group over a finite field'',
J. Fac. Sci. Univ. Tokyo Sect. I \textbf{10} (1964), 215-236

\bibitem{Iwa2} N. Iwahori,	
``Generalized Tits system (Bruhat decompostition) on $p$-adic 
semisimple groups'', 
pp. 71-83 in: \emph{Algebraic groups and discontinuous subgroups},
Proc. Sympos. Pure Math. \textbf{9},
American Mathematical Society, Providence RI, 1966

\bibitem{IwMa} N. Iwahori, H. Matsumoto,
``On some Bruhat decomposition and the structure,
of the Hecke rings of the $p$-adic Chevalley groups'',
Inst. Hautes \'Etudes Sci. Publ. Math \textbf{25} (1965), 5-48

\bibitem{Joh} B.E. Johnson,
\emph{Cohomology in Banach algebras},
Mem. Amer. Math. Soc. \textbf{127},
American Mathematical Society, Providence RI, 1972

\bibitem{Jul} P. Julg,
``$K$-th\'eorie \'equivariante et produits crois\'es'',
C.R. Acad. Sci. Paris \textbf{292} (1981), 629-632

\bibitem{Kar1} M. Karoubi,
``Espaces classifiants en $K$-th\'eorie'',
Trans. Amer. Math. Soc. \textbf{147} (1970), 75-115

\bibitem{Kar2} M. Karoubi,
``Homologie cyclique et $K$-th\'eorie'',
Ast\'erisque \textbf{149} (1987), 1-147

\bibitem{Kas} C. Kassel,
``Cyclic homology, comodules, and mixed complexes'',
J. Algebra \textbf{107.1} (1987), 195-216

\bibitem{Kat1} S.-I. Kato,
``Irreducibility of principal series representations
for Hecke algebras of affine type'',
J. Fac. Sci. Univ. Tokyo Sect. IA Math. 
\textbf{28.3} (1981), 929-943

\bibitem{Kat2} S.-I. Kato,
``A realization of irreducible representations
of affine Weyl groups'',
Indag. Math. \textbf{45.2} (1983), 193-201

\bibitem{Kat3} S.-I. Kato,
``Duality for representations of a Hecke algebra'',
Proc. Amer. Math. Soc. \textbf{119.3} (1993), 941-946

\bibitem{Kat4} Syu Kato,
``An exotic Deligne-Langlands correspondence for symplectic groups'',
arXiv:math.RT/0601155, 2006

\bibitem{Kat5} Syu Kato,
``On the geometry of exotic nilpotent cones'',
arXiv:math.RT/0607478, 2006

\bibitem{KaLu} D. Kazhdan, G. Lusztig,
``Proof of the Deligne-Langlands conjecture for Hecke algebras'',
Invent. Math. \textbf{87} (1987), 153-215

\bibitem{KNS} D. Kazhdan, V. Nistor, P. Schneider,
``Hochschild and cyclic homology of finite type algebras'',
Sel. Math. New Ser. \textbf{4.2} (1998), 321-359

\bibitem{KnVo} A.W. Knapp, D.A. Vogan,
\emph{Cohomological induction and unitary representations},
Princeton Mathematical Series \textbf{45},
Princeton University Press, Princeton NJ, 1995

\bibitem{Laf} V. Lafforgue,
``$K$-th\'eorie bivariante pour les alg\`ebres de Banach
et conjecture de Baum-Connes'',
Invent. Math. \textbf{149.1} (2002), 1-95

\bibitem{Lan} R.P. Langlands,
``On the classification of irreducible representations
of real algebraic groups'',
pp. 101-170 in: \emph{Representation theory and harmonic 
analysis on semisimple Lie groups},
Math. Surveys Monogr. \textbf{31},
American Mathematical Society, Providence RI, 1989

\bibitem{Lod} J.-L. Loday,
\emph{Cyclic homology 2nd ed.},
Mathematischen Wissenschaften \textbf{301},
Springer-Verlag, Berlin, 1997

\bibitem{LoQu} J.-L. Loday, D. Quillen,
``Cyclic homology and the Lie algebra homology of matrices'',
Comment. Math. Helvetici \textbf{59} (1984), 565-591

\bibitem{Lus1} G. Lusztig,
``Cells in affine Weyl groups'',
pp. 255-267 in: \emph{Algebraic groups and related topics},
Adv. Stud. Pure Math. \textbf{6},
North Holland, Amsterdam, 1985

\bibitem{Lus2} G. Lusztig,
``Cells in affine Weyl groups II'',
J. Algebra \textbf{109} (1987), 536-548

\bibitem{Lus3} G. Lusztig,
``Cells in affine Weyl groups III'',
J. Fac. Sci. Univ. Tokyo \textbf{34.2} (1987), 223-243

\bibitem{Lus4} G. Lusztig,
``Affine Hecke algebras and their graded version'',
J. Amer. Math. Soc \textbf{2.3} (1989), 599-635

\bibitem{Lus5} G. Lusztig,
``Classification of unipotent representations of simple 
$p$-adic groups'',
Int. Math. Res. Notices \textbf{11} (1995), 517-589

\bibitem{Lus6} G. Lusztig,
\emph{Hecke algebras with unequal parameters},
CRM Monograph Series \textbf{18},
American Mathematical Society, Providence RI, 2003

\bibitem{Mat2} H. Matsumoto,
``G\'en\'erateurs et relations des groupes de Weyl 
g\'en\'eralis\'es'',
C.R. Acad. Sci. Paris \textbf{258} (1964), 3419-3422

\bibitem{Mat1} H. Matsumoto,
\emph{Analyse harmonique dans les syst\`emes de Tits 
bornologiques de type affine},
Lecture Notes in Mathematics \textbf{590},
Springer-Verlag, Berlin - New York, 1977

\bibitem{Mey} R. Meyer,
``Homological algebra for Schwartz algebras of
reductive $p$-adic groups'',
pp. 263-300 in: \emph{Noncommutative geometry and number theory},
Aspects of Mathematics \textbf{E37},
Vieweg Verlag, Wiesbaden, 2006

\bibitem{Mic} E.A. Michael,
\emph{Locally multiplicatively-convex topological algebras},
Mem. Amer. Math. Soc. \textbf{11},
American Mathematival Society, Providence RI, 1952

\bibitem{Mis} P.A. Mischenko,
\emph{Invariant tempered distributions on the reductive $p$-adic
group $GL_n (F_p )$},
C.R. Math. Rep. Acad. Sci. Canada \textbf{4.2} (1982), 123-127

\bibitem{Mor} L. Morris,
``Tamely ramified intertwining algebras'',
Invent. Math. \textbf{114.1} (1993), 1-54

\bibitem{Nis2} V. Nistor,
``Higher index theorems and the boundary map
in cyclic cohomology'',
Doc. Math. \textbf{2} (1997), 263-295

\bibitem{Nis4} V. Nistor,
``A non-commutative geometry approach to the representation 
theory of reductive $p$-adic groups: Homology of Hecke 
algebras, a survey and some new results'',
pp. 301-323 in: \emph{Noncommutative geometry and number theory},
Aspects of Mathematics \textbf{E37},
Vieweg Verlag, Wiesbaden, 2006

\bibitem{Opd2} E.M. Opdam,
``A generating function for the trace of 
the Iwahori-Hecke algebra'',
Progr. Math. \textbf{210} (2003), 301-323

\bibitem{Opd3} E.M. Opdam,
``On the spectral decomposition of affine Hecke algebras'',
J. Inst. Math. Jussieu \textbf{3.4} (2004), 531-648

\bibitem{Opd4} E.M. Opdam,
``Hecke algebras and harmonic analysis'',
pp. 1227-1259 in: \emph{Proceedings of the International Congress
of Mathematicians - Madrid, August 22-30, 2006. Vol. II},
European Mathematical Society Publishing House, 2006

\bibitem{OsTe} H. Osaka, T. Teruya,
``Topological stable rank of inclusions of unital 
$C^*$-algebras'',
arXiv:math.OA/0311461, 2003

\bibitem{Phi1} N.C. Phillips,
\emph{Equivariant $K$-theory and freeness 
of group actions on $C^*$-algebras},
Lecture Notes in Mathematics \textbf{1274},
Springer-Verlag, Berlin, 1987

\bibitem{Phi2} N.C. Phillips,
``$K$-theory for Fr\'echet algebras'',
Int. J. Math. \textbf{2.1} (1991), 77-129

\bibitem{Ply1} R.J. Plymen,
``The reduced $C^*$-algebra of the $p$-adic group $GL (n)$'',
J. Funct. Anal. \textbf{72} (1987), 1-12

\bibitem{Ply2} R.J. Plymen,
``Reduced $C^*$-algebra for reductive $p$-adic groups'',
J. Funct. Anal. \textbf{88.2} (1990), 251-266

\bibitem{Ree} M. Reeder,
``Isogenies of Hecke algebras and a Langlands correspondence
for ramified principal series representations'',
Representation Theory \textbf{6} (2002), 101-126

\bibitem{Rie1} M.A. Rieffel,
``Dimension and stable rank in the $K$-theory of 
$C^*$-algebras'',
Proc. London Math. Soc. \textbf{46.2} (1983), 301-333

\bibitem{Rie2} M.A. Rieffel,
``The homotopy groups of the unitary groups of noncommutative tori'',
J. Operator Theory \textbf{17} (1987), 237-254

\bibitem{Roc} A. Roche,
``Types and Hecke algebras for principal series representations, 
of split reductive $p$-adic groups'',
Ann. Sci. \'Ecole Norm. Sup. \textbf{31.3} (1998), 361-413

\bibitem{Ros} J.M. Rosenberg,
``Appendix to: "Crossed products of UHF algebras by 
product type actions" by O. Bratteli'',
Duke Math. J. \textbf{46.1} (1979), 25-26

\bibitem{Sat} I. Satake,
``On a generalization of the notion of manifold'',
Proc. Nat. Acad. Sci. U.S.A. \textbf{42} (1956), 359-363

\bibitem{Schn} P. Schneider,
``The cyclic homology of $p$-adic reductive groups'',
J. f\"ur reine angew. Math. \textbf{475} (1996), 39-54

\bibitem{ScSt} P. Schneider, U. Stuhler,	
``Representation theory and sheaves on the Bruhat-Tits building'',
Publ. Math. Inst. Hautes \'Etudes Sci. \textbf{85} (1997), 97-191 

\bibitem{Schu} F. Schur,
``\"Uber die Darstellung der endlichen Gruppen
durch gebrochene linaere Subsitutionen'',
J. f\"ur reine angew. Math. \textbf{127} (1904), 20-50

\bibitem{Seg1} G.B. Segal,
``Classifying spaces and spectral sequences'',
Publ. Math. Inst. Hautes \'Etudes Sci. \textbf{34} (1968), 105-112

\bibitem{Seg2} G.B. Segal, 
``Equivariant $K$-theory'',
Publ. Math. Inst. Hautes \'Etudes Sci. \textbf{34} (1968), 129-151

\bibitem{Ser} J.-P. Serre, 
``G\'eom\'etrie alg\`ebrique et g\'eom\'etrie analytique'',
Ann. Inst. Fourier \textbf{6} (1955), 1-42

\bibitem{Shi} G. Shimura,
``Sur les int\'egrales attach\'ees aux formes automorphes'',
J. Math. Soc. Japan \textbf{11} (1959), 291-311

\bibitem{Sil1} A.J. Silberger,
\emph{Introduction to harmonic analysis on reductive $p$-adic groups},
Mathematical Notes \textbf{23},
Princeton University Press, Princeton NJ, 1979

\bibitem{Sil2} A.J. Silberger,
``The Langlands quotient theorem for p-adic groups'',
Math. Ann. \textbf{236.2} (1978), 95-104

\bibitem{Slo1} J. S\l omi\'nska,
``On the equivariant Chern homomorphism'',
Bull. Acad. Polon. Sci. Ser. Sci. Math. Astr. Phys. 
\textbf{24.10} (1976), 909-913

\bibitem{Slo2} K. Slooten,
\emph{A combinatorial generalization of the Springer 
correspondence for classical type},
Ph.D. Thesis, Universiteit van Amsterdam, 2003

\bibitem{Slo3} K. Slooten,
``Reducibility of induced discrete series representations
for affine Hecke algebras of type B'',
arXiv:math.RT/0511206, 2005

\bibitem{Sol} M.S. Solleveld, 
``Some Fr\'echet algebras for which the Chern character 
is an isomorphism'', 
$K$-theory \textbf{36} (2005), 275-290

\bibitem{Spr2} T.A. Springer,
\emph{Linear algebraic groups 2nd ed.},
Progress in Mathematics \textbf{9},
Birkh\"auser, Boston MA, 1998

\bibitem{Tak} M. Takesaki,
\emph{Theory of operator algebras I},
Springer-Verlag, New York, 1979

\bibitem{Tel} N. Teleman, 
``Microlocalisation de l'homologie de Hochschild'',
C.R. Acad. Sci. Paris \textbf{326} (1998), 1261-1264

\bibitem{Tit} J. Tits,
``Reductive groups over local fields'', pp. 29-69 in: 
\emph{Automorphic forms, representations and L-functions. Part 1},
Proc. Sympos. Pure Math. \textbf{33},
American Mathematical Society, Providence RI, 1979

\bibitem{Tou} J.C. Tougeron, 
\emph{Id\'eaux de fonctions diff\'erentiables},
Ergebnisse der Mathematik und ihrer Grenzgebiete \textbf{71}, 
Springer-Verlag, Berlin, 1972

\bibitem{Tsy} B.L. Tsygan,
``Homology of matrix Lie algebras over rings and the 
Hochschild homology'',
Russian Math. Surveys \textbf{38.2} (1983), 198-199

\bibitem{Vig} M.-F. Vign\'eras,
``On formal dimensions for reductive $p$-adic groups'',
pp. 225-266 in: 
\emph{Festschrift in honor of I. I. Piatetski-Shapiro 
on the occasion of his sixtieth birthday. Part I},
Israel Math. Conf. Proc. \textbf{2},
Weizmann, Jerusalem, 1990

\bibitem{Voi} C. Voigt,
``Chern character for totally disconnected groups'',
arXiv:math.KT/0608626, 2006

\bibitem{Wal} J.-L. Waldspurger,
``La formule de Plancherel pour les groupes $p$-adiques
(d'apr\`es Harish-Chandra)'',
J. Inst. Math. Jussieu \textbf{2.2} (2003), 235-333

\bibitem{Was1} A.J. Wassermann,
``Une d\'emonstration de la conjecture de Connes-Kasparov pour 
les groupes de Lie lin\'eaires connexes r\'eductifs'',
C.R. Acad. Sci. Paris \textbf{304} (1987), 559-562

\bibitem{Was2} A.J. Wassermann, 
``Cyclic cohomology of algebras of smooth functions on orbifolds'',
pp. 229-244 in: \emph{Operator algebras and applications. Vol. I},
London Math. Soc. Lecture Note Ser. \textbf{135},
Cambridge University Press, Cambridge, 1988

\bibitem{Wed} J.H.M. Wedderburn,
``On hypercomplex numbers'',
Proc. London Math. Soc. \textbf{6} (1908), 77-118

\bibitem{Wod} M. Wodzicki, 
``The long exact sequence in cyclic homology associated 
with an extension of algebras'',
C.R. Acad. Sci. Paris \textbf{306} (1988), 399-403

\bibitem{Zel} A.V. Zelevinsky,
``Induced representations of reductive $p$-adic groups II.,
On irreducible representations of $GL (n)$'',
Ann. Sci. \'Ecole Norm. Sup. \textbf{13.2} (1980), 165-210


\end{thebibliography}
